\documentclass[11pt]{flbook12}
\usepackage{fixltx2e}[2006/03/24]
\usepackage{times}
\usepackage{epsfig}
\usepackage{latexsym}
\usepackage{theorem}
\usepackage{amsfonts}
\usepackage{amssymb}
\usepackage{enumerate}
\usepackage{amsmath}
\usepackage{amscd}
\usepackage[all]{xy}
\usepackage{verbatim}
\usepackage{rotating}
\usepackage{isolatin1}
\usepackage{longtable} % for Johan's chapter 4
\usepackage{tikz}
\usetikzlibrary{arrows,automata}
\theoremstyle{change}
\theorembodyfont{\sl}
\newtheorem{thm}[subsection]{Theorem}
\newtheorem{prop}[subsection]{Proposition}
\newtheorem{lem}[subsection]{Lemma}
\newtheorem{cor}[subsection]{Corollary}

\newtheorem{supthm}[section]{Theorem}
\newtheorem{suplem}[section]{Lemma}
\newtheorem{supprop}[section]{Proposition}

\theorembodyfont{\rmfamily}
\newtheorem{defi}[subsection]{Definition}

\newtheorem{rem}[subsection]{Remark}

\newtheorem{suprem}[section]{Remark}

\newtheorem{exam}[subsection]{Example}

\makeatletter
\newenvironment{eqn}{\refstepcounter{subsection}
$$}{\leqno{\rm(\thesubsection)}$$\global\@ignoretrue}
\makeatother  



\makeatletter
\newenvironment{supeqn}{\refstepcounter{section}
$$}{\leqno{\rm(\thesection)}$$\global\@ignoretrue}
\makeatother





\newenvironment{prf}[1]{\trivlist
\item[\hskip \labelsep{\it
#1\hspace*{.3em}}]}{~\hspace{\fill}~$\square$\endtrivlist}

\renewenvironment{proof}{\begin{prf}{\bf Proof}}{\end{prf}}
\newcommand{\ol}{\overline}
\newcommand{\ul}{\underline}
\newcommand{\ZZ}{{\mathbb Z}}
\newcommand{\NN}{{\mathbb N}}
\newcommand{\QQ}{{\mathbb Q}}
\newcommand{\RR}{{\mathbb R}}
\newcommand{\CC}{{\mathbb C}}
\newcommand{\FF}{{\mathbb F}}
\newcommand{\HH}{{\mathbb H}}
\newcommand{\PP}{{\mathbb P}}
\newcommand{\GG}{{\mathbb G}}
\newcommand{\EE}{{\mathbb E}}
\renewcommand{\AA}{{\mathbb A}}
\renewcommand{\SS}{{\mathbb S}}
\newcommand{\TT}{\mathbb{T}}

\newcommand{\MM}{{\mathbb M}}
\newcommand{\BB}{{\mathbb B}}
\newcommand{\Qbar}{{\overline{\QQ}}}
\newcommand{\Zbar}{{\overline{\ZZ}}}
\newcommand{\Fbar}{{\overline{\FF}}}
\newcommand{\mmu}{{\mu}}
\newcommand{\Gal}{{\rm Gal}}
\newcommand{\Pic}{{\rm Pic}}
\DeclareMathOperator{\Spec}{Spec}

\newcommand{\rF}{\mathrm{F}}
\newcommand{\rH}{\mathrm{H}}

\newcommand{\rM}{\mathrm{M}}
\newcommand{\rN}{\mathrm{N}}

\newcommand{\rR}{\mathrm{R}}
\newcommand{\rS}{\mathrm{S}}
\newcommand{\eps}{\varepsilon}

\newcommand{\calO}{{\cal O}}
\newcommand{\calE}{{\cal E}}
\newcommand{\calL}{{\cal L}}

\newcommand{\calF}{{\cal F}}

\newcommand{\calC}{{\cal C}}
\newcommand{\calX}{{\cal X}}

\newcommand{\calV}{{\cal V}}
\newcommand{\calB}{{\cal B}}

\newcommand{\calA}{{\cal A}}

\newcommand{\divisor}{\mathrm{div}}
\newcommand{\fin}{\mathrm{fin}}
\newcommand{\Div}{\mathrm{Div}}

\newcommand{\GL}{{\rm GL}}
\newcommand{\PGL}{{\rm PGL}}
\newcommand{\SL}{{\rm SL}}
\newcommand{\om}{\underline{\omega}}
\newcommand{\Cusps}{{\rm Cusps}}
\newcommand{\lto}{\longrightarrow}
\newcommand{\into}{{\hookrightarrow}}
\newcommand{\onto}{\twoheadrightarrow}

\newcommand{\isomlto}{\;\tilde{\longrightarrow}\;}
\newcommand{\Hom}{{\rm Hom}}
\newcommand{\Isom}{{\rm Isom}}

\newcommand{\End}{{\rm End}}
\newcommand{\Aut}{{\rm Aut}}

\newcommand{\id}{{\rm id}}
\newcommand{\im}{{\rm im}}
\newcommand{\Sym}{{\rm Sym}}
\newcommand{\SU}{{\rm SU}}

\newcommand{\trace}{{\rm trace}}
\newcommand\tr{\mathop{\rm tr} \nolimits} % for Johan
\newcommand{\ld}{\langle}
\newcommand{\rd}{\rangle}
\newcommand{\Gm}{{\GG_{\rm m}}}

\newcommand{\abs}{{\rm abs}}

\newcommand{\Zhat}{\hat{\ZZ}}

\newcommand{\ab}{{\rm ab}}
\newcommand{\Supp}{{\rm Supp}}

\newcommand{\discr}{{\rm discr}}
\newcommand\disc{\mathop{\rm Disc} \nolimits} % for Johan
\newcommand{\rank}{{\rm rank}}

\newcommand{\et}{\mathrm{et}}
\newcommand{\sm}{\mathrm{sm}}

\DeclareMathOperator{\Vol}{Vol}

\DeclareMathOperator{\Cot}{Cot}

\newcommand{\Frob}{\mathrm{Frob}}
\newcommand{\new}{\mathrm{new}}

\newcommand{\compact}{\mathrm{c}}
\newcommand{\unr}{\mathrm{unr}}
\newcommand{\cusp}{\mathrm{cusp}}
\DeclareMathOperator{\li}{li}
\DeclareMathOperator{\Log}{Log}
\let\set\mathbb
\def\<<{\leavevmode
  \raise0.28ex\hbox{$\scriptscriptstyle\langle\!\langle$}\nobreak
  \hskip -.6pt plus.3pt minus.2pt\,}
\def\>>{\,\nobreak\hskip -.6pt plus.3pt minus.2pt
  \raise0.28ex\hbox{$\scriptscriptstyle\rangle\!\rangle$}}

\DeclareMathOperator{\Ker}{Ker}

\def\bbx{{\bf x}}
\def\bby{{\bf y}}
\def\bbz{{\bf z}}

\def\talpha{{\tilde \alpha}}
\def\tbeta{{\tilde \beta}}

\def\bD{{\bar D}}
\def\BOX{{\text{BOX}}}
\def\FAIL{{\text{FAIL}}}
\def\ZERO{{\text{ZERO}}}
\def\SQRT{{\text{SQRT}}}

\def\FQ{{\FF_Q}}
\def\Fp{{\FF_p}}

\def\Fq{{\FF_q}}
\def\Fqs{{\FF^*_q}}

\def\IM{\mathop{\rm{Im}}\nolimits }
\def\bff{{\bf f}}

\def\bFF{{\bf F}}
\def\bJ{{\bf J}}
\def\bz{{\rm z}}
\def\bv{{\bf v}}
\def\bx{{\rm x}}
\def\bbb{{\rm b}}
\def\bB{{\|}}
\def\by{{\rm y}}
\def\bzero{{{\bf 0}_g}}
\def\bzerogd{{{\bf 0}_{g_2}}}
\def\br{{\bf r}}
\def\bk{{\bf k}}
\def\bn{{\bf n}}
\def\bm{{\bf m}}
\def\ba{{\bf a}}
\def\bb{{\bf b}}
\def\bc{{\bf c}}
\def\bdeux{{{\bf 2}_g}}

\def\bun{{{\bf 1}_g}}
\def\bdeux{{{\bf 2}_g}}
\def\bungd{{{\bf 1}_{g_2}}}
\def\bq{{\bf q}}
\def\gH{{\mathfrak H}}
\def\gHb{{\mathfrak H}^\ast}
\def\gH{{\HH}}
\def\gHb{{\HH}^\ast}
\def\hgamma{{\vec{ \gamma}  }}
\def\Sym{\mathop{\rm{Sym}}\nolimits }
\def\CC{{\set C}}
\def\AA{{\set A}}
\def\EE{{\set E}}
\def\hT{{\hat T}}
\def\KK{{\set K}}
\def\LL{{\set L}}
\def\Div{{\rm Div}}
\def\End{\mathop{\rm End }}
\def\Ann{\mathop{\rm Ann }}
\def\Aut{\mathop{\rm Aut }}
\def\Fl{{\FF _l}}
\def\Gal{\mathop{\rm{Gal}}\nolimits }
\def\NN{{\set N}}

\def\PP{{\set P}}
\def\PP{{\set P}}

\def\Pic{\mathop{\rm{Pic}}\nolimits }
\def\QQ{{\set Q}}
\def\Qb{{\bar \QQ}}
\def\RR{{\set R}}
\def\Spec{\mathop{\rm{Spec}}\nolimits }
\def\TT{{\set T}}
\def\Id{{\rm{Id}}}
\def\Tr{\mathop{\rm{Tr}}\nolimits }
\def\Taylor{\mathop{\rm{Taylor}}\nolimits }
\def\GL{\mathop{\rm{GL}}\nolimits }
\def\ZZ{{\set Z}}
\def\FF{{\set F}}
\def\cB{{\mathcal B}}
\def\cG{{\mathcal G}}

\def\cA{{\mathcal A}}
\def\cD{{\mathcal D}}
\def\cS{{\mathcal S}}
\def\cF{{\mathcal F}}
\def\cH{{\mathcal H}}
\def\cJ{{\mathcal J}}
\def\cK{{\mathcal K}}

\def\cM{{\mathcal M}}
\def\cN{{\mathcal N}}
\def\cO{{\mathcal O}}
\def\cP{{\mathcal P}}
\def\cQ{{\mathcal Q}}
\def\cR{{\mathcal R}}
\def\cU{{\mathcal U}}
\def\cV{{\mathcal V}}

\def\mW{{\xi}}
\def\hmW{{\psi}}

\def\agot{{\mathfrak  a}}
\def\bgot{{\mathfrak  b}}
\def\KJ{{\mathfrak  j}}
\def\KJJ{{\mathfrak  J}}

\def\dgot{{\mathfrak  d}}

\def\mgot{{\mathfrak  m}}
\def\ngot{{\mathfrak  n}}

\def\pgot{{\mathfrak p}}

\def\gotB{{\mathfrak B}}
\def\gotC{{\mathfrak C}}
\def\gotD{{\mathfrak D}}

\def\vp{{\varphi}}
\def\hT{{\hat T}}

\newcommand{\GQl}{\Gal(\overline{\QQ}_\ell/\QQ_\ell)}
\newcommand{\GQ}{\Gal(\overline{\QQ}/\QQ)}
\newcommand{\GQp}{\Gal(\overline{\QQ}_p/\QQ_p)}

\def\sumprime_#1{\setbox0=\hbox{$\scriptstyle{#1}$}
  \setbox2=\hbox{$\displaystyle{\sum}$}
  \setbox4=\hbox{${}'\mathsurround=0pt$}
  \dimen0=.5\wd0 \advance\dimen0 by-.5\wd2
  \ifdim\dimen0>0pt
    \ifdim\dimen0>\wd4 \kern\wd4 \else\kern\dimen0\fi\fi
  \mathop{{\sum}'}_{\kern-\wd4 #1}}

\makeatletter
\def\eqalign#1{\null\,\vcenter{\openup\jot\m@th
  \ialign{\strut\hfil$\displaystyle{##}$&$\displaystyle{{}##}$\hfil
      \crcr#1\crcr}}\,}
\makeatother

\newcommand\Mat[4]{
  \left(
  \genfrac{}{}{0pt}{}{#1}{#3} 
  \thinspace 
  \genfrac{}{}{0pt}{}{#2}{#4}
  \right)
}

\begin{document}

\frontmatter	
\title{Computational aspects of modular forms and Galois representations}

\shorttitle{Computation of coefficients of modular forms}

\subtitle{how one can compute in polynomial time\\ 
the value of Ramanujan's tau at a prime}

\author{Jean-Marc Couveignes and \\ Bas Edixhoven, editors}
%% Here just the editors

\makehalftitle

\maketitle

\begin{bookepigraph}
Let $\tau\colon\NN\to\ZZ$ be defined by:
\[
\sum_{n\geq 0}\tau(n)q^n =  q\prod_{n\geq1}(1-q^n)^{24}
\quad\text{in $\ZZ[[q]]$}.
\]
Then we have, in $\ZZ/19\ZZ$:
\[
\begin{aligned}
\tau(10^{1000}+1357) & = \pm 4, \\
\tau(10^{1000}+7383) & = \pm 2, \\
\tau(10^{1000}+21567) & = \pm 3, \\
\tau(10^{1000}+27057) & = 0,\\
\tau(10^{1000}+46227) & = 0,\\
\tau(10^{1000}+57867) & = 0,\\
\tau(10^{1000}+64749) & = \pm 7,\\
\tau(10^{1000}+68367) & = 0,\\
\tau(10^{1000}+78199) & = \pm 8,\\
\tau(10^{1000}+128647) & = 0.
\end{aligned}
\]
\epigraphsource{Section~\ref{polytable},
  Lemma~\ref{johan_lem_trace_0}, and \cite{Magma}, \cite{SAGE}
  and~\cite{pari}.}
\end{bookepigraph}

\setcounter{tocdepth}{1}
\tableofcontents

\begin{thepreface}
This is a book about computational aspects of modular forms and the
Galois representations attached to them. The main result is the
following: Galois representations over finite fields attached to
modular forms of level one can, in almost all cases, be computed in
polynomial time in the weight and the size of the finite field. As a
consequence, coefficients of modular forms can be computed fast via
congruences, as in Schoof's algorithm for the number of points of
elliptic curves over finite fields. The most important feature of the
proof of the main result is that exact computations involving systems
of polynomial equations in many variables are avoided by
approximations and height bounds, i.e., bounds for the accuracy that
is necessary to derive exact values from the approximations.

The books authors are the two editors, Jean-Marc Couveignes and Bas
Edixhoven, together with Johan Bosman, Robin de Jong, and Franz
Merkl. Each chapter has its own group of authors.

Chapter~\ref{chap_intro} gives an introduction to the subject, precise
statements of the main results, and places these in a somewhat wider
context. Chapter~\ref{chap_mod_curves} provides the necessary
background concerning modular curves and modular
forms. Chapter~\ref{chap_first_descr} gives a first, informal
description of the algorithms.  These first three chapters should
allow readers without much background in arithmetic geometry to still
get a good idea of what happens in the book, skipping, if necessary,
some parts of Chapter~\ref{chap_mod_curves}.

Chapters~\ref{sec_intro_heights_and_ar} and~\ref{sec_couveignes_ZEROS}
provide the necessary background on heights and Arakelov theory, and
on algorithmic aspects of the computation with a desired accuracy of
the roots of complex polynomials and power series.

Chapters~\ref{chapcomput} and~\ref{chappgltau} are concerned with some
real computations of Galois representations attached to modular forms,
and end with a table dealing with all cases of forms of weight at
most~$22$ and finite fields of characteristic at most~$23$.

The main ingredients for the proof of the main result are established
in Chapters~\ref{chap_descr_X15l}, \ref{sec_appl_Arakelov},
\ref{subsec_merkl}, \ref{chap_bnd_height},
\ref{sec_couveignes_TORSION}, and~\ref{sec_couveignes_modp}. The
topics dealt with are, respectively: construction of suitable divisors
on modular curves, bounding heights using Arakelov theory, bounding
Arakelov invariants of certain modular curves, approximation of
divisors using complex numbers, and using finite fields.

The main result on the computation of Galois representations is proved
in Chapter~\ref{sec_comp_mod_l_rep}, where one finds a detailed
description of the algorithm and a rigorous proof of the complexity
bound.  

Chapter~\ref{chap_comp_coefs} contains the application of the main
result to the computation of coefficients of modular forms.

The Epilogue announces some work on generalisations and applications
that will be completed in the near future, as well as a direction of
further research outside the context of modular forms.

\section*{Acknowledgements}

This book started as a report, written in the context of a contract
(Contrat d'Études 04.42.217) between the University of Leiden and the
French CELAR (Centre Électronique de l'Armement). We thank the CELAR,
and, in particular, David Lubicz and Reynald Lercier, for this
financial support.

From 2005 until 2010, our project was supported financially in the
Netherlands by NWO, the Dutch organisation for scientific research
(Nederlandse organisatie voor Wetenschappelijk Onderzoek) in the form
of a VICI-grant to Bas Edixhoven.

As the contract with the CELAR stipulated, the report was published on
internet (arxiv, homepage), in May 2006, and our aim was to extract a
research article from it. We thank Eyal Goren for his suggestion to
expand it into this book instead.

Bas Edixhoven thanks Ren\'e Schoof for asking the question, in 1995,
that this book answers.

We thank Peter Bruin for carefully reading preliminary versions and
pointing out some errors, and Hendrik Lenstra for some suggestions and
references.

\vfill
\text{}
\section*{Dependencies between the chapters}
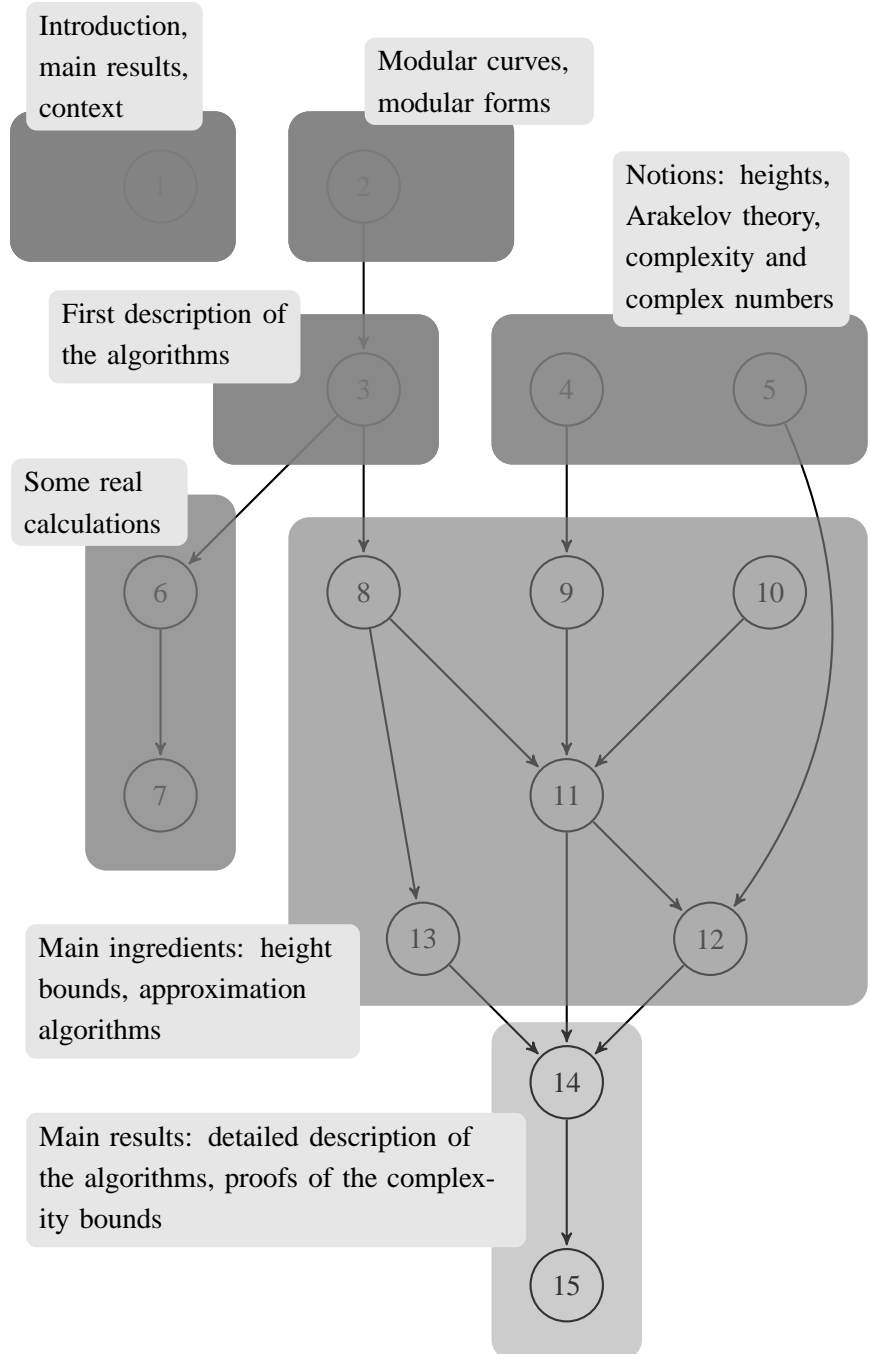
\vfill

\begin{tikzpicture}[auto,node distance=2.7cm,thick]

\tikzstyle{every state}=[fill=none,draw=black,text=black]
\tikzstyle{information text}=[rounded corners,fill=gray!20,inner sep=1ex]

\node[state] (1)                    {1};
\node[state] (2) [right of=1]       {2};
\node[state] (3) [below of=2]       {3};
\node[state] (4) [right of=3]       {4};
\node[state] (5) [right of=4]       {5};
\node[state] (8) [below of=3]       {8};
\node[state] (6) [left  of=8]       {6};
\node[state] (7) [below of=6]       {7};
\node[state] (9) [below of=4]       {9};
\node[state] (10)[right of=9]       {10};
\node[state] (11)[below of=9]       {11};
\node[state] (12)[below right of=11]{12};
\node[state] (13)[below left  of=11]{13};
\node[state] (14)[below left of=12]{14};
\node[state] (15)[below of=14] {15};

\begin{scope}[->,>=stealth',shorten >=1pt,]
\path (2) edge  node {} (3);
\path (3) edge  node {} (8);
\path (3) edge  node {} (6);
\path (4) edge  node {} (9);
\path (6) edge  node {} (7);
\path (8) edge  node {} (11);
\path (9) edge  node {} (11);
\path(10) edge  node {} (11);
\path(11) edge  node {} (12);
\path(11) edge  node {} (14);
\path (8) edge  node {} (13);
\path (5) edge  [bend left] node {} (12);
\path(12) edge  node {} (14);
\path(13) edge  node {} (14);
\path(14) edge  node {} (15);
\end{scope}

\begin{scope}[color=gray, rounded corners=8pt, fill opacity=0.4]
\fill (1) +(-2,-1) rectangle +(1,1);
\fill (2) +(-1,-1) rectangle +(2,1);
\fill (3) +(-2,-1) rectangle +(1,1);
\fill (4) +(-1,-1) rectangle +(4,1);
\fill (7) +(-1,-1) rectangle +(1,4);
\fill (8) +(-1,1) rectangle +(6.7,-5.5);
\fill (15) +(-1,-1) rectangle +(1,3.5);
\end{scope}

\draw
(1) +(-1.8,1.6)node[right,text width=2cm,style=information text]
{Introduction, main results, context};

\draw
(2) +(0,1.4)node[right,text width=2.7cm,style=information text]
{Modular curves, modular forms};

\draw
(3) +(-4.2,0.7)node[right,text width=3cm,style=information text]
{First description of the algorithms};

\draw
(4) +(0.6,2.0)node[right,text width=2.8cm,style=information text]
{Notions: heights, Arakelov theory, complexity and complex numbers};

\draw
(6) +(-2,1.2)node[right,text width=2cm,style=information text]
{Some real calculations};

\draw
(13) +(-5.3,-0.7)node[right,text width=4.1cm,style=information text]
{Main ingredients: height bounds, approximation algorithms};

\draw
(14) +(-7.2,-1.3)node[right,text width=6.3cm,style=information text]
{Main results: detailed description of the algorithms, proofs of
the complexity bounds};

\end{tikzpicture}

\vfill
\end{thepreface}

\mainmatter

\chapter{Introduction, main results, context }\label{chap_intro}

\author{B. Edixhoven}

\bigskip

\bigskip

% Authors: Bas and Jean-Marc.

\section{Statement of the main results}

As the final results in this book are about fast computation of
coefficients of modular forms, we start by describing the state of the
art in this subject.

A convenient way to view modular forms and their coefficients in this
context is as follows, in terms of Hecke algebras. For $N$ and $k$
positive integers, let $S_k(\Gamma_1(N))$ be the finite dimensional
complex vector space of cuspforms of weight~$k$ on the congruence
subgroup $\Gamma_1(N)$ of $\SL_2(\ZZ)$. Each $f$ in $S_k(\Gamma_1(N))$
has a power series expansion $f=\sum_{n\geq1}a_n(f)q^n$, a complex
power series converging on the unit disk. These $a_n(f)$ are the
coefficients of~$f$ that we want to compute, in particular for
large~$n$. For each positive integer~$n$ we have an endomorphism $T_n$
of $S_k(\Gamma_1(N))$, and we let $\TT(N,k)$ denote the
sub-$\ZZ$-algebra of $\End(S_k(\Gamma_1(N))$ generated by them. The
$\TT(N,k)$ are commutative, and free $\ZZ$-modules of rank the
dimension of $S_k(\Gamma_1(N))$, which is of polynomially bounded
growth in~$N$ and~$k$. For each $N$, $k$ and $n$ one has the identity
$a_n(f)=a_1(T_nf)$. The $\CC$-valued pairing between
$S_k(\Gamma_1(N),\CC)$ and $\TT(N,k)$ given by $(f,t)\mapsto a_1(tf)$
identifies $S_k(\Gamma_1(N),\CC)$ with the space of $\ZZ$-linear maps
from $\TT(N,k)$ to~$\CC$, and we can write $f(T_n)$ for~$a_n(f)$. All
together this means that the key to the computation of coefficients of
modular forms is the computation of the Hecke algebras $\TT(N,k)$ and
their elements~$T_n$. A modular form $f$ in $S_k(\Gamma_1(N))$ is
determined by the $f(T_i)$ with $i\leq
k{\cdot}[\SL_2(\ZZ):\Gamma_1(N)]/12$, hence if $T_n$ is known as a
$\ZZ$-linear combination of these $T_i$, then $f(T_n)$ can be computed
as the same $\ZZ$-linear combination of the~$f(T_i)$.

The state of the art in computing the algebras $\TT(N,k)$ can now be
summarised as follows.
\begin{quote}
\emph{
There is a deterministic algorithm, that on input positive
integers~$N$ and ~$k\geq2$, computes $\TT(N,k)$: it gives a
$\ZZ$-basis and the multiplication table for this basis, in running
time polynomial in $N$ and~$k$. Moreover, the Hecke operator $T_n$ can
be expressed in this $\ZZ$-basis in deterministic polynomial time
in~$N$, $k$ and~$n$.  
}
\end{quote}
We do not know a precise reference for this statement, but it is
rather obvious from the literature on calculations with modular forms
for which we refer to William Stein's book~\cite{Stein}, and in
particular to Section~8.10.2 of it. The algorithms alluded to above
use that $S_k(\Gamma_1(N))$, viewed as $\RR$-vector space, is
naturally isomorphic to the $\RR$-vector space obtained from the
so-called ``cuspidal subspace''
$\rH^1(\Gamma_1(N),\ZZ[x,y]_{k-2})_\cusp$ of the $\ZZ$-module
$\rH^1(\Gamma_1(N),\ZZ[x,y]_{k-2})$ in group cohomology. Here,
$\ZZ[x,y]_{k-2}$ is the homogeneous part of degree $k{-}2$ of the
polynomial ring $\ZZ[x,y]$ on which $\SL_2(\ZZ)$ acts via its standard
representation on $\ZZ[x,y]_1$. In this way,
$\rH^1(\Gamma_1(N),\ZZ[x,y]_{k-2})_\cusp$, modulo its torsion
subgroup, is a free $\ZZ$-module of finite rank that is a faithful
$\TT(N,k)$-module, and the action of the $T_n$ is described
explicitly. Algorithms based on this typically use a presentation of
$\rH^1(\Gamma_1(N),\ZZ[x,y]_{k-2})_\cusp$ in terms of so-called
``modular symbols'', and we call them therefore modular symbols
algorithms. The theory of modular symbols was developed by Birch,
Manin, Shokurov, Cremona, Merel,\ldots. It has led to many algorithms,
implementations and calculations, which together form the point of
departure for this book.

The computation of the element $T_n$ of $\TT(N,k)$, using modular
symbols algorithms, involves sums of a number of terms that grows at
least linearly in~$n$. If one computes such sums by evaluating and
adding the terms one by one, the computation of $T_n$, for $N$ and $k$
fixed, will take time at least linear in~$n$, and hence exponential in
$\log n$. The same is true for other methods for computing $T_n$ that
we know of: computations with $q$-expansions that involve
multiplication of power series, using linear combinations of theta
series, the ``graph method'' of Mestre and Oesterl\'e, and the
Lefschetz trace formula for correspondences, holomorphic or
not. Efforts to evaluate the encountered sums more quickly seem to
lead, in each case, again to the problem of computing coefficients of
modular forms. For example, the graph method leads to the problem of
computing quickly representation numbers of integer quadratic forms in
4 variables. In the case of the trace formula, there are maybe only
$O(\sqrt n)$ terms, but they contain class numbers of imaginary
quadratic orders, these numbers being themselves directly related to
coefficients of modular forms of half integral weight.

Let us now state one of the main results in this book,
Theorem~\ref{thm_comp_tn}.
\begin{quote}
\emph{
Assume that the generalised Riemann hypothesis (GRH) holds. There
exists a deterministic algorithm that on input positive integers~$n$
and~$k$, together with the factorisation of~$n$ into prime factors,
computes the element~$T_n$ of~$\TT(1,k)$ in running time polynomial
in $k$ and\/~$\log n$.
}
\end{quote}
The restriction to modular forms of level~$1$ in this result is there
for a technical reason. The result will certainly be generalised to
much more general levels; see the Epilogue at the end of this
book. The condition that the factorisation of $n$ into primes must be
part of the input is necessary because we do not have a polynomial
time algorithm for factoring integers. Vice versa,
see~Remark~\ref{rem_bach-charles} for evidence that factoring is not
harder than computing coefficients of modular forms.

Let us describe how the computation of Galois representations is used
for the computation of~$T_n$. Standard identities express $T_n$ in
terms of the $T_p$ for $p$ dividing~$n$. These $T_p$ are computed, via
the LLL basis reduction algorithm, from sufficiently many of their
images under morphisms $f$ from $\TT(1,k)$ to finite fields,
analogously to Schoof's algorithm for counting points of an elliptic
curve over a finite field. Indeed, for such an
$f\colon\TT(1,k)\to\FF$, with $p$ not the characteristic, $l$, say,
of~$\FF$, the image $f(T_p)$ is equal to the trace of
$\rho_f(\Frob_p)$, where $\rho_f\colon \Gal(\Qbar/\QQ)\to\GL_2(\FF)$
is the Galois representation attached to~$f$, and $\rho_f(\Frob_p)$ a
Frobenius element at~$p$. The representation $\rho_f$ is characterised
by the following three conditions: it is semi-simple, it is unramified
outside~$l$, and for all prime numbers $p\neq l$ one has:
\[
\trace(\rho_f(\Frob_p))=f(T_p)\quad\text{and}\quad
\det(\rho_f(\Frob_p))=p^{k-1}\quad \text{in\/ $\FF$.}
\]
It is \emph{the} main result of this book,
Theorem~\ref{thm_comp_rep_mod_l}, plus some standard computational
number theory, that enables us to compute $\rho_f(\Frob_p)$ in time
polynomial in $k$, $\#\FF$ and~$\log p$ (note the $\log$!). Under GRH,
existence of sufficiently many maximal ideals of small enough index is
guaranteed. We partly quote Theorem~\ref{thm_comp_rep_mod_l}.
\begin{quote}
\emph{ There is a deterministic algorithm that on input a positive
  integer~$k$, a finite field~$\FF$, and a surjective ring morphism
  $f$ from $\TT(1,k)$ to $\FF$ such that the associated Galois
  representation $\rho_f\colon \Gal(\Qbar/\QQ)\to\GL_2(\FF)$ is
  reducible or has image containing $\SL_2(\FF)$, computes $\rho_f$ in
  time polynomial in $k$ and\/~$\#\FF$.}
\end{quote}
By ``computing $\rho_f$'' we mean the following. Let $K_f\subset\Qbar$
be the finite Galois extension such that $\rho_f$ factors as the
natural surjection from $\Gal(\Qbar/\QQ)$ to $\Gal(K_f/\QQ)$, followed
by an injection into $\GL_2(\FF)$. Then to give $\rho_f$ means to give
$K_f$ as $\QQ$-algebra, in terms of a multiplication table with
respect to a $\QQ$-basis, together with a list of all elements of
$\Gal(K_f/\QQ)$, as matrices with coefficients in~$\QQ$, and, for each
$\sigma$ in $\Gal(K_f/\QQ)$, to give the corresponding element
$\rho_f(\sigma)$ of~$\GL_2(\FF)$.

Before we describe in more detail, in the next sections, some history
and context concerning our main results, we give one example and we
make some brief remarks. Many of these remarks are treated with more
detail further on.

The first non-trivial example is given by $k=12$. The space of
cuspidal modular forms of level one and weight $12$ is
one-dimensional, generated by the discriminant modular form $\Delta$,
whose coefficients are given by Ramanujan's $\tau$-function:
\[
\Delta = q\prod_{n\geq1}(1-q^n)^{24} = \sum_{n\geq 1}\tau(n)q^n 
= q-24q^2+252q^3+\cdots \quad\text{in $\ZZ[[q]]$}.
\]
In this case, the Hecke algebra $\TT(1,12)$ is the ring $\ZZ$, and,
for each $n$ in $\ZZ_{>0}$, we have $T_n=\tau(n)$. The results above
mean that:
\begin{quote}
\emph{for $p$ prime, Ramanujan's $\tau(p)$ can be computed in time
  polynomial in $\log p$.}
\end{quote}
For $l$ prime, let $\rho_l$ denote the Galois representation to
$\GL_2(\FF_l)$ attached to~$\Delta$. It was proved by Swinnerton-Dyer
that for $l$ not in $\{2,3,5,7,23,691\}$ the image of $\rho_l$
contains $\SL_2(\FF_l)$. This means that for all $l$ not in this short
list the representation $\rho_l$ has non-solvable image, and so cannot
be computed using computational class field theory. The classical
congruences for Ramanujan's $\tau$-function correspond to the $l$ in
the list above. Our results provide a generalisation of these
congruences in the sense that the number fields $K_l$ that give the
$\rho_l$ ``encode'' the $\tau(p)$ mod~$l$ in such a way that $\tau(p)$
mod~$l$ can be computed in time polynomial in $l$ and $\log p$, i.e.,
just the same complexity as in the case where one has explicit
congruences. 

More generally, we hope that non-solvable global field extensions
whose existence and local properties are implied by the Langlands
program can be made accessible to computation and so become even more
useful members of the society of mathematical objects. Explicit
descriptions of these fields make the study of global properties such
as class groups and groups of units possible. Certainly, if we only
knew the maximal abelian extension of~$\QQ$ as described by general
class field theory, then roots of unity would be very much welcomed.

The natural habitat for Galois representations such as the $\rho_f$
above is that of higher degree {\'e}tale cohomology with
$\FF_\ell$-coefficients of algebraic varieties over~$\Qbar$, together
with the action of~$\Gal(\Qbar/\QQ)$. Our results provide some
evidence that, also in interesting cases, such objects can be computed
in reasonable time. We stress that this question is not restricted to
varieties related to modular forms or automorphic forms. In fact,
thinking of elliptic curves, over~$\QQ$, say, knowing that these are
modular does not help for computing their number of points over finite
fields: Schoof's algorithm uses algebraic geometry, not modularity.

The problem of computing {\'e}tale cohomology with Galois action is
clearly related to the question of the existence of polynomial time
algorithms for computing the number of solutions in~$\FF_p$ of a fixed
system of polynomial equations over~$\ZZ$, when $p$ varies. Our
results treat this problem for the $11$-dimensional variety that gives
rise to~$\Delta$; see Section~\ref{sec_hist_cong_tau} for more details
and also for an explicit variety of dimension $19$ related to this.

The Epilogue at the end of this book describes a striking application
of a generalisation of our results to the problem of computing
representation numbers of the $\ZZ^{2k}$ equipped with the standard
inner product. This again is an example where only for small $k$ there
are explicit formulas, but where in general there (surely) exists an
algorithm that computes such numbers as quickly as if such formulas
did exist. Hence, from a computational perspective, such algorithms
form a natural generalisation of the finite series of formulas.

We very briefly describe the method by which we compute
the~$\rho_f$. Their duals occur in the higher degree {\'e}tale
cohomology of certain higher dimensional varieties, but no-one seems
to know how to compute with this directly. 

Via some standard methods in {\'e}tale cohomology (the Leray spectral
sequence, and passing to a finite cover to trivialise a locally
constant sheaf of finite dimensional $\FF_l$-vector spaces), or from
the theory of congruences between modular forms, it is well known that
the $\rho_f$ are realised by subspaces $V_f$ in the $l$-torsion
$J_l(\Qbar)[l]$ of the Jacobian variety $J_l$ of some modular
curve~$X_l$ defined over~$\QQ$. The field $K_l$ is then the field
generated by suitable ``coordinates'' of the points $x\in V_l\subset
J_l(\Qbar)[l]$. We are now in the more familiar situation of torsion
points on abelian varieties. But the price that we have paid for this
is that the abelian variety $J_l$ depends on~$l$, and that its
dimension, equal to the genus of~$X_l$, i.e., equal to
$(l-5)(l-7)/24$, grows quadratically with~$l$. This makes it
impossible to directly compute the $x\in V_l$ using computer algebra:
known algorithms for solving systems of non-linear polynomial
equations take time exponential in the dimension.

Instead of using computer algebra directly, Jean-Marc Couveignes
suggested that we use approximations and
height bounds. In its simplest form, this works as follows. Suppose
that $x$ is a rational number, $x=a/b$, with $a$ and $b$ in $\ZZ$
coprime. Suppose that we have an upper bound $M$ for
$\max(|a|,|b|)$. Then $x$ is determined by any approximation $y\in\RR$
of~$x$ such that $|y-x|<1/2M^2$, simply because for all $x'\neq x$
with $x'=a'/b'$, where $a'$ and $b'$ in $\ZZ$ satisfy
$\max(|a'|,|b'|)<M$, we have $|x'-x|=|(a'b-ab')/bb'|\geq 1/M^2$.

For the computation of~$K_l$, we consider the minimal polynomial $P_l$
in~$\QQ[T]$ of a carefully theoretically constructed generator
$\alpha$ of~$K_l$. We use approximations of all Galois conjugates
of~$\alpha$, i.e., of all roots of~$P_l$. Instead of working directly
with torsion points of $J_l$, we work with divisors on the
curve~$X_l$. Using this strategy, the problem of showing that $P_l$
can be computed in time polynomial in~$l$ is divided into two
different tasks. Firstly, to show that the number of digits necessary
for a good enough approximation of~$P_l$ is bounded by a fixed power
of~$l$. Secondly, to show that, given $l$ and~$n$, the coefficients of
$P_l$ can be approximated with a precision of $n$ digits in time
polynomial in~$n{\cdot}l$. The first problem is dealt with in
Chapters~\ref{sec_appl_Arakelov}, \ref{subsec_merkl},
and~\ref{chap_bnd_height}, using Arakelov geometry. The second problem
is solved in Chapters~\ref{sec_couveignes_TORSION}
and~\ref{sec_couveignes_modp}, in two ways: complex approximations
(numerical analysis), and approximations in the sense of reductions
modulo many small primes, using exact computations in Jacobians of
modular curves over finite fields. These five chapters form the
technical heart of this book. The preceding chapters are meant as an
introduction to them, or motivation for them, and the two chapters
following them give the main results as relatively straightforward
applications.

Chapters~\ref{chapcomput} and~\ref{chappgltau} stand a bit apart, as
they are concerned with some real computations of Galois
representations attached to modular forms. They use the method by
complex approximations, but do not use a rigorously proven bound for a
sufficient accuracy. Instead, the approximations provide good
\emph{candidates} for polynomials~$P_l$. The $P_l$ that are found have
the correct Galois group, and the right ramification
properties. Recent modularity results by Khare, Wintenberger and
Kisin, see \cite{Khare-Wintenberger1}, \cite{Khare-Wintenberger2}, and
\cite{Kisin1} and~\cite{Kisin2}, are then applied to \emph{prove} that
the candidates do indeed give the right Galois representations.

\section{Historical context: Schoof's algorithm}\label{sec_schoof}

The computation of Hecke operators from Galois representations and
congruences can be viewed as a generalisation of Schoof's method to
count points on elliptic curves over finite fields, see~\cite{Schoof1}
and~\cite{Schoof2}. René Schoof gave an algorithm to compute, for $E$
an elliptic curve over a finite field~$\FF_q$, the number $\#
E(\FF_q)$ of $\FF_q$-rational points in a time~$O((\log
q)^{5+\eps})$. His algorithm works as follows.

The elliptic curve is embedded, as usual, in the projective plane
$\PP^2_{\FF_q}$ as the zero locus of a Weierstrass equation, which, in
inhomogeneous coordinates, is of the form:
\[
y^2 + a_1xy +a_3y = x^3 + a_2 x^2 + a_4 x + a_6,
\]
with the $a_i$ in~$\FF_q$. We let $\FF_q\to\Fbar_q$ be an algebraic
closure. We let $\rF_q\colon E\to E$ denote the so-called
$q$-Frobenius. It is the endomorphism of $E$ with the property that
for all $(a,b)$ in the affine part of $E(\Fbar_q)$ given by the
Weierstrass equation above we have $\rF_q((a,b))=(a^q,b^q)$. The
theory of elliptic curves over finite fields says:
\begin{enumerate}
\item there is a unique integer $a$, called the \emph{trace} of $\rF_q$,
such that in the endomorphism ring of~$E$ one has $\rF_q^2-a\rF_q+q=0$;
\item $\# E(\FF_q) = 1-a+q$;
\item $|a|\leq 2 q^{1/2}$.
\end{enumerate}
So, computing $\# E(\FF_q)$ is equivalent to computing this
integer~$a$. Schoof's idea is now to compute $a$ modulo $l$ for small
prime numbers~$l$. If the product of the prime numbers $l$ exceeds
$4q^{1/2}$, the length of the interval in which we know $a$ to lie,
then the congruences modulo these $l$ determine $a$ uniquely. Analytic
number theory tells us that it will be sufficient to take all primes
$l$ up to approximately $(\log q)/2$.

Then the question is how one computes $a$ modulo~$l$. This should be
done in time polynomial in $\log q$ and~$l$. The idea is to use the
elements of order dividing~$l$ in~$E(\Fbar_q)$. We assume now that $l$
does not divide $q$, i.e., we avoid the characteristic of~$\FF_q$. For
each $l$, the kernel $E(\Fbar_q)[l]$ of multiplication by $l$ on
$E(\Fbar_q)$ is a two-dimensional vector space over $\FF_l$. The map
$\rF_q$ gives an endomorphism of $E(\Fbar_q)[l]$, and it follows that
the image of $a$ in $\FF_l$ is the unique element of $\FF_l$, also
denoted $a$, such that for each $v$ in $E(\Fbar_q)[l]$ we have
$a\rF_q(v) = \rF_q^2 v + qv$. We remark that the image of $a$ in
$\FF_l$ is the trace of the endomorphism of $E(\Fbar_q)[l]$ given
by~$\rF_q$, but this is not really used at this point.

To find this element $a$ of $\FF_l$, one proceeds as follows. We
suppose that $l\neq 2$. There is a unique monic element $\psi_l$ of
$\FF_q[x]$ of degree $(l^2-1)/2$, whose roots in $\Fbar_q$ are
precisely the $x$-coordinates of the $l^2-1$ non-zero elements in
$E(\Fbar_q)[l]$ (the rational function $x$ on $E$ is a degree two map
to $\PP^1_{\FF_q}$, which as such is the quotient for the
multiplication by $-1$ map on~$E$). One then lets $A_l$ be the
$\FF_q$-algebra obtained as:
\[
A_l := 
\FF_q[x,y]/(y^2 + a_1xy +a_3y -x^3 - a_2 x^2 - a_4 x - a_6,\psi_l(x)).
\]
The dimension of $A_l$ as $\FF_q$-vector space is $l^2-1$. An
equivalent description of $A_l$ is to say that it is the affine
coordinate ring of the subscheme of points of order~$l$ of~$E$. By
construction of~$A_l$, there is a tautological $A_l$-valued point $v$
in~$E(A_l)$ (its coordinates are the images of $x$ and $y$
in~$A_l$). Now to find the element $a$ of $\FF_l$ that we are looking
for one then tries one by one the elements $i$ in $0,\pm 1,\ldots,\pm
(l-1)/2$ until $i\rF_q(v) = \rF_q^2 v + qv$; then $i = a\mod l$. 

It is easy to see that all required computations can be done in time
$O((\log q)^{5+\eps})$ (using fast arithmetic for the elementary
operations, e.g., a multiplication in $A_l$ costs about $(l^2(\log
q))^{1+\eps}$ time; $l^2(\log q)$ is the number of bits needed to
store one element of~$A_l$).

For the sake of completeness, let us mention that shortly after the
appearance of Schoof's algorithm, Atkin and Elkies have added some
improvements to it, making it possible in certain cases to reduce the
dimension of the $\FF_q$-algebra from $l^2-1$ to linear in $l+1$
or~$l-1$. This improvement, called the Schoof-Atkin-Elkies (SEA)
algorithm, is important mainly for implementations. Its (average)
complexity is $O((\log q)^{4+\eps})$; for details, the reader is
referred to~\cite{Schoof2}.

\section{Schoof's algorithm described in terms of étale cohomology}
\label{sec_schoof_etale}
In order to describe Schoof's algorithm in the previous section, we
referred to the theory of elliptic curves over finite fields. But
there is a more general framework for getting information on the
number of rational points of algebraic varieties over finite fields:
cohomology, and Lefschetz's trace formula. Cohomology exists in many
versions. The version directly related to Schoof's algorithm is étale
cohomology with coefficients in~$\FF_l$. Standard references for étale
cohomology are \cite{SGA4}, \cite{SGA4.5}, \cite{SGA5}, \cite{Milne1},
\cite{Freitag-Kiehl}. The reader is referred to these references for
the notions that we will use below. We also recommend Appendix~C
of~\cite{Hartshorne1}.

For the sake of precision, let us say that we define the notion of
\emph{algebraic variety over a field $k$} to mean \emph{$k$-scheme
  that is separated and of finite type}. Attached to an algebraic
variety $X$ over a field $k$ there are \emph{étale cohomology groups
  with compact supports} $\rH^i_\compact(X_\et,\FF_l)$, for all $i\geq
0$ and for all prime numbers~$l$. Actually, the coefficients $\FF_l$
can be replaced by more general objects, sheaves of Abelian groups on
the étale site $X_\et$ of $X$, but we do not need this now. If $X$ is
a proper $k$-scheme, then the $\rH^i_\compact(X_\et,\FF_l)$ are equal
to the étale cohomology groups $\rH^i(X_\et,\FF_l)$ without condition
on supports.

If $k$ is separably closed then the $\rH^i_\compact(X_\et,\FF_l)$ are
finite dimensional $\FF_l$-vector spaces, zero for $i>2\dim(X)$. In
that case, they are the analog of the more easily defined cohomology
groups $\rH^i_\compact(X,\calF)$ for complex analytic varieties: the
derived functors of the functor that associates to a sheaf $\calF$ of
$\ZZ$-modules on $X$ equipped with its Archimedean topology its
$\ZZ$-module of global sections whose support is compact.

The construction of the $\rH^i_\compact(X_\et,\FF_l)$ is functorial
for proper morphisms: a proper morphism $f\colon X\to Y$ of algebraic
varieties over $k$ induces a pullback morphism $f^*$ from
$\rH^i_\compact(Y_\et,\FF_l)$ to $\rH^i_\compact(X_\et,\FF_l)$.

Let now $X$ be an algebraic variety over~$\FF_q$. Then we have the
$q$-Frobenius morphism $\rF_q$ from $X$ to itself, and, by extending
the base field from $\FF_q$ to $\Fbar_q$, from $X_{\Fbar_q}$ to
itself. This morphism $\rF_q$ is proper, hence induces maps:
\[
\rF_q^*\colon \rH^i_\compact(X_{\Fbar_q,\et},\FF_l) \lto 
\rH^i_\compact(X_{\Fbar_q,\et},\FF_l).
\]
Hence, for each $i$ in~$\ZZ$, the trace
$\trace(\rF_q^*,\rH^i_\compact(X_{\Fbar_q,\et},\FF_l))$ of the map
above is defined, and it is zero for $i<0$ and $i>2\dim(X)$. The set
of fixed points of $\rF_q$ on $X(\Fbar_q)$ is precisely the
subset~$X(\FF_q)$. The Lefschetz trace formula then gives the
following identity in~$\FF_l$:
\begin{eqn}\label{LTF}
\# X(\FF_q) = \sum_i (-1)^i
\trace(\rF_q^*,\rH^i_\compact(X_{\Fbar_q,\et},\FF_l)).
\end{eqn}
We can now say how Schoof's algorithm is related to étale
cohomology. We consider again an elliptic curve $E$ over a finite
field~$\FF_q$. We assume that $l$ does not divide~$q$. Then, as for any
smooth proper geometrically connected curve,
$\rH^0(E_{\Fbar_q,\et},\FF_l)=\FF_l$ and $\rF_q^*$ acts on it as the
identity, and $\rH^2(E_{\Fbar_q,\et},\FF_l)$ is one-dimensional and
$\rF_q^*$ acts on it by multiplication by~$q$, the degree
of~$\rF_q$. According to the trace formula~(\ref{LTF}), we have:
\[
\# E(\FF_q) = 1 - \trace(\rF_q^*,\rH^1(E_{\Fbar_q,\et},\FF_l)) + q.
\]
It follows that for the integer $a$ of the previous section, the
trace of Frobenius, we have, for all $l$ not dividing~$p$ the identity
in~$\FF_l$:
\[
a = \trace(\rF_q^*,\rH^1(E_{\Fbar_q,\et},\FF_l)).
\]
This identity is explained by the fact that there is a natural
isomorphism, compatible with the action of $\rF_q$:
\[
\rH^1(E_{\Fbar_q,\et},\FF_l) = E(\Fbar_q)[l].
\]
Let us describe how one constructs this isomorphism. On $E_\et$ we
have the short exact sequence of sheaves, called the Kummer sequence:
\[
0 \lto \mu_l \lto \Gm \lto \Gm \lto 0
\]
where the map on $\Gm$ is multiplication by~$l$ in the group law of
$\Gm$, i.e., taking $l$th powers. This short exact sequence gives an
exact sequence of cohomology groups after pullback
to~$E_{\Fbar_q,\et}$:
\begin{multline*}
\{1\}\lto \mu_l(\Fbar_q) \lto \Fbar_q^\times \lto \Fbar_q^\times \lto 
\rH^1(E_{\Fbar_q,\et},\mu_l) \lto \\
\lto \rH^1(E_{\Fbar_q,\et},\Gm) \lto \rH^1(E_{\Fbar_q,\et},\Gm) \lto \cdots
\end{multline*}
Just as for any scheme, one has:
\[
\rH^1(E_{\Fbar_q,\et},\Gm) = \Pic(E_{\Fbar_q})
\]
% reference ?????
It follows that
\[
\rH^1(E_{\Fbar_q,\et},\mu_l) = \Pic(E_{\Fbar_q})[l].
\]
Finally, using the exact sequence:
\[
0 \lto \Pic^0(E_{\Fbar_q}) \lto \Pic(E_{\Fbar_q})\stackrel{\deg}{\lto} 
\ZZ \lto 0
\]
and the fact that $E$ is its own Jacobian variety, i.e.,
$\Pic^0(E_{\Fbar_q})=E(\Fbar_q)$, we obtain:
\[
\rH^1(E_{\Fbar_q,\et},\mu_l) = \Pic(E_{\Fbar_q})[l] = 
\Pic^0(E_{\Fbar_q})[l] = E(\Fbar_q)[l].
\]
The choice of an isomorphism between $\mu_l(\Fbar_q)$ and $\FF_l$
gives us the desired isomorphism between $\rH^1(E_{\Fbar_q,\et},\FF_l)$
and~$E(\Fbar_q)[l]$. In fact, we note that by using the Weil pairing
from $E(\Fbar_q)[l]\times E(\Fbar_q)[l]$ to $\mu_l(\Fbar_q)$, we get an
isomorphism:
\[
\rH^1(E_{\Fbar_q,\et},\FF_l) = E(\Fbar_q)[l]^\vee
\]
that is more natural than the one used above; in particular, it does
not depend on the choice of an isomorphism $\FF_l\to\mu_l(\Fbar_q)$.

\section{Some natural new directions}

We have seen that the two-dimensional $\FF_l$-vector spaces that are
used in Schoof's algorithm for elliptic curves can also be seen as
étale cohomology groups. A natural question that arises is then the
following.
\begin{quote}
\emph{Are there other interesting cases where étale cohomology groups
can be used to construct polynomial time algorithms for counting
rational points of varieties over finite fields?}
\end{quote}
A more precise question is the following.
\begin{quote}
\emph{ Let $n$ and $m$ be in~$\ZZ_{\geq0}$, and let $f_1,\ldots,f_m$
  be in $\ZZ[x_1,\ldots,x_n]$. Is there an algorithm that on
  input a prime number $p$ computes $\#\{a\in\FF_p^n\,|\,\forall
  i:f_i(a)=0\}$ in time polynomial in $\log p$?  }
\end{quote}
We believe that the answer to this question is yes, and that the way
in which such an algorithm can work is to compute {\'e}tale
cohomology.

\subsection{Curves of higher genus}
The first step in the direction of this question was taken by Jonathan
Pila. In \cite{Pila1} he considered principally polarised Abelian
varieties of a fixed dimension, and curves of a fixed genus, and
showed that in those cases polynomial time algorithms for computing the
number of rational points over finite fields exist. In these cases,
the only relevant cohomology groups are in degree one, i.e., they are
of the form $\rH^1(X_{\Fbar_q,\et},\FF_l)$ with $X$ a smooth proper curve,
or an Abelian variety, over the field~$\FF_q$. As in Schoof's
algorithm, the way to deal with these cohomology groups is to view
them as $J(\Fbar_q)[l]$, the kernel of multiplication by $l$ on the
Abelian variety~$J$. In the case where $X$ is a curve, one lets $J$ be
the Jacobian variety of~$X$. 

As Pila makes use of explicit systems of equations for Abelian
varieties, his algorithm has a running time that is at least
exponential in the dimension of the Abelian variety, and hence, in the
case of curves, as a function of the genus of the curve.

The current state of affairs concerning the question of counting the
rational points of curves over finite fields seems still to be the
same: algorithms have a running time that is exponential in the
genus. As an illustration, let us mention that in
\cite{Adleman-Huang1} Adleman and Huang give an algorithm that
computes $\# X(\FF_q)$ in time $(\log q)^{O(g^2\log g)}$, where $X$ is
a hyperelliptic curve over $\FF_q$, and where $g$ is the genus of~$X$.

Recent progress in the case where the characteristic of the finite
fields $\FF_q$ is fixed, using so-called $p$-adic methods, will be
discussed in Section~\ref{sec_p-adic} below. In that case, there are
algorithms whose running time is polynomial in $g$ and~$\log q$.

\subsection{Higher degree cohomology, modular forms}
Another direction in which one can try to generalise Schoof's
algorithm is to varieties of higher dimension, where non-trivial
cohomology groups of degree higher than one are needed. In this
context, we would call the degree 2 cohomology group of a curve
trivial, because the trace of $\rF_q$ on it is~$q$.

More generally speaking, cohomology groups, but now with $l$-adic
coefficients, that are of dimension one are expected
%%% I think that we know this. If it is about H^i(X), X over F_q,
%%% then, by considering small affine opens of X, and alterations,
%%% reduce to smooth X that can be lifted with normal crossings
%%% divisor at infinity. Then use l-adic Hodge theory....%%%%%%%%%%%%
to have the property that the trace of $\rF_q$ can only be of the form
$q^n\zeta$, with $n$ an integer greater than or equal to zero, and
$\zeta$ a root of unity. This means that one-dimensional cohomology
groups are not so challenging. Indeed, it is the fact that for
elliptic curves over $\FF_p$ all integers in the Hasse interval
$[p+1-2p^{1/2},p+1+2p^{1/2}]$ can occur that makes the problem of
point counting very different from point counting on non-singular
quadric surfaces in~$\PP^3_{\FF_q}$, for example, where the outcome
can only be $q^2+2q+1$ or $q^2+1$.

It follows that the simplest case to consider is cohomology groups of
dimension two, in degree at least two, on which the action of $\rF_q$
is not given by a simple rule as in the one-dimensional case. Such
cohomology groups are provided by modular forms, as we will explain
later in Section~\ref{sec_modforms}. Let us just say for the moment,
that there is a direct relation with elliptic curves, via the concept
of \emph{modularity} of elliptic curves over~$\QQ$, that we will now
sketch.

Let $E$ be an elliptic curve over~$\QQ$, given by some Weierstrass
equation. Such a Weierstrass equation can be chosen to have its
coefficients in~$\ZZ$. A Weierstrass equation for $E$ with
coefficients in~$\ZZ$ is called \emph{minimal} if its
\emph{discriminant} is minimal among all Weierstrass equations for $E$
with coefficients in~$\ZZ$; this discriminant then only depends on $E$
and will be denoted~$\discr(E)$. In fact, two minimal Weierstrass
equations define isomorphic curves in $\PP^2_\ZZ$, the projective
plane over~$\ZZ$. In other words, $E$ has a Weierstrass minimal model
over~$\ZZ$, that will be denoted by~$E_\ZZ$. For each prime number
$p$, we let $E_{\FF_p}$ denote the curve over $\FF_p$ given by
reducing a minimal Weierstrass equation modulo~$p$; it is the fibre of
$E_\ZZ$ over~$\FF_p$. The curve $E_{\FF_p}$ is smooth if and only if
$p$ does not divide $\discr(E)$. The possible singular fibres have
exactly one singular point: an ordinary double point with rational
tangents, or with conjugate tangents, or an ordinary cusp. The three
types of reduction are called split multiplicative, non-split
multiplicative and additive, respectively, after the type of group law
that one gets on the complement of the singular point.  For each $p$
we then get an integer $a_p$ by requiring the following identity:
\[
p+1-a_p = \# E(\FF_p).
\]
This means that for all $p$, $a_p$ is the trace of $\rF_p$ on the
degree one étale cohomology of $E_{\Fbar_p}$, with coefficients in
$\FF_l$, or in $\ZZ/l^n\ZZ$ or in the $l$-adic numbers~$\ZZ_l$. For
$p$ not dividing $\discr(E)$ we know that $|a_p|\leq 2p^{1/2}$. If
$E_{\FF_p}$ is multiplicative, then $a_p=1$ or $-1$ in the split and
non-split case. If $E_{\FF_p}$ is additive, then $a_p=0$. We also
define, for each $p$ an element $\eps(p)$ in $\{0,1\}$ by setting
$\eps(p)=1$ for $p$ not dividing $\discr(E)$ and setting $\eps(p)=0$
for $p$ dividing $\discr(E)$. The \emph{Hasse-Weil $L$-function} of
$E$ is then defined as:
\[
L_E(s) = \prod_p L_{E,p}(s), \qquad
L_{E,p}(s) = (1 - a_pp^{-s} + \eps(p)pp^{-2s})^{-1}, 
\]
for $s$ in $\CC$ with $\Re(s)>3/2$ (indeed, the fact that $|a_p|\leq
2p^{1/2}$ implies that the product converges for such~$s$). To explain
this function more conceptually, we note that for all $p$ and for all
$l\neq p$ we have the identity:
\[
1 - a_pt + \eps(p)pt^2 = \det(1-t\rF_p^*,\rH^1(E_{\Fbar,\et},\QQ_l))
\]
The reader should notice that now we use étale cohomology with
coefficients in $\QQ_l$, the field of $l$-adic numbers, and not
in~$\FF_l$. The reason for this is that we want the last identity
above to be an identity between polynomials with integer coefficients,
and not with coefficients in~$\FF_l$.

The function $L_E$ was conjectured to have a holomorphic continuation
over all of~$\CC$, and to satisfy a certain precisely given functional
equation relating the values at $s$ and $2-s$. In that functional
equation appears a certain positive integer $N_E$ called the
\emph{conductor} of~$E$, composed of the primes $p$
dividing~$\discr(E)$ with exponents that depend on the behaviour of
$E$ at~$p$, i.e., on $E_{\ZZ_p}$. This conjecture on continuation and
functional equation was proved for semistable $E$ (i.e., $E$ such that
there is no $p$ where $E$ has additive reduction) by Wiles and
Taylor-Wiles, and in the general case by Breuil, Conrad, Diamond and
Taylor; see \cite{Edixhoven1} for an overview of this. In fact, the
continuation and functional equation are direct consequences of the
modularity of~$E$ that was proved by Wiles, Taylor-Wiles etc. (see
below). The weak Birch and Swinnerton-Dyer conjecture says that the
dimension of the $\QQ$-vector space $\QQ\otimes E(\QQ)$ is equal to
the order of vanishing of $L_E$ at~$1$.  Anyway, the function $L_E$
gives us integers $a_n$ for all $n\geq1$ as follows:
\[
L_E(s) = \sum_{n\geq 1} a_n n^{-s}, \quad \text{for $\Re(s) > 3/2$}.
\]
From these $a_n$ one can then consider the following function:
\[
f_E\colon \HH = \{\tau\in\CC\,|\,\Im(\tau)>0\}\to\CC, \quad 
\tau\mapsto\sum_{n\geq 1}a_n e^{2\pi i n\tau}.
\]
Equivalently, we have:
\[
f_E = \sum_{n\geq 1}a_nq^n, 
\quad \text{with}\quad q\colon \HH\to \CC, \quad \tau\mapsto e^{2\pi i\tau}.
\]
A more conceptual way to state the relation between $L_E$ and $f_E$ is
to say that $L_E$ is obtained, up to elementary factors, as the
\emph{Mellin transform} of~$f_E$:
\[
\int_0^\infty f_E(it)t^s\frac{dt}{t} = (2\pi)^{-s}\Gamma(s)L_E(s), 
\quad \text{for $\Re(s)>3/2$}.
\]
After all these preparations, we can finally state what the modularity
of $E$ means:
\begin{quote}
\emph{
$f_E$ is a modular form of weight two for the congruence subgroup
$\Gamma_0(N_E)$ of\/ $\SL_2(\ZZ)$.
}
\end{quote}
For some more details on the concept of modular forms we refer to
Section~\ref{sec_modforms}. At this moment, we just want to say that
the last statement means that $f_E$ has, as Mazur says in Singh's BBC
documentary on Wiles's proof of Fermat's Last Theorem, an enormous
amount of symmetry. This symmetry is with respect to the action of
$\GL_2(\QQ)^+$, the group of invertible 2 by 2 matrices with
coefficients in~$\QQ$ whose determinant is positive, on the upper half
plane~$\HH$. This symmetry gives, by Mellin transformation, the
functional equation of~$L_E$. Conversely, it had been proved
in~\cite{Weil1} by Weil that if sufficiently many twists of $L_E$ by
Dirichlet characters satisfy the conjectured holomorphic continuation
and functional equation, then $f_E$ is a modular form of the type
mentioned.

We now remark that Schoof's algorithm implies that, for $p$ prime, the
coefficient $a_p$ in the $q$-expansion of $f_E=\sum_{n\geq1}a_nq^n$
can be computed in time polynomial in $\log p$. One of the aims of the
research project described in this report is to generalise this last
fact to certain modular forms of higher weight. Before we give precise
definitions in Section~\ref{sec_modforms}, we will discuss a typical
case in the next section.

\section{More historical context: congruences for Ramanujan's 
$\tau$-function}
\label{sec_hist_cong_tau}
References for this section are the articles \cite{Serre1},
\cite{Swinnerton-Dyer1} and \cite{Deligne1} by Serre, Swinnerton-Dyer
and Deligne.

A typical example of a modular form of weight higher than two is the
\emph{discriminant} modular form, usually denoted~$\Delta$. One way to
view $\Delta$ is as the holomorphic function on the upper half plane
$\HH$ given by:
\begin{eqn}
\Delta = q\prod_{n\geq 1}(1-q^n)^{24},
\end{eqn}
where $q$ is the function from $\HH$ to $\CC$ given by $z\mapsto
\exp(2\pi i z)$. The coefficients in the power series expansion:
\begin{eqn}
\Delta = \sum_{n\geq 1}\tau(n)q^n
\end{eqn}
define the famous \emph{Ramanujan $\tau$-function}. 

To say that $\Delta$ is a modular form of weight $12$ for the group
$\SL_2(\ZZ)$ means that for all elements 
$(\begin{smallmatrix}a & b\\c & d\end{smallmatrix})$
of $\SL_2(\ZZ)$ the following identity holds for all $z$ in~$\HH$:
\begin{eqn}\label{eqn_delta_mod_form}
\Delta\left(\frac{az+b}{cz+d}\right) = (cz+d)^{12}\Delta(z),
\end{eqn}
which is equivalent to saying that the multi-differential form
$\Delta(z)(dz)^{\otimes 6}$ is invariant under the action
of~$\SL_2(\ZZ)$.  As $\SL_2(\ZZ)$ is generated by the elements
$(\begin{smallmatrix}1 & 1\\0 & 1\end{smallmatrix})$ and
$(\begin{smallmatrix}0 & -1\\1 & 0\end{smallmatrix})$, it suffices to
check the identity in~(\ref{eqn_delta_mod_form}) for these two
elements. The fact that $\Delta$ is $q$ times a power series in $q$
means that $\Delta$ is a \emph{cusp form}: it vanishes at
``$q=0$''. It is a fact that $\Delta$ is the first example of a
non-zero cusp form for $\SL_2(\ZZ)$: there is no non-zero cusp form
for $\SL_2(\ZZ)$ of weight smaller than $12$, i.e., there are no
non-zero holomorphic functions on $\HH$ satisfying
(\ref{eqn_delta_mod_form}) with the exponent $12$ replaced by a
smaller integer, whose Laurent series expansion in~$q$ is $q$ times a
power series. Moreover, the $\CC$-vector space of such functions of
weight $12$ is one-dimensional, and hence $\Delta$ is a basis of it.

The one-dimensionality of this space has as a consequence that
$\Delta$ is an eigenform for certain operators on this space, called
\emph{Hecke operators}, that arise from the action on~$\HH$ of
$\GL_2(\QQ)^+$, the subgroup of $\GL_2(\QQ)$ of elements whose
determinant is positive. This fact explains that the coefficients
$\tau(n)$ satisfy certain relations which are summarised by the
following identity of Dirichlet series (converging for $\Re(s)\gg 0$,
for the moment, or just formal series, if one prefers that):
\begin{eqn}
L_\Delta(s) := \sum_{n\geq 1} \tau(n)n^{-s} = 
\prod_p(1-\tau(p)p^{-s}+p^{11}p^{-2s})^{-1}.
\end{eqn}
These relations:
\begin{eqn}\label{eqn_tau_identities}
\begin{aligned}
& \tau(mn) = \tau(m)\tau(n) && \text{if $\gcd(m,n)=1$}\\
& \tau(p^n) = \tau(p^{n-1})\tau(p) - p^{11}\tau(p^{n-2}) 
&&\text{if $p$ is prime and $n\geq 2$}
\end{aligned}
\end{eqn}
were conjectured by Ramanujan, and proved by Mordell.  Using these
identities, $\tau(n)$ can be expressed in terms of the $\tau(p)$ for
$p$ dividing~$n$.

As $L_\Delta$ is the Mellin transform of $\Delta$, $L_\Delta$ is
holomorphic on~$\CC$, and satisfies the functional equation (Hecke):
\[
(2\pi)^{-(12-s)} \Gamma(12-s) L_\Delta(12-s) = 
(2\pi)^{-s} \Gamma(s) L_\Delta(s).
\]

The famous \emph{Ramanujan conjecture} states that for all primes $p$
one has the inequality:
\begin{eqn}\label{eqn_ramanujan_inequality}
|\tau(p)| < 2p^{11/2},
\end{eqn}
or, equivalently, that the complex roots of the polynomial
$x^2-\tau(p)x+p^{11}$ are complex conjugates of each other, and hence
are of absolute value~$p^{11/2}$. This conjecture was proved by
Deligne as a consequence of his article~\cite{Deligne1} and his proof
of the analog of the Riemann hypothesis in the Weil conjectures
in~\cite{Deligne2}. 

Finally, Ramanujan conjectured congruences for the integers $\tau(p)$
with $p$ prime, modulo certain powers of certain small prime
numbers. In order to state these congruences we define, for $n\geq 1$ and
$r\geq 0$:
\[
\sigma_r(n) := \sum_{1\leq d | n} d^r,
\]
i.e., $\sigma_r(n)$ is the sum of the $r$th powers of the positive
divisors of~$n$. We will now list the congruences that are given in
the first pages of~\cite{Swinnerton-Dyer1}:
\begin{align*}
&\tau(n)\equiv \sigma_{11}(n) &&\mod 2^{11} &&\text{if $n\equiv 1\mod
  8$}\\ 
&\tau(n)\equiv 1217\sigma_{11}(n) &&\mod 2^{13} &&\text{if $n\equiv 3\mod
  8$}\\
&\tau(n)\equiv 1537\sigma_{11}(n) &&\mod 2^{12} &&\text{if $n\equiv 5\mod
  8$}\\
&\tau(n)\equiv 705\sigma_{11}(n) &&\mod 2^{14} &&\text{if $n\equiv 7\mod
  8$}
\end{align*}

\begin{align*}
&\tau(n)\equiv n^{-610}\sigma_{1231}(n) &&\mod 3^6 &&\text{if $n\equiv 1\mod
  3$}\\ 
&\tau(n)\equiv n^{-610}\sigma_{1231}(n) &&\mod 3^7 &&\text{if $n\equiv 2\mod
  3$}
\end{align*}

\begin{align*}
&\tau(n)\equiv n^{-30}\sigma_{71}(n) &&\mod 5^3 &&\text{if $n$ is prime
  to $5$}
\end{align*}

\begin{align*}
&\tau(n)\equiv n\sigma_9(n) &&\mod 7 &&\text{if $n\equiv 0$, $1$, $2$
  or $4 \mod 7$}\\ 
&\tau(n)\equiv n\sigma_9(n) &&\mod 7^2 &&\text{if $n\equiv 3$, $5$ or 
  $6 \mod 7$}
\end{align*}

\begin{align*}
&\tau(p)\equiv 0 &&\mod 23 &&\text{if $p$ is prime and not a square mod $23$}\\ 
&\tau(p)\equiv 2 &&\mod 23 &&\text{if $p\neq 23$ is a prime of the form
    $u^2+23v^2$}\\ 
&\tau(p)\equiv -1 &&\mod 23 &&\text{for other primes $p\neq 23$}
\end{align*}

\begin{align*}
&\tau(n)\equiv \sigma_{11}(n)\mod 691
\end{align*}

The reader is referred to~\cite{Swinnerton-Dyer1} for the origin and
for proofs of these congruences. There, Swinnerton-Dyer remarks that
the proofs do little explain why such congruences occur. Serre
conjectured an explanation in~\cite{Serre1}. First of all, Serre
conjectured the existence, for each prime number~$l$, of a continuous
representation:
\begin{eqn}
\rho_l\colon \Gal(\Qbar/\QQ) \lto \Aut(V_l),
\end{eqn}
with $V_l$ a two-dimensional $\QQ_l$-vector space, such that $\rho_l$
is unramified at all primes $p\neq l$, and such that for all $p\neq l$
the characteristic polynomial of $\rho_l(\Frob_p)$ is given by:
\begin{eqn}
\det(1-x\Frob_p,V_l) = 1-\tau(p)x + p^{11}x^2.
\end{eqn}
To help the reader, let us explain what unramified at $p$ means, and
what the Frobenius elements $\Frob_p$ are. For $p$ prime, we let
$\QQ_p$ denote the topological field of $p$-adic numbers, and
$\QQ_p\to\Qbar_p$ an algebraic closure. The action of
$\Gal(\Qbar/\QQ)$ on the set $\Hom(\Qbar,\Qbar_p)$ of embeddings of
$\Qbar$ into $\Qbar_p$ is transitive, and each embedding induces an
injection from $\Gal(\Qbar_p/\QQ_p)$ into $\Gal(\Qbar/\QQ)$, the image
of which is called a decomposition group of $\Gal(\Qbar/\QQ)$
at~$p$. The injections from $\Gal(\Qbar_p/\QQ_p)$ into
$\Gal(\Qbar/\QQ)$ and the corresponding decomposition groups at $p$
obtained like this are all conjugated by the action of
$\Gal(\Qbar/\QQ)$. In order to go further we need to say a bit about
the structure of $\Gal(\Qbar_p/\QQ_p)$. We let $\QQ_p^\unr$ be the
maximal unramified extension of $\QQ_p$ in $\Qbar_p$, i.e., the
composite of all finite extensions $K$ of $\QQ_p$ in $\Qbar_p$ such
that $p$ is a uniformiser for the integral closure $O_K$ of $\ZZ_p$
in~$K$. We let $\ZZ_p^\unr$ be the integral closure of $\ZZ_p$ in
$\QQ_p^\unr$; it is a local ring, and its residue field is an
algebraic closure $\Fbar_p$ of~$\FF_p$. The sub-extension $\QQ_p^\unr$
gives a short exact sequence:
\begin{eqn}
I_p \into \Gal(\Qbar_p/\QQ_p) \onto \Gal(\Fbar_p/\FF_p).
\end{eqn}
The subgroup $I_p$ of $\Gal(\Qbar_p/\QQ_p)$ is called the inertia
subgroup. The quotient $\Gal(\Fbar_p/\FF_p)$ is canonically isomorphic
to $\Zhat$, the profinite completion of~$\ZZ$, by demanding that the
element $1$ of $\Zhat$ corresponds to the Frobenius element $\Frob_p$
of $\Gal(\Fbar_p/\FF_p)$ that sends $x$ to $x^p$ for each $x$
in~$\Fbar_p$. 

Let now $\rho_l$ be a continuous representation from $\Gal(\Qbar/\QQ)$
to $\GL(V_l)$ with $V_l$ a finite dimensional $\QQ_l$-vector
space. Each embedding of $\Qbar$ into $\Qbar_p$ then gives a
representation of $\Gal(\Qbar_p/\QQ_p)$ on~$V_l$. Different embeddings
give isomorphic representations because they are conjugated by an
element in the image of $\Gal(\Qbar/\QQ)$ under~$\rho_l$. We now
choose one embedding, and call the representation $\rho_{l,p}$ of
$\Gal(\Qbar_p/\QQ_p)$ on $V_l$ obtained like this the local
representation at $p$ attached to~$\rho_l$. This being defined,
$\rho_l$ is then said to be unramified at a prime $p$ if $\rho_{l,p}$
factors through the quotient $\Gal(\Qbar_p/\QQ_p) \to
\Gal(\Fbar_p/\FF_p)$, i.e., if $I_p$ acts trivially on~$V_l$. If
$\rho_l$ is unramified at~$p$, then we get an element
$\rho_l(\Frob_p)$ in $\GL(V_l)$. This element depends on our chosen
embedding of $\Qbar$ into $\Qbar_p$, but its conjugacy class under
$\rho_l(\Gal(\Qbar/\QQ))$ does not. In particular, we get a
well-defined conjugacy class in~$\GL(V_l)$, and so the characteristic
polynomial of $\rho_l(\Frob_p)$ is now defined if $\rho_l$ is
unramified at~$p$. 

Continuous representations such as $\rho_l$ can be reduced modulo
powers of $l$ as follows. The compactness of $\Gal(\Qbar/\QQ)$ implies
that with respect to a suitable basis of $V_l$ the representation
$\rho_l$ lands in $\GL_2(\ZZ_l)$, and hence gives representations to
$\GL_2(\ZZ/l^n\ZZ)$ for all $n\geq0$. This reduction of $\rho_l$
modulo powers of $l$ is not unique, but the semi-simplification of the
reduction modulo $l$ is well-defined, i.e., two reductions lead to the
same Jordan-Hölder constituents. According to Serre, the congruences
above would then be explained by properties of the image of~$\rho_l$.

For example, if the image of the reduction modulo~$l$ of~$\rho_l$ is
reducible, say an extension of two characters $\alpha$ and~$\beta$
from $\Gal(\Qbar/\QQ)$ to~$\FF_l^\times$, then one has the identity
in~$\FF_l$, for all $p\neq l$:
\begin{eqn}\label{5.8}
\tau(p) \equiv \alpha(\Frob_p) + \beta(\Frob_p).
\end{eqn}
The characters $\alpha$ and $\beta$ are unramified outside~$l$. By the
Kronecker-Weber theorem, the maximal Abelian subextension of $\QQ\to
\Qbar$ that is unramified outside $l$ is the cyclotomic extension
generated by all $l$-power roots of unity, with Galois
group~$\ZZ_l^\times$. It then follows that $\alpha=\chi_l^n$ and
$\beta=\chi_l^m$ for suitable $n$ and $m$, where $\chi_l$ is the
character giving the action on the $l$th roots of unity in~$\Qbar$:
for all $\sigma$ in $\Gal(\Qbar/\QQ)$ and for all $\zeta$ in
$\Qbar^\times$ with $\zeta^l=1$ one has
$\sigma(\zeta)=\zeta^{\chi_l(\sigma)}$. The identity~(\ref{5.8}) in
$\FF_l$ above then takes the form:
\begin{eqn}
\tau(p) = p^n + p^m \mod l, \quad\text{for all $p\neq l$},
\end{eqn}
which indeed is of the same form as the congruences mod $l$ for
$\tau(p)$ listed above. For example, the congruence mod~$691$
corresponds to the statement that the reduction modulo~$l$ of~$\rho_l$
contains the two characters $1$ and~$\chi_l^{11}$.

Deligne, in \cite{Deligne1}, proved the existence of the~$\rho_l$, as
conjectured by Serre, by showing that they occur in the degree one
$l$-adic étale cohomology of certain sheaves on certain curves, and in
the degree $11$ étale cohomology with $\QQ_l$-coefficients of a
variety of dimension~$11$. This last variety is, loosely speaking, the
$10$-fold fibred product of the universal elliptic curve. Deligne's
constructions will be discussed in detail in
Sections~\ref{sec_modforms} and~\ref{sec_galreps}. It should be said
that Shimura had already shown how to construct Galois representations
in the case of modular forms of weight two; in that case one does not
need étale cohomology, but torsion points of Jacobians of modular
curves suffice, see~\cite{Shimura1}.

At this point we give the following precise statement, relating
Ramanujan's $\tau$-function to point counting on an algebraic
variety~$C_{10}$ (more precisely, a quasi-projective scheme
over~$\ZZ$), for which one easily writes down a system of
equations. Moreover, the statement relates the weight of $\Delta$ to
the classical question in geometry on cubic plane curves passing
through a given set of points: up to $9$ points the situation is easy
and the count is given by a polynomial.

\begin{prop}\label{prop_tau_point-counting}
For $n\in\ZZ_{\geq0}$, $q$ a prime power and\/ $\FF_q$ a finite field
with $q$ elements, let $C_n(\FF_q)$ be the set of
$(C,P_1,\ldots,P_n)$, where $C$ is a smooth cubic in $\PP^2_{\FF_q}$,
and $P_i\in C(\FF_q)$. Then there are $f_0,\ldots,f_{10}\in\ZZ[x]$
such that for all $\FF_q$ and $n\leq 9$ one has $\#
C_n(\FF_q)/\#\PGL_3(\FF_q)=f_n(q)$, and for all prime numbers~$p$:
\[
\# C_{10}(\FF_p)\,/\,\#\PGL_3(\FF_p) = -\tau(p)+f_{10}(p).
\]
\end{prop}
\begin{proof}
For $n$ in $\ZZ_{\geq0}$ and $\FF_q$ a field with $q$ elements, let
$\calE_n(\FF_q)$ denote the category, and also its set of objects, of
$(E/\FF_q,P_1,\ldots,P_n)$, where $E/\FF_q$ is an elliptic curve, and
$P_i\in E(\FF_q)$; the morphisms are the isomorphisms $\phi\colon E\to
E'$ such that $\phi(P_i)=P_i'$. For each~$n$, the category
$\calE_n(\FF_q)$ has only finitely many objects up to isomorphism, and
one defines:
\[
\#\calE_n(\FF_q) =  \sum_{x\in\calE_n(\FF_q)}\frac{1}{\#\Aut(x)},
\]
where, in the sum, one takes one $x$ per isomorphism class. It is well
known (see~\cite{Deligne1}, \cite{Behrend}) that for $n\leq 9$ the
functions $q\mapsto \#\calE_n(\FF_q)$ are given by certain elements
$f_n$ in $\ZZ[x]$, and that there is an $f_{10}$ in $\ZZ[x]$ such that
for all prime numbers $p$ one has
$\#\calE_{10}(\FF_p)=-\tau(p)+f_{10}(p)$. In view of this, the claims
in Proposition~\ref{prop_tau_point-counting} are a consequence of the
following equality, for all $n\in\ZZ_{\geq0}$ and all prime
powers~$q$:
\begin{eqn}\label{eqn_tau_point-counting}
\text{for all $n\in\ZZ_{\geq0}$ and all~$\FF_q$:}\quad
\# C_n(\FF_q) = \#\PGL_3(\FF_q)\cdot\#\calE_n(\FF_q).
\end{eqn}
We prove~(\ref{eqn_tau_point-counting}) by comparing the subsets on
both sides in which the underlying curves are fixed.  

Let $n\in\ZZ_{\geq0}$ and $q$ a prime power. Let $F$ be a nonsingular
projective geometrically irreducible curve of genus one over~$\FF_q$,
and let $E_0$ be its Jacobian. Then $F$ is an $E_0$-torsor. By Lang's
theorem, Theorem~2 of~\cite{Lang2}, $F(\FF_q)$ is not empty.

Let $C_n(\FF_q)_F$ be the subset of $C_n(\FF_q)$ consisting of the
$(C,P_1,\ldots,P_n)$ with $C$ isomorphic to~$F$. The number of $C$ in
$\PP^2_{\FF_q}$ that are isomorphic to~$F$ is the number of embeddings
$i\colon F\to\PP^2_{\FF_q}$, divided by $\#\Aut(F)$. Such embeddings are
obtained from line bundles $\calL$ of degree~3 on~$F$, together with a
basis, up to~$\FF_q^\times$, of $\calL(F)$. Hence the number of
embeddings is $\#\PGL_3(\FF_q)\cdot\# E_0(\FF_q)$. The group
$\Aut(F)$ has the subgroup of translations, $E_0(\FF_q)$, with
quotient $\Aut(E_0)$. So we find:
\[
\# C_n(\FF_q)_F = \#\PGL_3(\FF_q)\cdot(\# E_0(\FF_q))^n/\#\Aut(E_0).
\]
On the other hand, let $\calE_n(\FF_q)_{E_0}$ be the full subcategory
of $\calE_n(\FF_q)$ with objects the $(E_0,P_1,\ldots,P_n)$, with
$P_i$ in~$E_0(\FF_q)$. The group $\Aut(E_0)$ acts on the set of
objects of $\calE_n(\FF_q)_{E_0}$, and this action \emph{is} the set
of morphisms in $\calE_n(\FF_q)_{E_0}$. This means that:
\[
\#\calE_n(\FF_q)_{E_0} = (\# E_0(\FF_q))^n/\#\Aut(E_0).
\]
Summing over the isomorphism classes of $F$ gives
(\ref{eqn_tau_point-counting}). 
\end{proof}
\begin{rem}
The polynomials $f_n$ mentioned in
Proposition~\ref{prop_tau_point-counting} have been computed by Carel
Faber and Gerard van der Geer. Their result is:
\begin{align*}
f_0 & = x,\\
f_1 & = x^2+x,\\
f_2 & = x^3+3x^2+x-1,\\
f_3 & = x^4+6x^3+6x^2-2x-3,\\
f_4 & = x^5+10x^4+20x^3+4x^2-14x-7,\\
f_5 & = x^6+15x^5+50x^4+40x^3-30x^2-49x-15,\\
f_6 & = x^7+21x^6+105x^5+160x^4-183x^2-139x-31,\displaybreak[1]\\
f_7 & = x^8+28x^7+196x^6+469x^5+280x^4-427x^3-700x^2 \\
    & \quad -356x-63,\\
f_8 & = x^9+36x^8+336x^7+1148x^6+1386x^5-406x^4-2436x^3\\
    & \quad -2224x^2-860x-127,\\
f_9 & = x^{10}+45x^9+540x^8+2484x^7+4662x^6+1764x^5-6090x^4\\
    & \quad -9804x^3-6372x^2-2003x -255,\\
f_{10} & = x^{11}+55x^{10}+825x^9+4905x^8+12870x^7+12264x^6\\
    & \quad -9240x^5-33210x^4-33495x^3-17095x^2-4553x-511.
\end{align*}
\end{rem}
We refer to Birch~\cite{Birch} for results on the distribution
of the number of rational points on elliptic curves over finite
fields, that also make $\tau(p)$ appear. 

In~\cite{Swinnerton-Dyer1}, Swinnerton-Dyer gives results, partly
resulting from his correspondence with Serre, in which the
consequences of the existence of the $\rho_l$ for congruences of
$\tau(p)$ modulo $l$ are explored. A natural question to ask is if
there are primes $l$ other than $2$, $3$, $5$, $7$, $23$ and $691$
modulo which there are similar congruences for~$\tau(p)$. 

For each $p\neq l$, $\tau(p)$ is the trace of $\rho_l(\Frob_p)$, and
the determinant of $\rho_l(\Frob_p)$ equals~$p^{11}$. Hence, a
polynomial relation between $\tau(p)$ and $p^{11}$, valid modulo some
$l^n$ for all $p\neq l$, is a relation between the determinant and the
trace of all $\rho_l(\Frob_p)$ in~$\GL_2(\ZZ/l^n\ZZ)$. But
Chebotarev's theorem (see~\cite{Lang1}, or~\cite{Cassels-Frohlich},
for example) implies that every element of the image of
$\Gal(\Qbar/\QQ)$ in $\GL_2(\ZZ/l^n\ZZ)$ is of the form $\Frob_p$ for
infinitely many~$p$. Hence, such a polynomial relation is then valid
for all elements in the image of $\Gal(\Qbar/\QQ)$ in
$\GL_2(\ZZ/l^n\ZZ)$. For this reason, the existence of non-trivial
congruences modulo $l^n$ as above for $\tau(p)$ depends on this
image. 

The image of $\Gal(\Qbar/\QQ)$ in $\ZZ_l^\times$ under
$\det\circ\rho_l$ is equal to the subgroup of $11$th powers
in~$\ZZ_l^\times$. To explain this, we note that $\det\circ\rho_l$ is
a continuous character from $\Gal(\Qbar/\QQ)$ to~$\ZZ_l^\times$,
unramified outside~$l$, and such that $\Frob_p$ is mapped to $p^{11}$
for all $p\neq l$; this implies that $\det\circ\rho_l$ is the $11$th
power of the $l$-adic cyclotomic character
$\chi_l\colon\Gal(\Qbar/\QQ)\to\ZZ_l^\times$, giving the action of
$\Gal(\Qbar/\QQ)$ on the $l$-power roots of unity.

In order to state the results in~\cite{Swinnerton-Dyer1}, one calls a
prime number $l$ \emph{exceptional (for $\Delta$)} if the image of
$\rho_l$, taking values in $\GL_2(\ZZ_l)$, does \emph{not}
contain~$\SL_2(\ZZ_l)$. For $l$ not exceptional, i.e., such that the
image of $\rho_l$ contains $\SL_2(\ZZ_l)$, the image of
$\Gal(\Qbar/\QQ)$ in $\GL_2(\ZZ_l)\times\ZZ_l^\times$, under
$(\rho_l,\chi_l)$, is the subgroup $H$ of elements $(g,t)$ such that
$\det(g)=t^{11}$. This subgroup $H$ maps surjectively to
$\FF_l\times\FF_l^\times$ under $(g,t)\mapsto(\trace(g),t)$, and
therefore there can be no congruence for $\tau(p)$ modulo~$l$ as
above.

The Corollary to Theorem~4 in~\cite{Swinnerton-Dyer1} states, among
others, that the list of primes that are exceptional for $\Delta$ is
$\{2,3,5,7,23,691\}$. The main tool that is used and that we have not
discussed is the theory of modular forms modulo~$l$, or, equivalently,
the theory of congruences modulo~$l$ between modular forms. As a
consequence, there are no similar congruences for $\tau(p)$ modulo
primes other than the ones listed above. The special form of the
congruences modulo~$23$ is explained by the fact that in that case the
image of $\Gal(\Qbar/\QQ)$ in $\GL_2(\FF_{23})$ is dihedral; in the
other cases the residual representation, i.e., the representation to
$\GL_2(\FF_l)$, is reducible. In the case $l=2$, Swinnerton-Dyer has
determined the image of $\Gal(\Qbar/\QQ)$ in $\GL_2(\ZZ_2)$ exactly:
see the appendix in~\cite{Swinnerton-Dyer1}. 

The direction in which we generalise Schoof's algorithm is to give an
algorithm that computes for prime numbers $l$ that are not exceptional
for $\Delta$ the field extension $\QQ\to K_l$ that corresponds to the
representation of $\Gal(\Qbar/\QQ)$ to $\GL_2(\FF_l)$ that comes from
$\Delta$. The field $K_l$ is given in the form~$\QQ[x]/(f_l)$. The
computation has a running time that is polynomial in~$l$. It is fair
to say that this algorithm makes the mod $l$ Galois representations
attached to $\Delta$ accessible to computation, at least
theoretically. As the field extensions that are involved are
non-solvable, this should be seen as a step beyond computational class
field theory, and beyond the case of elliptic curves, in the direction
to make the results of Langlands's program accessible to computations.

As a consequence, one can compute $\tau(p)$ mod~$l$ in time polynomial
in $\log p$ and $l$, by reducing $f_l$ as above mod~$p$ and some more
computations that will be described later (see
Section~\ref{sec_comp_tau}).  By doing this for sufficiently many $l$,
just as in Schoof's algorithm, one then gets an algorithm that
computes $\tau(p)$ in time polynomial in~$\log p$.

In Section~\ref{sec_comp_tn} the method used here is generalised to
the case of modular forms for $\SL_2(\ZZ)$ of arbitrary weight. The
main result there is Theorem~\ref{thm_comp_tn}.

\section{Comparison with $\protect\lowercase{p}$-adic methods}
\label{sec_p-adic}
Before we start seriously with the theory of modular forms and the
Galois representations attached to them in the next chapter, we
make a comparison between our generalisation of Schoof's algorithm and
the so-called $p$-adic methods that have been developed since 2000 by
Satoh~\cite{Satoh1}, Kedlaya~\cite{Kedlaya1} (see
also~\cite{Edixhoven4}),
%% Bas should put this preprint on arxiv and refer to there%%%%%%%%%%
Hubrechts, \cite{Hubrechts1}), Lauder and Wan~\cite{Lauder-Wan1},
\cite{Lauder-Wan2}, \cite{Lauder1} and~\cite{Lauder2}, Fouquet,
Gaudry, Gürel and Harley~\cite{Fouquet-Gaudry-Harley1},
\cite{Gaudry-Gurel1}, Denef and Vercauteren and
Castryk~\cite{Denef-Vercauteren1}, \cite{Castryck-Denef-Vercauteren1},
Mestre, Lercier and Lubicz~\cite{Lercier-Lubicz1}, Carls, Kohel and
Lubicz, \cite{Carls-Kohel-Lubicz}, \cite{Carls-Lubicz}, and Gerkmann,
\cite{Gerkmann1} and \cite{Gerkmann2}. Actually, we should notice that
such a method was already introduced in~\cite{Kato-Lubkin} in 1982,
but that this article seems to have been forgotten (we thank Fre
Vercauteren for having drawn our attention to this article).

In all these methods, one works with fields of small
characteristic~$p$, hence of the form $\FF_q$ with $q=p^m$ and $p$
fixed. All articles cited in the previous paragraph have the common
property that they give algorithms for computing the number of
$\FF_q$-rational points on certain varieties $X$ over~$\FF_q$, using,
sometimes indirectly, cohomology groups with $p$-adic coefficients,
whence the terminology ``$p$-adic methods''.

For example, Satoh~\cite{Satoh1} uses the canonical lift of ordinary
elliptic curves and the action of the lifted Frobenius endomorphism on
the tangent space, which can be interpreted in terms of the algebraic
de Rham cohomology of the lifted curve. Kedlaya~\cite{Kedlaya1} uses
Monsky-Washnitzer cohomology of certain affine pieces of hyperelliptic
curves. In fact, all cohomology groups used here are de Rham type
cohomology groups, given by complexes of differential forms on certain
$p$-adic lifts of the varieties in question. Just as an example, let
us mention that Kedlaya~\cite{Kedlaya1} gives an algorithm that for
fixed $p\neq 2$ computes the zeta functions of hyperelliptic curves
given by equations:
\[
y^2 = f(x),
\]
where $f$ has arbitrary degree, in time~$m^3\deg(f)^4$. The running
times of the other algorithms are all similar, but all have in common
that the running time grows at least linearly in~$p$ (or linear in
$O(p^{1/2})$, in~\cite{Harvey1}), hence exponentially in~$\log
p$. The explanation for this is that somehow in each case non-sparse
polynomials of degree at least linear in $p$ have to be manipulated.

Summarising this recent progress, one can say that, at least from a
theoretical point of view, the problem of counting the solutions of
systems of polynomial equations over finite fields of a fixed
characteristic~$p$ and in a fixed number of variables has been
solved. If $p$ is not bounded, then almost nothing is known about the
existence of polynomial time algorithms.

A very important difference between the project described here, using
étale cohomology with coefficients in~$\FF_l$, and the $p$-adic
methods, is that the Galois representations on $\FF_l$-vector spaces
that we obtain are \emph{global} in the sense that they are
representations of the absolute Galois group of the global
field~$\QQ$. The field extensions such as the $K_l=\QQ[x]/(f_l)$
arising from $\Delta$ discussed in the previous section have the
advantage that one can choose to do the required computations over the
complex numbers, approximating~$f_l$, or $p$-adically at some suitable
prime~$p$, or in $\FF_p$ for sufficiently many small~$p$. Also, as we
have said already, being able to compute such field extensions~$K_l$,
that give mod~$l$ information on the Frobenius elements at all
primes~$p\neq l$, is very interesting. On the other hand, the $p$-adic
methods force one to compute with $p$-adic numbers, or, actually,
modulo some sufficiently high power of~$p$, and it gives information
only on the Frobenius at~$p$. The main drawback of the étale
cohomology with $\FF_l$-coefficients seems to be that the degree of
the field extensions as $K_l$ to be dealt with grows exponentially in
the dimension of the cohomology groups; for that reason, we do not
know how to use étale cohomology to compute $\# X(\FF_q)$ for $X$ a
curve of arbitrary genus in a time polynomial in $\log q$ and the
genus of~$X$. Nevertheless, for modular curves, see the end of
Section~\ref{subsec_future}.

\chapter{Modular curves, modular forms, lattices, Galois representations}
\label{chap_mod_curves}

\author{B. Edixhoven}

\bigskip

\bigskip

%author: Bas.
\section{Modular curves} 
\label{sec_modcurves}
As a good reference for getting an overview of the theory of modular
curves and modular forms we recommend the article~\cite{Diamond-Im} by
Fred Diamond and John Im. This reference is quite complete as results
are concerned, and gives good references for the proofs of those
results. Moreover, it is one of the few references that treats the
various approaches to the theory of modular forms, from the classical
analytic theory on the upper half plane to the more modern
representation theory of adelic groups. Another good first
introduction could be the book \cite{Diamond-Shurman}. Let us also
mention that there is a forthcoming book~\cite{Conrad1} by Brian
Conrad, and also the information in the \emph{wikipedia} is getting
more and more detailed.

In this section our aim is just to give the necessary definitions and
results for what we need later (and we need at least to fix our
notation). Readers who want more details, or more conceptual
explanations are encouraged to consult \cite{Diamond-Im}.

\begin{defi}
For $n$ an integer greater than or equal to one we let $\Gamma(n)$ be
the kernel of the surjective morphism of groups
$\SL_2(\ZZ)\to\SL_2(\ZZ/n\ZZ)$ given by reduction of the coefficients
modulo~$n$, and we let $\Gamma_1(n)$ be the inverse image of the
subgroup of $\SL_2(\ZZ/n\ZZ)$ that fixes the element $(1,0)$
of~$(\ZZ/n\ZZ)^2$. Similarly, we let $\Gamma_0(n)$ be the inverse
image of the subgroup of $\SL_2(\ZZ/n\ZZ)$ that fixes the subgroup
$\ZZ/n\ZZ{\cdot}(1,0)$ of~$(\ZZ/n\ZZ)^2$. Hence the elements of
$\Gamma_0(n)$ are the $(\begin{smallmatrix}a & b\\ c & d
\end{smallmatrix})$ of $\SL_2(\ZZ)$ such that $c\equiv 0\mod n$, those
of $\Gamma_1(n)$ are the ones that satisfy the extra conditions
$a\equiv 1\mod n$ and $d\equiv 1\mod n$ and those of $\Gamma(n)$ are
the ones that satisfy the extra condition $b\equiv 0\mod n$.
\end{defi}
The group $\SL_2(\RR)$ acts on the upper half plane $\HH$ by
fractional linear transformations:
$$
\left(\begin{matrix}a & b\\ c & d  \end{matrix}\right)\cdot z = 
\frac{az+b}{cz+d}.
$$
The subgroup $\SL_2(\ZZ)$ of $\SL_2(\RR)$ acts discontinuously in the
sense that for each $z$ in $\HH$ the stabiliser $\SL_2(\ZZ)_z$ is
finite and there is an open neighbourhood $U$ of $z$ such that each
translate $\gamma U$ with $\gamma$ in $\SL_2(\ZZ)$ contains exactly
one element of the orbit $\SL_2(\ZZ){\cdot}z$ and any two translates
$\gamma U$ and $\gamma'U$ with $\gamma$ and $\gamma'$ in $\SL_2(\ZZ)$
are either equal or disjoint. This property implies that the quotient
$\SL_2(\ZZ)\backslash\HH$, equipped with the quotient topology and
with, on each open subset~$U$, the $\SL_2(\ZZ)$-invariant holomorphic
functions on the inverse image of~$U$, is a complex analytic manifold
of dimension one, i.e., each point of the quotient has an open
neighbourhood that is isomorphic to the complex unit disk. Globally,
the well-known $j$-function from $\HH$ to $\CC$ \emph{is} in fact the
quotient map for this action. One way to see this is to associate to
each $z$ in $\HH$ the elliptic curve $E_z:=\CC/(\ZZ+\ZZ z)$, and to
note that for $z$ and $z'$ in $\HH$ the elliptic curves $E_z$ and
$E_{z'}$ are isomorphic if and only if $z$ and $z'$ are in the same
$\SL_2(\ZZ)$-orbit, and to use the fact that two complex elliptic
curves are isomorphic if and only if their $j$-invariants are equal.

The quotient set $\Gamma(n)\backslash\HH$ can be identified with the
set of isomorphism classes of pairs $(E,\phi)$, where $E$ is a complex
elliptic curve and $\phi\colon (\ZZ/n\ZZ)^2\to E[n]$ is an isomorphism
of groups, compatible with the Weil pairing $E[n]\times E[n]\to
\mu_n(\CC)$ and the $\mu_n(\CC)$-valued pairing on $(\ZZ/n\ZZ)^2$ that
sends $((a_1,a_2),(b_1,b_2))$ to $\zeta_n^{a_1b_2-a_2b_1}$, where
$\zeta_n=e^{2\pi i/n}$.

The quotient $\Gamma_0(n)\backslash\HH$ is then identified with the
set of pairs $(E,G)$ where $E$ is a complex elliptic curve, and
$G\subset E$ a subgroup that is isomorphic
to~$\ZZ/n\ZZ$. Equivalently, we may view $\Gamma_0(n)\backslash\HH$ as
the set of isomorphism classes of $E_1\stackrel{\phi}{\to}E_2$, where
$\phi$ is a morphism of complex elliptic curves, and $\ker(\phi)$ is
isomorphic to~$\ZZ/n\ZZ$.

Finally, the quotient $\Gamma_1(n)\backslash\HH$ is then identified with the
set of pairs $(E,P)$ where $E$ is a complex elliptic curve, and $P$ is
a point of order~$n$ of~$E$. Explicitly: to each $z$ in $\HH$
corresponds the pair $(\CC/(\ZZ z+\ZZ),[1/n])$, where $[1/n]$ denotes
the image of $1/n$ in~$\CC/(\ZZ z+\ZZ)$. 

In order to understand that the quotients considered above are in fact
the complex analytic varieties associated with affine complex
algebraic curves, it is necessary (and sufficient!) to show that these
quotients can be compactified to compact Riemann surfaces by adding a
finite number of points, called \emph{the cusps}. As the quotient by
$\SL_2(\ZZ)$ is given by $j\colon \HH\to\CC$, it can be compactified
easily by embedding $\CC$ into~$\PP^1(\CC)$; the point $\infty$ of
$\PP^1(\CC)$ is called the cusp. Another way to view this is to note
that the equivalence relation on $\HH$ given by the action of
$\SL_2(\ZZ)$ identifies two elements $z$ and $z'$ with $\Im(z)>1$ and
$\Im(z')>1$ if and only if $z'=z+n$ for some $n$ in~$\ZZ$; this
follows from the identity, for all $(\begin{smallmatrix}a & b\\ c &
  d\end{smallmatrix})$ in $\SL_2(\RR)$ and $z$ in~$\HH$:
\begin{eqn}\label{eqn_transform_y}
\Im\left(\frac{az+b}{cz+d}\right) = \frac{\Im(z)}{|cz+d|^2}.
\end{eqn}
Indeed, if moreover $c\neq 0$, then:
\begin{eqn}
\frac{\Im(z)}{|cz+d|^2} \leq \frac{\Im(z)}{(\Im(cz))^2} = \frac{1}{c^2\Im(z)}.
\end{eqn}
Hence on the part ``$\Im(z)>1$'' of $\HH$ the equivalence relation
given by $\SL_2(\ZZ)$ is given by the action of $\ZZ$ by
translation. As the quotient for that action is given by the map
$q\colon \HH\to D(0,e^{-2\pi})^*$, $z\mapsto\exp(2\pi iz)$, where
$D(0,e^{-2\pi})^*$ is the open disk of radius $e^{-2\pi}$, centred
at~$0$, and with $0$ removed, we get an open immersion of
$D(0,e^{-2\pi})^*$ into $\SL_2(\ZZ)\backslash\HH$. The
compactification is then obtained by replacing $D(0,e^{-2\pi})^*$ with
$D(0,e^{-2\pi})$, i.e., by adding the centre back into the punctured
disk.

Let us now consider the problem of compactifying the other quotients
above. Let $\Gamma$ be one of the groups considered above, or, in
fact, any subgroup of finite index in~$\SL_2(\ZZ)$. We consider the
morphism $f\colon \Gamma\backslash\HH\to \SL_2(\ZZ)\backslash\HH=\CC$,
and our compactification $\PP^1(\CC)$ of~$\CC$. By construction, $f$
is proper (i.e., the inverse image of a compact subset of $\CC$ is
compact). Also, we know that ramification can only occur at points
with $j$-invariant $0$ or~$1728$. Let $D^*$ be the punctured disk
described above. Then $f\colon f^{-1}D^*\to D^*$ is an unramified
covering of degree $\#\SL_2(\ZZ)/\Gamma$ if $\Gamma$ does not contain
$-1$, and of degree $(\#\SL_2(\ZZ)/\Gamma)/2$ if $-1$ is
in~$\Gamma$. Up to isomorphism, the only connected unramified covering
of degree $n$, with $n\geq1$, of $D^*$ is the map $D_n^*\to D^*$, with
$D_n^*=\{z\in\CC\,|\,0<|z|<e^{-2\pi/n}\}$, sending $z\mapsto z^n$.  It
follows that $f^{-1}D^*$ is, as a covering of $D^*$, a disjoint union
of copies of such $D_n^*\to D^*$. Each $D_n^*$ has the natural
compactification $D_n:=\{z\in\CC\,|\,|z|<e^{-2\pi/n}\}$. We compactify
$\Gamma\backslash\HH$ by adding the origin to each punctured disk in
$f^{-1}D^*$. The points that we have added are called the cusps.  By
construction, the morphism $f\colon \Gamma\backslash\HH\to
\SL_2(\ZZ)\backslash\HH$ extends to the compactifications. It is a
fact that a compact Riemann surface can be embedded into some
projective space, using the theorem of Riemann-Roch, and that the
image of such an embedding is a complex algebraic curve. This means
that our quotients are, canonically, the Riemann surfaces associated
with smooth complex algebraic curves.

\begin{defi}
For $n\geq1$ we define $X(n)$, $X_1(n)$ and $X_0(n)$ to be the proper
smooth complex algebraic curves obtained via the compactifications of
$\Gamma(n)\backslash\HH$, $\Gamma_1(n)\backslash\HH$, and
$\Gamma_0(n)\backslash\HH$, respectively. The affine parts obtained by
removing the cusps are denoted $Y(n)$, $Y_1(n)$ and $Y_0(n)$.
\end{defi}
The next step in the theory is to show that these complex algebraic
curves are naturally defined over certain number fields. Let us start
with the $X_0(n)$ and~$X_1(n)$, which are defined over~$\QQ$. A simple
way to produce a model of $X_0(n)$ over~$\QQ$, i.e., an algebraic
curve $X_0(n)_\QQ$ over $\QQ$ that gives $X_0(n)$ via extension of
scalars via $\QQ\to\CC$, is to use the map:
$$
(j,j')\colon \HH \lto\CC\times\CC, \quad z\mapsto (j(z),j(nz)).
$$
This map factors through the action of $\Gamma_0(n)$, and induces a
map from $X_0(n)$ to $\PP^1\times\PP^1$ that is birational to its
image. This image is a curve in $\PP^1\times\PP^1$, hence the zero
locus of a bi-homogeneous polynomial often denoted~$\Phi_n$, the
minimal polynomial of $j'$ over~$\CC(j)$. One can then check, using
some properties of the $j$-function, that $\Phi_n$ has integer
coefficients. The normalisation of the curve in
$\PP^1_\QQ\times\PP^1_\QQ$ defined by $\Phi_n$ is then the desired
curve~$X_0(n)_\QQ$. As $\Phi_n$ has coefficients in~$\ZZ$, it even
defines a curve in $\PP^1_\ZZ\times\PP^1_\ZZ$ (here, one has to work
with schemes), whose normalisation $X_0(n)_\ZZ$ can be characterised
as a so-called \emph{coarse moduli space}. For this notion, and for
the necessary proofs, the reader is referred
to~\cite[II.8]{Diamond-Im}, to~\cite{Deligne-Rapoport} and
to~\cite{Katz-Mazur}. One consequence of this
statement is that for any algebraically closed field $k$ in which $n$
is invertible, the $k$-points of $Y_0(n)_\ZZ$ (the complement of the
cusps) correspond bijectively to isomorphism classes of
$E_1\stackrel{\phi}{\to}E_2$ where $\phi$ is a morphism of elliptic
curves over~$k$ of which the kernel is cyclic of order~$n$.

The notion of moduli space also gives natural models over $\ZZ[1/n]$
of $X_1(n)$ and $Y_1(n)$. For $n\geq4$ the defining property of
$Y_1(n)_{\ZZ[1/n]}$ is not hard to state. There is an elliptic curve
$\EE$ over $Y_1(n)_{\ZZ[1/n]}$ with a point $\PP$ in
$E(Y_1(n)_{\ZZ[1/n]})$ that has order $n$ in every fibre, such that
\emph{any} pair $(E/S,P)$ with $S$ a $\ZZ[1/n]$-scheme and $P$ in
$E(S)$ of order $n$ in all fibres arises by a \emph{unique} base
change:
\[
\xymatrix{
E \ar[r] \ar[d] & \EE \ar[d]\\
S \ar[r] & Y_1(n)_{\ZZ[1/n]}
}
\]
that is compatible with the sections $P$ and~$\PP$. The pair
$(\EE/Y_1(n)_{\ZZ[1/n]},\PP)$ is therefore called \emph{universal}. 

The moduli interpretation of $X(n)$ is a bit more complicated, because
of the occurrence of the Weil pairing on $E[n]$ that we have seen
above. The curve $X(n)$ has a natural model $X(n)_{\ZZ[1/n,\zeta_n]}$
over $\ZZ[1/n,\zeta_n]$. The complement of the cusps
$Y(n)_{\ZZ[1/n,\zeta_n]}$ then has an elliptic curve $\EE$ over it,
and an isomorphism $\phi$ between the constant group scheme
$(\ZZ/n\ZZ)^2$ and $\EE[n]$ that respects the pairings on each
side. The pair $(\EE/Y(n)_{\ZZ[1/n,\zeta_n]},\phi)$ is universal in
the same sense as above. We warn the reader that the notation $X(n)$
is also used sometimes for the moduli scheme for pairs $(E,\phi)$
where $\phi$ does not necessarily respect the pairings on the two
sides.

For $n$ and $m$ in $\ZZ_{>1}$ that are relatively prime we will
sometimes view $Y_1(nm)_{\ZZ[1/nm]}$ as the moduli space of triples
$(E/S,P_n,P_m)$, where $S$ is a scheme over $\ZZ[1/nm]$, $E/S$ an
elliptic curve, $P_n$ and $P_m$ in $E(S)$ that are everywhere (i.e.,
in every geometric fibre of $E/S$) of orders $n$ and~$m$,
respectively. Indeed, for such a triple, $P_n+P_m$ is everywhere of
order~$nm$, and the inverse construction starting with a point
$P_{nm}$ that is everywhere of order $nm$ is given by multiplying with
the two idempotents of $\ZZ/nm\ZZ$ corresponding to the isomorphism of
rings $\ZZ/nm\ZZ\to\ZZ/n\ZZ\times\ZZ/m\ZZ$.

\section{Modular forms}\label{sec_modforms}
Let us now turn our attention to modular forms. It will be enough for
us to work with modular forms for the congruence
subgroups~$\Gamma_1(n)$. Therefore, we restrict ourselves to that
case.

\begin{defi}\label{def_modform}
Let $n\geq 1$ and $k$ an integer. A (holomorphic) \emph{modular form} for
$\Gamma_1(n)$ is a holomorphic function $f\colon \HH\to\CC$ that
satisfies the following properties:
\begin{enumerate}
\item for all $(\begin{smallmatrix}a & b\\ c & d\end{smallmatrix})\in
\Gamma_1(n)$ and for all $z\in\HH$: 
\[
f((az+b)/(cz+d)) = (cz+d)^k f(z);
\]
\item $f$ is holomorphic at the cusps (see below for an explanation).
\end{enumerate}
A modular form is called a \emph{cuspform} if it vanishes at the cusps.
\end{defi}
We still need to explain the condition that $f$ is holomorphic at the
cusps. In order to do that, we first explain what this means at the
cusp~$\infty$. That cusp is the point that was added to the punctured
disk obtained by taking the quotient of $\HH$ by the unipotent
subgroup $(\begin{smallmatrix}1 & *\\ 0 & 1\end{smallmatrix})$, which
acts on $\HH$ by translations by integers. The coordinate of that disk
is $q$, the map that sends $z$ to $\exp(2\pi iz)$. Therefore, $f$
admits a Laurent series expansion in~$q$:
\begin{eqn}
f = \sum_{n\in\ZZ} a_n(f)q^n, \quad\text{called the $q$-expansion at~$\infty$}.
\end{eqn}
With this notation, $f$ is called holomorphic at~$\infty$ if $a_n(f)$
is zero for all $n<0$, and $f$ is said to vanish at $\infty$ if $a_n(f)$
is zero for all $n\leq 0$. 

To state this condition at the other cusps, we need some description
of the set of cusps. First, we note that $\PP^1(\CC)-\PP^1(\RR)$ is
the same as $\CC-\RR$, and therefore the disjoint union of $\HH$ and
its complex conjugate (which explains, by the way, that $\GL_2(\RR)^+$
acts by fractional linear transformations on~$\HH$). We can then
consider $\HH\cup\PP^1(\QQ)$ inside~$\PP^1(\CC)$, with the
$\SL_2(\ZZ)$ action on it. Then the subgroup $(\begin{smallmatrix}1 &
*\\ 0 & 1\end{smallmatrix})$ stabilises the point $\infty=(1:0)$
of~$\PP^1(\QQ)$, and $\infty$ can be naturally identified with the
origin that we added to the disk $D^*$ above, because $\infty$ is the
unique element of $\PP^1(\QQ)$ that lies in the closure of the inverse
image ``$\Im(z)>1$'' of $D^*$ in~$\HH$. Then, the images of the region
``$\Im(z)>1$'' under the action of elements of $\SL_2(\ZZ)$ correspond
bijectively to the elements of $\PP^1(\QQ)$ (note that $\SL_2(\ZZ)$
acts transitively on $\PP^1(\QQ)=\PP^1(\ZZ)$), and also to the maximal
unipotent subgroups of~$\SL_2(\ZZ)$ (i.e., the subgroups that consist
of elements whose eigenvalues are~$1$). It follows that we can
identify the set of cusps of $X_1(n)$ with
$\Gamma_1(n)\backslash\PP^1(\QQ)$, and that the images of the region
``$\Im(z)>1$'' under $\SL_2(\ZZ)$ give us punctured disks around the
other cusps.  Let $\gamma=(\begin{smallmatrix}a & b\\ c &
d\end{smallmatrix})$ be an element of $\SL_2(\ZZ)$. The conditions of
holomorphy and vanishing at the cusp $\gamma\infty = (a:c)$ are then
given in terms of the $q$-expansion of $z\mapsto (cz+d)^{-k}f(\gamma
z)$ at~$\infty$. The group $\gamma^{-1}\Gamma_1(n)\gamma$ contains the
group $(\begin{smallmatrix}1 & n*\\ 0 & 1\end{smallmatrix})$ (indeed,
$\Gamma_1(n)$ contains $\Gamma(n)$ and that one is normal
in~$\SL_2(\ZZ)$). Therefore, putting $q_n\colon \HH\to\CC$, $z\mapsto
\exp(2\pi iz/n)$, the function $z\mapsto (cz+d)^{-k}f(\gamma z)$ then
has a Laurent series expansion in $q_n$, and one asks that this
Laurent series is a power series (for holomorphy) or a power series
with constant term zero (for vanishing).

The space of modular forms of weight $k$ on $\Gamma_1(n)$ will be
denoted $M_k(\Gamma_1(n))$, and the subspace of cuspforms
by~$S_k(\Gamma_1(n))$. We define $M(\Gamma_1(n))$ to be the direct sum
over the $k$ in~$\ZZ$ of the~$M_k(\Gamma_1(n))$; it is a $\ZZ$-graded
commutative $\CC$-algebra under pointwise multiplication.

\begin{exam}\label{exam_eisenstein}
Some simple examples of modular forms for $\SL_2(\ZZ)$ are given by
Eisenstein series. For each even $k\geq 4$ one has the function~$E_k$:
\[
E_k\colon \HH\to\CC, \quad z\mapsto 
\frac{1}{2\zeta(k)}
\sum_{\stackrel{(n,m)\in\ZZ^2}{(n,m)\neq(0,0)}}\frac{1}{(n+mz)^k}.
\]
The $q$-expansions of these $E_k$ are given by:
\[\label{eqn_eis_series}
E_k = 1-\frac{2k}{B_k}\sum_{n\geq 1} \sigma_{k-1}(n)q^n,
\]
where the $B_k$ are the Bernoulli numbers defined by:
$$
\frac{te^t}{e^t-1} = \sum_{k\geq 0} B_k\frac{t^k}{k!},
$$
and where, as before, $\sigma_r(n)$ denotes the sum of the $r$th powers
of the positive divisors of~$n$. In particular, one has the formulas:
\begin{align*}
E_4 & = 1 + 240\sum_{n\geq 1}\sigma_3(n)q^n, \quad 
E_6 = 1 - 504\sum_{n\geq 1}\sigma_5(n)q^n, \\
\Delta & = \frac{E_4^3-E_6^2}{1728}.
\end{align*}

\begin{rem}\label{rem_bach-charles}
We note that, from a computational point of view, the coefficients
of $q^p$ with $p$ prime of the $E_k$ are very easy to compute, namely, up
to a constant factor they are $1+p^{k-1}$, but that computing the
$\sigma_{k-1}(n)$ for composite $n$ is equivalent to
factoring~$n$. This is a strong indication that, for computing
coefficients $a_n(f)$ of a modular form~$f$, there is a real difference
between the case where $n$ is prime and the case where $n$ is
composite. 

Indeed, Denis Charles and Eric Bach have shown that for $n=pq$ a
product of two distinct primes such that $\tau(n)\neq0$, one can
compute $p$ and $q$ from~$n$, $\tau(n)$ and $\tau(n^2)$ in time
polynomial in~$\log n$; see~\cite{Bach-Charles}.

The argument is very simple: one uses the identities
in~(\ref{eqn_tau_identities}) to compute the rational number
$\tau(p)^2/q^{11}$ and notes that the denominator is of the form $q^r$
with $r$ odd. According to a conjecture by Lehmer, $\tau(n)\neq0$ for
all $n\in\ZZ_{\geq 1}$. See Corollary~\ref{Lehmer}.
\end{rem}

The Eisenstein series $E_4$ and~$E_6$ generate the $\CC$-algebra
$M(\SL_2(\ZZ))$, and are algebraically independent:
\[
M(\SL_2(\ZZ)) = \CC[E_4,E_6].
\]
In particular, we have:
\[
\dim_\CC M_k(\SL_2(\ZZ)) = \#\{(a,b)\in\ZZ_{\geq0}^2\,|\,4a+6b=k\}.
\]
\end{exam}

The space $S_k(\Gamma_1(n))$ can be interpreted as the space of
sections of some holomorphic line bundle $\om^{\otimes k}(-\Cusps)$
on~$X_1(n)$, if $n\geq 5$ (for $n<4$ the action of $\Gamma_1(n)$
on~$\HH$ is not free, and for $n=4$ there is a cusp whose stabiliser
is not unipotent):
\begin{eqn}
S_k(\Gamma_1(n)) = \rH^0(X_1(n),\om^{\otimes k}(-\Cusps)), \quad
\text{if $n\geq 5$.}
\end{eqn}
This implies that the spaces $S_k(\Gamma_1(n))$ are finite
dimensional, and in fact zero if $k\leq 0$ because the line bundle in
question then has negative degree. The restriction to $Y_1(n)$ of the
line bundle giving the weight $k$ forms is given by dividing out the
action of $\Gamma_1(n)$ on $\CC\times\HH$ given by:
\begin{eqn}
(\begin{smallmatrix}a & b\\ c & d\end{smallmatrix})\colon
(x,z)\mapsto \left((cz+d)^k x,\frac{az+b}{cz+d}\right).
\end{eqn}
The extension of this line bundle over the cusps is then given by
decreeing that, at the cusp $\infty$, the constant section $1$ (which
is indeed invariant under the translations $z\mapsto z+n$) is a
generator for the bundle of holomorphic forms, and $q$ times $1$ is a
generator for the bundle of cusp forms.

The moduli interpretation for $Y_1(n)$ can be extended to the
holomorphic line bundles giving the modular forms as follows. Recall
that a point on $Y_1(n)$ is an isomorphism class of a pair $(E,P)$
with $E$ a complex elliptic curve and $P$ a point of order $n$
of~$E$. The complex line at $(E,P)$ of the bundle of forms of weight
$k$ is then $\om_{E}^{\otimes k}$, the $k$th tensor power of the dual
of the tangent space at $0$ of~$E$. In this way, a modular form $f$ of
weight $k$ for $\Gamma_1(n)$ can be described as follows: it is a
function that assigns to each $(E,P)$ an element $f(E,P)$ of
$\om_{E}^{\otimes k}$, varying holomorphically with $(E,P)$, and such
that it has the right property at the cusps (being holomorphic or
vanishing). The function $f$ has to be compatible with isomorphisms:
if $\phi\colon E\to E'$ is an isomorphism, and $\phi(P)=P'$, then
$f(E,P)$ has to be equal to $(\phi^*)^{\otimes k}f(E',P')$. In what
follows we will simply write $\phi^*$ for~$(\phi^*)^{\otimes k}$. 

The fact that $f$ should be holomorphic can be stated by
evaluating it on the family of elliptic curves that we have
over~$\HH$. Recall that to $z$ in $\HH$ we attached the pair
$(\CC/(\ZZ z+\ZZ),[1/n])$. Let us denote $x$ the coordinate of~$\CC$,
then $dx$ is a generator of the cotangent space at $0$ of this
elliptic curve. Then for $f$ a function as above, we can write:
\begin{eqn}
f((\CC/(\ZZ z+\ZZ),[1/n])) = F_f(z) {\cdot} (dx)^{\otimes k}, \quad
F_f\colon\HH\to\CC .
\end{eqn}
The function $F_f$ is then required to be holomorphic. The requirement
that $f$ is compatible with isomorphisms means precisely that $F_f$
transforms under $\Gamma_1(n)$ as in Definition~\ref{def_modform}
above.  The requirement that $f$ vanishes at the cusps is equivalent
to the statement that the Laurent expansions in $q^{1/n}\colon
z\mapsto \exp(2\pi iz/n)$ obtained by evaluating $f$ on all pairs
$(\CC/(\ZZ z+\ZZ),(az+b)/n)$, with $a$ and $b$ in $\ZZ$ such that
$(az+b)/n$ is of order $n$ are in fact power series with constant term
zero.

The spaces $S_k(\Gamma_1(n))$ are equipped with certain operators,
called Hecke operators and diamond operators. These operators arise
from the fact that for every element $\gamma$ of $\GL_2(\QQ)^+$ the
subgroups $\Gamma_1(n)$ and $\gamma\Gamma_1(n)\gamma^{-1}$ are
commensurable, i.e., their intersection has finite index in each of
them.  The diamond operators are then the simplest to describe. For
each $a$ in $(\ZZ/n\ZZ)^\times$, $Y_1(n)$ has the automorphism $\ld a\rd$
given by the property that it sends $(E,P)$ to $(E,aP)$. This action
is then extended on modular forms by:
\begin{eqn}
  (\ld a\rd f)(E,P) = f(E,aP).
\end{eqn}
Similarly, there are Hecke operators $T_m$ on $S_k(\Gamma_1(n))$ for
all integers $m\geq1$, defined by:
\begin{eqn}
  (T_mf)(E,P) = \frac{1}{m}\sum_\phi \phi^* f(E_\phi,\phi(P)),
\end{eqn}
where the sum runs over all quotients $\phi\colon E\to E_\phi$ of
degree $m$ such that $\phi(P)$ is of order~$n$. Intuitively, the
operator $T_m$ is to be understood as a kind of averaging operator
over all possible isogenies of degree~$m$. However, the normalising
factor $1/m$ is not equal to the inverse of the number of such
isogenies. Instead, this factor is there to make the Eichler-Shimura
isomorphism (see~(\ref{eqn_ES_iso})) $T_m$-equivariant.

Of course, each element $f$ of $S_k(\Gamma_1(n))$ is determined by its
$q$-expansion $\sum_{m\geq 1}a_m(f)q^m$ at the cusp~$\infty$. The
action of the Hecke operators can be expressed in terms of these
$q$-expansions (see \cite[(12.4.1)]{Diamond-Im}): 
\begin{eqn} \label{eqn_hecke_q-exp}
a_m(T_rf)=\sum_{\substack{0<d|(r,m)\\(d,n)=1}} d^{k-1}a_{rm/d^2}(\ld d\rd f), 
\end{eqn}
for $f$ in $S_k(\Gamma_1(n))$, $r$ and $m$ positive integers.

From this formula, a lot can be deduced. It can be seen that the $T_r$
commute with each other (but there are better ways to understand
this). The $\ZZ$-algebra generated by the $T_m$ for $m\geq 1$ and the
$\ld a\rd$ for $a$ in $(\ZZ/n\ZZ)^\times$ is in fact generated by the $T_m$
with $m\geq 1$, i.e., one does not need the diamond operators, and
also by the $T_p$ for $p$ prime and the $\ld a\rd$ with $a$
in~$(\ZZ/n\ZZ)^\times$ (see \cite[\S3.5]{Diamond-Im}). The multiplication
rules for the $T_m$ acting on $S_k(\Gamma_1(n))$ can be read off from
the formal identity (\cite[\S3.4]{Diamond-Im}):
\begin{eqn}\label{eqn_dirichlet_hecke}
\sum_{m\geq 1} T_m m^{-s} 
= \prod_p (1-T_pp^{-s}+p^{k-1}\ld p\rd p^{-2s})^{-1},
\end{eqn}
where $\ld p\rd$ is to be interpreted as zero when $p$ divides~$n$.
The fact that the Hecke and diamond operators commute means that they
have common eigenspaces. Taking $m=1$ in (\ref{eqn_hecke_q-exp})
gives:
\begin{eqn}\label{eqn_a1Tnf}
a_1(T_r f) = a_r(f).
\end{eqn}
It follows that if $f$ is a non-zero eigenvector for all $T_r$, then
$a_1(f)\neq 0$, so that we can assume that $a_1(f)=1$. Then, for all
$r\geq 1$, $a_r(f)$ is the eigenvalue for~$T_r$. In particular, this
means that the common eigenspaces for the $T_r$ are one-dimensional,
and automatically eigenspaces for the diamond operators. Eigenforms
with $a_1(f)=1$ are called \emph{normalised eigenforms}.

From (\ref{eqn_dirichlet_hecke}) above it follows that for a
normalised eigenform $f$ one has:
\begin{eqn}
\begin{aligned}
L_f(s) := & \sum_{m\geq 1}a_m(f)m^{-s} \\
= & \prod_p (1-a_p(f)p^{-s}+p^{k-1}\eps_f(p) p^{-2s})^{-1},   
\end{aligned}
\end{eqn}
where $\eps_f\colon(\ZZ/n\ZZ)^\times\to\CC^\times$ is the character via
which the diamond operators act on~$f$, with the convention that
$\eps_f(p)=0$ if $p$ divides~$n$. In particular, the $L$-function of a
modular form has such an Euler product expansion if and only if the
modular form is an eigenform for all Hecke operators.

An element of $S_k(\Gamma_1(n))$ that is a normalised eigenform for
all Hecke operators is called a \emph{newform} if the system of
eigenvalues $a_p(f)$, with $p$ not dividing~$n$, does not occur in a
level strictly smaller than~$n$, i.e., in some $S_k(\Gamma_1(m))$
with~$m<n$ (actually, we will see in a moment that one only needs to
consider the $m$'s dividing~$n$). The set of newforms in
$S_k(\Gamma_1(n))$ will be denoted $S_k(\Gamma_1(n))^\new$.

Now we want to recall briefly how one obtains a basis of
$S_k(\Gamma_1(n))$ in terms of the sets of newforms
$S_k(\Gamma_1(m))^\new$ for $m$ dividing~$n$. For details and
references to proofs, see~\cite[I.6]{Diamond-Im}. First of all, for
each $n$, $S_k(\Gamma_1(n))^\new$ is a linearly independent subset of
$S_k(\Gamma_1(n))$, hence finite. For $m$ dividing $n$ and for $d$
dividing $n/m$, we have a map $B_{n,m,d}\colon X_1(n)\to X_1(m)$,
whose moduli interpretation is that it maps $(E,P)$ to $(E/\ld
(n/d)P\rd,d'P)$, where $dd'=n/m$. For example, this means:
\begin{eqn}\label{eqn_degen_maps}
B_{n,m,d}\colon (\CC/(\ZZ z+\ZZ),1/n) \mapsto (\CC/(\ZZ zd +\ZZ),1/m),
\end{eqn}
which means that the cusp $\infty$ of $X_1(n)$ is mapped to the cusp
$\infty$ of~$X_1(m)$.  Each such map $B_{n,m,d}$ induces by pullback a
map:
\begin{eqn}
B_{n,m,d}^*\colon S_k(\Gamma_1(m)) \to S_k(\Gamma_1(n)).
\end{eqn}
In terms of $q$-expansions at the cusp $\infty$ we have, for $f$ in
$S_k(\Gamma_1(m))$:
\begin{eqn}\label{eqn_q-exp_basis}
B_{n,m,d}^*f = \sum_{r\geq 1}a_r(f)q^{dr},
\end{eqn}
i.e., the effect is just substitution of $q$ by~$q^d$. With these
definitions, we can describe a basis for $S_k(\Gamma_1(n))$:
\begin{eqn}\label{eqn_basis}
\coprod_{m|n}\coprod_{d|(n/m)} B_{n,m,d}^*S_k(\Gamma_1(m))^\new \quad
\text{is a basis of $S_k(\Gamma_1(n))$.}
\end{eqn}
In the case where $\Gamma_1(n)$ is replaced by $\Gamma_0(n)$, this kind
of basis is due to Atkin and Lehner. 

In the sequel, we will also make use of a (hermitian) inner product on
the $S_k(\Gamma_1(n))$: the Petersson scalar product. It is defined as
follows. For $f$ and $g$ in $S_k(\Gamma_1(n))$, viewed as functions on
$\HH$ as in Definition~\ref{def_modform} one has:
\begin{eqn}\label{eqn_pet_scal_prod}
\ld f,g\rd = 
\int_{\Gamma_1(n)\backslash\HH} f(z)\ol{g(z)}\,y^k\,\frac{dxdy}{y^2},
\end{eqn}
where the integral over $\Gamma_1(n)\backslash\HH$ means that one can
perform it over any fundamental domain. Indeed,
formula~\ref{eqn_transform_y} shows that the function $z\mapsto
f(z)\ol{g(z)}y^k$ is invariant under~$\Gamma_1(n)$. 

We also want to explain the definition of $\ld f,g\rd$ in terms of the
moduli interpretation of~$S_k(\Gamma_1(n))$, if $k\geq 2$. For
simplicity, let us suppose $n\geq 5$ now. Then $S_k(\Gamma_1(n))$ is
the space of global sections of $\om^{\otimes k}(-\Cusps)$
on~$X_1(n)$. Now we let $\Omega^1:=\Omega^1_{X_1(n)}$ denote the line
bundle of holomorphic differentials on~$X_1(n)$. Then there is an
isomorphism, named after Kodaira and Spencer:
\begin{eqn}
\mathrm{KS}\colon 
\om^{\otimes 2}(-\Cusps) \stackrel{\sim}{\lto} \Omega^1, 
\quad\text{Kodaira-Spencer isomorphism.}
\end{eqn}
Explicitly, for $f$ in $S_2(\Gamma_1(n))$, viewed as a
$\Gamma_1(n)$-invariant section of $\om^{\otimes 2}$ for the family of
elliptic curves over $\HH$ whose fibre at $z$ is $\CC/(\ZZ z+\ZZ)$ we
have:
\begin{eqn}\label{eqn_KS_isom}
\mathrm{KS}\colon f(dx)^{\otimes 2} \mapsto (2\pi i)^{-2} f \,\frac{dq}{q}.
\end{eqn}
Equivalently, for this family of elliptic curves, the Kodaira-Spencer
isomorphism sends $(dx)^{\otimes 2}$ to $(2\pi i)^{-2}(dq)/q$. Note
that indeed $(dx)^{\otimes 2}$ and $(dq)/q$ transform in the same way
under the action of~$\SL_2(\RR)$. We note that without $f$ being
required to vanish at the cusps, $\mathrm{KS}(f)$ could have poles of
order one at the cusps. The factor $(2\pi i)^{-2}$ is to make the
isomorphism compatible with the coordinates $t=\exp(2\pi ix)$ on
$\CC^\times/\ld\exp(2\pi i z)\rd$ (which is another way to write
$\CC/(\ZZ z+\ZZ)$), and the coordinate $q=\exp(2\pi i z)$ on the unit
disk. In those coordinates, that have a meaning ``over~$\ZZ$'', which
means that formulas relating them are power series (or Laurent series)
with integer coefficients, $\mathrm{KS}$ sends $((dt)/t)^{\otimes 2}$
to~$(dq)/q$.

For every complex elliptic curve, the one dimensional complex vector
space $\om_E$ has the inner product given by:
\begin{eqn}
\ld\alpha,\beta\rd = \frac{i}{2}\int_E \alpha\,\ol{\beta},
\end{eqn}
where we interpret $\alpha$ and $\beta$ as translation invariant
differential forms on~$E$. The factor $i/2$ is explained by the fact
that, for $z=x+iy$, one has $dx\,dy = (i/2)dz\,d\ol{z}$. Applying this
to the family of elliptic curves $\CC/(\ZZ z+\ZZ)$ over~$\HH$ gives an
inner product on the line bundle $\om$ on~$\HH$, and also on the line
bundle $\om$ on~$Y_1(n)$ (recall that we are supposing that $n\geq
5$). Taking tensor powers and duals, this induces inner products on
$\om^{\otimes k}$ for all~$k$. The Kodaira-Spencer
isomorphism~(\ref{eqn_KS_isom}) gives isomorphisms:
\begin{eqn}
\mathrm{KS}\colon \om^{\otimes k}(-\Cusps) \stackrel{\sim}{\lto}
\Omega^1\otimes\om^{\otimes (k-2)}.
\end{eqn}
For $f$ and $g$ in $S_k(\Gamma_1(n))$, now viewed as sections of
$\om^{\otimes k}(-\Cusps)$ over $X_1(n)$, one has:
\begin{eqn}
\ld f,g\rd = \frac{i}{2}\int_{X_1(n)}\ld
\mathrm{KS}(f),\mathrm{KS}(g)\rd,
\end{eqn}
where the inner product on the left hand side is the Petersson scalar
product~(\ref{eqn_pet_scal_prod}), and where for two local sections
$\omega_1\otimes\alpha_1^{\otimes (k-2)}$ and
$\omega_2\otimes\alpha_2^{\otimes (k-2)}$ of
$\Omega^1\otimes\om^{\otimes (k-2)}$ we have defined:
\begin{eqn}
\ld \omega_1\otimes\alpha_1^{\otimes (k-2)},
\omega_2\otimes\alpha_2^{\otimes (k-2)}\rd = 
\ld \alpha_1,\alpha_2\rd^{k-2}\omega_1\ol{\omega_2}.
\end{eqn}

The operators $T_m$ on $S_k(\Gamma_1(n))$ with $m$ relatively prime
to~$n$ are normal: they commute with their adjoint. As a consequence,
distinct newforms in $S_k(\Gamma_1(n))$ are orthogonal to each
other. On the other hand, the basis~(\ref{eqn_basis}) above of
$S_k(\Gamma_1(n))$ is not orthogonal if it consists of more than only
newforms.

\section{Lattices and modular forms}\label{sec_lat-forms}

Before we move on to Galois representations attached to modular
forms, we briefly discuss the relation between modular forms and
lattices. 

Let us consider a free $\ZZ$-module $L$ of finite rank~$n$, equipped
with a positive definite symmetric bilinear form $b\colon L\times
L\to\ZZ$. Then $L_\RR:=\RR\otimes L$ is an $\RR$-vector space of
dimension~$n$ on which $b$ gives an inner product, and hence $L$ is a
lattice in the euclidean space~$L_\RR$. For $m$ in~$\ZZ$ the
\emph{representation numbers} of $(L,b)$ are defined as:
\begin{eqn}
r_L(m) = r_{L,b}(m) := \#\{x\in L \;|\; b(x,x)=m \}.
\end{eqn}
In this situation, one considers the \emph{theta-function} attached
to~$(L,b)$:
\begin{eqn}\label{eqn_theta_defi}
\theta_L = \theta_{L,b} = \sum_{x\in L} q^{b(x,x)/2} = 
\sum_{m\geq0}r_L(m)q^{m/2}, \quad \HH\to\CC,
\end{eqn}
where $q^{1/2}\colon \HH\to\CC$ is the function 
$q^{1/2}\colon z\mapsto \exp(\pi iz)$.

If $(L,b)$ is the orthogonal direct sum of $(L_1,b_1)$ and
$(L_2,b_2)$ then we have:
\begin{eqn}\label{eqn_theta_of_sum}
\theta_L = \sum_{x\in L}q^{b(x,x)/2}=
\sum_{(x_1,x_2)\in L_1\times L_2} q^{b_1(x_1,x_1)/2+b_2(x_2,x_2)/2}
= \theta_{L_1}\cdot\theta_{L_2}.
\end{eqn}

The \emph{discriminant} of~$b$ is $\det(B)$, where $B$ is the matrix
of $b$ with respect to some basis of~$L$ (indeed, this determinant
does not depend on the choice of basis); we denote it by~$\discr(b)$.

We define the positive integer $N_L$ to be the exponent of the
cokernel of the map $\phi_b\colon L\to L^\vee$ given by~$b$, or,
equivalently, to be the denominator of~$B^{-1}$, where $B$ is the
matrix of~$b$ with respect to some basis $e$ of~$L$. The map
$N_L\phi_b^{-1}\colon L^\vee_\QQ\to L_\QQ$ restricts to a map
$N_L\phi_b^{-1}\colon L^\vee\to L$. Viewing $L$ as $(L^\vee)^\vee$ in
the usual way, this gives a positive definite symmetric bilinear form
$b'\colon L^\vee\times L^\vee\to\ZZ$. The matrix of this form with
respect to the basis $e^\vee$ dual to~$e$ is~$N_LB^{-1}$. Applying
this same construction to~$b'$ gives a $b''\colon L\times L\to\ZZ$
that is not necessarily equal to~$b$: one has $b=mb''$, with
$m\in\ZZ_{>0}$ and $b''$ primitive (i.e., the $\ZZ$-linear map
$L\otimes L\to\ZZ$ induced by~$b''$ is surjective). Poisson's
summation formula gives the following functional equation;
see~\cite[VII, \S6, Prop.~16]{Serre5}.
\begin{thm}\label{thm_theta_func_eqn}
Let $L$ be a free $\ZZ$-module, of finite rank~$n$, equipped with a
positive definite symmetric bilinear form $b\colon L\times
L\to\ZZ$. We have, with the notation as above, for all $z\in\HH$:
\[
\theta_{L,b}(-1/N_Lz) = 
\frac{(-N_Liz)^{n/2}}{\discr(b)^{1/2}} \theta_{L^\vee,b'}(z),
\]
where the square root of $-N_Liz$ is holomorphic in $z$ and positive
for $z\in \RR i$.
\end{thm}

The form~$b$ is called \emph{even} if $b(x,x)$ is even for all~$x$
in~$L$. Equivalently, $b$ is even if and only if the matrix $B$ of $b$
with respect to some basis of~$L$ has only even numbers on the
diagonal.

The form $b$ is called \emph{unimodular} if $\phi_b\colon L\to L^\vee$
is an isomorphism, or, equivalently, if $N_L=1$. In this case,
$\phi_b$ is an isomorphism from $(L,b)$ to $(L^\vee,b')$.

With this terminology, one has the following result,
see~\cite[Cor.~4.9.5]{Miyake1}, the proof of which has as main
ingredient the functional equation of Theorem~\ref{thm_theta_func_eqn}.
\begin{thm}\label{thm_theta-modular}
Let $L$ be a free $\ZZ$-module of finite rank~$n$, equipped with a
positive definite symmetric bilinear form $b\colon L\times L\to
\ZZ$. Assume that $n$ is even. Let $N_L$ be as defined above, and let
$\chi_L$ be the character given by:
\[
\chi_L\colon (\ZZ/N_L\ZZ)^\times\to\CC^\times, \quad 
(a\mod N_L)\mapsto \left(\frac{(-1)^{n/2}\discr(b)}{a}\right), 
\]
where the fraction denotes the Kronecker symbol.
\begin{enumerate}
\item The function $z\mapsto\theta_L(2z)$ is a (non-cuspidal) modular
form on $\Gamma_1(4N_L)$ of weight~$n/2$ and with character~$\chi_L$.
\item If $b$ is even then the function $\theta_L$ is a modular form on
$\Gamma_1(2N_L)$ of weight~$n/2$ and with character~$\chi_L$.
\item If both $b$ and $b'$ (see above for its definition) are even,
then the function $\theta_L$ is a modular form on~$\Gamma_1(N_L)$ of
weight~$n/2$ and with character~$\chi_L$.
\end{enumerate}
\end{thm}
This theorem says nothing about the case where $n$ is odd. In that
case, $\theta_L$ is a modular form of half-integral weight~$n/2$;
see~\cite[Cor.~4.9.7]{Miyake1}. For even unimodular forms, we have the
following corollary of Theorem~\ref{thm_theta-modular}.
\begin{cor}\label{cor_theta_even_unimod}
Let $L$ be a free $\ZZ$-module of finite rank~$n$, and let $b\colon
L\times L\to \ZZ$ be bilinear, symmetric, positive definite, even and
unimodular. Then $n$ is even, and $\theta_L$ is a modular form on
$\SL_2(\ZZ)$ of weight~$n/2$.
\end{cor}
\begin{proof}
As $b$ is unimodular, we have $N_L=1$. The fact that $n$ is even
follows from the fact that $b$ induces a non-degenerate alternating
bilinear form on $\FF_2\otimes L$. As $\phi_b\colon L\to L^\vee$ is an
isomorphism between $(L,b)$ and $(L^\vee,b')$, we have that $b'$ is
even as well. Theorem~\ref{thm_theta-modular} gives the conclusion.
\end{proof}
\begin{rem}
In fact, the rank of an even unimodular lattice is a multiple
of~$8$. This follows directly from Theorems~\ref{thm_theta_func_eqn}
and~\ref{thm_theta-modular} (see~\cite[Cor.~4.9.6]{Miyake1},
or~\cite[VII, \S6, Thm.~8]{Serre5}).
\end{rem}
Let us consider some examples.
\begin{exam}\label{exam_sum_of_n_squares}
For $n\in\ZZ_{\geq0}$ we consider $\ZZ^n$ with its standard inner
product. For $m$ in~$\ZZ$ we have:
\[
r_{\ZZ^n}(m) = \#\{x\in\ZZ^n\,|\,x_1^2+\cdots+x_n^2=m\},
\]
the number of ways in which $m$ can be written as a sum of $n$ squares
of integers. Theorem~\ref{thm_theta-modular} tells us that for even
$n$ the theta function $z\mapsto\theta_{\ZZ^n}(2z)$ is a modular form
on $\Gamma_1(4)$ of weight~$n/2$. According to
(\ref{eqn_theta_of_sum}) we have $\theta_{\ZZ^n} = \theta_\ZZ^n$, and
so all the functions $z\mapsto\theta_{\ZZ^n}(2z)$ are powers of the
modular form $\sum_{m\in\ZZ}q^{m^2}$ of weight~$1/2$ on~$\Gamma_1(4)$.
\end{exam}

\begin{exam}\label{exam_E8}
We consider the E8-lattice, i.e., $E8:=\ZZ^8$ equipped with the inner
product given by the Dynkin diagram~$E8$ with numbered vertices:
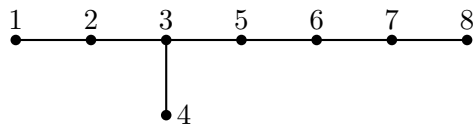
\[
\begin{tikzpicture}[thick]
\draw (0,0) -- (6,0);
\draw (2,0) -- (2,-1);
\fill (0,0) circle (0.07cm) node[anchor=south]{$1$}; 
\fill (1,0) circle (0.07cm) node[anchor=south]{$2$}; 
\fill (2,0) circle (0.07cm) node[anchor=south]{$3$}; 
\fill (2,-1) circle (0.07cm) node[anchor=west]{$4$}; 
\fill (3,0) circle (0.07cm) node[anchor=south]{$5$}; 
\fill (4,0) circle (0.07cm) node[anchor=south]{$6$}; 
\fill (5,0) circle (0.07cm) node[anchor=south]{$7$}; 
\fill (6,0) circle (0.07cm) node[anchor=south]{$8$}; 
\end{tikzpicture}
\]
that is, whose matrix with respect to the standard basis is:
\[
\left(
\begin{matrix}
2  & -1 &    &    &    &    &    &   \\
-1 & 2  & -1 &    &    &    &    &   \\
   & -1 & 2  & -1 & -1 &    &    &   \\
   &    & -1 &  2 &    &    &    &   \\
   &    & -1 &    & 2  & -1 &    &   \\
   &    &    &    & -1 & 2  & -1 &   \\
   &    &    &    &    & -1 & 2  & -1\\
   &    &    &    &    &    & -1 & 2
\end{matrix}
\right).
\]
The lattice $E8$ is unimodular and even, hence, by
Corollary~\ref{cor_theta_even_unimod}, $\theta_{E8}$ is a modular form on
$\SL_2(\ZZ)$ of weight~$4$, i.e., $\theta_{E8}$ is an element
of~$M_4(\SL_2(\ZZ))$. The dimension space of $M_4(\SL_2(\ZZ))$ is one,
with the Eisenstein series $E_4=1+240\sum_{n\geq 1}\sigma_3(n)q^n$ as
basis. Therefore, $\theta_{E8}$ is a constant times~$E_4$. Comparing
constant terms, we get:
\[
\theta_{E8} = E_4 = 1+240\sum_{n\geq 1}\sigma_3(n)q^n.
\]
\end{exam}

\begin{exam}\label{exam_Leech}
Let $L$ be the \emph{Leech lattice}. This lattice, which is of
rank~$24$, even and unimodular, is named after John Leech,
see~\cite{Leech2} and~\cite{Leech1}. Apparently, it had already been
discovered by Ernst Witt in 1940 (unpublished, see~\cite{Witt1}).
John Horton Conway showed in~\cite{Conway1} that $L$ is the only
non-zero even unimodular lattice of rank less than~$32$ with
$r_L(2)=0$; this also follows from Hans-Volker Niemeier's
classification of even unimodular lattices of rank~$24$
in~\cite{Niemeier1}.  According to Henry Cohn and Abhinav
Kumar~\cite{Cohn-Kumar1}, the Leech lattice gives the densest lattice
sphere packing in dimension~$24$.

Theorem~\ref{thm_theta-modular} above shows that the theta function
$\theta_L$ of the Leech lattice is a modular form of level~$1$ and
weight~$12$. The space of such modular forms is two-dimensional, with
basis the Eisenstein series~$E_{12}$ and the discriminant
form~$\Delta$, where:
\[
E_{12} = 1 + \frac{65520}{691}\sum_{m\geq 1}\sigma_{11}(m)q^m.
\]
Hence $\theta_L$ is a linear combination of $E_{12}$
and~$\Delta$. Comparing the coefficients of $q^m$ for $m=0$ and $m=1$
gives:
\[
\theta_L =  E_{12} - \frac{65520}{691}\Delta.
\]
\end{exam}

\section{Galois representations attached to eigenforms}\label{sec_galreps}
The aim of this section is to describe the construction of the Galois
representations attached to modular forms, that came up in the case
of $\Delta$ in Section~\ref{sec_hist_cong_tau}. Before giving the
construction, let us state the result, which is due, for $k=2$, to
Eichler and Shimura~\cite{Shimura1}, to Deligne~\cite{Deligne1} for
$k>2$, and to Deligne and Serre~\cite{Deligne-Serre} for $k=1$. See
Section~12.5 in~\cite{Diamond-Im}. A long account of the construction
in the case $k\geq 2$ will be given in the book~\cite{Conrad1}.
\begin{thm}\label{thm_exist_gal_rep}
Let $f$ be a normalised newform, let $n$ be its level, let $k$ be its
weight, and let $\eps\colon(\ZZ/n\ZZ)^\times\to\CC^\times$ be its
character. Then the subfield $K$ of $\CC$ generated over $\QQ$ by the
$a_n(f)$, $n\geq 1$, and the image of~$\eps$ is finite over~$\QQ$. For
every prime number $l$ and for any embedding $\lambda$ of $K$ into
$\Qbar_l$, there is a continuous two-dimensional representation
$V_\lambda$ over $\Qbar_l$ of $\Gal(\Qbar/\QQ)$ that is unramified
outside $nl$ and such that for each prime number $p$ not dividing $nl$
the characteristic polynomial of the Frobenius at~$p$ acting
on~$V_\lambda$ equals:
$$
\det(1-x\Frob_p,V_\lambda) = 1 - a_p(f)x + \eps(p)p^{k-1}x^2.
$$
\end{thm}
For $k\geq 2$ the representations $V_\lambda$ can be found in the
$l$-adic étale cohomology in degree $k-1$ of some variety of dimension
$k-1$, or in the cohomology in degree one of some sheaf on a curve, as
we will describe below. The determinant of the action of
$\Gal(\Qbar/\QQ)$ on~$V_\lambda$ is easily described. We let
$\chi_l\colon \Gal(\Qbar/\QQ)\to\ZZ_l^\times$ be the $l$-adic
cyclotomic character defined by $\sigma(z)=z^{\chi_l(\sigma)}$, for
all $\sigma$ in $\Gal(\Qbar/\QQ)$ and all $z$ in $\Qbar^\times$ of
$l$-power order. We let $\eps\colon\Gal(\Qbar/\QQ)\to K^\times$ be the
composition of the character $\eps\colon(\ZZ/n\ZZ)^\times\to K^\times$
with the mod~$n$ cyclotomic character
$\Gal(\Qbar/\QQ)\to(\ZZ/n\ZZ)^\times$ given by the action of
$\Gal(\Qbar/\QQ)$ on~$\mu_n(\Qbar)$. With these definitions, the
determinant of the action of $\Gal(\Qbar/\QQ)$ on~$V_\lambda$ is given
by the character~$\eps\chi_l^{k-1}$. As the image of $\Gal(\Qbar/\QQ)$
under the determinant of~$V_\lambda$ is infinite, its image in
$\GL(V_\lambda)$ is infinite. 

On the other hand, for $k=1$, the image of $\Gal(\Qbar/\QQ)$ in
$\GL(V_\lambda)$ is finite, and in fact all these representations when
$\lambda$ varies can be realised over some fixed finite extension
of~$\QQ$.  The proof of Theorem~\ref{thm_exist_gal_rep} by Deligne and
Serre in the case $k=1$ is quite different from the case $k\geq 2$:
the reductions to finite coefficient fields (see
Section~\ref{sec_red_to_jac}) can be constructed via congruences to
forms of weight~$2$, and then it is shown that these representations
can be lifted to characteristic zero. No direct construction of the
characteristic zero Galois representations for forms of weight one is
known. We remark that in the case $k=2$ the $V_\lambda$ occur in the
first degree étale cohomology with constant coefficients~$\QQ_l$ of
modular curves, hence can be constructed from $l$-power torsion points
of Jacobians of modular curves (in fact, of the modular
curve~$X_1(n)$).

The representation $V_\lambda$ is irreducible by a theorem of Ribet,
see Theorem~2.3 of~\cite{Ribet3}, and hence it is characterised by its
trace. As the Frobenius conjugacy classes at the primes not dividing
$nl$ are dense by Chebotarev's theorem, the representation $V_\lambda$
is unique up to isomorphism. Non-cuspidal eigenforms lead to Galois
representations that are reducible; as our interest lies in going
beyond class field theory, we do not discuss this case.

Let us now start the description of the construction, by Deligne, of
the representation $V_\lambda$ as in Theorem~\ref{thm_exist_gal_rep}
above in the case where $k\geq 2$. First, if $n< 5$, we replace $n$ by
say $5n$ and $f$ by a normalised Hecke eigenform in the 2-dimensional
$\CC$-vector space generated by $f(q)$ and~$f(q^5)$. Then $f$ is no
longer a newform, but it is an eigenform, which will be good enough,
and as $n\geq5$ we can view it as a section of the line bundle
$\om^{\otimes k}(-\Cusps)$ on the smooth complex projective
curve~$X_1(n)$. The eigenvalues at primes other than $5$ have not been
changed by this operation. As one can compute from the formulas in the
previous section, the two possible eigenvalues for $T_5$ on the space
generated by $f(q)$ and $f(q^5)$ are the two roots of the polynomial
$x^2-a_5(f)x + \eps(5)5^{k-1}$, i.e., the two eigenvalues of the
Frobenius element at $5$ attached to $f$ if $\lambda$ does not
divide~$5$. For a detailed computation for this, see Section~4 of
\cite{Coleman-Edixhoven}; that article also explains why one should
expect the two eigenvalues always to be distinct, and that this is a
theorem if $k=2$.

On $Y_1(n)$, we have a universal family $(\EE/Y_1(n),\PP)$ of elliptic
curves with a given point of order~$n$. Taking fibre-wise the
cohomology $\rH^1(\EE_s,\ZZ)$ gives us a locally constant sheaf on
$Y_1(n)$, denoted $\rR^1p_*\ZZ_\EE$ because it is the first higher
direct image of the constant sheaf $\ZZ_\EE$ on $\EE$ via the morphism
$p\colon \EE\to Y_1(n)$. The stalks of the locally constant sheaf
$\rR^1p_*\ZZ_\EE$ on $Y_1(n)$ are free $\ZZ$-modules of rank~$2$. More
concretely, the sheaf $\rR^1p_*\ZZ_\EE$ is obtained from the constant
sheaf $\ZZ^2$ on $\HH$ by dividing out the $\Gamma_1(n)$-action given
by:
\begin{eqn}
\gamma{\cdot}\left(\binom{n}{m},\tau\right) =
\left(\gamma{\cdot}\binom{n}{m},\gamma{\cdot}\tau\right) = 
\left(\binom{an+bm}{cn+dm},\frac{a\tau+b}{c\tau+d}\right),
\end{eqn}
where $\gamma=(\begin{smallmatrix} a & b\\ c & d
\end{smallmatrix})$. 

We will also use other locally constant sheaves on $Y_1(n)$ that are
obtained from $\rR^1p_*\ZZ_\EE$ by tensor constructions. The
classification of the irreducible representations of the algebraic
group $\GL_2$ over~$\QQ$ implies that these tensor constructions are
finite direct sums in which each term is a symmetric power of
$\rR^1p_*\ZZ_\EE$, tensored with a power of the determinant of
$\rR^1p_*\ZZ_\EE$. We define:
\begin{eqn}
\calF_k := \Sym^{k-2}(\rR^1p_*\ZZ_\EE),
\end{eqn}
where $\Sym^{k-2}$ denotes the operation of taking the $(k-2)$th
symmetric power. The sheaf $\calF_k$ is then obtained by dividing
out the $\Gamma_1(n)$-action on the constant sheaf $\Sym^{k-2}(\ZZ^2)$
on~$\HH$. It is useful to view $\ZZ^2$ as the $\ZZ$-submodule $\ZZ
x\oplus\ZZ y$ of the polynomial ring~$\ZZ[x,y]$. The grading $\ZZ[x,y]
= \oplus_i \ZZ[x,y]_i$ by the degree then gives the symmetric powers
of $\ZZ x\oplus\ZZ y$:
\begin{eqn}
\Sym^{k-2}(\ZZ^2) = \ZZ[x,y]_{k-2} = \bigoplus_{i+j=k-2} \ZZ x^iy^j.
\end{eqn}
We extend the sheaf $\calF_k$ to $X_1(n)$ by taking the direct image
via the open immersion $j\colon Y_1(n)\to X_1(n)$; this gives us
$j_*\calF_k$ on $X_1(n)$, again denoted~$\calF_k$. Outside the cusps,
$\calF_k$ is locally constant, with stalks free of rank $k-1$ as
$\ZZ$-modules. At the cusps, the stalks of $\calF_k$ are free of rank
one. At the cusp $\infty$ this follows from the fact that the
subring of invariants of $\ZZ[x,y]$ for the action of 
$(\begin{smallmatrix}1 & * \\ 0 & 1\end{smallmatrix})$
is~$\ZZ[x]$. At the other cusps it then follows by conjugating with a
suitable element of~$\SL_2(\ZZ)$. We note that for $k=2$ the sheaf
$\calF_k$ is the constant sheaf $\ZZ$ on~$X_1(n)$.

The \emph{Eichler-Shimura isomorphism} gives a relation between
modular forms and the cohomology of~$\calF_k$. One way to view this,
due to Deligne, is in terms of Hodge structures. More precisely, the
$\CC$-vector space $\CC\otimes \rH^1(X_1(n),\calF_k)$ carries a Hodge
decomposition:
\begin{eqn}\label{eqn_ES_iso}
\CC\otimes \rH^1(X_1(n),\calF_k) \isomlto
S_k(\Gamma_1(n))\oplus\ol{S_k(\Gamma_1(n))},
\end{eqn}
where the two terms on the right are of type $(k-1,0)$ and $(0,k-1)$,
respectively. The complex conjugation on the second term on the right
comes from the complex conjugation on the tensor factor $\CC$ on the
left.  A good reference for this decomposition and its properties
is~\cite{Bayer-Neukirch}; we will not go into details here. For an
account using group cohomology we refer to Section~12.2
of~\cite{Diamond-Im}. For $k=2$ all of this is quite easy. Via the
Kodaira-Spencer isomorphism (\ref{eqn_KS_isom}) it then is the
decomposition:
\begin{eqn}
\rH^1(X_1(n),\CC) =  \rH^0(X_1(n),\Omega^1)\oplus\ol{\rH^0(X_1(n),\Omega^1)}.
\end{eqn}
We should mention that instead of working with the sheaf $\calF_k$ on
the curve $X_1(n)$, one can also work with a constant sheaf on a
$k-1$-dimensional variety. As before, we let $(\EE,\PP)$ denote the
universal object over~$Y_1(n)$. Then we let $\EE^{k-2}$ denote the
$k-2$-fold fibre power of $\EE$ over~$Y_1(n)$; these are the simplest
cases of so-called Kuga-Sato varieties. The graded commutative algebra
structure on cohomology gives, for $s$ in~$Y_1(n)$, a map, equivariant
for the action of the symmetric group~$S_{k-2}$:
\[
\ZZ(\eps)\otimes\rH^1(\EE_s,\ZZ)\otimes\cdots\otimes\rH^1(\EE_s,\ZZ) 
\lto \rH^{k-2}(\EE^{k-2}_s,\ZZ),
\]
where $\ZZ(\eps)$ denotes the sign representation. Twisting this map
by $\ZZ(\eps)$ and taking co-invariants gives a map:
\begin{eqn}
\calF_{k,s} = \Sym^{k-2}(\rH^1(\EE_s,\ZZ)) \lto 
\rH^{k-2}(\EE^{k-2}_s,\ZZ)_\eps,
\end{eqn}
where the subscript $\eps$ means the largest quotient on which
$S_{k-2}$ acts via the sign representation.  In view of the Leray
spectral sequence for the cohomology $\rH(\EE^{k-2},\ZZ)$ of
$\EE^{k-2}$ in terms of the cohomology of the higher derived direct
images $\rH(Y_1(n),\rR p_*\ZZ_{\EE^{k-2}})$ it is then not so
surprising that $S_k(\Gamma_1(n))$ can be identified with a piece of
$\rH^{k-1}(\ol{\EE^{k-2}},\CC)$, where $\ol{\EE^{k-2}}$ is a certain
smooth projective model of~$\EE^{k-2}$ over~$X_1(n)$. Some details for
this can be found in \cite{Deligne1}, and more of them
in~\cite{Scholl1}, and probably still more in~\cite{Conrad1}. A very
explicit way to describe this identification is the map:
\begin{eqn}
\begin{aligned}
S_k(\Gamma_1(n)) &\lto \rH^{k-1}(\ol{\EE^{k-2}},\CC),\\
f & \mapsto (2\pi i)^{k-1} f d\tau\,dz_1\cdots dz_{k-2},
\end{aligned}
\end{eqn}
where $\tau$ is the coordinate on~$\HH$, and the $z_j$ are the
coordinates on the copies of $\CC$ using $\EE_\tau =
\CC/(\ZZ\tau+\ZZ)$. It is indeed easy to verify that the differential
form on the right is invariant under the actions of $\ZZ^{2(k-2)}$ and
$\SL_2(\ZZ)$, precisely because $f$ is a modular form of weight $k$
for~$\Gamma_1(n)$. The claim (proved in the references above) is that
it extends without poles over~$\ol{\EE^{k-2}}$. As it is a holomorphic
form of top-degree, it is automatically closed, and hence defines a
class in the de Rham cohomology of~$\ol{\EE^{k-2}}$, hence in
$\rH^{k-1}(\ol{\EE^{k-2}},\CC)$.

There are natural Hecke correspondences on
$\CC\otimes\rH^1(X_1(n),\calF_k)$ and on
$\rH^{k-1}(\ol{\EE^{k-2}},\CC)$, and the identification of
$S_k(\Gamma_1(n))$ as a piece of these cohomology groups is compatible
with these correspondences. Let now $f$ be our eigenform in
$S_k(\Gamma_1(n))$ as above. Then the Hecke eigenspace in
$\CC\otimes\rH^1(X_1(n),\calF_k)$ with the eigenvalues $a_m(f)$ for
$T_m$ is two-dimensional: the sum of the one-dimensional subspace $\CC
f$ in $S_k(\Gamma_1(n))$ and the one-dimensional subspace $\CC\ol{f'}$
in $\ol{S_k(\Gamma_1(n))}$, where $f'=\sum_{m\geq 1}\ol{a_m(f)}q^m$,
the Galois conjugate of $f$ obtained by letting complex conjugation
act on the coefficients of~$f$. This element $f'$ has eigenvalue
$\ol{a_m(f)}$ for $T_m$, hence $\ol{f'}$ has eigenvalue $a_m(f)$
again. The $(k-1)$-form corresponding to $\ol{f'}$ is 
$\ol{f'} d\ol{\tau}\,d\ol{z_1}\cdots d\ol{z_{k-2}}$, indeed a form of
type $(0,k-1)$. 

We let $\TT(n,k)$ or just $\TT$ denote the $\ZZ$-algebra
in $\End_\CC(S_k(\Gamma_1(n)))$ generated by the $T_m$ ($m\geq 1$) and
the~$\ld a\rd$ ($a$ in~$(\ZZ/n\ZZ)^\times$). The fact that the
Eichler-Shimura isomorphism~(\ref{eqn_ES_iso}) is equivariant for the
Hecke correspondences acting on both sides implies that the image of
$\rH^1(X_1(n),\calF_k)$ in $\CC\otimes\rH^1(X_1(n),\calF_k)$ is a
faithful $\TT(n,k)$-module. As this image is free of finite rank as
$\ZZ$-module, $\TT(n,k)$ is free of finite rank as $\ZZ$-module.

Let us for a moment drop the assumption that $n\geq5$. For $A$ a
subring of~$\CC$ and for $n\geq 1$, we let $M_k(\Gamma_1(n),A)$ be the
sub-$A$-module of $M_k(\Gamma_1(n))$ consisting of elements $g$ such
that $a_m(g)\in A$ for all $m\geq 0$. In particular,
$M_k(\Gamma_1(n),\ZZ)$ is the submodule of forms whose $q$-expansion
has all its coefficients in~$\ZZ$. Similarly, for $A$ a subring
of~$\CC$ and for $n\geq 1$, we let $S_k(\Gamma_1(n),A)$ be the
sub-$A$-module of $S_k(\Gamma_1(n))$ consisting of elements $g$ such
that $a_m(g)\in A$ for all $m\geq 1$. For example, $\Delta$ belongs to
$S_{12}(\SL_2(\ZZ),\ZZ)$. The $S_k(\Gamma_1(n),A)$ are
$\TT(n,k)$-submodules of $S_k(\Gamma_1(n))$; see Propositions~12.3.11
and~12.4.1 of~\cite{Diamond-Im}.

We have the following pairing between $\TT(n,k)$ and
$S_k(\Gamma_1(n),\ZZ)$:
\begin{eqn}\label{eqn_TS_pairing}
S_k(\Gamma_1(n),\ZZ) \times \TT(n,k) \lto \ZZ, \quad (g,t)\mapsto a_1(tg).
\end{eqn}
This pairing is perfect, in the sense that it identifies each side
with the $\ZZ$-linear dual of the other; this follows easily from the
identity~(\ref{eqn_a1Tnf}).  It follows that the $\ZZ$-dual
$S_k(\Gamma_1(n),\ZZ)^\vee$ of $S_k(\Gamma_1(n),\ZZ)$ is free of rank
one as $\TT(n,k)$-module.  See \cite[12.4.13]{Diamond-Im}. 

For any $\ZZ$-algebra $A$ we let $\TT_A=\TT(n,k)_A$ denote $A\otimes
\TT(n,k)$, and $\TT_A^\vee$ will denote the $A$-linear dual of~$\TT_A$. It
can be proved that $\TT_\QQ^\vee$ is free of rank one as
$\TT_\QQ$-module, i.e., that $\TT_\QQ$ is \emph{Gorenstein}. One proof
is by explicit computation, see Theorem~3.5 and Corollary~3.6
of~\cite{Parent1}.  Another, more conceptual proof, uses the Petersson
inner product, and a so-called Atkin-Lehner pseudo-involution
$w_{\zeta_n}$, to show that $S_k(\Gamma_1(n))^\vee$ is isomorphic as
$\TT_\CC$-module to $S_k(\Gamma_1(n))$ itself;
see~\cite[12.4.14]{Diamond-Im}. It follows that $S_k(\Gamma_1(n))$ is
free of rank one as $\TT_\CC$-module, and, if $n\geq5$, that
$\QQ\otimes\rH^1(X_1(n),\calF_k)$ and its dual
$\QQ\otimes\rH^1(X_1(n),\calF_k)^\vee$ are free of rank two as
$\TT_\QQ$-module. It is this freeness result that will lead to the
fact that the Galois representations we get are two-dimensional.

We assume again that $n\geq5$. The step from the cohomological
interpretation of modular forms, given, over the complex numbers, by
the Eichler-Shimura isomorphism~(\ref{eqn_ES_iso}), to two-dimensional
$l$-adic Galois representations is made by comparing the cohomology
groups above to their $l$-adic counterparts for the étale topology,
and noting that $p\colon \EE\to X_1(n)$ is naturally defined
over~$\ZZ[1/n]$ as we have seen at the end of
Section~\ref{sec_modcurves}. From now on we will denote by $X_1(n)$
this model over $\ZZ[1/n]$, and by $X_1(n)(\CC)$ the Riemann surface
given by~$X_1(n)$. For any $\ZZ[1/n]$-algebra $A$, $X_1(n)_A$ will
denote the $A$-scheme obtained from $X_1(n)$ by extending scalars via
$\ZZ[1/n]\to A$.

We let $\calF_{k,l}$ denote the sheaf of $\QQ_l$-vector spaces
$\QQ_l\otimes \calF_k$ on $X_1(n)$. Then we have a canonical
isomorphism:
\begin{eqn}
\rH^1(X_1(n)(\CC),\calF_{k,l}) = \QQ_l\otimes \rH^1(X_1(n)(\CC),\calF_k).
\end{eqn}
The sheaves $\calF_{k,l}$ can also be constructed on the étale site
$X_1(n)_\et$, by taking the first derived direct image of
the constant sheaf $\QQ_l$ on $\EE_\et$ under $p\colon
\EE\to Y_1(n)$, then the $(k-2)$th symmetric power of that and finally
the pushforward from $Y_1(n)$ to~$X_1(n)$. 

The usual comparison theorems (comparing cohomology for étale and
Archi\-me\-dean topology, and étale cohomology over various algebraically
closed fields) give:
\begin{eqn}
\begin{aligned}
\rH^1(X_1(n)(\CC),\calF_{k,l}) & = \rH^1(X_1(n)_{\CC,\et},\calF_{k,l})\\
& = \rH^1(X_1(n)_{\Qbar,\et},\calF_{k,l}).
\end{aligned}
\end{eqn}
We put:
\begin{eqn}
W_l := \rH^1(X_1(n)_{\Qbar,\et},\calF_{k,l})^\vee.
\end{eqn}
By the results and the comparisons above, $W_l$ is, as
$\TT_{\QQ_l}$-module, free of rank~2, and $\Gal(\Qbar/\QQ)$ acts
continuously on it. To be precise: an element $\sigma$ of
$\Gal(\Qbar/\QQ)$ acts as $((\id\times\Spec(\sigma^{-1}))^*)^\vee$,
which is indeed covariant in~$\sigma$. The fact that the Hecke
correspondences exist over $\QQ$ makes that the
$\Gal(\Qbar/\QQ)$-action on $W_l$ commutes with the Hecke
operators. The choice of a $\TT_{\QQ_l}$-basis of $W_l$ gives us a
representation:
\begin{eqn}
\rho_l\colon \Gal(\Qbar/\QQ) \lto \GL_2(\TT_{\QQ_l}).
\end{eqn}
Recall that we have fixed an eigenform $f$ in
$S_k(\Gamma_1(n),\CC)$. Sending a Hecke operator to its eigenvalue for
$f$ then gives us a morphism of rings:
\begin{eqn}
\phi_f\colon \TT \lto \CC.
\end{eqn}
We let $K(f)$ be the image of $\TT_\QQ$ under $\phi_f$; it is the
finite extension of $\QQ$ obtained by adjoining all coefficients
$a_m(f)$ of the $q$-expansion of~$f$. We now view $\phi_f$ as a
morphism from $\TT$ to~$K(f)$. The tensor product $\QQ_l\otimes K(f)$
is the product of the completions $K(f)_\lambda$, with $\lambda$
ranging through the finite places of $K(f)$ that divide~$l$. For each
such $\lambda$ we then get a morphism
$\phi_{f,\lambda}\colon\TT_{\QQ_l}\to K(f)_\lambda$, and a
representation:
\begin{eqn}
\rho_{f,\lambda}\colon \Gal(\Qbar/\QQ) \lto \GL_2(K(f)_\lambda).
\end{eqn}
These are the representations mentioned in
Theorem~\ref{thm_exist_gal_rep}. It may be useful to note that the
space on which the representation is realised is:
\begin{eqn}
V_{f,\lambda} := K(f)_\lambda\otimes_{\TT_{\QQ_l}} W_l.
\end{eqn}
The representations $\rho_{f,\lambda}$ are continuous by
construction. The sheaves $\calF_{k,l}$ on $X_1(n)_{\ZZ[1/nl]}$ are
``lisse'' away from the cusps, and tamely ramified at the cusps,
hence, by Proposition~2.1.9 of \cite[XIII, \S2]{SGA7.II},
$\rho_{f,\lambda}$ is unramified at all $p$ not dividing~$nl$.

In the case where $k=2$ the construction of the $\rho_{f,\lambda}$ is
much simpler, because then the sheaf $\calF_k$ is the constant
sheaf~$\ZZ$ on $X_1(n)(\CC)$. The use of étale cohomology can then be
replaced by Tate modules of the Jacobian variety of~$X_1(n)$. We let
$J:=J_1(n)$ be this Jacobian variety, actually an Abelian scheme
over~$\ZZ[1/n]$. Then we have: 
\begin{eqn}
W_l = \QQ\otimes\varprojlim_m J(\Qbar)[l^m].
\end{eqn}
The fact that for $p$ a prime not dividing $nl$ the characteristic
polynomial of $\rho_{f,\lambda}(\Frob_p)$ is as stated in
Theorem~\ref{thm_exist_gal_rep} is obtained by studying the reduction
modulo~$p$ of the Hecke correspondence~$T_p$, i.e., as a
correspondence on~$X_1(nl)_{\FF_p}$, compatibly with the
sheaf~$\calF_{k,l}$. For details we refer to Conrad's
book~\cite{Conrad1} and Deligne's article~\cite{Deligne1}. In the case
$k=2$ this result is known as the \emph{Eichler-Shimura congruence
relation}, expressing the endomorphism $T_p$ of $J_{\FF_p}$ as $F+\ld
p\rd V$, where $F$ denotes the Frobenius endomorphism, and $V$ its
dual, i.e, the endomorphism satisfying $FV=p=VF$ in
$\End(J_{\FF_p})$. For details in the case $k=2$ we refer to
Section~12.5 of~\cite{Diamond-Im}.

Now that we have sketched the construction of the $l$-adic Galois
representations attached to modular forms, we mention some more of
their properties, that are not mentioned in
Theorem~\ref{thm_exist_gal_rep} and in the remarks directly following
that theorem.

The fact that Deligne proved the Riemann hypothesis part of Weil's
conjectures in \cite{Deligne2} implies very precise bounds on the
coefficients of modular forms. The reason for that is that the roots
of the equation $x^2-a_p(f)x+\eps_f(p)p^{k-1}$ are eigenvalues of the
Frobenius at $p$ on the space
$\rH^{k-1}(\ol{\EE^{k-2}}_{\Fbar_p,\et},\QQ_l)$. We state these
bounds, called \emph{Ramanujan bounds}, in a theorem, due to Deligne
in the case $k\geq 2$, and to Deligne-Serre (\cite{Deligne-Serre}) in
the case $k=1$.
\begin{thm}
Let $f$ be a normalised newform, let $n$ be its level and $k$ its
weight. Then for $p$ not dividing $n$, we have:
\begin{eqn}
|a_p(f)| \leq 2{\cdot}p^{(k-1)/2}.
\end{eqn}
\end{thm}
A slightly weaker result than in the theorem above, stating that, for a
given $f$ as above, $|a_m(f)|=O(m^{k/2})$, can be obtained in a very
elementary way; see~\cite[Cor.~2.1.6]{Miyake1} (the idea is to use that
the function $z\mapsto |f(z)|(\Im(z))^k$ is bounded on~$\HH$ and to
view $a_m(f)$ as a residue).

Theorem~\ref{thm_exist_gal_rep} gives us information on the restriction
$\rho_{f,\lambda,p}$ of $\rho_{f,\lambda}$ to decomposition groups
$\Gal(\Qbar_p/\QQ_p)$ for $p$ not dividing~$nl$. Namely, the theorem
says that such restrictions $\rho_{f,\lambda,p}$ are unramified, and
it gives the eigenvalues of
$\rho_{f,\lambda,p}(\Frob_p)$. Unfortunately, it is not known if
$\rho_{f,\lambda,p}(\Frob_p)$ is semi-simple; see
\cite{Coleman-Edixhoven} for information on this.

We should note that also in the case that $p$ divides $nl$ almost
everything is known about~$\rho_{f,\lambda,p}$. For $p$ not dividing
$l$, this is the very general statement that the ``Frobenius
semi-simplification'' of $\rho_{f,\lambda,p}$ corresponds, via a
suitably normalised local Langlands correspondence, to a certain
representation $\pi_{f,p}$ of $\GL_2(\QQ_p)$ attached to~$f$. This
result is due, in increasing order of generality, to Langlands,
Deligne, and Carayol. For details on this the reader is referred
to~\cite{Carayol1}, which gives this result in the more general
context of Hilbert modular forms (i.e., $\QQ$ is replaced by a totally
real number field). The result for $p=l$ uses Fontaine's $p$-adic
Hodge theory, and is due to Saito (\cite{Saito1} for the case of
modular forms, and \cite{Saito2} for the case of Hilbert modular
forms).

% add a reference that the \rho_{\Delta,l} (l-adic!) cannot be obtained by
% tensor constructions from abelian varieties (Motivic question by
% Schappacher in a book of his (Lenny showed it to me)). Zie mijn
% dagboek. Don Blasius had already done this, reviewed by Ben Moonen.

\section{Galois representations over finite fields, and reduction to
  torsion in Jacobians}\label{sec_red_to_jac} 
We start this section by explaining how to pass from $l$-adic Galois
representations to Galois representations over finite fields.

Let $f=\sum a_mq^m$ be a (complex) normalised cuspidal eigenform for
all Hecke operators $T_m$, $m\geq1$, of some level $n\geq1$ and of
some weight $k\geq2$. As in Theorem~\ref{thm_exist_gal_rep} we have
the Galois representations $\rho_{f,\lambda}$, from $\Gal(\Qbar/\QQ)$
to~$\GL_2(\Qbar_l)$. It follows from the construction of those
representations that there is a finite subextension $\QQ_l\to E$ of
$\QQ_l\to\Qbar_l$ such that $\rho_{f,\lambda}$ takes its values
in~$\GL_2(E)$. (Actually, this can also be deduced from the continuity
alone; see the proof of Corollary~5 in~\cite{Dickinson1} for an
argument.) The question as to what the smallest possible $E$ is can be
easily answered. Such an $E$ must contain the traces $a_p(f)$ of the
$\rho_{f,\lambda}(\Frob_p)$ for all $p$ not dividing~$nl$. So let $K$
be the extension of $\QQ$ generated by the $a_p(f)$ with $p$ not
dividing~$n$, i.e., $K$ is the field of definition of the newform
corresponding to~$f$. Then $E$ can be taken to be $K_\lambda$, the
$l$-adic completion of $K$ specified by the embedding $\lambda$ of $K$
into~$\Qbar_l$ (see Section~12.5 in~\cite{Diamond-Im}).

Let now $\rho_{f,\lambda}\colon \Gal(\Qbar/\QQ)\to\GL_2(E)$ be a
realisation of $\rho_{f,\lambda}$ over~$E$ as above. As
$\rho_{f,\lambda}$ is semisimple (it is even irreducible), such a
realisation is unique up to isomorphism (because it is determined by
the traces). Let $O_E$ be the ring of integers in $E$, i.e., the
integral closure of $\ZZ_l$ in~$E$. As $\Gal(\Qbar/\QQ)$ is compact,
it stabilises some $O_E$ lattice in $E^2$ (in the set of lattices, the
orbits under $\Gal(\Qbar/\QQ)$ are finite, take the intersection, or
the sum, of the lattices in one orbit). This means that, after
suitable conjugation (choose such a lattice, and an $O_E$-basis of
it), $\rho_{f,\lambda}$ takes values in~$\GL_2(O_E)$. We let
$O_E\to\Fbar_l$ denote the morphism induced by the given embedding of
$E$ into~$\Qbar_l$ (we view $\Fbar_l$ as the residue field of the
subring of integers $\Zbar_l$ of~$\Qbar_l$). We can then define the
\emph{residual} Galois representation $\ol{\rho}_{f,\lambda}$ to be
the \emph{semi-simplification} of the composed representation
$\Gal(\Qbar/\QQ)\to\GL_2(O_E)\to\GL_2(\Fbar_l)$. Another choice of $E$
or of lattice or basis leads to an isomorphic~$\ol{\rho}_{f,\lambda}$,
but we note that without the operation of semi-simplification this
would not be true (see Chapter~III of~\cite{Serre3}).

Given $f$, all but finitely many of the $\ol{\rho}_{f,\lambda}$ are
irreducible. This was proved for $f$ of level one and with
coefficients in $\ZZ$ in Theorem~4 of~\cite{Swinnerton-Dyer1}. The
general case follows easily from Theorem~2.3 of~\cite{Faltings-Jordan},
which says that if $\ol{\rho}_{f,\lambda}$ is reducible with $l>k$ not
dividing~$n$, then
$\ol{\rho}_{f,\lambda}=\alpha\oplus\beta\ol{\chi}_l^{k-1}$ with
$\alpha$ and $\beta$ unramified outside~$n$, and
$\ol{\chi}_l\colon\Gal(\Qbar/\QQ)\to\FF_l^\times$ the mod~$l$
cyclotomic character. Moreover, the proof shows that the set of $l$
such that some $\ol{\rho}_{f,\lambda}$ is reducible can be bounded
explicitly. 

The next question that we want to answer is the following: what is the
smallest subfield of $\Fbar_l$ over which $\ol{\rho}_{f,\lambda}$ can
be realised? Just as for $\rho_{f,\lambda}$ itself, that subfield must
contain the traces of the $\ol{\rho}_{f,\lambda}(\Frob_p)$ for all $p$
not dividing~$nl$. That condition turns out to be sufficient, as we
will now show. So we let, in this paragraph, $\FF$ be the subfield of
$\Fbar_l$ that is generated by the images $\ol{a_p(f)}$ in $\Fbar_l$
of the $a_p(f)$ in~$\Zbar_l$. Then for any $\sigma$ in
$\Gal(\Fbar_l/\FF)$ the conjugate $\ol{\rho}_{f,\lambda}^\sigma$ of
$\ol{\rho}_{f,\lambda}$ and $\ol{\rho}_{f,\lambda}$ itself are both
semisimple and give the same characteristic polynomials as functions
on~$\Gal(\Qbar/\QQ)$. Therefore, by a theorem of Brauer-Nesbitt (see
Theorem~30.16 of~\cite{Curtis-Reiner1}), $\ol{\rho}_{f,\lambda}$ is
isomorphic to all its conjugates over~$\FF$. (A more general statement
of this kind is given in Exercise~1 of Section~18.2 of~\cite{Serre3}.)
The fact that $\Gal(\Fbar_l/\FF)$ is equal to $\Zhat$ then implies
that $\ol{\rho}_{f,\lambda}$ can be realised over~$\FF$. Let us give
an argument for that in terms of matrices, although a much more
conceptual argument would be to say that a ``gerbe over a finite field
is trivial''. Let $\sigma$ be the Frobenius element of
$\Gal(\Fbar_l/\FF)$, and let $s$ be an element of $\GL_2(\Fbar_l)$
such that for all $g$ in the image of $\ol{\rho}_{f,\lambda}$ we have
$\sigma(g) = sgs^{-1}$. Then take a $t$ in $\GL_2(\Fbar_l)$ such that
$s=\sigma(t)^{-1}t$. Then all $tgt^{-1}$ are in~$\GL_2(\FF)$. By
Brauer-Nesbitt, the realisation over~$\FF$ is unique.

For a discussion on possible images of $\ol{\rho}_{f,\lambda}$ we
refer the reader to the introduction of~\cite{Kiming-Verrill} (we note
however that for $f$ a ``CM-form'', i.e., a form for which all
$l$-adic Galois representations have dihedral image) infinitely many
of the $\ol{\rho}_{f,\lambda}$ can have dihedral image in
$\PGL_2(\Fbar_l)$). In particular, Theorem~2.1 of~\cite{Ribet1} states
that for $f$ not a CM-form only finitely many of the images of the
$\ol{\rho}_{f,\lambda}$ are exceptional
%%%%%exceptional should go in the index of terminology%%%%%%%%%%%%%%%%
in the sense that they are of order prime to~$l$.  See also
Theorem~\ref{thm_large_image} for the case where $n=1$. For $l>3$ such
that $\ol{\rho}_{f,\lambda}$ is irreducible and not exceptional, a
result of Dickson, see Chapter~XII of~\cite{Dickson}, or the proof of
Theorem~2.5 of~\cite{Ribet2}, says that the image
of~$\ol{\rho}_{f,\lambda}$ in~$\PGL_2(\Fbar_l)$ is, after suitable
conjugation, equal to $\PGL_2(\FF)$ or $\SL_2(\FF)/\{1,-1\}$ for some
finite extension~$\FF$ of~$\FF_l$. We note that this field~$\FF$ can
be smaller than the extension $\FF_l$ generated by the traces of
$\ol{\rho}_{f,\lambda}$ (indeed, twisting does not change the
projective image, but it can make the field generated by the traces
bigger).

%% Will we??? Still??? Look it up: search for {lem_transitive}
We will use later on the following lemma.
\begin{lem}\label{lem_transitive}
Let $l$ be a prime number, let $V$ be a two-dimensional $\FF_l$-vector
space, and let $G$ be a subgroup of~$\Aut(V)$ of order a multiple
of~$l$, and such that $V$ is irreducible as a representation
of~$G$. Then $G$ contains $\SL(V)$, and acts transitively
on~$V-\{0\}$. 
\end{lem}
\begin{proof}
Let $g_1$ be an element of $G$ of order~$l$. Then the kernel $L_1$ of
$g_1-\id_V$ is a line. As $V$ is irreducible, $L_1$ is not
$G$-invariant, hence we can take an element $h$ in~$G$ such that
$L_2:=hL_1$ is not~$L_1$. Then $g_2:=hg_1h^{-1}$ is of order~$l$ and
fixes~$L_2$. Let $e_1$ and $e_2$ be non-zero elements of $L_1$
and~$L_2$, respectively. Then with respect to the basis $e=(e_1,e_2)$
of~$V$, $g_1$ and $g_2$ are given by elementary matrices of the form
$(\begin{smallmatrix}1 & a\\ 0 & 1\end{smallmatrix})$ and
$(\begin{smallmatrix}1 & 0\\ b & 1\end{smallmatrix})$, respectively,
with $a$ and $b$ non-zero, and hence generate~$\SL(V)$.
\end{proof}
It follows that if $\ol{\rho}_{f,\lambda}$ takes values in
$\GL_2(\FF_l)$, and is irreducible and not exceptional, then $\im
\ol{\rho}_{f,\lambda}$ contains~$\SL_2(\FF_l)$, and therefore is the
subgroup of elements of $\GL_2(\FF_l)$ whose determinant is in the
image of the character~$\ol{\eps_f}{\cdot}\ol{\chi}_l^{k-1}$. In that
case, $\im\ol{\rho}_\lambda$ acts transitively on~$\FF_l^2-\{0\}$.

The properties of residual Galois representations that we have seen
above show that we do not need to define them via $l$-adic Galois
representations, but that we can start from maximal ideals in Hecke
algebras.
\begin{thm}\label{thm_repr_ass_to_mod_l_mod_form}
Let $n$ and $k$ be positive integers. Let $\FF$ be a finite field, and
$f\colon \TT(n,k)\to\FF$ a surjective morphism of rings. Then there is a
continuous semi-simple representation:
\[
\rho_f\colon\Gal(\Qbar/\QQ)\lto \GL_2(\FF)
\]
that is unramified outside $nl$, where $l$ is the characteristic
of~$\FF$, such that for all $p$ not dividing $nl$ we have, in~$\FF$:
\[
\trace(\rho(\Frob_p)) = f(T_p)\quad\text{and}\quad 
\det(\rho(\Frob_p))=f(\ld p\rd)p^{k-1}.
\]
Such a $\rho_f$ is unique up to isomorphism (i.e., up to conjugation).
\end{thm}
\begin{proof}
Let $n$, $k$, $\FF$ and $f$ be given. As $\TT:=\TT(n,k)$ is free of
finite rank as $\ZZ$-module, $\Spec(\TT)$ has only finitely many
irreducible components, each of which is one-dimensional and finite
over $\Spec(\ZZ)$. Therefore, the maximal ideal $\ker(f)$ of
$\TT_{\FF_l}$ is the specialisation of a maximal ideal $m$ of
$\TT_\QQ$. Let $K$ be the quotient $\TT_\QQ/m$. Then the quotient
morphism $\TT_\QQ\to K$ \emph{is} a normalised eigenform $\tilde{f}$ in
$S_k(\Gamma_1(n))_K$, and $\rho_f$ is the realisation over~$\FF$ of
the reduction of some~$\rho_{\tilde{f},\lambda}$.
\end{proof}
Let now $f$ be as in Theorem~\ref{thm_repr_ass_to_mod_l_mod_form}, and
let us suppose now that $\rho_f$ is irreducible. The construction of
$l$-adic Galois representations that we recalled in
Section~\ref{sec_galreps} implies that the dual of $\rho_f$ occurs in
$\rH^{k-1}(\ol{\EE^{k-2}}_{\Qbar,\et},\FF_l)$, as well as in
$\rH^1(X_1(n)_{\Qbar,\et},\ol{\calF}_{k,l})$, where $\ol{\calF}_{k,l}$
is defined as $\calF_{k,l}$ but with $\QQ_l$ replaced by~$\FF_l$. Let
us now assume that $k>2$. Then both these realisations are difficult
to deal with computationally. In the first representation the
difficulty arises from the degree $k-1$ étale cohomology; it seems to
be unknown how to deal explicitly with elements of such cohomology
groups. In the second representation, the elements of the cohomology
group are isomorphism classes of $\ol{\calF}_{k,l}$-torsors,
on~$X_1(n)_{\Qbar,\et}$. Such torsors can be described explicitly, as
certain covers of $X_1(n)_\Qbar$ with certain extra data. The set of
such torsors can probably be described by a system of polynomial
equations that can be written down in time polynomial in~$nl$ (think
of the variables as coefficients of certain equations for the
torsors). But the problem is that, apparently, there are no good
methods known to solve these systems of equations (the number of
variables grows too fast with~$l$ and the equations are not
linear). In fact, the \emph{satisfiability problem} SAT, which is
known to be NP-complete (Cook's theorem, see for
example~\cite{Moret1}, or wikipedia), is a special case of the problem
of deciding whether or not a polynomial system of equations over
$\FF_2$ has a solution over~$\FF_2$. We note that the description of
the set of torsors by a system of polynomial equations should also
work over suitable finite extensions of finite fields $\FF_p$, in time
polynomial in~$l\log p$.
%% See if this part does not overlap unnecessarily with the
%% introduction!!!!!!!!!!!!!!!!!!!!!!!!!!!!!!!!!

Another place where $\rho_f$ occurs is in $J_1(nl)(\Qbar)[l]$, i.e.,
in the $l$-torsion of the Jacobian of the modular curve with
level~$nl$, if $l+1\geq k$ and $l\not|n$. This means that at the cost
of increasing the level by a factor~$l$, we are reduced to dealing
with torsion points on Abelian varieties. Of course, the $l$-adic
representations $\rho_{\tilde{f},\lambda}$ attached to lifts of $f$
do not occur in the Jacobian of any curve, simply because the
Frobenius eigenvalues are Weil numbers of the wrong weight. What
happens here for $\rho_f$ is a ``mod~$l$ phenomenon'' having to do
with ``congruences'' between modular forms. Before we give a detailed
statement, let us explain why this happens (such explanations date
back at least to the 1960's; Shimura, Igusa, Serre,\ldots).

For simplicity, and only during this explanation, we assume that
$n\geq 5$. Then we have a universal elliptic curve with a given point
of order~$n$ over $\ZZ[1/nl]$-schemes: $(\EE/Y_1(n),\PP)$. We let
$p\colon\EE\to Y_1(n)$ denote the structure morphism. By definition,
we have:
\begin{eqn}
\ol{\calF}_{k,l} = \Sym^{k-2}\rR^1p_*\FF_l .
\end{eqn}
As explained at the end of Section~\ref{sec_schoof_etale}, we have a
natural isomorphism:
\begin{eqn}
\rR^1p_*\FF_l = \EE[l]^\vee .
\end{eqn}
And by the definition of $Y_1(nl)$, and the Weil pairing, we have an
exact sequence on $Y_1(nl)_\et$:
\begin{eqn}
0\lto \FF_l \lto \EE[l] \lto \mu_l \lto 0,
\end{eqn}
where $\FF_l$ and $\mu_l$ denote the corresponding constant
sheaves. It follows that the pullback of $\rR^1p_*\FF_l$ to
$Y_1(nl)_\et$ has a 2-step filtration with successive quotients
$\FF_l$ and $\mu_l^\vee$. Therefore, $\ol{\calF}_{k,l}$ has a
filtration in $k-1$ steps, with successive quotients $\FF_l^{\otimes
i}\otimes(\mu_l^\vee)^{\otimes j}=\mu_l^{\otimes -j}$, with $i+j=k-2$,
$i\geq0$, $j\geq0$. In particular, we get a map:
\begin{eqn}
\begin{aligned}
\rH^1(X_1(n)_{\Qbar,\et},\ol{\calF}_{k,l}) & \lto 
\rH^1(X_1(nl)_{\Qbar,\et},\ol{\calF}_{k,l}) \lto \\
 & \lto \rH^1(X_1(nl)_{\Qbar,\et},\FF_l) = J_1(nl)(\Qbar)[l]^\vee .    
\end{aligned}
\end{eqn}

This map explains that $\rho_f$ is likely to occur
in~$J_1(nl)(\Qbar)[l]$.  A better way to analyse this map is in fact
by studying the direct image of the constant sheaf $\FF_l$ via the map
$X_1(nl)\to X_1(n)$.  A recent detailed treatment of this method, and
precise results can be found in~\cite{Wiese1}.
% make the reference more precise....

Another way to show that $\rho_f$ occurs in~$J_1(nl)(\Qbar)[l]$ is to
study modular forms mod~$l$ of level $nl$ and of weight~$2$. This is
more complicated than the modular forms that we have seen before, as
it uses the study of the reduction mod $l$ of the modular
curve~$X_1(nl)$, which is not smooth. The study of these reductions
has its roots in Kronecker's congruence relation. The most complete
modern accounts of such material are given in the
article~\cite{Deligne-Rapoport} by Deligne and Rapoport and in the
book~\cite{Katz-Mazur} by Katz and Mazur. A construction of $\rho_f$
in $J_1(nl)(\Qbar)[l]$, following suggestions from Serre, was given by
Gross in~\cite{Gross1}.

We are now in a position to state the following theorem, that,
combining Gross's result with a so-called multiplicity one theorem,
gives us a useful realisation of~$\rho_f$. As it is nowadays customary
to say, it is a result due to ``many people'' (mainly Mazur, Ribet,
Gross (and Edixhoven for the multiplicity one part)).
\begin{thm}\label{thm_red_to_wt_2}
Let $n$ and $k$ be positive integers, $\FF$ a finite field and $l$ its
characteristic, and $f\colon\TT(n,k)\to\FF$ a surjective ring
morphism. Assume that $2<k\leq l{+}1$ and that the associated Galois
representation $\rho_f$ from $\Gal(\Qbar/\QQ)$ to $\GL_2(\FF)$ is
absolutely irreducible. Then there is a unique ring morphism
$f_2\colon\TT(nl,2)\to\FF$ such that for all $i\geq 1$ one has
$f_2(T_i)=f(T_i)$. The morphism $f_2$ is surjective. Let
$m_f=\ker(f_2)$, and let $V_f\subset J_1(nl)(\Qbar)$ denote the kernel
of~$m_f$, i.e., the $\FF$-vector space of elements $x$ in
$J_1(nl)(\Qbar)$ such that $tx=0$ for all $t$ in~$m_f$. Then $V_f$ is
a finite, non-zero, direct sum of copies of~$\rho_f$. If $k<l$ then
the multiplicity of $\rho_f$ in $V_f$ is one, i.e., $V_f$
realises~$\rho_f$. For all $a\in(\ZZ/nl\ZZ)^\times$, one has $f_2(\ld
a\rd)=f(\ld a\rd)a^{k-2}$, where we still denote by $a$ its images in
$\ZZ/n\ZZ$ and in~$\FF_l$.
\end{thm}
\begin{proof}
The existence of $f_2$ and the statement that $V_f$ is a successive
extenstion of copies of~$\rho_f$ are given in~\cite{Gross1} (see his
Proposition~11.8). In Section~6 of~\cite{Edixhoven3} it is proved,
applying results from~\cite{Boston-Lenstra-Ribet}, that $V_f$ is a
direct sum of copies of~$\rho_f$.  Case~1 of Theorem~9.2
of~\cite{Edixhoven3} gives the multiplicity one result.
\end{proof}
\begin{rem}
See~\cite{Wiese2}, Corollary~4.5, for a complete result on the
multiplicity one question for weights $k$ with $2\leq k\leq l{+}1$.
In particular, if $k=l$ and $\rho_f$ is unramified at~$l$ and
$\rho_f(\Frob_l)$ is scalar, then this multiplicity is not one.
\end{rem}
As we want to describe $V_f$ explicitly, we will need a bound on the
amount of Hecke operators needed to describe $\TT(nl,2)$ and its
ideal~$m_f$. We start by quoting a result of Jacob Sturm
(see~\cite{Sturm}).
\begin{thm}[Sturm]\label{thm_sturm}
Let $N\geq1$ be an integer, $\Gamma$ a subgroup of $\SL_2(\ZZ)$
containing $\Gamma(N)$. Let $N'$ be the ``width'' of the cusp~$\infty$
for~$\Gamma$, i.e., the positive integer defined by
$\Gamma\cap(\begin{smallmatrix}1 & \ZZ\\0 & 1\end{smallmatrix}) =
(\begin{smallmatrix}1 & N'\ZZ\\0 & 1\end{smallmatrix})$. Let $f$ be a
modular form on $\Gamma$ of weight~$k$, with coefficients in a
discrete valuation ring $R$ contained in~$\CC$. Let $F$ be the residue
field of~$R$, and suppose that the image $\sum a_nq^{n/N'}$ in
$F[[q^{1/N'}]]$ of the $q$-expansion of $f$ has $a_n=0$ for all $n\leq
k[\SL_2(\ZZ):\Gamma]/12$. Then $a_n=0$ for all~$n$, i.e., $f$ is
congruent to~$0$ modulo the maximal ideal of~$R$.
\end{thm}
This result of Sturm gives as a direct consequence a bound for up to
where one has to take $T_i$ so that one gets a system of generators of
the Hecke algebra as $\ZZ$-module, for a given level and weight. See
Section~9.4 of~\cite{Stein} for a detailed proof of Sturm's result,
and of this consequence. For convenience we also state and prove this
result in the precise context where we use it.
\begin{thm}\label{thm_gen_Hecke}
Let $N\geq 1$ and $k\geq 1$ be integers, and let $\TT(N,k)$ be the
Hecke algebra attached to $S_k(\Gamma_1(N))$, i.e., $\TT(N,k)$ is the
$\ZZ$-submodule of $\End_\CC(S_k(\Gamma_1(N)))$ generated by the
$T_n$, for $n\geq 1$, and the $\ld a\rd$, for $a$
in~$(\ZZ/N\ZZ)^\times$. Then $\TT(N,k)$ is generated, as $\ZZ$-module,
by the $T_i$ with $1\leq i\leq k[\SL_2(\ZZ):\Gamma_1(N)]/12$.
\end{thm}
\begin{proof}
Let $S$ be the $\ZZ$-module $S_k(\Gamma_1(N),\ZZ)$. Then
by~(\ref{eqn_TS_pairing}) we have isomorphisms of $\TT(N,k)$-modules:
$S=\TT(N,k)^\vee$, and $\TT(N,k)=S^\vee$. Now the result of Sturm
above says that for each prime number~$p$, the elements $T_i$, $1\leq
i\leq k[\SL_2(\ZZ):\Gamma_1(N)]/12$, generate the $\FF_p$-vector space
$\FF_p\otimes S^\vee$, and hence they generate
$\FF_p\otimes\TT(N,k)$. So, indeed, these $T_i$ generate $\TT(N,k)$ as
a $\ZZ$-module.
\end{proof}
We can now state a complement to Theorem~\ref{thm_red_to_wt_2}.
\begin{prop}\label{prop_compl_red_to_wt_2}
In the situation of Theorem~\ref{thm_red_to_wt_2}, the Hecke algebra
$\TT(nl,2)$ is generated, as $\ZZ$-module, by the Hecke operators
$T_i$ with $1\leq i\leq [\SL_2(\ZZ):\Gamma_1(nl)]/6$.
\end{prop}
We remark that, still in the same situation, giving generators of
$m_f$ is then a matter of simple linear algebra over $\FF_l$ in a
vector space of suitably bounded dimension.

We consider the particular case of a mod~$l$ eigenform of level one
and of weight~$k$, viewed as a ring morphism $f$ from $\TT(1,k)$ to a
finite extension of~$\FF_l$. Then we have the following result, that
states explicitly how the Galois representation attached to~$f$ is
realised in the Jacobian~$J_1(l)(\Qbar)$. Recall that $\TT(l,2)$, the
Hecke algebra acting on weight two cusp forms on $\Gamma_1(l)$, is
generated as $\ZZ$-module by the $T_j$ with $1\leq j\leq (l^2{-}1)/6$.

%% Generalised to more general forms.  
\begin{thm}\label{thm_red_to_wt_2_Delta}
Let $l$ be a prime number, let $k$ be an integer such that $2<k\leq
l{+}1$, and $f\colon \FF_l\otimes\TT(1,k)\to\FF$ a surjective ring
morphism with $\FF$ a finite field of characteristic~$l$, such that
the associated Galois representation
$\rho\colon\Gal(\Qbar/\QQ)\to\GL_2(\FF)$ is irreducible. Let
$f_2\colon \FF_l\otimes\TT(l,2)\to\FF$ be the morphism of rings such
that for all $m\in\ZZ_{\geq 1}$ we have $f_2(T_m)=f(T_m)$ (see
Theorem~\ref{thm_red_to_wt_2}). Let $(t_1,\ldots,t_r)$ be a system of
generators for $\ker(f_2)$. Let:
\[
V_f := 
\bigcap_{1\leq i\leq r} 
\ker(t_i,J_1(l)(\Qbar)[l]).
\]
Then $V_f$ is a $2$-dimensional $\FF$-vector space
realising~$\rho$. For $p\neq l$ prime, $\ld p\rd$ acts on $V_f$ as
multiplication by~$p^{k-2}$.

One can obtain a system of generators of $\ker(f_2)$ as follows. For
$i$ in $\{1,\ldots,(l^2-1)/6\}$, either $f(T_i)$ is an $\FF_l$-linear
combination of the $f(T_j)$ with $j{<}i$, or it is not. If it is not,
then let $t_i=0$. If it is, then pick one:
$f(T_i)=\sum_{j<i}a_{i,j}f(T_j)$, and let
$t_i=T_i-\sum_{j<i}a_{i,j}T_j$.
\end{thm}
\begin{proof}
Just as in the proof of Theorem~\ref{thm_red_to_wt_2}, we use
Theorem~9.2 of~\cite{Edixhoven3}, but this time in the case of level
one. Case~1 of that theorem deals with the $k$ that satisfy $2<
k<l$. Case~3 deals with the case $k=l$, because $\rho$, being
unramified outside~$\{l\}$, and being irreducible of dimension two,
\emph{is} ramified at~$l$. Case~4 deals with the case $k=l{+}1$,
because there are no nonzero cusp forms of weight two and level one.
\end{proof}
We also state the following definition and theorem here, because the
result, to be used later, is directly related to
Theorem~\ref{thm_red_to_wt_2}. The theorem is due, again, to ``many
people'', just as Theorem~\ref{thm_red_to_wt_2} itself.
\begin{defi}\label{defi_w_operator}
Let $N\geq 1$, and let $\ZZ[\zeta_N]$ be the subring of $\CC$
generated by a root of unity of order~$N$. To a pair
$(E/S/\ZZ[1/N,\zeta_N],P)$ consisting of an elliptic curve $E$ over a
$\ZZ[1/N,\zeta_N]$-scheme $S$, together with a point $P$ in $E(S)$
that is of order $N$ everywhere on~$S$, we associate another such pair
$(E'/S'/\ZZ[1/N,\zeta_N],P')$ as follows. Let $\beta\colon E\to E'$ be
the isogeny whose kernel is the subgroup of $E$ generated by~$P$. Let
$\beta^\vee\colon E'\to E$ be the dual of~$\beta$ (see Section~2.5
of~\cite{Katz-Mazur}).  Let $P'$ be the unique element of
$\ker(\beta^\vee)(S)$ such that $e_\beta(P,P')=\zeta_N$, where
$e_\beta$ is the perfect $\mu_N$-valued pairing between $\ker(\beta)$
and $\ker(\beta^\vee)$ as described in Section~2.8
of~\cite{Katz-Mazur}.  This construction induces an automorphism
$w_{\zeta_N}$ of the modular curve $X_1(N)_{\ZZ[1/N,\zeta_N]}$, called
an ``Atkin-Lehner pseudo-involution''.
\end{defi}

\begin{thm}\label{thm_Hecke_Gorenstein_etc}
In the situation of Theorem~\ref{thm_red_to_wt_2} the completion
$\TT_{m_{f,\lambda}}$ of $\TT$ at $m_{f,\lambda}$ is
\emph{Gorenstein}, i.e., the $\ZZ_l$-linear dual of
$\TT_{m_{f,\lambda}}$ is free of rank one as
$\TT_{m_{f,\lambda}}$-module. For all $r\geq 1$, the
$(\ZZ/l^r\ZZ)\otimes\TT_{m_{f,\lambda}}$-module
$J_1(nl)(\Qbar)[l^r]_{m_{f,\lambda}}$ is free of rank~$2$. 

For any $t$ in $\TT$ we have $t^\vee=wtw^{-1}$, where $t^\vee$ is the
dual of $t$ as endomorphism of the self-dual Abelian variety
$J_1(nl)_{\QQ(\zeta_{nl})}$, and where $w$ is the endomorphism of
$J_1(nl)_{\QQ(\zeta_{nl})}$ induced via Picard functoriality by the
automorphism $w_{\zeta_{nl}}$ of $X_1(nl)_{\QQ(\zeta_{nl})}$.

For $r\geq 0$, let $({\cdot},{\cdot})_r$ denote the Weil pairing on
$J_1(nl)(\Qbar)[l^r]$, and let $\ld {\cdot},{\cdot}\rd_r$ denote the
pairing defined by: 
\[
\ld x,y\rd_r = (x,w(y))_r .
\]
Then $\ld {\cdot},{\cdot}\rd_r$ is a perfect pairing on
$J_1(nl)(\Qbar)[l^r]$ for which the action of $\TT$ is
self-adjoint. As a consequence, $\ld {\cdot},{\cdot}\rd_r$ induces a
perfect pairing on~$J_1(nl)(\Qbar)[l^r]_{m_{f,\lambda}}$.
\end{thm}
\begin{proof}
See Sections~6.4 and~6.8 of~\cite{Edixhoven3}.
\end{proof}
\begin{rem}
See Corollary~4.2 of~\cite{Wiese2} for a proof that
$\TT_{m_{f,\lambda}}$ is \emph{not} Gorenstein if the multiplicity of
$\rho_f$ in $V_f$ is not one.
\end{rem}
The next result gives an effective criterion whether two modular forms
give isomorphic residual Galois representations.
\begin{prop}\label{prop_test_equality_res_repr}
Let $l$ be a prime number, $\FF$ a finite extension of~$\FF_l$, $k_1$
and $k_2$ in $\ZZ_{\geq0}$, and $f_1\colon\TT(1,k_1)\to\FF$ and
$f_2\colon\TT(1,k_2)\to\FF$ two morphisms of rings, and $i$
in~$\{0,\ldots,l{-2}\}$. Then $\rho_{f_1}$ and $\rho_{f_2}\otimes\chi_l^i$
are isomorphic if and only if $k_1=k_2{+}2i$ in $\ZZ/(l{-}1)\ZZ$ and
for all primes $p\neq l$ with $p\leq (l^2{-}1)/12$ we have
$f_1(T_p)=p^if_2(T_p)$.
\end{prop}
\begin{proof}
Assume first that $\rho_{f_1}$ and $\rho_{f_2}\otimes\chi_l^i$ are
isomorphic. Then we have
$\det\rho_{f_1}=\det(\rho_{f_2}\otimes\chi_l^i)$, hence
$\chi_l^{k_1-1}=\chi_l^{k_2-1+2i}$, hence $k_1=k_2{+}2i$ in
$\ZZ/(l{-}1)\ZZ$. For all primes $p\neq l$, we have
$f_1(T_p)=p^if_2(T_p)$ because they are the traces of the images under
$\rho_{f_1}$ and $\rho_{f_2}\otimes\chi_l^i$ of the Frobenius at~$p$.

Assume now that $k_1=k_2{+}2i$ in $\ZZ/(l{-}1)\ZZ$ and that for all primes
$p\neq l$ with $p\leq (l^2{-}1)/12$ we have
$f_1(T_p)=p^if_2(T_p)$. Then $\det\rho_{f_1}$ and
$\det(\rho_{f_2}\otimes\chi_l^i)$ are equal, hence it suffices to
prove that for all primes $p\neq l$ we have: $f_1(T_p)=p^if_2(T_p)$.

We will use some theory on ``Katz modular forms''; see Sections~2
and~3 of~\cite{Edixhoven3} for a short account. For $a$ in~$\ZZ$, we
denote by $M_a(1,\FF_l)$ the space of Katz modular forms of level one
and weight~$a$ over~$\FF_l$, and by $S_a(1,\FF_l)$ its subspace of
cuspidal forms. Our reason to use Katz modular forms
over~$\FF_l$ is that this gives us the Hasse invariant $A$ in
$M_{l-1}(1,\FF_l)$ and the operators $\theta\colon
M_a(1,\FF_l)\to S_{a+l+1}(1,\FF_l)$, for all $a\in\ZZ_{\geq 0}$, that,
on $q$-expansions, act as the differential
operator~$q{\cdot}d\!/\!dq$. See \cite[\S3]{Edixhoven3} for the
properties of $\theta$ that we will use.

The idea in what follows is to use $\theta$ to pass to eigenforms that
are annihilated by~$T_l$, and to pass to eigenforms of weight at most
$l^2{-}1$ by dividing by $A$ as many times as possible. Recall that
the $q$-expansion of $A$ is the constant~$1$.

We write $\theta^{l-1}f_1=A^{n_1}f_1'$, with $n_1$ maximal, and we let
$k_1'$ be the weight of~$f_1'$. Then $k_1'\leq l^2{-}1$ by Theorem~3.4
of~\cite{Edixhoven3} and the definition of $\theta$-cycles; note that
$l{+}1+(l{-}2)(l{+}1)=l^2{-}1$. Similarly, we write
$\theta^{l-1-i}f_2=A^{n_2}f_2'$, with $n_2$ maximal, and we let $k_2'$
be the weight of $f_2'$. Then $k_1'=k_2'$ in $\ZZ/(l{-1})\ZZ$, $f_1'$
and $f_2'$ are eigenforms, annihilated by~$T_l$, and with the same
eigenvalues for all $T_p$ with $p\leq (l^2{-}1)/12$ prime. This
implies that for all $m\leq (l^2{-}1)/12$ we have
$f_1'(T_m)=f_2'(T_m)$. If $k_1'\geq k_2'$, then
$f_1'=A^{(k_1'-k_2')/(l-1)}f_2'$ by Sturm's bound in this case: if the
difference were non-zero, then the order of vanishing at $\infty$
contradicts the degree of the line bundle of which it is a section. If
$k_2'\geq k_1'$, then $f_2'=A^{(k_2'-k_1')/(l-1)}f_1'$ for the same
reason. We conclude that for all primes $p\neq l$ we have
$f_1(T_p)=p^if_2(T_p)$. 
\end{proof}
\begin{rem}
Proposition~\ref{prop_test_equality_res_repr} can be generalised to
forms of higher level, by the similar trick of passing to a higher
level $n$ at which one has forms that gave the same Galois
representation, but with eigenvalue $0$ for all $T_p$ with $p$
dividing~$n$. 
\end{rem}
The next result gives some conditions under which the Galois
representation $\rho$ attached to a surjective ring morphism
$f\colon\TT(1,k)\to\FF$ has large image in the sense that is
contains~$\SL_2(\FF)$. It is an effective version of Theorem~5.1
of~\cite{Ribet4}. We will need such a result later on.
\begin{thm}\label{thm_large_image}
Let $k$ be a positive integer, $l$ a prime number with $l>6(k-1)$,
$\FF$ a finite field and $l$ its characteristic, and
$f\colon\TT(1,k)\to\FF$ a surjective morphism of rings such that the
associated Galois representation
$\rho\colon\Gal(\Qbar/\QQ)\to\GL_2(\FF)$ is irreducible. Then the
image of $\rho$ contains~$\SL_2(\FF)$, and is equal to the subgroup of
$\GL_2(\FF)$ of elements $g$ whose determinant is in the subgroup of
$k{-}1$th powers in~$\FF_l^\times$.
\end{thm}
\begin{proof}
As $\TT(1,k)=0$ for $k<12$, we have $k\geq12$, and hence $l>66$. As
$l>2$ and $\rho$ is odd, $\rho$ is absolutely irreducible. We also
have $k\leq l+1$.

We apply what is known about the restriction of $\rho$ to an inertia
subgroup $I$ at~$l$. We denote by $\psi$ and $\psi':=\psi^l$ the two
fundamental characters from $I$ to $\FF_{l^2}^\times$ of level~2 (the
tame quotient of $I$ is the projective limit of the $\FF_{l^n}^\times$
and the fundamental characters of level~$n$ to $\Fbar_l^\times$ are
those that are induced by ring morphisms $\FF_{l^n}\to\Fbar_l$). By
Theorems~2.5 (due to Deligne) and~2.6 (due to Fontaine)
in~\cite{Edixhoven3}, we have:
\begin{align*}
\rho|_I & = 
\left(\begin{matrix} \chi_l^{k-1} & * \\ 0 & 1\end{matrix}\right) 
\quad \text{if $f(T_l)\neq 0$, and}\\
\Fbar_l\otimes_\FF\rho|_I & = 
\left(\begin{matrix} \psi^{k-1} & 0 \\ 0 & {\psi'}^{k-1}\end{matrix}\right) 
\quad \text{if $f(T_l)=0$.}
\end{align*}
The classification of subgroups of $\GL_2(\FF)$ of order prime to~$l$
(see for example~\cite[\S2.5, Prop.~16]{Serre6}) says that the image
in $\PGL_2(\FF)$ of such a subgroup is either cyclic, dihedral or
isomorphic to $A_4$, $S_4$ or~$A_5$.  

As $\rho$ is absolutely irreducible, its image in $\PGL_2(\FF)$ cannot
be cyclic (note that the kernel of $\GL_2(\FF)\to\PGL_2(\FF)$ is the
center of~$\GL_2(\FF)$). 

Let us show that the projective image of $\rho$ cannot be~$A_4$, $S_4$
or~$A_5$. Assume that it is. Then the image of $\rho(I)$ in
$\PGL_2(\FF)$ is cyclic and has order at least $(l-1)/(k-1)$ (the
order of $\psi/\psi'$ is $l+1$). As we assume that $l-1$ is at least
$6(k-1)$, this image has an element of order at least~$6$, a
contradiction.

Let us show that the projective image of $\rho$ is not
dihedral. Assume that it is. Then the image of $\rho$ is contained in
the normaliser of a Cartan subgroup (i.e., the group of points of a
split or non-split maximal torus), and there is a quadratic
extension~$K$ of~$\QQ$ such that $\rho$ is the induction from
$\Gal(\Qbar/K)$ to $\Gal(\Qbar/\QQ)$ of a character of~$\Gal(\Qbar/K)$
that is not equal to its conjugate under~$\Gal(K/\QQ)$. As $\rho$ is
unramified outside~$l$, $K$ must be the quadratic extension of~$\QQ$
that is ramified precisely at~$l$. As $l-1>2(k-1)$ the description
above of $\rho|_I$ shows that there are precisely two lines
in~$\Fbar_l^2$ whose orbit under $\rho$ in $\PP^1(\Fbar_l)$ has order
at most~2: these are the coordinate axes (in the first case, the
extension must be split). But in the first case the characters on
these two lines are not conjugate under~$\Gal(K/\QQ)$, and in the
second case the action of $\rho$ on the set of these two lines is
not ramified. These contradictions show that the image of $\rho$
cannot be dihedral.

We conclude that the order of the image of $\rho$ is divisible
by~$l$. As $l>3$ a result of Dickson, see~\cite{Dickson}, Chapter~XII,
or rather the proof of Theorem~2.5 in~\cite{Ribet2}, says that the
image of~$\rho$ in~$\PGL_2(\FF)$ is, after suitable conjugation, equal
to $\PGL_2(\FF')$ or $\SL_2(\FF')/\{1,-1\}$ for some subfield $\FF'$
of~$\FF$. 

We claim that $\FF'=\FF$. Assume that it is not. We let $f'$ and
$\rho'$ be the conjugates of $f$ and $\rho$ by the Frobenius
automorphism of $\FF$ over~$\FF'$. Then $\rho$ and $\rho'$ are not
isomorphic because the traces of the image of $\rho$ generate~$\FF$
(use that $\TT(1,k)$ is generated as $\ZZ$-module by the $T_i$ with
$i\leq k/12<l$). But their projective representations to $\PGL_2(\FF)$
are equal. Hence $\rho'$ is a twist of~$\rho$ by some character
$\chi\colon\Gal(\Qbar/\QQ)\to\FF^\times$. As $\rho$ and $\rho'$ are
unramified outside~$l$, $\chi$ is unramified outside~$l$ and hence a
power of~$\chi_l$. But then we have $\theta^{i+1}f=\theta f'$, with
$\theta$ as in the proof of
Proposition~\ref{prop_test_equality_res_repr}. A look at the theta
cycles in Section~3 of~\cite{Edixhoven3} or Section~7
of~\cite{Jochnowitz} shows that then $k=(l+3)/2$ if $f(T_l)=0$, and
$k=(l+1)/2$ if $f(T_l)\neq 0$. This contradicts our assumption that
$l>6(k-1)$.

So the image $G$ of $\rho$ in $\PGL_2(\FF)$ contains the image
of~$\SL_2(\FF)$. Then, for each $a\in\FF$, $G$ contains elements of
the form $(\begin{smallmatrix}t & a\\0 & t\end{smallmatrix})$ and
$(\begin{smallmatrix}s & 0\\a & s\end{smallmatrix})$, for some $t$
and $s$ in~$\FF^\times$. Taking suitable powers, we conclude that $G$
contains all 
$(\begin{smallmatrix}1 & a\\0 & 1\end{smallmatrix})$ and
$(\begin{smallmatrix}1 & 0\\a & 1\end{smallmatrix})$, where $a$ ranges
through~$\FF$. These generate~$\SL_2(\FF)$. As
$\det\rho=\chi_l^{k-1}$ the last claim in the theorem follows.
\end{proof}
\begin{rem}
Eigenforms $f\colon\TT(1,k)\to\FF$ such that the projective image of
$\rho$ is $A_4$, $S_4$ or $A_5$ are related to complex modular forms
of weight one and level~$l$ or~$l^2$, see~\cite{Khare-Wintenberger1},
Theorem~10.1. There are tables of these. For example,
in~\cite{Basmaji-Kiming}, page~110, one finds an $A_5$-example with
$l=2083$, and an $S_4$-example with $l=751$. See also Section~4.3
of~\cite{Kilford-Wiese}. We note that for $f=\Delta$ the prime $23$
with $23-1=2(12-1)$ nicely illustrates one of the arguments that is
used in the proof above: $\rho$ is then dihedral. More generally, $f$
with $\rho$ dihedral come from class groups of imaginary quadratic
orders that are unramified outside~$l$.
\end{rem}

\chapter{First description of the algorithms}\label{chap_first_descr}

\author{B. Edixhoven and J.-M. Couveignes}

\bigskip

\bigskip

%author Jean-Marc and Bas.

% Gives the main ideas and motivates the chapters that follow

We put ourselves in the situation of Theorem~\ref{thm_red_to_wt_2},
and we ask how we can compute the Galois representation. More
explicitly, let $n$ and $k$ be positive integers, $\FF$ a finite field
and $l$ its characteristic, and $f\colon\TT(n,k)\to\FF$ a surjective
ring morphism. Assume that $2<k\leq l{+}1$, and that the associated
Galois representation $\rho\colon\Gal(\Qbar/\QQ)\to\GL_2(\FF)$ is
absolutely irreducible. Let $f_2\colon\TT(nl,2)\to\FF$ be the weight
two eigenform as in Theorem~\ref{thm_red_to_wt_2} and let
$m=\ker(f_2)$. Assume that the multiplicity of $\rho$ in
$V:=J_1(nl)(\Qbar)[m]$ is one, i.e., that $\rho$ is realised by~$V$.

We let $K\subset\Qbar$ be the field ``cut out by~$\rho$'', i.e., the
finite Galois extension of $\QQ$ contained in $\Qbar$ consisting of
the elements of $\Qbar$ that are fixed by all elements
in~$\ker(\rho)$. Then we have, by definition, the following
factorisation of~$\rho$:
\[
\rho\colon \Gal(\Qbar/\QQ) \onto \Gal(K/\QQ)
\;\into\; \GL_2(\FF).
\]

Our aim is then to compute such residual representations~$\rho$, in
time polynomial in~$n$, $k$ and~$\#\FF$. By computing $\rho$ we mean
giving $K$ as a $\QQ$-algebra, in the form of a monic polynomial in
$\QQ[T]$ that is the minimal polynomial of some generator~$t$ of~$K$,
and giving the elements $\sigma$ of $\Gal(K/\QQ)$ by giving their
matrices with respect to the $\QQ$-basis of $K$ consisting of the
first so many powers of~$t$, together with the element $\rho(\sigma)$
of~$\GL_2(\FF)$. Once given such an explicit description of $\rho$ it
becomes possible to compute $f(T_p)\in \FF$ in deterministic
 polynomial  time in $\log p$. Indeed this boils down to computing
the Frobenius endomorphism at $p$ for the algebra $A$. 
Chapter~\ref{chap_comp_coefs} explains how to do this.

It will be convenient for us to use the modern version of Galois
theory that says that the functor $A\mapsto \Hom_\QQ(A,\Qbar)$ is an
anti-equivalence from the category of finite separable $\QQ$-algebras
to that of finite discrete (continuous) $\Gal(\Qbar/\QQ)$-sets. An
inverse is given by the functor that sends $X$ to
$\Hom_{\Gal(\Qbar/\QQ)}(X,\Qbar)$, the $\QQ$-algebra of functions $f$
from $X$ to~$\Qbar$ such that $f(gx)=g(f(x))$ for all $g$ in
$\Gal(\Qbar/\QQ)$ and all $x$ in~$X$. Under this correspondence,
fields correspond to transitive $\Gal(\Qbar/\QQ)$-sets. 

As a first step towards the computation of~$\rho$ we let $A$ be the
$\QQ$-algebra corresponding to the $\Gal(\Qbar/\QQ)$-set~$V$. Before
we explain our strategy to compute~$A$, we sketch how one gets from
$A$ to $K$ and~$\rho$. The $\QQ$-algebra corresponding to $V\times V$
is $A\otimes A$. The addition map $V\times V\to V$ corresponds to a
morphism $A\to A\otimes A$, the co-addition. The $\FF^\times$-action
on~$V$ corresponds to an $\FF^\times$-action on~$A$. We will see later
that the co-addition and the $\FF^\times$-action on $A$ can be
computed by the same method by which $A$ will be computed. Viewing
$V\times V$ as $\Hom_\FF(\FF^2,V)$ gives a right-action by
$\GL_2(\FF)$ on $V\times V$, hence a left-action on $A\otimes A$. This
action can be expressed in the co-addition and the
$\FF^\times$-action. Let $B$ be the $\QQ$-algebra corresponding to the
subset $\Isom_\FF(\FF^2,V)$ of $\Hom_\FF(\FF^2,V)$. This factor $B$ of
$A\otimes A$ can be computed by linear algebra over~$\QQ$, using the
$\GL_2(\FF)$-action on $A\otimes A$. In terms of $V\times V$, one
removes the subset of $(v_1,v_2)$ that are linearly dependent, i.e.,
the point $(0,0)$ and the $\GL_2(\FF)$-orbit of $(V-\{0\})\times
\{0\}$. The field $K$ then corresponds to a $\Gal(\Qbar/\QQ)$-orbit in
$\Isom(\FF^2,V)$, hence is obtained by factoring $B$ as a product of
fields, using factoring algorithms, and choosing one of the
factors. See~\cite{LLL}, \cite{Lenstra_Arjen} and \cite{Landau1} for
the fact that such factoring can be done in polynomial time. The
equivalence between factoring algebras and polynomials is given
in~\cite{Lenstra_Hendrik_1}. Let $G\subset \GL_2(\FF)$ be the
stabiliser of the chosen factor~$K$. Then $G=\Gal(K/\QQ)$ and the
inclusion $\emph{is}$ a representation from $G$ to $\GL_2(\FF)$. Let
$\phi$ be in the chosen $\Gal(\Qbar/\QQ)$-orbit in
$\Isom_\FF(\FF^2,V)$. As this orbit is a right $G$-torsor on which
$\Gal(\Qbar/\QQ)$ acts, there is, for every $\sigma$ in
$\Gal(\Qbar/\QQ)$, a unique $g(\sigma)$ in $G$ such that
$\rho(\sigma)\circ\phi=\phi\circ g(\sigma)$. Note also that evaluation
at $\phi$ is an embedding of $K$ in~$\Qbar$, such that $\sigma$ in
$\Gal(\Qbar/\QQ)$ induces $g(\sigma)$ on~$K$. It follows that $\phi$
is an isomorphism between $\rho$ and the representation $g\colon
\Gal(\Qbar/\QQ)\onto\Gal(K/\QQ)=G\subset \GL_2(\FF)$.

We now turn to the question of how to compute the $\QQ$-algebra $A$
corresponding to~$V$. We wish to produce a generator of $A$, and its
minimal polynomial over~$\QQ$. This means that we must produce a
$\Qbar$-valued function $a$ on $V$ such that $a(\sigma
x)=\sigma(k(x))$ for all $x$ in $V$ and all $\sigma$
in~$\Gal(\Qbar/\QQ)$. Such a function is a generator of $A$ if and
only if it does not arise from a strictly smaller quotient of $V$ as
$\Gal(\Qbar/\QQ)$-set (such quotients correspond to subalgebras),
hence, equivalently, if and only if $a$ is injective. The minimal
polynomial over $\QQ$ of such a generator $a$ is given as follows:
\begin{supeqn}\label{eqn_min_pol_a}
P(T) = \prod_{x\in V} (T-a(x)).
\end{supeqn}

The question is now how to produce such a generator? A direct way
would be to compute the elements of $V$ in $J_1(nl)(\Qbar)$, by
writing down polynomial equations in a suitable coordinate system that
is defined over~$\QQ$, and solving them, using computer algebra. This
is essentially how Schoof's algorithm deals with elliptic
curves. However, the dimension of $J_1(nl)$ is quadratic
in~$l$. Writing down equations in polynomial time still seems
possible. But we do not know of a way of solving the equations in a
time that is not exponential in the dimension.

The
decisive idea is to use numerical computations to approximate the
coefficients of a minimal polynomial $P$ as above, in combination with
a bound on the \emph{height} of those coefficients. We recall that the
(standard, logarithmic) height of a rational number $a/b$, with $a$
and $b$ integers that are relatively prime, is $\log \max\{|a|,|b|\}$
(a variant would be $\log(a^2+b^2)$). This rational number $x=a/b$ is
known if we know an upper bound $h$ for its height, and an
approximation $y$ of it (in $\RR$, say), with
$|x-y|<e^{-2h}/2$. Indeed, if $x'=a'/b'$ also has height at most~$h$,
and $x'\neq x$, then:
\[
|x-x'| = \left|\frac{a}{b}-\frac{a'}{b'}\right| = 
\left|\frac{ab'-ba'}{bb'}\right| \geq \frac{1}{|bb'|} \geq e^{-2h}.
\]
We also note that there are good algorithms to deduce $x$ from such a
pair of an approximation~$y$ and a bound~$h$, for example by using
continued fractions, as we will now explain. 

In practice we will use rational approximations $y$ of~$x$. Every
rational number $y$ can be written uniquely as:
\[
[a_0,a_1,\ldots,a_n] = a_0+\cfrac{1}{a_1+
\cfrac{1}{\ddots\genfrac{}{}{0pt}{0}{}{a_{n-1}+\cfrac{1}{a_n}}}}
\]
where $n\in\ZZ_{\geq0}$, $a_0\in\ZZ$, $a_i\in\ZZ_{>0}$ for all $i>0$,
and $a_n>1$ if $n>0$. To find these $a_i$, one defines $a_0:=\lfloor
y\rfloor$ and puts $n=0$ if $y=a_0$; otherwise, one puts
$y_1:=1/(y-a_0)$ and $a_1=\lfloor y_1\rfloor$ and $n=1$ if $y_1=a_1$,
and so on. The rational numbers $[a_0,a_1,\ldots,a_i]$ with
$0\leq i\leq n$ are called the \emph{convergents} of the continued
fraction of~$y$. Then one has the following well known result,
see Theorem~184 from~\cite{Hardy-Wright}.
\begin{supprop}\label{prop_cont_frac}
Let $y$ be in $\QQ$, $a$ and $b$ in~$\ZZ$ with $b\neq 0$ and:
\[
\left|\frac{a}{b}-y\right|<\frac{1}{2b^2} .
\]
Then $a/b$ is a convergent of the continued fraction of~$y$. 
\end{supprop}

The question is now: how we are going to implement this method?
The basic idea in doing this is to not work on the Abelian
variety $J_1(nl)$ but rather on the product $X_1(nl)^g$ of copies
of~$X_1(nl)$, where $g$ is the genus of~$X_1(nl)$. To compare the
two, we first choose an effective divisor~$D_0=P_1+\cdots+P_g$ on
$X_1(nl)_\QQ$, and we consider the well-known map:
\begin{supeqn}\label{eqn_abel_jacobi}
\begin{aligned}
X_1(nl)^g & \lto J_1(nl), \\
(Q_1,\ldots,Q_g) & \mapsto [Q_1+\cdots+Q_g-D_0].  
\end{aligned}
\end{supeqn}

To understand the definition of this map rigorously, one must use the
interpretation of $X_1(nl)$ as its functor of points with values in
$\ZZ[1/nl]$-schemes, and that of $J_1(nl)$ as the degree zero part of the
relative Picard functor~$\Pic^0_{X_1(nl)/\ZZ[1/nl]}$. For the
necessary background on this, see Chapters~8 and~9 of~\cite{BLR1}. 
The divisor $D_0$ lives on~$X_1(nl)_\QQ$, and it extends uniquely over
$\ZZ[1/nl]$ to an effective relative Cartier divisor of degree $g$
on~$X_1(nl)$. The points $P_i$ of which $D_0$ is the sum need not be
rational over~$\QQ$. 

The inverse image of a point $x$ in $J_1(nl)(\Qbar)$ under the
map~(\ref{eqn_abel_jacobi}) can be described as follows. Let $\calL_x$
denote a line bundle of degree zero on $X_1(nl)_\Qbar$ that
corresponds to~$x$ ($x$ is an isomorphism class of such line
bundles). Then the inverse image of~$x$ is the set of
$(Q_1,\ldots,Q_g)$ such that $\calL_x$ has a rational section
whose divisor is $Q_1+\cdots+Q_g-D_0$, or, equivalently, the set of
$(Q_1,\ldots,Q_g)$ such that there is a non-zero section of
$\calL_x(D_0)$ with divisor $Q_1+\cdots+Q_g$.

When $x$ ranges over $J_1(nl)(\Qbar)$, the class of the $\calL_x(D_0)$
ranges over the set~$\Pic^g(X_1(nl)_\Qbar)$. The function on
$J_1(nl)(\Qbar)$ that assigns to $x$ the dimension $h^0(\calL_x(D_0))$
of the space of global sections of~$\calL_x(D_0)$ is semi-continuous in
the sense that for each $i$ the locus of $x$ where
$h^0(\calL_x(D_0))\geq i$ is closed (the condition $h^0(\calL_x(D_0))\leq
i$ need not be closed). On a non-empty open subset of $J_1(nl)(\Qbar)$
this value is one, as can be seen using the theorem of Riemann-Roch,
and Serre duality.  This means that for $x$ outside a proper closed
subset of $J_1(nl)(\Qbar)$, the inverse image in
$X_1(nl)^g(\Qbar)$ of $x$ consists of the $g$-tuples obtained by
permutation of coordinates of a single $(Q_1,\ldots,Q_g)$. Another
way to express this is to say that the map~(\ref{eqn_abel_jacobi})
above factors through the symmetric product $X_1(nl)^{(g)}$ and that
the map from $X_1(nl)^{(g)}$ to $J_1(nl)$ is birational (i.e., an
isomorphism on suitable non-empty open parts).

It is then reasonable to assume that we can take $D_0$ such that for
all $x$ in~$V$ there is, up to permutation of the coordinates, a
unique $Q=(Q_1,\ldots,Q_g)$ in $X_1(nl)^g(\Qbar)$ that is mapped to
$x$ via the map~(\ref{eqn_abel_jacobi}). On the other hand, on a curve
of high genus such as $X_1(nl)$ it is not clear how to make a large
supply of inequivalent effective divisors $D_0$ on~$X_1(nl)_\QQ$. We
will see later, in Theorem~\ref{thm_exist_D}, that we can indeed find
a suitable divisor, supported on the cusps, and defined
over~$\QQ(\zeta_l)$, on the~$X_1(5l)$, which will suffice for treating
almost all modular forms of level one.  

\begin{remark}
In situations where such a
cuspidal divisor cannot be found, one could try at random
$P_1,\ldots,P_g$ in~$X_1(nl)(L)$, corresponding to elliptic curves
lying in one isogeny class, with complex multiplications, for example
by~$\QQ(i)$. Then $L$ is a solvable Galois extension of~$\QQ$, so that
$K$ can be reconstructed from the compositum~$KL$. If one chooses
the~$P_i$ reasonably, the degree of~$L$ and the logarithm of the
discriminant of~$L$ are polynomial in~$l$. Another possibility is to
try to work with a divisor $D_0$ of degree smaller than~$g$, for
example a multiple of a rational cusp.
\end{remark}

Let us now assume that we have a divisor $D_0$ as described above. Then
we choose a non-constant function:
\[
f\colon X_1(nl)_\QQ \onto \PP^1_\QQ,
\]
that will have to satisfy some conditions that will be given in a
moment.

With these two choices, $D_0$ and $f$, and a choice of an integer~$m$,
we get an element $a_{D_0,f,m}$ of the $\QQ$-algebra $A$ corresponding
to~$V$ as follows. For $x$ in $V$ we let
$D_x=Q_{x,1}+\cdots+Q_{x,g}$ be the unique effective divisor of
degree $g$ such that:
\begin{supeqn}\label{eqn_def_D'_x}
x=[D_x-D_0].
\end{supeqn}
Note that indeed for $x=0$ we have $D_x=D_0$.  We assume that for all
$x\in V$ the divisor $D_x$ is disjoint from the poles of~$f$. Then,
for each $x$ in~$V$, we define:
\[
P_{D_0,f,x} = \prod_{i=1}^g\bigl(t-f(Q_{x,i})\bigr) 
\quad\text{in $\Qbar[t]$.}
\]
We then get an element $a_{D_0,f,m}$ of $A$ by evaluating the
$P_{D_0,f,x}$  at~$m$:
\begin{supeqn}\label{eqn_def_generator}
a_{D_0,f,m}\colon V\lto \Qbar, \quad
x\mapsto P_{D_0,f,x}(m) = \prod_{i=1}^g\bigl(m-f(Q_{x,i})\bigr).
\end{supeqn}
The condition that all the $D_x$ are disjoint from the locus of poles
of $f$ will not be guaranteed to hold later when we treat forms of
level one, but then it will be possible to omit the $Q_{x,i}$ at which
$f$ has a pole from the sum in~(\ref{eqn_def_generator}) ($f$ will
have its poles at certain cusps). For the moment, let us just assume
that this condition is satisfied. Then the $f_*D_x$, for $x$ in~$V$,
are effective divisors of degree~$g$ on~$\AA^1_\Qbar$.

We will choose $f$ in such a way that the $f_*D_x$ are distinct; we
assume now that this is so. Then there is an integer $m\geq0$ with
$m\leq g{\cdot}(\#\FF)^4$ such that $a_{D_0,f,m}$ is injective, and
hence a generator of~$A$: the polynomials $P_{D_0,f,x}$ are distinct
when $x$ varies, and $m$ must not be a root of any difference of two
of them.

Finally, we want to have control on the heights of the coefficients of
the minimal polynomial of~$a_{D_0,f,m}$, because these heights
determine the required precision of the approximations of those
coefficients that we must compute. The whole strategy depends on the
possibility to choose a divisor~$D_0$ and a function~$f$, such that,
when $n$, $k$ and~$\FF$ vary, those heights grow at most polynomially
in $n$, $k$ and~$\#\FF$. Using a great deal of machinery from Arakelov
theory, we will show (at least in the case $n=1$) that any reasonable
choices of~$D_0$ and~$f$ will lead to an at most polynomial growth of
those heights. Intuitively, and completely non-rigorously, one can
believe that this should work, because of the following argument. Our
$x$ are torsion points, so that their Néron-Tate height is zero. As
$x$ and $D_0$ determine~$D_x$, the height of $D_x$ should be not much
bigger than the height of~$D_0$. As we choose $D_0$ ourselves, it
should have small height. Finally, the height of $a_{D_0,f,m}$ should
be not much bigger than the sum of those of $f$ and $m$ and
the~$D_x$. Turning these optimistic arguments into rigorous statements
implies a lot of work that will be done in
Chapters~\ref{sec_appl_Arakelov}--\ref{chap_bnd_height}.  An important
problem here is that in Arakelov theory many results are available
that deal with a single curve over~$\QQ$, but in our situation we are
dealing with the infinitely many curves $X_1(nl)$ as $l$ varies.

A few words about the numerical computations involved. What we need is
that these can be done in a time that is polynomial in $n$ and $\#\FF$
and the number of significant digits that one wants for the
coefficients of the minimal polynomial~$P_{D_0,f,m}$
of~$a_{D_0,f,m}$. It is not at all obvious that this can be done, as
the genus of $X_1(nl)$ and hence the dimension of~$J_1(nl)$ are
quadratic in~$l$.

One way to do the computations is to use the complex uniformisations
of $X_1(nl)(\CC)$ and~$J_1(nl)(\CC)$. The Riemann surface
$X_1(nl)(\CC)$ can be obtained by adding finitely many cusps (the set
$\Gamma_1(nl)\backslash\PP^1(\ZZ)$) to the quotient
$\Gamma_1(nl)\backslash\HH$ (see Section~\ref{sec_modcurves}). This
means that $X_1(nl)(\CC)$ is covered by disks around the cusps, which
are well suited for computations (functions have $q$-expansions, for
example). In order to describe $J_1(nl)(\CC)$ as $\CC^g$ modulo a
lattice, we need a basis of the space of holomorphic differential
forms $\rH^0(X_1(nl)(\CC),\Omega^1)$. The basis that we work with is
the one provided by Atkin-Lehner theory, as given
in~(\ref{eqn_basis}); we write it as
$\omega=(\omega_1,\ldots,\omega_g)$. Then we have the following
complex description of the map~(\ref{eqn_abel_jacobi}):
\begin{supeqn}
\xymatrix{
X_1(nl)(\CC)^g \ar@{->>}[r] & J_1(nl)(\CC) \ar@{=}[r] 
& \CC^g/\Lambda \\
(Q_1,\ldots,Q_g) \ar@{|->}[r] 
& [\sum_{i=1}^gQ_i-\sum_{i=1}^gP_i] \ar@{=}[r] & 
\sum\limits_{i=1}^g\int\limits_{P_i}^{Q_i}\omega, 
}
\end{supeqn}
where $\Lambda$ is the period lattice with respect to this basis,
i.e., the image of $\rH_1(X_1(nl)(\CC),\ZZ)$ under integration of
the~$\omega_i$. This map can be computed up to any desired precision
by formal integration of power series on the disks mentioned
above. The coefficients needed from the power series expansions of the
$\omega_i$ can be computed using the method of modular symbols, as has
been implemented by William Stein in Magma (see his
book~\cite{Stein}).  We note that modular symbols algorithms can be
used very well to locate $V$ inside $l^{-1}\Lambda/\Lambda$, hence
in~$J_1(nl)(\CC)$.  A strategy to approximate a point
$Q_x=(Q_{x,1},\ldots,Q_{x,g})$ as above for a non-zero $x$ in~$V$ is
to lift the straight line that one can draw in $\CC^g/\Lambda$ from
$0$ to $x$ (within a suitable fundamental domain for~$\Lambda$) to a
path in $X_1(nl)(\CC)^g$ starting at~$(P_1,\ldots,P_g)$. In practice
this seems to work reasonably well, see Bosman's
Chapters~\ref{chapcomput} and~\ref{chappgltau}.  A theoretical
difficulty with this approach is that one needs to bound from below
the distance to the ramification locus of $X_1(nl)^{(g)}\to J_1(nl)$.
Chapter~\ref{sec_couveignes_TORSION} gets around this 
difficulty and provides a proven algorithm
for inverting the Jacobi map~(\ref{eqn_abel_jacobi}). 
The starting idea is to set $y=x/N$ for $N$ a large enough
integer. This $y$ is no longer an $l$-torsion point
but it is close to the origin in the
torus  $J_1(nl)(\CC)$,  and this helps finding a preimage
 $Q_y$ of $y$,  because the  behaviour of the Jacobi 
map~(\ref{eqn_abel_jacobi}) is well understood at least in the neighborhood 
of the origin. The divisor $Q_x$ we are looking for is such that
$Q_x-D_0$ and $N(Q_y-D_0)$ are linearly equivalent. So $Q_x$ can be computed
from $Q_y$ by repeated application of an explicit form of the
Riemann-Roch theorem. The resulting algorithm reduces to computing
approximations of the complex zeros of a great number of modular 
forms with level $5l$ and weight $4$. Chapter \ref{sec_couveignes_ZEROS}
explains how to approximate the complex
zeros of entire series. It also contains a reminder of the
necessary notions from computational complexity theory.

Another way to do the ``approximation'' is to compute the minimal
polynomial $P_{D_0,f,m}$ of $a_{D_0,f,m}$ modulo many small
primes~$p$. Indeed, the map~(\ref{eqn_abel_jacobi}) can be reduced
mod~$p$. In this case one has no analytic description of the curve and
its Jacobian, but one can make random points in $J_1(nl)(\FF_q)$ for a
suitable finite extension $\FF_p\to\FF_q$. Such random points can then
be projected, using Hecke operators, into~$V$. Elements of
$J_1(nl)(\FF_q)$ can be represented by divisors on $X_1(nl)_{\FF_q}$,
and all necessary operations can be done in polynomial time. This
approach is explained in detail in Chapter~\ref{sec_couveignes_modp}. 
In order to deduce a rational
number $x=a/b$ from the knowledge sufficiently many of its reductions
modulo primes $p$ not dividing $b$ we have the following well-known
result.

\begin{supprop}\label{prop_x_from_reductions}
Let $x=a/b$ be in~$\QQ$, with $a$ and $b$ in~$\ZZ$, relatively
prime. Let $M=\max\{|a|,|b|\}$. Let $S$ be a finite set of prime
numbers $p$ with $p$ not dividing~$b$, such that $\prod_{p\in
  S}p>2M^2$. For each $p$ in~$S$, let $x_p$ in~$\FF_p$ be the reduction
of~$x$, and let $L\subset \ZZ^2$ be the submodule of $(n,m)$ with the
property that for all $p$ in~$S$: $n-x_pm=0$ in~$\FF_p$. Then $(a,b)$
and $(-a,-b)$ are the shortest non-zero elements of~$L$ with respect
to the standard inner product on~$\RR^2$, and the lattice reduction
algorithm in dimension two, Algorithm~1.3.14 in~\cite{Cohen}, finds these
in time polynomial in~$\log M$.
\end{supprop}
\begin{proof}
The lattice reduction gives a shortest non-zero element, so it
suffices to show that, under the assumptions in the Proposition, the
two shortest non-zero elements of~$L$ are precisely $\pm(a,b)$. The
volume of $\RR^2/L$ is the index of $L$ in $\ZZ^2$, hence equals
$\prod_{p\in S}p$. Let $l_1$ be a shortest non-zero element
of~$L$. Then $\|l_1\|\leq\|(a,b)\|\leq\sqrt{2}M$. Let $l_2$ in $L$
be linearly independent of~$l_1$. Then:
\[
2M^2<\Vol(\RR^2/L)\leq \Vol(\RR^2/(\ZZ{\cdot}l_1+\ZZ{\cdot}l_2)) 
\leq \|l_1\|{\cdot}\|l_2\|\leq\sqrt{2}M\|l_2\|.
\]
Hence $\|l_2\|>\sqrt{2}M\geq\|(a,b)\|$. It follows that $(a,b)$ and
$l_1$ are linearly dependent, and hence $l_1=\pm(a,b)$.
\end{proof}

\begin{remark} In case one has a natural rigid analytic uniformisation at some
prime~$p$, one may want to use that. For the modular curves that we
are dealing with this is not the case, but the closely related Shimura
curves attached to quaternion algebras over~$\QQ$ do admit such
uniformisations at the primes where the quaternion algebra is ramified
(as was proved by Cerednik, Drinfeld, see~\cite{Boutot-Carayol}).
\end{remark}

\chapter{Short introduction to heights and Arakelov theory}
\label{sec_intro_heights_and_ar}

\author{B. Edixhoven and R. de Jong}

\bigskip

\bigskip

%authors: Bas and Robin

In Chapter~\ref{chap_first_descr} it has been explained how the
computation of the Galois representations $V$ attached to modular
forms over finite fields should proceed. The essential step is to
approximate the minimal polynomial $P$ of~(\ref{eqn_min_pol_a}) with
sufficient precision so that $P$ itself can be obtained. The topic
to be addressed now is to bound from above the precision that is
needed for this. This means that we must bound the heights of the
coefficients of~$P$. As was hinted to in Chapter~\ref{chap_first_descr},
we get such bounds using Arakelov theory, a tool that we discuss in
this section. It is not at all excluded that a direct approach to
bound the coefficients of~$P$ exists, thus avoiding the complicated
theory that we use. On the other hand, it is clear that the use of
Arakelov theory provides a way to split the work to be done in smaller
steps, and that the quantities occurring in each step are intrinsic in
the sense that they do not depend on coordinate systems or other
choices that one could make. We also want to point out that our method
does not depend on cancellations of terms in the estimates that we
will do; all contributions encountered can be bounded appropriately.

A good reference for a more detailed introduction to heights is Chapter~6
of~\cite{Cornell-Silverman}. Good references for the Arakelov theory
that we will use are~\cite{Faltings1} and~\cite{Moret-Bailly1}. A
general reference for heights in the context of Diophantine geometry
is~\cite{Bombieri-Gubler}.

\section{Heights on $\QQ$ and $\Qbar$}
\label{sec_heights}
The definition of the height of an element of~$\QQ$ has already been
given in Chapter~\ref{chap_first_descr}; for $x=a/b$ with $a$ and $b\neq0$
relatively prime integers, we have $h(x)=\log\max\{|a|,|b|\}$. We will
now give an equivalent definition in terms of absolute values
$|{\cdot}|_v$ on~$\QQ$ attached to all places~$v$ of~$\QQ$, the
finite places, indexed by the prime numbers, and the infinite place
denoted~$\infty$. 

The absolute value $|{\cdot}|_\infty$ is just the usual absolute value
on~$\RR$, restricted to~$\QQ$. We note that $\RR$ is the completion
of~$\QQ$ for~$|{\cdot}|_\infty$. For $p$ prime, we let $v_p$ be the
$p$-adic valuation:
\begin{eqn}
v_p\colon \ZZ \lto \ZZ\cup\{\infty\},
\end{eqn}
sending an integer to the maximal number of times that it can be divided
by~$p$. This valuation $v_p$ extends uniquely to~$\QQ$ subject to the
condition that $v_p(xy)=v_p(x)+v_p(y)$; we have
$v_p(a/b)=v_p(a)-v_p(b)$ for integers $a$ and $b\neq0$. We let
$|{\cdot}|_p$ denote the absolute value on~$\QQ$ defined by:
\begin{eqn}
|x|_p = p^{-v_p(x)},\quad |0|_p=0.
\end{eqn}
The completion of $\QQ$ with respect to~$|{\cdot}|_p$ is the locally
compact topological field~$\QQ_p$.  An important property of these
absolute values is that all together they satisfy the product formula:
\begin{eqn}
\prod_v |x|_v = 1, \quad \text{for all $x$ in~$\QQ^\times$}.
\end{eqn}
With these definitions, we have:
\begin{eqn}
h(x) = \sum_v \log\max\{1,|x|_v\}, \quad \text{for all $x\in\QQ$},
\end{eqn}
where $v$ ranges over the set of all places of~$\QQ$ (note that almost
all terms in the sum are equal to~$0$).

The height function on~$\QQ$ generalises as follows to number
fields. First of all, for a local field $F$ we define the natural
absolute value~$|{\cdot}|$ on it by letting, for $x$ in~$F^\times$,
$|x|_F$ be the factor by which all Haar measures on~$F$ are scaled by
the homothecy $y\mapsto xy$ on~$F$. For example, for $F=\CC$ we have
$|z|_\CC=z\ol{z}=|z|^2$, the square of the usual absolute value. Let
now $K$ be a number field. By a finite place of~$K$ we mean a maximal
ideal of~$O_K$. An infinite place of $K$ is an embedding of $K$
into~$\CC$, up to complex conjugation. For each place~$v$ of~$K$, let
$K_v$ be its completion at~$v$; as $K_v$ is a local field, we have the
natural absolute value $|{\cdot}|_v:=|{\cdot}|_{K_v}$ on~$K_v$ and
on~$K$. In this case, the product formula is true (this can be shown
easily by considering the adèles of~$K$, see Chapter~IV, Section~4,
Theorem~5 of~\cite{Weil2}). The height function on~$\QQ$ also
generalises to~$K$. For all for all $x$ in~$K$ we define:
\begin{eqn}
\begin{aligned}
h_K(x) := & \sum_v \log\max\{1,|x|_v\} = \\
= & \sum_{\text{$v$ finite}}\log\max\{1,|x|_v\} 
+ \sum_{\sigma\colon K\to\CC}\log\max\{1,|\sigma(x)|\}.
\end{aligned}
\end{eqn}
This function $h_K$ is called the \emph{height function of~$K$}. For
$K\to K'$ an extension of number fields, and for $x$ in~$K$, we have
$$h_{K'}(x)=(\dim_{K}K'){\cdot}h_K(x).$$ Therefore one has the
\emph{absolute} height function~$h$ on~$\Qbar$ defined by:
\begin{eqn}\label{eqn_height_on_qbar}
h\colon \Qbar\to\RR,\quad h(x) = \frac{h_K(x)}{\dim_\QQ K},
\end{eqn}
where $K\subset\Qbar$ is any number field that contains~$x$.

\section{Heights on projective spaces and on varieties}
\label{sec_heights_proj_sp}
For $n\geq0$ and for $K$ a number field, we define a height function
on the projective space $\PP^n(K)$ by:
\begin{eqn}\label{eqn_height_Pn}
\begin{aligned}
h_K((x_0:\cdots:x_n)) &:= \sum_v \log\max\{|x_0|,\ldots,|x_n|\},\\
h(x) &:= \frac{h_K(x)}{\dim_\QQ K},
\end{aligned}
\end{eqn}
where $v$ ranges through the set of all places of~$K$. We note that it
is because of the product formula that $h_K(x)$ is well-defined, and
that this definition is compatible with our earlier definition of the
height and absolute height on~$K$ if we view $K$ as the complement
of~$\infty$ in~$\PP^1(K)$. The functions $h$ on~$\PP^n(K)$ for
varying~$K$ naturally induce the absolute height function
on~$\PP^n(\Qbar)$.

A fundamental result, not difficult to prove, but too important to
omit here (even though we will not use it), is Northcott's finiteness
theorem.
\begin{thm}[Northcott]
\label{thm_Northcott}
Let $n$, $d$ and $C$ be integers. Then:
\[
\{x\in\PP^n(\Qbar)\;|\; \text{$h(x)\leq C$ and $\dim_\QQ(\QQ(x))\leq d$}\}
\]
is a finite set.
\end{thm}
For a proof the reader is referred to Chapter~6
of~\cite{Cornell-Silverman}, or to Section~2.4 of~\cite{Serre4}.

For any algebraic variety $X$ embedded in a projective space~$\PP^n_K$
over some number field~$K$, we get height functions~$h_K$ on~$X(K)$
and~$h$ on~$X(\Qbar)$ by restricting those from~$\PP^n$ to~$X$.

For later use, we include here some simple facts.  The height
functions on the projective spaces $\PP^n(\Qbar)$ are compatible with
embeddings as coordinate planes, for example by sending
$(x_0:\cdots:x_n)$ to $(x_0:\cdots:x_n:0)$, or to
$(0:x_0:\cdots:x_n)$.

For all $n\in\NN$, we view $\AA^n(\Qbar)$ as a subvariety of
$\PP^n(\Qbar)$, embedded in one of the $n+1$ standard ways by sticking
in a $1$ at the extra coordinate. For example, by sending
$(x_1,\ldots,x_n)$ to $(1:x_1:\cdots:x_n)$. This gives us, for
each~$n$, a height function $h\colon\AA^n(\Qbar)\to\RR$. These height
functions are also compatible with embeddings as coordinate
planes. For $n=1$ the height function on $\AA^1(\Qbar)=\Qbar$ is the
function in~(\ref{eqn_height_on_qbar}).
\begin{lem}\label{lem_sum_and_product}
Let $n\in\ZZ_{\geq 1}$, and $x_1,\ldots,x_n$ in~$\Qbar$. Then:
\begin{align*}
h(x_1\cdots x_n) &\leq\sum_{i=1}^n h(x_i)\, ,\quad
h(x_1+ \cdots + x_n) \leq \log n + \sum_{i=1}^n h(x_i) \, \\
h(x_i)& \leq h((x_1,\ldots,x_n)) \leq
h(x_1)+\cdots+h(x_n)\quad\text{for each~$i$}.
\end{align*}
\end{lem}
\begin{proof}
Let $n\in\ZZ_{\geq 1}$, and $x_1,\ldots,x_n$ in~$\Qbar$. Let
$K\subset\Qbar$ be a finite extension of~$\QQ$ containing the~$x_i$.
For the first inequality, we have:
\begin{align*}
h_K(x_1\cdots x_n) & = \sum_v \log\max\{1,|x_1\cdots x_n|_v\} \\
& \leq \sum_v\log(\max\{1,|x_1|_v\}\cdots\max\{1,|x_n|_v\}) \\
& = h_K(x_1)+\cdots+h_K(x_n)\, .
\end{align*}
For the second inequality, let $x=x_1+\cdots+x_n$. Then we have:
\begin{align*}
h_K(x) & = 
\sum_{\text{$v$ finite}} \log\max\{1,|\sum_i x_i|_v\} +
\sum_{\sigma\colon K\to\Qbar}\log\max\{1,|\sigma\sum_i x_i|\} \\
&\leq \sum_{\text{$v$ finite}} \log\max\{1,\max_i|x_i|_v\} + 
\sum_\sigma \log\max\{1,\sum_i|\sigma(x_i)|\}\\
&\leq \sum_{\text{$v$ finite}}\sum_i\log\max\{1,|x_i|_v\} + 
\sum_\sigma \max\{0,\log(n{\cdot}\max_i|\sigma(x_i)|)\} \\
&\leq \sum_{\text{$v$ finite},i}\log\max\{1,|x_i|_v\} + 
\sum_\sigma \log n 
+ \sum_{\sigma,i}\max\{0,\log|\sigma(x_i)|\} \\
&= (\dim_\QQ K){\cdot}\log n + \sum_i h_K(x_i)\, .
\end{align*}
For the third inequality, let $i$ be in~$\{1,\ldots,n\}$. We have:
\begin{eqnarray*}
h_K(x_1,\ldots,x_n) & =& h_K(1:x_1:\cdots:x_n) \\ &= &
\sum_v \log\max\{1,|x_1|_v,\ldots,|x_n|_v\} \\
& \geq& \sum_v \log\max\{1,|x_i|_v\} = h_K(x_i)\, .
\end{eqnarray*}
Finally, for the last inequality:
\begin{eqnarray*}
h_K(x_1,\ldots,x_n) & =& h_K(1:x_1:\cdots:x_n) \\&= &
\sum_v \log\max\{1,|x_1|_v,\ldots,|x_n|_v\} \\
& \leq &\sum_v \sum_i\log\max\{1,|x_i|_v\} = 
h_K(x_1)+\cdots+h_K(x_n)\, .
\end{eqnarray*}
\end{proof}

\begin{lem}\label{lem_bound_sym}
Let $d\geq1$ and $n\geq d$ be integers. Let $\Sigma_d$ denote the
elementary symmetric polynomial of degree~$d$ in $n$ variables. Let
$y_1,\ldots,y_n$ be in~$\Qbar$. Then we have:
\[
h\left(\Sigma_d(y_1,\ldots,y_n)\right) \leq 
n\log 2 + \sum_{1\leq i\leq n} h(y_i).
\]
\end{lem}
\begin{proof}
Let $K$ be the compositum of the fields $\QQ(y_i)$ for
$i=1,\ldots,n$. For each place $v$ of $K$, we let $|\cdot|_v$ be the
natural absolute value on $K_v$ and on $K$ as at the end of
Section~\ref{sec_heights}. By the triangle inequality we obtain, for
each place $v$ of~$K$:
\begin{align*}
|\Sigma_d(y_1,\ldots,y_n)|_v & \leq c(v,n) 
\max_{1\leq i_1 < \cdots < i_d \leq n} |y_{i_1} \cdots y_{i_d} |_v
\leq \\
& \leq c(v,n) \prod_{1\leq i\leq n} \max\{1,|y_i|_v\}\, ,
\end{align*}
where $c(v,n)=2^n$ if $v$ is Archimedean, and $c(v,n)=1$ if $v$ is
non-Archimedean. It follows that:
\[
\max \{ 1,|\Sigma_d(y_1,\ldots,y_n)|_v\} 
\leq c(v,n) \prod_{i=1}^n \max \{1,|y_i|_v \}\, .
\]
The proof of the lemma is finished by taking logarithms, summing over
the places $v$, and dividing by~$\dim_\QQ K$.
\end{proof}
\begin{lem}\label{lem_lower_bd_nonzero_elmt}
Let $x\neq0$ be in $\Qbar$, let $K=\QQ(x)$, and let $\sigma\colon
K\to\CC$. Then:
\[
|\sigma(x)| \geq e^{-(\dim_\QQ K){\cdot}h(x)}.
\]
\end{lem}
\begin{proof}
We have:
\[
\begin{aligned}
-\log|\sigma(x)| & = \log\left(|\sigma(x)|^{-1}\right) 
\leq \log\max\{1,|\sigma(x)|^{-1}\} \\
& \leq \sum_v \log\max\{1,|x|_v^{-1}\} = (\dim_\QQ K){\cdot}h(x^{-1}) \\
& = (\dim_\QQ K){\cdot}h(x),
\end{aligned}
\]
where the sum is over all places of~$K$.
\end{proof}

\begin{lem}\label{lem_bound_h_det}
Let $K$ be a number field, let $n$ be in $\ZZ_{\geq1}$, and let $a$ be
in $\rM_n(K)$. Then:
\[
h(\det(a))\leq \sum_{i,j}h(a_{i.j}) + \frac{1}{2}n\log n.
\]  
\end{lem}
\begin{proof}
Let $v$ be a finite place of~$K$. Then we have:
\[
\begin{aligned}
|\det(a)|_v & = \left|\sum_{s\in\rS_n}a_{1,s(1)}\cdots a_{n,s(n)}\right|_v 
\leq \max_s|a_{1,s(1)}|_v\cdots|a_{n,s(n)}|_v \\
& \leq \prod_{i,j}\max\{1,|a_{i,j}|_v\}.
\end{aligned}
\]
For $\sigma\colon K\to\CC$ we have, by Hadamard's inequality and the
comparison $\|{\cdot}\|\leq n^{1/2}\|{\cdot}\|_{\mathrm{max}}$
in~$\CC^n$ of the euclidean norm and the max-norm:
\[
\begin{aligned}
|\det\sigma(a)| & \leq\prod_j\|\sigma(a_j)\| 
= \prod_j \left(n^{1/2}\|\sigma(a_j)\|_{\mathrm{max}}\right) \\
& \leq n^{n/2}\prod_{i,j} \max\{1,|\sigma(a_{i,j})|\},
\end{aligned}
\]
where $a_j$ is the $j$th column of~$a$. Then we have:
\[
\begin{aligned}
h_K(\det(a)) & = \sum_v \log\max\{1,|\det(a)|_v\} +
\sum_\sigma\log\max\{1,|\det\sigma(a)|\} \\
& \leq \sum_v\sum_{i,j}\log\max\{1,|a_{i,j}|_v\} \\
& \quad + \sum_\sigma\left(\frac{n}{2}\log n +
\sum_{i,j}\log\max\{1,|\sigma(a_{i,j})|\}\right)\\
& = \frac{1}{2}(\dim_\QQ K){\cdot}n\log n + \sum_{i,j}h_K(a_{i,j}).
\end{aligned}
\]
Dividing by $\dim_\QQ K$ gives the result.
\end{proof}
\begin{lem}\label{lem_bound_h_sol_lin_eq}
Let $K$ be a number field, let $n$ be in $\ZZ_{\geq1}$, let $a$ be in
$\GL_n(K)$, and $y$ in~$K^n$. Let $x$ be the unique element in $K^n$
such that $ax=y$. Let $b$ be the maximum of all $h(a_{i,j})$
and~$h(y_i)$. Then we have, for all~$i$:
\[
h(x_i) \leq  2n^2b+n\log n.
\]
\end{lem}
\begin{proof}
We apply Cramer's rule: $x_i=\det(a(i))/\det(a)$, where
$a(i)$ in $\rM_n(K)$ is obtained by replacing the $i$th column
by~$y$. Lemma~\ref{lem_bound_h_det} gives us:
\[
h(\det(a)) \leq n^2b+\frac{1}{2}n\log n,
\quad h(\det(a(i))) \leq n^2b+\frac{1}{2}n\log n.
\]
Therefore: $h(x_i) = h(\det(a(i))/\det(a))\leq 2n^2b+n\log n$.
\end{proof}

\section{The Arakelov perspective on height functions}
We have just defined height functions $h_K$ and~$h$ on a variety~$X$
over a number field~$K$, embedded into some projective
space~$\PP^n_K$. Such an embedding determines a line bundle~$\calL$
on~$X$: the restriction of the line bundle $\calO(1)$ of~$\PP^n_K$
that corresponds to homogeneous forms of degree~$1$, in the variables
$x_0,\ldots,x_n$, say. The embedding of $X$ into~$\PP^n_K$ is given by
the global sections $s_0,\ldots,s_n$ of~$\calL$ obtained by
restricting the global sections $x_0,\ldots,x_n$ to~$X$. Now any
finite set of generating global sections $t_0,\ldots,t_m$ of~$\calL$
determines a morphism $f\colon X\to\PP^m_K$, inducing height functions
$h_{K,f}$ and~$h_f$ via pullback along~$f$. For $f$ and $f'$ two such
morphisms, the difference $|h_f-h_{f'}|$ is bounded on~$X(\Qbar)$ (see
Theorem~3.1 of Chapter~6 of~\cite{Cornell-Silverman}).  For this
reason, one usually associates to a line bundle~$\calL$ on a
variety~$X$ a class of height functions~$f_{\calL}$, i.e., an element
in the set of functions $X(\Qbar)\to\RR$ modulo bounded functions;
this map is then a morphism of groups on~$\Pic(X)$:
$f_{\calL_1\otimes\calL_2}\equiv f_{\calL_1}+f_{\calL_2}$. However, in
our situation, we cannot permit ourselves to work just modulo bounded
functions on each variety, as we have infinitely many curves $X_1(l)$
to deal with.

There is a geometric way to associate to a line bundle a specific
height function, not just a class of functions modulo bounded
functions. For this, the contributions from the finite as well as the
infinite places must be provided. Those from the finite places come
from a model of~$X$ over the ring of integers~$O_K$ of~$K$, i.e., an
$O_K$-scheme $X_{O_K}$ whose fibre over~$K$ is~$X$, together with a
line bundle~$\calL$ on $X_{O_K}$ whose restriction to~$X$ is the line
bundle that we had. The $O_K$-scheme $X_{O_K}$ is required to be
proper (e.g., projective). The contributions from the infinite places
are provided by a hermitian metric (or inner product) on~$\calL$, a
notion that we will briefly explain.

A hermitian metric on a locally free $\calO_X$-module of finite rank
$\calE$ consists of a hermitian metric $\ld{\cdot},{\cdot}\rd_x$ on
all $\CC$-vector spaces~$x^*\calE$, where $x$ runs through $X(\CC)$,
the set of $x\colon\Spec(\CC)\to X$. Each $x$ in $X(\CC)$ induces
a morphism $\Spec(\CC)\to\Spec(K)$, i.e., an embedding of $K$
into~$\CC$. Therefore, $X(\CC)$ is the disjoint union of the complex
analytic varieties~$X_\sigma$, indexed by the $\sigma\colon
K\to\CC$. A hermitian metric on $\calE$ consists of hermitian metrics
on all the holomorphic vector bundles $\calE_\sigma$ that $\calE$
induces on the~$X_\sigma$. The metrics to be used are required to be
continuous, i.e., for $U$ open in~$X$ and $s$ and~$t$ in~$\calE(U)$,
the function $x\mapsto \ld s(x),t(x)\rd_x$ on $U(\CC)$ must be
continuous. Actually, the metrics that we will use will live on
non-singular~$X$, and will be required to be smooth (infinitely
differentiable). Another condition that is usually imposed is a
certain compatibility between the metrics at a point~$x$ in~$X(\CC)$
and its complex conjugate~$\ol{x}$. We do not give this condition in
detail, but note that it will be fulfilled by the metrics that we will
use. It is also customary to denote a hermitian metric
$\ld{\cdot},{\cdot}\rd$ by its norm~$\|{\cdot}\|$, given by
$\|s\|^2=\ld s,s\rd$. Indeed, a suitable polarisation identity
expresses the hermitian metric in terms of its norm. A pair
$(\calE,\|{\cdot}\|)$ of a locally free $\calO_{X_{O_K}}$-module with
a hermitian metric~$\|{\cdot}\|$ is called a metrised vector bundle
on~$X_{O_K}$. Metrised vector bundles can be pulled back via morphisms
$f\colon W_{O_K}\to X_{O_K}$ between $O_K$-schemes of the type
considered.

An important example of the above is the case where $X=\Spec(K)$, just
a point, and $X_{O_K}=\Spec(O_K)$. A metrised line bundle
$(\calL,\|{\cdot}\|)$ then corresponds to an invertible $O_K$-module,
$L$, say, with hermitian metrics on the
$L_\sigma:=\CC\otimes_{\sigma,O_K}L$. The \emph{Arakelov degree}
of~$(\calL,\|{\cdot}\|)$ is the real number defined by:
\begin{eqn}\label{eqn_ar_degree}
\deg(\calL,\|{\cdot}\|) = 
\log\#(L/O_Ks) -\sum_{\sigma\colon K\to\CC}\log\|s\|_\sigma,
\end{eqn}
where $s$ is any non-zero element of~$L$ (independence of the choice
of~$s$ follows from the product formula). This definition should be
compared to that of the degree of a line bundle on a smooth projective
curve over a field: there one takes a rational section, and counts
zeros and poles. The first term in~(\ref{eqn_ar_degree}) counts the
zeros of~$s$ at the finite places. Interpreting this term in terms of
valuations, and then norms, at the finite places, then leads to the
second term which ``counts'' the ``zeros'' (or minus the ``poles'',
for that matter) at the infinite places. For a finite extension $K\to
K'$, and $(\calL,\|{\cdot}\|)$ on $\Spec(O_K)$ as above, the pullback
$(\calL',\|{\cdot}\|)$ to $\Spec(O_{K'})$ has degree $\dim_K K'$ times
that on~$\Spec(O_K)$.

We can now give the definition of the height given by a proper
$O_K$-scheme $X$ together with a hermitian line
bundle~$(\calL,\|{\cdot}\|)$. Let $x$ be in~$X(K)$. Then, by the
properness of~$X$ over~$O_K$, $x$ extends uniquely to an $O_K$-valued
point, also denoted~$x$, and one defines:
\begin{eqn}
h_K(x) := \deg x^*(\calL,\|{\cdot}\|).
\end{eqn}
The same method as the one use above can be applied to get an absolute
height $h\colon X(\Qbar)\to\RR$. For $K\to K'$ a finite extension,
each $x$ in $X(K')$ extends uniquely to an $x$ in~$X(O_{K'})$, and one
defines:
\begin{eqn}
h(x) := \frac{\deg x^*(\calL,\|{\cdot}\|)}{\dim_\QQ K'}.
\end{eqn}
It is not hard to verify that this height function~$h$ is in the class
(modulo bounded functions) that is attached to $X_K$ and $\calL_K$
(without metric); see Proposition~7.2 of Chapter~6
of~\cite{Cornell-Silverman}, or Theorem~4.5 of Chapter~V
in~\cite{Edixhoven-Evertse}.  In fact, for $X=\PP^n_{O_K}$, and
$\calL=\calO(1)$ with a suitable metric, the height $h$ just defined
is \emph{equal} to the one defined in~(\ref{eqn_height_Pn}).

\section{Arithmetic Riemann-Roch and intersection theory on
  arithmetic surfaces}
\label{subsec_ARR}

The context in which we are going to apply Arakelov theory is that of
smooth projective curves~$X$ over number
fields~$K$. In~\cite{Arakelov} Arakelov defined an intersection theory
on the \emph{arithmetic surfaces} attached to such curves, with the
aim of proving certain results, known in the case of functions fields,
in the case of number fields. The idea is to take a regular projective
model~$\calX$ over $B:=\Spec(O_K)$ of~$X$, and try to develop an
intersection theory on the surface~$\calX$, analogous to the theory
that one has when $K$ is a function field. If $K$ is a function field
over a finite field~$k$, say, one gets a projective surface $\calX$
over~$k$, fibred over the nonsingular projective curve~$B$ over~$k$
that corresponds to~$K$. On such a projective surface, intersecting
with principal divisors gives zero, hence the intersection pairing
factors through the Picard group of~$\calX$, the group of isomorphism
classes of invertible $\calO_\calX$-modules. In the number field case
one ``compactifies'' $B$ by formally adding the infinite places
of~$K$; the product formula then means that principal divisors have
degree zero. Instead of the Picard group of~$\calX$, one considers the
group of isomorphism classes of certain metrised line bundles
on~$\calX$, as defined above. In~\cite{Faltings1}, see also
Chapters~II, III and~I of~\cite{Szpiro1}, Faltings extended Arakelov's
work by establishing results such as a Grothendieck-Riemann-Roch
theorem in this context. Since then, Arakelov theory has been
generalised by Gillet and Soulé (see~\cite{Soule1} and
\cite{Faltings2}). Below, we will use the theory as given
in~\cite{Faltings1} and Chapter~II of~\cite{Szpiro1}. We start with
some preparations concerning Riemann surfaces. The aim of this
subsection is to give the arithmetic Riemann-Roch theorem as stated
and proved by Faltings.

Let $X$ be a compact Riemann surface of genus $g>0$. The space of
holomorphic differentials $\rH^0(X,\Omega_X^1)$ carries a natural
hermitian inner product: 
\begin{eqn}\label{eqn_nat_inner_pro}
(\omega,\eta) \mapsto \frac{i}{2} \int_X \omega \wedge \ol{\eta}.
\end{eqn}
Let $(\omega_1,\ldots,\omega_g)$ be an orthonormal basis with respect
to this inner product. This leads to a positive $(1,1)$-form $\mu$
on~$X$ given by:
\begin{eqn}\label{eqn_arakelov_1-1-form}
\mu = \frac{i}{2g} \sum_{k=1}^g \omega_k \wedge \ol{\omega_k},
\end{eqn}
independent of the choice of orthonormal basis. Note that $\int_X
\mu=1$.  We refer to \cite{Arakelov} for a proof of the following
proposition. Denote by $\mathcal{C}^\infty$ the sheaf of complex
valued $C^\infty$-functions on~$X$, and by $\mathcal{A}^1$ the sheaf
of complex $C^\infty$ $1$-forms on~$X$. Recall that we have a
tautological differential operator $d\colon \mathcal{C}^\infty \to
\mathcal{A}^1$. It decomposes as $d = \partial + \ol{\partial}$ where,
for any local $C^\infty$ function $f$ and any holomorphic local
coordinate $z$, with real and imaginary parts~$x$ and~$y$, one has
$\partial f = \frac{1}{2}(\frac{\partial f}{\partial x}
-i\frac{\partial f}{\partial y})\cdot dz$ and $\ol{\partial} f =
\frac{1}{2}(\frac{\partial f}{\partial x} +i\frac{\partial f}{\partial
y})\cdot d \ol{z}$.
\begin{prop}\label{prop_Ar-Gr-function}
For each $a$ in~$X$, there exists a unique real-valued $g_{a,\mu}$
in $\mathcal{C}^\infty(X - \{ a\})$ such that the following properties
hold:
\begin{enumerate}
\item we can write $g_{a,\mu} = \log |z - z(a)| + h$ in an open
neighbourhood of~$a$, where $z$ is a local holomorphic coordinate and
where $h$ is a $C^\infty$-function;
\item $\partial \ol{\partial} g_{a,\mu} = \pi i \mu$ on $X -
\{a\}$;
\item $\int_X g_{a,\mu} \mu = 0.$ 
\end{enumerate}
\end{prop}
We refer to $\mu$ and the $g_{a,\mu}$ as the Arakelov $(1,1)$-form and
the Arakelov-Green function, respectively. 
A fundamental property of the functions $g_{a,\mu}$ is that they give
an inverse to the map $\calC^\infty\to\calA^2$, $f\mapsto(-1/\pi
i)\partial\ol{\partial}f$, with $\calA^2$ the sheaf of complex
$C^\infty$ $2$-forms on~$X$, up to constants. For all $f$ in
$\calC^\infty(X)$ we have:
\begin{eqn}\label{eqn_green_inverts}
f(x) = 
\int_{y\in X} g(x,y) \frac{-1}{\pi i} (\partial\ol{\partial}f)y 
+\int_X f\,\mu_X\, .
\end{eqn}
For a proof of this see~\cite[pp.~393--394]{Faltings1},
or~\cite[Lemme~4]{Elkik}.

We note that Stokes' theorem implies $g_{a,\mu}(b)=g_{b,\mu}(a)$ for
all~$a$ and~$b$ in~$X$.  The Arakelov-Green functions determine
certain metrics, called admissible metrics, on all line
bundles~$\calO_X(D)$, where $D$ is a divisor on~$X$, as well as on the
holomorphic cotangent bundle~$\Omega^1_X$. To start, consider line
bundles of the form $\calO_X(a)$ with $a$ a point in~$X$ (the general
case with $D$ follows by taking tensor products). Let $s$ be the
tautological section of~$\calO_X(a)$, i.e. the constant
function~$1$. We define a smooth hermitian metric
$\|{\cdot}\|_{\calO_X(a)}$ on~$\calO_X(a)$ by putting $\log
\|s\|_{\calO_X(a)}(b) = g_{a,\mu}(b)$ for any $b$ in~$X$. By
property~2 of the Arakelov-Green function, the curvature form $(2\pi
i)^{-1}\partial\ol{\partial}\log(\|s\|^2)$ of $\calO_X(a)$ is equal
to~$\mu$.  To continue, it is clear that the functions $g_{a,\mu}$ can
be used to put a hermitian metric on the line bundle $\calO_{X\times
  X}(\Delta_X)$, where $\Delta_X$ is the diagonal on $X\times X$, by
putting $\log\|s\|(a,b)=g_{a,\mu}(b)$ for the tautological section~$s$
of $\calO_{X\times X}(\Delta_X)$. Restricting to the diagonal we have
a canonical adjunction isomorphism $\calO_{X \times
  X}(-\Delta_X)|_{\Delta_X} \isomlto \Omega_X^1$. We define a
hermitian metric $\|{\cdot}\|_{\mathrm{Ar}}$ on $\Omega_X^1$ by
insisting that this adjunction isomorphism be an isometry. It is
proved in~\cite{Arakelov} that this gives a smooth hermitian metric
on~$\Omega_X^1$, and that its curvature form is a multiple of~$\mu$.
From now on we will work with these metrics on~$\calO_X(P)$ and
$\Omega_X^1$ (as well as on tensor product combinations of them) and
refer to them as Arakelov metrics. Explicitly: for $D=\sum_P n_PP$ a
divisor on~$X$, we define $g_{D,\mu}:=\sum_P n_Pg_{P,\mu}$, and equip
$\calO_X(D)$ with the metric $\|{\cdot}\|$ for which $\log\|1\|(Q) =
g_{D,\mu}(Q)$, for all $Q$ away from the support of~$D$.  A metrised
line bundle $\calL$ in general is called admissible if, up to a
constant scaling factor, it is isomorphic to one of the admissible
bundles~$\calO_X(D)$, or, equivalently, if its curvature form
$\mathrm{curv}(\calL)$ is a multiple of~$\mu$. Note that then
necessarily we have $\mathrm{curv}(\calL) = (\deg \calL) \cdot \mu$ by
Stokes' theorem.
 
For any admissible line bundle~$\calL$, Faltings defines a certain
metric on the determinant of cohomology $\lambda(\calL) = \det
\rH^0(X,\calL) \otimes \det \rH^1(X,\calL)^\vee$ of the underlying
line bundle. This metric is the unique metric satisfying a set of
axioms. We recall these axioms (cf.~\cite{Faltings1}, Theorem~1):
(i)~any isometric isomorphism $\calL_1\isomlto\calL_2$ of admissible
line bundles induces an isometric isomorphism $\lambda(\calL_1) \isomlto
\lambda(\calL_2)$; (ii)~if we scale the metric on~$\calL$ by a
factor~$\alpha$, the metric on $\lambda(\calL)$ is scaled by a
factor~$\alpha^{\chi(\calL)}$, where $\chi(\calL)=\deg \calL -g +1$ is
the Euler-Poincaré characteristic of~$\calL$; (iii)~for any divisor
$D$ and any point $P$ on~$X$, the exact sequence $0 \to \calO_X(D-P)
\to \calO_X(D) \to P_*P^*\calO_X(D) \to 0$ induces an isometry
$\lambda(\calO_X(D)) \isomlto \lambda(\calO_X(D-P)) \otimes
P^*\calO_X(D)$; (iv)~for $\calL=\Omega^1_X$, the metric on
$\lambda(\calL) \cong \det \rH^0(X,\Omega^1_X)$ is defined by the
hermitian inner product $(\omega,\eta) \mapsto (i/2)\int_X
\omega \wedge \ol{\eta}$ on~$\rH^0(X,\Omega_X^1)$. In particular, for
an admissible line bundle~$\calL$ of degree~$g-1$, the metric on the
determinant of cohomology $\lambda(\calL)$ is independent of scaling.

It was proved by Faltings that we can relate the metric on the
determinant of cohomology to theta functions on the Jacobian
of~$X$. Let $\HH_g$ be the Siegel upper half space of complex
symmetric $g$-by-$g$-matrices with positive definite imaginary
part. Let $\tau$ in~$\HH_g$ be the period matrix attached to a
symplectic basis of $\rH_1(X,\ZZ)$ and consider the analytic Jacobian
$J_\tau(X) = \CC^g /(\ZZ^g + \tau \ZZ^g)$ attached to~$\tau$. On
$\CC^g$ one has a theta function 
$\vartheta(z;\tau) = \sum_{n\in\ZZ^g} 
\exp(\pi i\,{}^t\hspace{-0.1em}n \tau n 
+ 2\pi i\, {}^t\hspace{-0.1em}n z)$, 
giving rise to a reduced effective divisor~$\Theta_0$ and a line bundle
$\calO(\Theta_0)$ on~$J_\tau(X)$. Now consider on the other hand the
set $\mathrm{Pic}_{g-1}(X)$ of divisor classes of degree $g-1$
on~$X$. It comes with a canonical subset $\Theta$ given by the classes
of effective divisors. A fundamental theorem of Abel-Jacobi-Riemann
says that there is a canonical bijection $\mathrm{Pic}_{g-1}(X)\isomlto
J_\tau(X)$ mapping $\Theta$ onto~$\Theta_0$. As a result, we can equip
$\mathrm{Pic}_{g-1}(X)$ with the structure of a compact complex
manifold, together with a divisor $\Theta$ and a line
bundle~$\calO(\Theta)$.

The function $\vartheta$ is not well-defined on
$\mathrm{Pic}_{g-1}(X)$ or~$J_\tau(X)$.  We can remedy this by putting
\begin{eqn}\label{eqn_thetanorm}
\|\vartheta\|(z;\tau) = 
(\det \Im(\tau))^{1/4} \exp(-\pi\,{}^t\hspace{-0.1em}y
(\Im(\tau))^{-1} y)|\vartheta(z;\tau)|, 
\end{eqn}
with $y = \Im(z)$. One can check that $\|\vartheta\|$ descends to a
function on~$J_\tau(X)$. By our identification $\mathrm{Pic}_{g-1}(X)
\isomlto J_\tau(X)$ we obtain $\|\vartheta\|$ as a function
on~$\mathrm{Pic}_{g-1}(X)$. It can be checked that this function is
independent of the choice of~$\tau$.  Note that $\|\vartheta\|$ gives
a canonical way to put a metric on the line bundle $\calO(\Theta)$
on~$\mathrm{Pic}_{g-1}(X)$. For any line bundle $\calL$ of
degree~$g-1$ there is a canonical isomorphism from $\lambda(\calL)$ to
$\calO(-\Theta)[\calL]$, the fibre of $\calO(-\Theta)$ at the point
$[\calL]$ in $\mathrm{Pic}_{g-1}(X)$ determined by~$\calL$. Faltings
proves that when we give both sides the metrics discussed above, the
norm of this isomorphism is a constant independent of~$\calL$; he
writes it as~$\exp(\delta(X)/8)$. In more explicit terms, this means
that for any line bundle $\calL$ of degree $g{-}1$ on $X$ with
$h^0(\calL)=0$ (and hence $h^1(\calL)=0$) we have:
\begin{eqn}\label{eqn_deltaX}
\lambda(\calL)=\CC, \quad \|1\|_{\lambda(\calL)}^{-1} = 
\exp(\delta(X)/8){\cdot}\|\theta\|([\calL]).
\end{eqn}
The invariant $\delta(X)$ of~$X$ appears in the Noether formula, see
below.

We will now turn to intersections on an arithmetic surface. For us, an
arithmetic surface is a proper, flat morphism $p\colon\calX\to B$ with
$\calX$ a regular scheme, with $B$ the spectrum of the ring of
integers $O_K$ in a number field $K$, and with generic fibre a
geometrically connected and smooth curve~$X/K$. We say that $\calX$ is
of genus~$g$ if the generic fibre is of genus~$g$.  We will always
assume that $p$ is a semi-stable curve, unless explicitly stated
otherwise. After extending the base field if necessary, any
geometrically connected, smooth proper curve $X/K$ of positive genus
with $K$ a number field is the generic fibre of a unique semi-stable
arithmetic surface. 

An Arakelov divisor on~$\calX$ is a finite formal integral linear
combination of integral closed subschemes of codimension~$1$
of~$\calX$ plus a contribution $\sum_\sigma \alpha_\sigma \cdot
F_\sigma$ running over the complex embeddings of~$K$. Here
$\alpha_\sigma $ is a real number, and the symbols~$F_\sigma$
correspond to the compact Riemann surfaces~$X_\sigma$ obtained by base
changing $X/K$ to~$\CC$ via~$\sigma$.  We have an $\RR$-valued
intersection product $(\cdot,\cdot)$ for such divisors, respecting
linear equivalence. When we want to indicate which model $\calX$ is
used for this intersection product, we will use the
notation~$(\cdot,\cdot)_\calX$. The notion of principal divisor is
given as follows: let $f$ be a non-zero rational function in~$K(X)$,
then $(f)=(f)_\mathrm{fin}+(f)_\mathrm{inf}$ with $(f)_\mathrm{fin}$
the usual Weil divisor of~$f$ on~$\calX$, and with $(f)_\mathrm{inf}=
\sum_\sigma v_\sigma(f) \cdot F_\sigma$ with $v_\sigma(f) =
-\int_{X_\sigma} \log |f|_\sigma \mu_\sigma$. For a list of properties
of this intersection product we refer to~\cite{Arakelov},
\cite{Faltings1} or Chapter~II of~\cite{Szpiro1}.

It is proved in~\cite{Arakelov} that the group of linear equivalence
classes of Arakelov divisors is canonically isomorphic to the
group~$\widehat{\mathrm{Pic}}(\calX)$ of isometry classes of
admissible line bundles on~$\calX$. By an admissible line bundle
on~$\calX$ we mean the datum of a line bundle~$\calL$ on~$\calX$,
together with admissible metrics on the restrictions~$\calL_\sigma$
of~$\calL$ to the~$X_\sigma$. In particular we have a canonical
admissible line bundle $\omega_{\calX/B}$ whose underlying line bundle
is the relative dualising sheaf of~$p$. In many situations it is
convenient to treat intersection numbers from the point of view of
admissible line bundles. 

For example, if $P\colon B\to\calX$ is a section of~$p$, and $D$ is an
Arakelov divisor on~$\calX$, the pull-back $P^*\calO_\calX(D)$ is a
metrised line bundle on~$B$, and we have:
\begin{eqn}\label{eqn_def_int_sec}
(D,P)=\deg P^*\calO_\calX(D),
\end{eqn}
where the degree~$\deg$ of a metrised line bundle is as defined
in~(\ref{eqn_ar_degree}). As a second example, we mention that by
definition of the metric on~$\omega_{\calX/B}$, we have for each
section $P\colon B\to\calX$ of~$p$ an adjunction formula:
\begin{eqn}\label{eqn_adjunction_formula}
(P,P+\omega_{\calX/B})=0 \, .
\end{eqn}

For an admissible line bundle~$\calL$ on~$\calX$, we have the notion
of determinant of cohomology on~$B$, in this context denoted
by~$\det\rR p_*\calL$ (see Chapter~II of~\cite{Szpiro1}). By using the
description above for its metrisation over the complex numbers, we
obtain the determinant of cohomology on~$B$ as a metrised line
bundle. One of its most important features is a metrised Riemann-Roch
formula (cf.~\cite{Faltings1}, Theorem~3), also called
\emph{arithmetic Riemann-Roch formula}:
\begin{eqn}
\label{eqn_FRR}
\deg \det \rR p_*\calL = 
\frac{1}{2}(\calL,\calL \otimes \omega_{\calX/B}^{-1} ) +
\deg \det p_* \omega_{\calX/B} 
\end{eqn}
for any admissible line bundle~$\calL$ on~$\calX$.

The term $\deg \det p_* \omega_{\calX/B}$ is also known as the
\emph{Faltings height} of~$X$, the definition of which we will
now recall. We let $J_K$ be the Jacobian variety of~$X$, and $J$
its Néron model over~$B$. Then we have the locally free $O_K$-module
$\Cot_0(J):=0^*\Omega^1_{J/O_K}$ of rank~$g$, and hence the invertible
$O_K$-module of rank one:
\[
\omega_J:=\bigwedge^g 0^*\Cot_0(J).
\]
For each $\sigma\colon K\to\CC$ we have the scalar product on
$\CC\otimes_{O_K}\omega_J$ given by:
\[
\ld\omega|\eta\rd_\sigma = 
(i/2)^g(-1)^{g(g-1)/2}\int_{J_\sigma(\CC)}\omega\wedge\ol{\eta}.
\]
The Faltings height $h_K(X)$ is then defined to be the Arakelov
degree of this metrised line bundle:
\begin{eqn}
\label{def_Faltings_height}
h_K(X)=\deg(\omega_J),
\end{eqn}
and the absolute Faltings height (also called stable Faltings height)
$h_\abs(X)$ of~$X$ is defined as:
\begin{eqn}
\label{def_abs_Faltings_height}
h_\abs(X)=[K:\QQ]^{-1}\deg(\omega_J).
\end{eqn}
We remark that the stable Faltings height of $X$ does not change after
base change to larger number fields; that is why it is called
stable. Therefore, $h_\abs(X)$ can be computed from any model of $X$
over a number field as long as that model has stable reduction over
the ring of integers of that number field.

As $\calX\to B$ is semi-stable, a result of Raynaud
% give reference!!!!!!!!!!!!!!!!!!
gives that the connected component of~$0$ of~$J$ is the Picard scheme
$\Pic^0_{\calX/B}$, whose tangent space at~$0$ is
$\rR^1p_*\calO_\calX$. Therefore, $\Cot_0(J)$ is the same
as~$p_*\omega$, as locally free $O_K$-modules. A simple calculation
(see lemme~3.2.1 in Chapter~I of~\cite{Szpiro1}) shows that, with
these scalar products, $\omega_J$ and $\det p_*\omega$ are the same as
metrised $O_K$-modules. Therefore we have:
\begin{eqn}
\label{eqn_Faltings_height}
h_K(X) = \deg(\omega_J) = \deg \det p_*\omega .
\end{eqn}

One may derive from~(\ref{eqn_FRR}) the following projection formula:
let $E$ be a metrised line bundle on~$B$, and $\calL$ an admissible
line bundle on~$\calX$. Then the formula:
\begin{eqn}
\deg \det \rR p_* (\calL \otimes p^*E) = 
\deg \det \rR p_*\calL + \chi(\calL) \cdot \deg E
\end{eqn}
holds.  Here again $\chi(\calL)$ is the Euler-Poincaré characteristic
of~$\calL$ on the fibres of~$p$.

%JMCdebut

\chapter{Computing complex zeros  of polynomials and series}\label{sec_couveignes_ZEROS}

\author{J.-M. Couveignes}

\bigskip

\bigskip

%JMCdebut

The purpose of this chapter is twofold. We first want to prove 
Theorems~\ref{theorem:quadtree} and \ref{theorem:findingzeros}
 below about the complexity of computing complex
roots  of polynomials and zeros of power series.
The existence of a deterministic polynomial time
 algorithm for these purposes
 plays an important role in this book.
More importantly, we want to explain what it means for us to compute with
real or complex data  in polynomial time.
All the necessary concepts and algorithms already exist and 
are  provided partly by numerical analysis and partly 
by algorithmic complexity theory. However, the computational model of 
numerical analysis is not quite 
a Turing machine, but rather a real computer
with  floating point arithmetic. Such a computer makes
 rounding errors
at almost every step in the computation. In this context, it
is good enough  to  estimate   the 
conditioning of the problem and the stability of the used algorithm.
Statements about conditioning and stability
tend to be  local and qualitative. And this suffices
to identify and overcome most difficulties and design optimal methods.

Our situation however is quite different. We don't really care about
efficiency. Being  polynomial time is enough to us. On the other
hand, we want a rigorous, unconditional and fully  general proof that
the algorithms we use are polynomial time and return  a result
that is correct up to a small error that must  be bounded rigorously
in any case. For this reason,
 we shall not use floating point registers: we don't want to worry  about
the accumulation of rounding errors. We rather decompose the computation
in big blocks. Inside every block we only allow exact computations
(e.g. using integers or rational numbers). We also check that the function
computed by every such block is well conditioned  and we make a precise
statement for that. Finally, we need to control the accumulation of errors
in a chain of big blocks. But this shall not be too difficult because,
since the blocks are big and efficient enough, the general
organization of the algorithm is simple and involves few blocks.

In Section~\ref{sec:comp} we recall basic definitions in 
computational complexity theory.
Section~\ref{section:sqrt} deals with  the problem of computing square
roots. We illustrate  on this simple example what is expected from an
algorithm
in our context. The more general problem of computing complex roots of 
polynomials is treated in Section~\ref{section:roots}. Finally, we study
in Section~\ref{section:zerosseries} the problem of finding zeros of 
a converging power series.

\medskip

{\bf Notation}:  The symbol $\Theta$ in this chapter  stands for
a positive effective absolute constant. So any statement containing
this symbol becomes true if the symbol is replaced in every occurrence by some 
large enough real number.

\section{Polynomial time complexity classes}\label{sec:comp}

In this section we briefly recall  classical definitions from 
computational complexity theory.
Since we only need to
define the polynomial time complexity classes, we shall not go into
the details. We
refer the reader to Papadimitriou's
book \cite{papadimitriou} for a complete treatment of these
matters.

{\it Turing machines} are a theoretical model for computers. They
are finite automata (they have finitely many inner states)
but  they can write or read on  an infinite tape with
a tape head. A Turing machine
 can be defined by a transition table. For a given
inner state and current character read by the head, the transition
table provides the next inner state, which character to write on
the tape in place
of the current one, and how the head should move on the tape
(one step left, one step right, or no move at all).
See   \cite[Chapter  2]{papadimitriou} for a formal definition. 

A {\it decision problem} is a question that must be answered
by yes or no. For example deciding if an integer is prime.
The answer of a 
 {\it functional problem} is a more general
function of the question.
For example factoring an integer
is a functional problem.
If we want to solve a problem with a Turing machine, we write
the {\it input} on the tape, we run the Turing machine, and
we wait until it stops. We then  read the  {\it output}
on the tape. If the machine always stops and returns the correct
answer, we say that it solves the problem in question.
The time {\it complexity}  is the number of steps 
before  the Turing
machine has solved a given problem.
Such  a Turing
machine
is said to be {\it deterministic}  because
its behaviour only depends
on the input.
The {\it size} of the input is the number of bits required to
encode it. This is the space used  on the tape to write this
input. For example, the size of an integer is the number of bits
in its binary expansion.
A problem is said to be {\it deterministic
polynomial time} if there exists
a deterministic Turing machine that solves it in time
polynomial in the size of the input. The {\it class} of 
all functional  problems that can be solved
in 
deterministic polynomial time is denoted $\bf FP$ or $\bf FPTIME$.
The class of deterministic polynomial time decision problems
is denoted $\bf P$ or $\bf PTIME$.

There exist other models for complexity theory. For example one
may define {\it multitape Turing machines}. There also exist 
{\it random access machines}. All these models
lead to equivalent definitions of
the polynomial complexity classes.
%These latter  machines handle registers containing integers of arbitrary size, and they can perform arithmetic operationson these integers. They are 
An {\it algorithm} is   a sequence of elementary operations
and instructions. Any algorithm can be turned into
a  Turing machine, but this is  fastidious and rather useless since
conceptual description of the algorithm suffices to decide if
the number of elementary operations performed
by the algorithm
is polynomial in the size of
the input. If this is the case, we say that the algorithm
is deterministic polynomial time and we know that the 
corresponding problem is in $\bf PTIME$ or $\bf FPTIME$.

For example, if we want
to multiply two positive integers  $N_1$ and  $N_2$, then the size of the input $(N_1,N_2)$ is the number
of digits in $N_1$ and $N_2$ and this is  $\lceil 
\log_{10}(N_1+1)\rceil +\lceil \log_{10}(N_2+1)\rceil$. 
The number of elementary
operations required by the elementary school algorithm for multiplication is 
$\Theta \times \lceil 
\log_{10}(N_1+1)\rceil \times \lceil \log_{10}(N_2+1)\rceil$. 
The constant here depends on 
the (reasonable) definition we have chosen for what an elementary
operation is. We don't care about constants anyway. We say that the elementary school algorithm
is deterministic polynomial time.
There also exists a deterministic polynomial time algorithm
for Euclidean division (e.g. the elementary school one).
The extended Euclidean algorithm computes coefficients in B{\'e}zout's identity in
deterministic polynomial time also.
So addition, subtraction, multiplication and inversion in the ring $\ZZ/ \! N\ZZ$ can be performed in 
time polynomial in $\log N$. The  class  $a \bmod N$ in $\ZZ/\! N\ZZ$ is represented by
its smallest non-negative element. We denote it $a\% N$. This is the remainder in the Euclidean division
of $a$ by $N$.

A very important problem is exponentiation: given $a\bmod N$ with  $0\le a\le N-1$
and an integer $e\ge 1$, compute
$a^e\bmod N$.

Computing $a^e$ then reducing modulo $N$ is not a good idea because $a^e$ might be very large.
Another option would be to set  $a_1=a$ and compute $a_k=\left(a_{k-1}\times a\right)\%  N$ 
for $2\le k\le e$.
This requires $e-1$ multiplications and $e-1$ Euclidean divisions.
And we never deal with integers bigger than $N^2$.
The complexity of this method is thus $\Theta \times e\times (\log N)^2$ using elementary school algorithms.
It is well known however that we can do much better. We write
the  expansion of $e$ in base $2$
$$e=\sum_{0\le k\le K}\epsilon_k2^k$$
and we set $b_0=a$ and
$b_{k}=b_{k-1}^2 \% N$ for $1\le k\le K$. We then notice that
$$a^e\equiv \prod_{0\le k\le K}b_k^{\epsilon_k}\bmod N.$$
So we can compute $(a^e)\% N$ at the expense of $\Theta \times \log e$
multiplications and Euclidean  divisions between integers
$\le N^2$. The total number of elementary operations is thus
$\Theta\times \log e\times (\log N )^2$ with this method. 
So exponentiation in $\ZZ/N\ZZ$ lies in $\bf FPTIME$.
This is
an elementary but decisive result in algorithmic number theory.
The algorithm above is called {\it fast exponentiation} and it
makes sense in any group. We shall use it many times and in many
different
contexts.

A first interesting consequence is that  for $p$ an odd prime and 
$a$ an integer such that $1\le a\le p-1$,
we can compute the Legendre symbol
$$\left(\frac{a}{p}\right)\equiv a^{\frac{p-1}{2}}\bmod p$$
at the expense of $\Theta (\log p)^3$ elementary operations.
So testing  quadratic  residues is achieved in polynomial deterministic
time.
Assume now that we are interested in the following problem

\medskip 

\parbox{10cm}{\it 
Given an odd prime integer  $p$, find an integer 
$a$ such that $1\le a\le p-1$ and $a$
is not a square modulo $p$. \hfill  $({\star})$
}

\medskip 

This looks like a very easy problem because half of the non-zero
residues modulo $p$ are not squares. So we may just pick 
a random integer $a$  between $1$ and $p-1$ and  compute the
Legendre symbol $\left(\frac{a}{p}\right)=a^{\frac{p-1}{2}}$.
If the symbol  is $-1$  we output $a$. Otherwise we output
$\FAIL$. The probability of success is $1/2$ and  failing is not such a big problem
because we can rerun the algorithm: we just pick another random integer $a$.

This is a typical example of a {\it randomized Las Vegas} algorithm.
The behavior of the algorithm depends on the input of course,
but also on the result of some random choices. One has
to flip coins. A nice model for
such an algorithm would be a Turing machine that receives besides
the input, a long enough  (say infinite) one-dimensional array  $\cR$
consisting of $0$'s and 
$1$'s. 
Whenever the machine needs to flip a coin, she looks at the next 
entry in the array $\cR$. So the Turing machine does not need to flip coins: we provide
her with enough random data at the beginning.
We assume that the running time of the algorithm is bounded 
from above in terms of the size of the input only (this upper bound
should  not depend on the random  data $\cR$).
For \underline{each} input, we ask that the probability (on $\cR$) that
the Turing machine provides the correct answer is $\ge 1/2$. The random
data $\cR$ takes values in $\{0,1\}^{\NN}$. The measure on this latter set
is the limit of the uniform measures on $\{0,1\}^k$ when $k$ tends
to infinity.
If the Turing machine fails to return the correct answer, she should return $\FAIL$ instead.

We just proved that finding a non-quadratic residue modulo $p$ can be done 
in Las Vegas probabilistic polynomial time. There is no known algorithm that can be proven to solve
this problem in deterministic polynomial time. The class of
Las Vegas probabilistic polynomial time
 decision problems is denoted $\bf ZPP$.

\medskip 

We now consider another slightly more difficult problem

\medskip

\parbox{10cm}{\it 
Given an odd prime integer  $p$, find a generating set
$(g_i)_{1\le i\le I}$ for the cyclic group
$(\ZZ/p\ZZ)^*$. \hfill  $({\star}{\star} )$
}

\medskip 

We have a simple probabilistic algorithm for this problem. We compute
an integer $I$ such that
$$\log_2(3\log_2(p-1))\le I\le \log_2(3\log_2(p-1))+2 $$ and we pick $I$ random integers
$(a_i)_{1\le i\le I}$ in the interval $\left[ 1, p-1\right]$. The $a_i$
are uniformly distributed and pairwise independent. We set
$g_i=a_i\bmod p$ and we show that
the $(g_i)_{1\le i\le I}$ generate the group 
$(\ZZ/p\ZZ)^*$ with probability $\ge 2/3$.
Indeed, if  they don't, they must all lye in a maximal subgroup of
$(\ZZ/p\ZZ)^*$. The maximal subgroups 
of  $(\ZZ/p\ZZ)^*$ correspond to prime divisors of $p-1$.
Let $q$ be such a prime divisor. The probability that the
$(g_i)_{1\le i\le I}$ all lye in the subgroup of index $q$ is bounded from above
by $$\frac{1}{q^I}\le \frac{1}{2^I}$$
so the probability that the $(g_i)_{1\le i\le I}$ don't generate $(\ZZ/p\ZZ)^*$
is bounded from above by $2^{-I}$ times the number of prime divisors of $q-1$.
Since the latter is $\le \log_2(p-1)$, the probability of failure is 
$$\le \frac{\log_2(p-1)}{2^I}$$ and this
$\le 1/3$ by definition of $I$.

\medskip 

Note that here we have a new kind of probabilistic algorithm: the answer is correct with probability
$\ge 2/3$ but when the algorith fails, he may return a false answer. Such an algorithm (a Turing machine) is called
 {\it Monte Carlo} probabilistic. This is weaker than a Las Vegas algorithm. 
We just proved that problem $({\star}{\star})$ can be solved in Monte Carlo probabilistic polynomial time.
We don't know of any Las Vegas probabilistic polynomial time algorithm for this problem.
 The class of
Las Vegas probabilistic polynomial time
 decision problems is denoted $\bf BPP$.

In general, a Monte Carlo algorithm can be turned into a Las Vegas one provided  the answer can be checked efficiently, because
we then can force the algorithm to admit that he has failed.
Note also that if we set
$$I=f+\lceil \log_2(\log_2(p-1))\rceil $$
where $f$ is a positive integer; then the probability of failure in the algorithm above
is bounded from above by $2^{-f}$. So we can make this probability
 arbitrarily small at almost no cost.

The main purpose of this book is to prove 
statements about the complexity of computing coefficients
of modular forms. For example, Theorem~\ref{thm_computation_tau}
states that on input  a prime integer 
$p$, computing the Ramanujan function $\tau(p)$
can be done in deterministic polynomial time
in $\log p$. An important intermediate result is to prove
that one can compute some Galois representations modulo $l$
in
time polynomial in $l$. We shall
present two methods for computing such representations. 
Both  methods rely on computing  approximations. The first
method  computes  complex approximations and leads to a deterministic 
algorithm. This  is explained in Chapter~\ref{sec_couveignes_TORSION}
using  the main results in this Chapter~\ref{sec_couveignes_ZEROS}.
We also present in Chapter~\ref{sec_couveignes_modp} a probabilistic
method that relies on computations modulo small auxiliary
primes. The main  reason
why the latter methods are probabilistic is that
they require to find  
generating 
sets for  the Picard group of curves over finite fields. 
This is a generalization of problem
$({\star}{\star})$ and solving it in deterministic
polynomial time is out of reach at the moment.

\section{Computing the square  root of a positive real number}\label{section:sqrt}

In this section, we consider the following problem: 

\medskip 

\parbox{10cm}{\it 
Given a positive real number  $a$,  \\   compute  the  positive square
 root
$b=\sqrt a$
 of $a$. \hfill  $({\star}{\star}{\star})$
}

\medskip 

 We need an algorithm that
runs in deterministic polynomial time. This raises a few simple minded but
important questions about what should be called an algorithm  in this context.
In Section~\ref{subsection:turing} we try to formulate
problem $({\star}{\star}{\star})$ in a more precise way.  We  explain what is meant by an
algorithm
in this context, and which properties one would expect
from  such an algorithm. In Section~\ref{subsection:dicho} we present
the classical dichotomy algorithm and check that 
 it has polynomial
time complexity. The notions presented in this section are classical
and elementary and come from computational complexity 
 theory  \cite{papadimitriou} and
numerical  analysis \cite{Higham,Henrici}. The algorithms and methods we present are not 
original either, and they are far from optimal. We stress that our unique
goal here is to prove that a polynomial time algorithm (in 
a sense that can be made rigorous) exists for some
 classical computational
problem
regarding real  numbers.

\subsection{Turing machines and real numbers}\label{subsection:turing}

We  shall
use  deterministic Turing machines; as defined in  \cite[Chapter
  2]{papadimitriou} for example. 
There is something  annoying  with problem $({\star}{\star}{\star})$ however: both the input and the
output are real numbers. Both
 existing computing devices and  Turing machines only
deal with discrete  data. So they can't deal with
real numbers. We may imagine a  Turing machine or a computer
handling registers with real numbers  as in  \cite{BBS}. However, this
would not be of great use to us, because 
we plan to perform computations on real numbers as an intermediate step in the
computation of a discrete quantity:  our 
 basic idea is to compute an 
integer (having some arithmetic significance)
from a good enough real approximation of it. In the end,
we want a rigorous proof   that the discrete information we are interested
in can be computed by a deterministic Turing machine.
We need to prove that
a standard deterministic Turing machine can  efficiently and
safely 
compute with
real and complex numbers, or at least with approximations of them.
One  possible approach to this classical problem is interval arithmetic
as presented in \cite{hayes}. 
We shall follow a slightly different track, which is better
adapted to our situation. Our goal is 
to prove that a certain number of more or less elementary calculations
on complex numbers can be safely and efficiently performed (in a way that will be made more precise
soon) by an ordinary Turing machine. These calculations include root finding
of polynomials and power series,
 computation with divisors on modular curves, direct and
inverse Jacobi problems on these curves. In this section, the problem 
$({\star}{\star}{\star})$
will be used to illustrate a few simple ideas that will be applied more
systematically in the sequel.
The first question to be addressed  concerns the input.

\medskip 

{\it What is the input of problem $({\star}{\star}{\star})$ ?}

\medskip

Well, if a classical Turing machine is supposed to solve 
problem $({\star}{\star}{\star})$ it
cannot
be given the real number $a$ all at a time. That would be too big for
her. Instead of that,  we assume that the Turing machine is given a black box
$\BOX_a$. On
input a positive     integer $m$, the black box $\BOX_a$
returns a decimal fraction
$N_a\times 10^{-m_a}$ such that $\left| a-N_a\times 10^{-m_a}\right|
\le 10^{-m_a}$.
 If the
black box answers immediately, we will call it an {\it oracle} for $a$.
A more realistic situation is that the black box 
 answers in polynomial time. This means that
on input a positive integer $m_a$, the black box outputs the expected
numerator $N_a$  in time $\le A_a{m_a}^{d_a}$ where $A_a$ and 
$d_a$ are positive integers
depending on $a$ but not on $m_a$. We assume that a Turing machine calling to a
black box (or an oracle) must take the time to read and copy the integrality of
the oracle's answer. For example, a Turing machine with an oracle for $\pi$
cannot access the $10^{100}$-th digit without reading  the previous ones. 

In all the situations we shall be facing,
 there will be a Turing machine in the
black box. However, not every real number 
can be associated with such a Turing machine: the set of Turing
machines is countable and the set of real numbers is not. 
This is the theoretical reason for introducing black boxes there.

\medskip

{\it What should be  the output of a Turing machine
solving problem $({\star}{\star}{\star})$~?}

\medskip

Again, we don't expect the Turing machine to provide us
with the real $b=\sqrt a$ all at a time. We would be a bit embarrassed 
with it anyway. We rather expect the Turing machine, on input
a positive integer $m_b$ and a black box for $a$, will return 
a decimal fraction  $N_b\times 10^{-m_b}$ such that 
$|b-N_b\times 10^{-m_b}|\le 10^{-m_b}$. 
The square root Turing machine may 
call the black box  for  $a$ once or several times.

Altogether,
 the input of the square root Turing machine should consist of a black
box $\BOX_a$ 
for $a$ and a positive integer $m_b$ telling her the desired
absolute accuracy of the expected  result. And the output will be a decimal fraction
approximating $b$.

\medskip

{\it How do we define the complexity of a square root
Turing machine ? What does it mean for such a Turing 
machine to be polynomial time ?}

\medskip

Assume that we have a Turing machine SQRT
that computes square roots.
Assume  that the input of the square root machine SQRT consists of a
black box $\BOX_a$ for $a$ and a positive integer $m_b$ (the required
absolute accuracy of the result). 
%We cant define the  size of $a$.
We look for an upper bound for  the number of elementary operations 
performed by SQRT,
as a function  of $\log \max(a,1)$ and $m_b$. 
 Such a bound will
be called a {\it complexity} estimate for  SQRT. 
Notice that a call to the black box 
$\BOX_a$  will be counted as a single operation.

We assume that there exist two positive integers $A_\SQRT$ and
$d_\SQRT$ such that the complexity of SQRT is bounded above
by a polynomial 
$A_\SQRT (m_b+\log\max(a,1))^{d_\SQRT}$. 
Then the number of calls to 
$\BOX_a$ is certainly bounded by this number. And the absolute accuracy
required  from $\BOX_a$ cannot exceed this number either, otherwise
the machine SQRT would not even find the time to read the 
digits provided to her by $\BOX_a$. 
So by combining a black box $\BOX_a$ for the input
and the Turing machine SQRT,  we obtain a black box 
$\BOX_b$ for the output $b$. And if both SQRT  and $\BOX_a$
have polynomial time complexity, so is the resulting black box
$\BOX_b$.

\medskip

{\it Is problem $({\star}{\star}{\star})$ well conditioned ?}

\medskip

We have seen that the Turing machine SQRT cannot always access
the exact value
of the input $a$. Instead of that SQRT is
 provided with a black box that sends to her
approximations of $a$. We want to make sure that a good approximation
of $b= \sqrt a$ can be deduced from a good approximation of $a$.
Since the function $a\mapsto \sqrt a$ is $\frac{1}{2}$-Lipschitz
on the interval $[1,\infty]$, we have 

\begin{equation}\label{eq:condbig}
|\Delta b |= |\sqrt{a+\Delta a}-\sqrt{a}|\le 
\frac{1}{2}|\Delta a|
\end{equation}
\noindent  as soon as $a\ge 2$ and $|\Delta a|\le 1$.
So a small perturbation of the input  results in a small perturbation
of the expected output in that case. One says that the problem
is {\it well conditioned}. 
% multiplyit by a well chose power of $100$ before taking the square root.

Not every computational problem is well conditioned. For example,
computing the rounding function $a\mapsto \lceil a \rfloor$ 
is not well conditioned if one gets close to $\frac{1}{2}$ because
the function is not even continuous there.

We shall not need to formalize a definition of conditioning,  but
we shall check in several occasions that the function we want to
evaluate is $A$-Lipschitz for a reasonable constant $A$. A weaker
condition may suffice in some cases: for example, assume that we want
to compute a function $a\mapsto b$ and assume  that $-\log\min(1,|\Delta b|)$
is lower bounded by $(-\log\min(1,|\Delta a|))^{\frac{1}{e}}$ for
some fixed positive integer $e$. Then
the loss of accuracy is polynomial in some sense: one can obtain $m$
digits of $b$ from $\Theta m^e$ digits of $a$.

For example, if we consider the problem of computing
the square root of a positive real number  $a$, we notice that
the function $a\mapsto \sqrt a$ is not Lipschitz on $[0,+\infty[$ but
we have 

\begin{eqnarray}\label{eq:condsmall}
\nonumber |\Delta b |= |\sqrt{a+\Delta a}-\sqrt{a}|&=&\sqrt a\times \left|
\sqrt{1+\frac{\Delta a}{a}}
-1
\right|\\ &\le &
\sqrt a\times \sqrt{|\Delta a|}
\end{eqnarray}
\noindent whenever $|\Delta a |\le \min(a,a^2)$.

So for small values of $a$, we loose (no more than) 
half the absolute accuracy when
taking the  square root. This is enough for us to say that the problem is well
conditioned. 

\subsection{The dichotomy algorithm}\label{subsection:dicho}

Given a real interval $[M_1,M_2]$ and  a continuous function $$f :
[M_1,M_2]\rightarrow \RR$$ such that $f(M_1)f(M_2)<0$, the dichotomy algorithm finds an
approximation of a real  zero of $f$ in $[M_1,M_2]$.
We  use the dichotomy algorithm
to compute  the positive square root $b=\sqrt a$
of  a positive decimal  number
$a=N_a\times 10^{-m_a}$ where $N_a\ge 1$ and $m_a\ge 0$ are integers.
So we call $f : [0,\infty[ \rightarrow [0,\infty[ $  the map
$x\mapsto x^2-a$.
The algorithm below  only handles integers and decimal fractions.
 Let $h$ be the smallest integer such that
$10^h\ge N_a$.  If $h-m_a$ is even we set $M_2=10^{\frac{h-m_a}{2}}$,
 otherwise
we take $M_2=10^{\frac{h-m_a+1}{2}}$. We set $M_1=M_2/10$. We assume that 
we are given
also a positive integer $m_b$ (the required absolute  accuracy
of the result).

We use two registers $R_1$ and $R_2$ containing decimal fractions.
The initial value
of  $R_1$ is $M_1$ and the initial value of $R_2$ is $M_2$.

The algorithm then goes as follows:

\begin{enumerate}
\item If $f(\frac{R_1+R_2}{2})$ is zero or if
$|R_1-R_2|\le 10^{-m_b}$,   output $\frac{R_1+R_2}{2}$ and
 stop.
\item If $f(R_1)f(\frac{R_1+R_2}{2})>0$,   set $R_1:=\frac{R_1+R_2}{2}$ and
go to step $1$.
\item If $f(R_1)f(\frac{R_1+R_2}{2})<0$,  set $R_2:=\frac{R_1+R_2}{2}$
and go to step $1$.
\end{enumerate}

The algorithm ouputs a decimal fraction
 $\tilde b$ such that
$|b-\tilde b|\le 10^{-m_b}$. 
The loop is not executed  more
than $\Theta (m_b+h-m_a+1) $ times. 
So the denominator of $\tilde b$ is bounded above
by $10^{\Theta (m_a+m_b+h+1)}$ and the same holds for any intermediate result
in the course of the algorithm.
So the {\it complexity}  of the algorithm is polynomial in $m_a$, 
$\log N$ and $m_b$. And so is the size of the output.

\medskip

We now consider a few questions raised by this algorithm.

\medskip
 
{\it What to do  if the input $a$ is a real number rather than a decimal one ?}

\medskip

In that case, we assume that we are given a black box $\BOX_a$ for
$a$ and a positive  integer $m_b$ (the required absolute  accuracy of the result).
Our first task is to obtain from $\BOX_a$ a positive lower bound for 
$a$. So we first ask her for a decimal  approximation of
$a$ within $10^{-10}$. If the answer is zero, we 
ask her for a decimal  approximation of
$a$ within $10^{-20}$. If the answer is zero again, we
ask her for a decimal  approximation of
$a$ within $10^{-40}$. After a few steps, either we obtain a positive
lower bound for $a$ or we prove that $a$ is smaller than $10^{-2m_b}$.
In the later case, we output $\hat b=0$ which is a good enough 
approximation for $b=\sqrt a$.

The above shows that we can assume that
we know  the smallest 
positive integer $s$ such that $a$ is bigger than $10^{-s}$.

Inequalities   (\ref{eq:condsmall}) and (\ref{eq:condbig}) show
that if $l\ge 2s$ and $|\hat a -a|\le 10^{-\Theta l}$ then 
$|\sqrt{\hat a}-\sqrt a|\le 10^{-l}$   for some absolute constant
$\Theta$. So  we ask  $\BOX_a$ for a decimal approximation $\hat a$ of 
$a$ within $10^{-\Theta(m_b+2s+1)}$. 
We set $b=\sqrt a$ and $\hat b=\sqrt{\hat a}$ and we check that 
$|\hat b -b|\le 10^{-m_b-2s -1}$.
We 
 send the decimal number $\hat a$ to the 
dichotomy algorithm above and ask it for an approximation
of $\hat b =\sqrt{\hat a}$ within $10^{-m_b-1}$.

We obtain a decimal number $\widetilde{\hat b}$ such that
$|\widetilde{\hat b}-\hat b|\le 10^{-m_b-1}$.
We output $\widetilde{\hat b}$ and check that

\begin{equation}\label{eq:combine}
|\widetilde{\hat b}-b|\le |\widetilde{\hat b}-\hat b|+|\hat b-b|
\le 2\times 10^{-m_b-1}\le 10^{-m_b}.
\end{equation}

So we have designed an algorithm (a Turing machine) 
that computes the
positive square root of a positive number $a$ in time polynomial
in $\log (a+1)$ and the required accuracy. 

\medskip
 
{\it Is the above algorithm stable ?}

\medskip

People in numerical analysis  say that an algorithm is {\it stable}
 when the 
value output by the algorithm is not too far from the exact value.
As we just proved, the dichotomy algorithm
 can compute $b=\sqrt a$ within $10^{-m_b}$ in time
polynomial in $\log( a+1)$ and $m_b\ge 1$. So it can be said to be stable.
We shall allow ourselves to
use  this terminology; but we prefer to state
and prove clear and accurate complexity estimates like  the one above.

It is important  to make a distinction between
stability and conditioning. A {\it problem} can be said to be well conditioned.
An {\it algorithm} can be said to be stable. In the above proof that the
dichotomy algorithm for computing square roots is stable, we have
used the fact that the problem itself is well conditioned. This is
illustrated
by inequality  (\ref{eq:combine}).

\section{Computing the complex roots of a polynomial}\label{section:roots}

In this section, we consider the following problem: 

\medskip 

\parbox{10cm}{\it 
Given a degree $d\ge 2$ 
unitary polynomial with complex coefficients
 $$P(x) = x^d+\sum_{0\le k\le d-1} a_kx^k,$$   compute  the  complex
roots  of $P(x)$. \hfill  $({\star}{\star}{\star}{\star})$
}

\medskip

The input of problem  $({\star}{\star}{\star}{\star})$ consists of an integer $d\ge 2$
and a black box $\BOX_P$
for the coefficients of $P(x)$.
 On
input a positive     integer $m$ and an index $k$ such 
that $0\le k\le d-1$, the black box $\BOX_P$
returns a  decimal fraction

$$(N_1+N_2 i)\times 10^{-m-1}\text{ such that }
|a_k- (N_1+N_2i )\times
10^{-m-1} |\le  10^{-m}$$
where $i=\sqrt{-1}\in \CC$.

A Turing machine ROOTS solving problem $({\star}{\star}{\star}{\star})$  should be given 
also a positive  integer  $m_Z$ telling her  the required accuracy 
of the result. Let  $Z=[z_1]+[z_2]+\cdots +[z_d]$  be the divisor
of $P$. This is the formal sum of complex roots, counted with
multiplicities. This is an effective divisor of degree $d$.
On input  a positive  integer $m_Z$ and a
black box for the coefficients of $P(x)$, the machine ROOTS
is expected to return an approximation
$\hat Z=[\hat z_1]+[\hat z_2]+\cdots +[\hat z_d]$ of $Z$ within $10^{-m_Z}$.
This means that there should exist a permutation of the indices
$\tau \in \cS_d$ such that $|z_{\tau (k)}-\hat z_k|\le 10^{-m_Z}$
for every $1\le k\le d$.

The rest of this section is devoted to proving the following theorem.

\begin{thm}[Computing roots of polynomials]\label{theorem:quadtree}
There exists a deterministic algorithm that on input 
 a degree $d$  unitary    
polynomial $$P(x)=
x^d+\sum_{0\le k\le d-1} a_kx^k$$ in $\CC[x]$ and a positive integer  $m_Z$,
computes an approximation 
$$\hat Z=[\hat z_1]+[\hat z_2]+\cdots +[\hat z_d]$$
of the divisor 
$$Z=[z_1]+[z_2]+\cdots +[z_d]$$ of $P$,  within  $10^{-m_Z}$.
This means that there exists a permutation of the indices
$\tau \in \cS_d$ such that $|z_{\tau (k)}-\hat z_k|\le 10^{-m_Z}$
for every $1\le k\le d$.
 The running time
is polynomial in $d$, $\log H$
 and the required accuracy $m_Z\ge 1$. Here
$H$ is the smallest integer bigger than the absolute value
 of  all coefficients 
in $P(x)$.
\end{thm}

In Section~\ref{subsection:condroots} we prove that the problem
is well conditioned in some sense. In Section~\ref{subsection:quadtree}
we recall the principles of Weyl's Quadtree algorithm.
We recall in  Section~\ref{subsection:pan}
that there  exists an  exclusion function that
is   sharp enough and easy to compute.
 This finishes the proof of Theorem~\ref{theorem:quadtree}.

\subsection{Conditioning}\label{subsection:condroots}

Our first concern is to check that the problem of finding roots
is well conditioned  in some sense. 
We first need to define {\it clusterings} of roots. Let $\epsilon$ be a positive real
number. An $\epsilon$-clustering for $P(x)$ consists of a positive integer
$K$ and a  pair $(c_k,n_k)$ for every $1\le k\le K$ such that
the following conditions are satisfied:

\begin{itemize}
\item $c_k$ is a complex
number and $n_k$ is a positive integer  for every $1\le k\le K$.
\item if $k_1\not =k_2$ then $|c_{k_1}-c_{k_2}|>2\epsilon$.
\item There
are $n_k$ roots of $P(x)$, counting 
multiplicities,  in the open disk $D(c_k,\epsilon/2)$ of center
$c_k$ and radius $\epsilon/2$.
\item There
are $n_k$ roots of $P(x)$ in the open disk $D(c_k,\epsilon)$, counting
multiplicities.
\item $\sum_{1\le k\le K}n_k=d$.
\end{itemize}

So we want to squeeze  the complex roots of $P(x)$ into small disks that are
distant  enough from each other.
Note that there may not exist an 
$\epsilon$-clustering for every $\epsilon$. Problems may occur
when the distance between two roots of $P$ is close to
$\epsilon$.
However, for every positive $\epsilon$,  there exists an $\epsilon '$
such that 

\begin{equation}\label{eq:cluster}
2^{-d^2}\epsilon \le \epsilon '\le \epsilon
\end{equation}
 and an 
$\epsilon'$-clustering for $P(x)$.  Indeed, we consider the interval
$$S=[-d^2\log 2+\log\epsilon,\log\epsilon]$$ and for every pair $(z,z')$
of distinct roots of $P(x)$ we remove the interval $$[\log |z-z'|-\log 2,\log
  |z-z'|+\log 2]$$ to $S$. The resulting set $T$ is not empty because
there are at most $d(d-1)/2$ pairs of distinct roots. Any $\epsilon '$
such that $\log \epsilon '$  belongs to $T$ is fine.

\medskip 

Now,  call $Z=[z_1]+[z_2]+\cdots +[z_d]$
the divisor of  $P$.
We call  $A$ the smallest integer 
bigger than the absolute values  of the coefficients of
$P(x)$.
Let $\Delta(x)$ be a polynomial of degree $\le d-1$. Let $\delta$
be the maximum of the absolute values of the coefficients of $\Delta(x)$.
We want to
compare
the roots of $P$ and the roots of $P+\Delta$.
 %We call $M$ the {\it naive height} of $P$.
The absolute values of the roots of $P$ are $\le dA$.
We fix an  $\epsilon \le 1$. We know that
 there exists an $\epsilon'$-clustering $(c_k,n_k)_{1\le k\le
  K}$ of the roots of $P$ for some  $2^{-d^2}\epsilon 
\le \epsilon ' \le \epsilon$. 
Let $z$ be a complex number such that
$|z-c_k|=\epsilon '$
for some $1\le k\le K$.  The absolute value of $P(z)$ is lower bounded by
$(\epsilon ')^d2^{-d}$.  On the other hand, the absolute value of $\Delta(z)$  is upper
bounded 
by $d\delta (dA+2\epsilon ')^{d-1}$. So if  

\begin{equation}\label{eq:rouche}
\delta <  \left(\frac{\epsilon '}{2(dA+2)} \right)^d
\end{equation}
\noindent  we deduce from Rouch{\'e}'s theorem that 
$P+\Delta$  has $n_k$ roots inside $D(c_k,\epsilon ')$. As a consequence, the
roots
$\hat z_1$, $\hat z_2$, \ldots, $\hat z_d$ of $P+\Delta$ can be indexed in such a way that
$|z_j-\hat z_{j}|\le 2\epsilon '\le 2\epsilon$ for every 
$1\le j \le d$. Roughly speaking, the meaning of inequality (\ref{eq:rouche}) is that 
when passing from coefficients to roots, the accuracy is divided  by (no 
significantly more
than) the degree $d$ of the polynomial $P$. We thus have proven the following
lemma.

\begin{lem}[Conditioning of the roots]\label{lemma:condroot}
There exists a positive constant $\Theta$ such that the following is true.
Let $d\ge 1$ be an integer and let $P(x)$ be a degree $d$ unitary polynomial with
complex coefficients. Let $A$ be the smallest integer bigger than the absolute values
of the coefficients of $P$.  
Let $\Delta(x)$ be a degree $d-1$
polynomial
with coefficients bounded above by $\delta$ in absolute value.
Let $\epsilon $ be the unique positive real such that
\begin{equation}\label{eq:condroot}
\log \delta=d\log\epsilon -d\log A-\Theta d^3.
\end{equation}

Assume that $\epsilon \le 1$. 
Let $Z=[z_1]+[z_2]+\cdots +[z_d]$ be the divisor of
$P$ and let $\hat Z=[\hat z_1]+[\hat z_2]+\cdots +[\hat z_d]$
be the divisor of  $P(x)+\Delta(x)$.
There exists a permutation of the indices
$\tau \in \cS_d$ such that $|z_{\tau (k)}-\hat z_k|\le  \epsilon$
for every $1\le k \le d$.
\end{lem}

This lemma tells us that if we are looking for an approximation
within $10^{-m_Z}$ of the divisor  $Z=[z_1]+[z_2]+\cdots +[z_d]$ 
of a unitary polynomial $P(x)\in \CC[x]$ given by a blackbox $\BOX_P$, we
may replace $P$ by a good enough approximation of it having e.g.
decimal coefficients.
Indeed, we first compute    the smallest integer $A$
bigger than the absolute values of all  coefficients of
$P(x)$. We set $\epsilon=10^{-m_Z-1}$ and let $\delta$ be the real
number  given
by Equation~(\ref{eq:condroot}). Let $m_P$ be the smallest integer
such that $10^{-m_P}\le \delta/d$. We call  the black box $\BOX_P$ and
obtain for every $0\le k\le d-1$ a 
decimal fraction $\hat a_k$ with denominator $10^{m_P+1}$
such that $|\hat a_k-a_k|\le 10^{-m_P}\le \delta/d$.
Let $\hat P=x^d+\sum_{0\le k\le d-1}\hat a_kx^k$.
We compute the discriminant of $\hat P$. If it is zero, we
replace $\hat P$ by $\hat P +10^{-m_P}$ (we increase the constant
term by $10^{-m_P}$.) If the discriminant is zero again, we
add $10^{-m_P}$ to the constant term again. We go on like that
until the discriminant of $\hat P$ is non-zero. This process
stops after $d$ steps at most (seen as a polynomial in the indeterminate
$a_0$ the discriminant has degree $d-1$, so it cannot cancel $d$
times.) In the end we obtain a unitary  polynomial $\hat P$ with decimal 
coefficients $\hat a_k$ such that $|\hat a_k-a_k|\le \delta$
and the discriminant of $\hat P$ is non zero. 
If $\hat Z=[\hat z_1]+[\hat z_2]+\cdots +[\hat z_d]$ is the divisor
of $\hat P$, there exists a permutation  of the indices
$\tau \in \cS_d$ such that $|z_{\tau (k)}-\hat z_k|\le  \epsilon$
for every $1\le k\le d$.
The coefficients of $\hat P$ are decimal fractions
with denominator  $10^{m_P+1}$
where $m_P$ is an integer such that 

\begin{equation}\label{eq:f}
m_P \le  \Theta d(m_Z+  d^2+\log A).
\end{equation}
 In order to approximate the roots of $P$ within $10^{-m_Z}$
it  suffices to approximate
the roots of $\hat P$ within $10^{-m_Z-1}$.
So in the sequel we shall assume that we are given a unitary polynomial with
coefficients in $\ZZ[ i ,\frac{1}{10}]$ having no multiple root. 

\subsection{Weyl's Quadtree algorithm}\label{subsection:quadtree}

We now describe a simple-minded variant of the celebrated Weyl's Quadtree
algorithm to compute complex roots of polynomials.
Let $$P(x)= x^d+\sum_{0\le k\le d-1} a_kx^k$$
 be a degree $d\ge 1$ unitary polynomial with
coefficients in $\ZZ[i,\frac{1}{10}]$.
So every coefficient $a_k$ is a fraction
$N_k\times 10^{-m_P-1}$ with $N_k=N_{1,k}+N_{2,k}i$ and 
$N_{1,k}$, $N_{2,k}$ are in $\ZZ$. 
%We denote by $N$ the maximum of the  modules $|N_i|$ of the numerators. 
We  assume that $a_0$ is not  zero.
Let $g$ be the smallest positive integer such that the distance
between any two distinct roots of $P$ is at least $10^{-g}$.
  For every complex
number $z$ we denote by $r(z)$ the distance between $z$ and the closest
root of $P$. We denote by $R(z)$ the distance between $z$ and the 
furthest 
root of $P$. Computing $r(z)$ seems difficult unless one already knows
the roots of $P$. However, we assume that we can
compute for every $z$ in $\ZZ[ i ,\frac{1}{10}]$
a decimal $\rho(z)$ such that $\rho(z)\le r(z)\le 1.01\times \rho(z)$.
Such a $\rho$ is called an {\it exclusion function}. We shall give
in Section~\ref{subsection:pan} an example of such an exclusion
function that can be computed efficiently.

We first construct 
 a square $\cQ$ in the complex plane $\CC$ that is large enough
to contain all roots of $P$. We  take for $\cQ$
a   square with center the origin
and side length $$s=2dA$$ where $A$
is the smallest integer bigger than the absolute values of all coefficients
in $P$.

Now we divide $\cQ$ into four squares of side lenght $s/2$: the top
left square $\cQ_1$, the top right square $\cQ_2$, the bottom left
square $\cQ_3$ and the bottom right square $\cQ_4$.
For each $1\le k\le 4$ we evaluate the exclusion function
$\rho$ at the center $c_k$ of $\cQ_k$. If $\rho(c_k)$ is bigger
than  the  half   diagonal  $s/(2\sqrt 2)$ of $\cQ_k$ then we know that there is no
root of $P$ in $\cQ_k$. So we erase this square.

Next we consider all those squares that have not been erased and we
divide each of them into four smaller squares with side length 
$s/4$. We evaluate the exclusion function at the center of every such
square. If the value of the exclusion function if bigger than the
half diagonal $s/(4\sqrt 2)$ of the square in question, we erase it.

We go on like that,  dividing all remaining squares  into four smaller
ones at each step.
The number of remaining squares is never bigger than $4d$. The reason is
that a given root of $P$ cannot compromise more than
$4$ squares at a time (such a situation would occur if the root in question
were  very close to the intersection of four contiguous squares).

After $n$ steps, the side length of the remaining squares is
$s/2^n$. If $s/2^n$ is much smaller than the minimum distance between
two roots of $P$,  then there remains exactly $d$ groups of
contiguous squares, and they each contain a single root of $P$.
So the number of steps is bounded by a constant
times  

\begin{equation}\label{eq:nbetapes}
\log (dA)+m_Z+g 
\end{equation}
\noindent  where $A$
is the smallest integer bigger than the absolute values of all coefficients
in $P$,  and $m_Z$ is the required accuracy of the result, and $10^{-g}$
is a lower bound for the distance between  any two roots of $P$.

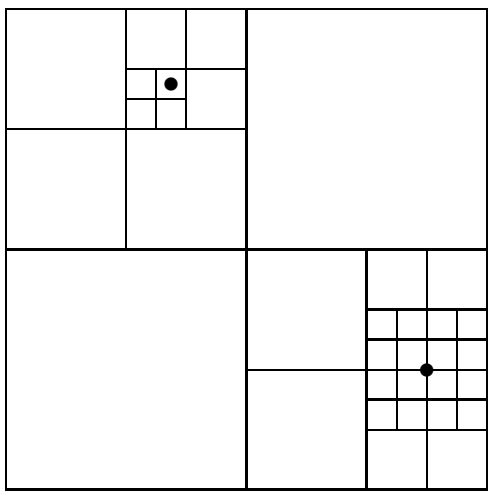
\begin{figure}[h]
\begin{center}
\setlength{\unitlength}{1657sp}%
\begingroup\makeatletter\ifx\SetFigFontNFSS\undefined%
\gdef\SetFigFontNFSS#1#2#3#4#5{%
   \reset@font\fontsize{#1}{#2pt}%
   \fontfamily{#3}\fontseries{#4}\fontshape{#5}%
   \selectfont}%
\fi\endgroup%
\begin{picture}(7330,7330)(2186,-8676)
{\color[rgb]{0,0,0}\thicklines
\put(8551,-6811){\circle*{180}}
}%
{\color[rgb]{0,0,0}\put(4726,-2536){\circle*{180}}
}%
{\color[rgb]{0,0,0}\put(2251,-1411){\line( 1, 0){7200}}
\put(9451,-1411){\line( 0,-1){7200}}
\put(9451,-8611){\line(-1, 0){7200}}
\put(2251,-8611){\line( 0, 1){7200}}
}%
{\color[rgb]{0,0,0}\put(5851,-1411){\line( 0,-1){7200}}
}%
{\color[rgb]{0,0,0}\put(2251,-5011){\line( 1, 0){7200}}
}%
{\color[rgb]{0,0,0}\put(4051,-1411){\line( 0,-1){3600}}
}%
{\color[rgb]{0,0,0}\put(2251,-3211){\line( 1, 0){3600}}
}%
{\color[rgb]{0,0,0}\put(5851,-6811){\line( 1, 0){3600}}
}%
{\color[rgb]{0,0,0}\put(4951,-1411){\line( 0,-1){1800}}
}%
{\color[rgb]{0,0,0}\put(4051,-2311){\line( 1, 0){1800}}
}%
{\color[rgb]{0,0,0}\put(8551,-5011){\line( 0,-1){3600}}
}%
{\color[rgb]{0,0,0}\put(7651,-5011){\line( 0,-1){3600}}
}%
{\color[rgb]{0,0,0}\put(7651,-5911){\line( 1, 0){1800}}
}%
{\color[rgb]{0,0,0}\put(7651,-7711){\line( 1, 0){1800}}
}%
{\color[rgb]{0,0,0}\put(4501,-2311){\line( 0,-1){900}}
}%
{\color[rgb]{0,0,0}\put(4051,-2761){\line( 1, 0){900}}
}%
{\color[rgb]{0,0,0}\put(8101,-5911){\line( 0,-1){1800}}
}%
{\color[rgb]{0,0,0}\put(9001,-5911){\line( 0,-1){1800}}
}%
{\color[rgb]{0,0,0}\put(7651,-6361){\line( 1, 0){1800}}
}%
{\color[rgb]{0,0,0}\put(7651,-7261){\line( 1, 0){1800}}
}%
\end{picture}%
\end{center}
\caption{Weyl's Quadtree method}\label{figure:pan}
\end{figure}

We illustrate this process on Figure~\ref{figure:pan}. The two roots
are represented by two bullets.
%The blue squares are erased first. Then the green ones, the purple
%ones and red ones in this order. 

Note that if the discriminant of $P$ is non-zero, its absolute value is a least
$10^{-(m_P+1)(2d-1)}$ 
so the distance between any two distinct roots is at least
$10^{-(m_P+1)(2d-1)}(2dA)^{-d(d-1)}$. We deduce

$$g\le \Theta d(m_P+d^2 + d\log A).$$

Combining this with the estimates in  Equations~(\ref{eq:nbetapes}) and
(\ref{eq:f}) we deduce that
the number of steps in  Weyl's Quadtree algorithm is

\begin{equation}\label{eq:nbetapes2}
\le \Theta d^2(\log (A)+ d^2 +m_Z)
\end{equation}
\noindent where $d$ is the degree of the polynomial, $A$ the smallest
integer bigger than the coefficients and $m_Z$ the required accuracy
for the roots of $P(x)$.

To finish the proof of Theorem~\ref{theorem:quadtree}  it remains to prove
that there exists an exclusion  function $\rho$ that can be evaluated quickly
enough.
This is done in the next Section~\ref{subsection:pan}.

\subsection{Buckholtz inequalities}\label{subsection:pan}

In this section we recall useful inequalities due to Buckholtz
\cite{buch1, buch2} and
we explain how to deduce a nice exclusion function following Pan \cite{pan}.
Let $P(x)= x^d+\sum_{0\le k\le d-1} a_kx^k$ be a degree $d\ge 2$
unitary polynomial with complex coefficients.  Assume that 
$a_0$ is not  zero. Let $Z=[z_1]+[z_2]+\cdots +[z_d]$ 
be the divisor of $P$.  For every integer $v$ (positive or
negative) we call

$$\nu_v=z_1^v+\dots+z_d^v$$
\noindent  the $v$-th power sum. Buckholtz has shown the following  inequality

\begin{equation}\label{eq:buckholtz}
\frac{1}{5}
\max_{1\le v\le d}|z_v|\le \max_{1\le v \le d}{\left( \frac{|\nu_v|}{d}
  \right)}^{\frac{1}{v}}\le \max_{1\le v\le d}|z_v|.
\end{equation}

One can easily deduce a more general and sharper statement. If $M\ge 1$
is an integer then

\begin{equation}\label{eq:buch2}
{5}^{-\frac{1}{M}}
\max_{1\le v\le d}|z_v|\le \max_{1\le v \le d}{\left( \frac{|\nu_{Mv}|}{d}
  \right)}^{\frac{1}{Mv}}\le \max_{1\le v\le d}|z_v|.
\end{equation}

Recall that for any $z\in \CC$ we call $r(z)$ 
the distance between $z$ and the closest
root of $P$ and $R(z)$ the distance between $z$ and the 
furthest 
root of $P$. 
From Equation~(\ref{eq:buch2}) we deduce an estimate for 
$R(0)$.

\begin{equation}\label{eq:buch3}
\max_{1\le v \le d} 
\left( \frac{|\nu_{Mv}|}{d}
\right)^{\frac{1}{Mv}}
\le R(0) \le 5^{\frac{1}{M}}\max_{1\le v \le d} 
\left( \frac{|\nu_{Mv}|}{d}
\right)^{\frac{1}{Mv}}.
\end{equation}

If we apply inequality (\ref{eq:buch3}) to the reciprocal polynomial of $P(x)$
we obtain

\begin{equation}\label{eq:turan2}
\frac{ 5^{-\frac{1}{M}}}{ \max_{1\le v \le d} 
\left( \frac{| \nu_{-Mv}|}{d}
\right)^{\frac{1}{Mv}}}
\le r(0) \le \frac{1}{\max_{1\le v \le d} 
\left( \frac{| \nu_{-Mv}|}{d}
\right)^{\frac{1}{Mv}}}.
\end{equation}

The power sum $\nu_v$ can be  computed using Newton formulae.  
Assume  that $P(x)$ has coefficients $a_k=
N_k\times 10^{-m_P}$ with $N_k=N_{k,1}+N_{k,2} i $ and 
$N_{k,1}$, $N_{k,2}$ are in $\ZZ$.  Let $A$
be  the smallest integer bigger than the absolute values of the
coefficients $a_k$.
Then computing $\nu_v$  takes
time polynomial in $v$, $m_P$ and $\log A$. 

Extracting the  $v$-th power of a
decimal fraction  $N\times 10^{-m}$ can be done
in time polynomial in $v$, $\log N$, $m$ and the required accuracy,  using a
dichotomy algorithm as in Section~\ref{section:sqrt}.
If $M\ge 200$ then the approximation factor $5^{\frac{1}{N}}$ is smaller than
$1.01$. So we obtain a very sharp estimate for $r(0)$.

For any $z\in \CC$ we can apply inequality (\ref{eq:turan2}) to  
the polynomial
$P(x+z)$ and obtain a good approximation of  $r(z)$.
Assume  that $P(x)$ has coefficients $a_i=
N_i\times 10^{-m_P}$ with $N_i=N_{i,1}+N_{i,2}i$.  
Assume  that $z=(z_{1}+z_{2}i )\times 10^{-m_z}$ and
$z_1$, $z_2$ are in $\ZZ$.
Let $A$
be smallest integer bigger than    the absolute values $|a_i|$.
Then, using Equation~(\ref{eq:turan2}) for $N=200$ 
and the change of variable $x\mapsto x+z$, we can
compute a decimal number  $\rho(z)$ such that

$$\rho(z)\le r(z)\le 1.01\times \rho(z)$$
\noindent and this takes 
time polynomial in $d$, $m_P$,  $\log A$, 
$m_z$, and $$\log\max(|z_1|,|z_2|,1).$$

In the course Weyl's Quadtree algorithm, $m_z$ is   the number
of steps which is bounded in (\ref{eq:nbetapes2}).
Both $|z_1|$ and $|z_2|$ are bounded by $dA10^{m_z}$. An upper
bound for $m_P$ is given in (\ref{eq:f}). 

So the calculation
of any value of the exclusion function that is required in the
course of Weyl's  algorithm takes polynomial time in $d$,
$\log A$ and the required accuracy $m_Z$.
This finishes the proof of Theorem~\ref{theorem:quadtree}.

\section{Computing the zeros of  a power series}\label{section:zerosseries}

Given a power series $f(x)=f_0+\sum_{k\ge 1}f_kx^k\in \CC[[x]]$
 with positive radius of
  convergence $R$, we may want to compute the zeros of $f$ inside the open
  disk $D(0,R)$. However, the number of such zeros may very well be
  infinite. So we should restrict to a smaller disc $D(0,r)$ with
$0<r<R$. 
Then it makes sense to wonder how many  zeros there are in the disk
$D(0,r)$ and try to compute approximations of these  zeros. However, the
  problem
of counting zeros in such a  disk is 
  not a well conditioned problem,  because if $f$ has a zero $z$ of absolute value
  exactly 
$r$, then an infinitesimal perturbation of $f$ may push this zero inside or
outside $D(0,r)$. So we should allow the algorithm to choose a $r'$ that is
very close to $r$ and compute the number of zeros inside $D(0,r')$. Then it
makes sense to ask for approximations of these zeros.
The  input of a  Turing machine  $\ZERO$ computing zeros of power series would consist 
of a black box $\BOX_f$ for the coefficients of $f(x)$. On input an integer $K\ge 0$
and an integer $m\ge 1$, the black box $\BOX_f$  returns for every $0\le  k\le K$
a decimal fraction
$(N_{1,k}+N_{2,k} i)\times 10^{-m-1}$  such that  $|f_k- (N_{1,k}+N_{2,k}i)\times
10^{-m-1} |\le  10^{-m}$. Note that the Turing machine $\ZERO$ is not allowed to ask
about one  coefficient individually. In particular, she needs at least time
$\Theta k$  to receive  any  information about coefficient $f_k$ from the black
box $\BOX_f$. 
Unfortunately this is not enough for her to compute an approximation of the
zeros,  because there may appear a huge coefficient $f_k$ very far away in the
development of $f$. If we want $\ZERO$ to be able to compute zeros without
knowing all coefficients of the power series,  we should at least provide her
with an upper bound for the coefficients in $f(x)$.  We introduce the
following definition.

\begin{defi}[Type of a power series]\label{def:type}
Let $A\ge 1$ be a real and $n\ge 1$ an integer. A power series
$f(x)=f_0+\sum_{k\ge 1}f_kx^k$
is said to be of type $(A,n)$ if for every $k\ge 0$ we have

$$|f_k|\le A(k+1)^n.$$
\end{defi}

The radius of convergence of a power
 series of type $(A,n)$ is at least $1$. Elementary results about series of
 type
$(A,n)$ are collected  in Section~\ref{subsection:elemAn}. The reader is
 advised 
to read quickly the results in this section
 before going further.
 We shall assume that the Turing machine $\ZERO$ is given two integers
$a\ge 1$ and $n\ge 1$ such that the power series $f(x)$ is of type $(\exp(a),n)$.

There is still one difficulty to overcome. For any finite set $F\subset \CC$
we may consider the polynomial $P_F(x)=\prod_{z\in F}(x-z)$. If we divide
$P_F(x)$
by its $L^\infty$ norm we obtain a polynomial $Q_F(x)$ having all its coefficient
bounded by $1$ in absolute value. We may regard
$Q_F$ as a power series of type $(1,1)$. Since $F$ can be  arbitrarily large, we
deduce that we cannot  bound the complexity of finding zeros of a power series
in a given disk $D(0,r)$
just in terms of the type $(A,n)$ of the series. However, we guess that a
power series  of type $(1,1)$
having too many zeros inside  a small compact set contained in  its disk of
convergence, must be small everywhere on this compact set. Also we shall assume
that the Turing machine $\ZERO$ is given a lower bound for the maximum of $f$ 
on $D(0,r)$. More precisely, we assume that $1/2\le r < 1$ and
provide the Turing machine  with
a  positive  integer $\mu$ such that 
$|f(z)|>  \exp(-\mu)$ for {\it at least one}  $z$ in $D(0,1/2)$.
%$-\log \min(1,\max_{z\in D(0,r)}|f(z)|)\le \mu $.

In this section, we shall prove the following theorem.

\begin{thm}[Counting zeros of a power series]\label{theorem:countingzeros}
Let $$f(x)=f_0+\sum_{k\ge 1}f_kx^k\in \CC[[x]]$$ be a power series.
Let $a\ge 1$ be an integer and set $A=\exp(a)$.
Let $n\ge 1$ be an integer such that $f$ is of type $(A,n)$. 
Let $r$ be a real number such that $1/2\le r<1$.  Let $\mu$ be a positive
 integer such 
that there exists at least one $z$ in $D(0,1/2)$ such that 
$|f(z)|> \exp(-\mu)$.

The number of zeros of $f$ in the closed disk $\bD (0,r)$ 
is bounded by a polynomial in $n$, $a$, $\mu$ and $(1-r)^{-1}$.
More precisely, there exists an absolute constant $\Theta$ such that the number
of  zeros of $f$ in the closed disk $\bD (0,r)$  is bounded above
by 

$$\frac{\Theta(n^2+\mu+a)^2}{(1-r)^{13}}.$$
\end{thm}

We shall prove also the existence of an algorithm that computes
approximations of these zeros.

\begin{thm}[Approximating the  zeros of a power series]\label{theorem:findingzeros}
There is   a deterministic algorithm that on input an integer 
$m\ge 1$, an
integer $a\ge 1$, an integer $n\ge 1$, an integer $\mu\ge 1$, 
a rational number $r=1-1/o$ where $o\ge 2$ is
an integer,  and a power series  $$f(x)=f_0+\sum_{k\ge 1}f_kx^k\in
\CC[[x]]$$
of type $(\exp(a),n)$, such that $|f(z)|> \exp(-\mu)$ for at least
one $z$ in $D(0,1/2)$ returns

\begin{itemize}
\item A rational number $r'$ such that $|r-r'|\le 10^{-m}$,
\item The number  $J$ of the zeros  of $f(x)$ in the closed disk $\bD(0,r')$,
or equivalently the degree $J$ of the
 divisor  $Z=[z_1]+[z_2]+\cdots +[z_J]$  of $f$ restricted to 
the closed disk $\bD(0,r')$, 

\item Assuming $J\ge 1$, for every $1\le j\le J$,  a decimal  $$\hat z_j= (M_{1,j}+M_{2,j}i)\times
  10^{-m-1}$$  where  $M_{1,j}$ and $M_{2,j}$ are integers, such that 
$$\hat Z=[\hat z_1]+[\hat z_2]+\cdots +[\hat z_J]$$ approximates 
the divisor  $Z=[z_1]+[z_2]+\cdots +[z_J]$ within $10^{-m}$. More precisely, there
 exists a permutation of the indices
$\tau \in \cS_J$ such that $|z_{\tau (j)}-\hat z_j|\le 10^{-m}$ 
for every $1\le j\le J$.
\end{itemize}

The algorithm runs in time polynomial in $a$, $n$, $\mu$, $o=(1-r)^{-1}$,
and $m$.
\end{thm}

Here is the main idea in the proof of these theorems. For every positive
integer $u$ we write $f(x)=P_u(x)+R_u(x)$ where $P_u(x)=\sum_{0\le k\le
  u-1}f_kx^k$  is the principal part and $R_u(x)=\sum_{k\ge u}f_kx^k$
is the remainder term of order $u$. We expect  that if $u$ is large enough,
the roots of $P_u(x)$ in $D(0,r)$ sharply approximate the zeros  of $f(x)$ 
in $D(0,r)$. It would then suffice to apply the algorithm and theorem in
Section~\ref{section:roots}.

Our first concern will be to prove that the problem of finding zeros
of a power series is well conditioned in some sense: a small perturbation
does not affect too much the zeros. We cannot use Lemma~\ref{lemma:condroot}
about the conditioning of roots
of polynomials,  and we cannot adapt its proof either,  because both the statement
and the  proof involve  the degree of the polynomial. Instead of that, we
shall
first study the conditioning  of every  zero independently.
 In
Section~\ref{subsection:newtonpoly} we introduce the
Newton polynomial of a power series and we show how it can help
locating the zeros in the neighborhood of the origin.
 In Section~\ref{subsection:zero0}
we  deduce  that 
if $f(0)$ is very small, then $f(x)$ has a zero in a small neighborhood 
of $0$. This proves that a zero at the origin is well conditioned.
In Section~\ref{subsection:analyticcont}  we 
state and prove a simple effective
version of the analytic continuation theorem. 
In the next Section~\ref{subsection:anyzero}, we deduce
that  any zero of $f(x)$ is well conditioned in some sense.
A global statement about the conditioning of all
zeros is stated and proved in Section~\ref{subsection:globalzero}.
In the next Section~\ref{subsection:proofthzeros}
we use this result to  finish the 
proof of  Theorem~\ref{theorem:countingzeros}.
We give the proof of Theorem~\ref{theorem:findingzeros} and the corresponding algorithm in
Section~\ref{subsection:proofalgo}.
Section~\ref{subsection:elemAn} is devoted to the
proof of a few  elementary results
about power series of type $(A,n)$.

\subsection{The Newton polygon of a power
  series}\label{subsection:newtonpoly}

Let $f(x)=1+\sum_{k\ge 1}f_kx^k\in \CC[[x]]$  be a power series such that
$f(0)=1$. We assume that the radius $R$ of convergence is positive (it might
be infinite). Let $r$ be a real such that $0<r<R$. We try to mimic
non-archimedean analysis and relate the modulus of zeros of $f(x)$
and the slopes of its Newton polygon.
Let $d$ be the degree of $f(x)$ in the indeterminate $x$. So 
$0\le d\le \infty$ and most of the time $d=\infty$. 
The {\it Newton cloud} of $f(x)$ is the subset of $\RR\times \RR$ 
consisting  of couples $(k,-\log|f_k|)$ for all $k\ge 0$ such that 
$f_k\not =0$.  Let $\Phi$ be the set of all affine functions
$\phi : [0,d]\rightarrow \RR$ whose graph stays below the 
Newton cloud: for every $k\ge 0$ one has
$\phi(k)\le -\log|f_k|$.  The {\it Newton function} is 
a function $\cN : [0,d]\rightarrow \RR$ defined as the sup 
of all $\phi$ in $\Phi$. For every $t\in [0,d]$

$$\cN(t)=\sup_{\phi\in \Phi}\phi(t).$$

This is a convex function on $[0,d]$. It is continuous
and piecewise affine. Indeed,
it is affine on the  interval $[k,k+1]$ for every $0\le k\le d-1$.
The {\it Newton polygon} is the set $\cP=\{(x,y)\in [0,d]\times \RR \, |\, 
\cN(x)\le y\}$.
The {\it vertices} of the Newton polygon are  the points 
$(k,\cN(k))$   for all $k$ where $\cN$ is not differentiable (such a
$k$ must be an integer). This includes the 
point $(0,0)$, and the point $(d,\cN(d))$ also when $d$ is finite.
The Newton polygon has a vertical {\it edge} $](0,+\infty),(0,0)]$.
If $d$ is finite and non-zero, there is another vertical edge
$[(d,\cN(d)),(d,+\infty)[$. The remaining edges connect successive 
vertices of $\cP$.
A {\it slope} of the Newton polygon is a slope of one of its
edges. 
A {\it supporting line} is a line $L$ such that $\cP$ is entirely contained
in one of the two closed half planes determined by $L$, and 
$L$ contains at least one point of $\cP$. For every element $\alpha$
 in $[-\infty,\log R[$  there is a single supporting line $L_\alpha$
with slope
$\alpha$. If the Newton polygon has finitely many edges, then there also exists
a supporting line with slope $\log R$.
If $\alpha$ is a slope of the Newton polygon, then 
$L_\alpha$ contains the corresponding edge of $\cP$.
 
Assume  that $\alpha \in \RR$ is not a  slope of the Newton polygon. The supporting line  
$L_\alpha$ contains a single vertex $v=(k,\cN(k))$ of $\cP$.  Let
$\alpha^-\in [-\infty,+\infty[$ be the slope of the edge just before $v$.
Let 
$\alpha^+\in ]-\infty,+\infty]$ be the slope of the edge just after
$v$. The vector $(-\alpha,1)$ is orthogonal to $L_\alpha$. The Newton
polygon
$\cP$ is contained in the closed half plane $H=\{(x,y)| 
(-\alpha,1).(x,y)\ge (-\alpha,1).(k,\cN(k))\}$. We deduce that for $z$ a complex number such that $\log |z|=\alpha$, the
term $f_kz^k$ dominates the power series in some sense. Indeed, for every
$j\ge 0$ we have 

\begin{eqnarray*}
-\log |f_jz^j|=-\log |f_j|-j\alpha &= &(-\alpha,1).(j,-\log |f_j|)\\
&\ge& (-\alpha,1).(k,-\log |f_k|)
\end{eqnarray*}
\noindent so $f_kz^k$ has a bigger absolute value than any other term in the power
series.

Assume now that $j>k$. The point $(j,-\log |f_j|)$ is above the supporting line $L_{\alpha^+}$
with slope $\alpha^+$. So $|f_j|\le |f_k|\exp(-(j-k)\alpha^+)$.
Since  $\log |z|=\alpha$,  the  term $f_jz^j$ has absolute value  bounded above by
$|f_k||z|^k$ times  $\exp(-(j-k)(\alpha^+-\alpha))$.
The sum $\sum_{j>k} |f_j||z|^j$ is thus bounded above by $|f_k||z|^k$ times  $\frac{a}{1-a}$ where
$a=\exp(-(\alpha^+-\alpha))$.

If $j<k$, then the point $(j,-\log |f_j|)$ is above the supporting line $L_{\alpha^-}$
with slope $\alpha^-$. So $|f_j|\le |f_k|\exp(-(k-j)\alpha^-)$.
Since  $\log |z|=\alpha$,  the  term $f_jz^j$ has absolute value  bounded above by
$|f_k||z|^k$ times  $\exp(-(k-j)(\alpha-\alpha^-))$.
The sum $\sum_{j<k} |f_j||z|^j$ is thus bounded above by $|f_k||z|^k$ times  $\frac{b}{1-b}$ where
$b=\exp(-(\alpha-\alpha^-))$. 

If $\frac{a}{1-a}+\frac{b}{1-b}<1$ then the power series $f(x)$ is dominated
by the term $f_kz^k$ for $|z|=\exp(\alpha)$. %So it does not cancel on the
                                %circle with center $0$ and radius
                                %$\exp(\alpha)$. 
We deduce from Rouch{\'e}'s theorem that $f(x)$ has exactly $k$ zeros inside
$D(0,\exp(\alpha))$
counting multiplicities. This is the case in particular if both $a$ and $b$
are smaller than $\frac{1}{3}$.
We deduce the following lemma.

\begin{lem}[Newton polygon]\label{lemma:Newtonpolygon}
Let  $f(x)=1+\sum_{k\ge 1}f_kx^k$ be a power series. Let $R>0$ be its radius
of convergence. Let  $z \in D(0,R)$ be a z{e}ro of
$f(x)$. There exists a slope  $\sigma$ of the Newton polygon such that  $|\log |z| - \sigma| \le \log 3$.

Let  $\cP_3$ be the set obtained by removing to  $]-\infty, \log R[$ all the
     intervals
 $[\sigma-\log 3, \sigma +\log 3]$
where  $\sigma$ is any slope of the  Newton polygon.

If  a real $\alpha$ belongs to  $\cP_3$, then there is a unique vertex $(k,\cN(k))$
lying on the supporting line with slope $\alpha$. The  power series $f(x)$ has exactly 
 $k$ z{e}ros, counting multiplicities,  inside the open disk
$D(0,\exp(\alpha))$; and none on its boundary.

Let  $\cP_4$ be the set obtained by removing to  $]-\infty, \log R[$ all the
     intervals
 $[\sigma-\log 4, \sigma +\log 4]$
where  $\sigma$ is any slope of the  Newton polygon.

If  a real $\alpha$ belongs to  $\cP_4$, then there is a unique vertex $(k,\cN(k))$
lying on the supporting line with slope $\alpha$. For any complex number $z$ such
that $|z|=\exp(\alpha)$, one has  $|f(z)|\ge \frac{1}{3}\times |f_k|\times {|z|^k}$.

\end{lem}

\subsection{Conditioning of a zero at the origin}\label{subsection:zero0}

In this section, we prove that a power series taking a very small
value at the origin must have a very small zero. More precisely we prove
the following lemma.

\begin{lem}[The smaller zero of a power series]\label{lemma:pluspetitzero}
Consider  a power series 
 $F(x)=F_0+\sum_{k\ge 1}F_kx^k\in \CC[[x]]$
of type $(A,n)$ where $A\ge 1$ is a real and  $n\ge 1$ is an integer.
Assume that $|F_0|<1$ and $$\sqrt{-\log |F_0|}\ge \Theta(n^2+\log A).$$
Then \underline{at least one} of the following two statements holds true

\begin{itemize}
\item $F(x)$ has a zero $\xi$ such that $-\log |\xi| \ge  \sqrt{-\log |F_0|}$,
\item for every $z$ in $D(0,1/2)$ the absolute value of $F(z)$ is upper bounded 
by $\exp(-0.031\sqrt{-\log |F_0|})$.
\end{itemize}
\end{lem}

\bigskip 

Before giving the proof of this lemma, let us sketch the main idea in
this proof: if $F$ is not small everywhere,  there
must be a coefficient in this series that is not too
small. Since $F_0$ is very small, the first finite slope
in the Newton polygon must be very negative. But the lack
of roots in the neighborhood of $0$ and Lemma~(\ref{lemma:Newtonpolygon}) force the Newton polygon to be
smooth at the beginning: the many first slopes must keep
very small. So the series has many huge coefficients, and 
it cannot be of type $(A,n)$. A contradiction.

\bigskip 

So let $F(x)$ be a power series of type $(A,n)$ where $A\ge 1$ is a real and
$n\ge 1$ is an integer. If $F_0=0$ then the first condition  in the conclusion
of the lemma
is met and we are
done. So we assume that $F_0\not= 0$ and we set $f(x)=F(x)/F_0$.
If $F(x)$ is constant, then the second condition in the conclusion
of the lemma
is met and we are
done. So  we assume that $f(x)$ is not constant. The Newton polygon of $f(x)$
has at least one finite slope. We call $\sigma_0=-\infty$,
$\sigma_1$, $\sigma_2$, \ldots \, the successive slopes of the Newton 
polygon. If $\sigma_1\ge 0$ then all coefficients $f_i$ are 
upper bounded by $1$ in absolute value and the second condition in the conclusion
of the lemma is met. So we may assume that $\sigma_1$ is negative.
Recall the definition of $\cP_3$ in Lemma~\ref{lemma:Newtonpolygon}
and consider the intersection $\cP_3\cap ]\sigma_1,0[$. If this
intersection is not empty, we call $\rho$ its infimum and set
$r=\exp( \rho )$. If the intersection is empty, we set $\rho=0$ and
$r=1$.

We assume that the first condition in the conclusion of 
Lemma~\ref{lemma:pluspetitzero} is not met, and we show that in that case the second
condition holds true. 

Since the first condition in the conclusion of 
Lemma~\ref{lemma:pluspetitzero} is not met, we know that  all zeros of $f(x)$ have absolute value
bigger than $\exp(- \sqrt{-\log |F_0|})$. Using 
Lemma~\ref{lemma:Newtonpolygon} we deduce that 

\begin{equation}\label{eq:logr}
-\rho=-\log r\le  \sqrt{-\log |F_0|}.
\end{equation}

The interval $]\sigma_1, \log r[$ does not meet $\cP_3$. So
it is contained in the union of the intervals
$[\sigma_1, \sigma_1+\log 3]$, $[\sigma_2-\log 3,\sigma_2+\log 3]$,
$[\sigma_3-\log 3,\sigma_3+\log 3]$,
$[\sigma_4-\log 3,\sigma_4+\log 3]$, \ldots \,

We deduce that all these intervals should adjust to tile all of
$[\sigma_1, \log r ]$. So $\sigma_2\le \sigma_1+2\log 3$,
$\sigma_3\le \sigma_2+2\log 3$, \ldots, $\sigma_k\le \sigma_{k-1}
+2\log 3$ as long as $\sigma_{k-1}<  \log r-\log 3$.
So $\sigma_k\le \sigma_1+2(k-1)\log 3$ for every 
$k$ from $1$ to $\ell$
where 

\begin{equation}\label{eq:defl}
\ell =\left\lfloor \frac{\log r-\sigma_1}{2\log 3}\right\rfloor.
\end{equation}

We deduce by summation that

\begin{equation}\label{eq:supnl}
\cN(\ell)\le \sigma_1+\sigma_2+\cdots+\sigma_{\ell }\le \ell
\sigma_1+\ell(\ell-1)\log 3.
\end{equation}
\noindent where $\cN$ is the Newton function of the power series $f(x)$.

From (\ref{eq:defl}) we deduce that
$\ell \le \frac{\log r-\sigma_1}{2\log 3}$ and  $\sigma_1\le \log r
-2\ell\log 3$. Reporting in (\ref{eq:supnl}) we obtain

\begin{equation}\label{eq:supnl2}
\cN(\ell)\le \ell \log r-\ell^2\log 3.
\end{equation}

On the other hand, the power series $F(x)$ is of type $(A,n)$, so $f(x)$
is of type $(A/|F_0|,n)$. So the Newton function $k\mapsto \cN(k)$
is lower bounded by the convex function $k\mapsto -\log A+\log |F_0|-n\log(k+1)$. We deduce 

\begin{equation}\label{eq:infnl}
\cN(\ell)\ge -\log A+\log |F_0|-n\log(\ell+1)\ge -\log A+\log |F_0|-n\ell.
\end{equation}

From (\ref{eq:supnl2}) and (\ref{eq:infnl}) we deduce 

\begin{equation*}
\ell^2\log 3-(n+\log r)\ell+\log |F_0|-\log A\le 0.
\end{equation*}

In case  $\ell \ge n+\log r$ we deduce from the inequality above that $\ell^2(\log
3-1)\le \log A-\log |F_0|$. So

\begin{eqnarray}\label{eq:supl}
\nonumber \ell&\le &\max(n+\log r,\sqrt{\frac{\log A-\log |F_0|}{\log 3 -1}})\\
&\le& \max(n,\sqrt{\frac{\log A-\log |F_0|}{\log 3 -1}}).
\end{eqnarray}

We have made the assumption  that

\begin{equation}\label{eq:hyplogF}
\sqrt{-\log |F_0|}\ge \Theta(\log A+n^2).
\end{equation}

From (\ref{eq:supl}) and (\ref{eq:hyplogF}) we deduce

\begin{equation}\label{eq:supl2}
\ell\le 3.2\sqrt{-\log |F_0|}.
\end{equation}

From (\ref{eq:defl}) we know that $\ell\ge  \frac{\log r -\sigma_1}{2\log 3}
-1$ so

\begin{equation}\label{eq:infs1}
\sigma_1\ge -2\ell \log 3-2\log 3 +\log r.
\end{equation}

From (\ref{eq:supl2}) and (\ref{eq:infs1}) and (\ref{eq:logr})  we deduce

\begin{equation}\label{eq:infs12}
\sigma_1\ge -9\sqrt{-\log |F_0|}.
\end{equation}

So for every $k\ge 0$ we have

\begin{equation}\label{eq:supfk}
|F_k|\le |F_0|\exp(-k\sigma_1)\le |F_0|\exp(9k\sqrt{-\log |F_0|}).
\end{equation}

This estimate will show very useful for small values of $k$. For bigger
values of $k$, we use the fact that the power series $F(x)$ is of type
$(A,n)$.  

Now we take some $z\in D(0,1/2)$ and we try to bound $|F(z)|$.
We set 

\begin{equation}\label{eq:defu}
u=\left\lfloor \frac{\sqrt{-\log |F_0|}}{10}  \right\rfloor
\end{equation}
and we check that 
$u\ge \frac{4n^2}{(\log 2)^2}$ 
because of (\ref{eq:hyplogF}).

We write $F(x)=P_u(x)+R_u(x)$ where $P_u(x)=\sum_{0\le k\le
  u-1}F_kx^k$  is the principal part and $R_u(x)=\sum_{k\ge u}F_kx^k$
is the remainder term of order $u$. We now bound $|P_u(z)|$  and 
$|R_u(z)|$ separately.

On the one hand, using Equation~(\ref{eq:supfk}) 

\begin{eqnarray}\label{eq:supPuz}
\nonumber |P_u(z)|&\le &\sum_{0\le k\le u-1 } 2^{-k}
\exp(\log |F_0| + 9(u-1)\sqrt{-\log |F_0|})\\ &\le& \exp(\frac{\log |F_0|}{10}).
\end{eqnarray}

On the other hand, using Lemma~\ref{lemma:majreste}

\begin{equation}\label{eq:supRuz}
|R_u(z)|\le {n!A}{2^{-\frac{u}{2}+n+1}}\le 2^{-\frac{\sqrt{-\log |F_0|}}{21}}
\end{equation}
\noindent because of (\ref{eq:hyplogF}).

Altogether we have

$$|F(z)|\le |P_u(z)|+|R_u(z)|\le \exp(\frac{\log
  |F_0|}{10})+2^{-\frac{\sqrt{-\log |F_0|}}{21}}\le 2^{-\frac{\sqrt{-\log
      |F_0|}}{22}}$$
\noindent using (\ref{eq:supPuz}) and  (\ref{eq:supRuz}),  then
(\ref{eq:hyplogF}).
This proves that the second condition in the conclusion of 
Lemma~\ref{lemma:pluspetitzero} is met. So the proof of this lemma is finished.
\hfill $\Box$

\subsection{Analytic continuation of a power series of type
  $(A,n)$}\label{subsection:analyticcont}
In this section  we consider a power series $f(x)=\sum_{k\ge 0}f_kx^k$ 
of type $(A,n)$ with $A\ge 1$ a real and $n\ge 1$ an integer.
Let $c$ be a complex number with absolute value smaller than $1$ and let 
$r$ be a real such that $0<r<1-|c|$. The disk $D(c,r)$ is contained 
in $D(0,1)$. We want to prove that if $|f(z)|$ is  very small  
for every $z\in D(c,r)$ then $|f(z)|$ is small for every 
$z\in D(0,1/2)$. 
We shall need the following definition.

\begin{defi}[Balanced disk]
Let $D(c,r)$ be a disk contained in the unit disk $D(0,1)$. We say that 
$D(c,r)$ is {\it balanced} if $r=\frac{1-|c|}{2}$. The meaning of this
condition is that the distance between $D(c,r)$ and the unit circle
is equal to the radius of $D(c,r)$.  If this condition is met,
the circle $D(c,r)$ is denoted $D_c$.

Assume further that $|c|>\frac{1}{5}$ and set $c'=(|c|-\frac{r}{2})\times
\frac{c}{|c|}$. Let $D_{c'}$ be the balanced disk with center $c'$. The radius
of $D_{c'}$ is $r'=\frac{5}{4}r$ and $1-|c'|=\frac{5}{4}(1-|c|)$. We say
that $D_{c'}$ is the {\it son} of $D_c$.  

If $|c|\le \frac{1}{5}$, then the  son  of $D_c$ is defined to be
$D(0,1/2)$.
\end{defi}

\begin{figure}[ht]
\begin{center}
\[
\begin{tikzpicture}[thick]
\draw (-0.5,0) -- (9,0);
\draw (0,-0.1) -- (0,0.1)node[anchor=south]{$0$};
\draw (8,-0.1) -- (8,0.1)node[anchor=south]{$1$};
\draw (4.8,-0.1) -- (4.8,0.1);
\draw (4.8,0) circle (1.6 cm) node[anchor=south west]{$c$};
\draw (4,-0.1) -- (4,0.1);
\draw (4,0) circle (2 cm) node[anchor=south]{$c'$};
\draw (4,0) circle (0.8 cm);
\end{tikzpicture}
\]
\end{center}
\caption{Son of a balanced disk}\label{figure:cercles}
\end{figure}
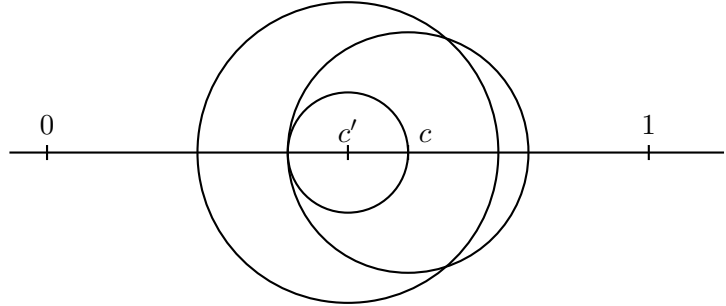

The following lemma states that a power series of type $(A,n)$ that
is very small  on a balanced disk, must be small also on the son of
this disk.

\begin{lem}[Analytic continuation]\label{lemma:analyticcont}
Let $f(x)=\sum_{k\ge 0}f_kx^k$ be a power series of type $(A,n)$ where $A\ge 1$ is
a real and $n\ge 1$ is an integer. Let $c$ be a complex number with
absolute value smaller than $1$ and let $D_c$ be the balanced disk with center
$c$. Call $r=\frac{1-|c|}{2}$  the radius of $D_c$.
Let $D_{c'}$ be the son of $D_c$. Let 
$\epsilon$ be a real  number in $]0,1[$ such that
$-\log \epsilon \ge \Theta(\log A+n(n+|\log r|))$.
 Assume that $|f(z)|\le \epsilon$ for every $z\in D_c$.
Then $|f(z)|\le \epsilon ^{\frac{1}{3}}$ for every $z\in D_{c'}$.
\end{lem}

In order to prove this lemma we observe that the disk $D(c',r/2)$ 
is contained in $D_c$. So the absolute value of $f$ is bounded by $\epsilon$
in $D(c',r/2)$. Cauchy's integral formula then gives an upper bound 
for the successive derivatives of $f$ at $c'$.

\begin{equation}\label{eq:supf(i)}
|f^{(k)}(c')|=\left|\frac{k!}{(2i\pi )}\int_{|\zeta|=r/2}
\frac{f(c' + \zeta)}{\zeta^{k+1}}d\zeta \right|\le \epsilon\frac{2^{k}k!}{r^{k}}.
\end{equation}

In order to bound $|f(z)|$ for $z\in D_{c'}$ we 
call $T(x)=\sum_{k\ge 0}  \frac{f^{(k)}(c')}{k !}x^k$
the Taylor expansion of $f$ at $c'$. We
choose an integer $u\ge 0$
and write $T(x)=P_{c',u}(x)+R_{c',u}(x)$ where $P_{c',u}(x)=\sum_{0\le k\le
  u-1}   \frac{f^{(k)}(c')}{k !} x^k$  
is the principal part and $R_{c',u}(x)=\sum_{k\ge u} \frac{f^{(k)}(c')}{k !} x^k$
is the remainder term of order $u$. 
Any $z$ in the son $D_{c'}$ can be written $z=c'+t$ whith 
$|t|<r'=\frac{5r}{4}$. We now bound $|P_{c',u}(t)|$  and 
$|R_{c',u}(t)|$ separately.

On the one hand, using Equation~(\ref{eq:supf(i)}) we obtain

\begin{eqnarray}\label{eq:majorprincipal}
|P_{c',u}(t)| &\le & \sum_{0\le k \le u-1} \left| \frac{f^{(k)}(c')}{k !}t^k\right|
\nonumber\\ &\le &   \sum_{0\le k \le u-1} \left(\frac{5r}{4}\right)^{k}\frac{2^{k}\epsilon}{r^{k}}\le \epsilon \left(\frac{5}{2}\right)^{u}.
\end{eqnarray}

On the other hand, we set $t=(1-|c'|)y$ and $z=c'+t=c'+(1-|c'|)y$.
The remainder $R_{c',u}((1-|c'|)y)$ is nothing but the remainder
of order $u$ of the refocused\footnote{The refocused series is defined in
  section \ref{subsection:elemAn}}  power series $y\mapsto f_{c'}(y)=
f(c'+(1-|c'|)y)$.  Since $z=c'+y(1-|c'|)$ belongs to the balanced
disk $D_{c'}$ we know that $y$ belongs to $D(0,1/2)$, the balanced
disk with center $0$.  So $|y|<1/2$. We now apply the refocusing
Lemma~\ref{lemma:refocusing} together  with Lemma~\ref{lemma:majreste}.

The refocusing Lemma~\ref{lemma:refocusing} tells  us that
$f_{c'}$ has type $(A',n+1)$ where
$A'\le n!A\left( \frac{2e}{1-|c'|}\right)^{n+2}\le n!A\left(
\frac{e}{r}\right)^{n+2}$. Lemma \ref{lemma:majreste} applied
to the refocused series $f_{c'}$ then says that

\begin{equation}\label{eq:majrestec'}
|R_{c',u}(t)|\le A''(u+1)^{n+1}2^{-u}
\end{equation}
\noindent where 

\begin{equation}\label{eq:majA''}
A''\le (n+1)!n! A\left(\frac{2e}{r}\right)^{n+2}\le A\exp(K n(n+|\log r|))
\end{equation}
\noindent for some positive constant $K$.

We set $u=\left\lceil  \frac{|\log \epsilon |}{2\log \frac{5}{2}}  \right\rceil$.
From (\ref{eq:majorprincipal}) we deduce 

\begin{eqnarray}\label{eq:majPu}
\nonumber \log |P_{c',u}(t)|&\le &
\log \epsilon +u\log\frac{5}{2}\le 
\log \epsilon +\frac{|\log\epsilon|}{2}+ \log\frac{5}{2}\\ &\le &
0.49 \log \epsilon 
\end{eqnarray}
\noindent  using  the hypothesis 

\begin{equation}\label{eq:hypeps}
-\log \epsilon \ge \Theta(\log A+n(n+|\log r|)).
\end{equation}

As far as the remainder is concerned we obtain from (\ref{eq:majrestec'})
and (\ref{eq:majA''})  that $\log |R_{c',u}(t)|$ is

\begin{equation}\label{eq:sumR}
\le \log A+ K n(n+|\log r|)+(n+1)\log (u+1)-u\log 2.
\end{equation}

We show that the negative term $-u\log 2$ dominates this sum.
First of all, from the definition of $u$ and
the hypothesis (\ref{eq:hypeps}) 
we deduce that $$\log A+Kn(n+|\log r |)\le 0.01\times u\log 2.$$
The same hypothesis (\ref{eq:hypeps})  implies
that $$u\ge \left(\frac{n+1}{0.01\times \log 2}\right)^2$$
so $0.01\times u\log 2\ge (n+1)\sqrt u\ge
(n+1)\log(u+1)$. Inequality (\ref{eq:sumR}) then implies

\begin{equation}\label{eq:supRfin}
\log |R_{c',u}(t)|\le -u\times 0.98\log 2\le \frac{0.98\log 2}{2\log\frac{5}{2}}\log
\epsilon\le 0.37 \log \epsilon .
\end{equation}

From (\ref{eq:majPu}) and (\ref{eq:supRfin}) we deduce that if  $z\in D_{c'}$
and $z=c'+t$ then

\begin{eqnarray}\label{majfin}
\nonumber \log |f(z)|&=&\log |P_{c',u}(t)+R_{c',u}(t)|\\
\nonumber  &\le &\log
\left( 2\max(|P_{c',u}(t)|,|R_{c',u}(t)|) \right)\\ &\le &\log 2+0.37\log \epsilon
\le \frac{\log \epsilon}{3}.
\end{eqnarray}

This finishes the proof of Lemma~\ref{lemma:analyticcont}. 
We notice that the exponent $\frac{1}{3}$ in the conclusion
of this lemma could be replaced by any real smaller than $\frac{\log 2}{\log
  \frac{5}{2}}$.\hfill $\Box$
\subsection{Conditioning of any zero}\label{subsection:anyzero}

In  this section we combine
Lemma~\ref{lemma:pluspetitzero} and Lemma~\ref{lemma:analyticcont} to
prove that any zero of a power series is well conditioned 
in some sense.
In other words, we prove that
 a power series taking a very
small value  at some point $c$
inside  its disk of convergence,  must have a zero that is very close to $c$.
More precisely we prove the following lemma.

\begin{lem}[Conditioning of a  zero]\label{lemma:anyzero}
Let $$f(x)=f_0+\sum_{k\ge 1}f_kx^k\in \CC[[x]]$$ be a power series 
of type $(A,n)$ where $A\ge 1$ is a real and  $n\ge 1$ is an integer.
Let $c\in D(0,1)$ be a complex number with absolute value smaller than $1$.
Assume that $|f(c)|<1$ and  

\begin{equation}\label{eq:condanyzero}
(1-|c|)^{6}\times \sqrt{-\log |f(c)|}\ge \Theta(n^2+\log A).
\end{equation}

Then \underline{at least one} of the following two statements holds true

\begin{itemize}
\item $f(x)$ has a zero $\xi$ such that $-\log |\xi-c| \ge  \sqrt{-\log |f(c)|}$,
\item for every $z$ in $D(0,1/2)$ the absolute value of $f(z)$ is upper bounded 
by $$\exp \left( -\frac{(1-|c|)^{5}\times \sqrt{-\log |f(c)|}}{100}\right).$$
\end{itemize}
\end{lem}

Indeed, let $f_c(y)=f(c+y(1-|c|))$ be the refocused series of $f$ at $c$.
According to Lemma~\ref{lemma:refocusing}, it has type $(A',n+1)$ where 

\begin{eqnarray}\label{eq:majA'2}
\log(A')&\le &\log\left(n!A\left( \frac{2e}{1-|c|}\right)^{n+2}\right)\le 
K(n^2+\log A+n|\log(1-|c|)|)\nonumber \\
&\le &K(n^2+\log A+\frac{n}{1-|c|})
\end{eqnarray}
\noindent where $K$ is a  positive constant.

We now apply Lemma~\ref{lemma:pluspetitzero} to $f_c$.
From inequality  (\ref{eq:majA'2}) and the hypothesis
(\ref{eq:condanyzero}) we deduce  that $f_c$ satisfies  the
hypothesis  in Lemma~\ref{lemma:pluspetitzero}. 
We deduce  that one at least 
of the two following conditions is met:

\medskip 

$\bullet$ Either $f_c$ has a zero $y_0$ such that $\log |y_0|\le -\sqrt{-\log |f(c)|} $,
in which case we call $\xi=c+y_0(1-|c|)$ the corresponding zero of $f$. And 
we are done, because $|\xi-c|\le |y_0|$.

\medskip 

$\bullet$ Otherwise, $|f_c(y)|$ is bounded by $\exp(-0.031\times \sqrt{-\log |f(c)|}\, )$
for any $y$ in $D(0,1/2)$. Equivalently, $f(z)$ is bounded by
$$\exp(-0.031\times \sqrt{-\log |f(c)|}\, )$$
for any $z$ in the balanced disk $D_c$ with center $c$ and radius $(1-|c|)/2$.
In that case, we apply Lemma~\ref{lemma:analyticcont} several times.
Indeed, we set 

$$w=\left\lceil \frac{-\log (1-|c|)}{\log \frac{5}{4}}\right\rceil $$
\noindent  and apply
 $w$ times  Lemma~\ref{lemma:analyticcont}. 

We conclude that
for every $z\in D(0,1/2)$, the absolute value of $f(z)$  is bounded by
$$\exp(-{0.031\times 3^{-w}\times \sqrt{-\log |f(c)|}}\, ).$$ Since

$$3^w\le 3(1-|c|)^{-\frac{\log 3}{\log\frac{5}{4}}}\le 3(1-|c|)^{-5}$$
\noindent 
we find that $f$ is bounded in absolute value
by $$\exp(-{0.01\times (1-|c|)^5 \times \sqrt{-\log |f(c)|}}\, )$$ on $D(0,1/2)$. 

Of course, we must check that the hypothesis of Lemma~\ref{lemma:analyticcont}
are satisfied every time we apply it. This is the case if 

$$(1-|c|)^5 \times \sqrt{-\log |f(c)|}\ge \Theta(\log A+n(n+|\log (1-|c|)|).$$

We conclude using %inequality (\ref{eq:majA'2}),
 the hypothesis
(\ref{eq:condanyzero})  and the fact
that $$|\log (1-|c|)|\le (1-|c|)^{-1}.$$\hfill $\Box$

\subsection{Global conditioning of zeros }\label{subsection:globalzero}

In this section, we apply Lemma~\ref{lemma:anyzero} to several
zeros at a time and evaluate the conditioning of the divisor of zeros
of a power series. The following lemma states that the zeros of a power
series can be approximated by the zeros of its principal part.

\begin{lem}[Global conditioning of  the  zeros]\label{lemma:allzeros}
Let $$f(x)=f_0+\sum_{k\ge 1}f_kx^k\in \CC[[x]]$$ be a power series 
of type $(A,n)$ where $A\ge 1$ is a real and  $n\ge 1$ is an integer.
Let $r$ and $\rho$ be two real numbers such that
$1/2\le r<1$ and $0<\rho<1$.
Assume

\begin{equation}\label{eq:condallzeros}
(1-r)^{6}\times |\log \rho |\ge \Theta(n^2+\log A).
\end{equation}

Let $u$ be an integer such that
 $$u\ge \frac{16(\log \rho)^2}{|\log r|}$$ and
$u\le \frac{16(\log \rho)^2}{1-r}$.
Let $v$ be an integer such that the two following conditions
hold true

\begin{itemize}
\item $1\le v\le u$,
\item For any integer $w$ such that $v\le w$ and $w \le  u-1$, the coefficient $f_w$
is bounded in absolute value by $\exp(-6(\log \rho)^2)$.
\end{itemize}

Note in particular that $v=u$ is fine.

Then for every $z$ in $D(0,r)$, the remainder $R_v(z)$ of order $v$ is bounded
in  absolute value by $\exp(-4(\log \rho)^2)$

\begin{equation}\label{eq:borneRv}
|R_v(z)|\le \exp(-4(\log \rho)^2),
\end{equation}
\noindent and  \underline{at least one} of the following two statements holds true

\begin{itemize}
\item for every $z$ in $D(0,1/2)$ the absolute value of $f(z)$ is upper bounded 
by $$\exp \left( 0.01 \times (1-r)^{5}\times \log \rho \right),$$
\item there exists a positive real $R$ such that $r-\rho\le R\le r$
and such that, inside the disk $D(0,R)$, the zeros of $f(x)$
are approximated within $\rho$ 
by the zeros of its  principal part $P_v(x)
=\sum_{0\le k\le v-1}f_kx^k$. In particular, the
number $J$ of these zeros (counting multiplicities)
is at most  $v-1$. More precisely,  let 
$Z=[z_1]+[z_2]+\cdots +[z_J]$
be  the divisor of $f(x)$ restricted to $D(0,R)$, and let
 $\hat Z=[\hat z_1]+[\hat z_2]+\cdots +[\hat z_J]$ be 
 the divisor
of $P_v(x)$ restricted to $D(0,R)$. Then, either $J=0$ or
there exists a permutation  of the indices
$\tau \in \cS_J$ such that $|z_{\tau (j)}-\hat z_j|\le  \rho$
for every $1\le j\le J$.
\end{itemize}

\medskip 

We can be a bit more explicit about the real $R$ above. Any $R$ such that
$R\in [r-\rho,r]$ and $|R-|z||>2\rho^2$ for every root $z$ of $P_v(x)$ is fine.
\end{lem}

\bigskip
The meaning of  this lemma is that the zeros of $f(x)$ are
well approximated by the zeros of a well chosen truncation
of $f(x)$. The lemma tells us where the series should be truncated,
depending on the required absolute  accuracy: if we want $m$ digits of accuracy,
we need to consider the $m^2$ first terms in the series.
One may wonder why we have introduced the integer 
$v$. The point is that if for the given integer $u$, the coefficient
$f_u$ is very small, then we can ignore it and truncate a bit
earlier. This freedom will be used to ensure that the leading coefficient
in $P_v(z)$ is not too small. Then, we can normalize $P_v(z)$ and apply
 Theorem~\ref{theorem:quadtree}
of Section~\ref{section:roots}. It is important to normalize because
Theorem~\ref{theorem:quadtree} only applies to unitary polynomials.

\bigskip

Now let us prove Lemma~\ref{lemma:allzeros}.

Let $\epsilon = \exp(-4(\log \rho)^2)$. Using Lemma~\ref{lemma:majreste}
we prove that $R_u(z)$ is bounded in absolute value by $\epsilon/2$ for any 
$z$ in $D(0,r)$. Indeed, the hypothesis  $u\ge \frac{4n^2}{(\log |z|)^2}$ in
Lemma~\ref{lemma:majreste}
results from Equation~(\ref{eq:condallzeros}) and the definition of $u$. Equation~(\ref{supRukappa})
then 
gives

\begin{eqnarray*}
 \log |R_u(z)| &\le &n^2+\log A+\frac{u\log r}{2}-(n+1)\log(1-r)\\
 &\le &  n^2+\log A-{8(\log \rho)^2}-(n+1)\log(1-r)
\end{eqnarray*}
from the
  definition of $u$.

Using (\ref{eq:condallzeros}) we deduce 

$$\log |R_u(z)| \le -7(\log \rho)^2 \le -\log 2+\log \epsilon.$$

We check that 

\begin{eqnarray}\label{eq:boundu}
u&\le& \frac{16(\log \rho)^2}{{1-r}}\\\nonumber
&\le & \Theta(\log \rho)^{2+\frac{1}{6}}\text{ because of   (\ref{eq:condallzeros}).}
\end{eqnarray}

Now, we prove that $R_v(z)$ is bounded in absolute value by $\epsilon$. 
Indeed

\begin{eqnarray*}
|R_v(z)| &\le & |R_u(z)|+(u-v)\exp(-6(\log \rho)^2)\\
&\le &  \epsilon/2+\Theta  (\log \rho)^{2+\frac{1}{6}}  \exp(-6(\log
\rho)^2)\text{ using equation (\ref{eq:boundu}), }\\
&\le & \epsilon/2+\exp(-5(\log \rho)^2) \text{  because of   (\ref{eq:condallzeros}),}\\
&\le & \epsilon.
\end{eqnarray*}

The principal part $P_v(x)$ is a degree $v-1$ polynomial. It has no more
than $u-1$ zeros. 
From Equation~(\ref{eq:boundu}) we deduce that $u\rho ^2\le \Theta\rho^2 (\log\rho)^{2+\frac{1}{6}}$.
So $4u\rho^2 \le \rho$ using   (\ref{eq:condallzeros}).
In particular $4u\rho^2$ is   smaller than $r$ and smaller
than $\rho$.
The interval  $[r-\rho,r]$ then consists of positive reals and must contain
at least one  real number $R$ such that $|R-|z||>2\rho^2$ for every root $z$ of
$P_v(x)$. We fix any such $R$.

Let 
$\cU$ be the union of all
 open disks $D(\xi,2\rho^2)$ where $\xi$ is any root of $P_v(x)$ lying in
$\bar D(0,R)$.
We check that $\cU$ is contained in $D(0,R)$.
Let $\cD=\bar D(0,R)-\cU$.  It is clear that $P_v(x)$ has no root in $\cD$.

\medskip 
$\bullet$ We first
assume that $|P_v(z)|>\epsilon $ for every $z$ on the boundary of $\cD$. We
apply  Rouch{\'e}'s theorem to $\bar D(0,R)$ and deduce that
 $f(x)$ and $P_v(x)$  have the same number of zeros
with  absolute value $\le R$.

Now let $\cV$ be any  connected component  of $\cU$. The boundary
of $\cV$ is contained in the boundary of $\cD$. Applying Rouch{\'e}'s
theorem to $\bar\cV$ we deduce that 
 $f(x)$ and $P_v(x)$  have the same number of zeros
in $\bar \cV$. Since $\bar\cV$ is the union of a number $< v$
of closed
disks of radius $2\rho^2$, we know that its diameter is $\le 4u\rho^2$.
We deduce that the zeros of $P_v(x)$ inside $\bar D(0,R)$ approximate
the zeros  of $f(x)$ inside $\bar D(0,R)$ within $4u\rho^2\le \rho$. More
precisely,
 let  $Z=[z_1]+[z_2]+\cdots +[z_J]$
be  the divisor of $f(x)$ restricted to $D(0,R)$. Then 
  the divisor  $\hat Z$
of $P_v(x)$ restricted to $D(0,R)$
has degree $J$ also. Assume that  $J\ge 1$ and let 
 $\hat Z=[\hat z_1]+[\hat z_2]+\cdots +[\hat z_J]$ be  the divisor
of $P_v(x)$ restricted to $D(0,R)$. Then
 there exists a permutation  of the indices
$\tau \in \cS_J$ such that $|z_{\tau (j)}-\hat z_j|\le  \rho$
for every $1\le j\le J$. So the lemma is proven in that case.

\medskip 

$\bullet$ We now 
assume that $|P_v(\zeta)|\le \epsilon $ for some
 $\zeta$ on the boundary of $\cD$. Since $|\zeta|\le r$ and
 $P_v(x)$ has no
zero in  $\bar D(\zeta,\rho^2)$, we deduce, using Lemma~\ref{lemma:anyzero},
that $P_v(z)$ is bounded by 

$$\exp(-0.01\times (1-r)^5\times |\log \rho^2|)
=\exp(0.02\times (1-r)^5\times \log \rho )$$
\noindent  for every
$z\in D(0,1/2)$. So for every such $z$  we have 

\begin{eqnarray*}
|f(z)|&\le &|P_v(z)|+|R_v(z)|
\le \exp(0.02\times (1-r)^5\times \log \rho) + \epsilon\\&\le& \exp(0.01\times
(1-r)^5\times \log \rho)
\end{eqnarray*}
\noindent  and the lemma is proven in that case. \hfill $\Box$

\medskip 

In the next two Sections~\ref{subsection:proofthzeros}
and \ref{subsection:proofalgo}
we use Lemma~\ref{lemma:allzeros} 
to prove Theorem~\ref{theorem:countingzeros} and Theorem~\ref{theorem:findingzeros}.

\subsection{Counting zeros of power series}\label{subsection:proofthzeros}

We now can give a proof of Theorem~\ref{theorem:countingzeros}. 

Let $K$ be a large enough positive constant and let $0<\rho<1$ 
be the unique real number such that 

\begin{equation}\label{eq:defrho}
|\log \rho|\times (1-r)^6 = K(\mu+n^2+\log A).
\end{equation}

We set $r'=(1+r)/2$ and check that $1-r'=(1-r)/2$. We apply 
Lemma~\ref{lemma:allzeros} to the series $f(x)$ over the disk $D(0,r')$.

It is clear that condition (\ref{eq:condallzeros}) is satisfied provided $K$
is large enough. We note also that 

$$u\le \frac{16(\log \rho)^2}{1-r'} =
\frac{16K^2(\mu+ n^2+\log
  A)^2}{(1-r')^{13}}.$$

We set $v=u$ and observe that, provided $K$ is large enough,
 the first conclusion in Lemma~\ref{lemma:allzeros} is not compatible
with the definition of $\rho$ given by Equation~(\ref{eq:defrho}).
Therefore, the second conclusion in Lemma~\ref{lemma:allzeros}
must hold true: there exists an $R$ such that $r'-\rho\le R \le r'$
and the number of zeros of $f(x)$ inside $D(0,R)$ is bounded  by
$u$.  If $K$ is large enough then $\rho$ is smaller than $r'-r=(1-r)/2$.
So $D(0,r)\subset D(0,R)$ and  the number of zeros in $D(0,r)$
is bounded by $u$.
 \hfill $\Box$
\subsection{An algorithm for finding zeros of power series}\label{subsection:proofalgo}

In this section we describe the algorithm announced by Theorem~\ref{theorem:findingzeros} for computing the zeros of a power series $f(x)$.
The input of the algorithm is described in the statement of 
Theorem~\ref{theorem:findingzeros} and at the beginning of 
Section~\ref{section:zerosseries}.
We are given a black box $\BOX_f$ for the coefficients of $f(x)$. We are given
two integers $a\ge 1$ and $n\ge 1$ such that $f(x)$ is of type $(A,n)$ where
$A=\exp(a)$. We are
given also an integer $\mu\ge 1$ such that there exists at least one $z$ in
$D(0,1/2)$
such that $|f(z)|>\exp(-\mu)$. We don't need the value of this $z$. Knowing
its existence suffices. Finally, we are given two integers $m\ge 1$ and $o\ge
2$
and we are supposed to compute approximations  within $10^{-m}$ 
for  the zeros of $f(x)$ in some  disk $D(0,r')$ where
$|r'-r|\le 10^{-m}$ and $r=1-1/o$.

\medskip

Let $K$ be a large enough constant. We assume that $K$ is an integer.
We set 

\begin{equation}\label{eq:defk'}
m'=K(m+ (\mu+n^2+a)o^6)\text{ and }
\rho'=\exp(-m').
\end{equation}

We apply Lemma~\ref{lemma:allzeros} to the series $f(x)$ over the disk $D(0,r)$ with the  accuracy
$\rho'$.
It is clear that condition (\ref{eq:condallzeros}) is satisfied provided $K$
is large enough. It is clear also that the first conclusion
in this lemma is not compatible with the definitions of 
$\mu$ and $m'$. So the second  conclusion must hold true. We take for $u$
the value 

$$u =  \frac{16(\log \rho')^2}{1-r}
= 16 o K^2(m+(\mu+ n^2+a)o^6)^2.$$

We first show that there exists at least one integer $w$ such that $1\le w\le u$ and
$|f_{w-1}|\ge 10^{-6(m')^2}$.
Otherwise, for any $z\in D(0,1/2)$ we would have 

\begin{eqnarray*}
|f(z)|&\le & |P_{u}(z)|+|R_u(z)|\\
&\le & u10^{-6(m')^2}+  \exp(-4(\log \rho ')^2)\text{ using equation
  (\ref{eq:borneRv})}\\
&\le &\exp(-\mu) 
\end{eqnarray*} 
 using the definition of $\rho'$ in
  (\ref{eq:defk'})  and assuming $K$ is large enough.
But the later inequality  contradicts the hypothesis on $\mu$ in 
Theorem~\ref{theorem:findingzeros}.

So let $v$ be an  integer  such that $1\le v\le u$ and
$|f_{v-1}|\ge 0.5\times 10^{-6(m')^2}$ and
$|f_{w}|\le  10^{-6(m')^2}$ for all  $w$ such that
$w\ge v$ and $w\le u-1$. Let $R$ be any
real such that 
$R\in [r-\rho',r]$ and $|R-|z||>2(\rho')^2$ for every root $z$ of $P_v(x)$.
Then, inside the disk $D(0,R)$, the zeros of $P_v(x)$ approximate the zeros of 
$f(x)$ within $\rho'$.
If $K$ is  large enough then $\rho' <0.5\times 10^{-m}$.

We call $P(x)=P_v(x)/f_{v-1}$ the unitary polynomial associated with
$P_v(x)$. Its coefficients are bounded above by $2Av^n10^{6(m')^2}$. So we can compute approximations $(\alpha_j)_{1\le j\le v-1}$
of the roots of $P(x)$ within
$(\rho')^2$ using the  algorithm provided by Theorem~\ref{theorem:quadtree}
of Section~\ref{section:roots}. We assume that the $\alpha_j$ are sorted by
increasing absolute value.
We consider the interval
$[r-\rho',r]$ and we remove to it all intervals 
$[|\alpha_j|-3(\rho')^2, |\alpha_j|+3(\rho')^2]$.
 Let  $r'$ be a rational number in this set. We count the $\alpha_j$
that belong to  the disk $D(0,r')$. Assume that there are $J$ of them. We output
the divisor $[\alpha_1]+[\alpha_2]+\cdots+[\alpha_J]$.

This leads to the following algorithm:

\begin{enumerate}
\item Set \\
$m'=K(m+o^6(\mu+n^2+a))$ and $u= 16K^2o(m+o^6(\mu+n^2+a))^2$
 and look for an 
integer $v$ such that \\$1\le v\le u$ and
$|f_{v-1}|\ge 0.5\times 10^{-6(m')^2}$ and
$|f_{w}|\le  10^{-6(m')^2}$ for all  $w$ such that
$w\ge v$ and $w\le u-1$.
\item Using the algorithm provided by Theorem~\ref{theorem:quadtree},
compute approximations $(\alpha_j)_{1\le j\le v-1}$ of the roots of
$P(x)=P_v(x)/f_{v-1}$ within $\exp(-2m')$.
\item Pick a rational number $r'$ in the set 
\\$[r-\exp(-m'),r]-\bigcup_{j}[|\alpha_j|-3\exp(-2m'), |\alpha_j|+3\exp(-2m')].$
\item Output $r'$ and those $\alpha_j$ that have absolute value $\le r'$.
\end{enumerate}

\subsection{Power series of type $(A,n)$}\label{subsection:elemAn}

In this section, we review simple elementary results about power series
in one variable.
Recall Definition~\ref{def:type} of the {\it type} of a power series in one
variable and let $f(x)=\sum_{k\ge 0}f_kx^k$ be a power series of
type $(A,n)$ with $A\ge 1$ a real number and $n\ge 1$ an integer.
Let $u\ge 0$ be an integer. We
write $f(x)=P_u(x)+R_u(x)$ where $P_u(x)=\sum_{0\le k\le
  u-1}f_kx^k$  is the principal part and $R_u(x)=\sum_{k\ge u}f_kx^k$
is the remainder term of order $u$. We first want
to  bound $|R_u(z)|$  for $z\in D(0,1)$.

\begin{eqnarray}
|R_u(z)|&=& \left| \sum_{k\ge u}  f_kz^k\right| \le   \sum_{k\ge u}
|f_k||z|^k \nonumber \\
&\le & A\sum_{k\ge u}(k+1)^n|z|^k \nonumber  \\
&\le&A|z|^u\sum_{k\ge 0}(k+u+1)^{n}|z|^k \nonumber  \\
&\le&A|z|^u(u+1)^n\sum_{k\ge 0}(k+1)^{n}|z|^k \nonumber \\
&\le&A|z|^u(u+1)^n\frac{n!}{(1-|z|)^{n+1}} \label{eq:supRu} 
\end{eqnarray}

We set 

\begin{equation}\label{eq:defkappa}
\kappa=\frac{n!A|z|^{\frac{u}{2}}}{(1-|z|)^{n+1}}
\end{equation}
\noindent  and we show that
$|R_u(z)|\le \kappa$ provided $u\ge \frac{4n^2}{(\log |z|)^2}$.
Indeed, if $u\ge \frac{4n^2}{(\log |z|)^2}$ then 
$n\log(u+1)\le n\sqrt u\le \frac{u\left|\log|z|\right|}{2}$ so

$$u\log|z| +n\log(u+1)\le \frac{u\log|z|}{2}.$$

Using (\ref{eq:supRu}) we deduce 

$$\log |R_u(z)|\le \log A+\frac{u\log|z|}{2}+\log\frac{n!}{(1-|z|)^{n+1}}=\log
\kappa.$$

We thus have  proven the following lemma.

\begin{lem}[Bounding the remainder of a power series]\label{lemma:majreste}
Let $f(x)$ be a power series of type $(A,n)$ where
$A\ge 1$ is a real number and $n\ge 1$ is an integer. Let $z$ be a complex
number such that $|z|<1$. Let $u\ge 0$ be an integer and let $R_u(x)
=\sum_{k\ge u}f_kx^k$
be the remainder of order $u$ of $f(x)$.
We have

$$|R_u(z)|\le A|z|^u(u+1)^n\frac{n!}{(1-|z|)^{n+1}}$$
\noindent and if $u\ge   \frac{4n^2}{(\log |z|)^2}$ then 

\begin{equation}\label{supRukappa}
|R_u(z)|\le \frac{n!A|z|^{\frac{u}{2}}}{(1-|z|)^{n+1}}
\end{equation}

\end{lem}

Note that if we set $u=0$ in (\ref{eq:supRu}) we obtain

\begin{equation}\label{eq:supfz}
|f(z)|\le \frac{n!A}{(1-|z|)^{n+1}}
\end{equation}

\medskip 

Now let $f(x)$ be a series of type $(A,n)$ and
let $c$ be a complex number such that $|c|<1$. We
set $f_c(x)=f(c+x(1-|c|))$ and call $f_c$ the {\it refocused}
series of $f(x)$ at $c$. 
We want to bound the coefficients of the power series $f_c(x)$.
These coefficients are related to the successive derivatives
of $f$ at $c$. Let $k\ge 0$ be an integer. We set

$$\rho=\frac{(k+1)+(n+1)|c|}{k+n+2}.$$

From Cauchy's integral  formula

\begin{eqnarray}\label{eq:plusCauchy}
|f^{(k)}(c)|&=&
\left|\frac{k !}{2i\pi }\int_{|\zeta|=\rho}\frac{f(\zeta)}{ (\zeta-c)^{k+1}}
d\zeta \right|\nonumber\\
&\le & n!k!A\frac{(k+n+2)^{k+n+2}}{(k+1)^{k+1}(n+1)^{n+1}} \frac{1}{(1-|c|)^{n+k+2}}
\end{eqnarray}
\noindent using (\ref{eq:supfz}).
 
We notice that

\begin{eqnarray}\label{eq:supBelyi}
\frac{(k+n+2)^{k+n+2}}{(k+1)^{k+1}(n+1)^{n+1}}&\le&
\left(1+\frac{n+1}{k+1}\right)^{k+1}\left(1+\frac{k+1}{n+1}\right)^{n+1}\nonumber
\\
&\le& \exp(n+1)2^{n+1}(k+1)^{n+1}.
\end{eqnarray}

From (\ref{eq:supfz}) and (\ref{eq:supBelyi}) we deduce that the refocused
power series $f_c(x)$ has  type 

$$(n!A(1-|c|)^{-n-2}\exp(n+1)2^{n+1},n+1).$$

We   thus have proven  the following lemma.

\begin{lem}[Refocusing a power series in one variable]\label{lemma:refocusing}
Let $$f(x)\in \CC[[x]]$$ be a power series of type $(A,n)$ where $A\ge 1$ 
is a
real
and $n\ge 1$ is an integer. Let $c$ be a complex number with absolute value
 smaller
than $1$. The refocused series $f_c(x)$ is defined as $f_c(x)=f(c+x(1-|c|))$. 
It is a power series of type $(n!A(1-|c|)^{-n-2}\exp(n+1)2^{n+1},n+1)$.
\end{lem}

%JMCfin

\chapter{Computations with modular forms and Galois
  representations}\label{chapcomput} 

\author{J. Bosman}

\bigskip

\bigskip

%% Chapter 2 (without 2.3.1 and 2.3.2) of Johan's thesis.
In this chapter we will discuss several aspects of the practical side
of computating with modular forms and Galois representations.  We
start by discussing computations with modular forms and from there on
work towards the computation of polynomials associated with modular
Galois representations.  Throughout this chapter, we will denote the
space of cusp forms of weight $k$, group $\Gamma_1(N)$ and character
$\varepsilon$ by~$S_k(N,\varepsilon)$.

\section{Modular symbols}
Modular symbols provide a way of doing symbolic calculations with
modular forms, as well as the homology of modular curves.  In this
section our intention is to give the reader an idea of what is going
on rather than a complete and detailed account of the material.  For
more details and further reading on the subject of modular symbols,
the reader could take a look at~\cite{manin}, \cite{Shokurov1}
and~\cite{merel}.  A computational approach to the material can be
found in~\cite{Stein-thesis} and~\cite{Stein}.
 
\subsection{Definitions}
Let $A$ be the free abelian group on the symbols $\{\alpha,\beta\}$
with $\alpha,\beta\in\PP^1(\QQ)$. Consider the subgroup $I\subset A$
generated by all elements of the forms:
\[
\{\alpha,\beta\}+\{\beta,\gamma\}+\{\gamma,\alpha\},
\quad
\{\alpha,\beta\}+\{\beta,\alpha\},
\quad\text{and}\quad
\{\alpha,\alpha\}.
\]
We define the group:
\[
\MM_2 := (A/I)/\text{torsion}\index{$M_k(\H)$@$\MM_k$} 
\]
as the quotient of $A/I$ by its torsion subgroup.  By a slight abuse
of notation, we will denote the class of $\{\alpha,\beta\}$ in this
quotient also by $\{\alpha,\beta\}$.  We have an action $\GL^+_2(\QQ)$
on $\MM_2$ by:
\[
\gamma\{\alpha,\beta\} := \{\gamma\alpha,\gamma\beta\},
\]
where $\gamma$ acts on $\PP^1(\QQ)$ by fractional linear
transformations.

For $k\geq 2$, we consider also the abelian group
$\ZZ[x,y]_{k-2}\subset\ZZ[x,y]$ of homogeneous polynomials of degree
$k-2$ and we let matrices in $\GL^+_2(\QQ)$ with integer coefficients
act on it on the left by:
\[
\Mat{a}{b}{c}{d}P(x,y) := 
P(dx-by,-cx+ay).
\]
We define:
\[
\MM_k := \ZZ[x,y]_{k-2}\otimes\MM_2,\index{$M_k(\H)$@$\MM_k$}
\]
and we equip $\MM_k$ with the component-wise action of integral
matrices in $\GL^+_2(\QQ)$ (that is $\gamma(P\otimes\alpha) =
\gamma(P)\otimes\gamma(\alpha)$).
\begin{defi}
Let $k\geq 2$ be an integer. Let $\Gamma\subset\SL_2(\ZZ)$ be a subgroup 
of finite index and let $I\subset \MM_k$ be the subgroup generated by all 
elements of the form $\gamma x-x$ with $\gamma\in\Gamma$ and $x\in\MM_k$. 
Then we define the space of \emph{modular symbols}\index{modular symbols} 
of weight $k$ for $\Gamma$ to be the quotient of $\MM_k/I$ by its torsion 
subgroup and we denote this space by 
$\MM_k(\Gamma)$:\index{$M_k(\H)$@$\MM_k$}
\[
\MM_k(\Gamma):=(\MM_k/I)/\text{torsion}.
\]
\end{defi}
In the special case $\Gamma=\Gamma_1(N)$, which we will mostly be
interested in, $\MM_k(\Gamma)$ is called the space of modular symbols
of weight $k$ and level $N$. The class of $\{\alpha,\beta\}$ in
$\MM_k(\Gamma)$ will be denoted by $\{\alpha,\beta\}_\Gamma$ or, if no
confusion exists, by $\{\alpha,\beta\}$.

The group $\Gamma_0(N)$ acts naturally on $\MM_k(\Gamma_1(N))$ and
induces an action of $(\ZZ/N\ZZ)^\times$ on $\MM_k(\Gamma_1(N))$.  We
denote this action by the diamond symbol~$\langle
d\rangle$\index{$D$@$\langle d\rangle$}.  The operator $\langle
d\rangle$ on $\MM_k(\Gamma_1(N))$ is called a \emph{diamond
  operator}\index{diamond operator}.  This leads to the notion of
modular symbols with character.
\begin{defi}
Let $\varepsilon\colon (\ZZ/N\ZZ)^\times\to\CC^\times$ be a Dirichlet
character.  Denote by ${\ZZ[\varepsilon]\subset\CC}$ the subring
generated by all values of $\varepsilon$. Let $I$ be the
$\ZZ[\varepsilon]$-submodule of
$\MM_k(\Gamma_1(N))\otimes\ZZ[\varepsilon]$ generated by all elements
of the form $\langle d\rangle x-\varepsilon(d)x$ with
$d\in(\ZZ/N\ZZ)^\times$ and $x\in\MM_k(\Gamma_1(N))$. Then we define
the space $\MM_k(N,\varepsilon)$\index{$M_k(\H)$@$\MM_k$} of modular
symbols\index{modular symbols} of weight $k$, level $N$ and character
$\varepsilon$ as the $\ZZ[\varepsilon]$-module:
\[
\MM_k(N,\varepsilon):=
\big(\MM_k(\Gamma_1(N))\otimes\ZZ[\varepsilon]/I\big)/\text{torsion}.
\]
\end{defi}
We denote the elements of $\MM_k(N,\varepsilon)$ by
$\{\alpha,\beta\}_{N,\varepsilon}$ or simply by $\{\alpha,\beta\}$. If
$\varepsilon$ is trivial, then we have
$\MM_k(N,\varepsilon)\cong\MM_k(\Gamma_0(N))$.

Let $\BB_2$ be the free abelian group on the symbols $\{\alpha\}$ with
$\alpha\in\PP^1(\QQ)$, equipped with action of $\SL_2(\ZZ)$ by
$\gamma\{\alpha\}=\{\gamma\alpha\}$ and define
$\BB_k$ as $\ZZ[x,y]_{k-2}\otimes\BB_2$ with component-wise
$\SL_2(\ZZ)$-action.  Elements of $\BB_k$\index{$B_{k,\q}$@$\BB_k$}
are called \emph{boundary modular symbols}\index{modular
  symbols!boundary}.  For a subgroup $\Gamma<\SL_2(\ZZ)$ of finite
index, we define $\BB_k(\Gamma)$ as:
\[
\BB_k(\Gamma):= (\BB_k/I)/\text{torsion}
\]
where $I$ is the subgroup of $\BB_k$ generated by all elements $\gamma
x - x$ with $\gamma\in\Gamma$ and $x\in\BB_k$. We define:
\[
\BB_k(N,\varepsilon) := 
\left(\left(\BB_k(\Gamma_1(N))\otimes\ZZ[\varepsilon]\right)/I\right)/
\text{torsion},
\]
where $I$ is the $\ZZ[\varepsilon]$-submodule of
$\BB_k(\Gamma_1(N))\otimes\ZZ[\varepsilon]$ generated by the elements
$\gamma x - \varepsilon(\gamma)x$ with~$\gamma\in\Gamma_0(N)$.

We have \emph{boundary homomorphisms}:
\[
\delta\colon \MM_k(\Gamma)\to\BB_k(\Gamma)\quad\text{and}\quad
\delta\colon \MM_k(N,\varepsilon)\to\BB_k(N,\varepsilon)
\]
defined by:
\[
\delta\left(P\otimes\{\alpha,\beta\}\right)=
P\otimes\{\beta\}-P\otimes\{\alpha\}. 
\]
The spaces of \emph{cuspidal modular symbols}\index{modular
  symbols!cuspidal}, denoted by
$\SS_k(\Gamma)$\index{$S_k(\H)$@$\SS_k$} and
$\SS_k(N,\varepsilon)$\index{$S_k(\H)$@$\SS_k$} respectively are
defined as the kernel of $\delta$.

\subsection{Properties}\label{modsymprops}
One can interpret the symbol $\{\alpha,\beta\}$ as a smooth path in
$\HH\cup\PP^1(\QQ)$ from the cusp $\alpha$ to the cusp $\beta$, lying
in $\HH$ except for the endpoints $\alpha$ and $\beta$. It can be
shown that this interpretation induces an isomorphism:
\[
\MM_2(\Gamma)\cong
H_1(X_\Gamma,\text{cusps},\ZZ).
\]
Here the homology is taken of the topological pair $(X_1(N),\text{cusps})$.
We also get an isomorphism:
\[
\SS_2(\Gamma)\cong H_1(X_\Gamma,\ZZ).
\]
So we immediately see that there is a perfect pairing:
\[
\left(\SS_2(\Gamma(N))\otimes\CC\right)
\times
\left(
S_2(\Gamma(N))\oplus \ol{S}_2(\Gamma(N))
\right)\to\CC
\]
defined by:
\[
\left(\{\alpha,\beta\},f\oplus g\right)
\mapsto
\int_\alpha^\beta\left(
f\frac{dq}{q} + g\frac{d\ol{q}}{\ol{q}}
\right).
\]
More generally, there is a pairing
\begin{eqn}\label{modsympairing}
\MM_k(\Gamma_1(N))
\times
\left(
S_k(\Gamma_1(N))\oplus \ol{S}_k(\Gamma_1(N))
\right)\to\CC
\end{eqn}
defined by:
\[
\left(P\otimes\{\alpha,\beta\},f\oplus g\right)
\mapsto
2\pi i\int_\alpha^\beta
\left(
f(z)P(z,1)dz
-
g(z)P(\ol{z},1)d\ol{z}
\right),
\]
which becomes perfect if we restrict and then tensor the left factor
to $\SS_k(\Gamma(N))\otimes\CC$. This pairing induces a pairing:
\[
\left(\MM_k(N,\varepsilon)\right)
\times 
\left(S_k(N,\varepsilon)\oplus\ol{S}_k(N,\varepsilon)\right)\to\CC
\]
which becomes perfect when the left factor is restricted and then
tensored to
${\SS_k(N,\varepsilon)\otimes_{\ZZ[\varepsilon]}\CC}$. From now on we
will denote all these pairings with the notation:
\[
\left(x,f\right)\mapsto \langle x,f\rangle.
\]

\subsection{The star involution}\index{star involution}
On the spaces $\MM_k(\Gamma_1(N))$ and $\MM_k(N,\varepsilon)$ we have
an involution~$\iota^*$:
\[
\iota^*(P(x,y)\otimes\{\alpha,\beta\}):=-P(x,-y)\otimes\{-\alpha,-\beta\},
\]
which is called the \emph{star involution}. It preserves cuspidal
subspaces.  We define $\SS_k(\Gamma_1(N))^+$ and
$\SS_k(\Gamma_1(N))^-$ subspaces of $\SS_k(\Gamma_1(N))$ where
$\iota^*$ acts as $+1$ and $-1$ respectively and we use similar
definitions for $\SS_k(N,\varepsilon)^\pm$.  It can be shown that the
pairing~(\ref{modsympairing}) induces perfect pairings:
\[
(\SS_k(\Gamma_1(N))^+\otimes\CC)\times S_k(\Gamma_1(N))\to\CC
\]
and:
\[
(\SS_k(\Gamma_1(N))^-\otimes\CC)\times \ol{S}_k(\Gamma_1(N))\to\CC
\]
and similarly for the spaces with character. This allows us to work
sometimes in modular symbols spaces of half the dimension of the full
cuspidal space.

\subsection{Hecke operators}\index{Hecke operator}
%% Left out because it made only sense in my thesis and not here:
% Hecke operators on modular symbols are defined in a similar way as
% on modular forms  
% (see Subsection \ref{subsecmodformhecke}). 
Let $k{\geq}2$ and $N{\geq}1$ be given. Then for
$\gamma\in\GL_2^+(\QQ)\cap\text{M}_2(\ZZ)$ we define an operator
$T_\gamma$ on $\MM_k(\Gamma_1(N))$ by letting
$\gamma_1,\ldots,\gamma_r$ be double coset representatives for
$\Gamma_1(N)\setminus \Gamma_1(N)\gamma\Gamma_1(N)$ and putting
\begin{eqn}\label{Tgammamodsym} %\indT
T_\gamma(x):=\sum_{i=1}^r \gamma_i x
\quad\text{for $x\in\MM_k(\Gamma_1(N))$}.
\end{eqn}
It follows from \cite[Theorem 4.3]{Shokurov1} that this operator is
well-defined.  For a prime number $p$ we put $T_p=T_\gamma$ for
$\gamma=\Mat{1}{0}{0}{p}$ and for positive integers $n$ we define
$T_n$ by means of the formal identity~\eqref{eqn_dirichlet_hecke}.
The operators $T_n$ are called Hecke operators.%\indT

The Hecke operators preserve the subspace $\SS_k(\Gamma_1(N))$ and
induce an action on the spaces $\MM_k(N,\varepsilon)$ and
$\SS_k(N,\varepsilon)$. Furthermore, from \cite[Theorem
  4.3]{Shokurov1} one can conclude that the diamond and Hecke
operators are self-adjoint with respect to the pairings defined in the
previous subsection: one has
\begin{eqn}\label{Tadjmodsym}
\langle Tx,f\rangle=\langle x,Tf\rangle
\end{eqn}
for any modular symbol $x$, cusp form $f$ and diamond or Hecke
operator $T$ for which this relation is well-defined.  Here, for
anti-holomorphic cusp forms we define the Hecke action by
$T\ol{f}=\ol{Tf}$.  Also, the Hecke operators commute with the star
involution $\iota^*$.

In conclusion, we have seen how we can write cusp forms spaces as the
dual of modular symbols spaces. The computation of Hecke operators on
these modular symbols spaces would enable us to compute $q$-expansions
of cusp forms: $q$-coefficients of newforms can be computed once we
can compute the eigenvalues of Hecke operators. But because
of~(\ref{Tadjmodsym}) this reduces to the computation of the
eigenvalues of Hecke operators on modular symbols spaces. In
computations one often works with the spaces
$\SS_k(N,\varepsilon)^+\otimes_{\ZZ[\varepsilon]}\QQ(\varepsilon)$
because these have smaller dimension than
$\SS_k(\Gamma_1(N))\otimes\QQ$.  Since we also know how all cusp forms
arise from newforms of possibly lower level (see~\eqref{eqn_basis}),
this allows us to compute the $q$-expansions of a basis for the spaces
$S_k(\Gamma_1(N))$ and $S_k(N,\varepsilon)$.  For precise details on
how these computations work, please read \cite[Chapter~9]{Stein}.

\subsection{Manin symbols}\index{Manin symbols}
If we want to do symbolic calculations with modular symbols, then the
above definitions are not quite applicable since the groups of which
we take quotients are not finitely generated. The \emph{Manin symbols}
enable us to give finite presentations for the spaces of modular
symbols.

First we need some definitions and lemmas. For a positive integer $N$
we define a set:
\[
E_N := \left\{(c,d)\in (\ZZ/N\ZZ)^2 : \gcd(N,c,d)=1\right\}.
\]
Define the following equivalence relation on~$E_N$:
\[
(c,d)\sim (c',d')
\ \stackrel{\text{def}}{\Longleftrightarrow}\ 
\text{there is an } a\in(\ZZ/N\ZZ)^\times\text{ with }(c,d)=(ac',ad')
\]
and the denote the quotient by $P_N$:
\begin{eqn}\label{defPN}
P_N := E_N/\sim.
\end{eqn}
The following lemma is easily verified:
\begin{lem}\label{gamma10sl2}
Let $N$ be a positive integer. Then the maps
\begin{align*}
&\Gamma_1(N)\setminus\SL_2(\ZZ)\to E_N:
\ol{\Mat{a}{b}{c}{d}}\mapsto (\ol{c},\ol{d})\quad\text{and}\\
&\Gamma_0(N)\setminus\SL_2(\ZZ)\to P_N:
\ol{\Mat{a}{b}{c}{d}}\mapsto \ol{(c,d)}
\end{align*}
are well-defined and bijective.
\end{lem}
This lemma enables us to write down an explicit set of coset
representatives for the orbit spaces $\Gamma_1(N)\setminus\SL_2(\ZZ)$
and $\Gamma_0(N)\setminus\SL_2(\ZZ)$.  The following lemma provides us
a first step in reducing the set of generators for the spaces of
modular symbols:
\begin{lem}
Each space $\MM_2(\Gamma_1(N))$ or $\MM_2(N,\varepsilon)$ is generated
by the symbols $\{a/c,b/d\}$ with $a,b,c,d\in\ZZ$ and $ad-bc=1$, where
in this notation a fraction with denominator equal to zero denotes the
cusp at infinity.
\end{lem}
\noindent
Calculating the continued fraction expansion at each cusp in $\QQ$
gives us immediately an algorithm to write a given element of $\MM_2$
in terms of the generators in the lemma. Furthermore, note that:
\[
\left\{\frac{a}{c},\frac{b}{d}\right\}=\Mat{a}{b}{c}{d}\{\infty,0\},
\]
so that we can write each element of $\MM_2$ as a sum of
$\gamma\{\infty,0\}$ with $\gamma\in\SL_2(\ZZ)$.

Let's consider the space $\MM_2(\Gamma_1(N))$. As we saw, it is
generated by the elements $\gamma\{\infty,0\}$ where $\gamma$ runs
through $\SL_2(\ZZ)$. Now, two matrices $\gamma$ define the same
element this way if they are in the same coset of the quotient
$\Gamma_1(N)\setminus\SL_2(\ZZ)$. According to Lemma~\ref{gamma10sl2}
such a coset can be uniquely identified with a pair
$(c,d)\in(\ZZ/N\ZZ)^2$. The corresponding element in
$\MM_2(\Gamma_1(N))$ is also denoted by $(c,d)$. This element $(c,d)$
is called a \emph{Manin symbol}. Clearly, there are only a finite
number of Manin symbols so we now know a finite set of generators for
$\MM_2(\Gamma_1(N))$.

For arbitrary $k$ we define the Manin symbols in $\MM_k(\Gamma_1(N))$
as the symbols of the form $P\otimes (c,d)$ where $P$ is a monomial in
$\ZZ[x,y]_{k-2}$ and $(c,d)$ a Manin symbol in
$\MM_2(\Gamma_1(N))$. In this case as well there are finitely many
Manin symbols and they generate the whole space.

In the modular symbols spaces with a given character~$\varepsilon$ we
have, for all $\gamma\in\Gamma_0(N)$, that
${\gamma(\alpha)=\varepsilon(\alpha)}$. Now for each element of $P_N$
we choose according to Lemma~\ref{gamma10sl2} a corresponding element
$\gamma\in\SL_2(\ZZ)$ and hence an element in $\MM_2(N,\varepsilon)$,
which we call again a Manin symbol. Note that this Manin symbol
depends on the choice of $\gamma$, but because of the relation
$\gamma(x)=\varepsilon(x)$ these chosen Manin symbols always form a
finite set of generators for $\MM_2(N,\varepsilon)$ as a
$\ZZ[\varepsilon]$-module. Likewise, $\MM_k(N,\varepsilon)$ is
generated by elements $P\otimes (c,d)$ with $P$ a monomial in
$\ZZ[x,y]_{k-2}$ and $(c,d)$ a Manin symbol in~$\MM_2(N,\varepsilon)$.

If we want to do symbolic calculations, then besides generators we
also need to know the relations between the Manin symbols. For
$\MM_k(\Gamma_1(N))$ one can do the following.
\begin{prop}\label{mkg1pres}
Let $N$ be a positive integer and let $A$ be the free abelian group on
the Manin symbols of the space $\MM_k(\Gamma_1(N))$. Let $I\subset A$
be the subgroup generated by the following elements:
\begin{align*}
P(x,y)\otimes(c,d) &\,+\, P(-y,x)\otimes(-d,-c),\cr
P(x,y)\otimes(c,d) &\,+\, 
P(-y,x-y)\otimes(-d,-c-d)\cr &\,+\, 
P(-x+y,-x)\otimes (-c-d,-c),\cr
P(x,y)\otimes(c,d) &\,-\, P(-x,-y)\otimes(c,d),
\end{align*}
where $P(x,y)\otimes (c,d)$ runs through all Manin symbols.  Then
$\MM_k(\Gamma_1(N))$ is naturally isomorphic to the quotient of $A/I$
by its torsion subgroup.
\end{prop}
\noindent
For the modular symbols spaces $\MM_k(N,\varepsilon)$ we have a
similar proposition.
\begin{prop}
Let $N$ and $\varepsilon$ be given. Let $A$ be the free
$\ZZ[\varepsilon]$-module on the Manin symbols of
$\MM_k(N,\varepsilon)$. Let $I\subset A$ be the submodule generated by
the elements given in Proposition~\ref{mkg1pres} plus for each
$n\in(\ZZ/N\ZZ)^\times$ the elements:
\[
P(x,y)\otimes\ol{(nc,nd)}\,-\,\varepsilon(n)P(x,y)\otimes\ol{(c,d)}.
\]
Then $\MM_k(N,\varepsilon)$ is naturally isomorphic to the quotient of
$A/I$ by its torsion submodule.
\end{prop}
These presentations enable us to perform symbolic calculations very
efficiently.

A remark on the computation of Hecke operators is in order here.  The
formula~(\ref{Tgammamodsym}) does not express the Hecke action on
Manin symbols in terms of Manin symbols. However, one can use other
formulas to compute Hecke operators.  The following theorem, due to
Merel, allows us to express Hecke operators more directly in terms of
Manin symbols:
\begin{thm}[see {\cite[Theorem 2]{merel}}]
On the spaces $\MM_k(\Gamma_1(N))$ and $\MM_k(N,\varepsilon)$ the
Hecke operator $T_n$ satisfies the following relation:
\[
T_n(P(x,y)\otimes(u,v))
=
\sumprime_{{a>b\geq 0\atop d>c\geq 0}\atop ad-bc=n}
P(ax+by,cx+dy)\otimes(au+cv,bu+dv),
\]
where the prime in the notation for the sum means that terms with
$\gcd(N,au+cv,bu+dv)\not=1$ have to be omitted.
\end{thm}

One would also like to express $\SS_k(\Gamma_1(N))$ and
$\SS_k(N,\varepsilon)$ in terms of the Manin symbols.  The following
proposition will help us.
\begin{prop}[see {\cite[Proposition 4]{merel}}]
Let integers $N\geq 1$ and $k\geq 2$ be given. Define an equivalence
relation on the vector space $\QQ[\Gamma_1(N)\setminus \QQ^2]$ by:
\[
[\ol{\lambda x}]\sim \mathrm{sign}(\lambda)^k[\ol{x}]
\quad\mbox{for $\lambda\in\QQ^\times$ and $x\in\QQ^2$}.
\]
Then the map:
\[
\mu\colon \BB_k(\Gamma_1(N))\to\QQ[\Gamma_1(N)\setminus \QQ^2]\,/\!\sim
\]
given by:
\[
\mu\colon P\otimes\left\{\frac{a}{b}\right\}\mapsto 
P(a,b)\left[\ol{\left(a\atop b\right)}
\right]\quad\mbox{($a$, $b$ coprime integers)}
\]
is well-defined and injective.
\end{prop}
The vector space $\QQ[\Gamma_1(N)\setminus \QQ^2]\,/\!\sim$ is finite
dimensional. The above proposition shows that $\SS_k(\Gamma_1(N))$ is
the kernel of $\mu\delta$, which is a map that can be computed in
terms of Manin symbols. The computation of $\SS_k(N,\varepsilon)$ can
be done in a similar way, see~\cite[Section 8.4]{Stein}.

\section{Intermezzo: Atkin-Lehner
  operators}\label{atleop}\index{Atkin-Lehner operator} 
In the rest of this chapter, we will be using the Atkin-Lehner
operators on $S_k(\Gamma_1(N))$ from time to time. This section
provides a brief treatment of the properties that we need.  The main
reference for this material is~\cite{Atkin-Li}.

Let $Q$ be a positive divisor of $N$ such that $\gcd(Q,N/Q)=1$. Let
$w_Q\in\GL_2^+(\QQ)$ be any matrix of the form:
\begin{eqn}\label{wQmatrix}
w_Q = \Mat{Qa}{b}{Nc}{Qd}
\end{eqn}
with $a,b,c,d\in\ZZ$ and $\det(w_Q)=Q$. The assumption $\gcd(Q,N/Q)=1$
ensures that such a $w_Q$ exists. A straightforward verification shows
$f|_kw_Q\in S_k(\Gamma_1(N))$. Now, given $Q$, this $f|_kw_Q$ still
depends on the choice of $a,b,c,d$. However, we can use a
normalisation in our choice of $a,b,c,d$ which will ensure that
$f|_kw_Q$ only depends on $Q$. Be aware of the fact that different
authors use different normalisations here. The one we will be using is:
\begin{eqn}\label{atlenorm}
a\equiv 1\bmod N/Q,\quad b\equiv 1 \bmod Q,
\end{eqn}
which is the normalisation used in~\cite{Atkin-Li}.  We define:
\begin{eqn}\label{defWQ}\index{$W_Q$}
\begin{aligned}
W_Q(f) :&= Q^{1-k/2}f|_kw_Q \\
&=\frac{Q^{k/2}}{(Ncz+Qd)^k}f\left(\frac{Qaz+b}{Ncz+Qd}\right),
\end{aligned}
\end{eqn}
which is now independent of the choice of $w_Q$ and call $W_Q$ an
\emph{Atkin-Lehner operator}.  In particular we have:
\[
W_N(f) = \frac{1}{N^{k/2}z^k}f\left(\frac{-1}{Nz}\right).
\]

An unfortunate thing about these Atkin-Lehner operators is that they
do not preserve the spaces $S_k(N,\varepsilon)$. But we can say
something about it. Let $\varepsilon\colon
(\ZZ/N\ZZ)^\times\to\CC^\times$ be a character and suppose that $f$ is
in $S_k(N,\varepsilon)$. By the Chinese Remainder Theorem, one can write
$\varepsilon$ in a unique way as
$\varepsilon=\varepsilon_Q\varepsilon_{N/Q}$ such that $\varepsilon_Q$
is a character on $(\ZZ/Q\ZZ)^\times$ and $\varepsilon_{N/Q}$ is a
character on $(\ZZ/(N/Q)\ZZ)^\times$.  It is a fact that:
\[
W_Q(f)\in S_k(N,\ol{\varepsilon}_Q\varepsilon_{N/Q}).
\]
Also, there is a relation between the $q$-expansions of $f$
and~$W_Q(f)$:
\begin{thm}\label{thmWQqexp}
Let $f$ be a newform in $S_k(N,\varepsilon)$. Take $Q$ dividing $N$
with $\gcd(Q,N/Q)=1$. Then:
\[
W_Q(f)=\lambda_Q(f)g
\]
with $\lambda_Q(f)$ an algebraic number in $\CC$ of absolute value 1 and
$g$  a newform in $S_k(N,\ol{\varepsilon}_Q\varepsilon_{N/Q})$. Suppose
now that $n$ is a positive integer and write $n=n_1n_2$ where $n_1$
consists only of prime factors dividing $Q$ and $n_2$ consists only of
prime factors not dividing $Q$. Then we have:
\[
a_n(g)=
\varepsilon_{N/Q}(n_1)
\ol{\varepsilon}_Q(n_2)
\ol{a_{n_1}(f)}
a_{n_2}(f).
\]
\end{thm}
The number $\lambda_Q(f)$ in the above theorem is called a
\emph{pseudo-eigenvalue}\index{pseudo-eigenvalue} for the Atkin-Lehner
operator.  In some cases there exists a closed expression for it.  In
the notation of the following theorem, $g(\chi)$ denotes the Gauss sum
of a Dirichlet character $\chi$ of conductor~$N(\chi)$:
\begin{eqn}\label{Gausssum}
g(\chi):=\sum_{\nu\in(\ZZ/N(\chi)\ZZ)^\times}\chi(\nu)
\exp\left(\frac{2\pi i\nu}{N(\chi)}\right).
\end{eqn}
\begin{thm}\label{lambdaqthm}
Let $f\in S_k(N,\varepsilon)$ be a newform and suppose $q$ is a prime
that divides $N$ exactly once. Then we have:
\[
\lambda_q(f)=
\left\{\begin{array}{ll}
g(\varepsilon_q)q^{-k/2}\ol{a_q(f)}
&\text{if $\varepsilon_q$ is non-trivial,}\cr
-q^{1-k/2}\ol{a_q(f)}&\mbox{if $\varepsilon_q$ is trivial.}
\end{array}\right.
\]
%Here, $g(\varepsilon_q)$ is the Gauss sum of $\varepsilon_q$.
\end{thm}
\begin{thm}[see {\cite[Theorem 2]{Asai1}}]\label{thmAsai}
Let $f\in S_k(N,\varepsilon)$ be a newform with $N$ square-free. For
$Q\mid N$ we have:
\[
\lambda_Q(f)=\varepsilon(Qd-\frac{N}{Q}a)\prod_{q\mid
  Q}\varepsilon(Q/q)\lambda_q(f). 
\]
Here, $a$ and $d$ are defined by~(\ref{wQmatrix}). Moreover, this
identity holds without any normalisation assumptions on the entries of
$w_Q$, as long as we define $\lambda_q(f)$ by the formula given in
Theorem~\ref{lambdaqthm}.
\end{thm}

\section{Basic numerical evaluations}
In this section we will describe how to perform basic numerical
evaluations, such as the evaluation of a cusp form at a point in $\HH$
and the evaluation of an integral of a cusp form between to points in
$\HH\cup\PP^1(\QQ)$. Again, the focus will be on performing actual
computations.

\subsection{Period integrals: the direct method}\label{subsecperint1}
In this subsection we will stick to the case $k=2$, referring to
\cite[Chapter~10]{Stein} for a more general approach (see also
\cite[Section~2.10]{cremona} for a treatment of~$\Gamma_0(N)$). So fix
a positive integer $N$ and an $f\in S_2(\Gamma_1(N))$. Our goal is to
efficiently evaluate the integral pairing $\langle x,f\rangle$ for
$x\in\SS_2(\Gamma_1(N))$.

Let us indicate why it suffices to look at newforms $f$. Because
of~\eqref{eqn_basis}, it suffices to look at $f=\alpha_d(f')$ with
$f'\in S_k(\Gamma_1(M))$ a newform for some $M\mid N$ and $d\mid
N/M$. By \cite[Theorem~4.3]{Shokurov1} we have:
\[
\langle x,f\rangle=\langle x,\alpha_d(f')\rangle = 
d^{1-k}\left\langle\Mat{d}{0}{0}{1}x,f'\right\rangle
\]
so that computing period integrals for $f$ reduces to computing period
integrals of the newform $f'$.

Let us now make the important remark that for each $z\in\HH$ we can
numerically compute $\int_\infty^z fdq/q$ by formally integrating the
$q$-expansion of~$f$:
\begin{eqn}\label{intser}
\int_\infty^zf\frac{dq}{q}=
\sum_{n\geq 1}\frac{a_n(f)}nq^n\quad\mbox{where $q=\exp(2\pi iz)$}.
\end{eqn}
The radius of convergence of this series is $1$ and the coefficients
are small (that is, estimated by $\tilde{O}(n^{(k-3)/2})$). So if $\Im
z\gg 0$ then we have $|q|\ll 1$ and the series converges rapidly. To
be more concrete, for $\Im z>M$ we have $|q^n|<\exp(-2\pi Mn)$ so if
we want to compute $\int_\infty^z fdq/q$ to a precision of $p$
decimals, we need to compute about $\frac{p\log 10}{2\pi M}\approx
0.37\frac{p}{M}$ terms of the series.

To compute a period integral we remark that for any
$\gamma\in\Gamma_1(N)$ and any $z\in\HH\cup\PP^1(\QQ)$ any continuous,
piecewise smooth path $\delta$ in $\HH\cup\PP^1(\QQ)$ from $z$ to
$\gamma z$, the homology class of $\delta$ pushed forward to
$X_1(N)(\CC)$ depends only on $\gamma$
\cite[Proposition~1.4]{manin}. Let us denote this homology class by:
\[
\{\infty,\gamma\infty\}\in\SS_2(\Gamma_1(N))\cong H_1(X_1(N)(\CC),\ZZ)
\]
and remark that all elements of $H_1(X_1(N)(\CC),\ZZ)$ can be written
in this way. As we also have $S_2(\Gamma_1(N))\cong
H^0(X_1(N)_\CC,\Omega^1)$, this means we can calculate
$\int_{\{\infty,\gamma\infty\}}f\frac{dq}{q}$ by choosing a smart path
in $\HH\cup\PP^1(\QQ)$:
\[
\int_\infty^{\gamma\infty}f\frac{dq}{q}\,=\,
\int_{z}^{\gamma z}f\frac{dq}{q}\,=\,
\int_{\infty}^{\gamma z}f\frac{dq}{q} \,-\, \int_{\infty}^zf\frac{dq}{q}.
\]
If we write $\gamma=\Mat{a}{b}{c}{d}$ then a good choice for $z$ is:
\[
z=-\frac{d}{c} + \frac{i}{|c|}.
\]
In this case we have $\Im z=\Im\gamma z=1/|c|$ so in view
of~(\ref{intser}), to compute the integral to a precision of $p$
decimals we need about $\frac{pc\log 10}{2\pi}\approx 0.37 pc$ terms
of the series.

Another thing we can use is the Hecke compatibility
from~(\ref{Tadjmodsym}). Put:
\[
W_f := \left(\SS_2(\Gamma_1(N))/I_f\SS_2(\Gamma_1(N))\right)\otimes\QQ,
\]
where $I_f$ is the Hecke ideal belonging to $f$, i.e. the kernel of
the map $\TT\to\CC$ that sends $T_n$ to $a_n$ for all $n$ (here, as
usual, $\TT$ denotes the Hecke algebra attached to
$S_2(\Gamma_1(N))$).  The space $W_f$ has the structure of a vector
space over $(\TT/I_f)\otimes\QQ\cong K_f$ of dimension $2$. This means
that computing any period integral of $f$, we only need to precompute
$2$ period integrals. So one tries to find a $K_f$-basis of $W_f$
consisting of elements $\{\infty,\gamma\infty\}$ where
$\gamma\in\Gamma_1(N)$ has a very small $c$-entry. In practice it
turns out that we do not need to search very far.

\subsection{Period integrals: the twisted method}\label{subsecperint2}
In this subsection we have the same set-up as in the previous
subsection.  There is another way of computing period integrals for
$f\in S_2(\Gamma_1(N))$ which sometimes beats the method described in
the previous subsection. The method described in this subsection is
similar to \cite[Section 2.11]{cremona} and makes use of winding
elements and twists.

The \emph{winding element}\index{winding element} of
$\MM_2(\Gamma_1(N))$ is simply defined as the element $\{\infty,0\}$
(some authors define it as $\{0,\infty\}$, this is just a matter of
sign convention). Integration over this element is easy to perform
because we can break up the path in a very neat way:
\begin{multline*}
\int_\infty^0f\frac{dq}{q} =\,
\int_\infty^{i/\sqrt{N}}f\frac{dq}{q}\,+\,\int_{i/\sqrt{N}}^0f\frac{dq}{q} =
\\ 
=\, \int_\infty^{i/\sqrt{N}}f\frac{dq}{q}\,+\,\int_{i/\sqrt{N}}^\infty
W_N(f)\frac{dq}{q} 
=\,\int_\infty^{i/\sqrt{N}}(f-W_N(f))\frac{dq}{q}.
\end{multline*}
Now, choose an odd prime $\ell$ not dividing $N$ and a primitive Dirichlet
character ${\chi\colon \ZZ\to\CC}$ of conductor $\ell$.  If $f\in
S_k(\Gamma_1(N))$ is a newform then $f\otimes\chi$ is a newform in
$S_k(\Gamma_1(N\ell^2))$, where:
\[
f\otimes\chi = 
\sum_{n\geq 1} a_n(f)\chi(n)q^n.
\]
The following formula to express $\chi$ as a linear combination of
additive characters is well-known:
\[
\chi(n)=\frac{g(\chi)}{\ell}
\sum_{\nu=1}^{\ell-1}\ol{\chi}(-\nu)
\exp\left(\frac{2\pi i \nu n}{\ell}\right),
\]
where $g(\chi)$ is the Gauss sum of $\chi$ (see~(\ref{Gausssum})).  It
follows now immediately that:
\[
f\otimes\chi \,=\,
\frac{g(\chi)}{\ell}
\sum_{\nu=1}^{\ell-1}\chi(-\nu)f\left(z+\frac\nu\ell\right)
\,=\,
\frac{g(\chi)}{\ell}
\sum_{\nu=1}^{\ell-1}\chi(-\nu)\,f\left|\Mat{\ell}{\nu}{0}{\ell}\right..
\]
For $f\in S_2(\Gamma_1(N))$ we now get the following useful formula
for free:
\begin{eqn}\label{intfchi}
\left\langle
\{\infty,0\},
f\otimes\chi
\right\rangle
=
\frac{g(\chi)}\ell
\left\langle
\sum_{\nu=0}^{l-1}
\chi(-\nu)\left\{\infty,\frac{\nu}{\ell}\right\},
f
\right\rangle.
\end{eqn}
The element 
$\sum_{\nu=0}^{l-1}\chi(-\nu)\left\{\infty,\frac{\nu}{\ell}\right\}$ of 
$\MM_k(\Gamma_1(N))\otimes\ZZ[\chi]$ or of some other modular symbols space 
where it is well-defined is called a 
\emph{twisted winding element}\index{winding element} or, more precisely the 
\emph{$\chi$-twisted winding element}. Because of formula~(\ref{intfchi}), we
can calculate the pairings of newforms in $S_2(\Gamma_1(N))$ with twisted 
winding elements quite efficiently as well.

We can describe the action of the Atkin-Lehner operator $W_{N\ell^2}$ on 
$f\otimes\chi$:
\[
W_{N\ell^2}(f\otimes\chi)=
\frac{g(\chi)}{g(\ol{\chi})}
\varepsilon(\ell)\chi(-N)
\lambda_N(f)
\tilde{f}\otimes\ol{\chi},
\]
where $\tilde{f}=\sum_{n\geq 1}\ol{a_n(f)}q^n$ (see for example
\cite[Section 3]{Atkin-Li}).  So in particular we have the following
integral formula for a newform $f\in S_2(N,\varepsilon)$:
\begin{eqn}\label{intfchi2}
\begin{aligned}
& \int_\infty^0f\otimes\chi\frac{dq}q =
\int_\infty^{i/(\ell\sqrt{N})}
(f\otimes\chi-W_{N\ell^2}(f\otimes\chi))\frac{dq}{q} = \\
&\quad = \int_\infty^{i/(\ell\sqrt{N})}
\left(f\otimes\chi-\frac{g(\chi)}{g(\ol\chi)}\chi(-N)\varepsilon(\ell)
\lambda_N(f)\tilde{f}\otimes\ol\chi\right)\frac{dq}q.
\end{aligned}
\end{eqn}
So to calculate:
\[
\left\langle\sum_{\nu=0}^{l-1}\chi(-\nu)\left\{\infty,
\frac{\nu}{\ell}\right\},f\right\rangle
\]
we need to evaluate the series~(\ref{intser}) at $z$ in $\HH$ with
$\Im z=1/(\ell\sqrt{N})$ which means that for a precision of $p$
decimals we need to sum about $\frac{p\ell\sqrt{N}\log 10}{2\pi}\approx
0.37p\ell\sqrt{N}$ terms of the series. In the spirit of the previous
subsection, we try several $\ell$ and $\chi$, as well as the untwisted
winding element $\{\infty,0\}$, until we can make a $K_f$-basis for
$W_f$. It follows from \cite[Theorems~1 and~3]{Shimura2} that we can
always find such a basis. Also here, it turns out that in practice we
do not need to search very far. The method that requires the least
amount of of $q$-expansion terms is preferred.

\subsection{Computation of $q$-expansions at various
  cusps}\index{q-expansion@$q$-expansion!at cusp} 
The upper half plane $\HH$ is covered by neighbourhoods of the
cusps. If we want to evaluate a cusp form $f\in S_k(\Gamma_1(N))$ or
an integral of a cusp form at a point in such a neighbourhood then it
is useful to be able to calculate the $q$-expansion of $f$ at the
corresponding cusp. We shall mean by this the following: A cusp $a/c$
can be written as $\gamma\infty$ with $
\gamma=\Mat{a}{b}{c}{d}\in\SL_2(\ZZ)$. Then a $q$-expansion of $f$ at
$a/c$ is simply the $q$-expansion of $f|_k\gamma$. This notation is
abusive, since it depends on the choice of $\gamma$. The $q$-expansion
will be an element of the power series ring $\CC[[q^{1/w}]]$ where $w$
is the width of the cusp $a/c$ and $q^{1/w}=\exp(2\pi iz/w)$.

If the level $N$ is square-free this can be done
symbolically. However, for general $N$ it is not known how to do this,
but we shall give some attempts that do at least give numerical
computations of $q$-expansions. We use that we can compute the
$q$-expansions of newforms in $S_k(\Gamma_1(N))$ at $\infty$ using
modular symbols methods.

\subsubsection*{The case of square-free $N$}
The method we present here is due to Asai~\cite{Asai1}.  Let $N$ be
square-free and let $f\in S_k(\Gamma_1(N))$ be a newform of character
$\varepsilon$. The main reason for being able to compute
$q$-expansions at all cusps in this case is because the group
generated by $\Gamma_0(N)$ and all $w_Q$ (see~(\ref{wQmatrix})) acts
transitively on the cusps, something that is not true when $N$ is not
square-free.

So let $\gamma=\Mat{a}{b}{c}{d}\in\SL_2(\ZZ)$ be given. Put:
\[
c'=\frac{c}{\gcd(N,c)},\quad\text{and}\quad Q=\frac{N}{\gcd(N,c)}.
\]
Let $r\in\ZZ$ be such that $d\equiv cr\bmod Q$ and define
$b',d'\in\ZZ$ by:
\[
Qd'=d-cr\quad\text{and}\quad b'=b-ar.
\]
Then we have:
\[
\Mat{a}{b}{c}{d}
=\Mat{Qa}{b'}{Nc'}{Qd'}
\Mat{Q^{-1}}{rQ^{-1}}{0}{1}.
\]
Theorems~\ref{thmWQqexp} and~\ref{thmAsai} tell us how
$\Mat{Qa}{b'}{Nc'}{Qd'}$ acts on $q$-expansions. The action of
$\Mat{Q^{-1}}{rQ^{-1}}{0}{1}$ on $q$-expansions is simply:
\[
\sum_{n\geq 1} a_nq^n\mapsto Q^{1-k}\sum_{n\geq 1}a_n\zeta_Q^{rn}q^{n/Q}
\quad\mbox{with $\zeta_Q=\exp(\frac{2\pi i}{Q})$}.
\]
This shows how the $q$-expansion of $f|_k\gamma$ can be derived from
the $q$-expansion of~$f$.

Let us now explain how to do it for oldforms as well. By induction
and~\eqref{eqn_basis} we may suppose $f=\alpha_p(f')$ with $p\mid N$
prime, $f'\in S_k(\Gamma_1(N/p))$ and that we know how to compute the
$q$-expansions of $f'$ at all the cusps. Let $\gamma=\Mat{a}{b}{c}{d}$
be given. Then we have:
\[
f|_k\gamma=p^{1-k}f'\big|_k\Mat{p}{0}{0}{1}\gamma=
p^{1-k}f'\big|_k\Mat{pa}{pb}{c}{d}.
\] 
We will now distinguish on two cases: $p\mid c$ and $p\nmid c$. If
$p\mid c$ then we have a decomposition:
\[
\Mat{pa}{pb}{c}{d}=\Mat{a}{pb}{c/p}{d}\Mat{p}{0}{0}{1}
\] 
and we know how both matrices on the right hand side act on $q$-expansions. 
If $p\nmid c$, choose $b',d'$ with $pad'-b'c=1$. Then we have:
\[
\Mat{pa}{pb}{c}{d}=\Mat{pa}{b'}{c}{d'}\beta
\]
with $\beta\in\GL_2^+(\QQ)$ upper triangular, so also in this case we
know how both matrices on the right hand side act on $q$-expansions.

\subsubsection*{The general case}
In a discussion with Peter Bruin, the author figured out an attempt to
drop the assumption that $N$ be square-free and compute $q$-expansions
of cusp forms numerically in this case. The idea is to generalise the
$W_Q$ operators from Section~\ref{atleop}.

So let $N$ be given. Let $Q$ be a divisor of $N$ and put
$R=\gcd(Q,N/Q)$. Let $w_Q$ be any matrix of the form:
\[
w_Q=\Mat{RQa}{b}{RNc}{Qd}
\quad\mbox{with $a,b,c,d\in\ZZ$}
\]
such that $\det w_Q=QR^2$ (the conditions guarantee us that such
matrices do exist). One can then verify:
\[
\Gamma_1(NR^2)<w_Q^{-1}\Gamma_1(N) w_Q,
\]
so that slashing with $w_Q$ defines a linear map:
\[
S_k(\Gamma_1(N))\oplus \ol{S}_k(\Gamma_1(N))
\stackrel{|w_Q}{\longrightarrow}
S_k(\Gamma_1(NR^2))\oplus\ol{S}_k(\Gamma_1(NR^2))
\]
which is injective since the slash operator defines a group action on
the space of all functions $\HH\to\CC$.

On the other hand, $w_Q$ defines an operation on $\MM_k$ which can be
shown to induce a linear map:
\[
w_Q\colon \SS_k(\Gamma_1(NR^2))\otimes\QQ\to\SS_k(\Gamma_1(N))\otimes\QQ
\]
that satisfies the following compatibility with respect to the
integration pairing between modular symbols and cusp forms (see
\cite[Theorem 4.3]{Shokurov1}):
\begin{eqn}\label{wQpairing}
\langle w_Qx,f\rangle=\langle x,f|_kw_Q\rangle.
\end{eqn}
Let $(x_1,\ldots,x_r)$ and $(y_1,\ldots,y_s)$ be bases of
$\SS_k(\Gamma_1(N))\otimes\QQ$ and of
${\SS_k(\Gamma_1(NR^2))\otimes\QQ}$ respectively. Then one can write
down a matrix $A$ in terms of these basis that describes the map $w_Q$
since we can express any symbol $P\otimes\{\alpha,\beta\}$ in terms of
Manin symbols. The matrix $A^t$ then defines the action of $w_Q$ in
terms of the bases of the cusp forms spaces that are dual to
$(x_1,\ldots x_r)$ and $(y_1,\ldots,y_s)$.

Now, let $(f_1,\ldots,f_r)$ be a basis of
$S_k(\Gamma_1(N))\oplus\ol{S}_k(\Gamma_1(N))$ and let
$(g_1,\ldots,g_s)$ be a basis of
$S_k(\Gamma_1(NR^2))\oplus\ol{S}_k(\Gamma_1(NR^2))$ (for instance we
could take bases consisting of eigenforms for the Hecke operators away
from $N$).  Define matrices:
\[
B:=\left(\langle x_i,f_j\rangle\right)_{i,j}
\quad\text{and}\quad
C:=\left(\langle y_i,g_j\rangle\right)_{i,j}.
\]
These can be computed numerically as the entries are period
integrals. Then the matrix $C^{-1}A^tB$ describes the map
$\cdot|_kw_Q$ in terms of the bases $(f_1,\ldots,f_r)$ and
$(g_1,\ldots,g_s)$. Hence if we can invert $C$ efficiently, then we
can numerically compute the $q$-expansion of $f|_kw_Q$ with $f\in
S_k(\Gamma_1(N))$.

Let now a matrix $\gamma=\Mat{a}{b}{c}{d}\in\SL_2(\ZZ)$ be given. Put:
\[
c':=\gcd(N,c)\quad\text{and}\quad Q:=N/c'.
\] 
Because of $\gcd(c/c',Q)=1$ we can find $\alpha\in(\ZZ/Q\ZZ)^\times$
with $\alpha c/c'\equiv 1\mod Q$. If we lift $\alpha$ to
$(\ZZ/N\ZZ)^\times$ then we have $\alpha c\equiv c'\mod N$. Let now
$d'\in\ZZ$ be a lift of $\alpha d\in (\ZZ/N\ZZ)^\times$. Because
$\alpha c$ and $\alpha d$ together generate $\ZZ/N\ZZ$ we have
$\gcd(c',d')=1$ and so we can find $a',b'\in\ZZ$ that
satisfy $a'd'-b'c'=1$. According to Lemma~\ref{gamma10sl2}, we have:
\[
\gamma=\gamma_0\Mat{a'}{b'}{c'}{d'}
\quad\mbox{with $\gamma_0\in\Gamma_0(N)$}.
\]
Put $R=\gcd(c',Q)$. Then we have $\gcd(NR,Q^2Ra')=QR\gcd(c',Qa')$ and
hence $\gcd(NR,Q^2Ra')=QR^2$, so there exist $b'',d''\in\ZZ$ with:
\[
w_Q:=\Mat{QRa'}{b''}{NR}{Qd''}
\]
having determinant $QR^2$. One can now verify that we have
$\Mat{a'}{b'}{c'}{d'}=w_Q\beta$ with $\beta\in\GL_2^+(\QQ)$ upper
triangular.  So in the decomposition:
\[
\gamma=\gamma_0w_Q\beta
\]
we can compute the slash action of all three matrices on the right
hand side in terms of $q$-expansions, hence also of $\gamma$.

In conclusion we see that in this method we have to increase the level
and go to $S_k(\Gamma_1(NR^2))$ for the square divisors $R^2$ of $N$
to compute $q$-expansions of cusp forms in $S_k(\Gamma_1(N))$ at
arbitrary cusps.

%BLABLA interessante vraag om te bekijken is hoe wQ met Hecke- en
%diamantoperatoren commuteert. 

\subsection{Numerical evaluation of cusp forms}\label{subsecevalmodform}
For $f\in S_k(\Gamma_1(N))$ and a point $P\in\HH$ we wish to compute
$f(P)$ to a high numerical precision. Before we do this let us say
some words on how $P$ should be represented. The transformation
property of modular forms implies that representing $P$ as $x+iy$ with
$x,y\in\RR$ is not a good idea, as this would be numerically very
unstable when $P$ is close to the real line.  Instead, we represent
$P$ as:
\begin{eqn}\label{pointinh}
P=\gamma z
\quad\text{with $\gamma\in\SL_2(\ZZ),$\, $z=x+iy,$\, $x\ll\infty$\,
  and\, $y\gg 0$}. 
\end{eqn}
For instance, one could demand that $z$ be in the standard fundamental
domain:
\[
\calF:=\{z\in\HH :\, |\Re z|\leq 1/2\,\text{ and }\,|z|\geq 1\}
\]
for $\SL_2(\ZZ)$ acting on $\HH$, although this is not strictly
necessary.

So let $P=\gamma z$ be given, with
$\gamma=\Mat{a}{b}{c}{d}\in\SL_2(\ZZ)$ and $\Im z > M$, say. Let
$w=w(\gamma)$ be the width of the cusp $\gamma\infty$ with respect to
$\Gamma_1(N)$. To compute $f(P)$ we make use of a $q$-expansion of $f$
at $\gamma\infty$:
\[
f(P)= (cz+d)^k (f|_k\gamma)(z) =
(cz+d)^k\sum_{n\geq 1} a_nq^{n/w}
%\quad\mbox{where $q^{1/w}=\exp(2\pi iz/w)$}.
\]
The radius of convergence is $1$ and the coefficients are small
(estimated by $\tilde{O}(n^{(k-1)/2})$). So to compute $f(P)$ to a
precision of $p$ decimals we need about $\frac{pw\log 10}{2\pi
  M}\approx 0.37\frac{pw}{M}$ terms of the $q$-expansion of
$f|_k\gamma$.

Of course, we have some freedom in choosing $\gamma$ and $z$ to write
down $P$.  We want to find $\gamma$ such that $P=\gamma z$ with $\Im
z/w(\gamma)$ as large as possible. In general, one can always write
$P=\gamma z$ with $z\in\mathcal{F}$ so one obtains:
\begin{eqn}\label{maxImz}
\max_{\gamma\in\SL_2(\ZZ)}\frac{\Im \gamma^{-1}P}{w(\gamma)}\geq 
\frac{\sqrt{3}}{2N}.
\end{eqn}
We see that in order to calculate $f(P)$ to a precision of $p$
decimals it suffices to use about $\frac{pN\log
  10}{\sqrt{3}\pi}\approx 0.42pN$ terms of the $q$-expansions at each
cusp. Although for many points $P$ there is a better way of writing it
as $\gamma z$ in this respect than taking $z\in\mathcal{F}$, it seems
hard to improve the bound $\frac{\sqrt{3}}{2N}$ in general.

We wish to adjust the representation sometimes from $P=\gamma z$ to
$P=\gamma'z'$ where $\gamma'\in\SL_2(\ZZ)$ is another matrix, for
instance because during our calculations $\Re z$ has become too large
or $\Im z$ has become too small (but still within reasonable
bounds). We can make $\Re z$ smaller by putting $z':=z-n$ for
appropriate $n\in\ZZ$ and putting
$\gamma':=\gamma\Mat{1}{n}{0}{1}$. Making $\Im z$ larger is rather easy
as well.  We want to find $\gamma''=\Mat{a}{b}{c}{d}\in\SL_2(\ZZ)$
such that:
\[
\Im\gamma''z=\frac{\Im z}{|cz+d|^2}
\]
is large. This simply means that we have to find a small vector $cz+d$
in the lattice $\ZZ z+\ZZ$, something which can be done easily if $\Re
z\ll\infty$ and $\Im z\gg 0$. If $c$ and $d$ are not coprime we can
divide both by their greatest common divisor to obtain a smaller
vector. The matrix $\gamma''$ can now be completed and we put
$z':=\gamma''z$ and $\gamma':=(\gamma'')^{-1}$.

\subsection{Numerical evaluation of integrals of cusp
  forms}\label{subsecevalint} 
In this subsection we will describe for $f\in S_2(\Gamma_1(N))$ and
$P\in\HH$ how to evaluate the integral $\int_\infty^P fdq/q$. As in
the previous subsection, we assume $P$ to be given by means
of~(\ref{pointinh}). The path of integration will be broken into two
parts: first we go from $\infty$ to a cusp $\alpha$ near $P$ and then
we go from $\alpha$ to~$P$.

\subsubsection*{Integrals over paths between cusps}
The pairing~(\ref{modsympairing}) gives a map:
\[
\Theta\colon \MM_2(\Gamma_1(N))\to 
\Hom_{\CC}\left(S_2(\Gamma_1(N)),\CC\right),
\]
which is injective when restricted to $\SS_2(\Gamma_1(N))$.  The image
of $\Theta$ is a lattice of full rank, hence the induced map:
\[
\SS_2(\Gamma_1(N))\otimes\RR
\to
\Hom_{\CC}\left(S_2(\Gamma_1(N)),\CC\right)
\]
is an isomorphism. In particular we obtain a map:
\[
\Phi\colon \MM_2(\Gamma_1(N))\to\SS_2(\Gamma_1(N))\otimes\RR,
\]
which is an interesting map to compute if we want to calculate
integrals of cusp forms along paths between cusps. The map $\Phi$ is
called a \emph{period mapping}\index{period mapping}.

The Manin-Drinfel'd theorem (see~\cite[Corollary 3.6]{manin}
and~\cite[Theorem~1]{Drinfeld1}) tells us that
$\im(\Phi)\subset\SS_2(\Gamma_1(N))\otimes\QQ$. This is equivalent to
saying that each degree 0 divisor of $X_1(N)$ which is supported on
cusps defines a torsion point of $J_1(N)$. The proof given
in~\cite{Drinfeld1} already indicates how to compute $\Phi$ with
symbolic methods: let $p$ be a prime that is $1\bmod N$. Then the
operator $p+1-T_p$ on $\MM_2(\Gamma_1(N))$ has its image in
$\SS_2(\Gamma_1(N))$. The same operator is invertible on
$\SS_2(\Gamma_1(N))\otimes\QQ$. So we simply have:
\[
\Phi = (p+1-T_p)^{-1}(p+1-T_p),
\]
where the rightmost $p{+}1{-}T_p$ denotes the map from
$\MM_2(\Gamma_1(N))$ to $\SS_2(\Gamma_1(N))$ and the leftmost
$p{+}1{-}T_p$ denotes the invertible operator on
$\SS_2(\Gamma_1(N))\otimes\QQ$. For other methods to compute $\Phi$,
see~\cite[Section~10.6]{Stein}. So we can express the integral of
$fdq/q$ between any two cusps $\alpha$ and $\beta$ in terms of period
integrals, which we have already seen how to compute:
\[
\int_\alpha^\beta f\frac{dq}{q} = \langle\Phi(\{\alpha,\beta\}),f\rangle.
\]

\subsubsection*{Integrals over general paths}
We can imitate the previous subsection pretty much. 
Write $P\in\HH$ as $P=\gamma z$ with $\gamma\in\SL_2(\ZZ)$ such that 
$\Im z/w(\gamma\infty)$ is as large as possible. Then we have 
\begin{eqn}\label{intgenpath}
\begin{aligned}
\int_\infty^P f\frac{dq}{q}&=
\int_\infty^{\gamma\infty} f\frac{dq}{q} + 
\int_{\gamma\infty}^{\gamma z} f\frac{dq}{q} \nonumber \\ &=
\int_\infty^{\gamma\infty} f\frac{dq}{q} + 
\int_\infty^z (f|_2\gamma)\frac{dq}{q}.
\end{aligned}
\end{eqn}
The integral $\int_\infty^{\gamma\infty} f\frac{dq}{q}$ is over a path
between two cusps so we can compute it by the above discussion and the
integral $\int_\infty^z (f|_2\gamma)\frac{dq}{q}$ can be computed
using the $q$-expansion of $f|_2\gamma$:
\[
\int_\infty^z (f|_2\gamma)\frac{dq}{q}=w\sum_{n\geq 1} \frac {a_n}{n}q^{n/w},
\]
where $w=w(\gamma)$, $q^{1/w}=\exp(2\pi iz/w)$ and $f|_2\gamma=\sum
a_nq^{n/w}$.  Because of~(\ref{maxImz}), computing about $\frac{pN\log
  10}{\sqrt{3}\pi}\approx 0.42pN$ terms of the series should suffice
to compute $\int_\infty^P f\frac{dq}{q}$ for any $P\in\HH$.

Note also that we can use formula~(\ref{intgenpath}) to compute the
pseudo-eigenvalue $\lambda_Q(f)$ by plugging in $\gamma=w_Q$ and a $z$
for which both $\Im z$ and $\Im w_Qz$ are high and for which
$\int_\infty^zW_q(f)dq/q$ is not too close to zero.
 
%% qwerty - Comment just put here to find this spot.

\section{Applying numerical calculations to Galois representations}
\label{secexplnumcomput}
Let $f$ be a newform (of some level and weight) and let
$\lambda\mid\ell$ be a prime of its coefficient field.  From
Section~\ref{sec_red_to_jac} we know that a residual Galois
representation $\ol\rho=\ol\rho_{f,\lambda}$ is attached to the pair
$(f,\lambda)$. The fixed field $K_\lambda$ of $\ker(\ol\rho)$ in
$\Qbar$ is a number field.  The results from
Chapter~\ref{sec_comp_mod_l_rep} point out that we know that computing
$\ol\rho$ essentially boils down to computing a polynomial that has
$K_\lambda$ as splitting field.  In this section we describe how
numerical calculations can be used to compute such a polynomial. We
will follow ideas from Chapter~\ref{chap_first_descr}.

Theorem~\ref{thm_red_to_wt_2} shows that we can reduce this problem to
the case of a form of weight $2$ in most interesting cases.  Hence we
will assume that $f$ is a newform in $S_2(\Gamma_1(N))$.  Assume that
the representation $\ol\rho_{f,\lambda}$ is absolutely irreducible and
let $\TT$ be the Hecke algebra acting on $J_1(N)$.  There is a
subspace $V_\lambda$ of $J_1(N)(\Qbar)[\ell]$ on which both $\TT$ and
$\GQ$ act, such that the action of $\GQ$ defines
$\ol\rho_{f,\lambda}$.

\subsection{Approximation of torsion points}
The Jacobian $J_1(N)_\CC$ can be described as follows. Pick a basis
$f_1,\ldots,f_g$ of $S_2(\Gamma_1(N))$. Put:
\[
\Lambda:=\left\{\int_\gamma(f_1,\ldots,f_g)\frac{dq}{q} : 
[\gamma]\in H_1(X_1(N)(\CC),\ZZ)\right\}\subset\CC^g.
\]
This is a lattice in $\CC^g$ of full rank. By the Abel-Jacobi theorem
we have an isomorphism:
\begin{align*}
J_1(N)(\CC)&\stackrel{\sim}{\longrightarrow}\CC^g/\Lambda, \\
\big[\sum_i\left([Q_i]-[R_i]\right)\big]&\mapsto 
\sum_i\int_{R_i}^{Q_i}(f_1,\ldots,f_g)\frac{dq}{q}.
\end{align*}
Let now a divisor $\sum_{i=1}^g[R_i]$ on $X_1(N)$ be given.
Identifying $J_1(N)(\CC)$ with $\CC^g/\Lambda$ in this way, we get a
birational morphism
\begin{align*}
\phi\colon \Sym^gX_1(N)(\CC) &\to\CC^g/\Lambda,\\
(Q_1,\ldots,Q_g) &\mapsto
\sum_{i=1}^g\int_{R_i}^{Q_i}(f_1,\ldots,f_g)\frac{dq}{q}.
\end{align*}
The homology group $H_1(X_1(N)(\CC),\ZZ)$ is canonically isomorphic to
the modular symbols space $\SS_2(\Gamma_1(N))$.  The period lattice
$\Lambda$ can thus be computed numerically using the methods from
Subsections~\ref{subsecperint1} and~\ref{subsecperint2}.  Since we can
compute the action of $\TT$ on $\SS_2(\Gamma_1(N))\cong\Lambda$, we
can write down the points in
$\frac{1}{\ell}\Lambda/\Lambda\subset\CC^g/\Lambda$ that correspond to
the points of $V_\lambda$.  The aim is now to compute the divisors on
$X_1(N)_\CC$ that map to these points along $\phi$.  In our
computations, we assume without proof that $V_\lambda$ lies beneath
the good locus of $\phi$, i.e.\ the map $X_1(N)^g\to\CC^g/\Lambda$
induced by $\phi$ is \'etale above $V_\lambda$.

We start calculating with a small precision. Let a non-zero $P$ in
$V_\lambda(\CC)\subset\CC^g/\Lambda$ be given.  First we try out a lot
of random points $Q{=}(Q_1,\ldots,Q_g)$ in~$X_1(N)(\CC)^g$. Here, each
$Q_i$ will be written as $Q_i=\gamma_i w_i$, with $\gamma_i$ in a set
of representatives for $\Gamma_1(N)\setminus\SL_2(\ZZ)$ and
$w_i\in\mathcal{F}$.  We can compute $\phi(Q)$ using methods from
Subsection~\ref{subsecevalint}.  We work with the point $Q$ for which
$\phi(Q)$ is closest to $P$.  If we in fact already know some points
$Q$ with $\phi(Q)$ approximately equal to a point in $V_\lambda(\CC)$,
then we could also take one of those points as a starting point $Q$ to
work with.

The next thing to do is adjust $Q$ so that $\phi(Q)$ comes closer to
$P$.  We'll make use of the Newton-Raphson approximation method.  Let
$\phi'\colon \HH^g\to\CC^g/\Lambda$ be the function defined by:
\[
\phi'(z_1,\ldots z_g)=\phi(\gamma_1z_1,\ldots\gamma_gz_g).
\]
We observe that for a small vector $h=(h_1,\ldots h_g)\in\CC^g$ we have:
\[
\phi'(w_1+h_1,\ldots,w_g+h_g) = \phi(Q)+ hD + O(\| h\|^2)
\]
with:
\[
D=\left.\left(
\begin{array}{ccc}
\frac{\partial \phi'_1}{\partial z_1} & \cdots 
& \frac{\partial\phi'_g}{\partial z_1} \cr
\vdots & \ddots & \vdots \cr
\frac{\partial \phi'_1}{\partial z_g} & \cdots 
& \frac{\partial\phi'_g}{\partial z_g} 
\end{array}
\right)\right|_{(w_1,\ldots,w_g)}
.
\]
From the definition of $\phi$ we can immediately deduce:
\[
\frac{\partial\phi'_i}{\partial z_j}(w_1,\ldots,w_g) = 
2\pi i\cdot (f_i|_2\gamma_j)(w_j),
\]
where we apologise for the ambiguous~$i$.  We can thus compute the
matrix $D$ using the methods of
Subsection~\ref{subsecevalmodform}. Now choose a small vector
$v=(v_1,\ldots,v_g)\in\CC^g$ such that $\phi(Q)+v$ is closer to $P$
than $\phi(Q)$ is. For example, $v$ can be chosen among all vectors of
a bounded length so that $\phi(Q)+v$ is closest to $P$. If we write:
\[
h = vD^{-1},
\]
then we expect $\phi'(w_1{+}h_1,\ldots,w_g{+}h_g)$ to be approximately
equal to $\phi(Q){+}v$. If this is not the case, then we try the same
thing with a smaller~$v$. It could be that this still fails, for
instance because we are too close to the bad locus of the
map~$\phi$. In that case, we start with a new random point~$Q$.

We repeat the above adjustments until we are (almost) as close as we
can get considering our calculation precision.  It might happen that
the $w_i$ become too wild, i.e. $|\Re w_i|$ becomes too large or $\Im
w_i$ becomes too small.  If this is the case we adjust the way we
write $Q_i$ as $\gamma_iw_i$ using the method described in
Subsection~\ref{subsecevalmodform}.  We can always replace the
$\gamma_i$ then by a small matrix in the same coset of
$\Gamma_1(N)\setminus\SL_2(\ZZ)$.

Once we have for each $P\in V_\lambda-\{0\}$ a point $Q$ such that
$\phi(Q)$ is approximately equal to $P$, we can start increasing the
precision.  We double our calculation precision and repeat the above
adjustments ($\phi(Q){+}v$ will in this case be equal to~$P$). We repeat
this a few times until we have very good approximations.

\subsection{Computation of polynomials}
Now, we will choose a function in $h\in\QQ(X_1(N))$ and evaluate it at
the components of the points in $\phi^{-1}(V_\lambda)$.  With the
discussion of Chapter~\ref{chap_first_descr} in mind, we want $h$ to take
values of small height. Since $h$ multiplies heights of points roughly
by $\deg(h)$, we want to find a function of small degree.  Take any
$k$ and a basis $h_1,\ldots, h_n$ of $S_k(\Gamma_1(N))$ such that the
$q$-expansions of the $h_i$ lie in $\ZZ[[q]]$ and such that the
exponents of the first non-zero terms of these $q$-expansions form a
strictly increasing sequence. We propose to use
$h=W_N(h_{n-1})/W_N(h_n)$ as a function to use (assuming $n\geq
2$). Remember from Section~\ref{sec_modforms} that $S_k(\Gamma_1(N))$
is the space of global sections of the line bundle
$\mathcal{L}=\omega^{\otimes k}(-\text{cusps})$ on $X_1(N)$, base
changed to~$\CC$.  Be aware of the fact that the cusp $\infty$ of
$X_1(N)$ is not defined over $\QQ$, but the cusp $0$ is.  Since we
demand the $q$-expansions to have rational coefficients, the sections
$W_N(h_1),\ldots,W_N(h_n)$, with $W_N$ an Atkin-Lehner operator, are
defined over $\QQ$ and they have increasing order at~$0$. One can now
verify that for $h=W_N(h_{n-1})/W_N(h_n)$ we have:
\[
\deg(h)\leq\deg(\mathcal{L})-v_\infty(h_{n-1})\leq\deg(\mathcal{L})-
\dim H^0(\mathcal{L})+2 \leq g+1.
\]
For $k=2$ and $g\geq 2$ we have $\mathcal{L}\cong\Omega^1(X_1(N))$ and
we get $g$ as an upper bound for $\deg(h)$. Using methods from
Subsection~\ref{subsecevalmodform}, we can evaluate $h$ numerically.
The author is not aware of a sophisticated method for finding a
function $h\in\QQ(X_1(N))$ of minimal degree in general; this minimal
degree is called the \emph{gonality} of the curve $X_1(N)$. Published
results on these matters seem to either be limited to $X_0(N)$ or to
concern only \emph{lower} bounds for the gonality of modular curves,
see for example~\cite{Abramovich1}, \cite[Chapter~3]{Baker-thesis}
or~\cite{Poonen1}.

Now put, for $P\in V_\lambda(\CC)-\{0\}$:
\[
\alpha_P=\sum_{i=1}^gh(Q_i),\quad
\text{where $\phi(Q_1,\ldots,Q_g)=P$}.
\]
We work out the product in:
\[
P_\lambda(x):=\prod_{P\in V_\lambda(\CC)-\{0\}} \left(x-\alpha_P\right) 
= \sum_{k=0}^n a_kx^k,
\quad\mbox{where $n=\deg P_\lambda$}.
\]
The coefficients $a_k$ are rational numbers that we have computed
numerically.  Since the height of $P_\lambda$ is expected to be not
too large, the denominators of the $a_k$ should have a relative small
common denominator. The LLL algorithm can be used to compute integers
$p_0,\ldots,p_{n-1}, q$ such that $|p_k-a_kq|$ is small for all $k$,
see \cite[Proposition 1.39]{LLL}.  If the sequence $(a_k)$ is
arbitrary, then we'll be able to find $p_k$ and $q$ such that
$|p_k-a_kq|$ is roughly of order $q^{-1/n}$ for each $k$, but not much
better than that. So if it happens that we find $p_k$ and $q$ with
$|p_k-a_kq|$ much smaller than $q^{-1/n}$ for all $k$, then we guess
that $a_k$ is equal to $p_k/q$. If we cannot find such $p_k$ and $q$
then we will double the precision and repeat all the calculations
described above.

Heuristically, the calculation precision that is needed to find the
true value of $a_k$ is about
$(1+1/n)\cdot\mathrm{height}(P_\lambda)/\log(10)$ decimals.  Another
way of finding rational approximations of the $a_k$ is by
approximating them using continued fractions. For this method, the
precision needed to find the true value of $a_k$ would be about
$2\cdot\mathrm{height}(P_\lambda)/\log(10)$ decimals.

Since the degree of $P_\lambda$ will be quite large, we won't be able
to do many further calculations with it. In particular it may be hard
to verify whether all the guesses we made were indeed correct.
Instead, we will look at the following variant. If $\mathfrak{m}$ is
the Hecke ideal of $f\bmod\lambda$, then $V_\lambda$ is a vector space
over $\TT/\mathfrak{m}$. The representation $\ol\rho_{f,\lambda}$
induces an action $\tilde\rho_\lambda$ of $\GQ$ on the set
$\PP(V_\lambda)$ of lines in $V_\lambda$. We can attach a polynomial
$\tilde{P}_\lambda$ to this projectivised representation
$\tilde\rho_\lambda$, analogously to the way this was done for
$\ol\rho$. This polynomial will have a smaller degree than
$P_\lambda$.  We put:
\[
\tilde{P}_\lambda(x)=\prod_{L\in\PP(V_\ell)}
\big(x\,\,-\!\!\sum_{P\in L-\{0\}}\alpha_P\big)=\sum_{k=0}^mb_kx^k,
\quad\mbox{where $m=\deg\tilde{P}_\lambda$}.
\]
As above, if the calculation precision is sufficient we can use
lattice reduction algorithms to compute the exact values of the $b_k$.

\subsection{Reduction of polynomials}
Although the polynomial $\tilde{P}_\lambda$ will not have a very huge
height, its height is still too large to do any useful computations
with it.  The first step in making a polynomial of smaller height
defining the same number field is computing the maximal order of that
number field.  Let $q$ be the common denominator of the coefficients
and put $p_k=b_kq$.  Consider the polynomial:
\[
Q(x)=q\cdot\tilde{P}_\lambda(x)=qx^m+p_{m-1}x^{m-1}+\cdots+p_0.
\]
We make ourselves confident that we correctly computed $Q(x)$
(although we won't prove anything at this point yet). For instance, we
verify that $Q(x)$ is irreducible and that its discriminant has the
prime factors of $N\ell$ in it.  We can also compute for several
primes $p$ not dividing $\disc(Q(x))$ the decomposition type of $Q(x)$
modulo~$p$ and verify that it could be equal to the cycle type of
$\tilde\rho(\Frob_p)$.  If not, we again double the precision and
repeat the above calculations.

Let now $\alpha$ be a root of $\tilde{P}_\lambda(x)$ and write down
the order:
\[
\calO := \ZZ+
\sum_{k=1}^{m-1}\left(\ZZ\cdot\sum_{j=0}^{k-1}a_{m-j}\alpha^{k-j}\right),
\]
which is an order that is closer to the maximal order than
$\ZZ[q\alpha]$ (see \cite[Subsection 2.10]{Lenstra_Hendrik_1}).  Being
confident in the correctness of $Q(x)$, we know where the number field
$K$ defined by it ramifies and thus we can compute its maximal order
(see \cite[Section 6 and Theorems 1.1 and~1.4]{Buchmann-Lenstra}).
Having done this, we embed $\calO_K$ as a lattice into $\CC^m$ in the
usual way and we use the LLL algorithm to compute a basis of small
vectors in $\calO_K$.  We can then search for an element of small
length in $\calO_K$ that generates $K$ over $\QQ$. Its defining
polynomial $\tilde{P}'_\lambda$ will have small coefficients. See
also~\cite{Cohen-Diaz-y-Diaz}.

In the computation of the polynomials $P_\lambda$ and
$\tilde{P}'_\lambda$ we made several guesses and assumptions that we
cannot prove to be correct.  In Chapter~\ref{chappgltau}, we work out
in special cases how we can use established parts of Serre's
conjecture to prove afterwards for polynomials of the style
$\tilde{P}'_\lambda$ that they indeed belong to the modular Galois
representations that we claim they belong to.  See~\cite{Bosman1} for
another example of this.  In the unlikely case that such tests may
fail we can of course make adjustments like choosing another function
$h$ or another divisor to construct~$\phi$.

\subsection{Further refinements}
The Jacobian $J_1(N)$ has large dimension (for $N$ prime this
dimension is $(N{-}5)(N{-}7)/24$).  It could be that our newform $f$
is an element of $S_2(\Gamma)$ with
$\Gamma_1(N)\lneq\Gamma<\Gamma_0(N)$. In that case we work with the
curve $X_\Gamma$, which is given its $\QQ$-structure by defining it as
a quotient of $X_1(N)$. The Jacobian $J_\Gamma$ of $X_\Gamma$ is
isogenous to an abelian subvariety of $J_1(N)$ that contains
$V_\lambda$, so this works perfectly well.

In the case $\Gamma=\Gamma_0(N)$ we can sometimes go a step
further. The operator $W_N$ on $X_0(N)$, sending $z$ to $-1/Nz$, is
defined over $\QQ$.  If $f$ is invariant under $W_N$, one can work
with the curve $X_0^+(N):=X_0(N)/\langle W_N\rangle$. Its Jacobian
$J_0^+(N)$ is isogenous to an abelian subvariety of $J_1(N)$ that
contains $V_\lambda$, so also here it works. Some words on the
computation of the homology of $X_0^+(N)$ are in order.  The action of
$W_N$ on $X_0(N)$ induces an action on $H_1(X_0(N)(\CC),\ZZ)$ and on
$H_1(X_0(N)(\CC),\text{cusps},\ZZ)$.  Since paths between cusps on
$X_0^+(N)(\CC)$ lift to paths between cusps on $X_0(N)(\CC)$ we have a
surjection:
\[
H_1(X_0(N),\text{cusps},\ZZ)\twoheadrightarrow 
H_1(X_0^+(N)(\CC),\text{cusps},\ZZ).
\]
The kernel of this map consists of the elements $[\gamma]$ in
$H_1(X_0(N),\text{cusps},\ZZ)$ satisfying ${W_N([\gamma])=-[\gamma]}$.
So modular symbols methods allow us to compute
$H_1(X_0^+(N)(\CC),\text{cusps},\ZZ)$ as a quotient of
$\MM_2(\Gamma_0(N))$.  Let $\BB_2^+(\Gamma_0(N))$ be the free abelian
group on the cusps of $X_0^+(N)(\CC)$ and define:
\[
\delta\colon H_1(X_0^+(N)(\CC),\text{cusps},\ZZ)\to \BB_2^+(\Gamma_0(N)),
\quad\{\alpha,\beta\}\mapsto \{\beta\}-\{\alpha\}.
\]
Then $H_1(X_0^+(N)(\CC))=\ker(\delta)$.

\chapter{Polynomials for projective representations of level one
  forms}\label{chappgltau} 

\author{J. Bosman}

\bigskip

\bigskip

%% Chapter 4 from Johan's thesis.

\section{Introduction}
In this chapter we explicitly compute mod-$\ell$ Galois
representations attached to modular forms. To be precise, we look at
cases with $\ell\leq 23$ and the modular forms considered will be cusp
forms of level $1$ and weight up to~$22$.  We present the result in
terms of polynomials associated with the projectivised representations.
As an application, we will improve a known result on Lehmer's
non-vanishing conjecture for Ramanujan's tau function
(see~\cite[p.\ 429]{Lehmer-vanishing}).
%% can the space \ of the previous line be broken???

To fix a notation, for any $k\in\ZZ$ satisfying $\dim
S_k(\SL_2(\ZZ))=1$ we will denote the unique normalised cusp form in
$S_k(\SL_2(\ZZ))$ by $\Delta_k$.  We will denote the coefficients of
the $q$-expansion of $\Delta_k$ by $\tau_k(n)$:
\[
\Delta_k(z)=\sum_{n\geq 1} \tau_k(n)q^n\in S_k(\SL_2(\ZZ)).
\]
From $\dim S_k(\SL_2(\ZZ))=1$ it follows that the numbers $\tau_k(n)$
are integers.  For every $\Delta_k$ and every prime $\ell$ there is a
continuous representation:
\[
\ol{\rho}_{\Delta_k,\ell}\colon \GQ\to\GL_2(\FF_\ell)
\] 
such that for every prime $p\not=\ell$ we have that the characteristic
polynomial of $\ol{\rho}_{\Delta_k,\ell}(\Frob_p)$ is congruent to
$X^2-\tau_k(p)X+p^{k-1}\bmod\ell$.  For a summary on the exceptional
representations $\ol{\rho}_{\Delta_k,\ell}$ and the corresponding
congruences for $\tau_k(n)$, see~\cite{Swinnerton-Dyer1}.

\subsection{Notational conventions}
Throughout this chapter, for every field $K$ we will fix an algebraic
closure $\overline{K}$ and all algebraic extension fields of $K$ will
be regarded as subfields of $\overline{K}$. Furthermore, for each
prime number $p$ we will fix an embedding
$\overline{\QQ}\hookrightarrow\overline{\QQ}_p$ and hence an embedding
$\GQp\hookrightarrow\GQ$, whose image we call $D_p$.  We will use
$I_p$ to denote the inertia subgroup of~$\GQp$.

All representations (either linear or projective) in this chapter will
be \emph{continuous}.  For any field $K$, a linear representation
$\rho\colon G{\to}\GL_n(K)$ defines a projective representation
$\tilde\rho\colon G{\to}\PGL_n(K)$ via the canonical map
$\GL_n(K){\to}\PGL_n(K)$.  We say that $\tilde\rho\colon
G{\to}\PGL_n(K)$ is \emph{irreducible} if the induced action of $G$ on
$\PP^{n-1}(K)$ fixes no proper subspace. So for $n=2$ this means that
every point of $\PP^1(K)$ has its stabiliser subgroup not equal
to~$G$.

\subsection{Statement of results}\label{statres}

\begin{thm}\label{PGLprop}
For every pair $(k,\ell)$ occurring the table in
Section~\ref{polytable}, let the polynomial $P_{k,\ell}$ be defined as
in that same table.  Then the splitting field of each $P_{k,\ell}$ is
the fixed field of $\operatorname{Ker}(\tilde\rho_{\Delta_k,\ell})$
and has Galois group $\PGL_2(\FF_\ell)$.  Furthermore, if
$\alpha\in\overline{\QQ}$ is a root of $P_{k,\ell}$ then the subgroup
of $\Gal(\overline{\QQ}/\QQ)$ fixing $\alpha$ corresponds via
$\tilde\rho_{\Delta_k,\ell}$ to a subgroup of $\PGL_2(\FF_\ell)$
fixing a point of~$\PP^1(\FF_\ell)$.
\end{thm}
For completeness we also included the pairs $(k,\ell)$ for which
$\rho_{k,\ell}$ is isomorphic to the action of $\GQ$ on the
$\ell$-torsion of an elliptic curve.  These are the pairs in the table
in Section~\ref{polytable} with $\ell=k-1$, as there the
representation is the $\ell$-torsion of $J_0(\ell)$, which happens to
be an elliptic curve for $\ell\in\{11,17,19\}$.  A simple calculation
with division polynomials \cite[Chapter II]{Lang-elliptic} can be used
to treat these cases. In the general case, one has to work in the more
complicated Jacobian variety $J_1(\ell)$, which has dimension $12$ for
$\ell=23$ for instance.

We can apply Theorem~\ref{PGLprop} to verify the following result.
\begin{cor}\label{Lehmer}
The non-vanishing of $\tau(n)$ holds for all $n$ such that:
\[
n<22798241520242687999\approx 2\cdot 10^{19}.
\]
\end{cor}
In~\cite{Jordan-Kelly}, the non-vanishing of $\tau(n)$ was verified 
for all $n$ such that:
\[
n<22689242781695999\approx 2\cdot 10^{16}.
\] 

To compute the polynomials, the author used the approach described in
Section~\ref{secexplnumcomput}.  After the initial computations some
of the polynomials had coefficients of almost 2000 digits, so
reduction techniques were absolutely necessary.  The used algorithms
do not give a proven output, so we have to concentrate on the
verification.  We will show how to verify the correctness of the
polynomials in Section~\ref{secproof} after setting up some
preliminaries about Galois representations in
Section~\ref{secgalrep}. In Section~\ref{secLehmer} we will point out
how to use Theorem~\ref{PGLprop} in a calculation that verifies
Corollary~\ref{Lehmer}. All the calculations were perfomed using {\sc
  Magma} (see~\cite{Magma}).

\section{Galois representations}\label{secgalrep}
This section will be used to state some results on Galois
representations that we will need in the proof of
Theorem~\ref{PGLprop}.
\subsection{Liftings of projective representations}
Let $G$ be a topological group, let $K$ be a topological field and let
$\tilde\rho\colon G\to\PGL_n(K)$ be a projective representation. Let
$L$ be an extension field of $K$.  By a \emph{lifting} of $\tilde\rho$
over $L$ we shall mean a representation $\rho\colon G\to\GL_n(L)$ that
makes the following diagram commute:
\[
\SelectTips{cm}{10}
\xymatrix{
G \ar [r]^-{\tilde\rho} \ar [d]_{\rho} & \PGL_n(K) \ar @{^(->}[d]\cr
\GL_n(L) \ar @{->>}[r] & \PGL_n(L)
}
\]
where the maps on the bottom and the right are the canonical ones. If
the field $L$ is not specified then by a lifting of $\tilde\rho$ we
shall mean a lifting over~$\overline{K}$.

An important theorem of Tate arises in the context of liftings. For
the proof we refer to \cite[Section 6]{Serre9}.  Note that in the
reference representations over $\CC$ are considered, but the proof
works for representations over arbitrary algebraically closed fields.

\begin{thm}[Tate]\label{tate2}
Let $\tilde\rho\colon \Gal(\ol{\QQ}/\QQ)\to\PGL_n(K)$ be a projective
representation of $\Gal(\ol{\QQ}/\QQ)$ over a field~$K$.  Then for
each prime number~$p$, there exists a lifting $\rho'_p\colon
D_p\to\GL_n(\overline{K})$ of $\tilde\rho|_{D_p}$.  If these liftings
$\rho'_p$ have been chosen so that all but finitely many of them are
unramified, then there is a unique lifting $\rho\colon
\Gal(\overline{\QQ}/\QQ)\to\GL_n(\overline{K})$ such that for all
primes $p$ we have:
\[
\rho|_{I_p}=\rho'_p|_{I_p}.
\]
\end{thm}

\begin{lem}\label{localunramified}
Let $p$ be a prime number and let $K$ be a field. Suppose that we are
given a projective representation ${\tilde\rho_p\colon
  \GQp\to\PGL_n(K)}$ that is unramified.  Then there exists a lifting
$\rho_p\colon \GQp\to\GL_n(\overline{K})$ of $\tilde\rho_p$ that is
unramified as well.
\end{lem}
\begin{proof}
Since $\tilde\rho$ is an unramified representation of $\GQp$, it
factors through $\Gal(\overline{\FF}_p/\FF_p)\cong\hat{\ZZ}$ and is
determined whenever we know the image of
${\Frob_p\in\Gal(\overline{\FF}_p/\FF_p)}$. By continuity, this image
is an element of $\PGL_n(K)$ of finite order, say of order $m$.  If we
take any lift $F$ of $\tilde\rho(\Frob_p)$ to $\GL_n(K)$ then we have
$F^m=a$ for some $a\in K^\times$.  So $F':=\alpha^{-1}F$, where
$\alpha\in\overline{K}$ is any $m$-th root of $a$, has order $m$ in
$\GL_n(\overline{K})$.  Hence the homomorphism
$\GQp\to\GL_n(\overline{K})$ obtained by the composition:
\[
\SelectTips{cm}{10}
\entrymodifiers={+!!<0pt,\fontdimen22\textfont2>}
\xymatrixcolsep{.25pc}
\xymatrix{
\GQp \ar @{->>}[rr]
&&
\Gal(\overline{\FF}_p/\FF_p) \ar [rrr]^-\sim
&&&
\hat{\ZZ} \ar @{->>}[rr]
&&
\ZZ/m\ZZ \ar [rrr]^-{1\mapsto F'}
&&&
\GL_n(\overline{K})
}
\]
lifts $\tilde\rho$ and is continuous as well as unramified.
\end{proof}

\subsection{Serre invariants and Serre's conjecture}
Let $\ell$ be a prime. A Galois representation $\rho\colon
\GQ\to\GL_2(\overline{\FF}_\ell)$ has a \emph{level} $N(\rho)$ and a
\emph{weight} $k(\rho)$.  The definitions were introduced by Serre
(see \cite[Sections 1.2 \& 2]{Serre8}). Later on, Edixhoven found an
improved definition for the weight, see \cite[Section 4]{Edixhoven3}.
The definitions agree in the cases of our interest, but in the general
formulation of Theorem~\ref{MoonTaguchi} later on, Edixhoven's
definition applies.  The level $N(\rho)$ is defined as the
prime-to-$\ell$ part of the Artin conductor of $\rho$ and equals $1$
if $\rho$ is unramified outside~$\ell$.  The weight is defined in
terms of the local representation $\rho|_{D_\ell}$; its definition is
rather lenghty so we will not write it out here. When we need results
about the weight we will just state them.  Let us for now mention that
one can consider the weights of the twists $\rho\otimes\chi$ of a
representation $\rho\colon \GQ\to\GL_2(\overline{\FF}_\ell)$ by a
character $\chi\colon \GQ\to\overline{\FF}_\ell^\times$. If one
chooses $\chi$ so that $k(\rho\otimes\chi)$ is minimal, then we always
have $1\leq k(\rho\otimes\chi)\leq\ell+1$ and we can in fact choose
our $\chi$ to be a power of the $\operatorname{mod} \ell$ cyclotomic
character.

Serre conjectured \cite[Conjecture 3.2.4]{Serre8} that if $\rho$ is
irreducible and odd, then $\rho$ belongs to a modular form of level
$N(\rho)$ and weight $k(\rho)$. Oddness here means that the image of a
complex conjugation has determinant $-1$.  A proof of this conjecture
in the case $N(\rho)=1$ has been published by Khare, building on ideas
of himself and Wintenberger:
\begin{thm}
[Khare \& Wintenberger, {\cite[Theorem~1.1]{Khare1}}]\label{serreconj}
Let $\ell$ be a prime number and let $\rho\colon
\Gal(\overline{\QQ}/\QQ)\to\GL_2(\overline{\FF}_\ell)$ be an odd
irreducible representation of level $N(\rho)=1$. Then there exists a
modular form $f$ of level $1$ and weight $k(\rho)$ which is a
normalised eigenform and a prime $\lambda\mid\ell$ of $K_f$ such that
$\rho$ and $\ol{\rho}_{f,\lambda}$ become isomorphic after a suitable
embedding of $\FF_\lambda$ into~$\overline{\FF}_\ell$.
\end{thm}

\subsection{Weights and discriminants}
If a representation $\rho\colon \GQ\to\GL_2(\overline{\FF}_\ell)$ is
wildly ramified at $\ell$ it is possible to relate the weight to
discriminants of certain number fields.  In this subsection we will
present a theorem of Moon and Taguchi, \cite[Theorem 3]{Moon-Taguchi},
on this matter and derive some results from it that are of use to us.
\begin{thm}[Moon \& Taguchi]
\label{MoonTaguchi} Consider a wildly ramified representation
from $\rho\colon\GQl\to\GL_2(\overline{\FF}_\ell)$. Let $\alpha\in\ZZ$
be such that $k(\rho\otimes\chi_\ell^{-\alpha})$ is minimal where
$\chi_\ell\colon \GQl\to\FF_\ell^\times$ is the mod $\ell$ cyclotomic
character.  Put $\tilde{k}=k(\rho\otimes\chi_\ell^{-\alpha})$, put
$d=\gcd(\alpha,\tilde{k}-1,\ell-1)$ and define $m$ in $\ZZ$ by letting
$\ell^m$ be the wild ramification degree of
$K:=\overline{\QQ}_\ell^{\operatorname{Ker}(\rho)}$ over
$\QQ_\ell$. Then we have:
\[
v_\ell(\mathcal{D}_{K/\QQ_\ell})=
\left\{
\begin{array}{ll}
1+\frac{\tilde{k}-1}{\ell-1}-\frac{\tilde{k}-1+d}{(\ell-1)\ell^m}
&\text{if\, $2\leq\tilde{k}\leq\ell$,} \cr
2+\frac{1}{(\ell-1)\ell}-\frac{2}{(\ell-1)\ell^m}
&\text{if\, $\tilde{k}=\ell+1$,} 
\end{array}
\right.
\]
where $\mathcal{D}_{K/\QQ_\ell}$ denotes the different of $K$ over
$\QQ_\ell$ and $v_\ell$ is normalised by $v_\ell(\ell)=1$.
\end{thm}
We can simplify this formula to one which is useful in our case.  In
the following corollaries $v_\ell$ denotes a valuation at a prime
above $\ell$ that is normalised by $v_\ell(\ell)=1$.
\begin{cor}\label{cordisc1}
Let $\tilde\rho\colon \Gal(\overline{\QQ}/\QQ)\to\PGL_2(\FF_\ell)$ be
an irreducible projective representation that is wildly ramified at
$\ell$.  Take a point in $\PP^1(\FF_\ell)$, let
${H\subset\PGL_2(\FF_\ell)}$ be its stabiliser subgroup and let $K$ be
the number field defined as:
\[
K = \overline{\QQ}^{\tilde\rho^{\,-1}(H)}.
\]
Then the $\ell$-primary part of $\disc(K/\QQ)$ is related to the
minimal weight $k$ of the liftings of $\tilde\rho$ by the following
formula:
\[
v_\ell(\disc(K/\QQ))=k+\ell-2.
\]
\end{cor}
\begin{proof}
Let $\rho$ be a lifting of $\tilde\rho$ of minimal weight.  Since
$\rho$ is wildly ramified, after a suitable conjugation in
$\GL_2(\overline{\FF}_\ell)$ we may assume:
\begin{eqn}\label{rhoIl}
\rho|_{I_\ell}=\Mat{\chi_\ell^{k-1}}{*}{0}{1},
\end{eqn}
where $\chi_\ell\colon I_\ell\to\FF_\ell^\times$ denotes the mod
$\ell$ cyclotomic character; this follows from the definition of
weight.  The canonical map
$\GL_2(\overline{\FF}_\ell)\to\PGL_2(\overline{\FF}_\ell)$ is
injective on the subgroup $\Mat{*}{*}{0}{1}$, so the subfields of
$\overline{\QQ}_\ell$ cut out by $\rho|_{I_\ell}$ and
$\tilde\rho|_{I_\ell}$ are equal, call them $K_2$.  Also, let
$K_1\subset K_2$ be the fixed field of the diagonal matrices in
$\operatorname{Im}\rho|_{I_\ell}$.  We see from~(\ref{rhoIl}) that in
the notation of Theorem~\ref{MoonTaguchi} we can put $\alpha=0$, $m=1$
and $d=\gcd(\ell-1,k-1)$.  So we have the following diagram of field
extensions:
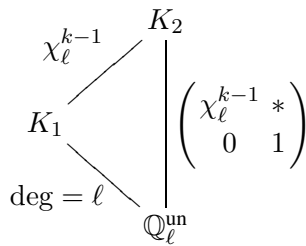
\[
\xymatrix{
& K_2 \ar @{-}[dd]^{\displaystyle{\Mat{{\chi_\ell^{k-1}}}{*}{0}{1}}}\cr 
K_1 \ar @{-}[ru]^{\displaystyle{\chi_\ell^{k-1}}} \cr
& \QQ_\ell^{\text{un}}\ar @{-}[lu]^{\displaystyle{\deg=\ell}}
}
\]
The extension $K_2/K_1$ is tamely ramified of degree $(\ell-1)/d$
hence we have:
\[
v_\ell(\mathcal{D}_{K_2/K_1})
=\frac{(\ell-1)/d-1}{(\ell-1)\ell/d} 
=\frac{\ell-1-d}{(\ell-1)\ell}.
\]
Consulting Theorem~\ref{MoonTaguchi} for the case $2\leq k \leq \ell$
now yields
\begin{align*}
v_\ell(\mathcal{D}_{K_1/\QQ_\ell^\text{un}})
&=
v_\ell(\mathcal{D}_{K_2/\QQ_\ell^\text{un}})
-
v_\ell(\mathcal{D}_{K_2/K_1})
\cr&=
1+\frac{k-1}{\ell-1}-\frac{k-1+d}{(\ell-1)\ell}
-
\frac{\ell-1-d}{(\ell-1)\ell}
=
\frac{k+\ell-2}{\ell}
\end{align*}
and also in the case $k=\ell+1$ we get:
\[
v_\ell(\mathcal{D}_{K_1/\QQ_\ell^\text{un}})
=
2+\frac{1}{(\ell-1)\ell}-\frac{2}{(\ell-1)\ell}
-
\frac{\ell-2}{(\ell-1)\ell}
=
\frac{k+\ell-2}{\ell}.
\]

Let $L$ be the number field
$\overline{\QQ}^{\operatorname{Ker}(\tilde\rho)}$. From the
irreducibility of $\tilde\rho$ and the fact that
$\operatorname{Im}\tilde\rho$ has an element of order $\ell$ it
follows that the induced action of $\GQ$ on $\PP^1(\FF_\ell)$ is
transitive and hence that $L$ is the normal closure of $K$ in
$\overline{\QQ}$. This in particular implies that $K/\QQ$ is wildly
ramified.  Now from $[K:\QQ]=\ell+1$ it follows that there are two
primes in $K$ above $\ell$: one is unramified and the other has
inertia degree $1$ and ramification degree $\ell$.  From the
considerations above it now follows that any ramification subgroup of
$\Gal(L/\QQ)$ at $\ell$ is isomorphic to a subgroup of
$\Mat{*}{*}{0}{1}\subset\GL_2(\overline{\FF}_\ell)$ of order
$(\ell-1)\ell/d$ with $d\mid \ell-1$.  Up to conjugacy, the only
subgroup of index $\ell$ is the subgroup of diagonal matrices.  Hence
$K_1$ and $K_{\lambda_2}^\text{un}$ are isomorphic field extensions of
$\QQ_\ell^\text{un}$, from which:
\[
v_\ell(\disc(K/\QQ))=
v_\ell(\disc(K_1/\QQ_\ell^\text{un}))=\ell\cdot
v_\ell(\mathcal{D}_{K_1/\QQ_\ell^\text{un}})=k+\ell-2. 
\]
follows.
\end{proof}
\begin{cor}\label{cordisc2}
Let $\tilde\rho\colon \Gal(\overline{\QQ}/\QQ)\to \PGL_2(\FF_\ell)$ be
an irreducible projective representation and let $\rho$ be a lifting
of $\tilde\rho$ of minimal weight.  Let $K$ be the number field
belonging to a point of $\PP^1(\FF_\ell)$, as in the notation of
Corollary~\ref{cordisc1}.  If $k\geq 3$ is such that:
\[
v_\ell(\disc(K/\QQ))=k+\ell-2
\]
holds, then we have $k(\rho)=k$. 
\end{cor}
\begin{proof}
From $v_\ell(\disc(K/\QQ))=k+\ell-2\geq\ell+1$ it follows that
$\tilde\rho$ is wildly ramified at $\ell$ so we can apply
Corollary~\ref{cordisc1}.
\end{proof}

\section{Proof of the theorem}\label{secproof}
To prove Theorem~\ref{PGLprop} we need to do several verifications. We
will derive representations from the polynomials $P_{k,\ell}$ and
verify that they satisfy the conditions of
Theorem~\ref{serreconj}. Then we know there are modular forms attached
to them that have the right level and weight and uniqueness follows
then easily.

First we we will verify that the polynomials $P_{k,\ell}$ from the
table in Section~\ref{polytable} have the right Galois group. The
algorithm described in \cite[Algorithm 6.1]{Geissler-Klueners} can be
used perfectly to do this verification; proving
$A_{\ell+1}\not<\Gal(P_{k,\ell})$ is the most time-consuming part of
the calculation here.  It turns out that in all cases we have:
\begin{eqn}\label{rightgalois}
\Gal(P_{k,\ell})\cong \PGL_2(\FF_\ell).
\end{eqn}
That the action of $\Gal(P_{k,\ell})$ on the roots of $P_{k,\ell}$ is
compatible with the action of $\PGL_2(\FF_\ell)$ follows from the
following lemma.
\begin{lem}
Let $\ell$ be a prime and let $G$ be a subgroup of
$\PGL_2(\FF_{\ell})$ of index $\ell{+}1$.  Then $G$ is the stabiliser
subgroup of a point in $\PP^1(\FF_\ell)$.  In particular, any
transitive permutation representation of $\PGL_2(\FF_\ell)$ of degree
$\ell{+}1$ is isomorphic to the standard action on $\PP^1(\FF_\ell)$.
\end{lem}
%\begin{proof} It suffices show that $G$ is contained in the
%stabiliser subgroup of a point: $G$ cannot be strictly contained in
%it because of its index.  From
%$\#\PGL_2(\FF_\ell)=(\ell-1)\ell(\ell+1)$ it follows that $G$
%has order $(\ell-1)\ell$. If $G$ has an orbit of length $\ell$ in
%$\PP^1(\FF_\ell)$ then it also has an orbit of length
%$1$ and we are done. So suppose that $G$ has no orbit of length
%$\ell$. Note that $G$ contains an $\ell$-Sylow subgroup of
%$\PGL_2(FF_\ell)$. But an $\ell$-Sylow subgroup of
%$\PGL_2(\FF_\ell)$ is conjugate to $\{[\Mat{1}{*}{0}{1}]\}$
%and has an orbit of lenght $\ell$. This implies that $G$ has to act
%transitively on $\PP^1(\FF_\ell)$ but all its orbits
%have length dividing $\ell-1$.  The statement about the permutation
%representation is now immediate.  \end{proof}
\begin{proof}
This follows from \cite[Proof of Theorem~6.25]{Suzuki}.
\end{proof}
So now we have shown that the second assertion in
Theorem~\ref{PGLprop} follows from the first one.

Next we will verify that we can obtain representations from this that
have the right Serre invariants.  Let us first note that every
automorphism of the group $\PGL_2(\FF_\ell)$ is an inner automorphism.
This implies that for every $P_{k,\ell}$, two
isomorphisms~(\ref{rightgalois}) define isomorphic representations
$\GQ{\to}\PGL_2(\FF_\ell)$ via composition with the canonical map
$\GQ{\twoheadrightarrow}\Gal(P_{k,\ell})$. In other words, every
$P_{k,\ell}$ gives a projective representation $\tilde\rho\colon
\GQ{\to}\PGL_2(\FF_\ell)$ that is well-defined up to isomorphism.

Now, for each $(k,\ell)$ in the table in Section~\ref{polytable}, the
polynomial $P_{k,\ell}$ is irreducible and hence defines a number
field:
\[
K_{k,\ell}:=\QQ[x]/(P_{k,\ell}),
\]
whose ring of integers we will denote by $\mathcal{O}_{k,\ell}$.  It
is possible to compute $\mathcal{O}_{k,\ell}$ using the algorithm from
\cite[Section 6]{Buchmann-Lenstra} (see also
\cite[Theorems~1.1~\&~1.4]{Buchmann-Lenstra}), since we know what kind
of ramification behaviour to expect.  In all cases it turns out that
we have:
\[
\disc(K_{k,\ell}/\QQ)=(-1)^{(\ell-1)/2}\ell^{k+\ell-2}.
\]
We see that for each $(k,\ell)$ the representation
$\tilde\rho_{k,\ell}$ is unramified outside $\ell$.  From
Lemma~\ref{localunramified} it follows that for each $p\not=\ell$, the
representation $\tilde\rho_{k,\ell}|_{\GQp}$ has an unramified
lifting.  Above we saw that via $\tilde\rho_{k,\ell}$ the action of
$\GQ$ on the set of roots of $P_{k,\ell}$ is compatible with the
action of $\PGL_2(\FF_\ell)$ on $\PP^1(\FF_\ell)$, hence we can apply
Corollary~\ref{cordisc2} to show that the minimal weight of a lifting
of $\tilde\rho_{k,\ell}$ equals $k$.  Theorem~\ref{tate2} now shows
that every $\tilde\rho_{k,\ell}$ has a lifting $\rho_{k,\ell}$ that
has level $1$ and weight $k$. From
$\operatorname{Im}\tilde\rho_{k,\ell}=\PGL_2(\FF_\ell)$ it follows
that each $\rho_{k,\ell}$ is absolutely irreducible.

To apply Theorem~\ref{serreconj} we should still verify that
$\rho_{k,\ell}$ is odd in each case.  I thank Robin de Jong for
pointing out that this is immediate: Since the weight of
$\rho_{k,\ell}$ is $k$, we have
$\det\rho_{k,\ell}|_{I_\ell}=\chi_\ell^{k-1}|_{I_\ell}$ where
$\chi_\ell$ is the mod $\ell$ cyclotomic character. Now,
$\det\rho_{k,\ell}|_{D_p}$ is unramified for $p\not=\ell$ and hence
$\det\rho_{k,\ell}$ must equal $\chi_\ell^{k-1}$ on all of $\GQ$
(apply Theorem~\ref{tate2} with $n=1$ for instance). But then, since
$k$ is even, we have that $\det\rho_{k,\ell}$ evaluated at a complex
conjugation equals $-1$ hence $\rho_{k,\ell}$ is odd.

So now that we have verified all the conditions of
Theorem~\ref{serreconj} we remark as a final step that all spaces of
modular forms $S_k(\SL_2(\ZZ))$ involved here are $1$-dimensional. So
the modularity of each $\rho_{k,\ell}$ implies immediately the
isomorphism $\rho_{k,\ell}\cong\ol{\rho}_{\Delta_k,\ell}$, hence also
$\tilde\rho_{k,\ell}\cong\tilde{\rho}_{\Delta_k,\ell}$ , which
completes the proof of Theorem~\ref{PGLprop}.

\section{Proof of the corollary}\label{secLehmer}
If $\tau$ vanishes somewhere, then the smallest positive integer $n$
for which $\tau(n)$ is zero is a prime.  This was observed by Lehmer
\cite[Theorem~2]{Lehmer-vanishing} and can also be seen using the
following argument: Suppose $n$ is the smallest positive integer with
$\tau(n)=0$.  From the multiplicative property of $\tau$ it follows
that $n$ is a power of a prime~$p$.  If $\tau(p)\not=0$ then from
$|\tau(p)|<p^6$ and the recursion for $\tau(p^r)$ it follows that
$v_p(\tau(p^r))=r\cdot v_p(\tau(p))$ for all $r$, so $\tau(p^r)$ would
never be zero.

Using results on the exceptional representations for $\tau(p)$, Serre
pointed out \cite[Section 3.3]{Serre7} that if $p$ is a prime number
with $\tau(p)=0$ then $p$ can be written as:
\[
p=hM-1
\]
with 
\begin{align*}
&M=2^{14}3^75^3691=3094972416000,\\
&\left(\frac{h+1}{23}\right)=1\quad\text{and}
\quad h\equiv 0,30 \text{ or } 48\bmod 49.
\end{align*}
In fact $p$ is of this form if and only if $\tau(p)\equiv 0\bmod
23\cdot 49\cdot M$ holds.  Knowing this, we will do a computer search
on these primes $p$ and verify whether $\tau(p)\equiv 0\bmod\ell$ for
$\ell\in\{11,13,17,19\}$.  To do that we will use the following lemma.
\begin{lem}\label{johan_lem_trace_0}
Let $K$ be a field of characteristic not equal to $2$. Then the
following conditions on $M\in\GL_2(K)$ are equivalent:
\begin{enumerate}
\item[(1)]
$\tr M=0$.
\item[(2)] For the action of $M$ on $\PP^1(K)$, there are $0$
  or $2$ orbits of length $1$ and all other orbits have length $2$.
\item[(3)]
The action of $M$ on $\PP^1(K)$ has an orbit of length $2$.
\end{enumerate}
\end{lem}
\begin{proof}
We begin with verifying (1)$\,\Rightarrow\,$(2).  Suppose $\tr
M=0$. Matrices of trace $0$ in $\GL_2(K)$ have distinct eigenvalues in
$\overline{K}$ because of $\operatorname{char}(K)\not=2$.  It follows
that two such matrices are conjugate if and only if their
characteristic polynomials coincide.  Hence $M$ and $M':=
\Mat{0}{1}{-\det M}{0}$ are conjugate so without loss of generality we
assume $M=M'$.  Since $M^2$ is a scalar matrix, all the orbits of $M$
on $\PP^1(K)$ have length $1$ or $2$.  If there are at least $3$
orbits of length $1$ then $K^2$ itself is an eigenspace of $M$ hence
$M$ is scalar, which is not the case. If there is exactly one orbit of
length $1$ then $M$ has a non-scalar Jordan block in its Jordan
decomposition, which contradicts the fact that the eigenvalues are
distinct.

The implication 
(2)$\,\Rightarrow\,$(3) is trivial so that leaves proving 
(3)$\,\Rightarrow\,$(1).
Suppose that $M$ has an orbit of length $2$ in
$\PP^1(K)$. After a suitable conjugation, we may assume that
this orbit is $\{[{1 \choose 0}],[{0 \choose 1}]\}$. But this means
that $M\sim\Mat{0}{a}{b}{0}$ for certain $a,b\in K$ hence $\tr M=0$.
\end{proof}
In view of the above lemma it follows from Theorem~\ref{PGLprop} that
for $\ell\in\{11,13,17,19\}$ and $p\not=\ell$ we have $\tau(p)\equiv
0\bmod\ell$ if and only if the prime $p$ decomposes in the number
field $\QQ[x]/(P_{12,\ell})$ as a product of primes of degree $1$ and
$2$, with degree $2$ occurring at least once.  For
$p\nmid\disc(P_{12,\ell})$, which is a property that all primes $p$
satisfying Serre's criteria possess, we can verify this condition by
checking whether $P_{12,\ell}$ has an irreducible factor of degree $2$
over $\FF_p$. This can be easily checked by verifying:
\[
\ol{x}^{p^2}=
\ol{x}\quad\text{and}\quad\ol{x}^p\not= \ol{x}
\quad\text{in}\quad \FF_p[x]/(\ol{P}_{12,\ell}).
\]
Having done a computer search, it turns out that the first few
primes satisfying Serre's criteria as well as $\tau(p)\equiv 0\bmod
11\cdot 13\cdot 17\cdot 19$ are
\begin{align*}
&22798241520242687999,\ 60707199950936063999,\\ 
&93433753964906495999,\ 102797608484376575999.
\end{align*}
\begin{remark}
The unpublished paper~\cite{Jordan-Kelly} in which Jordan and Kelly
obtained the previous bound for the verification of Lehmer's
conjecture seems to be unfindable.  Kevin Buzzard asked me the
question what method they could have used.  If we weaken the above
search to using only the prime $\ell=11$ we obtain the same bound as
Jordan and Kelly did. So our speculation is that they searched for
primes $p$ satisfying Serre's criteria as well as $\tau(p)\equiv
0\bmod 11$. This congruence can be verified using an elliptic curve
computation, as was already remarked in Subsection~\ref{statres}.
\end{remark}

\section{The table of polynomials}\label{polytable}
In this section we present the table of polynomials that is referred
to throughout this chapter.
%\begin{center}
% The table should be numbered as a subsection, preferably
\begin{longtable}{|c|l|}
\caption*{Polynomials belonging to projective modular representations}
\cr

\hline{}&{}\cr \multicolumn{1}{|c|}{$(k,\ell)$} & 
\multicolumn{1}{c|}{$P_{k,\ell}$}\cr{}&{}\cr \hline 
\endfirsthead

\multicolumn{2}{c}%
{{\tablename\ continued from previous page}} \cr
\hline{}&{}\cr \multicolumn{1}{|c|}{$(k,\ell)$} &
\multicolumn{1}{c|}{$P_{k,\ell}$}\cr{}&{}\cr\hline 
\endhead

%\multicolumn{2}{|r|}{{Continued on next page}} \cr \hline
\multicolumn{2}{r}{\footnotesize{Continued on next page}}
\endfoot

\hline
\endlastfoot

$(12,11)$&$x^{12} - 4x^{11} + 55x^9 - 165x^8 + 264x^7 - 341x^6 + 330x^5
$\\*&${} - 165x^4 - 55x^3 + 99x^2 - 41x - 111$\cr
\hline
$(12,13)$&$x^{14} + 7x^{13} + 26x^{12} + 78x^{11} + 169x^{10} + 52x^9 - 702x^8
$\\*&${}- 1248x^7 + 494x^6 + 2561x^5 + 312x^4 - 2223x^3 
$\\*&${}+ 169x^2 + 506x - 215$\cr
\hline
$(12,17)$&$x^{18} - 9x^{17} + 51x^{16} - 170x^{15} + 374x^{14} -
578x^{13} $\\*&${}+ 493x^{12}
- 901x^{11}+ 578x^{10} - 51x^9 + 986x^8 + 1105x^7 $\\*&${}+ 476x^6
+ 510x^5 + 119x^4 + 68x^3 + 306x^2 + 273x + 76$
\cr\hline
$(12,19)$&$x^{20} - 7x^{19} + 76x^{17} - 38x^{16} - 380x^{15} + 114x^{14} 
$\\*&${}+ 1121x^{13}
- 798x^{12} - 1425x^{11} + 6517x^{10} + 152x^9 $\\*&${}- 19266x^8 
- 11096x^7
 + 16340x^6 + 37240x^5 + 30020x^4 $\\*&${}- 17841x^3 - 47443x^2 
- 31323x - 8055$\cr
\hline
$(16,17)$&$x^{18} - 2x^{17} - 17x^{15} + 204x^{14} - 1904x^{13} 
+ 3655x^{12} $\\*&${}+ 5950x^{11}
 - 3672x^{10} - 38794x^9 + 19465x^8 + 95982x^7 $\\*&${}- 280041x^6
 - 206074x^5 + 455804x^4 + 946288x^3 $\\*&${}- 1315239x^2 + 606768x 
- 378241$\cr
\hline
$(16,19)$&$x^{20} + x^{19} + 57x^{18} + 38x^{17} + 950x^{16} +
4389x^{15} $\\*&${}+ 20444x^{14}
 + 84018x^{13} + 130359x^{12} - 4902x^{11} $\\*&${}- 93252x^{10} 
+ 75848x^9
 - 1041219x^8 - 1219781x^7 $\\*&${}+ 3225611x^6 + 1074203x^5 
- 3129300x^4 - 2826364x^3 $\\*&${}+ 2406692x^2 + 6555150x 
- 5271039$\cr
\hline
$(16,23)$&$x^{24} + 9x^{23} + 46x^{22} + 115x^{21} - 138x^{20} - 
1886x^{19} $\\*&${}+ 1058x^{18} 
+ 59639x^{17} + 255599x^{16} + 308798x^{15} $\\*&${}- 1208328x^{14}
 - 6156732x^{13} - 10740931x^{12} $\\*&${}+ 2669403x^{11} 
+ 52203054x^{10}
 + 106722024x^9 $\\*&${}+ 60172945x^8  - 158103380x^7- 397878081x^6
$\\*&${} - 357303183x^5 + 41851168x^4 + 438371490x^3 $\\*&${}+ 484510019x^2
 + 252536071x + 55431347$\cr
\hline
$(18,17)$&$x^{18} - 7x^{17} + 17x^{16} + 17x^{15} - 935x^{14} +
799x^{13} $\\*&${} + 9231x^{12}
- 41463x^{11} + 192780x^{10} + 291686x^9 $\\*&${}- 
390014x^8 + 6132223x^7  -3955645x^6 
+ 2916112x^5 $\\*&${}+ 45030739x^4 - 94452714x^3 
+ 184016925x^2 $\\*&${}- 141466230x + 113422599$\cr
\hline
$(18,19)$&$x^{20} + 10x^{19} + 57x^{18} + 228x^{17} - 361x^{16} -
3420x^{15} $\\*&${}+ 23446x^{14}
 + 88749x^{13} - 333526x^{12} - 1138233x^{11} $\\*&${}+ 1629212x^{10} 
+ 13416014x^9 + 7667184x^8 $\\*&${}- 208954438x^7 + 95548948x^6
 + 593881632x^5 $\\*&${}- 1508120801x^4 - 1823516526x^3 
+ 2205335301x^2 $\\*&${}+ 1251488657x - 8632629109$\cr    
\hline
$(18,23)$&$x^{24} + 23x^{22} - 69x^{21} - 345x^{20} - 483x^{19} -
6739x^{18} $\\*&${}+ 18262x^{17}
 + 96715x^{16} - 349853x^{15} + 2196684x^{14} $\\*&${}- 7507476x^{13}
 + 59547x^{12} + 57434887x^{11} $\\*&${}- 194471417x^{10} + 545807411x^9 
 + 596464566x^8 $\\*&${}- 9923877597x^7 + 33911401963x^6 
- 92316759105x^5 $\\*&${}+ 157585411007x^4 - 171471034142x^3  
 $\\*&${}+ 237109280887x^2 - 93742087853x  + 97228856961$\cr
\hline
$(20,19)$&$x^{20} - 5x^{19} + 76x^{18} - 247x^{17} + 1197x^{16} -
8474x^{15} $\\*&${}+ 15561x^{14}
 - 112347x^{13} + 325793x^{12} - 787322x^{11} $\\*&${}+ 3851661x^{10}
 - 5756183x^9 + 20865344x^8 $\\*&${}- 48001353x^7 + 45895165x^6 
 - 245996344x^5 $\\*&${}+ 8889264x^4 - 588303992x^3 - 54940704x^2 
$\\*&${}- 538817408x + 31141888$\cr 
\hline
$(20,23)$&$x^{24} - x^{23} - 23x^{22} - 184x^{21} - 667x^{20} -
5543x^{19} $\\*&${}- 22448x^{18}
 + 96508x^{17} + 1855180x^{16} $\\*&${}+ 13281488x^{15} + 66851616x^{14} 
+ 282546237x^{13} $\\*&${}+ 1087723107x^{12} + 3479009049x^{11} 
+ 8319918708x^{10} $\\*&${}+ 8576048755x^9 - 19169464149x^8 
- 111605931055x^7 $\\*&${}- 227855922888x^6 - 193255204370x^5 
$\\*&${}+ 176888550627x^4 + 1139040818642x^3 $\\*&${}+ 1055509532423x^2 
 + 1500432519809x $\\*&${}+ 314072259618$\cr 
\hline
$(22,23)$&$x^{24} - 2x^{23} + 115x^{22} + 23x^{21} + 1909x^{20} +
22218x^{19} $\\*&${}+ 9223x^{18}
 + 121141x^{17} + 1837654x^{16} - 800032x^{15} $\\*&${}+ 9856374x^{14}
 + 52362168x^{13} - 32040725x^{12} $\\*&${}+ 279370098x^{11}
 + 1464085056x^{10} + 1129229689x^9 $\\*&${}+ 3299556862x^8  
 + 14586202192x^7 + 29414918270x^6 $\\*&${}+ 45332850431x^5 
 - 6437110763x^4 - 111429920358x^3 $\\*&${}- 12449542097x^2 
 + 93960798341x - 31890957224$\cr
\end{longtable}
%\end{center}

%% The main "blocks" for the algorithms and for bounding their complexity

%% -- Decription X1(5l) : notation+important technicalities.
\chapter{Description of $X_1(\protect\lowercase{5l})$}\label{chap_descr_X15l}

\author{B. Edixhoven}

\bigskip

\bigskip

%author: Bas.

%% maybe improve the title of the chapter.... 
%% and rearrange the contents of the sections that follow. Maybe move
%% the explicit description of Y_1(5l) over \ZZ[1/5l] to the chapter
%% where modular curves are introduced...

\section{Construction of a suitable cuspidal divisor 
on~$X_1(\protect\lowercase{5l})$}
\label{sec_constr_D}
In this section we put ourselves in the situation of
Theorem~\ref{thm_red_to_wt_2_Delta}: $l$ is a prime number, $k$ is an
integer such that $2<k\leq l{+}1$, and $f$ is a surjective ring
morphism $\TT(1,k)\to\FF$ a with $\FF$ a finite field of
characteristic~$l$, such that the associated Galois representation
$\rho\colon\Gal(\Qbar/\QQ)\to\GL_2(\FF)$ is irreducible. We let $V$
denote the two-dimensional $\FF$-vector space in $J_1(l)(\Qbar)[l]$
that realises~$\rho$.

As explained in Chapter~\ref{chap_first_descr}, we would like to have an
effective divisor $D_0$ on $X_1(l)_\QQ$ of degree the genus of
$X_1(l)$ such that for all non-zero $x$ in the submodule $V$ of
$J_1(l)(\Qbar)$ we have $h^0(\calL_x(D_0))=1$. It would be nice to
have a cuspidal divisor (i.e., a divisor supported on the cusps) with
this property. The first complication is that the cusps are not all
rational over~$\QQ$: half of them have the maximal real subfield
of~$\QQ(\zeta_l)$ as field of definition. Moreover, even working with
all the cusps, we have not succeeded to find a cuspidal divisor $D_0$
with the desired properties. On the other hand, below we will give
explicitly a cuspidal divisor $D_0$ on the
curve~$X_1(5l)_{\QQ(\zeta_l)}$ that has the property that
$h^0(\calL_x(D_0))=1$ for each $x$ in $J_1(5l)(\Qbar)$ that
specialises to $0$ at some place of $\Qbar$ over~$l$. In particular,
$D_0$ has the required property for $V$ embedded in $J_1(5l)(\Qbar)$
in an arbitrary way, provided that the image of $\rho$
contains~$\SL_2(\FF)$ (this will be shown in
Section~\ref{sec_setup_tau}). We have chosen to work with $X_1(5l)$,
but the same method will work for modular curves corresponding to some
level structure if the prime to $l$ part of the level structure is
fine, and of genus zero.

For the rest of this section, our assumptions are the following: $l$
is a prime number, not equal to~$5$. We let
$X:=X_1(5l)_{\QQ(\zeta_l)}$ over~$\QQ(\zeta_l)$. The genus of~$X$ is
$(l-2)^2$. References for facts about $X$ that we use can be found
in~\cite{Gross1}, and also in~\cite{Edixhoven3}; they are derived from
results in~\cite{Deligne-Rapoport} and in~\cite{Katz-Mazur}.

The curve $X_0(5)_\QQ$ has $2$ cusps, both $\QQ$-rational, called $0$
and~$\infty$ (after the points of $\PP^1(\QQ)$ of which they
come). The cusp $\infty$ has as moduli interpretation the degenerate
elliptic curve (generalised elliptic curve in the terminology
of~\cite{Deligne-Rapoport}): the $1$-gon, equipped with the unique
subgroup of order~$5$ of~$\Gm$. The cusp~$0$ corresponds to the
$5$-gon, equipped with a subgroup of order~$5$ that meets all $5$
components. The group $\FF_5^\times$ acts (as diamond operators) on
$X_1(5)$, with quotient~$X_0(5)$; the subgroup $\{\ld\pm 1\rd\}$ acts
trivially, and the quotient $\FF_5^\times/\{\ld\pm 1\rd\}$ by this
subgroup acts faithfully. The inverse images of $0$ and $\infty$ both
consist of two cusps. Those over $0$ are $\QQ$-rational (the subgroup
of order~$5$ of the $5$-gon is the constant groupscheme~$\ZZ/5\ZZ$),
whereas those over~$\infty$ are conjugated over~$\QQ(\sqrt{5})$. We
fix one $\QQ$-rational cusp $c$ of~$X_1(5)$.

The group $\FF_l^\times$ acts faithfully, and in fact, even freely,
on~$X$. The set of cusps of $X$ over the cusp $c$ of $X_1(5)$ form two
$\FF_l^\times$-orbits, corresponding to the type of degenerate
elliptic curve that they correspond to: $5$-gon or $5l$-gon. The orbit
corresponding to the $5$-gon consists of points rational over
$\QQ(\zeta_l)$, all conjugates of each other. The orbit corresponding
to the $5l$-gon consists of $\QQ$-points.

We let $J$ denote the Jacobian of~$X$. What we want is an effective
divisor $D_0$ of degree $g$ on $X$ (with $g$ the genus of~$X$),
supported on the cusps over~$c$, such that for all $x$ in $J(\Qbar)$
that specialise to $0$ at some place of $\Qbar$ over~$l$ we have
$h^0(X_\Qbar,\calL_x(D_0))=1$. For the notion of specialisation we use
Néron models; the reader is referred to~\cite{BLR1} for this
notion. For $x$ in $J(K)$ with $K\subset\Qbar$ a finite extension of
$\QQ(\zeta_l)$, and $\lambda$ a place of $\Qbar$ over~$l$, we say that
$x$ specialises to~$0$ at~$\lambda$ if $x$, viewed as element of
$J_{O_K}(O_K)$, with $J_{O_K}$ the Néron model of $J$ over~$O_K$,
specialises to $0$ at the place of $K$ given by~$\lambda$. For
$K\subset K'$ a finite extension we have $J(K)\subset J(K')$, hence we
can also view $x$ as element of~$J(K')$. The notion of $x$
specialising to zero is the same for $K$ and $K'$, because
$J_{\ZZ[\zeta_l]}$ is semi-stable at~$l$.

The moduli interpretation of~$X$ gives a semi-stable model
$X_{\ZZ[\zeta_l,1/5]}$ over $\ZZ[\zeta_l,1/5]$, described
in~\cite{Gross1} for example. A result of Raynaud identifies the
connected component of the Néron model $J_{O_K}$ with the connected
component of the Picard scheme of~$X_{O_K}$ (see Section~9.5
of~\cite{BLR1}).  This means that for $x$ in $J(K)$ specialising to
$0$ at $\lambda$ the line bundle $\calL_x$ on $X_K$ associated with $x$
can be extended uniquely over the local ring $O_{K,\lambda}$ to a line
bundle $\calL_x$ on $X_{O_{K,\lambda}}$ such that the restriction
$\ol{\calL_x}$ of $\calL_x$ to the special fibre $X_{\FF_\lambda}$ is
trivial. The divisor $D_0$ on $X_K$ extends, by taking the Zariski
closure, to an effective Cartier divisor on~$X_{O_{K,\lambda}}$.

We now note that $h^0(X_K,\calL_x(D_0))$ is at least one, by
Riemann-Roch, and that $h^0(X_K,\calL_x(D_0))$ is at most
$h^0(X_{\FF_\lambda},\ol{\calL_x}(\ol{D_0}))$ by semi-continuity of
cohomology of coherent sheaves.  As
$\ol{\calL_x}(\ol{D_0}))=\calO(\ol{D_0})$, it now suffices to take $D_0$
such that $h^0(X_{\Fbar_l},\calO(\ol{D_0}))=1$. We do this by looking at
the geometry of~$X_{\Fbar_l}$. As $\ZZ[\zeta_l]$ has a unique morphism
to~$\Fbar_l$, the curve $X_{\Fbar_l}$ does not depend on~$K$. The
scheme of cusps of $X_{\ZZ_l[\zeta_l]}$ is finite étale over
$\ZZ_l[\zeta_l]$, hence the cusps lying over $c$ specialise
injectively to~$X_{\FF_l}$.

The curve $X_{\Fbar_l}$ is the union of two irreducible components,
$X_1$ and $X_2$, say, both isomorphic to the Igusa curve of level $l$
over $X_1(5)$ over~$\Fbar_l$, that meet transversally in the set
$\Sigma$ of supersingular points. We will take $\ol{D_0}=D_1+D_2$, with
$D_1$ on $X_1$ and $D_2$ on~$X_2$; note that the cusps are disjoint
from~$\Sigma$, so $\ol{D_0}$ lies in the smooth locus
of~$X_{\Fbar_l}$. 

In order to simplify the notation, we let $X$ and $D_0$
denote~$X_{\Fbar_l}$ and~$\ol{D_0}$, from now on, in this section. We
let $\Omega_X$ be the dualising sheaf on~$X$ (see Section~8
of~\cite{Gross1}, or~\cite{Mazur-Ribet1}): it is the invertible
$\calO_X$-module obtained by gluing $\Omega^1_{X_1}(\Sigma)$ and
$\Omega^1_{X_2}(\Sigma)$ along~$\Sigma$ via the residue maps at the
points of $\Sigma$ on~$X_1$ and minus the residues maps at the points
of $\Sigma$ on~$X_2$.  By Riemann-Roch, what we want is that
$h^1(X,\calO(D_0))=0$, and hence, by Serre duality, that
$h^0(X,\Omega_X(-D_0))=0$. In other words, an element of
$\rH^0(X,\Omega_X)$ that vanishes on $D_0$ must be zero. Restriction
to~$X_1$ gives a short exact sequence:
\begin{multline*}
0 \to \rH^0(X_2,\Omega^1_{X_2/\Fbar_l}(-D_2)) \to \rH^0(X,\Omega_X(-D_0)) 
\to \\
\to \rH^0(X_1,\Omega^1_{X_1/\Fbar_l}(\Sigma-D_1)) \to 0.
\end{multline*}
Hence, it suffices to take $D_1$ such that
$\rH^0(X_1,\Omega^1_{X_1/\Fbar_l}(\Sigma-D_1))=0$ and $D_2$ such that
$\rH^0(X_2,\Omega^1_{X_2/\Fbar_l}(-D_2))=0$. Let us now first see of
which degrees $d_1$ and $d_2$ we want to take $D_1$ and~$D_2$. Let
$g_1$ and $g_2$ be the genera of $X_1$ and $X_2$ (note: they are
equal). Then we have that $g=g_1+g_2+\#\Sigma -1$, and $g=\deg D_0 =
d_1+d_2$. It can be shown in several ways that the degree of the sheaf
$\ul{\omega}$ on $X_1(5)_{\FF_l}$ is one. Either by explicit
computation, using the equations of Proposition~\ref{prop_Y_1_5}, or
by the following argument. The curve $X_1(5)_{\FF_l}$ over~$\FF_l$ has
genus zero. The Kodaira-Spencer isomorphism on~$X_1(5)_{\FF_l}$, from
$\ul{\omega}^{\otimes 2}\to\Omega^1(\Cusps)$, plus the fact that the
divisor of cusps has degree~$4$, give that the degree of
$\ul{\omega}^{\otimes 2}$ is~$2$.  Therefore, the Hasse invariant,
%%% give reference, to Katz-Mazur for example!!!!!!!!!!!!!!
being a global section of~$\ul{\omega}^{\otimes l{-}1}$, has exactly
$l{-}1$ zeros on $X_1(5)$ over $\Fbar_l$ and therefore we have:
\begin{eqn}
\#\Sigma = l-1.
\end{eqn}
Applying Hurwitz's formula to the covering $X_1\to X_1(5)_{\Fbar_l}$,
which is totally ramified over $\Sigma$ and unramified outside it,
gives:
\[
2g_1-2 = -2(l-1) + (l-1)(l-2),
\]
and hence:
\begin{eqn}
g_1 = \frac{1}{2}(l-2)(l-3). 
\end{eqn}
This implies that we want to take:
\begin{eqn}
\begin{aligned}
d_1 & = g_1 + \#\Sigma -1 = \frac{1}{2}(l-1)(l-2), \\
d_2 & = g_2 = \frac{1}{2}(l-2)(l-3).
\end{aligned}
\end{eqn}
Now we use equations to compute with. We choose a coordinate $z$ on
$X_1(5)_{\Fbar_l}$, i.e., an isomorphism from $X_1(5)_{\Fbar_l}$
to~$\PP^1_{\Fbar_l}$, such that $z(\Sigma)$ does not contain $0$ or
$\infty$ and such that $z^{-1}0$ is the (rational) cusp $0$
of~$X_1(5)_{\Fbar_l}$. Let $f$ be the monic polynomial in $z$ whose
zeros are the elements of~$\Sigma$, each with multiplicity one. Then
$X_1$ and $X_2$ are both isomorphic to the cover of $X_1(5)_{\Fbar_l}$
given by the equation $y^{l-1}=f$, by the following argument.  The
complete local rings of $X_1$ at the points of~$\Sigma$, with their
$\FF_l^\times$-actions, are all isomorphic to each other because, by a
theorem of Serre and Tate, these can all be described in terms of the
deformation theory of one $l$-divisible group over~$\Fbar_l$;
see~\cite[\S5.2--5.3]{Katz-Mazur}. There is a general theory of cyclic
possibly ramified covers such as $X_1\to X_1(5)_{\Fbar_l}$, based on
the decomposition of $\calO_{X_1}$ as $\calO$-module on
$X_1(5)_{\Fbar_l}$ for the $\FF_l^\times$-action.
%%% Include a reference. Liu's book???
It shows that the cover $X_1\to X_1(5)_{\Fbar_l}$ is the cover of
$l{-}1$th roots of the global section~$1$ of the invertible sheaf
$\calO(\Sigma)$ on $X_1(5)_{\Fbar_l}$, in an invertible $\calO$-module
$\calL$ with a given isomorphism
$\calO(\Sigma)\to\calL^{\otimes(l-1)}$, where the $\FF_l^\times$-action
may have been changed by an automorphism of~$\FF_l^\times$. As
$X_1(5)_{\Fbar_l}$ has genus zero, we can take $\calL$ to
be~$\calO_{X_1}(z^{-1}\infty)$. In fact, Section~12.8
of~\cite{Katz-Mazur} shows that $X_1$ is obtained from
$X_1(5)_{\Fbar_l}$ by extracting the $l{-}1$th root of the Hasse
invariant, in~$\ul{\omega}$.

We compute a basis of $\rH^0(X_1,\Omega^1_{X_1/\Fbar_l}(\Sigma))$. On
$X_1$ we have:
\begin{eqn}
-y^{l-2}dy = f'dz.
\end{eqn}
Hence $(dz)/y^{l-2}=-(dy)/f'$ is a generating section of
$\Omega^1_{X_1/\Fbar_l}$ on the affine part given by our
equation. Hence $(dz)/y^{l-1}$ is generating section of
$\Omega^1_{X_1/\Fbar_l}(\Sigma)$ on the affine part, and it is
$\FF_l^\times$-invariant. At each point of $X_1$ over the point where $z$
has its pole, both $z$ and $y$ have a simple pole, and $(dz)/y^{l-1}$
has order $-2+l-1 = l-3$. So we have a basis:
\begin{eqn}\label{eqn_basis_diff_forms_1}
\rH^0(X_1,\Omega^1_{X_1/\Fbar_l}(\Sigma)) = 
\bigoplus_{i+j\leq l-3}\Fbar_l z^iy^j \cdot (dz)/y^{l-1}.
\end{eqn}
Note that this agrees with the fact that $d_1=\frac{1}{2}(l-1)(l-2)$. 

We can now say how to choose $D_1$. At each of the $l{-}1$ points
where $z$ has a zero we must give a multiplicity. In the coordinate
system given by $z$ and~$y$, these points are the ones of the form
$(0,b)$ with $b\in \Fbar_l^\times$ satisfying $b^{l-1}=f(0)$. Here is
how we choose~$D_1$: just distribute the multiplicities
$(0,1,\ldots,l-2)$ over these points. Then one sees that any linear
combination of our basis elements that vanishes on $D_1$ is zero, as
follows. At all points of~$D_1$, $z$ has a simple zero. Let $\omega$
be an element of $\rH^0(X_1,\Omega^1_{X_1/\Fbar_l}(\Sigma))$, with
coordinates $\lambda_{i,j}$ in the
basis~(\ref{eqn_basis_diff_forms_1}). Assume that $\omega$ vanishes
on~$D_1$ (taking multiplicities into account). As there are $l{-}2$
points in $D_1$ with multiplicity~$>0$, the polynomial
$\sum_j\lambda_{0,j}y^j$, being of degree $\leq l{-}3$, must be
zero. As there are $l{-}3$ points in $D_1$ with multiplicity~$>1$, the
polynomial $\sum_j\lambda_{1,j}zy^j$, being of degree $\leq l{-}4$,
must be zero. And so on.

Now we do~$D_2$. A basis is the following:
\begin{eqn}
\rH^0(X_2,\Omega^1_{X_2/\Fbar_l}) = 
\bigoplus_{i+j\leq l-4} \Fbar_l z^iy^j \cdot (dz)/y^{l-2}.
\end{eqn}
We note that this agrees with $g_2=(l-3)(l-2)/2$. So, for $D_2$, just
distribute the multiplicities $(0,0,1,\ldots,l-3)$ over the
points where $z$ has a zero. The same argument as the one we used for
$D_1$ shows that any $\omega$ in $\rH^0(X_2,\Omega^1_{X_2/\Fbar_l})$
that vanishes on~$D_2$ is zero.

We summarise our results. As the action of $\FF_5^\times$ permutes the
two $\QQ$-rational cusps of~$X_1(5)$, our arguments above work for
both of them.
\begin{thm}\label{thm_exist_D}
Let $l$ be a prime number not equal to~$5$. Let $c$ be one of the two
$\QQ$-rational cusps of~$X_1(5)$. Then the cusps of $X_1(5l)$ over $c$
are $\QQ(\zeta_l)$-rational, and consist of two $\FF_l^\times$-orbits,
on which $\FF_l^\times$ acts freely. Let $D_1$ be a divisor on
$X_1(5l)_{\QQ(\zeta_l)}$ obtained by distributing the multiplicities
$(0,1,\ldots,l-2)$ over one of these two orbits. Let $D_2$ be the
divisor obtained by distributing the multiplicities
$(0,0,1,\ldots,l-3)$ over the other orbit. Then $D_0:=D_1+D_2$ has
degree equal to the genus of $X_1(5l)$ and has the property that for
any $\Qbar$-point $x$ of the Jacobian of $X_1(5l)$ that specialises
to~$0$ at some place over~$l$ we have
$h^0(X_1(5l)_\Qbar,\calL_x(D_0))=1$.
\end{thm}

\section{The exact setup for the level one case}
\label{sec_setup_tau}

In Chapter~\ref{chap_first_descr} we described our strategy for computing
the residual Galois representations attached to a fixed
newform. That strategy depends on properties of divisors $D_0$ and
functions~$f$ to be chosen, on modular curves of varying level. These
$D_0$ and $f$ must satisfy a number of conditions. In general we do
not know how to choose divisors~$D_0$ of which we can prove, without a
computer computation, that they have the required property. This is
the main reason for which we will now restrict ourselves to just the
case of modular forms of level one. 

The aim of this section is to describe exactly our strategy for
computing the residual representations $V$ in the situation of
Theorem~\ref{thm_red_to_wt_2_Delta}: $l$ is a prime number, $k$ is an
integer such that $2<k\leq l{+}1$, and $f$ is a surjective ring
morphism $\TT(1,k)\to\FF$ a with $\FF$ a finite field of
characteristic~$l$, such that the associated Galois representation
$\rho\colon\Gal(\Qbar/\QQ)\to\GL_2(\FF)$ is irreducible, under the
extra hypothesis that the image of $\rho$ contains~$\SL_2(\FF)$. By
Theorem~\ref{thm_large_image}, this hypothesis holds when $\rho$ is
irreducible and $l\geq 6k-5$. We let $V$ denote the two-dimensional
$\FF$-vector space in $J_1(l)(\Qbar)[l]$ that realises~$\rho$.

Theorem~\ref{thm_exist_D} gives us a divisor~$D_0$
on~$X_1(5l)_{\QQ(\zeta_l)}$ that we want to use. Therefore, we want to
embed $V$ into~$J_1(5l)(\Qbar)[l]$.

Let $\pi\colon X_1(5l)\to X_1(l)$ be the standard map (i.e., the one
that forgets the $5$-part of the level structure, the one denoted
$B_{5l,l,1}$ in Section~\ref{sec_modforms}). Then the degree of $\pi$
is $5^2-1=24$, which is prime to~$l$. This implies that $\pi_*\pi^*$
is multiplication by $24$ on~$J_1(l)$, and that $\pi^*$ is injective
on $J_1(l)(\Qbar)[l]$. We have a projector:
\begin{eqn}
\frac{1}{24}\pi^*\pi_* \colon J_1(5l)(\Qbar)[l] 
\onto \pi^* J_1(l)(\Qbar)[l] \subset J_1(5l)(\Qbar)[l].
\end{eqn}
We will consider $V$ embedded in $J_1(5l)(\Qbar)[l]$ via its embedding
into $J_1(l)(\Qbar)[l]$, followed by~$\pi^*$. 
\begin{prop}\label{prop_uniqueD'}
Let $l$ be a prime number, let $k$ be an integer such that $2<k\leq
l{+}1$, and $f$ a surjective ring morphism $\TT(1,k)\to\FF$ a with
$\FF$ a finite field of characteristic~$l$, such that the image of the
associated Galois representation
$\rho\colon\Gal(\Qbar/\QQ)\to\GL_2(\FF)$ contains~$\SL_2(\FF)$. We let
$V$ denote the pullback as above of the two-dimensional $\FF$-vector
space in $J_1(l)(\Qbar)[l]$ that realises~$\rho$.  Let $D_0$ be a
divisor on $X_1(5l)_{\QQ(\zeta_l)}$ as given in
Theorem~\ref{thm_exist_D}. Then, for every $x$ in~$V_l$, we have
$h^0(X_1(5l)_\Qbar,\calL_x(D_0))=1$.
\end{prop}
\begin{proof}
In view of Theorem~\ref{thm_exist_D}, it suffices to show that for
each $x$ in $V$ there is a place of $\Qbar$ over~$l$ at which $x$
specialises to~$0$. The notion of specialisation is explained in
Section~\ref{sec_constr_D}. Let $J_{\ZZ[\zeta_l]}$ denote the Néron
model of $J:=J_1(5l)$ over~$\ZZ[\zeta_l]$. Then $V$ is the group of
$\Qbar$-points of an $\FF$-vector space scheme $\calV_{\QQ(\zeta_l)}$
in $J_{\QQ(\zeta_l)}$. Let $\calV$ be the Zariski closure of
$\calV_{\QQ(\zeta_l)}$ in~$J_{\ZZ[\zeta_l]}$. Then it is shown in
Section~12 of ~\cite{Gross1} and in Section~6 of~\cite{Edixhoven3}
that $\calV_{\ZZ_l[\zeta_l]}$ is finite locally free over
$\ZZ_l[\zeta_l]$, and that the dimension as $\FF$-vector space scheme
of the local part of $\calV_{\ZZ_l[\zeta_l]}$ is $1$ if
$f(T_l)\neq 0$ and $2$ if~$f(T_l)=0$. We note that it does not
matter if we take Zariski closure in~$J_1(l)$ or in~$J_1(5l)$, as
$\pi^*$ gives a closed immersion of the $l$-torsion of $J_1(l)$ over
$\ZZ_l$ into that of~$J_1(5l)$.

This means that at each place of~$\Qbar$ over~$l$ there is a non-zero
$x$ in $V$ that specialises to~$0$. Under our assumptions, the image
of $\Gal(\Qbar/\QQ(\zeta_l))$ acting on~$V$ is~$\SL(V)$. Hence
$\Gal(\Qbar/\QQ(\zeta_l))$ acts transitively on~$V-\{0\}$. Hence for
each $x$ in~$V-\{0\}$ there is at least one place of~$\Qbar$ over~$l$
where $x$ specialises to~$0$.
\end{proof}
The fact that our divisor $D_0$ lives on~$X_1(5l)_{\QQ(\zeta_l)}$, and
not on~$X_1(5l)_\QQ$, forces us to work over $\QQ(\zeta_l)$, and not
over~$\QQ$, as in Chapter~\ref{chap_first_descr}. 

We let $X_l$ denote~$X_1(5l)_{\QQ}$ and $g_l$ its genus, and we let
$A_{l,\QQ(\zeta_l)}$ denote the $\QQ(\zeta_l)$-algebra that
corresponds to the $\Gal(\Qbar/\QQ(\zeta_l))$-set~$V$. In order to
explain the notation $A_{l,\QQ(\zeta_l)}$, we note that this
$\QQ(\zeta_l)$-algebra is obtained from the $\QQ$-algebra $A_l$ (that
corresponds to the $\Gal(\Qbar/\QQ)$-set $V$) by extension of
scalars.

Proposition~\ref{prop_uniqueD'} gives that for each $x$ in $V$ there
is a unique effective divisor $D_x=\sum_{i=1}^{g_l} Q_{x,i}$ of
degree~$g_l$ on $X_{l,\Qbar}$ such that $x=[D_x-D_0]$
in~$J_l(\Qbar)$. Note that, for $x=0$, it is indeed true that
$D_x=D_0$, hence the two notations are consistent. The following
observation will make the exposition in
Chapter~\ref{sec_couveignes_TORSION} somewhat easier. As for each $x$
in $V$ specialises to $0$ at some place of $\Qbar$ above~$l$, the
divisor $D_x$ specialises to the cuspidal divisor $D_0$ at such a
place, and hence none of all $Q_{x,i}$ can be a CM-point, in
particular:
\begin{eqn}\label{eqn_jnotramifiedatQxi}
\text{for all $x$ and $i$:}\quad j(Q_{x,i})\not\in \{0,1728\}.
\end{eqn}
The uniqueness of $D_x$ implies that:
\begin{eqn}
D_{gx}=gD_x, \quad
\text{for all $x$ in $V$ and $g$ in $\Gal(\Qbar/\QQ(\zeta_l))$}.
\end{eqn}
We write each $D_x$ as:
\begin{eqn}
D_x = D_x^\fin + D_x^\cusp,
\end{eqn}
where $D_x^\cusp$ is supported on the cusps of $X_{l,\Qbar}$ and where
$D_x^\fin$ is disjoint from the cusps. The next lemma shows that
$D_x^\fin$ determines~$x$, and its proof uses only that $\rho$ is
absolutely irreducible, not that its image contains~$\SL_2(\FF)$.
\begin{lem}\label{lem_x_from_Dx_fin}
In this situation, the map from $V$ to the set of effective divisors
on $X_{l,\Qbar}$ that sends $x$ to $D_x^\fin$ is injective.
\end{lem}
\begin{proof}
Suppose that it is not. We take $x_1$ and $x_2$ in~$V$, distinct, such
that $D_{x_1}^\fin = D_{x_2}^\fin$.  Then the element $x_1-x_2$ in $V$
is nonzero and is represented by the cuspidal divisor
$D_{x_1}-D_{x_2}$. The cusps of $X_l$ are rational over
$\QQ(\zeta_{5l})$. Hence $x_1-x_2$ gives an injection $\FF\to V$ of
representations of $\Gal(\Qbar/\QQ(\zeta_{5l}))$, where $\FF$ has
trivial action. But that gives, by adjunction of induction and
restriction, a nonzero map from the regular representation of
$\Gal(\QQ(\zeta_{5l})/\QQ)$ over $\FF$ to~$V$, necessarily surjective
because $V$ is irreducible. But then the image of
$\rho\colon\Gal(\Qbar/\QQ)\to\GL_2(\FF)$ is abelian. As $\rho$ is odd,
the two eigenspaces in~$V$ of any complex conjugation then
decompose~$V$, in contradiction with the irreducibility of~$V$.
\end{proof}
As the cusps of $X_{l,\Qbar}$ form a $\Gal(\Qbar/\QQ)$-stable subset
of~$X_{l,\Qbar}$ we have:
\begin{eqn}
\text{for all $g$ in $\Gal(\Qbar/\QQ(\zeta_l))$:}\quad
D_{gx}^\fin = gD_x^\fin, \quad D_{gx}^\cusp = gD_x^\cusp.
\end{eqn}
Hence the map that sends $x$ in~$V$ to $D_x^\fin$ is
$\Gal(\Qbar/\QQ(\zeta_l))$-equivariant. 

We will now produce a suitable function $f_l\colon X_l\to\PP^1_\QQ$,
in order to push the set of $\{D_x^\fin\;|\;x\in V\}$ injectively and
$\Gal(\Qbar/\QQ(\zeta_l))$-equivariantly to the set
$\{f_{l,*}D_x^\fin\;|\;x\in V\}$ of divisors on~$\AA^1_\Qbar$.

We start by giving an explicit description of the curve $Y_1(5)$
over~$\ZZ[1/5]$. In order to do that, we determine a universal triple
$(E/S,P)$ where $E/S$ is an elliptic curve over an arbitrary scheme,
and $P$ in $E(S)$ is everywhere of order~$5$, i.e., for every
$\Spec(A)\to S$ with $A$ non-zero, the image of $P$ in~$E(A)$ has
order~$5$. The base of this triple is the open part $Y_1(5)'$ of the
model $Y_1(5)$ over~$\ZZ$ (constructed in Chapter~8
of~\cite{Katz-Mazur}) where the $P$ has order~$5$ (i.e., $Y_1(5)'$ is
the complement of the irreducible component of $Y_1(5)_{\FF_5}$ where
the point~$P$ generates the kernel of Frobenius). The equation of this
universal triple can also be found on page~7 of Tom Fisher's thesis,
see~\cite{Fisher1}. 

\begin{prop}\label{prop_Y_1_5}
Let $E/S$ be an elliptic curve, and $P\in E(S)$ a point that is
everywhere of order~$5$. Then $(E/S,P)$ arises via a unique base
change from the following triple:
\[
\left\{
\begin{aligned}
&E: y^2 + (b+1)xy + by = x^3 + bx^2\\
&Y_1(5)' =  \Spec(\ZZ[b,1/\discr(E)]), \quad \discr(E) =
-b^5(b^2+11b-1)\\
&P = (0,0).
\end{aligned}
\right.
\]
The $j$-invariant of $E$ is given by:
\[
j(E) =-(b^4 + 12b^3 + 14b^2 - 12b + 1)^3/b^5(b^2+11b-1).
\]
\end{prop}
\begin{proof}
Our proof is modelled on Section~2.2 of~\cite{Katz-Mazur}; basic
properties of Weierstrass equations for elliptic curves are used
without being mentioned.

Let $(E/S,P)$ be given, with $P$ everywhere of order~$5$. Choose a
parameter~$t$ at~$0$, up to order~$2$, i.e., a trivialisation
of~$\omega_{E/S}$. Note: we are working locally on~$S$, here; in the
end, as we will succeed in making things unique, our construction will
be global. Note: $t$ is unique up to $t'=ut$, with $u\in R^\times$
($S=\Spec(R)$ now).

Choose $x$ a global function on $E-0(S)$ such that
$x=t^{-2}+\cdots$. Then $x$ is unique up to $x'=x+a$, $a\in R$. Make
$x$ unique by demanding that $x(P)=0$ (this is alright because $0$
and~$P$ are disjoint).

Choose $y=t^{-3}+\cdots$ regular on~$E-0(S)$. Then $y$ is unique up to
$y'=y+ax+b$. Make $y$ unique by demanding that $y(P)=0$ and that the
tangent of~$E$ at~$P$ is the line given by the
equation~$y=0$. (Indeed, use $b$ (uniquely) to get $y(P)=0$, then note
that the tangent at~$P$ is nowhere the line given by $x=0$ because $P$
is nowhere annihilated by~$2$).

The equation for $E$ is of the form:
\[
y^2 + a_1 xy + a_3 y = x^3 + a_2 x^2,
\]
because the coefficients usually called $a_4$ and~$a_6$ are zero. We
also see that $a_3$ is a unit because $E$ is smooth at~$(0,0)$. The
coefficient~$a_2$ is a unit, because $P$ is nowhere annihilated
by~$3$.

Now we try to get rid of $u$ (the ambiguity in the choice of~$t$). If
$t'=u^{-1}t$, then $a_i'=t^ia_i$, hence we can make $t$ unique by
demanding that $a_2=a_3$. We do that, and then we have the following
equations.

The elliptic curve $E$ and the point $P$ are given by:
\[
y^2 + a xy + b y = x^3 + b x^2,\quad P=(0,0).
\]

We have $5{\cdot}P\equiv 5{\cdot}0$, hence there is a unique $f$ on
$E-0(S)$ of the form:
\[
f=xy+\alpha y+\beta x^2+\gamma x+\delta
\]
such that the divisor of~$f$ is $5{\cdot}P-5{\cdot}0$. As $x$ and $y$
have order one and two, respectively, at~$P$, we have
$\gamma=\delta=0$. The function $f$ with divisor $5{\cdot}P-5{\cdot}0$
is given by:
\[
f = xy + \alpha y + \beta x^2.
\]
Here we know that $b$, $\alpha$ and~$\beta$ are in~$R^\times$, because
$v_P(x)=1$ and $v_P(y)=2$ everywhere on~$S$. Now we have to compute
what it means that $v_P(f)=5$. This means that the intersection
multiplicity of the two curves $E$ and~$V(f)$ at~$(0,0)$ is~$5$. A
systematic way to compute that is to do successive blow-ups; that
works nicely, but we will not do that here. A much faster way to do
the computation is to take suitable linear combinations of the
equations for~$E$ and $f$ directly. One finds the equations:
\[
\beta = -\alpha, \quad a = \alpha^{-1}b+1,\quad
\alpha=1. 
\]
\end{proof}
The reason we give such a detailed description of $Y_1(5)$ is that it
gives us functions on all the $Y_1(5l)$, at least over $\ZZ[1/5l]$, as
stated in the following proposition.

\begin{prop}\label{prop_bxy}
Let $l\neq 5$ be prime. Let $E/Y_1(5)'$ be the elliptic curve given in
Proposition~\ref{prop_Y_1_5}. Then $Y_1(5l)$ and $E[l]-\{0\}$ agree
over~$\ZZ[1/5l]$: for $S$ a $\ZZ[1/5l]$-scheme and $Q$ in
$(E[l]-\{0\})(S)$ we get, by pullback from $Y_1(5)'$, an elliptic curve
over $S$ with an $S$-valued point $P_5$ that is everywhere of
order~$5$, and an $S$-valued point $P_l$ that is everywhere of
order~$l$. In particular, the functions $b$, $x$ and~$y$ on
$E[l]-\{0\}$ give functions $b_l$, $x_l$ and $y_l$ on $Y_1(5l)$
over~$\ZZ[1/5l]$ that generate its coordinate ring.
\end{prop}
\begin{proof}
This is standard. The construction above gives a morphism, over
$\ZZ[1/5l]$, from $E[l]-\{0\}$ to~$Y_1(5l)$. Conversely, an elliptic
curve over $S$ with such points $P_5$ and $P_l$ gives a point of
$E[l]-\{0\}$ by the universality of $Y_1(5)'$ and the fact that
$E[l]$ is finite étale over $Y_1(5)'$ away from
characteristic~$l$. The second statement follows from the fact that
$E[l]-\{0\}$ is a closed subscheme of the affine scheme $E-\{0\}$.
\end{proof}
The functions $b_l$, $x_l$ and $y_l$ on $Y_1(5l)$ over~$\ZZ[1/5l]$
have the following moduli interpretations. For any $\ZZ[1/5l]$-algebra $A$,
and any $Q$ in $Y_1(5l)(A)$, a point corresponding to a triple
$(E,P_5,P_l)$ with $E$ an elliptic curve over~$A$, $P_5$ in $E(A)$ a
point that is everywhere of order~$5$ and $P_l$ in $E(A)$ a point that
is everywhere of order~$l$, there are unique elements $b_l(Q)$,
$x_l(Q)$ and $y_l(Q)$ in~$A$ such that $(E/A,P_5)$ is uniquely
isomorphic to the pair given by:
\[
y^2+(b_l(Q)+1)xy +b_l(Q)y = x^3+b_l(Q)x^2, \quad P_5=(0,0).
\]
Then, in these coordinates, we have: 
\[
P_l = (x_l(Q),y_l(Q)).
\]
Similarly, we define regular functions $x'_l$ and $y'_l$ on
$Y_1(5l)_{\ZZ[1/5l]}$ by the condition that, in the coordinates above,
we have:
\[
l^{-1}P_5+P_l = (x'_l(Q),y'_l(Q)),
\]
where $l^{-1}P_5$ is the unique point $Q$ in $E[5](A)$ with $lQ=P_5$.
 
We note that the pair of functions $(b_l,x'_l)$ embeds
$Y_1(5l)_{\ZZ[1/5l]}$ in the affine plane~$\AA^2_{\ZZ[1/5l]}$. Indeed,
assume that $Q$ and $Q'$ are in $Y_1(5l)(A)$, corresponding to
$(E,P_5,P_l)$ and $(E',P_5',P_l')$, with $b_l(Q)=b_l(Q')$ and
$x'_l(Q)=x'_l(Q')$. Then, by Proposition~\ref{prop_Y_1_5}, $(E,P_5)$
is uniquely isomorphic to $(E',P_5')$, and so we simply consider them
to be equal. Then, $l^{-1}P_5+P_l$ and $l^{-1}P_5+P_l'$ have the same
$x$-coordinate. Hence, locally on $\Spec(A)$,
$l^{-1}P_5+P_l=\pm(l^{-1}P_5+P_l')$. Multiplying by $l$ we see that
the sign cannot be a minus.

Using the functions~$b_l$ and~$x'_l$, we can now say how we will
choose the function~$f_l$. We return to the situation right after
Lemma~\ref{lem_x_from_Dx_fin}. 

For $x\in V$, let $d_x$ be the degree of~$D_x^\fin$, and let us write
$D_x$ as a sum of points in~$X_l(\Qbar)$ as follows:
\[
D_x = \sum_{i=1}^g Q_{x,i}, \quad\text{with}\quad
D_x^\fin = \sum_{i=1}^{d_x} Q_{x,i}, \quad 
D_x^\cusp = \sum_{i=d_x+1}^g Q_{x,i}.
\]
We note that $d_0=0$, as $D_0$ is a cuspidal divisor, and that for all
non-zero $x$ in $V$ the $d_x$ are equal, as they are permuted
transitively by $\Gal(\Qbar/\QQ(\zeta_l))$.
%%%%%%%%say_that_our_construction_does_not_use_this???????????????????????

The set $S$ of points in~$\AA^2(\Qbar)$ consisting of the
$(b_l(Q_{x,i}),x'_l(Q_{x,i}))$, with $x$ in $V$ and $i$ in
$\{1,\ldots,d_x\}$ has at most $g{\cdot}(\#\FF)^2$ elements. We want
to project $S$ injectively into $\AA^1(\Qbar)$ with a map of the form
$(a,b)\mapsto a+nb$ for a suitable integer~$n$. As there are at most
$g^2{\cdot}(\#\FF)^4$ pairs of distinct elements in~$S$, at most that
number of integers~$n$ is excluded. Hence there exists an integer $n$
with $0\leq n\leq g^2{\cdot}(\#\FF)^4$ such that the function
$f_l:=b_l+nx'_l$ has the required property that the $f_{l,*}D_x^\fin$,
for $x\in V$, are all distinct.

Let $f_l$ be such a function. For each $x$ in~$V$, $f_{l,*}D_x^\fin$
gives us a polynomial $P_{D_0,f_l,x}$ with coefficients in $\Qbar$
given by:
\[
P_{D_0,f_l,x}(t) = \prod_{i=1}^{d_x} (t-f_l(Q_{x,i})) 
\quad \text{in $\Qbar[t]$.} 
\]
Vice versa, each $P_{D_0,f_l,x}$ gives us the divisor
$f_{l,*}D_x^\fin$ by taking the roots, with multiplicity. Therefore,
the map that sends $x$ to $P_{D_0,f_l,x}$ is injective, and
$\Gal(\Qbar/\QQ(\zeta_l))$-equivariant. 

The next step is to ``encode'' each $P_{D_0,f_l,x}$ in a single
element of~$\Qbar$, respecting the action of
$\Gal(\Qbar/\QQ(\zeta_l))$. We do this by evaluating at a suitable
integer~$m$, i.e., by sending $P_{D_0,f_l,x}$ to
$P_{D_0,f_l,x}(m)$. For a given $m$, this map is injective if and only
if for any distinct $x_1$ and $x_2$ in $V$, $m$ is not a root of the
difference of $P_{D_0,f_l,x_1}$ and $P_{D_0,f_l,x_2}$. Each of these
differences has at most $g$ roots, and as there are less than $(\#\FF)^4$
such differences, there are at most $g{\cdot}(\#\FF)^4$ integers to
avoid. So there is a suitable $m$ with $0\leq m\leq
g{\cdot}(\#\FF)^4$. Composing our maps, we obtain a generator for the
$\QQ(\zeta_l)$-algebra $A_{l,\QQ(\zeta_l)}$ associated with~$V$:
\[
a_{D_0,f_l,m}\colon V\to \Qbar, \quad x\mapsto P_{D_0,f_l,x}(m).
\]
We let:
\begin{eqn}\label{eqn_min_pol_a_2}
P_{D_0,f_l,m} := \prod_{x\in V} (T-a_{D_0,f_l,m}(x)) 
\quad\text{in $\QQ(\zeta_l)[T]$}
\end{eqn}
be the minimal polynomial over~$\QQ(\zeta_l)$ of~$a_{D_0,f_l,m}$.

%% -- Arakelov I (original results that hold for any curve)
\chapter{Applying Arakelov theory}
\label{sec_appl_Arakelov}

\author{B. Edixhoven and R. de Jong}

\bigskip

\bigskip

%authors: Bas and Robin

In this chapter we start applying Arakelov theory in order to derive a
bound for the height of the coefficients of the
polynomials~$P_{D_0,f_l,m}$ as in~(\ref{eqn_min_pol_a_2}). We proceed
in a few steps. The first step, taken in
Section~\ref{sec_height_and_intersection}, is to relate the height of
the $b_l(Q_{x,i})$ as in Section~\ref{sec_setup_tau} to intersection
numbers on~$X_l$. The second step, taken in
Section~\ref{sec_control_diff}, is to get some control on the
difference of the divisors $D_0$ and $D_x$ as in~(\ref{eqn_def_D'_x}).
Certain intersection numbers concerning this difference are bounded in
Theorem~\ref{maininequality}, in terms of a number of invariants in
the Arakelov theory on modular curves~$X_l$. These invariants will
then be bounded in terms of~$l$ in Sections~\ref{sec_bnd_height_X},
\ref{sec_bnd_theta}, and~\ref{sec_bnd_green_f}. Finally, in
Section~\ref{sec_final_estimates}, the height of the coefficients of
the~$P_{D_0,f_l,m}$ will be bounded. In this chapter, we do our best
to formulate the most important results, Theorem~\ref{thm_projection},
Theorem~\ref{mainequality}, and Theorem~\ref{maininequality} in the
context of curves over number fields, i.e., outside the context of
modular curves.

\section{Relating heights to intersection numbers}
\label{sec_height_and_intersection}
We pick up the notation as at the end of Section~\ref{sec_setup_tau},
so we have a modular curve~$X_l =X_1(5l)_\QQ$ with $l>5$, non-constant
morphisms $b_l$ and $x'_l\colon X_l \to \PP^1_\QQ$, and certain
divisors $D_x^\fin = \sum_{i=1}^{d_x} Q_{x,i}$ on~$X_{l,\Qbar}$ that
have support outside the cusps.  It is our objective in this
subsection to link the absolute height $h(b_l(Q_{x,i}))$ of the
algebraic numbers $b_l(Q_{x,i})$ to certain quantities coming from
Arakelov intersection theory. The height of $x'_l(Q_{x,i})$ will be
bounded in terms of $h(b_l(Q_{x,i}))$ in
Section~\ref{bds_height_x_and_y_in_height_b}. The final estimates for
these heights, depending only on~$l$, will be done in
Section~\ref{sec_final_estimates}.

\begin{thm}\label{thm_bound_hbQ}
Let $x$ be in $V$, and let $i$ be in~$\{1,\ldots,d_x\}$.  Let $K$ be a
number field containing $\QQ(\zeta_{5l})$ and such that $Q_{x,i}$ is
defined over~$K$. Let $\calX$ be the minimal regular model of $X_l$
over~$K$. Then we have:
\begin{multline*}
h(b_l(Q_{x,i}))   \leq \frac{1}{[K:\QQ]} 
\left( (Q_{x,i},b_l^*\infty)_\calX +
l^2  \sum_\sigma \sup_{X_\sigma} g_{\sigma} \right. \\
\left. + \frac{1}{2}\sum_\sigma\int_{X_\sigma}\log(|b_l|^2+1)\mu_{X_\sigma}\right)
+\frac{1}{2}\log 2\, .
\end{multline*}
Here $\sigma$ runs through the embeddings of $K$ into~$\CC$ and
$g_\sigma$ is the Arakelov-Green function on~$X_{l,\sigma}$.
\end{thm}  
In the next chapters we shall derive bounds that are polynomial in $l$
for all terms in the above estimate.

Theorem~\ref{thm_bound_hbQ} will be derived from
Theorem~\ref{thm_projection} below, which states a fairly general
result. We start with a definition. Let $K$ be a number field and
consider~$\PP^1_{O_K}$. Let $\infty$ denote the $O_K$-point $(1:0)$
of~$\PP^1_{O_K}$. For any section $P$ in $\PP^1_{O_K}(O_K)$ we define
by $(P,\infty)_{\PP^1}$ the degree (see~(\ref{eqn_ar_degree})) of
$P^*O_{\PP^1}(\infty)$, where $O_{\PP^1}(\infty)$ has the Fubini-Study
metric, i.e. the metric $\| \cdot \|_{\PP^1}$ given by:
\begin{eqn}\label{eqn:fub-study-metric}
\| 1 \|_{\PP^1} (x_0 : x_1) := 
\frac{ |x_1| }{(|x_0|^2 + |x_1|^2)^{1/2}}
\end{eqn}
over~$\PP^1_\CC$. Here $1$ is the tautological section
of~$O_{\PP^1}(\infty)$.

\begin{thm}\label{thm_projection}
Let $X$ be a geometrically irreducible, smooth and complete curve of
positive genus over a number field $K$ and let $\calX$ be a proper
semi-stable model of $X$ over the ring of integers $O_K$
of~$K$. Suppose that we have a non-constant morphism $f\colon
X\to\PP^1_K$ and a $K$-rational point $Q$ of~$X$ with $f(Q) \neq
\infty$.  Assume the following: the Zariski closure of
$\Supp(\divisor(f)_+) \cup \Supp(\divisor(f)_-)$ in $\calX$ is étale
over~$O_K$. For any closed point $s$ of~$\Spec(O_K)$, denote by
$m_s(f)$ the supremum of the multiplicities of $\divisor(f)_-$
on~$\calX$ along the irreducible components of the fibre at $s$
of~$\calX$.  Then we have the inequality:
\begin{multline*}
(f(Q),\infty)_{\PP^1} \leq (Q,f^*\infty)_{\calX}   +
\deg f \sum_\sigma \sup_{X_\sigma} g_{\sigma} \\
+ (1/2)\sum_\sigma \int_{X_\sigma}\log(|f|^2+1)\mu_{X_\sigma} 
+ \sum_s m_s(f) \log \#k(s) \, .
\end{multline*}
Here the first sum runs over the embeddings of~$K$ into~$\CC$, and the
last sum runs over the closed points of~$\Spec(O_K)$.
\end{thm}
\begin{proof}
Note that the locus of indeterminacy of $f$ on $\calX$ consists of
finitely many closed points. This implies that there exists a blow-up
$\tilde{\calX} \to \calX$ of $\calX$ such that $f$ extends to a
regular map $f\colon\tilde{\calX}\to\PP^1_{O_K}$.  For any such
$\tilde{\calX}$ we have by construction:
\begin{align*}
(f(Q),\infty)_{\PP^1} & = 
\deg f(Q)^*(O_{\PP^1}(\infty),\|\cdot\|_{\PP^1}) \\
&=\deg Q^*f^*(O_{\PP^1}(\infty),\|\cdot\|_{\PP^1})\\ 
&=\deg Q^*(O_{\tilde{\calX}}(f^*\ol{\infty}), \| \cdot \|_{\PP^1}) \, ,
\end{align*}
where we write $\ol{\infty}$ to emphasise that $f^*\ol{\infty}$ is the
inverse image under $f$ of $\infty(\Spec(O_K))$, and not the Zariski
closure $\ol{f^*\infty}$ of the inverse image under $f$ of
$\infty(\Spec(K))$.  If we let $\| \cdot \|_X$ denote the canonical
Arakelov metric on $O_X(f^* \infty)$ then we can write:
\begin{align*}
\deg Q^*(O_{\tilde{\calX}}(f^*\ol{\infty}), \| \cdot \|_{\PP^1}) & = 
\deg Q^*(O_{\tilde{\calX}}(f^*\ol{\infty}), \| \cdot \|_X \cdot
\frac{\| \cdot \|_{\PP^1}}{\| \cdot \|_X}) \\
& = \deg Q^*(O_{\tilde{\calX}}(f^*\ol{\infty}), \| \cdot \|_X) \\
& \qquad - \sum_\sigma
\log((\frac{\|\cdot\|_{\PP^1}}{\|\cdot\|_X})(Q_\sigma)) \\
& =  (Q,f^*\ol{\infty})_{\tilde{\calX}} 
+ \sum_\sigma
\log((\frac{\|\cdot\|_X}{\|\cdot\|_{\PP^1}})(Q_\sigma))  \, .
\end{align*}
A bound for $\log((\|\cdot\|_X/\|\cdot\|_{\PP^1})(Q_\sigma))$
follows by testing on the tautological section~$1$, giving:
\[
\log\|1\|_X(Q_\sigma)-\log\|1\|_{\PP^1}(Q_\sigma) =
g_\sigma(f^* \infty,Q_\sigma) + \frac{1}{2} \log(|f(Q_\sigma)|^2+1)\, .
\] 
Applying Proposition~\ref{prop_peter_bruin} below at this point shows
that:
\begin{multline*}
g_\sigma(f^* \infty,Q_\sigma) + \frac{1}{2} \log(|f(Q_\sigma)|^2+1)
\leq (\deg f)\sup_{X_\sigma} g_{\sigma} \\
+  \frac{1}{2}\int_{X_\sigma} \log(|f|^2+1) \,\mu_{X_\sigma}\, .
\end{multline*}
This accounts for the second and third terms in the bound of the
theorem. We are finished once we prove that
$(Q,f^*\ol{\infty})_{\tilde{\calX}} -
(Q,\ol{f^*\infty})_\calX$ is bounded by $\sum_s
m_s(f) \log \# k(s)$ for a particular choice of~$\tilde{\calX}$.
(The usual projection formula shows that in fact
$(Q,f^*\ol{\infty})_{\tilde{\calX}}$ is independent of the choice
of~$\tilde{\calX}$.)  On any $\tilde{\calX}$ we write $f^*\ol{\infty}$
as a sum $f^*\ol{\infty} = (f^*\ol{\infty})_\mathrm{hor} +
(f^*\ol{\infty})_\mathrm{vert}$ of a horizontal and a vertical
part. Note that $(f^*\ol{\infty})_\mathrm{hor} = \ol{f^*\infty}$,
with the Zariski closure now taken in~$\tilde{\calX}$. Since the local
intersection multiplicities of $Q$ and $\ol{f^*\infty}$ do not go up
when passing from $\calX$ to~$\tilde{\calX}$, we have
$(Q,(f^*\ol{\infty})_\mathrm{hor})_{\tilde{\calX}} =
(Q,\ol{f^*\infty})_{\tilde{\calX}} \leq (Q,\ol{f^*\infty})_\calX$ and
thus we are reduced to proving that
$(Q,(f^*\ol{\infty})_\mathrm{vert})_{\tilde{\calX}}$ is bounded from
above by $\sum_s m_s(f) \log \# k(s)$ for a particular choice
of~$\tilde{\calX}$.

We exhibit a specific blow-up, and we calculate which multiplicities
$f$ acquires along the irreducible components of the vertical fibres
of this blow-up.  Note that the locus of indeterminacy of $f$ on
$\calX$ consists precisely of the closed points of $\calX$ where an
irreducible component $C$ from the zero divisor $\divisor(f)_+$ of $f$
on $\calX$ and an irreducible component $C'$ from its polar divisor
$\divisor(f)_-$ meet.  Now since by assumption the Zariski closure of
$\Supp(\divisor(f)_+) \cup \Supp(\divisor(f)_-)$ in $\calX$ is étale
over~$O_K$, this can only happen when at least one of~$C$, $C'$ is
vertical. In such points where this happens we have to perform a
sequence of successive blowings-up until a component arises with
multiplicity~$0$ for~$f$, so that the components with positive
multiplicities and the components with negative multiplicities are
separated from each other.

We begin by observing that it will cause no harm if we pass to a
finite extension $K \to K'$. Indeed, both the left hand side and the
right hand side of the inequality that we wish to prove get multiplied
by $[K':K]$ if we do this. Here is why: for the terms
$(f(Q),\infty)_{\PP^1}$ and
$(Q,f^*\infty)$ the scaling by a factor $[K':K]$
follows from general properties of the Arakelov intersection product,
cf.~\cite{Faltings1}, p.~404 for example.
%% The statement there is not precise enough to my (Bas) liking; as
%% stated, it is just wrong. So if there is a better reference I am
%% for it.
(Note that it is understood that over~$K'$, intersection products are
taken on the minimal resolution of the pullback of the model~$\calX$.)
That the term $\sum_\sigma \sup_{X_\sigma} g_\sigma$ scales by a
factor $[K':K]$ is obvious. Finally, fix a closed point~$s$
of~$\Spec(O_K)$ and let $s'$ be any closed point of $\Spec(O_{K'})$
above it. Denoting by~$e_{s'}$ the ramification index of~$s'$ over~$s$
and by~$f_{s'}$ the degree of the residue field extension of $s'$
over~$s$, we see that for any $s'$ above~$s$, the integer $m_s(f)$
gets multiplied by~$e_{s'}$, and the number $\log \#k(s)$ gets
multiplied by~$f_{s'}$. Using that $\sum_{s'} e_{s'} f_{s'} = [K':K]$,
the sum running over the closed points $s'$ above~$s$, we see finally
that also the term $\sum_s m_s(f) \log \# k(s)$ gets multiplied by
$[K':K]$.

Starting with $\calX$ over~$O_K$, we first do the following. Let $x$
be a closed point on $\calX$ that is the intersection of a vertical
component~$C$ and a horizontal component~$C'$ having non-zero
multiplicities $m$ and~$m'$ for $f$ that have different sign. After
blowing up in~$x$, we obtain an exceptional divisor~$E$ whose
multiplicity for~$f$ is $m+m'$. We have two distinguished points
on~$E$, one lying on the strict transform of~$C$, and one lying on the
strict transform of~$C'$. At exactly one of them there is a sign
change for the multiplicities, or $m+m'=0$. If a sign change happens
at the double point lying on the strict transform of~$C'$, then we repeat
the process. If a sign change happens at the double point lying on the
strict transform of~$C$ or if $m+m'=0$, we stop, and continue with a
new point~$x'$, if available.

We end up with a blow-up $\calX'\to\calX$ such that an intersection of
two irreducible components~$C$, $C'$ that have different sign in
$\divisor(f)$ on $\calX'$ only occurs for~$C$, $C'$ both vertical. For
this, we did not yet need to extend the ground field~$K$. In order to
continue, we note the following. Suppose that we have a closed
point~$x$ on the model~$\calX'$ of $X$ over $O_K$ which is a double
point of a vertical fibre, and two irreducible components~$C$, $C'$ of
that vertical fibre pass through~$x$, having non-zero multiplicities
$m$ and $m'$ for $f$ that differ in sign. Assume that $K\to K'$ is a
Galois extension that ramifies over the image of~$x$ in~$\Spec(O_K)$,
with a ramification index~$e$ that is a multiple of~$m-m'$. In passing
to the minimal resolution~$\tilde{\calX}$ of~$\calX'_{O_{K'}}$, the
point~$x$ gets replaced by a chain of $e-1$ projective lines of
self-intersection~$-2$. The multiplicities of~$f$ along these
components change in $e$ steps from~$em$ to~$em'$, so that the steps
are $m-m'$ and a multiplicity~$0$ will appear somewhere, because
$m-m'$ is a divisor of~$em$.

Thus we see how we can reach our goal: take a Galois extension $K\to
K'$ that ramifies as specified above over the images in~$\Spec(O_K)$
of the double points where components meet with a different sign
for~$f$. (This is always possible.)  By our remarks above, it suffices
to prove the inequality over~$K'$. By construction, the morphism~$f$
extends over the model~$\tilde{\calX}$ that arises in this
way. Moreover, it follows from the construction that for $s'$ a closed
point of $\Spec(O_{K'})$ and $s$ its image in $\Spec(O_K)$ we have
$m_{s'}(f)\leq e_{s'}m_s(f)$. Hence the sum of the local intersection numbers
$(Q,(f^*\ol{\infty})_\mathrm{vert})_{s'}$ for all
$s'$ over~$s$ is bounded
from above by~$[K':K] m_s(f)\log\# k(s)$. This is what we needed to prove.
\end{proof}

\begin{prop}\label{prop_peter_bruin}
Let $f\colon X\to\PP^1$ be a finite morphism of Riemann surfaces with
$X$ connected and of positive genus.  Consider on $X-f^{-1}\infty$ the
function:
\[ 
h(x) = g(f^*\infty,x) + \frac{1}{2} \log(|f(x)|^2+1) \, . 
\]
Then $h$ extends uniquely to an element of $\calC^\infty(X)$, also
denoted~$h$.  For all $x\in X$ we have:
\[ 
h(x)\leq (\sup_X g) \deg f + \frac{1}{2}\int_X \log(|f|^2+1) \,\mu_X
\, .
\]
\end{prop}
\begin{proof}
Let us first show that $h$ extends to a $C^\infty$ function on~$X$. In
fact, as the beginning of the proof of Theorem~\ref{thm_projection}
indicates, $h$ is the logarithm of the function $x\mapsto
\|{\cdot}\|_X(x)/\|{\cdot}\|_{\PP^1}(x)$ that gives the quotient of
two metrics on $O_X(f^*\infty)$: the Arakelov metric $\|{\cdot}\|_X$
and the pullback of the Fubini-Study metric~$\|{\cdot}\|_{\PP^1}$
(see~\ref{eqn:fub-study-metric}). But then $h$ is
in~$\calC^\infty(X)$. We have, by~(\ref{eqn_green_inverts}):
\[
h(x) = 
\int_{y\in X} -g(x,y) \frac{1}{\pi i} (\partial\ol{\partial}h)y 
+\int_X h\,\mu_X \,.
\]
For $s$ any local holomorphic generator of $O_X(f^*\infty)$ we have:
\[
\frac{1}{\pi i}\partial\ol{\partial} h = 
\frac{1}{2\pi i}\partial\ol{\partial}\log(\|s\|_X^2) 
-\frac{1}{2\pi i}\partial\ol{\partial}\log(\|s\|_{\PP^1}^2)
=(\deg f) \mu_X - f^* \mu_{\PP^1} \, , 
\]
where $\mu_{\PP^1}$ is the curvature form of the Fubini-Study metric
on~$O_{\PP^1}(\infty)$. Substituting this in the previous equality,
and using that for all $x$ in $X$ we have $\int_{y\in
X}g(x,y)\mu_X(y)=0$, we get:
\[
h(x) = \int_{y\in X}g(x,y)(f^*\mu_{\PP^1})y + \int_X \frac{1}{2}
\log(|f|^2+1)\,\mu_X\,. 
\]
As $\mu_{\PP^1}$ defined as the curvature form of the Fubini-Study
metric on $O_{\PP^1}(\infty)$, we have
$\int_{\PP^1}\mu_{\PP^1}=\deg(O_{\PP^1}(\infty))=1$. As the metric is
invariant under the transitive action of~$\SU_2$, $\mu_{\PP^1}$ is
everywhere positive. In fact, one can compute that $\mu_{\PP^1} =
(i/2\pi)\, dz\, d\ol{z} /(1+|z|^2)^2$.  We end up with:
\[
h(x) \leq (\sup_X g) \deg f + \frac{1}{2} \int_X \log(|f|^2+1)\mu_X 
\]
as required.
\end{proof}

\begin{proof}
[Proof of Theorem~\ref{thm_bound_hbQ}] In order to simplify our
notation, we drop the subscript $l$ in~$b_l$. For $P$ in $\PP^1(O_K)$
we put:
\[
h'(P) = (P,\infty)/[K:\QQ]. 
\]
Then we have $h(P) \leq
h'(P)+(1/2)\log 2$, as for all $x\in\CC^2$ we have
$|x_1|^2+|x_2|^2\leq 2\max\{|x_1|,|x_2|\}^2$.  In order to bound
$h'(b(Q_{x,i}))$ from above we want to apply
Theorem~\ref{thm_projection}. It follows from the definition of the
morphism~$b$ that both the zero divisor $\divisor(b)_+$ and the polar
divisor $\divisor(b)_-$ of~$b$ on~$X_l$ have as their support only
$K$-rational closed points, namely, cusps.  In particular, we never
have $b(Q_{x,i})=\infty$, by construction of~$D''_x$. We have also
seen that the Zariski closure in $\calX$ of $\Supp(\Cusps)$ is étale
over~$O_K$ (as $O_K$-valued points the cusps are disjoint), and hence
the same holds for the Zariski closure in~$\calX$ of
$\Supp(\divisor(b)_+) \cup \Supp(\divisor(b)_-)$.
Theorem~\ref{thm_projection} now gives us that:
\begin{multline*}
[K:\QQ]{\cdot}h'(b(Q_{x,i}))\leq (Q_{x,i},b^*\infty)_\calX +
\deg b \sum_\sigma \sup_{X_\sigma} g_{\sigma} \\
+ (1/2)\sum_\sigma\int_{X_\sigma} \log(|b|^2+1)\mu_{X_\sigma}
+ \sum_s m_s(b) \log \# k(s)
\end{multline*}
with $m_s(b)$ the supremum of the multiplicities of $\divisor(b)_-$ on
$\calX$ along the irreducible components of the fibres of $\calX$
at~$s$.  We are done if we can prove that $\deg b$ is at most~$l^2$,
and that $m_s(b)=0$ for all~$s$.  The definition of~$b$ shows directly
that its degree is $l^2{-}1$ (it is the degree of the natural morphism
from $X_1(5l)$ to~$X_1(5)$).  Let us now show that $m_s(b)=0$ for
all~$s$.  For this we evidently need information on the
divisor~$\divisor(b)$ on~$\calX$.

We start with working over~$\QQ(\zeta_{5l})$. From the discussion in
Section~\ref{sec_constr_D} we recall that there is a fine moduli
scheme $Y_1(5l)_{\ZZ[\zeta_{5l}]}$ over $ \ZZ[\zeta_{5l}]$ of elliptic
curves with balanced level structure (terminology
from~\cite{Katz-Mazur}). Let $E(5l)\to Y_1(5l)_{\ZZ[\zeta_{5l}]}$ be
the universal elliptic curve and let~$P_5$, $P_l$ be the tautological
points of order~$5$ and~$l$. From Proposition~\ref{prop_Y_1_5} we
recall the elliptic curve $E\to Y_1(5)'$ with
$Y_1(5)'=\Spec(\ZZ[b,1/\discr(E)])$. The elliptic curve $E(5l) \to
Y_1(5l)_{\ZZ[\zeta_{5l}, 1/5]}$ arises from $E \to Y_1(5)'$ by a
unique base change $Y_1(5l)_{\ZZ[\zeta_{5l}, 1/5]}\to Y_1(5)'$. This
gives the regular function $b$ on~$Y_1(5l)_{\ZZ[\zeta_{5l}, 1/5]}$. As
$b$ is invertible on $Y_1(5)'$ it is invertible
on~$Y_1(5l)_{\ZZ[\zeta_{5l}, 1/5]}$. We conclude that~$\divisor(b)$
on~$X_1(5l)_{\ZZ[\zeta_{5l}]}$ is a certain linear combination of the
irreducible components of the closed subschemes~$\Cusps$ and
$X_1(5l)_{\FF_5[\zeta_l]}$.  In order to find this linear combination,
we examine the multiplicities of $b$ along the irreducible components
that we have isolated.

We start with the multiplicities along the irreducible components of
the divisor~$\Cusps$. It is sufficient to study the situation
over~$\CC$, and here we can make a beginning by looking at
$\divisor(b)$ on~$X_1(5)_\CC$.  From the equations in
Proposition~\ref{prop_Y_1_5} we obtain that over the cusp~$0$
of~$X_0(5)_\CC$ lie two cusps, say $c_1$ and~$c_2$, with~$c_1$, say,
corresponding to the $5$-gon with the tautological point of order~$5$
being on a component adjacent to the connected component of~$0$, and
the other, $c_2$, corresponding to the $5$-gon with the tautological
point of order~$5$ being on a component that is not adjacent to the
connected component of~$0$. We have $\divisor(b) = \pm (c_1 - c_2)$
on~$X_1(5)_\CC$; we could compute the exact sign but that is not
important for us. The divisor of $b$ on $X_1(5l)_\CC$ is obtained by
pulling back its divisor on $X_1(5)_\CC$ via the forgetful map
$X_1(5l)_\CC \to X_1(5)_\CC$. Hence, pulling back the divisor
$c_1-c_2$ we get plus or minus the divisor of $b$ on $X_1(5l)_\CC$;
the multiplicities are just the ramification indices above the cusps
$c_1$ and~$c_2$.  Since these are in $\{1,l\}$, we obtain that the
multiplicities of~$b$ along the irreducible components of $\Cusps$ are
just~$1$ or~$l$ in absolute value.

Next we calculate the multiplicities of $b$ along the irreducible
components of $X_1(5l)_{\FF_5[\zeta_l]}$. The structure of a connected
component of the scheme $X_1(5l)_{\FF_5[\zeta_l]}$ is as follows: it
consists of two irreducible components, one on which $P_5$ has
order~$1$, and one having an open part where $P_5$ has
order~$5$. These two irreducible components intersect (transversally)
in the supersingular points. 

We denote by $\Gamma$ the union of the irreducible components
over~$\FF_5$ on which $P_5$ has order~$1$. The construction of the
scheme~$Y_1(5)'$ immediately gives us a forgetful map
$X_1(5l)_{\ZZ[\zeta_{5l}]} - \Gamma - \Supp(\Cusps) \to
Y_1(5)'$. Since $b$ is invertible on~$Y_1(5)'$, the same holds for $b$
along the irreducible components of $X_1(5l)_{\FF_5[\zeta_l]}$, except
possibly for the irreducible components in~$\Gamma$. But the
multiplicity of $b$ along such an irreducible component is then also
zero, as can be seen by the following argument. Let $C \cup C'$ be a
connected component of $X_1(5l)_{\FF_5[\zeta_l]}$, with the
irreducible component~$C$ corresponding to $P_5$ having order~$1$. All
the horizontal components of $\divisor(b)$ on
$X_1(5l)_{\ZZ[\zeta_{5l}]}$ specialise to~$C'$. We know that $b$ has
multiplicity~$0$ along $C'$ and hence it restricts to a non-trivial
rational function, also denoted~$b$, on~$C'$. The degree of $b$ on
$C'$ is zero, or equivalently $m (C',C) +
(C',\divisor(b)_\mathrm{hor}) =0$, where $m$ is the multiplicity of
$b$ along~$C$. Now, since $(C',\divisor(b)_\mathrm{hor})$ is zero and
$(C',C)$ isn't, we get $m=0$.

All in all we conclude that the absolute values of the multiplicities
of the irreducible components in~$\divisor(b)$ on
$X_1(5l)_{\ZZ[\zeta_{5l}]}$ are bounded by a constant times~$l$, and
that all multiplicities along irreducible components of fibers over
closed points $s$ of $\Spec(\ZZ[\zeta_{5l}])$ are zero. In particular, for
all closed points $s$ we have $m_s(b)=0$. This implies in fact
that the rational function $b$ on $X_1(5l)_{\ZZ[\zeta_{5l}]}$ extends
to a morphism to~$\PP^1_\ZZ$. As this is a useful fact, we record
it in a Proposition.
%%% put this in a better place??????
This completes the proof of Theorem~\ref{thm_bound_hbQ}.
\end{proof}
\begin{prop}
Let $l>5$ be prime a prime number. The rational function $b_l$ on
$X_1(5l)_{\ZZ[\zeta_{5l}]}$ from Proposition~\ref{prop_bxy} extends to
a morphism to~$\PP^1_\ZZ$.
\end{prop}

\section{Controlling $D_{\protect\lowercase{x}}-D_0$}
\label{sec_control_diff}
In this subsection, the hypotheses are as follows (unless stated
otherwise). We let $K$ be a number field, $O_K$ its ring of integers,
$B:=\Spec(O_K)$, $p\colon \calX\to B$ a regular, split semi-stable
curve over~$B$ whose generic fibre~$X\to\Spec K$ is geometrically
irreducible and of genus $g\geq 1$. We let $D$ be the closure in
$\calX$ of an effective divisor of degree~$g$ (also denoted~$D$)
on~$X$. We let $x$ be a $K$-rational torsion point of the Jacobian
of~$X$, i.e., a torsion element of~$\Pic(X)$, which has the property
that there is a unique effective divisor $D_x$ on~$X$ such that
$x=[D_x-D]$. Finally, we let $P\colon B\to\calX$ be a section of~$p$,
i.e., an element of~$\calX(B)$.

We denote by $\Phi_{x,P}$ the unique finite vertical \emph{fractional}
divisor~$\Phi$ (i.e., with rational coefficients that are not
necessarily integral) on $\calX$ such that $(D_x-D-\Phi,C)=0$ for all
irreducible components $C$ of fibres of~$p$, and such that $P(B)$ is
disjoint from the support of~$\Phi$. It is not difficult to see that a
$\Phi$ satisfying the first condition exists and that it is unique up
to adding multiples of fibers of~$p$ (the intersection pairing
restricted to the divisors with support in a fibre is negative
semi-definite); see Lemme~6.14.1 of~\cite{Moret-Bailly1}. The second
condition removes the ambiguity of adding multiples of fibres.

We denote by $\delta_s$ the number of singular points in the geometric
fibre at a closed point $s$ of~$B$.
\begin{thm} \label{mainequality} 
The $O_B$-module $\rR^1 p_* O_\calX (D_x)$ is a torsion module
on~$B$, and we have:
\begin{multline*}
(D_x,P) + \log \# \rR^1 p_* O_\calX(D_x) +
\frac{1}{8} (\omega_{\calX/B},\omega_{\calX/B}) + 
\frac{1}{8} \sum_s \delta_s \log \# k(s) \\ 
= (D,P) - \frac{1}{2}(D+\Phi_{x,P},D+\Phi_{x,P}-\omega_{\calX/B}) + 
\frac{1}{2} \deg \det p_* \omega_{\calX/B}   \\ 
+\sum_\sigma \int_{X_\sigma} \log \| \vartheta \|
(D_x^\sigma  - Q) \cdot \mu_\sigma (Q) + \frac{g}{2} [K:\QQ] \log(2\pi).
\end{multline*}
Here $s$ runs over the closed points of~$B$, and $\sigma$ runs through
the complex embeddings of~$K$.
\end{thm}
We derive Theorem \ref{mainequality} from three lemmas. For the moment
we work in~$\QQ\otimes_\ZZ\widehat{\mathrm{Pic}}(\calX)$.
\begin{lem} \label{firstlemma}
The admissible line bundles:
\[
O_\calX(D_x{-}D) \otimes p^*P^*O_\calX(D_x{-}D)^\vee  
\quad\text{and}\quad O_\calX(\Phi_{x,P}) 
\]
are numerically equivalent. That is, for any admissible line bundle
$F$ on $\calX$ we have:
\[ 
(O_\calX(D_x-D) \otimes p^*P^*O_\calX(D_x-D)^\vee ,F) = 
(O_\calX(\Phi_{x,P}),F).
\]
\end{lem}
\begin{proof} 
In this proof we just write $\Phi$ for~$\Phi_{x,P}$. We denote the
first line bundle in the lemma by~$\Psi$. Since $D_x-D$ is torsion,
there is a positive integer $N$ such that $\Psi^{\otimes N}$ is
trivial on the generic fibre as a classical line bundle (that is,
without taking the metrics into account). We have a canonical
isomorphism $P^* \Psi \isomlto O_B$ on~$B$. Combining, we find that
$\Psi^{\otimes N}$ has a rational section~$s$ with $\divisor_\calX(s)$
vertical and with $P^*s \mapsto 1$. The latter condition implies that
$P$ intersects to zero with $\divisor_\calX(s)$ for the Arakelov
intersection product. On the other hand, as
$p^*P^*O_\calX(D_x-D)^\vee$ is trivial on the fibres of~$p$ over
finite places of~$B$, we have $(N(D_x-D)-\divisor_\calX(s),C)=0$ for
all irreducible components~$C$ of fibres of~$p$. Hence in fact $\Phi =
\frac{1}{N} \divisor_\calX(s)$. To prove the lemma, it suffices now to
prove that $\Psi^{\otimes N} \isomlto O_\calX(\divisor_\calX(s))$ given
by $s \mapsto 1$, with $1$ the tautological section, is an
isometry. Because of admissibility, it suffices to check that this is
so when restricted to $P$; but here we get the canonical isomorphism
$P^* \Psi \isomlto O_B$. This is indeed an isometry by the definition
of~$\Psi$.
\end{proof}
\begin{lem} \label{secondlemma}
Let $X$ be a compact Riemann surface of genus $g\geq 1$. Let $D$ be an
effective divisor on $X$ of degree $g$ satisfying $h^0(D)=1$. Then
the determinant of cohomology $\lambda(O_X(D))$ of $D$ is identified
with $H^0(X,O_X(D))$. Further, the formula:
\[
\log \| 1 \| + \frac{\delta(X)}{8} + \int_X \log \| \vartheta
\|(D-Q) \cdot \mu_X(Q) =0 
\]
holds for the length (with respect to Faltings' metrisation of the
determinant of cohomology) of the tautological section~$1$ of
$H^0(X,O_X(D))$.
\end{lem}
\begin{proof}
Since $h^0(D)=1$, $H^0(X,O_X(D))=\CC$. Therefore, the set of points
$Q$ on $X$ such that $h^0(D-Q) >0$ is the support of~$D$. Let $Q$ be
a point outside the support of~$D$. Then $h^0(D-Q)=0$. According to
the axioms for the metrisation of the determinant of cohomology, the
exact sequence:
\[
0 \to O_X(D-Q) \to O_X(D) \to Q_*Q^* O_X(D) \to 0
\]
gives rise to an isometry:
\begin{multline*}
\lambda(O_X(D)) \isomlto \lambda(O_X(D-Q))
\otimes Q^*O_X(D) \cong \\
\cong O(-\Theta)[O_X(D-Q)] \otimes Q^*O_X(D).
\end{multline*}
Taking the norm on left and right of a tautological section we obtain,
using~(\ref{eqn_deltaX}):
\[
\|1\| = \exp(-\delta(X)/8) \cdot \| \vartheta \|(D-Q)^{-1} \cdot G(D,Q),
\]
where $\log G(D,Q)=g_{D,\mu_X}(Q)$.  Taking logarithms and then
integrating against $\mu_X(Q)$ gives the result.
\end{proof}
\begin{lem} (Noether formula) \label{noether}
We have: 
\begin{multline*}
12 \, \deg \det p_* \omega_{\calX/B} =
(\omega_{\calX/B},\omega_{\calX/B}) + \sum_s \delta_s \log \# k(s) \\
+ \sum_\sigma \delta(X_\sigma) - 4g [K:\QQ] \log (2\pi) \, ,
\end{multline*}
the first sum running over the closed points of~$B$, the second sum
running over the complex embeddings of~$K$.
\end{lem}
\begin{proof} See~\cite{Faltings1} and~\cite{Moret-Bailly2}.
\end{proof}
\begin{proof}
[Proof of Theorem~\ref{mainequality}] We first show that $\rR^1 p_*
O_\calX(D_x)$ is a torsion module. As it is a coherent $O_B$-module,
it suffices to show that it is zero on the generic point of~$B$, i.e.,
that $\rH^1(X,O_X(D_x))$ is zero. By Riemann-Roch, we have 
$h^0(D_x)-h^1(D_x)=1$. By definition of~$D_x$, we have
$O_X(D_x)\cong\calL_x(D)$. And by construction of~$D$, we have
$h^0(\calL_x(D))=1$. This shows that $h^1(D_x)=0$.

Let us now prove the identity in Theorem~\ref{mainequality}. We start
by noting that, by~(\ref{eqn_def_int_sec}):
\[
(D_x-D,P)=\deg P^*O_\calX(D_x-D).
\]
By Lemma~\ref{firstlemma}, $O_\calX(D_x) \otimes
p^*P^*O_\calX(D_x-D)^\vee$ and $O_\calX(D+\Phi_{x,P})$ are numerically
equivalent. The Riemann-Roch theorem then gives:
\begin{multline*}
\deg \det\rR p_*(O_\calX(D_x) \otimes p^*P^*O_\calX(D_x-D)^\vee) = \\
= \frac{1}{2} (D+\Phi_{x,P},D+\Phi_{x,P}-\omega_{\calX/B}) + 
\deg \det p_*\omega_{\calX/B} \, .
\end{multline*}
By the projection formula for the determinant of cohomology we can
write the left-hand side as:
\begin{multline*}
\deg \det\rR p_*(O_\calX(D_x) \otimes p^*P^*O_\calX(D_x-D)^\vee) = \\
= \deg \det\rR p_* O_\calX(D_x) - \deg P^*O_\calX(D_x-D) \, 
\end{multline*}
Since $p_*O_\calX(D_x)$ is canonically trivialised by the
function~$1$, the term $\deg \det\rR p_* O_\calX(D_x)$ can be
computed as:
\[
\deg \det\rR p_* O_\calX(D_x) = 
-\sum_\sigma \log \| 1 \|_\sigma - \log \# \rR^1 p_* O_\calX(D_x),
\]
where for each complex embedding $\sigma$, the norm $\| 1 \|_\sigma$
is the length of the tautological section~$1$ of
$\lambda(O_{X_\sigma}(D_x))=H^0(O_{X_\sigma}(D_x))$. By Lemma
\ref{secondlemma} we can then write:
\begin{align*}
\deg \det\rR p_* O_\calX(D_x) = & \sum_\sigma \int_{X_\sigma} 
\log \|\vartheta \|_\sigma(D_x^\sigma - Q) \cdot \mu_\sigma(Q) \\
& + \sum_\sigma \delta(X_\sigma)/8 - \log \# \rR^1 p_* O_\calX(D_x) \, .
\end{align*}
Combining everything gives:
\begin{align*}
(D_x-D,P) = & -\frac{1}{2}(D+\Phi_{x,P},D+\Phi_{x,P}- \omega_{\calX/B}) - 
\deg \det p_* \omega_{\calX/B} \\ 
&+ \sum_\sigma \int_{X_\sigma} 
\log \| \vartheta \|(D_x^\sigma - Q) \cdot \mu_\sigma(Q) \\
& + \sum_\sigma \delta(X_\sigma)/8 - \log \# \rR^1 p_* O_\calX(D_x) \, .
\end{align*}
We obtain the required formula upon eliminating $\sum_\sigma
\delta(X_\sigma)/8$ with the Noether formula, Lemma~\ref{noether}.
\end{proof}
\begin{thm}\label{maininequality} 
We have an upper bound:
\[
\begin{aligned}
(D_x,P) + \log \# \rR^1 p_* O_\calX(D_x) \leq 
& -\frac{1}{2} (D,D-\omega_{\calX/B}) \\
& + 2g^2\sum_{s \in B} \delta_s \log \#k(s) \\
& + \sum_\sigma \log \|\vartheta\|_{\sigma,\sup} + 
\frac{g}{2}[K:\QQ]\log(2\pi) \\ 
&+ \frac{1}{2} \deg \det p_* \omega_{\calX/B} + (D,P) \, ,
\end{aligned} 
\] 
where $s$ runs through the closed points of~$B$, and where the
supnorm $\|\vartheta\|_{\sigma,\sup}$ is taken
over~$\Pic^{g-1}(X_\sigma)$.
\end{thm}
The required upper bound follows directly from
Theorem~\ref{mainequality} by using Lemma~\ref{removingPhi} below and
the fact that $(\omega_{\calX/B},\omega_{\calX/B}) \geq 0$ (cf.~Thm.~5
of~\cite{Faltings1}).
\begin{lem} \label{removingPhi}
We have an upper bound:
\begin{multline*}
-\frac{1}{2} (D+\Phi_{x,P}, D+\Phi_{x,P}-\omega_{\calX/B}) \leq \\
\leq -\frac{1}{2}(D,D-\omega_{\calX/B}) 
+ 2g^2 \sum_{s \in B} \delta_s \log \#k(s) \, ,
\end{multline*}
with $s$ running through the closed points of~$B$.
\end{lem}
\begin{proof}
In this proof we just write $\Phi$ for~$\Phi_{x,P}$. By the definition
of $\Phi$, we have $(D_x-D-\Phi,\Phi)=0$, or in other words,
$(\Phi,\Phi)=(D_x-D,\Phi)$.  Using this we can write:
\[
-\frac{1}{2}(D+\Phi,D+\Phi-\omega_{\calX/B}) = 
-\frac{1}{2}(D,D-\omega_{\calX/B}) + 
\frac{1}{2}(\Phi,\omega_{\calX/B}-D-D_x)\, . 
\] 
We write $\Phi=\sum_C \Phi(C) \cdot C$, and for any finite fibre~$F_s$
of~$p$ we put $A_s := \sup_C |\Phi(C)|$ with $C$ running through the
irreducible components of~$F_s$.  Since $\omega_{\calX/B},D$ and
$D_x$ intersect any irreducible component $C$ with non-negative
multiplicity, we find
\begin{align*}
\frac{1}{2}(\Phi ,\omega_{\calX/B}-D-D_x) & \leq
\frac{1}{2} \left( \sum_s A_s F_s,\omega_{\calX/B}+D+D_x \right) \\
& \leq 2g \sum_s A_s \log \#k(s) \, .
\end{align*}
We are going to prove that $A_s \leq g\delta_s$, and then we are done.
So let $s$ be a finite place of~$B$. Let $S_0$ be the set of
irreducible components of~$F_s$, and let $S_1$ be the set of double
points on~$F_s$. Let $\Gamma_s$ be the dual graph of~$F_s$ (thus, the
set of vertices of~$\Gamma_s$ corresponds to~$S_0$, the set of edges
corresponds to~$S_1$, and the graph is defined by the incidence
relations). Choose an orientation on~$\Gamma_s$. This gives rise to
the usual source and target maps $s$ and $t\colon S_1 \to S_0$.
Consider the boundary and coboundary maps $d_* = t_* - s_*\colon
\QQ^{S_1} \to \QQ^{S_0}$ and $d^* = t^* - s^* : \QQ^{S_0} \to
\QQ^{S_1}$. Then $d_* d^*\colon \QQ^{S_0} \to \QQ^{S_0}$ is given by
minus the intersection matrix of~$F_s$. In particular, the map
$d_*d^*$ sends $\Phi$ to the map $u\colon C\mapsto
-(\Phi,C)=(D-D_x,C)$.  The kernel of~$d_* d^*$ consists exactly of
the constant functions, and the image consists of the orthogonal
complement of the constant functions. Now consider the graph
$\Gamma_s$ as an electric circuit, where each edge has a resistance of
1 Ohm. By Ohm's law and by spelling out the maps~$d_*$ and~$d^*$ we
see that if we let at each vertex~$C$ a current of $u(C)$ Ampère enter
the circuit, subject to the condition that $\sum_C u(C)=0$, the
potentials $\varphi(C)$ at each vertex $C$ will be given, up to
addition of a constant function, by a solution of the equation $d_*d^*
\varphi = u$. Hence $\Phi$ is the potential corresponding to the
current $C\mapsto (D-D_x,C)$, normalised by the condition that
$\Phi(C_P)=0$ with $C_P$ the component that $P$ specialises to. We
must bound the $|\Phi(C)|$ for $C$ varying over~$S_0$. The worst case
that may happen is that $\Gamma_s$ is a chain, with $D'$ and $D$
specialising entirely to the beginning and end point, respectively. In
this case, the biggest potential difference is $g \cdot (\# S_0-1)$ in
absolute value, so that we arrive at $|\Phi(C)| \leq g \cdot (\#
S_0-1)$.  Now note that $\Gamma_s$ is connected and that $X/K$ has
split semi-stable reduction. This gives $\# S_0 -1 \leq \delta_s$ and
hence $|\Phi(C)| \leq g \delta_s$, as required.
\end{proof}

%% -- Merkl
\chapter{An upper bound for Green functions on Riemann surfaces}
\label{subsec_merkl}

\author{F. Merkl}

\bigskip

\bigskip

% author: Franz Merkl%%%%%%%%%%%%%%%%%%%%%%%%%%%%%%%%%%%%%%%%%%%%%%%%%
We begin with explaining the setup and the results of this subsection.
Let $X$ be a compact Riemann surface, endowed with a 2-form
$\mu\geq 0$ that fulfils $\int_X\mu=1$.  Let $\ast$ denote rotation by
$90^\circ$ in the cotangential spaces (with respect to the holomorphic
structure); in a coordinate $z=x+iy$ 
this means $\ast dx=dy$, $\ast dy=-dx$ and, equivalently, 
$\ast dz=-idz$, $\ast d\bar
z=id\bar z$.  In particular, the Laplace operator on real $C^\infty$
functions on $X$ can be written as $d{\ast}d=2i\partial\bar\partial$.

For $a,b\in X$, let $g_{a,b}\colon X-\{a,b\}\to\mathbb{R}$ 
be the (unique) solution on $X$ (in
the sense of distributions) of the following 
differential equation:
\[
d{\ast} d g_{a,b}=\delta_a-\delta_b\quad\text{on $X$}
\]
with the normalising condition:
\[
\int_{X-\{a,b\}} g_{a,b}\,\mu=0.
\]
Note that $g_{a,a}=0$.  The distributional differential equation for
$g_{a,b}$ is equivalent to the following two more elementary
conditions.  Firstly, $g_{a,b}$ is a real-valued harmonic function on
$X-\{a,b\}$.  Secondly, it has logarithmic singularities near $a$ and
near $b$ of the following type: for any local coordinate $z$ near~$a$,
the function $P\mapsto g_{a,b}(P)-(2\pi)^{-1}\log|z(P)-z(a)|$ extends
to a harmonic function in a neighbourhood of~$a$, and for any local
coordinate $w$ near~$b$, the function $P\mapsto
g_{a,b}(P)+(2\pi)^{-1}\log|w(P)-w(b)|$ extends to a harmonic function
in a neighbourhood of~$b$.  The existence of such a function $g_{a,b}$
is shown in the theorem in paragraph II.4.3 on page~49
in~\cite{Farkas-Kra}.  However, for $a$ close to~$b$, the proof of
Lemma~\ref{merkl_lemma0} below also shows the existence of $g_{a,b}$
as a by-product.  The difference of any two solutions of the
differential equation for $g_{a,b}$ extends to a global harmonic
function on $X$ and thus is a constant. Hence, the normalising
condition for $g_{a,b}$ determines the function $g_{a,b}$ uniquely.

Now, for $x\in X- \{a\}$, let
\[
g_{a,\mu}(x) :=\int_{b\in X-\{x\}} g_{a,b}(x)\,\mu(b).
\]
Then we have:
\[
d{\ast} d g_{a,\mu}(x)=\delta_a-\mu
\]
in the sense of distributions, and:
\[
\int_{X-\{a\}} g_{a,\mu}\,\mu=0.
\]

We consider an atlas  of $X$ consisting of $n$ local coordinates:
\[
z^{(j)}\colon U^{(j)}\to \CC, \quad j=1,\ldots, n,
\]
such that each range $z^{(j)}[U^{(j)}]$ contains the closed unit disk.
For any radius $0<r\le 1$ and $j\in\{1,\ldots,n\}$, we define the
disk:
\[
U_r^{(j)}=\{P\in U^{(j)}\;|\; |z^{(j)}(P)|<r\}\, .
\]
We fix a radius $0<r_1<1$ once and for all. Our aim is to prove the
following result.

\begin{supthm}\label{merkl_thm1}
Assume that the open sets $U_{r_1}^{(j)}$ with $j$ in $\{1,\ldots,n\}$
cover~$X$. Next, assume that $c_1$ is a positive real number such that
for all $j$ in $\{1,\ldots,n\}$ we have:
\[
\mu\le c_1|dz^{(j)}\wedge d\bar z^{(j)}| \quad \text{on $U_1^{(j)}$}.
\]
Finally, assume that for all $j$ and $k$ in $\{1,\ldots,n\}$:
\[
\sup_{U_1^{(j)}\cap U_1^{(k)}}
\left|
\frac{dz^{(j)}}
{dz^{(k)}}
\right|\le M
\]
holds with some constant $M\ge 1$. Then for some positive constants
$c_7$, $c_9$, $c_{10}$ and~$c_{11}$, depending only on~$r_1$, we have,
for all $a$ in~$X$:
\begin{supeqn}\label{merkl_thm1_ineq1}
g_{a,\mu} \le n(c_{10}+c_1 c_{11}+c_7\log M) + \frac{\log 2}{2\pi}
\end{supeqn}
and, for all $j$ such that $a\in U_{r_1}^{(j)}$:
\begin{supeqn}\label{merkl_thm1_ineq2}
\begin{aligned}
& \lim_{x \to a} 
\left| 
g_{a,\mu}(x) - \frac{1}{2\pi}\log |z^{(j)}(x) - z^{(j)}(a) | 
\right| \le \\ 
& \qquad\le n(c_{10}+c_1 c_{11}+c_7\log M) + \frac{\log M}{2\pi} + c_9.
\end{aligned}
\end{supeqn}
\end{supthm}
We start by considering just one coordinate $z=z^{(j)}$ for a
fixed~$j$.  To simplify the notation in this section, we drop the
superscript $(j)$ in $U=U^{(j)}$, $z=z^{(j)}$, and so on. We fix three
radii $0<r_1<r_2<r_3<1$ once and for all.  The radii $r_2$ and $r_3$
should depend only on~$r_1$; e.g.  $r_2=(2r_1+1)/3$, $r_3=(r_1+2)/3$
is an admissible choice.  Furthermore, we fix a partition of unity:
let $\chi\colon X\to[0,1]$ be a $C^\infty$ function which is compactly
supported in the interior of $U_1$ with $\chi=1$ on $\ol{U_{r_2}}$,
and set $\chi^c=1-\chi$.  More specifically, we take
$\chi=\tilde\chi(|z|)$ on $U_1$ with a smooth function
$\tilde\chi\colon\RR\to[0,1]$ such that $\tilde\chi(r)=0$ for $r\ge
1-\epsilon$ with some $\epsilon>0$, and $\tilde\chi(r)=1$ for $r\le
r_2$.  The shape function $\tilde \chi$ may be taken independently of
$X$ and the choice of the coordinate $z$, only depending on~$r_2$.

We shall use the 2-norm of a (real valued) 1-form $\omega$ over a
measurable set $Y\subseteq X$ defined by:
\[
\|\omega\|_Y :=\left(\int_Y\omega\wedge\ast \omega\right)^{1/2}
=
\left(2i\int_Y\omega_{1,0}\wedge\omega_{0,1}\right)^{1/2},
\]
where $\omega=\omega_{1,0}+\omega_{0,1}$ is the decomposition of
$\omega$ in its components in $T_{(1,0)}X$ and~$T_{(0,1)}X$.  In the
case $Y=X$, we just write $\|\omega\|$ for~$\|\omega\|_X$.

Given $a$ and $b$ in $U_{r_1}$, we define the following function, having
logarithmic singularities in $a$ and~$b$:
\[
f_{a,b} :=\frac{1}{2\pi}
\log\left|
\frac{(z-z(a))(\ol{z(a)}z-1)}{(z-z(b))(\ol{z(b)}z-1)}\right|
\quad\text{on $U_1-\{a,b\}$.}
\]
Note that the singularities at $1/\ol{z(a)}$ and
$1/\ol{z(b)}$ do not lie within the unit disk. We note that:
\[
d{\ast} df_{a,b}=\delta_a-\delta_b
\]
holds on $U_1$ in the sense of distributions, and that $f_{a,b}$
fulfils Neumann boundary conditions on~$\partial U_1$.  
One can see this as follows.
The meromorphic function on $U_1$ given by:
\[
P\mapsto \frac{(z(P)-z(a))(\ol{z(a)}z(P)-1)}
{(z(P)-z(b))(\ol{z(b)}z(P)-1)}
=
\frac{(z(P)-z(a))(\frac{1}{z(P)}-\ol{z(a)})}{(z(P)-z(b))
(\frac{1}{z(P)}-\ol{z(b)})}
\]
takes positive real values on~$\partial U_1$.
Let $\Log$ denote the principal branch of
the logarithm. The function
\[
q_{a,b}(P):=
\frac{1}{2\pi}\Log\frac{(z(P)-z(a))(\ol{z(a)}z(P)-1)}
{(z(P)-z(b))(\ol{z(b)}z(P)-1)}
\]
is defined and holomorphic for $P$ in a neighbourhood of $\partial
U_1$, and it takes real values for $P\in \partial U_1$.  As a
consequence, the directional derivative of the imaginary part $\Im
q_{a,b}$ tangential to $\partial U_1$ vanishes on $\partial U_1$.
Using holomorphy, this implies that the directional derivative of the
real part $\Re q_{a,b}$ in normal direction to $\partial U_1$ vanishes
also on $\partial U_1$.  Using $\Re q_{a,b}(P)=f_{a,b}(P)$ for $P$ in
a neighbourhood of $U_1$, this proves the claimed Neumann boundary
conditions for~$f_{a,b}$.

Finally, for $a\in U_{r_1}$ and $P\in U_1-\{a\}$, we set:
\[
\label{merkl_def_la}
l_a(P):=\frac 1{2\pi} \chi(P) \log|z(P)-z(a)|,
\]
extended by~$0$ to~$X-\{a\}$. 

Our first step in the proof of Theorem \ref{merkl_thm1} is the
following key lemma.
\begin{suplem}
\label{merkl_lemma0}
For $a$ and $b$ in $U_{r_1}$, the supremum $\sup_X|g_{a,b}-l_a+l_b|$ is
bounded by a constant $c_2=c_4+c_1 c_5$, with $c_4,c_5$ depending only
on~$r_1$.
\end{suplem}
\begin{suprem} Note that $g_{a,b}-l_a+l_b$ has removable singularities
at $a$ and $b$, since the logarithmic singularities cancel.  The
constant $c_2$ is uniform in the choice of the Riemann surface~$X$,
and uniform in the choice of $a,b\in U_{r_1}$. The choice of the
coordinate $z$ influences $c_2$ only via the dependence of $c_1$ on
the choice of~$z$.  The radii $r_2$, $r_3$ and the shape function
$\tilde\chi$ are viewed as $r_1$-dependent parameters; this is why we
need not emphasise in the lemma that $c_2$ also depends on these
quantities.
\end{suprem}  

\begin{proof}
(of Lemma~\ref{merkl_lemma0}) We define the 2-form:
\[
u_{a,b} :=d{\ast} d(\chi^cf_{a,b}) \quad\text{on $U_1- \ol{U_{r_1}}$}
\]
and extend it by~$0$ to the whole surface~$X$.  Note that $u_{a,b}$ is
supported in $\ol{U_1}- U_{r_2}$, since $\chi^c$ varies only there,
and since $f_{a,b}$ is harmonic.  Consider the following variational
principle on square integrable 1-forms~$\omega$. We want to minimise
$\|\omega\|^2$ with the constraint:
\[
d{\ast}\omega = u_{a,b}
\]
in the sense of distributions.  Writing the constraint with test
functions, we see that the minimisation problem is taken over the
following closed affine linear subspace of $L^2(X,T^\ast X)$:
\[
V= \{\omega\in L^2(X,T^\ast X):\;-\int_Xdg\wedge \ast\omega =
\int_X gu_{a,b} \;\text{for all $g\in C^\infty(X)$}\}.
\]
The space $V$ is nonempty, since $\tilde\omega_{a,b}\in V$ holds for
the following $1$-form:
\[
\tilde{\omega}_{a,b}=
\left\{
\begin{array}{ll}
d (\chi^c f_{a,b})
&
\text{on $U_1- U_{r_1}$,}
\\0&\text{otherwise.}
\end{array}
\right.
\]
Indeed, using Stokes' theorem, we have:
\[
\begin{aligned}
-\int_Xdg\wedge \ast \tilde{\omega}_{a,b} 
&= -\int_{U_1} dg\wedge \ast \tilde{\omega}_{a,b} \nonumber \\
&=- \int_{\partial U_1} g\, {\ast} \tilde{\omega}_{a,b}
+\int_{U_1} g\,d{\ast} \tilde{\omega}_{a,b}.
\end{aligned}
\]
The first summand in the last expression vanishes by the Neumann boundary
conditions of $\chi^c f_{a,b}=f_{a,b}$ on~$\partial U_1$, and the
second summand equals:
\[
\int_{U_1} g\,d{\ast} \tilde{\omega}_{a,b}=\int_X g u_{a,b}
\]
by the definition of~$u_{a,b}$.

Our minimisation problem has a unique solution $\omega_{a,b}\in V$.
It fulfils:
\begin{supeqn}\label{eqn_merkl_1}
\int_X \omega_{a,b}\wedge \sigma=0
\quad\mbox{for all closed $C^\infty$ 1-forms $\sigma$.}
\end{supeqn}
Indeed: if $d\sigma=0$, then $\omega_{a,b}+t\ast \sigma\in V$ holds
for all $t\in\RR$, since $\omega_{a,b}\in V$ and
$d{\ast}(\ast\sigma)=-d\sigma=0$.  Thus:
\[
0= \left.
\frac d{dt}\|\omega_{a,b}+t\ast \sigma\|^2\right|_{t=0}
=
-2\int_X \omega_{a,b}\wedge \sigma.
\]
In particular:
\[
\int_X \omega_{a,b}\wedge dg=0
\]
for all $g\in C^\infty(X)$, i.e.\ $d\omega_{a,b}=0$ in the sense of
distributions.  Since $d{\ast} \omega_{a,b}=u_{a,b}$ and
$d\omega_{a,b}=0$, we get that $\omega_{a,b}$ is smooth.  This follows
from (hypo-)elliptic regularity, as treated in Corollary~4.1.2 on
page~101 in~\cite{Hoermander}. Precisely speaking, this corollary
treats only the case of a single partial differential
equation. However, as is mentioned at the end of Section~4.0 on
page~97 of the reference, the extension of the result to systems of
partial differential equations with as many equations as unknowns, as
needed here, follows trivially.  Then equation (\ref{eqn_merkl_1})
implies that $\omega_{a,b}$ is exact:
\[
\omega_{a,b}=d\tilde{g}_{a,b}
\]
for some $\tilde g_{a,b}\in C^\infty(X)$; see for example
\cite{Forster}, Corollary~19.13.  We normalise $\tilde g_{a,b}$ such
that:
\begin{supeqn}\label{eqn_merkl_2}
\int_X \tilde{g}_{a,b}\mu=0,
\end{supeqn}
to make it uniquely determined.

We set 
\[
h_{a,b}:=\tilde g_{a,b} +\chi f_{a,b}-\int_X \chi f_{a,b}\mu.
\]

We are now going to prove that $d{\ast} d h_{a,b}=\delta_a-\delta_b$.
We claim that $d{\ast} d (\chi f_{a,b})= -u_{a,b}+\delta_a-\delta_b$
holds.  We prove this equality separately on the three sets $X-
\Supp\chi$, $X- \Supp\chi^c$, and $U_1- \ol{U_{r_1}}$, which
cover~$X$.  The claimed equality holds on $X- \Supp\chi$, because both
sides vanish there.  It holds also on $X- \Supp\chi^c$, because there
$u_{a,b}=0$ and $d{\ast} d (\chi f_{a,b})=d{\ast} d f_{a,b}=
\delta_a-\delta_b$ are valid.  Finally, on $U_1- \ol{U_{r_1}}$, the
function $f_{a,b}$ is harmonic, which implies that $d{\ast} d (\chi
f_{a,b})=-d{\ast} d (\chi^c f_{a,b})=-u_{a,b}$ on this annulus, which
neither contains $a$ nor~$b$. Thus the claim $d{\ast} d (\chi
f_{a,b})= -u_{a,b}+\delta_a-\delta_b$ holds in all cases.

Combining this with the fact $d{\ast} d \tilde g_{a,b}= u_{a,b}$,
we conclude 
\[
d{\ast} dh_{a,b}=d{\ast} d \tilde g_{a,b}+d{\ast} d (\chi f_{a,b})
=u_{a,b}-u_{a,b}+\delta_a-\delta_b
\]
and thus
\[
d{\ast} dh_{a,b}=\delta_a-\delta_b.
\]

Furthermore, using the normalisation $\int_X\mu=1$ and $\int_X\tilde
g_{a,b}\mu=0$, we observe
\[
\int_X h_{a,b}\mu=0.
\]
Because $g_{a,b}$ is uniquely characterised by its properties $d{\ast}
dg_{a,b}=\delta_a-\delta_b$ and $\int_X g_{a,b}\mu=0$, we conclude
$g_{a,b}=h_{a,b}$.  Thus, we have shown
\[
g_{a,b}
=
\tilde g_{a,b} +\chi f_{a,b}-\int_X \chi f_{a,b}\mu.
\]

The function:
\[
g_{a,b}^{(1)} = \tilde g_{a,b}+\chi f_{a,b}
\]
is harmonic on $X-\{a,b\}$, and:
\[
g_{a,b}^{(2)}= \tilde g_{a,b}-\chi^c f_{a,b}
\]
is harmonic on~$U_1$; in particular both functions are harmonic on the
annulus $A:=U_1- \ol{U_{r_2}}$.  Now for every harmonic function $g$
on~$A$, we have a bound:
\[
\max_{\partial U_{r_3}}{g}-\min_{\partial U_{r_3}}{g}
\le c_3 \|dg\|_{A}
\]
with some positive constant $c_3$ depending only on $r_2$ and~$r_3$;
note that the circle $\partial U_{r_3}$ is relatively compact in the
annulus~$A$. We bound $\|dg_{a,b}^{(2)}\|_A$ from above:
\[
\|dg_{a,b}^{(2)}\|_A\le \|d\tilde g_{a,b}\|_A+\|d(\chi^c f_{a,b})\|_A.
\]
We estimate the first summand as follows, using that
$\omega_{a,b}=d\tilde g_{a,b}$ solves the above variational problem:
\[
\|d\tilde g_{a,b}\|_A \le \|d\tilde g_{a,b}\| = \|\omega_{a,b}\|\le 
\|\tilde \omega_{a,b}\| = \|\tilde \omega_{a,b}\|_A = 
\|d(\chi^c f_{a,b})\|_A;
\]
we used that $\tilde\omega_{a,b}$ is supported in~$A$.  Thus we have:
\[ 
\|dg_{a,b}^{(2)}\|\le 2\|d(\chi^c f_{a,b})\|_A,
\]
which is bounded by a constant, uniformly in $a$ and $b$ in~$U_{r_1}$.

This also allows us to estimate~$g_{a,b}^{(1)}$: on~$A$, we know
$g_{a,b}^{(1)}=g_{a,b}^{(2)}+f_{a,b}$, hence,
\[
\|dg_{a,b}^{(1)}\|_A\le\|dg_{a,b}^{(2)}\|_A+\|df_{a,b}\|_A
\le 2\|d(\chi^c f_{a,b})\|_A+\|df_{a,b}\|_A.
\]
Both summands on the right hand side are bounded by constants, only
depending on $r_1$ and~$r_2$, but uniformly in $a$ and $b$
in~$U_{r_1}$.  To summarise, we have shown that:
\[
\max_{\partial U_{r_3}}{g_{a,b}^{(j)}}-\min_{\partial U_{r_3}}{g_{a,b}^{(j)}}
\]
($j=1,2$) are uniformly bounded by a constant depending only on $r_1$,
$r_2$, and~$r_3$.  However, $g_{a,b}^{(1)}$ is harmonic on
$X-U_{r_3}$, and $g_{a,b}^{(2)}$ is harmonic on $\ol{U_{r_3}}$,
which both have the same boundary~$\partial U_{r_3}$.  Thus, by the
maximum principle:
\[
\max_{X- U_{r_3}}{g_{a,b}^{(1)}}-\min_{X- U_{r_3}}{g_{a,b}^{(1)}} =
\max_{\partial U_{r_3}}{g_{a,b}^{(1)}} -
\min_{\partial U_{r_3}}{g_{a,b}^{(1)}}
\]
and:
\[
\max_{\ol{U_{r_3}}}{g_{a,b}^{(2)}} -
\min_{\ol{U_{r_3}}}{g_{a,b}^{(2)}} =
\max_{\partial U_{r_3}}{g_{a,b}^{(2)}} -
\min_{\partial U_{r_3}}{g_{a,b}^{(2)}}.
\]
Furthermore, $\max_{X- U_{r_3}}|\chi f_{a,b}|$ and
$\max_{\ol{U_{r_3}}}|\chi^c f_{a,b}|$ are bounded, uniformly in $a$
and $b$ in~$U_{r_1}$, by a constant only depending on $r_1$ and~$r_3$.
Using $\tilde g_{a,b}=g_{a,b}^{(1)}-\chi f_{a,b}$ on $X- U_{r_3}$ and
$\tilde g_{a,b}=g_{a,b}^{(2)}+\chi^c f_{a,b}$ on $\ol{U_{r_3}}$, we
conclude that $\max_X\tilde g_{a,b}-\min_X\tilde g_{a,b}$ is bounded
on $X= (X- U_{r_3}) \cup \ol{U_{r_3}}$ by a constant $c_6$ only
depending on the radii $r_1$, $r_2$ and~$r_3$.  Using the
normalisation condition~(\ref{eqn_merkl_2}), we know that:
\[
\max_X\tilde g_{a,b}\ge 0\ge \min_X\tilde g_{a,b}
\]
holds; thus:
\[
\max_X|\tilde g_{a,b}|\le \max_X\tilde g_{a,b}-\min_X\tilde g_{a,b}
\]
is also bounded by the same constant.

From this we get a bound for:
\[
g_{a,b}-\chi f_{a,b} = \tilde g_{a,b} -\int_X \chi f_{a,b}\mu.
\]
Indeed, we estimate:
\[
\left|\int _X\chi f_{a,b}\mu\right| \le 
\int_{U_1}|f_{a,b}|\mu \le c_1\int_{U_1}|f_{a,b}\,dz\wedge d \bar z|,
\]
which is uniformly bounded for $a,b\in U_{r_1}$ by a constant $c_1c_5
$ with $c_5$ depending only on~$r_1$; note that the logarithmic
singularities are integrable.  Combining the bounds for $\max_X|\tilde
g_{a,b}|$ and $\left|\int _X\chi f_{a,b}\mu\right|$, we conclude that
$\sup_X|g_{a,b}-\chi f_{a,b}|$ is bounded by a constant $c_6+c_1 c_5$
with $c_6,c_5$ depending on~$r_1$.  Since:
\[
\sup_X|\chi f_{a,b} -l_a+l_b|
=
\frac{1}{2\pi}
\sup_{U_1}\left|\chi \log
\left|\frac{\ol{z(a)}z-1}{\ol{z(b)}z-1}\right|\right|
\]
is bounded, uniformly in $a,b\in U_{r_1}$ and~$X$, the key lemma
follows (with $c_4$ being the sum of $c_6$ and the uniform upper bound
last mentioned).
\end{proof}

\begin{proof}
(of Theorem~\ref{merkl_thm1}) Since we now work with varying
coordinates, we include again the superscript coordinate index~$(j)$
in the coordinate~$z^{(j)}$, its domain~$U^{(j)}$, but also in
$U_r^{(j)}$, $\chi^{(j)}$, and~$l_a^{(j)}$.
\begin{suplem}
\label{merkl_lemma1}
Consider two coordinates $z^{(j)}$ and~$z^{(k)}$, with $k$ and $j$ in
$\{1,\ldots,n\}$.  Assume that $x$ is in $U^{(j)}_{r_1}\cap U^{(k)}_{r_1}$
and that $y$ is in $U^{(j)}_{r_2}$ with
$|z^{(j)}(y)-z^{(j)}(x)|<(r_2-r_1)/M$. Then $y$ is in~$U^{(k)}_{r_2}$.
\end{suplem}
\begin{proof}
The intersection $U^{(j)}_{r_1}\cap U^{(k)}_{r_1}$ is an open
neighbourhood of~$x$.  Assume that there exists
$y\in\ol{U}_{r_2}^{(j)}$ with $|z^{(j)}(y)-z^{(j)}(x)|<(r_2-r_1)/M$
and $y\notin U^{(k)}_{r_2}$.  Then there is also such a point $y$ with
minimal distance $|z^{(j)}(y)-z^{(j)}(x)|$ from~$x$, since
$\ol{U}_{r_2}^{(j)}- U^{(k)}_{r_2}$ is compact.  For this point~$y$,
we conclude $y\in\partial U_{r_2}^{(k)} \subseteq \ol{U}_{r_2}^{(k)}$,
and the straight line from $x$ to $y$ in the $z^{(j)}$-coordinate is
contained in $\ol{U}_{r_2}^{(j)}\cap \ol{U}_{r_2}^{(k)}$.  By the mean
value theorem, we conclude $|z^{(k)}(y)-z^{(k)}(x)|\le M
|z^{(j)}(y)-z^{(j)}(x)|<r_2-r_1$, hence $|z^{(k)}(y)|<r_2$, since
$|z^{(k)}(x)|\le r_1$.  This contradicts~$y\in\partial U_{r_2}^{(k)}$.
\end{proof}
We choose a smooth partition of unity $\phi^{(j)}:X\to[0,1]$,
$j=1,\ldots,n$, such that $\phi^{(j)}$ is supported
in~$U_{r_1}^{(j)}$.  For $a\in X$, we set:
\[
h_a := \sum_j \phi^{(j)}(a) l_a^{(j)}.
\]
\begin{suplem} \label{merkl_lemma2}
Let $a\in U_{r_1}^{(k)}$, $y\in X$, $y\neq a$.  Then we have:
\[
l_a^{(k)}(y)\le \frac{\log 2}{2\pi}.
\]
\end{suplem}
\begin{proof}
This follows immediately from the definition of the
function~$l_a^{(k)}$, since $|z^{(k)}(y)-z^{(k)}(a)|\le 2$ whenever
$y\in\Supp(\chi^{(k)})$.
\end{proof}
\begin{suplem}
\label{merkl_lemma3}
For all $a,b\in X$ we have the inequality:
\[  
\sup_X|g_{a,b}-h_a+h_b| \le n(c_{10}+c_1 c_5+c_7\log M)
\]
with constants~$c_{10}$, $c_5$, and $c_7$ depending only on~$r_1$.
\end{suplem}
\begin{proof}
We first show for $a\in U_{r_1}^{(k)}\cap U_{r_1}^{(j)}$ that:
\[
\sup_X|l_a^{(k)}-l_a^{(j)}| \le
\frac{1}{2\pi} [\log M +|\log (r_2-r_1)|+\log 2].
\]
To prove this, let $y\in X$.  We distinguish 3 cases in order to prove
that $l_a^{(k)}(y)-l_a^{(j)}(y)$ is bounded from above by the right
hand side.

\paragraph{\bf case 1:} $y\in U_1^{(j)}$ with
$|z^{(j)}(y)-z^{(j)}(a)|< (r_2-r_1)/M$.  In particular, we have
$|z^{(j)}(y)|< |z^{(j)}(a)|+(r_2-r_1)/M\le r_2$ (recall that $M \geq
1$), hence $a,y\in U_{r_2}^{(j)}$.  Consequently, the straight line
$[a,y]^{(j)}$ from $a$ to $y$ in the $z^{(j)}$-coordinate is contained
in~$U_{r_2}^{(j)}$.  Then Lemma~\ref{merkl_lemma1} implies that
$[a,y]^{(j)}\subseteq U_{r_2}^{(k)}$. Using
$\chi^{(j)}(y)=\chi^{(k)}(y)=1$, we conclude by the mean value theorem
that:
\[
l_a^{(k)}(y)-l_a^{(j)}(y) =
\frac 1{2\pi}
\log\left|
\frac{z^{(k)}(y)-z^{(k)}(a)}{z^{(j)}(y)-z^{(j)}(a)}
\right|
\le
\frac{\log M}{2\pi},
\]
which is bounded by the right hand side.

\paragraph{\bf case 2:}
$y\notin U_1^{(j)}$. Then $l_a^{(j)}(y)=0$,
and we conclude, using Lemma~\ref{merkl_lemma2}, that:
\[
l_a^{(k)}(y)-l_a^{(j)}(y)=l_a^{(k)}(y)\le \frac{\log 2}{2\pi}.
\]

\paragraph{\bf case 3:}
$y\in U_1^{(j)}$ and $|z^{(j)}(y)-z^{(j)}(a)|\ge (r_2-r_1)/M$;
thus:
\[
l_a^{(k)}(y)-l_a^{(j)}(y)\le \frac{\log 2}{2\pi}-l_a^{(j)}(y)
\le \frac{1}{2\pi}(\log 2- \chi^{(j)}(y)\log[(r_2-r_1)/M]),
\]
which is also bounded by the right hand side.

The upper bound for $l_a^{(j)}(y)-l_a^{(k)}(y)$ in our claim is
obtained by exchanging~$j$ and~$k$. Thus the claim is proven.

We conclude:
\begin{supeqn}\label{merkl_eqn1}
\begin{aligned}
|h_a-l_a^{(j)}| & \le 
\sum_{k}\phi^{(k)}(a) |l_a^{(k)}-l_a^{(j)}| \le \\ 
& \le \frac{1}{2\pi}\left(\log M +|\log (r_2-r_1)|+\log 2\right).
\end{aligned}
\end{supeqn}
Combining this with Lemma~\ref{merkl_lemma0}, we conclude for $a,b\in
U_{r_1}^{(j)}$:
\[
\begin{aligned}
|g_{a,b}-h_a+h_b|
&\le 
|g_{a,b}-l_a^{(j)}+l_b^{(j)}|+|h_a-l_a^{(j)}|+|h_b-l_b^{(j)}| \\
&\le  c_{10}+c_1 c_5+c_7\log M
\end{aligned}
\]
with some constants $c_{10}, c_5, c_7$ depending only on~$r_1$ (a
possible choice is $c_7=(\log M)/\pi$ and
$c_{10}=(|\log(r_2-r_1)|+\log 2)/\pi+c_4$).

Finally, for general $a,b\in X$, we choose a finite sequence of points
$a=a_0,a_1,\ldots, a_m=b$ in~$X$ and indices $j_1,\ldots, j_m$ with
$m\le n$ and $a_{i-1}, a_i\in U_{r_1}^{(j_i)}$ for all~$i=1,\ldots,
m$.  Using:
\[
g_{a,b}=\sum_{i=1}^m g_{a_{i-1},a_i},
\]
we get by estimating:
\[
|g_{a,b}-h_a+h_b| \le
\sum_{i=1}^m |g_{a_{i-1},a_i}-h_{a_{i-1}}+h_{a_i}| \le 
n(c_{10}+c_1 c_5+c_7\log M) 
\]
the claim of the lemma.
\end{proof}

We define:
\[
h_\mu(x) :=\int_{b\in X} h_b(x)\,\mu(b),\qquad (x\in X).
\]
\begin{suplem}\label{merkl_lemma4}
We have:
\[
\sup_X |h_\mu|\le n c_1 c_8,
\]
with some universal constant~$c_8$. Furthermore, we have:
\[
\sup_{\substack{b,x\in X\\ b\neq x}}h_b(x)\le \frac{\log 2}{2\pi} .
\]
\end{suplem}
\begin{proof}
We observe first that for all $w\in \CC$ with $|w|\le 1$ the integral:
\[
\frac{1}{2\pi}\int_{|z|\le 1}|\log |z-w||\,|dz\wedge d\bar z| 
\]
is bounded from above by a universal constant~$c_8$.  We conclude that
for all $x\in X$ we have:
\[
\begin{aligned}
\int_{b\in U_{r_1}^{(j)}} |l_b^{(j)}(x)|\phi^{(j)}(b)\,\mu(b)
& \le 
c_1 \int_{b \in U_{r_1}^{(j)}} 
|l_b^{(j)}(x)|\,|dz^{(j)}\wedge d\ol{z^{(j)}}| \\
& \le c_1 c_8.
\end{aligned}
\]
Let $x\in X$. We get the first estimate:
\[
|h_\mu(x)|\le 
\sum_{j=1}^n \int_{U_1^{(j)}}|l_b^{(j)}(x)|\phi^{(j)}(b)\,\mu(b) 
\le n c_1 c_8.
\]
Finally, the second estimate follows from Lemma~\ref{merkl_lemma2}:
\[
h_b= \sum_{j=1}^n \phi^{(j)} l_b^{(j)}\le \frac{\log 2}{2\pi} ,
\]
as required.
\end{proof}
\begin{supprop}\label{merkl_prop1}
For some positive constants $c_{10}$, $c_7$, and~$c_{11}$ that depend
only on~$r_1$ we have, uniformly in $a$ and $x \neq a$ on~$X$:
\[
|g_{a,\mu}(x)- h_a(x)| \le n(c_{10} +c_1 c_{11} +c_7 \log M) .
\] 
\end{supprop}
\begin{proof}
Indeed, averaging Lemma~\ref{merkl_lemma3} over~$b$ with respect
to~$\mu$, we obtain:
\[
\sup_X|g_{a,\mu}-h_a+h_\mu| \le n(c_{10}+c_1 c_5+c_7\log M).
\]
By Lemma~\ref{merkl_lemma4}, one has $|h_\mu|\le n c_1 c_8$. Combining
gives what we want (we can take $c_{11}= c_5+c_8$).
\end{proof}
\begin{supprop}\label{merkl_prop2}
Let $c_{10}$, $c_7$, and~$c_{11}$ be as
in~Proposition~\ref{merkl_prop1}, and let $a$ be in~$X$. Then
$\lim_{x\to a}|g_{a,\mu}(x)- h_a(x)|$ exists, and we have:
\[
\lim_{x\to a}|g_{a,\mu}(x)- h_a(x)| \le n(c_{10}+c_1 c_{11}+c_7\log M).
\]
\end{supprop}
\begin{proof}
The functions $g_{a,\mu}$ and $h_a$ have the same logarithmic
singularity at~$a$; hence the limit exists. The estimate then follows
from Proposition~\ref{merkl_prop1}.
\end{proof}
We can now finish the proof of Theorem~\ref{merkl_thm1}. We have
seen in~(\ref{merkl_eqn1}) that:
\[
|h_a-l_a^{(j)}|\le 
\frac{1}{2\pi}\left(\log M +|\log (r_2-r_1)|+\log 2\right) .
\]
Combining this with Proposition~\ref{merkl_prop2} and using the
definition of~$l_a^{(j)}$ gives the second estimate of the theorem. As
to the first estimate, using:
\[
g_{a,\mu}\le h_a+|g_{a,\mu}-h_a| 
\]
we obtain it by applying the upper bound for $h_a$ in
Lemma~\ref{merkl_lemma4} and the upper bound for $|g_{a,\mu} - h_a|$
in Proposition~\ref{merkl_prop1}. This ends the proof of
Theorem~\ref{merkl_thm1}.
\end{proof}

%% -- Arakelov II (more about modular curves)
\chapter{Bounds for Arakelov invariants of modular curves}
\label{chap_bnd_height}

\author{B. Edixhoven and R. de Jong}

\bigskip

\bigskip

% authors: Bas and Robin

In this chapter, we give bounds for all quantities on the right hand
side in the inequality in Theorems~\ref{thm_bound_hbQ}
and~\ref{maininequality}, in the context of the modular curves
$X_1(5l)$ with $l>5$ prime, using the upper bounds for Green functions
from the previous chapter. The final estimates are given in the last
section.

\section{Bounding the height of~$X_1(\protect\lowercase{pl})$}
\label{sec_bnd_height_X}
As before, for $l>5$ prime, we let $X_l$ be the modular
curve~$X_1(5l)$, over a suitable base that will be clear from the
notation. We let $g_l$ denote the genus of~$X_l$; we have $g_l>1$. A
model $X_{l,\ZZ}$ is given by \cite{Katz-Mazur}, as well as a model
$X_{l,\ZZ[\zeta_{5l}]}$ that is semi-stable; see
Chapter~\ref{chap_descr_X15l}. The aim of this section is to prove a
suitable bound for the stable Faltings height of~$X_l$
(see~\ref{def_Faltings_height}). We will in fact give such a bound for
the modular curves $X_1(pl)$ with $p$ and $l$ distinct primes. Before
we get to that, we prove some intermediate results, that will also be
important in the next section.

\begin{lem} \label{lem_bnd_coefficients}
Let $N\geq1$ be an integer, and let:
\[
\calB_2(N) := 
\coprod_{M|N}\coprod_{d|(N/M)} B_{N,M,d}^*S_2(\Gamma_1(M))^\new
\]
be the basis of $S_2(\Gamma_1(N))$ obtained from newforms of levels
dividing~$N$ as explained in~(\ref{eqn_basis}).  Let $f=\sum_{n\geq
  1}a_n(f) q^n$ be an element of~$\calB_2(N)$. Then we have for all
$n\geq1$:
\[
|a_n(f)| \leq 2n.
\]
\end{lem}
\begin{proof}
As $a_n(B_{N,M,d}^*f)=a_{n/d}(f)$ (see~\ref{eqn_q-exp_basis}), it
suffices to treat the case that $f$ is a newform of some level~$M$
dividing~$N$.  We use the Weil bounds on the $a_p(f)$ for all
primes~$p$. We recall from Section~1.8 of~\cite{Deligne-Serre} that we
have an equality of formal Dirichlet series:
\[
\sum_{n\geq 1} a_n(f)n^{-s} = 
\prod_{p|M}(1-a_p(f)p^{-s})^{-1}
\prod_{p\nmid M} (1-\alpha_p p^{-s})^{-1}(1-\beta_p p^{-s})^{-1}
\]
with the following properties. For $p\nmid M$ we have $|\alpha_p| =
|\beta_p| = \sqrt{p}$. For $p|M$ we have:
\[
\left\{
\begin{aligned}
a_p(f) & = 0&  &\text{if $p^2|M$},\\
a_p(f) & = 0&  &
\text{if $\eps_f$ factors through $(\ZZ/(M/p)\ZZ)^\times$},\\
|a_p(f)| & = p^{1/2}&  &
\text{if $\eps_f$ does not factor through $(\ZZ/(M/p)\ZZ)^\times$},\\
|a_p(f)| & = 1&  &
\text{if $p^2|M$ and $\eps_f$ factors through $(\ZZ/(M/p)\ZZ)^\times$}.
\end{aligned}
\right.
\]
Using that:
\[
(1-a_p(f) p^{-s})^{-1} = \sum_{k\geq 0}a_p(f)^kp^{-sk},
\]
and that:
\[
(1-\alpha_p p^{-s})^{-1}(1-\beta_p p^{-s})^{-1} = 
(\sum_{k\geq 0}\alpha_p^kp^{-sk}) (\sum_{k\geq 0}\beta_p^kp^{-sk})
\]
we find that for arbitrary $n$ we have $|a_n(f)| \leq
\sigma_{0,M}(n)\sqrt{n}$, where $\sigma_{0,M}(n)$ is the number of
positive divisors of~$n$ that are prime to~$M$, and a simple estimate
leads to $|a_n(f)| \leq 2n$.
\end{proof}
The following lemma states a very well known lower bound for the
Petersson norm of a normalised cuspform.
% reference ??
\begin{lem}\label{eqn_lowerb_pet_norm}
Let $N\geq1$ and let $\omega=fdq/q$ be the holomorphic $1$-form on
$X_1(N)(\CC)$ attached to a cusp form $f=\sum_na_n(f)q^n$
in~$S_2(\Gamma_1(N))$ with $a_1(f)=1$. Then we have:
\[
\|\omega\|^2 = \frac{i}{2} \int_{X_1(N)} \omega \wedge
\ol{\omega} \geq \pi e^{-4\pi}.
\]
\end{lem}
\begin{proof}
We have $\omega=\sum_{n\geq1} a_n(f) q^n dq/q$ in the coordinate
$q=e^{2\pi i z}$, where $z$ is the standard coordinate on the upper
half plane~$\HH$. If we let $x$ and $y$ be the real and imaginary
parts of~$z$ we have:
\[
\frac{i}{2} \omega \wedge \ol{\omega} = 4\pi^2|f|^2\, dx\,dy
\]
Let $F$ be the region in~$\HH$ given by the conditions $|x|<1/2$ and
$y>1$.  Then:
\[
\begin{aligned} 
\|\omega\|^2 & \geq \int_F 4\pi^2 |f(z)|^2 dx dy \\
& = 4\pi^2 \sum_{m,n\geq 1} a_m(f) \ol{a_n(f)} 
\int_{-1/2}^{1/2} e^{2\pi i (m-n)x} \int_1^\infty e^{-2\pi(m+n)y}dy \\ 
&= 4\pi^2 \sum_{n\geq 1} |a_n(f)|^2 e^{-4\pi n}/{4\pi n} \, .
\end{aligned} 
 \]
From the first term (note that $a_1(f)=1$) we obtain $\|\omega\|^2
\geq \pi e^{-4\pi}$.
\end{proof}
We now specialise to a slightly less special case than our
curves~$X_l$: the curves $X_1(pl)$ with $p$ and~$l$ two distinct prime
numbers. We call an Atkin-Lehner basis for $\Omega^1(X_1(pl))$ any
basis of $\Omega^1(X_1(pl))$ given by an ordering of the set
$\calB_2(pl)$. We start by describing, in a notation that is slightly
different from the one used in~(\ref{eqn_degen_maps}), the degeneracy
maps that are used for the definition of~$\calB_2(pl)$. This time, we
call them source and target maps:
\[
\left\{
\begin{aligned}
s_l  \colon X_1(pl) \to X_1(p), \quad & (E,P,L)\mapsto (E,P)\\
t_l   \colon X_1(pl) \to X_1(p), \quad & (E,P,L)\mapsto (E/\ld L\rd,P)\\
s_p   \colon X_1(pl) \to X_1(l), \quad & (E,P,L)\mapsto (E,L)\\
t_p   \colon X_1(pl) \to X_1(l), \quad & (E,P,L)\mapsto (E/\ld P\rd,L)
\end{aligned}
\right.
\]
where $(E,P,L)$ denotes an elliptic curve~$E$ with a point~$P$ of
order~$p$ and a point~$L$ of order~$l$. Note that $s_l$ and $t_l$ have
degree $l^2-1$, and that $s_p$ and $t_p$ have degree~$p^2-1$.  For any
integer $M\geq1$ we denote by $\Omega^1(X_1(M))^\new$ the set of
holomorphic 1-forms in $\Omega^1(X_1(M))$ of the form $fdq/q$ with $f$
in~$S_2(X_1(M))^\new$. Our next goal is to get information on the Gram
matrix of an Atkin-Lehner basis of~$\Omega^1(X_1(pl))$. As described
above, the contribution to $\Omega^1(X_1(pl))$ of each $f$ in
$S_2(\Gamma_1(pl))^\new$ is the subspace $\CC fdq/q$. The contribution
of an $f$ in~$S_2(\Gamma_1(p))^\new$ is the 2-dimensional space
generated by $s_l^*fdq/q$ and $t_l^*fdq/q$, and, of course, each $f$
in $S_2(\Gamma_1(l))^\new$ contributes the 2-dimensional space
generated by $s_p^*fdq/q$ and $t_p^*fdq/q$.

\begin{lem} \label{innerproducts}
For $f$ in~$S_2(\Gamma_1(l))^\new$ and $\omega=fdq/q$ we have:
\[
\begin{aligned}
\ld s_p^* \omega, s_p^* \omega \rd & = (p^2-1) \|\omega\|^2 \, , \\
\ld t_p^* \omega, t_p^* \omega \rd & = (p^2-1) \|\omega\|^2 \, , \\
\ld s_p^* \omega, t_p^* \omega \rd & = (p-1)\ol{a_p(f)} \|\omega\|^2 .
\end{aligned}
\]
We have similar equalities with $p$ and $l$ switched.
\end{lem}
\begin{proof} The first two equalities are clear. As to the latter,
note first that:
\begin{multline*}
\ld s_p^*\omega,t_p^*\omega\rd = 
\frac{i}{2}\int_{X_1(pl)}s_p^*\omega\wedge\ol{t_p^*\omega} = 
\frac{i}{2}\int_{X_1(l)}s_{p,*}(s_p^*\omega\wedge\ol{t_p^*\omega}) = \\
=\frac{i}{2}\int_{X_1(l)}\omega\wedge\ol{s_{p,*}t_p^*\omega} \, .
\end{multline*}
Next note that $s_p\colon X_1(pl)\to X_1(l)$ and $t_p\colon X_1(pl)\to
X_1(l)$ factor through the forget map $X_1(pl)\to X_1(l;p)$ where the
latter curve corresponds to the moduli problem $(E,P,G)$ with $P$ of
order~$l$ and $G$ a subgroup of order~$p$. This forget map has degree
$p-1$, and the correspondence on $X_1(l)$ induced by $X_1(l;p)$ is the
standard Hecke correspondence~$T_p$. We find that $s_{p,*}t_p^*\omega
= (p-1)T_p^*\omega$. By the standard relation between eigenvalues and
$q$-coefficients
% reference ???
we have $T_p^*\omega = a_p(f)\omega$, so finally:
\[
\ld s_p^*\omega,t_p^*\omega\rd =
\frac{i}{2}\int_{X_1(l)}\omega\wedge\ol{(p-1)T_p^*\omega} = 
(p-1)\ol{a_p(f)}\|\omega\|^2
\]
as required. 
\end{proof}
\begin{cor} \label{gram} 
Let $p$ and $l$ be two distinct primes. The structure of the Gram
matrix $(\ld \omega_i , \omega_j \rd)_{i,j}$ of holomorphic 1-forms
attached to an Atkin-Lehner basis for $\Omega^1(X_1(pl))$ is as
follows. Two subspaces attached to distinct elements of the union of
$S_2(\Gamma_1(pl))^\new$, $S_2(\Gamma_1(l))^\new$ and
$S_2(\Gamma_1(p))^\new$ are orthogonal to each other, hence the Gram
matrix decomposes into blocks corresponding to these subspaces. The
contribution of an element $f$ in~$S_2(\Gamma_1(pl))^\new$ is the
1-by-1 block~$\|fdq/q\|^{2}$. The contribution of an element $f$
in~$S_2(\Gamma_1(l))^\new$ is the 2-by-2 block:
\[
(p-1)\|fdq/q\|^2
\left(
\begin{matrix}
p+1 & \ol{a_p(f)} \\ a_p(f) & p+1
\end{matrix}
\right) ,
\]
where the norm $\|fdq/q\|^2$ is taken on $X_1(l)$). The contribution
of an element $f$ in~$S_2(\Gamma_1(p))^\new$ is the 2-by-2 block:
\[
(l-1)\|fdq/q\|^2
\left(
\begin{matrix}
l+1 & \ol{a_l(f)} \\ a_l(f) & l+1
\end{matrix}
\right) ,
\]
where the norm $\|fdq/q\|^2$ is taken on~$X_1(l)$.
\end{cor} 
\begin{cor} \label{determinantgram}
The determinant of the Gram matrix of the holomorphic 1-forms
attached to an Atkin-Lehner basis for $\Omega^1(X_1(pl))$ is bounded
below by $(\pi e^{-4\pi})^g$.
\end{cor}
\begin{proof}
By the Weil-Ramanujan-Deligne bounds (or, in this case, the Weil
bounds, as the weight of the modular forms here is two), the
determinant of a 2-by-2 block as in Corollary~\ref{gram} is bounded
below by $\|fdq/q\|^4$. We obtain our corollary by invoking
Lemma~\ref{eqn_lowerb_pet_norm}.
\end{proof}
\begin{cor} \label{expressionmu}
The Arakelov (1,1)-form $\mu$ on $X_1(pl)$ is given by:
\begin{multline*} 
\mu = \frac{i}{2g}\sum_\omega 
\frac{\omega\wedge\ol{\omega}}{\|\omega\|^2} + 
\frac{i}{2g}\sum_\omega \left(
\frac{(p+1)s_p^*(\omega\wedge\ol{\omega}) + 
(p+1)t_p^*(\omega\wedge\ol{\omega})}
{(p-1)\|\omega\|^2((p+1)^2-|a_p(f_\omega)|^2)} \right. \\
\left. - \frac{a_p(f_\omega)s_p^*\omega\wedge\ol{t_p^*\omega} + 
\ol{a_p(f_\omega)}t_p^*\omega\wedge\ol{s_p^*\omega}}
{(p-1)\|\omega\|^2((p+1)^2-|a_p(f_\omega)|^2)}
\right)\\
+ 
\frac{i}{2g}\sum_\omega \left(
\frac{(l+1)s_l^*\omega\wedge\ol{s_l^*\omega} + 
(l+1)t_l^*\omega\wedge\ol{t_l^*\omega}}
{(l-1)\|\omega\|^2((l+1)^2-|a_l(f_\omega)|^2)} \right. \\
\left. - \frac{a_l(f_\omega)s_l^*\omega\wedge\ol{t_l^*\omega} + 
\ol{a_l(f_\omega)}t_l^*\omega\wedge\ol{s_l^*\omega}}
{(l-1)\|\omega\|^2((l+1)^2-|a_l(f_\omega)|^2)}
\right)\, ,
\end{multline*} 
with the first sum running over $\Omega^1(X_1(pl))^\new$, 
the second sum running over $\Omega^1(X_1(l))^\new$, the third sum
running over $\Omega^1(X_1(p))^\new$, and where $f_\omega$ is defined
by $\omega=f_\omega dq/q$. 
\end{cor}
\begin{proof}
Consider first an arbitrary compact Riemann surface $X$ and let
$\omega=(\omega_1,\ldots,\omega_g)$ be an arbitrary basis
of~$\Omega^1(X)$. Let $a:=\ld \omega,\omega\rd$ be the $g$-by-$g$
matrix given by $a_{i,j} = \ld\omega_i,\omega_j\rd$. Note that
$\ol{a}=a^t$. Let $b=a^{t,-1}$, the inverse of the transpose
of~$a$. Then we claim that the Arakelov $(1,1)$-form on $X$ can be
written as:
\[
\mu = \frac{i}{2g}\sum_{i,j} b_{i,j}\omega_i\wedge\ol{\omega_j}.
\]
To see this, note that for $\omega$ an orthonormal basis this is the
correct expression, and that changing to $\omega'=\omega{\cdot}g$ with
any invertible~$g$ does not change~$\mu$, as one may directly
calculate. 

In our case, the basis $(\omega_1,\ldots,\omega_g)$ that we take is an
Atkin-Lehner basis. Using Corollary~\ref{gram} one obtains the
expression that we gave.
\end{proof}
We remark that Abbes and Ullmo have determined the Arakelov
$(1,1)$-form on $X_0(n)$ for all square free $n\geq1$ such that
$X_0(n)$ has genus at least one in~\cite{Abbes-Ullmo1}. It should not
be hard to generalise their result to~$X_1(n)$ for square free~$n$.

Now we arrive at the main result of this section. We recall that the
Faltings height of a curve, and its stable or absolute version, have
been briefly described in~(\ref{def_abs_Faltings_height}).
\begin{thm} \label{thm_height}
For the stable Faltings height of~$X_1(pl)$, for distinct prime
numbers~$p$ and~$l$, one has:
\[
h_\abs(X_1(pl)) = O((pl)^2\log(pl)).
\]
\end{thm}
\begin{proof} 
This proof is an adaptation of an argument in Section~5
of~\cite{Coleman-Edixhoven} where the case $X_0(p)$ with $p$ prime was
treated. We may and do assume that $X_1(pl)$ has genus at least one.

We start with a general observation. For $X_K$ a curve over a number
field, and $K\to L$ a finite extension, we claim that:
\[
[L:\QQ]^{-1}h_L(X_L) \leq [K:\QQ]^{-1}h_K(X_K).
\]
This inequality simply results from the fact that for the Néron models
of the Jacobians the identity morphism on the generic fibres extends
to a morphism:
\[
(J_{O_K})_{O_L} \lto J_{O_L}.
\]

For $n$ a positive integer, we let $X_\mu(n)_\QQ$ denote the modular
curve corresponding to elliptic curves with an embedding
of~$\mu_n$. The reason for considering this variant of~$X_1(n)$ is
that the cusp $\infty$ of $X_\mu(n)$ is $\QQ$-rational. Of
course, $X_1(n)_\QQ$ and $X_\mu(n)_\QQ$ become isomorphic
over~$\QQ(\zeta_{n})$, and therefore we have, for all~$n$:
\[
h_\abs(X_1(n)_\QQ) = h_\abs(X_\mu(n)_\QQ).
\]
For more details about these $X_\mu(n)$ we refer to sections~9.3
and~12.3 of~\cite{Diamond-Im}. 

The general observation above gives:
\[
h_\abs(X_\mu(n)_\QQ) \leq h_\QQ(X_\mu(n)_\QQ). 
\]
Because of this, it suffices to establish the bound of the
theorem for the $h_\QQ(X_\mu(pl)_\QQ)$. 

Let $p$ and $l$ be given, and let $X$ be the model over~$\ZZ$ of
$X_\mu(pl)\QQ$ obtained by normalisation of the $j$-line~$\PP^1_\ZZ$
in the function field of~$X_\mu(pl)\QQ$. As $X$ is proper over~$\ZZ$,
the $\QQ$-rational point $\infty$ extends to an element~$\infty$
in~$X(\ZZ)$, which is known to lie in the open part $X^\sm$ of $X$
where the structure morphism to~$\Spec(\ZZ)$ is smooth,
see~\cite{Diamond-Im}. 
%% make the reference more precise.
In terms of the Tate curve over~$\ZZ((q))=\ZZ((1/j))$, the cusp
$\infty$ is the immersion of~$\mu_n$, over~$\ZZ$, in the $n$-torsion
of the Tate curve (see Sections~8.6--8.11 of~\cite{Katz-Mazur}).

We let $J$ be the Néron model over~$\ZZ$ of the Jacobian of the
curve~$X_\QQ$. Then, by the defining property, the embedding of
$X_\QQ$ into~$J_\QQ$ that sends $\infty$ to~$0$ extends to a morphism
from $X^\sm$ to~$J$. This morphism induces via pullback of
differential forms a morphism from $\Cot_0(J)$ to $S(\ZZ)$, the
sub-$\ZZ$-module of $\Omega^1(X_\QQ)$ of forms whose $q$-expansion
at~$\infty$ has coefficients in~$\ZZ$ (see
around~(\ref{eqn_TS_pairing})). As $\Cot_0(J)$ and $S(\ZZ)$ are both
$\ZZ$-structures on~$\Omega^1(X(\CC))$, we have (see
around~(\ref{def_Faltings_height})):
\[
\begin{aligned}
h_\QQ(X_\QQ) & = \deg(\bigwedge^g 0^*\Cot_0(J)) = 
-\log\Vol((\RR\otimes\Cot_0(J))/\Cot_0(J)) \\
& \leq -\log\Vol(\RR\otimes S(\ZZ)/S(\ZZ)),
\end{aligned}
\]
where the volume form on $\RR\otimes\Cot_0(J)$ comes from integration
over~$J(\CC)$, and that on $\RR\otimes S(\ZZ)$ from integration
over~$X(\CC)$. 

Let $\TT\subset\End(J)$ be the Hecke algebra, generated by all~$T_i$,
$i\geq 1$, and the $\ld a\rd$, $a$ in~$(\ZZ/pl\ZZ)^\times$. We have a
perfect pairing (see~(\ref{eqn_TS_pairing})): 
\[
\TT \times S(\ZZ) \to \ZZ, \quad (t,\omega) \mapsto a_1(t\omega).
\]
Using the duality we can write: 
\[
-\log \Vol(\RR\otimes S(\ZZ)/S(\ZZ)) = \log\Vol(\RR\otimes\TT/\TT)
\]
where the volume form on $\RR\otimes\TT$ is dual to the one
on~$\RR\otimes S(\ZZ)$.  Now consider an Atkin-Lehner basis
$(\omega_1,\ldots,\omega_g)$ of~$\Omega^1(X)_\CC$.  Let $\Vol'$ denote
the volume with respect to the volume form on $\RR\otimes\TT$ induced
by the one on $\CC\otimes S(\ZZ)$ for which the basis
$(\omega_1,\ldots,\omega_g)$ is an orthonormal basis. Then we have:
\[
\log\Vol(\RR\otimes \TT/\TT) = 
\log \Vol'(\RR \otimes \TT/\TT) - 
\frac{1}{2}\log\det(\ld \omega,\omega \rd) \, , 
\] 
where $\ld \omega,\omega\rd$ is the matrix whose $(i,j)$-coefficient
is~$\ld \omega_i,\omega_j \rd$.  By Corollary~\ref{determinantgram} we
have:
\[
- \log \det(\ld\omega,\omega\rd) \leq g(4\pi-\log\pi) = O((pl)^2).
\]
%% reference for g=O((pl)^2) ??????????????????
It remains to bound $\log \Vol'(\RR \otimes \TT/\TT)$.  Let $\Gamma$
be the set of integers $i \geq 1$ such that there exists an $\omega$
in $\Omega^1(X(\CC))$ with a zero of exact order $i-1$ at~$\infty$.
Then $\Gamma$ is the set of integers $i\geq 1$ such that
$h^0(X(\CC),\Omega^1(-i\infty))$ is strictly less (and hence exactly
one less) than $h^0(X(\CC),\Omega^1((-i+1)\infty)$. As
$h^0(X(\CC),\Omega^1)=g$, and $h^0(X(\CC),\Omega^1(-(2g-1)\infty))=0$,
there are exactly $g$ such integers, and we can write
$\Gamma=\{i_1,\ldots,i_g\}$ with:
\[
1= i_1 < i_2 < \ldots < i_g \leq 2g-1.
\]
Under the pairing between $\TT$ and~$S(\ZZ)$, each Hecke operator
$T_i$ is sent to the element $\omega\mapsto a_i(\omega)$ of the dual
of~$S(\ZZ)$, where the $a_i(\omega)$ are given by the $q$-expansion:
\[
\omega = \sum_{i\geq 1}a_i(\omega)q^i\, (dq/q) 
= \sum_{i\geq1}a_i(\omega) q^{i-1}\,dq.
\]
It follows that the elements $T_{i_1},\ldots,T_{i_g}$ of the free
$\ZZ$-submodule $\TT$ are linearly independent. Hence $\TT'$, the
submodule of~$\TT$ generated by these $T_{i_j}$ has finite index. We
thus find:
\[
\log\Vol(\RR\otimes \TT/\TT) \leq 
\log \Vol'(\RR \otimes \TT'/\TT') -
\frac{1}{2}\log\det(\ld\omega,\omega\rd) \, .
\]
Now we have $g=r_1+2r_2$, where $r_1$ is the number of elements of our
basis $(\omega_1,\ldots,\omega_g)$ of $\CC\otimes\TT^\vee$ that are
fixed by the complex conjugation. We let:
\[
\phi\colon \RR\otimes\TT \lto \RR^{r_1}\times\CC^{r_2}\times \CC^{r_2} 
\lto \RR^{r_1}\times\CC^{r_2}
\]
be the map obtained from our basis (each $\omega_i$ gives $t\mapsto
a_1(t\omega_i)$), composed with the projection. We view
$\RR^{r_1}\times\CC^{r_2}$ as $\RR^g$ by decomposing each factor $\CC$
as $\RR\oplus\RR i$. Then we have:
\[
\Vol'(\RR\otimes \TT/\TT') 
= 2^{r_2}|\det(\phi(T_{i_1}),\ldots,\phi(T_{i_g}))|.
\]
By construction, each $\phi(T_{i_j})_k$ is the real or imaginary part
of some~$a_{i_j}(\omega_l)$. Hence by Lemma~\ref{lem_bnd_coefficients}
we have:
\[
|\phi(T_{i_j})_k| \leq 2i_j.
\]
We obtain:
\[
|\det(\phi(T_{i_1}),\ldots,\phi(T_{i_g}))| \leq 
\prod_{j=1}^g ( 2i_j \sqrt{g} ) \leq (4g^2)^g.
\]
Hence, finally:
\[
\log \Vol'(\RR \otimes \TT/\TT) \leq 
r_2\log 2 + g(\log 4 + 2\log g).
\]
Noting that $r_2\leq g/2$ and that $g=O((pl)^2)$ completes our proof.
\end{proof}

\section{Bounding the theta function on 
$\Pic^{\protect\lowercase{g}-1}(X_1(\protect\lowercase{pl}))$}
\label{sec_bnd_theta}
The aim of this section is to give a bound for the supnorm of the
theta function that occurs in Theorem~\ref{maininequality}. 

\begin{thm}\label{thm_bound_theta} 
For $X=X_1(pl)$, with $p$ and $l$ distinct primes for which the genus
of $X_1(pl)$ is at least one, we have $\log\|\vartheta\|_{\sup} =
O((pl)^6)$ .
\end{thm}
We start with two lemmas, which are possibly of independent interest.
\begin{lem}\label{boundimtau}
Let $X=V/\Lambda$ be a principally polarised complex Abelian variety and let
$H\colon V \times V \to \CC$ be its Riemann form. Let
$\lambda_1,\ldots,\lambda_{2g}$ be the successive minima of the lattice
$\Lambda$, with norm defined by $\|x\|^2 =  H(x,x)$. 
Let $(e_1,\ldots,e_{2g})$ be a symplectic basis of~$\Lambda$,
i.e., a basis with respect to which the matrix of the symplectic form
$\Im(H)$ is, in $g$ by $g$ block form, equal to 
$(\begin{smallmatrix}0 & 1\\-1 & 0\end{smallmatrix})$. The one has:
\[
\sqrt{\det \Im(\tau)} \leq 
\frac{(2g)!}{2^g} \frac{V_{2g}}{V_g} \lambda_{g+1} \cdots \lambda_{2g} \, ,
\]
where $\tau$ is the period matrix in~$\HH_g$ corresponding to
$(e_1,\ldots,e_{2g})$. Here $V_n$ denotes the volume of the unit ball in 
$\RR^n$ with its standard euclidean inner product.
\end{lem}
\begin{proof}
We consider the lattice
$M=\ZZ{\cdot}e_1\oplus\cdots\oplus\ZZ{\cdot}e_g $ in the real
subvector space $W=\RR{\cdot}e_1\oplus \cdots\oplus\RR{\cdot}e_g$
of~$V$. Denote by $\mu_1,\ldots,\mu_g$ the successive minima of~$M$,
where the norm is given by restricting $H$ to~$W$ (note that $M$ is
isotropic for the symplectic form, so that $H$ takes real values
on~$W$). We have $(H(e_i,e_j))_{i,j} = ( \Im(\tau) )^{-1}$, so that
the volume (with respect to the inner product on $W$ given by $H$) of
$W/M$ is equal to $(\det \Im(\tau))^{-1/2}$, and hence by Minkowski's
second fundamental inequality:
\[
\mu_1 \cdots \mu_g \leq 2^g \frac{ (\det \Im(\tau))^{-1/2} }{V_g } \,.
\] 
On the other hand we have:
\[
\mu_1 \cdots \mu_g \geq \lambda_1 \cdots \lambda_g =
(\lambda_1\cdots\lambda_{2g})\cdot(\lambda_{g+1}\cdots\lambda_{2g})^{-1}
\]
and since the volume of $V/\Lambda$ is 1 we obtain by Minkowski's first
fundamental inequality:
\[
\lambda_1\cdots\lambda_{2g} \geq \frac{2^{2g}}{(2g)!}\frac{1}{V_{2g}} \, .
\]
Combining we find a lower bound:
\[
\mu_1 \cdots \mu_g \geq 
\frac{2^{2g}}{(2g)!} \frac{1}{V_{2g}} 
\cdot(\lambda_{g+1} \cdots \lambda_{2g})^{-1} \, .
\]
Combining this with the upper bound for $\mu_1 \cdots \mu_g$ we obtain
the required formula.
\end{proof}
\begin{lem} \label{generators}
Let $N \geq 3$ be an integer.  The group $\Gamma_1(N)$ is generated by
its elements whose entries are bounded from above in absolute value
by~$N^6/4$.
\end{lem}
\begin{proof}
We first note the following: let $G$ be a group, and let $S\subset G$
be a set of generators. Let $X$ be a transitive $G$-set and let $x$ be
in~$X$. For each $y$ in~$X$, let $g_y$ be an element of~$G$ such that
$g_yx=y$; we demand that~$g_x=1$. Then the $g_{sy}^{-1}sg_y$, for $s$
in~$S$ and $y$ in~$X$, form a system of generators for the stabiliser
$G_x$ of~$x$.  To see this, first replace $S$ by $S\cup S^{-1}$. Let
$g$ be in~$G_x$. Write $g=s_n\cdots s_1$ with $s_i$ in~$S$. Then we
can write:
\[
g = s_n\cdots s_1 = 
g_{s_ny_n}^{-1}s_ng_{y_n}\cdots g_{s_1y_1}^{-1}s_1g_{y_1}
\quad\text{with $y_i=(s_{i-1}\cdots s_1)x$}.
\]
The equality holds because $s_iy_i=y_{i+1}$, and $g_{y_1}=1$ and
$g_{s_ny_n}=1$.  Now we apply this to our case. We take
$G=\SL_2(\ZZ)$, and we take $X$ to be the subset of $(\ZZ/N\ZZ)^2$
consisting of the elements of order~$N$. This is a transitive
$G$-set. We let $x=(1,0)$; then $G_x$ is identified
with~$\Gamma_1(N)$.  Let $S$ be the set consisting of
$(\begin{smallmatrix}1 & 1\\0 & 1\end{smallmatrix})$ and
$(\begin{smallmatrix}1 & 0\\1 & 1\end{smallmatrix})$ and their
inverses. Then $S$ generates~$G$.  We now apply the previous argument
to find generators of~$\Gamma_1(N)$.  So let
$y=(\ol{a},\ol{b})$ be in~$X$.  Then, thinking of
$\ZZ/N\ZZ$ as the product of its local rings, we see that there is a
$u$ in $\ZZ$ with $|u|\leq N/2$ and $\ol{a+bu}$
in~$(\ZZ/N\ZZ)^\times$. Put $a_1:=a+bu$ and $b_1:=b$. Then:
\[
\left(\begin{matrix}a_1 \\b_1\end{matrix}\right) = 
\left(\begin{matrix}1 & u\\ 0 & 1\end{matrix}\right)
\left(\begin{matrix}a \\b \end{matrix}\right).
\]
Next, there is a $v$ in $\ZZ$ with $|v|\leq N/2$ and $b_1+a_1v=1$
mod~$N$, i.e.:
\[
\left(\begin{matrix}a_1 \\1\end{matrix}\right) = 
\left(\begin{matrix}1 & 0\\ v & 1\end{matrix}\right)
\left(\begin{matrix}a_1 \\b_1 \end{matrix}\right)\bmod N.
\]
Finally, let $w$ be in $\ZZ$ with $|w|\leq N/2$ and with image $\ol{a_1}$
in~$\ZZ/N\ZZ$. Then one has:
\[{\scriptsize 
\left(\begin{matrix}a \\b \end{matrix}\right) = 
\left(\begin{matrix}1 & -u\\ 0 & 1\end{matrix}\right)
\left(\begin{matrix}1 & 0\\ -v & 1\end{matrix}\right)
\left(\begin{matrix}1 & w\\ 0 & 1\end{matrix}\right)
\left(\begin{matrix}0 & -1\\ 1 & 0\end{matrix}\right)
\left(\begin{matrix}1 \\0 \end{matrix}\right) \bmod N.}
\]
Writing out, we have:
\[  { \scriptsize 
\left(\begin{matrix}1 & -u\\ 0 & 1\end{matrix}\right)
\left(\begin{matrix}1 & 0\\ -v & 1\end{matrix}\right)
\left(\begin{matrix}1 & w\\ 0 & 1\end{matrix}\right)
\left(\begin{matrix}0 & -1\\ 1 & 0\end{matrix}\right)
 =
\left(\begin{matrix} w-u(1-vw)& -1-uv\\ 1-vw & v\end{matrix}\right)}
\]
so that we can put:
\[
g_y = g_{(a,b)} = 
\left(\begin{matrix} w-u(1-vw)& -1-uv\\ 1-vw & v\end{matrix}\right).
\]
The absolute values of the coefficients of $g_{(a,b)}$ are smaller
than~$N^3/4$, if $N \geq 3$, and the lemma follows. 
\end{proof}
\begin{proof}
[Proof of Theorem~\ref{thm_bound_theta}] 
Recall that $\|\vartheta\|(z;\tau)$ is given by:
\[
\|\vartheta\|(z;\tau) = 
(\det \Im(\tau))^{1/4} \exp(-\pi\,{}^t\hspace{-0.1em}y
(\Im(\tau))^{-1}  y)  |\vartheta(z;\tau)| , 
\]
where $y = \Im(z)$ and where $\tau$ is a period matrix in the Siegel
upper half plane $\HH_g$ corresponding to~$X$. We first deal with the
factor $\det \Im(\tau)$ and for this we invoke
Lemma~\ref{boundimtau}. Choose once more an Atkin-Lehner basis
$(\omega_1,\ldots,\omega_g)$ for~$\Omega^1(X)$. Using the dual basis
in $\Omega^1(X)^\vee$ we write:
\[
J(X)= \CC^g/\Lambda\, ,
\]
where:
\[
\Lambda=\mathrm{Image}\left(H_1(X,\ZZ) \to \CC^g \colon \gamma \mapsto
\int_\gamma (\omega_1,\ldots,\omega_g)\right).
\]
The polarisation form for $J(X)$ is given by:
\[
(z,w)\mapsto
{}^t\hspace{-0.1em}z\cdot(\ld\omega,\omega\rd)^{-1}_{i,j}\cdot\ol{w}.
\]
Denote by $\| \cdot \|_P$ the corresponding norm on~$\CC^g$. We also
consider the standard hermitian inner product on $\CC^g$, which is
just $(z,w) \mapsto {}^t\hspace{-0.1em}z \cdot \ol{w}$.  Here we denote the
corresponding norm by $\| \cdot \|_E$.  From the next two lemmas we
obtain:
\[
\left( \lambda_{g+1} \cdots \lambda_{2g} \right)^2 \leq
\left( g e^{4\pi} (pl)^{46} / \pi \right)^g
\]
and hence, by Lemma~\ref{boundimtau}, the estimate:
\[
\log (\det \Im(\tau)) = O((pl)^2\log(pl)).
\]
\begin{lem}
The lattice $\Lambda$ is generated by the subset of its elements $x$
that satisfy $\|x\|_E^2\leq g{\cdot}(pl)^{46}$.
\end{lem}
\begin{proof}
For the moment put $N=pl$. Following the natural surjections:
\[
\Gamma_1(N) \onto \Gamma_1(N)^\ab = H_1(Y_1(N),\ZZ) \onto
H_1(X_1(N),\ZZ)  
\]
we see that any generating set for $\Gamma_1(N)$ gives a generating
set for $H_1(X_1(N),\ZZ)$. We are going to take generators for
$\Gamma_1(N)$ 
as given by Lemma~\ref{generators}. In particular, the absolute values
of their coefficients are bounded by~$N^6/4$. We have to see now what
this implies for $\|x\|_E$ for corresponding elements $x$
of~$\Lambda$.  Concretely, choose a 
$g=(\begin{smallmatrix}a & b\\c & d\end{smallmatrix})$ 
in~$\Gamma_1(N)$. The image in $\Lambda$ can be given as follows:
choose any $z$ in~$\HH$ and any path in $\HH$ from $z$ to~$gz$. This
gives us a loop in~$X_1(N)$, and the class of that loop is the image
of $g$ in~$H_1(X_1(N),\ZZ)$. In order to get to $\Lambda$ we compute
the periods of $(\omega_1,\ldots,\omega_g)$ around this loop. We want
to get bounds for these periods. In order to do this, note first that
we can assume that $c \neq 0$. Indeed, the $g$ with $c=0$ are
unipotent, hence have trivial image in~$H_1(X_1(N),\ZZ)$. Now we make
the following choices. First, we want to take a $z$ in $\HH$ with
$\Im(z) =\Im(gz)$. Since, for all $z$ in~$\HH$, $\Im(gz) =
|cz+d|^{-2}\Im(z)$, the condition that $\Im(z) =\Im(gz)$ is equivalent
to $|cz+d|=1$. We choose $z = -d/c + i/|c|$, i.e., such that $cz+d=\pm
i$, depending on the sign of~$c$. Second, the path that we take is the
straight line from $z$ to~$gz$.  We have $|c|\geq N$ because $g$ is
in~$\Gamma_1(N)$. Using furthermore that the absolute values of the
coefficients of $g$ are bounded by $N^6/4$ we get $|gz-z| \leq |az+b|
+ |z| \leq N^{11}/10$ (where we have used that $N\geq 2$).  For the
period of an element $\omega=fdq/q$ of the Atkin-Lehner basis we
obtain from this that:
\[
\left| \int_z^{gz} \omega \right| = 
\left|\int_z^{gz}f{\cdot}(dq)/q \right| =
\left|\int_z^{gz}f{\cdot}2\pi i\,dw \right| \leq 
\frac{2\pi N^{11}}{10} \|f\|, 
\]
where $w$ denotes the standard coordinate of~$\HH$ (i.e., the
inclusion map into~$\CC$), and $\|f\|$ the supnorm of $f$ on the
straight line from $z$ to~$gz$; recall that $q=\exp(2\pi iw)$. When
writing $f=\sum_{n\geq 1}a_n(f)q^n$ we have $|a_n(f)|\leq 2n$ by
Lemma~\ref{lem_bnd_coefficients}.  Noting furthermore that that we
have $\Im z=|c|^{-1} \geq 4/N^6$ it follows that:
\[
\|f\| \leq 2\sum_{n\geq 1}n e^{-8\pi N^{-6}n} = \frac{2r}{(1-r)^2}, 
\quad\text{where}\quad r = e^{-8\pi /N^6}.
\]
Hence:
\[
\|f\| \leq (N^6/4\pi)^2, \quad\text{and}\quad 
\left|\int_z^{gz}f(dq)/q \right| \leq 
\frac{2\pi}{10} N^{11} \left(\frac{N^6}{4\pi}\right)^2 \leq N^{23}.
\]
This means that all $g$ coordinates of our element $x$ of~$\CC^g$ are,
in absolute value, at most~$N^{23}$. Hence $\|x\|_E^2$, being the sum
of the squares of these coordinates, is at
most~$g{\cdot}N^{46}=g{\cdot}(pl)^{46}$.
\end{proof}
\begin{lem}
For any $x$ in $\CC^g$ we have the estimate:
\[
\|x\|^2_P \leq (e^{4\pi}/\pi) \|x\|_E^2 .
\] 
\end{lem}
\begin{proof}
By Lemma~\ref{gram} the matrix $(\ld\omega,\omega\rd)^{-1}$ is
almost diagonal, having in fact diagonal elements $1/\|\omega\|^2$
corresponding to newforms $\omega$ on~$X_1(pl)$, and 2-by-2 blocks
corresponding to newforms $\omega$ on~$X_1(l)$ and~$X_1(p)$. The
2-by-2 block corresponding to a newform $\omega$ on~$X_1(l)$ is:
\[
\frac{1}{(p-1) \|\omega\|^2 ((p+1)^2-|a_p(\omega)|^2)} 
\left( \begin{array}{cc} p+1 &
-\ol{a_p(\omega)} \\ -a_p(\omega) & p+1 \end{array} \right) \, ,
\]
where the norm $\|\omega\|^2$ is taken on~$X_1(l)$. 

The 2-by-2 block corresponding to a newform $\omega$ on~$X_1(p)$ is:
\[
\frac{1}{(l-1) \|\omega\|^2 ((l+1)^2-|a_l(\omega)|^2)} 
\left( \begin{array}{cc} l+1 &
-\ol{a_l(\omega)} \\ -a_l(\omega) & l+1 \end{array} \right) \, ,
\]
where the norm $\|\omega\|^2$ is taken on~$X_1(p)$.

A short calculation shows that for any $(z_1,z_2)$ in~$\CC^2$ one has:
\[
\frac
{
\left(\begin{array}{cc}z_1 & z_2\end{array}\right) 
\left( \begin{array}{cc} p+1 &
-\ol{a_p(\omega)} \\ -a_p(\omega) & p+1 \end{array} 
\right) \left( \begin{array}{c}
\ol{z_1} \\ \ol{z_2} \end{array} \right)
}
{(p-1)((p+1)^2-|a_p(\omega)|^2)} 
\leq \left( |z_1|^2 + |z_2|^2 \right)\, ,
\]
and similarly for $X_1(p)$, so that all in all one gets:
\[
\|(z_1,\ldots,z_g)\|_P^2 \leq 
\frac{ |z_1|^2 }{ \| \omega_1 \|^2 } + \cdots + 
\frac{ |z_g|^2 }{ \|\omega_g \|^2} \, .
\]  
The lemma follows by the lower bound from Lemma~\ref{eqn_lowerb_pet_norm}.
\end{proof}
Next we consider the factor 
$\exp(-\pi\,{}^t\hspace{-0.1em}y (\Im(\tau))^{-1}y)|\vartheta(z;\tau)|$.  
Since in our previous estimates the choice of
$\tau$ was irrelevant, it will cause no loss of generality here if we
restrict to $\tau$ lying in the so-called Siegel fundamental
domain~$F_g$, which is the set of matrices $\tau=x+iy$ satisfying the
conditions:
\begin{quote}
\begin{enumerate}
\item for each entry $x_{ij}$ of $x$ one has $|x_{ij}| \leq \frac{1}{2}$, 
\item for all $\gamma$ in $ \mathrm{Sp}(2g,\ZZ)$ one has $\det
\Im(\gamma{\cdot}\tau) \leq \det \Im(\tau)$,
\item $y$ is Minkowski reduced, i.e.,  for each $\xi =
(\xi_1,\ldots,\xi_g)$ in $\ZZ^g$ and each $i$ such that
$\xi_i,\ldots,\xi_g$ are non-zero, one has 
$\xi\, y\, {}^t\hspace{-0.1em}\xi \geq y_{ii}$ and moreover, for each
$1 \leq i \leq g-1$ one has $y_{i,i+1} \geq 0$.
\end{enumerate}
\end{quote}
It is well known that $F_g$ contains at least one representative from
each $\mathrm{Sp}(2g,\ZZ)$-orbit on~$\HH_g$.  We claim that for $\tau$
in $F_g$ the estimate:
\[
\exp(-\pi\, {}^t\hspace{-0.1em}y (\Im(\tau))^{-1} y) |\vartheta(z;\tau)|
\leq 2^{3g^3+5g} 
\]
holds, for all $z$ in~$\CC^g$.  Thus, this factor gives us a
contribution~$O((pl)^6)$.  In order to prove the estimate, write $y =
\Im(z) = (\Im(\tau)) \cdot b$ with $b$ in~$\RR^g$. Then it is easy to
see that:
\[
\exp(-\pi\, {}^t\hspace{-0.1em}y(\Im(\tau))^{-1}y)|\vartheta(z;\tau)| 
\leq \sum_{n \in \ZZ^g} 
\exp(-\pi\, {}^t\hspace{-0.1em}(n+b)(\Im(\tau)) (n+b)) \, .
\]
Since the $\Im(\tau)$ are Minkowski reduced we have, for any~$m$
in~$\RR^g$ (cf. \cite{Igusa1}, V~\S 4):
\[
{}^t\hspace{-0.1em}m \Im(\tau)m \geq 
c(g) \sum_{i=1}^g m_i^2 (\Im(\tau))_{ii}, \quad
c(g)=\left(\frac{4}{g^3}\right)^{g-1}\left(\frac{3}{4}\right)^{g(g-1)/2}.
\]
Moreover, we have $(\Im(\tau))_{i,i} \geq \sqrt{3}/2$ for
$i=1,\ldots,g$. From this we derive:
\begin{multline*}
\sum_{n \in \ZZ^g} \exp(-\pi\, {}^t\hspace{-0.1em}(n+b)
(\Im(\tau))  (n+b))  
\leq \\
\leq \sum_{n \in \ZZ^g} \exp\left(-\sum_{i=1}^g 
\pi c(g) (n_i+b_i)^2 (\Im(\tau))_{i,i} \right) \leq \\
\leq \prod_{i=1}^g \sum_{n_i \in \ZZ} \exp( - \pi c(g) (n_i +
 b_i)^2 (\Im(\tau))_{i,i} ) \leq \\
\leq \prod_{i=1}^g \frac{2}{1-\exp(-\pi c(g) (\Im(\tau))_{ii})} \leq \\
\leq 2^g \left( 1 + \frac{2}{\sqrt{3} \pi c(g)} \right)^g \, .
\end{multline*}
From this and the formula for $c(g)$ the required estimate follows and
the proof of Theorem~\ref{thm_bound_theta} is finished.
\end{proof}

\section{Upper bounds for Arakelov Green functions on 
the curves~$X_1(\protect\lowercase{pl})$}
\label{sec_bnd_green_f}
The aim of this section is to give an upper bound for the Arakelov
Green functions on the curves $X_1(pl)$ that will enable us to bound from
above the contributions of the intersection numbers in the right hand
side of the inequality in Theorem~\ref{maininequality}. As the
$X_l(\CC)$ are compact, it is clear that for each $l$ such an upper
bound exists, but we need such upper bounds that grow as most as a
power of~$l$.

In order to establish such upper bounds we will use Franz Merkl's
result on Green functions on arbitrary Riemann surfaces given in
Chapter~\ref{subsec_merkl}. 

Instead of using the result of Merkl for our work we could certainly
also have used recent work by Jorgenson and Kramer
in~\cite{Jorgenson-Kramer1}. The results of Jorgenson and Kramer date
back to the same time as those of Merkl (early Spring 2004). We chose
to use Merkl's results because his approach is more elementary, and we
had the details earlier than those of Jorgenson and Kramer.

The following theorem gives a suitable upper bound for the
Arakelov-Green functions~$g_{a,\mu}$ (see~(\ref{eqn_arakelov_1-1-form})
and Proposition~\ref{prop_Ar-Gr-function}) on the modular
curves~$X_1(pl)$ with $p$ and $l$ distinct primes.
\begin{thm} \label{appl_merkl_thm1}
There is a real number $c$ such that for all pairs of distinct prime
numbers $p$ and~$l$ for which the genus of $X_1(pl)$ is at least one
and for all distinct $a$ and~$b$ on $X_1(pl)(\CC)$ we have:
\[ 
g_{a,\mu}(b) \leq c{\cdot}(pl)^6.
\]
Let $\infty$ denote the cusp $\infty$ on~$X_1(pl)$, and let $q$ be the
standard local coordinate around~$\infty$ given by the map
$\tau\mapsto \exp(2\pi i\tau)$ from the region $\Im\tau>1$
in~$\HH$ to~$\CC$. Then we have:
\[
\left|\log\|dq\|_\mathrm{Ar}(\infty)\right| = O((pl)^6),
\]
where $\|{\cdot}\|_\mathrm{Ar}$ denotes the Arakelov metric
on~$\Omega^1$ (see Section~\ref{subsec_ARR}).
\end{thm}
\begin{proof} 
We write for the moment $N$ for~$pl$. We will apply
Theorem~\ref{merkl_thm1}, but we will carry out the estimates on the
more symmetrical modular curve $X(N)$ which for us is $\Gamma(N)
\backslash (\HH \cup \PP^1(\QQ))$. Let $h \colon X(N) \to X_1(N)$ be
the canonical map; it has degree~$N$. We let $\mu$ denote the Arakelov
$(1,1)$-form on~$X_1(N)$, and we define $\mu'=h^*\mu/N$. The
characterising properties of Green functions directly imply that:
\begin{eqn}\label{appl_merkl_eqn0}
h^*g_{a,\mu} = \sum_{h(b)=a}g_{b,\mu'},
\end{eqn}
where the $b$ are counted with multiplicity.

As in Section~\ref{subsec_merkl} we fix a constant~$r_1$ with
$0<r_1<1$; we take $r_1:=3/4$.  We need to construct an atlas with
charts $z^{(j)}\colon U^{(j)} \to \CC$ for $X(N)$ with all
$z^{(j)}(U^{(j)})$ containing the closed unit disk and with the
$U_{r_1}^{(j)}$ covering~$X(N)$.

We start with a construction of a local coordinate $z\colon U \to \CC$
in a neighbourhood of the standard cusp~$\infty$.  Since
$\SL_2(\ZZ/N\ZZ)$ acts transitively on the set of cusps of~$X(N)$,
this construction will suffice to give the full atlas. Our initial
coordinate is induced by the map $z$ from $\HH$ to $\CC$ that sends
$\tau$ to~$e^{2\pi i\tau/N}$. As the following lemma is valid for all
integers $n\geq1$, we state it in that generality, and will apply it
with $n:=N$.
\begin{lem}
Let $n$ be in~$\ZZ_{\geq1}$. The subset in~$\HH$ given by the
conditions $-1/2\leq\Re \tau <n-1/2$ and $\Im \tau > 1/n$ is mapped
injectively to~$X(n)(\CC)$.
\end{lem}
\begin{proof}
First we note that for $(\begin{smallmatrix}a & b\\c &
d\end{smallmatrix})$ in~$\SL_2(\ZZ)$ and for $\tau$ in~$\HH$ we
have:
\[
\Im \left(\frac{a\tau+b}{c\tau+d}\right) = \frac{\Im(\tau)}{|c\tau+d|^2}.
\]
Let us call $D$ the set of $\tau\in\HH$ that satisfy the two
conditions of the lemma: 
\[
D = \{\tau\in\HH\;|\; 
\text{$-1/2\leq\Re(\tau)<n-1/2$ and $\Im(\tau)>1/n$}\}.
\]
Let $\tau$ be in~$D$, and let
$\gamma=(\begin{smallmatrix}a&b\\c&d\end{smallmatrix})$ be
in~$\Gamma(n)$, such that $\gamma\tau\neq\tau$. If $c=0$ then, as
$ad=1$, we have $a=d$ and $\tau'=\tau\pm b$ with $b$ a non-zero
multiple of~$n$, and so $\tau'$ is not in~$D$. If $c\neq0$ then we
have $|c|\geq n$ because~$n|c$, and:
\begin{align*}
\Im\left(\frac{a\tau+b}{c\tau+d}\right) & =
\frac{\Im(\tau)}{|c\tau+d|^2} 
\leq \frac{\Im(\tau)}{(\Im(c\tau))^2} \leq 
\frac{\Im(\tau)}{n^2(\Im(\tau))^2} = \\
& = \frac{1}{n^2\Im(\tau)} < \frac{1}{n^2{\cdot}(1/n)} = \frac{1}{n},
\end{align*}
hence $\gamma\tau$ is not in~$D$.
\end{proof}
In particular, the region of $\tau$ with $-1/2\leq \Re\tau < N-1/2$
and $\Im\tau>1/2$ is mapped injectively into $X(N)$ to give an open
neighbourhood $U$ of~$\infty$. We could replace the condition
``$\Im\tau>1/2$'' by ``$\Im\tau>1/N$'', but that would not make the
work to be done significantly easier. The map $\tau\mapsto e^{2\pi
i\tau/N}$ gives an isomorphism:
\begin{eqn}\label{appl_merkl_eqn1}
z\colon U\lto D(0,e^{-\pi/N})\subset\CC.
\end{eqn}

The region of $\tau$ with $-1/2\leq\Re\tau<N-1/2$ and $\Im\tau >3/4$
gives an open neighbourhood~$V$ of~$\infty$, contained in~$U$, such
that the translates of~$V$ under~$\SL_2(\ZZ/N\ZZ)$ cover~$X(N)$ (note
that $3/4<\sqrt{3}/2$). The image of $V$ under $z$ is the
disk~$D(0,e^{-3\pi/2N})$. However, the quotient of the radii
$e^{-\pi/N}$ and $e^{-3\pi/2N}$ tends to~$1$ as $N$ tends to infinity,
hence we cannot work with these disks directly.

What we do instead is the following. We define a new coordinate
$z':=e^{3\pi/2N}z$ to get $z'V=D(0,1)$. Then
$z'U=D(0,e^{\pi/2N})$. Let $\eps:=\eps(N):=e^{\pi/2N}-1$ be the
difference between the two new radii. Then $\eps>\pi/2N$. We can
choose $O(\eps^{-2})$ open disks $D(a,\eps)$ with centre $a$ in
$D(0,1)$, such that the union of the $D(a,\eps/2)$ contains $D(0,1)$;
we let $A$ denote the set of these~$a$. The $D(a,\eps)$ are contained
in~$z'U=D(0,1+\eps)$.
%% Hier een mooi plaatje maken???

The group $\SL_2(\ZZ/N\ZZ)$ acts transitively on the set of cusps
of~$X(N)$. For each cusp~$c$, we choose a $g_c$ in~$\SL_2(\ZZ/N\ZZ)$
such that $c=g_c\infty$. The open sets of our atlas for $X(N)$ are
then the $U^{(a,c)}$ with $a$ in~$A$ and $c$ a cusp, defined by:
\[
U^{(a,c)} := g_c\cdot (z')^{-1}D(a,\eps).
\]
The required coordinates $z^{(a,c)}$ on the $U^{(a,c)}$ are defined by the
composition of isomorphisms:
\[
\xymatrix{
z^{(a,c)}\colon U^{(a,c)} \ar[r]^(0.6){g_c^{-1}} &
U^{(a,\infty)} \ar[r]^{z'} & D(a,\eps) \ar[r]^{-a} & 
D(0,\eps) \ar[r]^{{\cdot}3/2\eps} & D(0,\frac{3}{2}).
}
\]
Indeed, the images $z^{(a,c)}U^{(a,c)}$ contain the unit disk, and
$U^{(a,c)}_{r_1}$ corresponds via $z'\circ g_c^{-1}$ to the subdisk
$D(a,\eps/2)$ of $D(a,\eps)$, hence the $U^{(a,c)}_{r_1}$
cover~$X(N)$. The exact number $n=n(N)$ of $U^{(a,c)}$ is the
cardinality of $A$ times the number of cusps, hence~$n=O(N^4)$. We
choose a numbering of $A\times\{\text{cusps}\}$ with the integers
$\{1,\ldots,n\}$, and we will denote our charts as: 
\begin{eqn}\label{appl_merkl_eqn2}
z^{(j)}\colon U^{(j)}\to D(0,3/2)\subset\CC .
\end{eqn}

\begin{lem}\label{appl_merkl_lem0}
For the local coordinates $z^{(j)}\colon U^{(j)}\to\CC$ on $X(N)$ that
we have just defined, the following holds. For all $j$ and~$k$
in~$\{1,\ldots,n\}$ we have:
\[
\sup_{U_1^{(j)}\cap U_1^{(k)}}\left|\frac{dz^{(j)}}{dz^{(k)}}\right| 
\leq M,
\]
with $M=6$.
\end{lem}
\begin{proof}
Let $j$ and~$k$ be in $\{1,\ldots,n\}$. If $j$ and~$k$ arise from the
same cusp, then $z^{(j)}$ and $z^{(k)}$ differ by a translation, hence
$dz^{(j)}/dz^{(k)}=1$. Now suppose that $j$ and~$k$ arise from two
distinct cusps. We may suppose then, by acting with an element
of~$\SL_2(\ZZ/N\ZZ)$, that $k$ arises from the standard
cusp~$\infty$. Let $x$ denote the cusp that $j$ arises from. The
coordinate $z^{(j)}$ is then obtained as above from an element $g_x$
of $\SL_2(\ZZ)$ that sends $\infty$ to~$x$. Let us write
$g_x^{-1}=(\begin{smallmatrix}a & b\\c & d\end{smallmatrix})$. Note
that $c\neq 0$, hence $|c|\geq 1$. Let $z$ be a point in~$\HH$ with
$-1/2\leq \Re z < N-1/2$ that maps to an element in~$U_1^{(j)}\cap
U_1^{(k)}$. Then we know that $1/2<\Im z<1$ because disks given by
$\Im z>1$ around different cusps do not meet at all. Likewise, we then
know that:
\[
\frac{1}{2} < \Im\left(\frac{az+b}{cz+d}\right) = 
\frac{\Im(z)}{|cz+d|^2} < 1.
\]
Hence, as $(\Im z)/|cz+d|^2 \leq (\Im z)/c^2(\Im z)^2$, we have $\Im z
< 2/c^2$ which gives in fact $|c|=1$. Under these conditions, we
estimate:
\begin{multline*}
\left| \log \left|
\frac{d\exp(2\pi i\frac{az+b}{cz+d}/N)}{d\exp(2\pi i z/N)}
\right|\right| = 
\left| \log \left|
\frac{\exp(2\pi i\frac{az+b}{cz+d}/N)d\frac{az+b}{cz+d}}{\exp(2\pi i z/N)dz}
\right|\right| = \\
= \left|\log\left| \exp(2\pi i\frac{az+b}{cz+d}/N) \right| -
\log\left|\exp(2\pi i z/N)\right| - \log|cz+d|^2\right| \leq \\
\leq 4\pi/N + 4\pi/N + \log 4.
\end{multline*}
So indeed, for $N\geq 6$, we can take $M=6$. Some explanations are
perhaps in order here: as $\Im z$ and $\Im (az+b)/(cz+d)$ are between
$1/2$ and~$2$, $|\exp(2\pi i z/N)|$ and $|\exp(2\pi i
\frac{az+b}{cz+d}/N)|$ are between $\exp(-4\pi/N)$ and
$\exp(-\pi/N)$. As $\Im (az+b)/(cz+d)=(\Im z)/|cz+d|^2$, we see that
$|cz+d|^2$ is between $1/4$ and~$4$.
\end{proof}

Our next task is to produce a suitable bound, as in
Theorem~\ref{merkl_thm1}, of the type $\mu\le c_1|dz^{(j)}\wedge d\bar
z^{(j)}|$. We start with a bound for $\mu$ on disks around~$\infty$
on~$X_1(pl)$.

\begin{lem}\label{appl_merkl_lem1}
Let $r$ be a real number such that $0<r<1$. We map $D(0,r)$ to
$X_1(pl)$ by sending $q\neq0$ to $(\CC^\times/q^\ZZ,\zeta_{pl})$. The
image of this map is the image in $X_1(pl)$ of the region in~$\HH$
defined by the condition ``$\Im\tau>-(\log r)/2\pi$'', plus the
cusp~$\infty$. We still denote by $\mu$ the $(1,1)$-form on $D(0,r)$
induced by~$\mu$. Then we have, on~$D(0,r)$:
\[
\mu \leq \frac{28e^{4\pi}}{\pi}\frac{1}{(1-r)^4}\cdot \frac{i}{2}dqd\ol{q} .
\]
\end{lem}
\begin{proof}
We first bound, on the disk~$D(0,r)$, and for a newform~$f$, the
functions $\sum_{n\geq 1}a_n(f)q^{n-1}$ and $\sum_{n\geq
  1}a_n(f)q^{nl-1}$. We have, for $|q|<r$ by
Lemma~\ref{lem_bnd_coefficients}:
\[
|\sum_{n\geq 1}a_n(f)q^{n-1}| \leq \sum_{n\geq1}
|a_n(f)|r^{n-1}  \leq  2\sum_{n\geq 1}
nr^{n-1} = \frac{2}{(1-r)^2} \, ,
\]
and next:
\[
|\sum_{n\geq 1}a_n(f)q^{nl-1}| 
\leq 2\sum_{n\geq 1} nr^{nl-1} = \frac{2r^{l-1}}{(1-r^l)^2}
\quad\text{for $|q|<r$}.
\]
Now recall from Corollary~\ref{expressionmu} that for $\mu$ we have the
expression:
\begin{multline*} 
\mu = \frac{i}{2g}\sum_\omega 
\frac{\omega\wedge\ol{\omega}}{\|\omega\|^2} + 
\frac{i}{2g}\sum_\omega \left(
\frac{(p+1)s_p^*(\omega\wedge\ol{\omega}) + 
(p+1)t_p^*(\omega\wedge\ol{\omega})}
{(p-1)\|\omega\|^2((p+1)^2-|a_p(f_\omega)|^2)} \right. \\
\left. - \frac{a_p(f_\omega)s_p^*\omega\wedge\ol{t_p^*\omega} + 
\ol{a_p(f_\omega)}t_p^*\omega\wedge\ol{s_p^*\omega}}
{(p-1)\|\omega\|^2((p+1)^2-|a_p(f_\omega)|^2)}
\right)\\
+ 
\frac{i}{2g}\sum_\omega \left(
\frac{(l+1)s_l^*\omega\wedge\ol{s_l^*\omega} + 
(l+1)t_l^*\omega\wedge\ol{t_l^*\omega}}
{(l-1)\|\omega\|^2((l+1)^2-|a_l(f_\omega)|^2)} \right. \\
\left. - \frac{a_l(f_\omega)s_l^*\omega\wedge\ol{t_l^*\omega} + 
\ol{a_l(f_\omega)}t_l^*\omega\wedge\ol{s_l^*\omega}}
{(l-1)\|\omega\|^2((l+1)^2-|a_l(f_\omega)|^2)}
\right)\, ,
\end{multline*} 
the first sum running over~$\Omega^1(X_1(pl))^\new$, the second sum
running over $\Omega^1(X_1(l))^\new$, and the third sum running
over~$\Omega^1(X_1(p))^\new$.  We bound the different terms from the
above expression for~$\mu$.  The contribution of an $\omega$ in
$\Omega^1(X_1(pl))^\new$ gives, for $|q|<r$:
\begin{align*}
\frac{i}{2g\|\omega\|^2} \omega\wedge\ol{\omega} & = 
\frac{1}{g\|\omega\|^2} 
\left|\sum_{n\geq 1}a_n(f_\omega)q^{n-1}\right|^2
\cdot \frac{i}{2}dqd\ol{q} \\
& \leq \frac{1}{g} \frac{e^{4\pi}}{\pi} \left(\frac{2}{(1-r)^2}\right)^2
\cdot \frac{i}{2}dqd\ol{q} \\
& = \frac{1}{g} \frac{e^{4\pi}}{\pi}\frac{4}{(1-r)^4}\cdot
\frac{i}{2}dqd\ol{q} 
\end{align*}
The contribution of an element $\omega$ of~$\Omega^1(p)^\new$ is:
\begin{multline*}
\frac{i}{2g}\sum_\omega \left(
\frac{(l+1)s_l^*\omega\wedge\ol{s_l^*\omega} + 
(l+1)t_l^*\omega\wedge\ol{t_l^*\omega}}
{(l-1)\|\omega\|^2((l+1)^2-|a_l(f_\omega)|^2)} \right. \\
\left. - \frac{a_l(f_\omega)s_l^*\omega\wedge\ol{t_l^*\omega} + 
\ol{a_l(f_\omega)}t_l^*\omega\wedge\ol{s_l^*\omega}}
{(l-1)\|\omega\|^2((l+1)^2-|a_l(f_\omega)|^2)}
\right) \leq \\
\leq \frac{1}{g} \frac{e^{4\pi}}{\pi(l-1)^3}
\left(\frac{4(l+1)}{(1-r)^4}+
\frac{4(l+1)r^{2(l-1)}}{(1-r^l)^4}+\right. \\
\left. \frac{16\sqrt{l}r^{l-1}}{(1-r)^2(1-r^l)^2}\right)
\cdot \frac{i}{2}dqd\ol{q} .
\end{multline*}
Here one uses the Weil bounds on~$a_l(f_\omega)$.
Symmetrically (in $p$ and~$l$), the contribution to $\mu$ of an
element $\omega$ of~$\Omega^1(l)^\new$ is:
\begin{multline*}
\frac{i}{2g}\sum_\omega \left(
\frac{(p+1)s_p^*\omega\wedge\ol{s_p^*\omega} + 
(p+1)t_p^*\omega\wedge\ol{t_p^*\omega}}
{(p-1)\|\omega\|^2((p+1)^2-|a_p(f_\omega)|^2)} \right. \\
\left. - \frac{a_p(f_\omega)s_p^*\omega\wedge\ol{t_p^*\omega} + 
\ol{a_p(f_\omega)}t_p^*\omega\wedge\ol{s_p^*\omega}}
{(p-1)\|\omega\|^2((p+1)^2-|a_p(f_\omega)|^2)}
\right) \leq \\
\leq \frac{1}{g} \frac{e^{4\pi}}{\pi(p-1)^3}
\left(\frac{4(p+1)}{(1-r)^4}+
\frac{4(p+1)r^{2(p-1)}}{(1-r^p)^4}+\right. \\
\left. \frac{16\sqrt{p}r^{p-1}}{(1-r)^2(1-r^p)^2}\right)
\cdot \frac{i}{2}dqd\ol{q} .
\end{multline*}
Now we sum all contributions up, over the elements of
$\Omega^1(X_1(pl))^\new$, $\Omega^1(X_1(p))^\new$,
and~$\Omega^1(X_1(l))^\new$.  We get for $|q|<r$:
\begin{multline*}
\mu \leq \frac{4e^{4\pi}}{\pi}
\left(
\frac{1}{(1-r)^4} +
\frac{1}{(1-r)^4} + \frac{r^{2(l-1)}}{(1-r^l)^4} +
\frac{r^{l-1}}{(1-r)^2(1-r^l)^2} + \right.\\
\left.\frac{1}{(1-r)^4} + \frac{r^{2(p-1)}}{(1-r^p)^4} +
\frac{r^{p-1}}{(1-r)^2(1-r^p)^2}
\right)
\cdot \frac{i}{2}dqd\ol{q} ,
\end{multline*}
and finally:
\[
\mu \leq \frac{28e^{4\pi}}{\pi}
\frac{1}{(1-r)^4}
\cdot \frac{i}{2}dqd\ol{q}, \quad \text{for $|q|<r$.}
\]
\end{proof}

Our next step is to consider the disks $gU$ in $X(pl)$, where $g$ is
in $\SL_2(\ZZ/pl\ZZ)$ and where $U$ is as in~(\ref{appl_merkl_eqn1}).

\begin{lem}\label{appl_merkl_lem2}
Let $g$ be in $\SL_2(\ZZ/pl\ZZ)$ and let $z\colon U\to
D(0,e^{-\pi/pl})$ be as in~(\ref{appl_merkl_eqn1}). Let $z_g:=z\circ
g^{-1}\colon gU\to D(0,e^{-\pi/pl})$. Then we have, for the
restriction to~$gU$ of the pullback $h^*\mu$ of $\mu$ along $h\colon
X(N)\to X_1(N)$:
\[
(h^*\mu)|_{gU} \leq c(pl)^4|dz_gd\ol{z_g}|,
\]
with $c$ independent of $p$ and~$l$.
\end{lem}
\begin{proof}
To prove this, we consider the map $h\circ g \circ z^{-1}$ from
$D(0,e^{-\pi/pl})$ to~$X_1(pl)$ and the pullback of $\mu$
to~$D(0,e^{-\pi/pl})$. We observe that $\mu$ is invariant under all
automorphisms of~$X_1(pl)$. This applies in particular to the diamond
operators and the Atkin-Lehner pseudo-involutions (defined
in~(\ref{defi_w_operator})). As the group generated by these
automorphisms permutes the cusps of $X_1(pl)$ transitively, we can
take such an automorphism $\alpha$ such that $\alpha\circ h\circ g
\circ z^{-1}$ sends $D(0,e^{-\pi/pl})$ to a disk around the
cusp~$\infty$, where we can then apply
Lemma~\ref{appl_merkl_lem1}. The pullbacks of $\mu$ via $h\circ g
\circ z^{-1}$ and $\alpha\circ h\circ g\circ z^{-1}$ are the same. We
are also free to replace the coordinate $z$ by $\zeta z$ with $\zeta\in\CC$
such that $|\zeta|=1$. 

The map $h\circ g \circ z^{-1}$ sends a point $0\neq q\in
D(0,e^{-\pi/pl})$ to the point of $X_1(pl)$ corresponding to
$(\CC^\times/q^{pl\ZZ},\zeta_p^a(q^l)^b,\zeta_l^c(q^p)^d)$ for certain
$a$ and $b$ in $\FF_p$ and $c$ and $d$ in $\FF_l$ depending
on~$g$. After replacing $h$ with $h$ composed with a suitable diamond
operator, and $z$ by $\zeta z$ with $\zeta$ a suitable element
of~$\mu_{pl}(\CC)$, we are in one of four cases, that we will treat
one by one.

In the first case, $q$ is mapped to
$(\CC^\times/q^{pl\ZZ},\zeta_{pl})$. Then the map $h\circ g \circ
z^{-1}$ factors as the cover $D(0,e^{-\pi/pl})\to D(0,e^{-\pi})$ of
degree~$pl$ sending $q$ to~$q^{pl}$, followed by the map of
Lemma~\ref{appl_merkl_lem1} that sends $q^{pl}$ to
$(\CC^\times/q^{pl\ZZ},\zeta_{pl})$. Then we have, on
$D(0,e^{-\pi/pl})$:
\begin{align*}
h^*\mu & \leq 
\frac{28e^{4\pi}}{\pi}\frac{1}{(1-e^{-\pi})^4}\cdot 
\frac{i}{2}d(q^{pl})d(\ol{q^{pl}}) \\ 
& \leq
\frac{28e^{4\pi}}{\pi}\frac{1}{(1-e^{-\pi})^4}(pl)^2{\cdot}
\frac{i}{2}dqd\ol{q}.
\end{align*}

In the second case, $q$ is mapped to
$(\CC^\times/q^{pl\ZZ},q^l,\zeta_l)$. In this case, we compose it with
the pseudo-involution $w_{\zeta_p}$, which brings us to the point
$(\CC^\times/q^{l\ZZ},\zeta_p,\zeta_l)$. The map then factors as the
$l$th power map from $D(0,e^{-\pi/pl})$ to $D(0,e^{-\pi/p})$, followed
by the map of Lemma~\ref{appl_merkl_lem1} composed with a suitable
diamond operator. We find:
\begin{align*}
h^*\mu & \leq 
\frac{28e^{4\pi}}{\pi}\frac{1}{(1-e^{-\pi/p})^4}\cdot 
\frac{i}{2}d(q^l)d(\ol{q^l}) \\
& \leq
\frac{28e^{4\pi}}{\pi}\frac{1}{(1-e^{-\pi/p})^4}l^2{\cdot}
\frac{i}{2}dqd\ol{q}.
\end{align*}

The third case is obtained by interchanging the roles of $p$ and~$l$,
so we will not make it explicit.

In the fourth case, $q$ is mapped to $(\CC^\times/q^{pl\ZZ},q)$. We
compose with the pseudo-involution $w_{\zeta_{pl}}$, which brings us
to $(\CC^\times/q^{\ZZ},\zeta_{pl})$. This is the map of
Lemma~\ref{appl_merkl_lem1}. We find:
\[
h^*\mu \leq 
\frac{28e^{4\pi}}{\pi}\frac{1}{(1-e^{-\pi/pl})^4}{\cdot}
\frac{i}{2}dqd\ol{q}.
\]

In these four cases, we see that the factor in front of
$(i/2)dqd\ol{q}$ in the upper bound for~$h^*\mu$ on $D(0,e^{-\pi/pl})$
is~$O((pl)^4)$. This gives the required estimate on~$gU$, as $dq$ on
$D(0,e^{-\pi/pl})$ corresponds to $dz_g$ on~$gU$.
\end{proof}

\begin{lem}\label{appl_merkl_lem3}
For the local coordinates $z^{(j)}$ and the real $(1,1)$-form $\mu'$
on $X(pl)$ as defined in~(\ref{appl_merkl_eqn2}) we have:
\[
\mu' \leq c_1 |dz^{(j)} \wedge d\ol{z}^{(j)}|
\]
on $U_1^{(j)}$ with $c_1 = c_1(pl)=O(pl)$.
\end{lem}
\begin{proof}
First of all, we have, by definition: $\mu'=(1/pl)h^*\mu$. The
definition of the $z^{(j)}$ (see~\ref{appl_merkl_eqn2}) plus the
definitions $z'=e^{3\pi/2pl}z$ and $\eps=e^{\pi/2pl}-1$ give:
\[
dz = \frac{2\eps}{3}e^{-3\pi/2pl}{\cdot}dz^{(j)}.
\]
If $pl$ gets large, then the factor $e^{3\pi/2pl}$ tends to~$1$  and 
for $\eps$ we have $\eps=(\pi/2pl)(1+O(1/pl))$. Combining all this with
Lemma~\ref{appl_merkl_lem2} finishes the proof.
\end{proof}
We can now finish the proof of Theorem~\ref{appl_merkl_thm1}. We apply
Theorem~\ref{merkl_thm1} on $X(pl)$ with the $(1,1)$-form $\mu'$. Then
we have $n=O((pl)^4)$, $M=6$ (Lemma~\ref{appl_merkl_lem0}) and
$c_1=O(pl)$ (Lemma~\ref{appl_merkl_lem3}). We obtain
from~(\ref{merkl_thm1_ineq1}) that there exists a constant $c$ such
that $g_{b,\mu'}(b')\leq c{\cdot}(pl)^5$ for all distinct primes $p$
and~$l$, and all distinct $b$ and~$b'$ on~$X(pl)$. For distinct $a$
and~$a'$ on $X_1(pl)$ we then have (see~(\ref{appl_merkl_eqn0})):
\[
g_{a,\mu}(a') = \sum_{h(b)=a}g_{b,\mu'}(h(a')) \leq c{\cdot}(pl)^6. 
\]
The statement that $\log \|dz^{(j)}\|_\mathrm{Ar}(P)= O((pl)^6)$
follows from the inequality~(\ref{merkl_thm1_ineq2}) in
Theorem~\ref{merkl_thm1}. Indeed, if locally we write $g_{a,\mu}$ as
$g_{a,\mu} = \log |z-z(a)| + f$ then $f(a)=-\log
\|dz\|_\mathrm{Ar}(a)$.
\end{proof} 

\section{Bounds for intersection numbers on~$X_1(\protect\lowercase{pl})$}
\label{bds_int_nmbrs}
In this section, we will bound the intersection numbers occurring
in the right hand side of the inequality in
Theorem~\ref{maininequality}, in the situation described in
Section~\ref{sec_setup_tau}.

\begin{thm}\label{thm_bds_int_nmbrs}
Let $p$ and $l$ be two distinct prime numbers, both at least~$5$, and
let $\calX$ be the semistable model over $B:=\Spec\ZZ[\zeta_{pl}]$
provided by~\cite{Katz-Mazur}. 
%%% Here, refer to the description given in this book.
For two cusps $P$ and $Q$ (possibly equal) in $\calX(B)$ we have:
\[
(P,P)\leq 0,\quad\text{and}\quad  |(P,Q)| = O((pl)^7) .
\]
For a cuspidal effective divisor $D$ of degree~$g$ on~$\calX$ we have:
\[
|(D,D-\omega_{\calX/B})| = O((pl)^{11}) .
\]
\end{thm}
\begin{proof}
As $p$ and $l$ are at least~$5$, the genus of $X_1(pl)$ is at least
two.  By the adjunction formula (see~(\ref{eqn_adjunction_formula})).
we have $-(P,P)=(P,\omega_{\calX/B})$, and by \cite{Faltings1},
Theorem~5, we have $(P,\omega_{\calX/B})\geq 0$, hence $(P,P)\leq0$.

Let us now derive an upper bound for $(P,\omega_{\calX/B})$. As the
automorphism group of~$\calX$ over~$B$ preserves the Arakelov
intersection product on~$\calX$ and acts transitively on the cusps, it
suffices to do this for the standard cusp~$\infty$. The Fourier
expansion at~$\infty$ of the rational function $j$ on~$\calX$ is of
the form $j=1/q + f$ with $f\in\ZZ[[q]]$. Therefore, $1/j$ is regular
in a neighbourhood of~$\infty$, and has a zero of order one
along~$\infty$. It follows that $d(1/j)$ generates
$\infty^*\omega_{\calX/B}$, and $d(1/j)=dq$
in~$\infty^*\omega_{\calX/B}$. By definition of the Arakelov
intersection product and Theorem~\ref{appl_merkl_thm1} we then have:
\[
(\infty,\omega_{\calX/B})=
-[\QQ(\zeta_{pl})):\QQ]\log\|dq\|_\mathrm{Ar}(\infty) = O((pl)^7).
\]
We now know $|(P,P)|=O((pl)^7)$ for all cusps $P$ in~$\calX(B)$. We
will now show that $|(P,Q)|=O((pl)^7)$. By the Theorem of
Manin-Drinfeld, see~\cite{Drinfeld1}, the image of the divisor $P-Q$
in $J_1(pl)(\QQ(\zeta_{pl}))$ is of finite order. Let $\Phi$ be a
vertical fractional divisor such that for any irreducible component $C$
of a fibre of $\calX$ over~$B$ we have $(P-Q-\Phi,C)=0$. By
\cite{Hriljac} or Theorem~4 of \cite{Faltings1} we have
$(P-Q-\Phi,P-Q-\Phi)=0$. Equivalently, we have:
\[
2(P,Q)=(P,P)+(Q,Q)-(P-Q,\Phi).
\]
The term $(P-Q,\Phi)$ can be dealt with by the method used in the
proof of Lemma~\ref{removingPhi}. We work this out in this special
situation. We make $\Phi$ unique by demanding that its support is
disjoint from~$P$. Of course, this does not change the number
$(P-Q,\Phi)$, but it makes it easier to talk about~$\Phi$. The support
of $\Phi$ is contained in the reducible fibers. 
%% Here we should refer to a place where the model is described and
%% the number of ss points is given.
These are exactly the fibers in the characteristics~$p$ and~$l$. Let
us estimate the contribution at the prime~$p$. We have to sum over the
maximal ideals of the $\FF_p$-algebra
$\FF_p[x]/(x^{l-1}+\cdots+1)$. All residue fields of this algebra are
isomorphic to a finite extension~$\FF$ of~$\FF_p$, and the number of
them is $(l-1)/\dim_{\FF_p}\FF$. Let $\calX_\FF$ be the fibre at one
of these residue fields, and let $\Phi_\FF$ be the part of $\Phi$ that
has support in~$\calX_\FF$. Then $\calX_\FF$ is the union of two
irreducible components $I$ and~$I'$, with transversal intersection at
the supersingular points. The number of $\Fbar_p$-valued supersingular
points is given by:
\[
s:=\#\calX_\FF(\Fbar_p)^\mathrm{s.s.} = 
\#\left(X_1(l)(\Fbar_p)^\mathrm{s.s.}\right) = \frac{(p-1)(l^2-1)}{24}.
\]
The degree on $I$ of the restriction to it of $\calO_\calX(I')$
is~$s$, and, symmetrically, $\deg_{I'}\calO_\calX(I)=s$. As
$I+I'=\calX_\FF$ is a principal divisor in a neighbourhood
of~$\calX_\FF$, the restrictions of $\calO_\calX(I+I')$ to $I$ and
$I'$ are trivial, we have $\deg_I\calO_\calX(I)=-s$ and
$\deg_{I'}\calO_\calX(I')=-s$. It follows that $\Phi_\FF$ is one of
the following fractional divisors: $\Phi_\FF=0$ if $P$ and $Q$
specialise to points on the same irreducible component of~$\calX_\FF$;
$\Phi_\FF=(1/s){\cdot}I$ if $P$ specialises to a point on $I'$ and $Q$
to a point on~$I$; $\Phi_\FF=(1/s){\cdot}I'$ if $P$ specialises to a
point on $I$ and $Q$ to a point on~$I'$. If we denote by
$(P-Q,\Phi)_\FF$ the contribution to $(P-Q,\Phi)$ at the
fibre~$\calX_\FF$, we have:
\[
|(P-Q,\Phi)_\FF| \leq (2/s){\cdot}\log\#\FF . 
\]
Summing this over the residue fields of $\FF_p[x]/(x^{l-1}+\cdots+1)$
gives, for the contribution $(P-Q,\Phi)_p$ to $(P-Q,\Phi)$ at~$p$:
\[
|(P-Q,\Phi)_p| \leq
 \frac{l{-}1}{\dim_{\FF_p}(\FF)}{\cdot}\frac{2}{s}{\cdot}\log\#\FF
= \frac{48\log p}{(p-1)(l+1)}. 
\]
%%%% include a picture of the dual graph ????
Likewise, we have, for the contribution at~$l$ to~$(P-Q,\Phi)$:
\[
|(P-Q,\Phi)_l| \leq \frac{48\log l}{(l-1)(p+1)}. 
\]
So, finally:
\[
|(P-Q,\Phi)| \leq 
\frac{48\log l}{(l-1)(p+1)} + \frac{48\log p}{(p-1)(l+1)}.
\]
The estimate $|(P,Q)|=O((pl)^7)$ now follows.

To get to the second statement of the theorem, note that
\begin{align*}
(D,D-\omega_{\calX/B}) & =
(D,D+\omega_{\calX/B})-2(D,\omega_{\calX/B}) \\
& = \sum_{k \neq l} (P_k,P_l) - 2(D,\omega_{\calX/B}) ,
\end{align*}
where $D=P_1+\cdots +P_g$, with repetitions allowed.  By our previous
estimates, we get
$|(D,D-\omega_{\calX/B})|=O(g^2(pl)^7)=O((pl)^{11})$.
\end{proof}
We will also need a lower bound for the intersection number of two
distinct points on~$X_1(pl)$.
\begin{thm}\label{thm_lower_bnd_gr_function}
There is an integer~$c$ such that for all pairs of distinct primes $p$
and $l$ such that $X_1(pl)$ has genus at least one, for any extension
$K$ of $\QQ(\zeta_{5l})$ and for $P$ and $Q$ distinct points in
$X_1(pl)(K)$ we have:
\[
\frac{1}{[K:\QQ]}(P,Q) \geq c (pl)^6, 
\]
where $(P,Q)$ is the Arakelov intersection number of $P$ and $Q$ on
the minimal regular model of $X_1(pl)$ over~$O_K$.
\end{thm}
\begin{proof}
We have:
\[
(P,Q) = (P,Q)_\fin + (P,Q)_\infty,
\]
with $(P,Q)_\fin$ the contribution from the finite places of~$K$, and
$(P,Q)_\infty$ the contribution from the infinite places. As $P\neq
Q$, we have $(P,Q)_\fin\geq 0$. On the other hand, we have:
\[
(P,Q)_\infty = \sum_\sigma -g_\sigma(P_\sigma,Q_\sigma).
\]
By Theorem~\ref{appl_merkl_thm1} we have:
\[
g_\sigma(P_\sigma,Q_\sigma) \leq c{\cdot}(pl)^6
\]
for some absolute constant~$c$. This finishes the proof.
\end{proof}

\section{A bound for
  $\protect\lowercase{h}(\protect\lowercase{x}'_{\protect\lowercase{l}}(Q))$  
in terms 
of~$\protect\lowercase{h}(\protect\lowercase{b_{\protect\lowercase{l}}}(Q))$}
\label{bds_height_x_and_y_in_height_b}
In this section we do what was promised at the beginning of
Section~\ref{sec_height_and_intersection}, by stating and proving the
following proposition and a corollary.

\begin{prop}\label{prop_height_in_X} 
There is a real number $c$ such that the following holds. Let $b$ be
in~$\Qbar$, such that $b^5(b^2+11b-1)\neq0$. Let $(u,v)$ in~$\Qbar^2$
be a torsion point on the elliptic curve $E_b$ over~$\Qbar$ given
by the equation:
\begin{eqn}
y^2+(b+1)xy + by = x^3 + bx^2 \, ,
\end{eqn}
i.e., on the fibre at $b$ of the universal elliptic curve with a point
of order~$5$ given in Proposition~\ref{prop_Y_1_5}. Then the absolute
heights $h(u)$ and $h(v)$ are bounded from above by $c+14h(b)$.
\end{prop}
\begin{proof}
Let $b$, $u$ and $v$ in $\Qbar$ be as in the statement of the
proposition.  We will now invoke known bounds for the difference
between the Weil height and the N\'eron-Tate height on elliptic curves
over number fields. Such bounds are given, for example,
in~\cite{Demjanenko1}, \cite{Zimmer},
and~\cite{Silverman2}. In~\cite{Silverman2} bounds are given for
elliptic curves given by general Weierstrass equations, but under the
assumption that the coefficients are algebraic integers. Therefore,
for us it seems better to use the following bound in Zimmer's theorem
on page~40 of~\cite{Zimmer}.
\begin{thm}[Zimmer]\label{thm_zimmer}
Let $E$ be an elliptic curve over~$\Qbar$ given by a Weierstrass
equation $y^2=x^3+Ax+B$, and let $P\in E(\Qbar)$. Let $h$ denote the
absolute Weil height on~$\PP^2(\Qbar)$, and $\hat{h}$ the absolute
Néron-Tate height on~$E$ attached to~$h$. Then one has:
\begin{multline*}
-2^{-1}(2^{-1}h(1:A^3:B^2)+7\log 2) \leq h(P)-\hat{h}(P) \leq \\
\leq 2^{-1}h(1:A^3:B^2)+6\log 2\, .
\end{multline*}
\end{thm}
So, we must compare our plane elliptic curve~$E_b$ with one given by a
standard Weierstrass equation. We put:
\[
v_1:=v+((b+1)u+b)/2,\quad u_1:=u+(b+(b+1)^2/4)/3\, .
\]
Then $(u_1,v_1)$ is a point on the elliptic curve $E'_b$ given by a
Weierstrass equation $y^2=x^3+Ax+B$, with $A$ and~$B$ polynomials
in~$b$, with coefficients in~$\QQ$, of degrees at most~$4$ and~$6$,
respectively. We note that $A$ and $B$ depend only on~$b$, not
on~$(u,v)$. Using Lemma~\ref{lem_sum_and_product} and writing $A$ as
$A_0+b(A_1+b(A_2+b(A_3+bA_4)))$, we see that there is a real
number~$c_1$, such that for all~$b$ we have $h(A)\leq
c_1+4h(b)$. Similarly, there is a $c_2$ such that $h(B)\leq
c_2+6h(b)$. Therefore, there is a $c_3$ such that $h(1:A^3:B^2)\leq
c_3+24h(b)$. Zimmer's theorem~\ref{thm_zimmer}, plus the fact that the
N\'eron-Tate height of torsion points is zero, tells us that there is
a $c_4$ such that for all~$b$ and for all torsion points $(u,v)$
on~$E_b$, we have:
\[
h(u_1,v_1)\leq c_4+12h(b)\, .
\]
Expressing $u$ and $v$ in $u_1$ and $v_1$, and using again
Lemma~\ref{lem_sum_and_product}, we get a real number $c_5$ such that
for all $b$ and all $(u,v)$ as in the proposition we are proving, we
have:
\[
h(u)\leq c_5+14h(b),\quad h(v)\leq c_5+14h(b)\, .
\]
This ends the proof of Proposition~\ref{prop_height_in_X}.
\end{proof}
\begin{cor}\label{cor_heightx'l_of_Q}
There is a real number $c$ such that for each~$l$ and each $Q_{x,i}$
as in the beginning of Section~\ref{sec_height_and_intersection} we
have:
\[
 h(x'_l(Q_{x,i}))\leq c+14h(b_l(Q_{x,i})).
\]
\end{cor}
\begin{proof}
By definition (Section~\ref{sec_setup_tau}), $x'_l(Q_{x,i})$ is the
$x$-coordinate of a point of order~$5l$ on the fibre at $b_l(Q_{x,i})$
of the elliptic curve~$E$ in
Proposition~\ref{prop_Y_1_5}. Proposition~\ref{prop_height_in_X} gives
the result.
\end{proof}

\section{An integral over $X_1(\protect\lowercase{5l})$}
In this section we will give an upper bound for the integral appearing
in the estimate in Theorem~\ref{thm_bound_hbQ}. We recall the
situation: $l$ is a prime number with $l>5$, and $b$ is the regular
function on $Y_1(5)_\QQ$ given by Proposition~\ref{prop_Y_1_5}. We
also view $b$ as a regular function $b_l$ on~$Y_1(5l)_\QQ$ via pullback
along the map $Y_1(5l)_\QQ\to Y_1(5)_\QQ$ that sends $(E/S,P_5,P_l)$
to $(E/S,P_5)$, where $S$ is any $\QQ$-scheme, $E/S$ an elliptic
curve, $P_5$ in $E(S)$ everywhere of order~$5$, and $P_l$ in $E(S)$
everywhere of order~$l$. 

\begin{prop}\label{prop_bound_integral}
There exist real numbers $A$ and $B$ such that for all primes $l>5$:
\[
\int_{X_1(5l)(\CC)} \log(|b_l|^2 +1) \, \mu \leq A+B\cdot l^6\, ,
\]
where $\mu$ is the Arakelov $(1,1)$-form. 
\end{prop}
\begin{proof} 
In order to simplify the notation in the proof, we will let $X_1(5l)$
denote the Riemann surface $X_1(5l)_\QQ(\CC)$ of complex points of the
curve $X_1(5l)_\QQ$ over~$\QQ$, and we will drop de subsscript~$l$
in~$b_l$. We will denote points of $Y_1(5l)$ as triples $(E,P_5,P_l)$
with $E$ a complex elliptic curve, with points $P_5$ of order~$5$ and
$P_l$ of order~$l$. Similarly, we will denote points of $Y_1(5)$ by
pairs $(E,P)$, and points on the $j$-line just by elliptic curves.

The equations in Proposition~\ref{prop_Y_1_5} show that the rational
function $b$ on $X_1(5)$ has exactly one pole, that it is of order
one, and that at that point, the function~$j\colon X_1(5)\to\PP^1$ has
a pole of order~$5$. The region in $\HH$ consisting of the $\tau$ with
$\Im(\tau)>1$ gives an embedding of the disk $D(0,e^{-2\pi})$ into the
$j$-line, sending $q\neq 0$ to~$\CC^\times/q^\ZZ$. The inverse image
of this disk under $j\colon X_1(5)\to\PP^1$ consists of $4$ disks, one
around each cusp. The two disks around the cusps where $j$ is ramified
are given by the embeddings $D(0,e^{-2\pi/5})\to X_1(5)$, sending
$q\neq 0$ to $(\CC^\times/q^{5\ZZ},q)$ and to
$(\CC^\times/q^{5\ZZ},q^2)$. As the integral in the proposition that
we are proving does not change if we replace $b$ by its image under a
diamond operator $\ld a\rd$ with $a\in\FF_5^\times$, we may and do
suppose that $b$ has its pole at the center of the punctured disk
\[
U = \{ (\CC^\times/q^{5\ZZ}, q) : q\in D(0,e^{-2\pi/5})^*\}
\subset X_1(5)\, .
\]
The integral of $\log(|b|^2 +1)\mu$ over the complement of the inverse
image of $U$ in~$X_1(5l)$ is bounded by the supremum of $\log(|b|^2
+1)$ on $X_1(5)-U$. This upper bound does not depend on~$l$. Hence it
is enough to prove that the integral of $\log(|b|^2 +1)\mu$ over the
inverse image of $U$ in $X_1(5l)$ is bounded by $B \cdot l^6$, for a
suitable~$B$. Now this inverse image of $U$ is a union of $2(l-1)$
punctured disks $U_c$ and $V_d$ with $c$ and $d$ running
through~$\FF_l^\times$ given by:
\[
U_c=\{(\CC^\times/q^{5\ZZ}, q, \zeta_l^c) : q\in D(0,e^{-2\pi/5})^*\}
\]
and:
\[
V_d=\{(\CC^\times/q^{5l\ZZ}, q^l, q^{5d}) : q\in D(0,e^{-2\pi/5l})^*\}\, .
\] 
For each $c$ in $\FF_l^\times$ the map $X_1(5l)\to X_1(5)$ restricts
to the isomorphism $U_c\to U$ given by $q\mapsto q$. For each $d$ in
$\FF_l^\times$, the restriction $V_d\to U$ is given by $q\mapsto q^l$.

Around the standard unramified cusp of $X_1(5l)$ we have the punctured
disk:
\[
W_l = \{ (\CC^\times/q^\ZZ,\zeta_5,\zeta_l) : q\in D(0,e^{-2\pi/5})^*\}.
\]
By applying a suitable element from the group of automorphisms of
$X_1(5l)$ generated by the Atkin-Lehner pseudo-involutions and the
diamond operators we can establish isomorphisms of the $U_c$
with~$W_l$. The point $(\CC^\times/q^{5\ZZ},q,\zeta_l^c)$ in $U_c$ is
then first sent to $(\CC^\times/q^{5\ZZ},q,\zeta_l)$ by~$\ld
c^{-1}\rd$, and then to $(\CC^\times/q^\ZZ,\zeta_5,\zeta_l)$ by the
Atkin-Lehner pseudo-involution that divides out by the group generated
by the point of order~$5$ (see Definition~\ref{defi_w_operator}).  The
coordinate $q$ on $W_l$ is then identified in this way with the
coordinate $q$ on~$U_c$. We observe that $\mu$ is invariant under each
automorphism of~$X_1(5l)$. Lemma~\ref{appl_merkl_lem1} gives a real
number~$C_1$ such that for all~$l$ and~$c$, the positive real
$(1,1)$-form $\mu$ on $U_c$ can be estimated from above by:
\[
\mu|_{U_c} \leq C_1 \cdot \frac{i}{2} dq d\ol{q} \, .
\]
Similarly, for each $d$ in~$\FF_l^\times$ a suitable automorphism of
$X_1(5l)$ maps $V_d$ to the punctured disk~$W'_l$:
\[
W'_l = \{ (\CC^\times/q^\ZZ,\zeta_5,\zeta_l) : 
q\in D(0,e^{-2\pi/5l})^*\}\, .
\]
Lemma~\ref{appl_merkl_lem1} gives a real number~$C_2$ such that for
all~$l$ and~$d$, the positive real $(1,1)$-form $\mu$ on $V_d$ can be
estimated from above by:
\[
\mu|_{V_d} \leq 
C_2\cdot\frac{1}{(1-e^{-2\pi/5l})^4}\cdot\frac{i}{2} dq d\ol{q} \, .
\]
Now the function $qb$ on $U$ extends to a holomorphic function on a
disk containing~$U$, hence $|qb|$ is bounded on~$U$. Hence there is a
real number $C_3>1$ such that $|b|^2+1 \leq C_3{\cdot}|q^{-2}|$ on~$U$. It
follows that on all $U_c$ we have $|b|^2+1 \leq C_3{\cdot}|q^{-2}|$, and
that on all $V_d$ we have $|b|^2+1 \leq C_3{\cdot}|q^{-2l}|$ (note that
under $V_d\to U$, $(\CC^\times/q^{5l\ZZ},q^l,q^{5d})$ is sent to
$(\CC^\times/q^{5l\ZZ},q^l)$). We remark that
$1/(1-e^{-x})=x^{-1}(1+O(x))$ as $x$ tends to~$0$ from above. Hence
there is a $C_4\in\RR$ such that for all~$l$:
\[
\frac{1}{(1-e^{-2\pi/5l})^4} \leq C_4{\cdot}l^4\, .
\]
We get, for all~$l$:
\begin{multline*} 
\int_{b^{-1}U} \log(|b|^2+1) \, \mu  \leq  
(l-1)C_1\int_{D(0, e^{-2\pi/5})} \log(C_3{\cdot}|q^{-2}|) \cdot
\frac{i}{2} dq d\ol{q} \\
  + (l-1)C_2 \frac{1}{(1-e^{-2\pi/5l})^4} \int_{D(0,e^{-2\pi/5l})}
\log(C_3{\cdot}|q^{-2l}|)\cdot\frac{i}{2} dq d\ol{q} \leq\\
 \leq  \int_{|z| < 1}\left( lC_1\log(C_3{\cdot}|z|^{-2})
  + C_2C_4 l^5 \log(C_3{\cdot}|z|^{-2l})\right)\frac{i}{2}dz d\ol{z} \\ 
 \leq  \int_{|z| < 1}\left( C_1l\log(C_3{\cdot}|z|^{-2})
  + C_2C_4l^6\log(C_3{\cdot}|z|^{-2})\right) \frac{i}{2} dz d\ol{z} \\
 =  (C_1l+C_2C_4l^6){\cdot}(\pi\log C_3 + \pi)\, .
\end{multline*} 
This finishes the proof of Proposition~\ref{prop_bound_integral}
\end{proof}

\section{Final estimates of the Arakelov contribution}
\label{sec_final_estimates}
We will now put the estimates of the preceding sections together, in
the situation of Section~\ref{sec_setup_tau}. We briefly recall this
situation. We have a prime number $l>5$, and $X_l$ denotes the modular
curve~$X_1(5l)$, over~$\QQ$, and $g_l$ its genus. The Jacobian variety
of $X_l$ is denoted by~$J_l$. In $J_l(\Qbar)[l]$ we have the
$\Gal(\Qbar/\QQ)$-module~$V$ that realises the representation
$\rho$ from $\Gal(\Qbar/\QQ)$ to $\GL_2(\FF)$ attached to a surjective
ring morphism $f\colon\TT(1,k)\onto\FF$ such that the image of $\rho$
contains~$\SL_2(\FF)$. We have an effective divisor~$D_0$
on~$X_{l,\QQ(\zeta_l)}$, of degree~$g_l$, supported on the cusps. For
every $x$ in~$V$ there is a unique effective divisor
$D_x=Q_{x,1}+\cdots+Q_{x,g_l}$ of degree~$g_l$ such that $x$ is the
class of $D_x-D_0$. We have written $D_x=D_x^\fin+D_x^\cusp$, where
$D_x^\cusp$ is the part of $D_x$ supported on the cusps. The numbering
of the $Q_{x,i}$ is such that $D_x^\fin=Q_{x,1}+\cdots+Q_{x,d_x}$. We
have morphisms $b_l$ and $x'_l$ from $X_{l,\QQ}$ to~$\PP^1_\QQ$ that,
seen as rational functions, have their poles contained in the set of
cusps of~$X_l$.

The following proposition gives upper bounds for the absolute heights
of the algebraic numbers $b_l(Q_{x,i})$ and~$x'_l(Q_{x,i})$,
polynomial in~$l$. The height function~$h$ used here is as defined
in~(\ref{eqn_height_on_qbar}). The proof of the proposition combines
the involved arguments of the previous sections.
\begin{prop}\label{prop_arakelov_contrib_1}
There is an integer~$c$, such that for all $x$ in $V$ as above, and
for all $i\in\{1,\ldots,d_x\}$, the absolute heights of~$b_l(Q_{x,i})$
and~$x'_l(Q_{x,i})$ are bounded from above by~$c{\cdot}l^{12}$.
\end{prop}
\begin{proof}
We will just write $b$ for~$b_l$. Let $V$ be as in the
proposition. Theorem~\ref{thm_bound_hbQ} shows that for all $x\in V$,
all $i\in\{1,\ldots,d_x\}$ and all number fields $K$ containing
$\QQ(\zeta_{5l})$ over which $Q_{x,i}$ is rational, we have:
\begin{eqn}\label{eqn_arak_contrib_1.1}
\begin{aligned}
h(b(Q_{x,i})) \leq & \frac{1}{[K:\QQ]} 
\left( (Q_{x,i},b^*\infty)_\calX +
l^2  \sum_\sigma \sup_{X_\sigma} g_{\sigma} + 
\vphantom{\frac{1}{2}\sum_\sigma\int_{X_\sigma}}\right.\\
& \quad \left. 
+ \frac{1}{2}\sum_\sigma\int_{X_\sigma}\log(|b|^2+1)\mu_{X_\sigma}\right) 
+ \frac{1}{2}\log 2\, .    
\end{aligned}
\end{eqn}
Here $\calX$ is the minimal regular model of $X_l$ over $B:=\Spec
O_K$. 

Let us first concentrate on the second and third terms of the right
hand side
of~(\ref{eqn_arak_contrib_1.1}). Theorem~\ref{appl_merkl_thm1} gives
an integer $c_1$ such that for all $\sigma$ we have $\sup_{X_\sigma}
g_\sigma \leq c_1{\cdot}l^6$, uniformly in
all~$l$. Proposition~\ref{prop_bound_integral} gives an integer $c_2$
such that for all $\sigma$ we have $(1/2)\int_{X_{\sigma}}
\log(|b|^2+1)\mu_{X_\sigma}\leq c_2{\cdot}l^6$, uniformly in all~$l$.
Hence, for all $x\in V$, all $i\in\{1,\ldots,d_x\}$ and all number
fields $K$ containing $\QQ(\zeta_{5l})$ over which $Q_{x,i}$ is
rational, we have:
\begin{eqn}\label{eqn_arak_contrib_1.2}
h(b(Q_{x,i})) \leq \frac{1}{[K:\QQ]}(Q_{x,i},b^*\infty)_\calX +
c_1{\cdot}l^8 + c_2{\cdot}l^6\, ,
\end{eqn}
uniformly in all~$l$ and~$V$.

We now concentrate on the first term. We recall that $b^*\infty$ is an
effective cuspidal divisor on~$X_l$. Let $x$ be in~$V$. We apply
Theorem~\ref{maininequality}, where (in the notation of that theorem)
we assume that all $Q_{x,i}$ ($i\in\{1,\ldots,g_l\}$) are $K$-rational
and that $K$ contains~$\QQ(\zeta_{5l})$, and that $P$ is a cusp. We
also use the obvious fact that $\log\#\rR^1p_*O_\calX(D_x)$ is
nonnegative. That gives:
\begin{multline*}
\frac{(D_x,P)}{[K:\QQ]}  \leq 
\frac{1}{[K:\QQ]} \left( -\frac{1}{2} (D_0,D_0-\omega_{\calX/B}) + 
2g_l^2\sum_{s \in B} \delta_s \log \#k(s) \right. \\
+ \sum_\sigma \log \|\vartheta\|_{\sigma,\sup} + 
\frac{g_l}{2}[K:\QQ]\log(2\pi) \\ 
+ \left. \frac{1}{2} \deg \det p_* \omega_{\calX/B} + (D_0,P)\right) \, .
\end{multline*} 
Theorem~\ref{thm_bds_int_nmbrs}, applied with
$B=\Spec(\ZZ[\zeta_{5l}])$, gives that:
\[
\frac{1}{[K:\QQ]}|(D_0,D_0-\omega_{\calX/B})| = O(l^{10}), 
\]
and that:
\[
\frac{1}{[K:\QQ]} |(D_0,P)| = O(l^8),
\]
as $D_0$ is an effective cuspidal divisor of degree $g_l$ and
$g_l=O(l^2)$.
By Theorem~\ref{thm_bound_theta} we have:
\[
\frac{1}{[K:\QQ]} \sum_\sigma \log \|\vartheta\|_{\sigma,\sup} 
= O(l^6).
\]
By Theorem~\ref{thm_height} we have:
\[
\frac{1}{[K:\QQ]} \deg \det p_* \omega_{\calX/B} = O(l^2\log l).
\]
Finally, we have:
\[
\frac{g_l^2}{[K:\QQ]}\sum_{s \in B} \delta_s \log \#k(s) = O(l^6), 
\]
by the following argument. The only non-trivial contributions come
from $s$ over~$5$ and over~$l$. The total contribution at~$5$ is
independent of which extension~$K$ of $\QQ(\zeta_5)$ we use, and for
$\QQ(\zeta_5)$ there is only one $s$ over~$5$, $k(s)=\FF_5$, and
$\delta_s$ equals the number of supersingular points
in~$X_1(l)(\Fbar_5)$, which is~$O(l^2)$. The contribution from~$l$ can
be computed over~$\QQ(\zeta_l)$. Then there is one~$s$, and
$k(s)=\FF_l$, and $\delta_s$ is the number of supersingular points
in~$X_1(5)(\Fbar_l)$, which is~$O(l)$ (see Section~\ref{sec_constr_D}).

Putting these last estimates together, we get that there is an integer
$c_3$ such that for all~$l$ and~$x$ we have:
\begin{eqn}\label{eqn_upperbnd_D'P}
\frac{1}{[K:\QQ]}(D_x,P) \leq c_3{\cdot}l^{10}.
\end{eqn}
Now $D_x$ is the sum of $g_l$ points~$Q_{x,i}$. In order to get upper
bounds for the individual $(Q_{x,i},P)$ we need a lower bound for these.
Theorem~\ref{thm_lower_bnd_gr_function} gives us a lower bound if
$Q_{x,i}\neq P$ and Theorem~\ref{thm_bds_int_nmbrs} gives us one if
$P=Q_{x,i}$. Putting these together, we get an integer $c_4$
such that for all $l$, $x$ and~$i$:
\[
\frac{1}{[K:\QQ]}(Q_{x,i},P) \geq c_4{\cdot}l^6.
\]
The last two estimates together imply that there is an integer $c_5$
such that for all~$l$, $x$ and~$i$ we have:
\[
\frac{1}{[K:\QQ]}(Q_{x,i},P) \leq c_5{\cdot}l^{10}.
\]
As $b^*\infty$ is an effective cuspidal divisor on~$X_l$, of
degree~$O(l^2)$, the previous inequality implies that there is an
integer $c_6$ such that for all $l$, $x$ and
$i\in\{1,\ldots,d_x\}$ we have:
\[
\frac{1}{[K:\QQ]}(Q_{x,i},b^*\infty)_\calX
\leq c_6{\cdot}l^{12}.
\]
This finishes the proof concerning the height
of~$b_l(Q_{x,i})$. Corollary~\ref{cor_heightx'l_of_Q} then finishes
the proof.
\end{proof}
We recall, from the end of Section~\ref{sec_setup_tau}, that, in the
situation as described in the beginning of this section, we take a
linear combination $f_l:=b_l+nx'_l$ with $0\leq n\leq g_l^2(\#\FF)^4$,
such that under the map $f_l\colon X_{l,\Qbar}\to\PP^1_\Qbar$ the
divisors~$D_x^\fin=Q_{x,1}+\cdots+Q_{x,d_x}$, for $x\in V$, have
distinct images~$f_{l,*}D_x^\fin$. Suppose that $f_l$ is any such
linear combination. The $f_{l,*}D_x^\fin$ are then distinct effective
divisors of degree~$d_x$ on~$\AA^1_{\Qbar}$.

For each $x$ in~$V$, we get a polynomial $P_{D_0,f_l,x}$ with
coefficients in $\Qbar$ given by:
\[
P_{D_0,f_l,x}(t) = \prod_{i=1}^{d_x} (t-f_l(Q_{x,i})) 
\quad \text{in $\Qbar[t]$,} 
\]
and the map that sends $x$ to $P_{D_0,f_l,x}$ is injective, and
$\Gal(\Qbar/\QQ(\zeta_l))$-equivariant.

We have seen that there is an integer $m$ with $0\leq m\leq
g{\cdot}(\#\FF)^4$ such that the map:
\[
a_{D_0,f_l,m}\colon V\to \Qbar, \quad x\mapsto 
P_{D_0,f_l,x}(m) = \prod_{i=1}^{d_x} (m-f_l(Q_{x,i}))
\]
is injective and hence a generator of the $\QQ(\zeta_l)$-algebra
$A_{l,\QQ(\zeta_l)}$ associated with~$V$. Assume that $m$ is such an
integer. The following theorem gives our final upper bound
for the absolute height of the coefficients of the minimal polynomial:
\begin{eqn}\label{eqn_min_pol_a_3}
P_{D_0,f_l,m} = \prod_{x\in V}(T-a_{D_0,f_l,m}(x)) = \sum_jP_jT^j 
\quad\text{in $\QQ(\zeta_l)[T]$}  
\end{eqn}
of $a_{D_0,f_l,m}$ over~$\QQ(\zeta_l)$.

\begin{thm}\label{thm_arakelov_contrib_2}
There exists an integer $c$  such that for all $l$, $V$, $D_0$, $f_l$
and $m$ as above we have, for all~$P_j$ as in~(\ref{eqn_min_pol_a_3}):
\[
h(P_j)\leq c{\cdot}l^{14}{\cdot}(\#\FF)^2.
\]
\end{thm}
\begin{proof}
Let $c_1$ be an integer as given by
Proposition~\ref{prop_arakelov_contrib_1}. Let $l$, $V$, $D_0$, $f_l$
and $m$ be as in the theorem.  For each $x\in V$ and each $i$ in
$\{1,\ldots,d_x\}$ we have, by
Proposition~\ref{prop_arakelov_contrib_1}, the definition of~$f_l$,
and Lemma~\ref{lem_sum_and_product}:
\begin{eqn}\label{eqn_ineq_hflQxi}
\begin{aligned}
h(f_l(Q_{x,i})) & = h(b_l(Q_{x,i})+nx'_l(Q_{x,i})) \\
& \leq \log 2 + h(b_l(Q_{x,i})) + h(nx'_l(Q_{x,i})) \\
& \leq \log 2 + h(b_l(Q_{x,i})) + h(n) + h(x'_l(Q_{x,i})) \\
& \leq \log 2 + c_1{\cdot}l^{12} + \log\left(l^2{\cdot}(\#\FF)^4\right) 
+ c_1{\cdot}l^{12} \\
& \leq 2c_1{\cdot}l^{12} + 2\log l + 4\log(\#\FF) + \log 2 \\
& \leq c_2{\cdot}l^{12} + 4\log(\#\FF),
\end{aligned}
\end{eqn}
for $c_2=c_1+1$.

In order to simplify the notation during the rest of this proof, we
write $a$ for $a_{D_0,f_l,m}$ and $P$ for~$P_{D_0,f_l,m}$. Then we
have, for each $x\in V$ (using $d_x\leq l^2$
and~(\ref{eqn_ineq_hflQxi})):
\begin{eqn}\label{eqn_h_ax}
\begin{aligned}
h(a(x)) & = h\left(\prod_{i=1}^{d_x}\left(m-f_l(Q_{x,i})\right)\right)
\leq \sum_{i=1}^{d_x}h(m-f_l(Q_{x,i})) \\
& \leq \sum_{i=1}^{d_x}\left(\log 2 + h(m) + h(f_l(Q_{x,i}))\right) \\
& \leq d_x{\cdot}\left(\log 2 + \log\left(l^2(\#\FF)^4\right) +
c_2{\cdot}l^{12} + 4\log(\#\FF)\right)\\
& \leq c_3{\cdot}l^{14} + 8l^2\log(\#\FF),
\end{aligned}
\end{eqn}
where $c_3=c_2+1$.

Let $j$ be in $\{0,\ldots,\#V\}$. Then $P_j$ is, up to a sign, the
value of the elementary symmetric polynomial of degree $\#V-j$
evaluated in the $a(x)$, where $x$ ranges
through~$V$. Lemma~\ref{lem_bound_sym}, together
with~(\ref{eqn_h_ax}), gives us:
\[
\begin{aligned}
h(P_j) & \leq \#V{\cdot}\log 2 + \sum_{x\in V}h(a(x))\\
& \leq \#V{\cdot}\left(\log 2 + c_3{\cdot}l^{14}+8l^2\log(\#\FF)\right)\\
& \leq c{\cdot}(\#\FF)^2{\cdot}l^{14},
\end{aligned}
\]
where $c=c_3+1$, and where we have used that $\log\#\FF\leq l\log
l$ since $\FF$ is a quotient of~$\TT(1,k)$ and $k\leq l+1$. 
\end{proof}

A last consequence of all Arakelovian estimates is the following upper bound
for the term $\log\#\rR^1p_*O_\calX(D_x)$ in Theorem~\ref{maininequality}.
\begin{thm}\label{thm_boundR1p*}
There is an integer~$c$ such that for all $l>5$ prime and all~$x$ in
$V$ we have:
\[
\frac{1}{[K:\QQ]}\log\#\rR^1p_*O_\calX(D_x) \leq c{\cdot}l^{10}.
\] 
\end{thm}
\begin{proof}
We have already seen, in the proof of
Proposition~\ref{prop_arakelov_contrib_1}, that the right hand side of
the inequality in Theorem~\ref{maininequality}, divided by~$[K:\QQ]$,
is bounded from above by a constant times~$l^{10}$. We have also seen
in Theorem~\ref{thm_lower_bnd_gr_function} that the term
$(D_x,P)/[K:\QQ]$ on the left hand side is bounded from below by a
constant times~$l^8$ (recall that $D_x$ is of degree~$O(l^2)$). This
proves the inequality.
\end{proof}
This upper bound will be very useful for us, as the next
interpretation shows. Recall that $X_l$ has good reduction
over~$\ZZ[1/5l]$. 
\begin{thm}\label{thm_exist_good_places}
There is an integer $c$ with the following property. Let $l$, $X_l$,
$V$ and $D_0$ as in the beginning of this section. A prime
number~$p\not|5l$ is said to be \emph{$V$-good} if for all $x$ in~$V-\{0\}$
the following two conditions are satisfied:
\begin{enumerate}
\item at all places $v$ of $\Qbar$ over~$p$ the specialisation
$D_{x,{\Fbar_p}}$ at $v$ is the unique effective divisor on the
reduction~$X_{l,\Fbar_p}$ such that the difference with $D_{0,\Fbar_p}$
represents the specialisation of~$x$;
\item the specialisations of the non-cuspidal part $D_x^\fin$ of $D_x$
at all~$v$ above~$p$ are disjoint from the cusps.
\end{enumerate}
Then we have:
\[
\sum_{\text{$p$ not $V$-good}} \log p \leq c l^{12}{\cdot}(\#\FF)^2.
\]
\end{thm}
\begin{proof}
First of all, a prime number~$p$ satisfies conditions~(1) and~(2) for
all~$x$ in $V-\{0\}$ if and only if it satisfies them one of them, as
$\Gal(\Qbar/\QQ(\zeta_l))$ acts transitively on $V-\{0\}$ by
assumption.

We take $K$ to be the extension of $\QQ(\zeta_l)$ that corresponds to
the transitive $\Gal(\Qbar/\QQ(\zeta_l))$-set $V-\{0\}$; this is the
field of definition of one~$x$ in~$V-\{0\}$. Then
$[K:\QQ]=O(l{\cdot}(\#\FF)^2)$. We define $S(V)$ to be the image
in~$\Spec(\ZZ[1/5l])$ of the support in $\Spec(O_K)$ of the finite
$O_K$-module $\rR^1p_*O_\calX(D_x)$. By Theorem~\ref{thm_boundR1p*}
we have:
\[
\log\# (\ZZ[1/5l]\otimes\rR^1p_*O_\calX(D'_x)) = O(l^{11}{\cdot}(\#\FF)^2),
\]
and hence:
\[
\sum_{p\in S(V)} \log p = O(l^{11}{\cdot}(\#\FF)^2).
\]
We claim that the $p$ in $S(V)$ are precisely the primes
$p\not\in\{5,l\}$ such that condition~(1) is not satisfied for~$x$. To
see this, we first note that for a morphism $O_K\to\Fbar_p$ the
canonical map from $\Fbar_p\otimes_{O_K}\rR^1p_*O_\calX(D_x)$ to
$\rH^1(X_{l,\Fbar_p},O_{\calX_{\Fbar_p}}(D_{x,\Fbar_p}))$ is an
isomorphism (see Theorem~III.12.1 of~\cite{Hartshorne1}; base change
and cohomology in top dimension commute). The divisor $D_{x,\Fbar_p}$ is
the unique effective divisor in its linear equivalence class if and
only if $h^0(X_{l,\Fbar_p},O_{\calX_{\Fbar_p}}(D_{x,\Fbar_p}))=1$, which, by
Riemann-Roch, is equivalent to
$h^1(X_{l,\Fbar_p},O_{\calX_{\Fbar_p}}(D_{x,\Fbar_p}))=0$.

Now we let $T(V)$ denote the set of primes $p\not\in\{5,l\}$ such that
at least one specialisation of $D_x^\fin$ at a place of~$K$ above $p$
is not disjoint from the cusps. Taking into account that
$[K:\QQ]=O(l{\cdot}(\#\FF)^2)$, equation~(\ref{eqn_upperbnd_D'P})
gives us an upper bound:
\[
(D_x,P) = O(l^{11}{\cdot}(\#\FF)^2)
\]
As the degree of $D_x^\cusp$ is at most~$O(l^2)$,
Theorem~\ref{thm_bds_int_nmbrs} gives us:
\[
|(D_x^\cusp,P)| = O(l^7{\cdot}(\#\FF)^2), \quad\text{hence}\quad 
(D_x^\fin,P) = O(l^{11}{\cdot}(\#\FF)^2).
\]
As the divisor $\Cusps$ has degree~$O(l)$, we have:
\[
(D_x^\fin,\Cusps) = O(l^{12}{\cdot}(\#\FF)^2).
\]
The intersection number $(D_x^\fin,\Cusps)$ is the sum of
$(D_x^\fin,\Cusps)_\fin$, the contribution of the finite places, and
$(D_x^\fin,\Cusps)_\infty$, the contribution of the infinite
places. We have:
\[
(D_x^\fin,\Cusps)_\infty = \sum_{i,P,\sigma} -g_\sigma(Q_{x,i},P),
\]
where the sum is taken over the~$i$ with $1\leq i\leq d_x$, over the
cusps~$P$ and the $\sigma\colon K\to \CC$. 
Then Theorem~\ref{appl_merkl_thm1} gives the upper bound:
\[
(D_x^\fin,\Cusps)_\fin = (D_x^\fin,\Cusps) - (D_x^\fin,\Cusps)_\infty 
= O(l^{12}{\cdot}(\#\FF)^2).
\]
By the definition of our set $T(V)$ and the definition of
$(D_x^\fin,\Cusps)_\fin$ we get:
\[
\sum_{p\in T(V)} \log p \leq (D_x^\fin,\Cusps)_\fin 
= O(l^{12}{\cdot}(\#\FF)^2).
\]
The proof of the theorem is then finished by noticing that the set of
primes $p\not\in\{5,l\}$ that are not $V$-good is precisely the union
of $S(V)$ and~$T(V)$.
\end{proof}

%JMCdebut
%% -- Computing V_\ell over the complex numbers
\chapter{Approximating $V_f$ over the complex numbers}\label{sec_couveignes_TORSION}

\author{J.-M. Couveignes}%JMCdebut

\bigskip

\bigskip

In this chapter, we address the problem of computing 
torsion divisors on modular curves  with an application
to the explicit calculation of modular representations.
We assume we  are given an even integer $k > 2$,  a prime integer $l>6(k-1)$, a finite field $\FF$
with characteristic $l$, and a ring epimorphism $f : \TT(1,k)\rightarrow  \FF$. We want to compute
the associated Galois representation $\rho_f : \Gal(\Qbar/\QQ)\rightarrow \GL_2(\FF)$.  This
representation lies in the jacobian of the  modular curve $X_1(l)$. Indeed,  let $f_2 :  \TT(l,2)
\rightarrow \FF$ be the unique ring homomorphism such that $f_2(T_m)=f(T_m)$ for every positive integer $m$. Let 
$V_f$ be the subgroup of $J_1(l)(\Qbar)[l]$ cut out by the kernel of $f_2$. This is a dimension $2$
vector space over $\FF$. Given a finite generating set $(t_1, \ldots, t_r)$ with $r=(l^2-1)/6$ as in Theorem~\ref{thm_red_to_wt_2_Delta} we can rewrite
$V_f$ as a finite intersection 

\begin{equation}\label{eq:defVltor}
V_f = \bigcap_{1\le i \le r} \ker\left(t_i, J_1(l)(\Qbar)[l]\right).
\end{equation}

 Then $V_f$ realizes  $\rho_f$ and we may
write $\rho_f$ as a morphism

\[
        \rho_f :   \Gal(\Qbar/\QQ) \to \GL (V_f).
\]

We will assume that the image of $\rho_f$ contains $\SL  (V_f)$. Otherwise $\rho_f$ would be reducible according to Theorem~\ref{thm_large_image}. 
And the reducible case is treated in Section~\ref{sec_red_to_irr}.

We want to compute the  splitting field $K_f$  of $\rho_f$ 
as an extension of $\QQ$. 
Computing $V_f$ using general algorithms from computer algebra, 
like Buchberger's
algorithm, seems difficult in this situation  because
$V_f$ is defined as a subset of the $l$-torsion subgroup
inside $J_1(l)$. A naive algebraic description of $V_f$
would lead us to write down an equation for $J_1(l)[l]$; something 
similar to an  $l$-division
polynomial for $J_1(l)$. The degree of such 
a polynomial would be the cardinality of the group of torsion points, that is
 $l^{2g(X_1(l))}$ where the genus $g(X_1(l))$ of $X_1(l)$
grows quadratically 
in $l$. Such a degree  is far too large for us: we are looking for an
algorithm with  polynomial
time complexity in $l$.

We describe below   an algorithm
 for computing  elements in $V_f$. This is a deterministic  algorithm
and the running time is polynomial  in $l$.
We shall work with
the jacobian $J_1(5l)$ rather than with $J_1(l)$. 
We set $$X=X_1(5l)$$
and we denote by $g$ the genus
of the latter curve. We note that
the conditions above imply $k\ge 4$,  $l\ge 19$ and  
 $g\ge 289$.
Using the map 
$B_{5l,l,1} : X_1(5l) \rightarrow X_1(l)$
defined  in Section~\ref{sec_modforms}, and the 
associated  morphism $B_{5l,l,1}^* : J_1(l) \rightarrow J_1(5l)$
between Jacobian varieties, we can see the
  subspace
$V_f$ as an $\Fl$-subspace inside the $l$-torsion
of the jacobian  $$J=J_1(5l).$$%\index{$J$, the jacobian $J_1(5l)$}
To avoid confusions we shall call 
$W_f\subset J_1(5l)$ the image of $V_f\subset J_1(l)$ by
$B_{5l,l,1}^*$. We call $W_f$ the {\it Ramanujan subspace} associated 
with $f$.
Elements in $J$ are represented by
divisors on the curve $X$. 
For every class in $W_f$, we   compute  
a sharp enough
 approximation  for some  divisor in this class.
It will be explained in Chapter~\ref{sec_comp_mod_l_rep} how one can
compute  the splitting field of this divisor from such a
complex approximation,  using the upper bounds  for the naive height of torsion
 divisors on modular curves proven in Section~\ref{sec_final_estimates}.

\smallskip

This chapter is organized as follows.
In Section~\ref{sec_torsion_points} we recall how points on $X$ are 
represented using standard coordinates
taking values in the complex unit disk. 
Section~\ref{section:latticeperiods} recalls state of the art algorithms for computing
the lattice of periods $\Lambda$ of the jacobian $J$ of the  modular curve $X$. 
In Section~\ref{sec:modfunjmc}
we describe an algebraic model
for  $X_1(5l)$ and we relate it
to the analytic model  $\Gamma_1(5l)\backslash \gHb$.

The four next sections collect useful intermediate results.
Section~\ref{sec_torsion_ps} provides explicit inequalities relating coefficients and
values of converging power series.
In Section~\ref{section:JW} we prove formal identities relating
Jacobians and Wronskians,  that are necessary 
for the local study of the Jacobi integration map.
Section~\ref{section:quantjacob} collects simple quantitative facts
about the Jacobi integration map.
Section~\ref{sec_torsion_QJ}  relates several natural norms on the space
of parabolic modular forms of weight $4$.

A
point on $J$  can be represented
in two different ways. We may  consider it as a class $x+\Lambda$ in $\CC^g/\Lambda=J(\CC)$. We may also fix a degree $g$
divisor $\Omega$ on $X$ and represent an element in $J(\CC)$ by a divisor $Q-\Omega$ in the corresponding linear
equivalence class, where $Q$ is an effective degree $g$ divisor on $X$. 
In Sections~\ref{section:opjac}  and \ref{sec_torsion_arith} we adopt this latter point of view and we show that it is
very convenient for computational purposes and leads to polynomial time algorithms.
Unfortunately, the points  we are interrested in (the $x$ belonging to the Ramanujan subspace $W_f$) are rather
difficult to characterize and compute in this form. However,
assuming the divisor  $\Omega$ has been chosen correctly
(e.g. we take for $\Omega$ the divisor $D_0$ constructed
 in Section~\ref{sec_constr_D}), to every $x$ in $W_f$ there corresponds a unique divisor $Q_x$ such that 
$Q_x-\Omega$ lies in the class represented by $x+\Lambda\in J(\CC)$. Such a $Q_x$ will be  called a {\it Ramanujan divisor}.
Computing $x\in \CC^g$ is not too difficult because 
the defining equations of $V_f$ given in Equation~(\ref{eq:defVltor}) become linear in the analytic model
$\CC^g/\Lambda$. The difficulty then is to compute $Q_x$ once given $x$. This is a typical example of the inverse
Jacobi problem. Section~\ref{sec_torsion_invjac} provides a partial general solution for this inverse Jacobi 
problem: it explains how, given
$\Omega$ and some $x+\Lambda$, one can find a divisor $Q$ such that $Q-\Omega$ lies in the corresponding linear
equivalent  class. Since we are only working with approximations, we must control  the error made in computing
$Q$ from $x$. The output divisor $Q'$ is hopefully close to $Q$ but most likely not equal to it.
We call $x'$ the image of $Q'$ by the Jacobi integration map. 
Statements in Section~\ref{sec_torsion_invjac} control the difference between $x$ and $x'$. We can't  hope a much better result
in full generality
since, in general, the divisor $Q$ is not even unique,
because  $x$ could lye in the singular locus  of the Jacobi map.
To relate the distance between $Q'$ and $Q$ and the distance between $x'$ and $x$,  we need some information
about the local behavior of the Jacobi map at $Q-\Omega$. Using results from Arakelov
theory proven in Section~\ref{sec_final_estimates}, we show in Sections~\ref{section:stabalg},
\ref{sec_torsion_heights} and \ref{section:BXg} that when $x$ lies in $W_f$ then  the error on $x$ and the error on $Q_x$
are nicely related.  This finishes the proof of the main Theorem~\ref{theorem:torsion_main} in this chapter.
The last Section~\ref{section:finalC} provides a more algebraic variant of
this theorem.

In this chapter, we shall use several times the main statements  in Chapter~\ref{sec_couveignes_ZEROS} and in particular Lemma~\ref{lemma:anyzero} 
and Theorem~\ref{theorem:findingzeros}. These statements 
basically say that it is possible to compute efficiently 
sharp approximations of zeros of power series, provided we don't
prospect    near  the boundary of the disk of convergence. In particular 
these zeros are well
conditioned: they are not dramatically affected by a small 
perturbation of the series.

We suggest that the reader look at the first
pages and main statements in Chapter~\ref{sec_couveignes_ZEROS} before going further in this chapter.
%He or she may also  skim through sections \ref{subsection:type},\ref{section:refoc} and \ref{section:boundrem} below.

%An index can be found at the end of this chapter.

\begin{remark} The symbol $\Theta$ in this chapter stands for
a positive effective absolute constant. So any statement containing
this symbol becomes true if the symbol is replaced in every occurrence by some 
large enough real number.
\end{remark}

\begin{remark} In this chapter the letter $i$ stands for
the square root of $-1$ in $\CC$ having positive imaginary part.
\end{remark}

\section{Points, divisors and coordinates on $X$}\label{sec_torsion_points}

In this section we recall how points, functions, forms and divisors
are represented on a modular curve. 
We denote by $\cF$\index{$\cF$, the fundamental domain for the action
of $SL_2(\ZZ)$ on $\gH$}
the classical fundamental domain
for the action of $\SL_2(\ZZ)$ on the Poincar{\'e} upper half plane $\gH$.
We set  $D=\bar \cF \cup \infty$.
We set 

$$T=\left( \begin{array}{cc} 1&1\\0&1  \end{array} \right).$$\index{$T$, the translation $z\mapsto z+1$}

For every positive integer $w$ we set

$$D_w=\bigcup_{0\le k\le w-1}T^k(D).$$\index{$D_w$, the union of the $w$ first translates of $D$}

We denote by $F_w$\index{$F_w\subset D(0,1)$, is mapped by $\mu_\gamma$
onto a $2$-cell on $X$ containing the cusp $\gamma(\infty)$} the image of $D_w$ by the map 
$z\mapsto \exp(2iw^{-1}\pi z)$. This is a compact subset of the open disk
$D(0,1)\in \CC$. It is even contained in $D(0,\exp(-\pi/w))$.
Any cusp on $X(\CC)$ can be written $\gamma(\infty)$
for some   $\gamma \in \SL_2(\ZZ)$. These $\gamma$
can be chosen once for all with entries 

\begin{equation*}%\label{eq:bornegamma}
\le l^\Theta
\end{equation*}
in absolute value. 
We denote by $\Xi$\index{$\Xi\subset \SL_2(\ZZ)$, parametrizes the cusps} the set of all these chosen $\gamma$.
The set $\Xi$ parametrizes the cusps of $X$. We assume  that the identity
belongs to $\Xi$. It parametrizes the cusp $\infty$ itself.
We write the topological space $X(\CC)$  as a union

\begin{equation*}
X(\CC)=\bigcup_{\gamma \in \, \Xi}\gamma (D_{w_\gamma})=
\bigcup_{\gamma \in  \, \Xi}\, \bigcup_{0\le k\le w_\gamma -1}\gamma(T^k(D))
\end{equation*}
where $w_\gamma$ is the width of the
cusp $\gamma(\infty)$.
We say that the $\gamma T^k$ arising in this
union form a  {\it standard system of right cosets representatives}
for  $\Gamma_1(5l)$ in $\SL_2(\ZZ)$.

Every point on $X$ is represented by
a complex number $z$ in  $\gamma(D_{w_\gamma}) \subset \gH$ for
some  $\gamma$ in $\Xi$. 
But $\gamma^{-1}(z)\in D_{w_\gamma}$ will often be  more convenient.
And 

\begin{equation}\label{eq:qgamma}\index{$q_\gamma$, a local parameter at the
cusp $\gamma(\infty)$}
q_\gamma=\exp(2iw_\gamma^{-1}\pi \gamma^{-1}(z))
\end{equation}
 is even 
more convenient.
So most of the time, a complex point on $X$ will be given
as a pair $(\gamma, q)$ where $\gamma \in \Xi$
and $q\in F_{w_\gamma}\subset D(0,\exp(-\pi/w_\gamma))$ 
is the  value of $q_\gamma$ at this point.
We set 

$$D_\gamma = \gamma (D_{w_\gamma} ) \subset \Gamma_1(5l)\backslash \gHb=X_1(5l).$$\index{$D_\gamma$, a $2$-cell on $X$ containing $\gamma(\infty)$}

When $\gamma$ is the identity $\Id$, we sometimes write
$q$
instead of $q_\Id$. The parameter $q$ traditionally   plays a more
important role. Functions and forms on $X_1(5l)$ are often identified
with  their $q$-expansions. 
The function field $\CC(X_1(5l))$ can thus be identified with 
a subfield of the field of Puiseux  series $\CC\{\{q\}\}$.  In
particular, we have an action of $\Aut(\CC)$ on modular functions and
we can define the field
$\QQ(X_1(5l))$ of $\QQ$-rational functions on $X_1(5l)$ to be the field
of functions having $q$-expansion  with $\QQ$-rational  coefficients.

We define the distance $$d_\gamma(P_1,P_2)$$\index{$d_\gamma(P_1,P_2)$ a distance on $D_\gamma\subset X$}
between two points $P_1=(\gamma,q_1)$
and $P_2=(\gamma,q_2)$ in  $D_\gamma$ to be the modulus $|q_2-q_1|$ of the
 difference
of their $q_\gamma$ coordinates. Of course $d_\gamma$ extends to
$D(0,1)\supset F_{w_\gamma}$. 

Given a complex number $q$ in the open disk $D(0,1)$, let 
$z\in \gH \cup \infty$ be such that
$$\exp(2iw_\gamma^{-1}\pi\gamma^{-1}(z))=q.$$ 
In the special case  $q=0$ we set $z = \gamma(\infty)$.
Such a $z$ may not be  unique. But two such $z$'s are mapped
onto each other by some power of $\gamma \times T^{w_\gamma}\times \gamma^{-1}$. Since the latter lies in $\Gamma_1(5l)$ we have
defined a map  $$\mu_\gamma : D(0,1)\rightarrow X_1(5l).$$\index{$\mu_\gamma
: D(0,1)\rightarrow X_1(5l)$, the map associated with  $\gamma$}
This is the  parameterization associated with $\gamma$. It sends 
 $F_{w_\gamma}$ onto $D_\gamma$. Any form (resp. function) on $X_1(5l)$
can be lifted to $D(0,1)$ along the map $\mu_\gamma$. For example,
Klein's modular function $\KJ (z)$\index{$\KJ$, Klein's modular function} is usually given as a function of 
 $q=q_\Id$.  There exists a Laurent series $\bJ (x)$ in
the indeterminate  $x$ such that $\KJ (z)=\bJ(q)$. Further

\begin{equation}\label{eq:744}
\bJ (x)=\frac{1}{x}+744+\sum_{k\ge 1}c(k)x^k
\end{equation}
where the coefficients $c(k)$ are rational integers.
It can be checked easily that the expansion of $\KJ$ at the
cusp $\gamma(\infty)$ is given by 

\begin{equation}\label{eq:KJJ}
\KJ (z)=\bJ(q_\gamma^{w_\gamma})
\end{equation}
\noindent where $w_\gamma$ is the width of the cusp $\gamma(\infty)$. It is
a consequence of a famous theorem by Petersson and Rademacher that 
the coefficient
$c(k)$ is  bounded from above by 
\begin{equation}\label{eq:Petersson}
\Theta^{\sqrt{k}}.
\end{equation}

Now let $d$ be a positive integer and let
 $Q=Q_1+Q_2+\dots+Q_{d}$ be a degree $d$ effective divisor on $X$.
Let $\epsilon$ be a non-negative  real number.
We say that $Q$ is $\epsilon$-{\it simple}\index{$\epsilon$-simple divisor} if the following conditions hold true.

\begin{enumerate}
\item  
For every integer $k$  such that $1\le k\le d$, the point $Q_k$
belongs to $D_{\gamma_k}$ for a unique $\gamma_k$ in $\Xi$. So
  $Q_k=(\gamma_k,q_k)$. We ask  that
$q_k$ lies in the interior of $F_{w_k}$ where
$w_k$ is the width of the cusp $\gamma_k(\infty)$.
\item The distance between $q_k$ and the boundary of $F_{w_k}$ is 
$>\epsilon$.
\item If $1\le k_1 < k_2 \le d$ and $\gamma_{k_1}=\gamma_{k_2}=\gamma$, we write $Q_{k_1}=(\gamma,q_{1})$ and $Q_{k_2}=
(\gamma,q_2)$ and we ask that  $|q_2-q_1|>\epsilon$.
\end{enumerate}

Not every divisor $Q$  is $\epsilon$-simple but if $\epsilon <1/(d\Theta)$ there 
exists an $\epsilon$-simple
 divisor $Q'=Q'_1+Q'_2+\dots+Q'_{d}$
such that for every $1\le k\le d$
we have $Q'_k=(\gamma_k,q'_k)$
and $|q'_k-q_k|\le \Theta d\epsilon$.

\section{The lattice of periods}\label{section:latticeperiods}

This section is devoted to the explicit calculation of the lattice
of periods of $X$. All the algorithms in this section are detailed
in the two books by Cremona \cite{cremona} and Stein \cite{Stein} and
in Bosman's thesis \cite{Bosman3}. See also Chapter~\ref{chapcomput}.

We  need a 
 complex analytic description of the torus
$J(\CC)$ as $\CC^g/\Lambda$ where $\Lambda$ is the lattice of periods.
We  first compute an explicit description of the first
group in singular cohomology $$H_1^{\rm sing}(X_1(5l),\ZZ).$$
Using Manin-Shokurov theory we find a basis ${\cB_1^{\rm sing}}$\index{$\cB_1^{\rm sing}$, a basis for the singular
homology}
 for this
$\ZZ$-module. 
Every element in this basis is an integer  linear combination
of  Manin symbols 

\begin{equation}\label{eq:sumManin}
\sum_{\gamma} c_\gamma \{\gamma(0),\gamma(\infty)\}.
\end{equation}

The $\gamma$ in the sum (\ref{eq:sumManin}) runs over
the  standard system of right cosets representatives
for  $\Gamma_1(5l)$ in $\SL_2(\ZZ)$. 
The integer coefficients $c_\gamma$ can be chosen
to be 

\begin{equation*}%\label{eq:bornecgamma}
\le \exp(l^\Theta)
\end{equation*}
in absolute value. 

We also need a basis $\cB^1_{\rm DR}$\index{$\cB^1_{\rm DR}$, 
a basis for $\cH^1=H^1_{\rm DR}(X_1(5l))$}
of the space of holomorphic differentials
$$\cH^1=H^1_{\rm DR}(X_1(5l))$$
or equivalently a basis of the space 
$S_2(\Gamma_1(5l))$ of weight
two  cusp forms.   We shall use the standard basis made
of normalized newforms of level $5l$ together
with normalized newforms of level $l$ lifted
to level $5l$ by the two degeneracy maps.

Let $f=\sum_{k\ge v}f_kq^k$ be a form in this basis.
The $q$-valuation $v$ of $f$ is $1$ or $5$.
The first non-zero coefficient  $f_v$
 in the $q$-expansion of $f$
is $1$.
The  coefficients $f_k$ 
in the $q$-expansion of $f$ are algebraic integers. The modulus
 of $f_k$  is $k^\Theta$. One can compute an approximation
of $f_k$ within $\exp(-m)$ in deterministic polynomial time $(klm)^\Theta$.

The action of  Atkin-Lehner involutions $w_5$ and $w_l$ is expressed in the
 basis $\cB^1_{\rm DR}$ by 
theorem 2 of \cite{Asai1}. The action of the diamond operators is known also 
because  every 
form in $\cB^1_{\rm DR}$ is an eigenform for the Hecke algebra $\TT^{(5l)}$
generated by the operators $T_n$ for $n$ prime to $5l$.

We  also need the expansion of every form $f(q)$
in $\cB_{\rm DR}^1$
at every cusp $\gamma(\infty)$.
More precisely, $f(q)q^{-1}dq$ should be rewritten as 
$h(q_\gamma)q_\gamma^{-1}dq_\gamma$ for 
every $\gamma$ in $\Xi$.
Since the level
$5l$ is squarefree, the group generated by the Atkin-Lehner 
involutions and the diamond operators
acts simply transitively on the cusps. So there is an automorphism
in this group that sends $\infty$ to $\gamma(\infty)$.
This automorphism can be represented by a matrix $W_\gamma$ in $\GL_2(\QQ)$
having integer entries
as explained in  Section~\ref{atleop} of Chapter~\ref{chapcomput}.
 If
$$\gamma=\left( \begin{array}{cc} a&b\\c&d  \end{array} \right)$$
then the 
width $w_\gamma$ of the cusp $\gamma(\infty)$ is 
$$w_\gamma=\frac{5l}{\gcd(5l,c)}.$$
Let $r$ be the unique integer in $[0,w_\gamma[$
such that $d\equiv cr \bmod w_\gamma$. Set $b'=b-ar$
and $d'=(d-cr)/w_\gamma$ and $c'=c/\gcd(5l,c)$. Then
\begin{equation}\label{eq:Wgamma}
W_\gamma=\left( \begin{array}{cc} aw_\gamma&b'\\5lc'&w_\gamma d'  \end{array} \right)
\end{equation}
 and
the  product 
\begin{equation}\label{eq:Wgamma2}
W^{-1}_\gamma \times 
 \gamma= \left( \begin{array}{cc} w_\gamma^{-1}&rw_\gamma^{-1}\\0&1 \end{array} \right)
\end{equation}
 fixes $\infty$ and it acts on Fourier
expansions like the substitution $q\mapsto \zeta q^{1/w_\gamma} $ for some
root of unity $\zeta=\exp(\frac{2ri\pi}{w_\gamma})$. 
Since the action of $W_\gamma$ on forms
is known, we can
compute the expansion of all forms in $\cB_{\rm DR}^1$ at all
cusps in deterministic polynomial time $(klm)^\Theta$ where $k$
is the $q$-adic accuracy and $m$ the complex absolute
accuracy of coefficients.

Once we have computed a basis for both the singular homology $H_1^{\rm sing}(X_1(5l),\ZZ)$
and the de~Rham cohomology $H^1_{\rm DR}(X_1(5l))$ we can compute the lattice $\Lambda$\index{$\Lambda$, the lattice of periods}
of periods. Since we are given a basis $\cB^1_{\rm DR}$ of holomorphic differentials, the lattice
$\Lambda$ is well defined inside $\CC^{\cB^1_{\rm DR}}$ as the image of the integration
map $H_1^{\rm sing}(X_1(5l),\ZZ)\rightarrow \CC^{\cB^1_{\rm DR}}$ sending a cycle $c$ onto the vector
$(\int_{c}\omega)_{\omega\in \cB^1_{\rm DR}}$. The image of the basis $\cB_1^{\rm sing}$ by the integration map
is a basis $\cB_{\rm per}$\index{ $\cB_{\rm per}$, a basis of the lattice of periods} of the lattice $\Lambda$ of periods.
The so-called matrix of periods has entries $\int_{c}\omega$ where $c$ is a cycle in the basis
${\cB_1^{\rm sing}}$ of  $H_1^{\rm sing}(X_1(5l),\ZZ)$
and $\omega=f(q)q^{-1}dq$ is a form in the basis $\cB^1_{\rm DR}$ of $H^1_{\rm DR}(X_1(5l))$. Computing these 
periods reduces to evaluating integrals of the form $\int_{\alpha}^\beta f(q)q^{-1}dq$ where
$\alpha$ and $\beta$ are two  cusps. We first cut this integral in two pieces
$\int_{\alpha}^\beta f(q)q^{-1}dq=\int_{i}^\beta f(q)q^{-1}dq - \int_{i}^\alpha f(q)q^{-1}dq$.
Since the group
generated by Atkin-Lehner  involutions  and diamond operators
acts transitively on the cusps, we can reduce to the computation of integrals
of the form
$\int_{\alpha}^\infty f(q)q^{-1}dq$ where $\alpha=(a+bi)/c$ 
and $a$, $b$  and $c$ are integers bounded by $l^\Theta$ in absolute value.
Since the coefficient $f_k$  in the $q$-expansion of $f$ are bounded by $k^\Theta$ we
can compute  approximations of the entries  in the period matrix within $\exp(-m)$
in deterministic polynomial time $(lm)^\Theta$.

We note that the $L^\infty$ norm on $\CC^{\cB^1_{\rm DR}}$ induces a distance 
$d_J$\index{$d_J$,  a distance on the jacobian $J$} on
the quotient $\CC^{\cB^1_{\rm DR}}/\Lambda=J(\CC)$

\begin{equation}\label{eq:distJ}
d_J(x+\Lambda,y+\Lambda)=\min_{z\in \Lambda} |x-y-z|_\infty.
\end{equation}

 This distance will be useful when evaluating rounding
errors in the course of numerical computations.
We denote by 
$$\phi : X\rightarrow J$$
the Jacobi integration map. This map is well defined once we have
chosen a degree $1$ divisor on $X$ as origin. For any $\gamma$
in $\Xi$, the restriction of $\phi$ to
$D_\gamma$ 
is Lipschitz with constant $l^\Theta$ according
to Equation~(\ref{eq:lipjac}). More precisely, if $P_1=(\gamma,q_1)$ and 
$P_2=(\gamma,q_2)$
are two points in $D_\gamma$ then

\begin{equation*}%\label{eq:lipphi}
d_J(\phi(P_2),\phi(P_1))\le l^\Theta \times |q_2-q_1|. 
\end{equation*}

For every positive
integer $k$, 
we also denote by $\phi$ the integration map $\phi : X^k\rightarrow J$\index{$\phi : X^k\rightarrow J$, the Jacobi integration map}. 
We denote
by 
$$\phi' : \Div(X)\rightarrow J$$
the map induced by $\phi$ on the group of divisors
on  $X$. The restriction of $\phi'$ to the subgroup
$\Div^0(X)$ of degree $0$ divisors is independent of
the origin we have chosen.

\section{Modular functions}\label{sec:modfunjmc}

Since we plan to compute the splitting field of some very special
divisors on the modular curve $X=X_1(5l)$ we must be able to evaluate
some well chosen  modular functions of weight $0$ and
level $5$ or $5l$ at a given
point $z\in \HH$. In this section we describe algebraic models
for $X_1(5)$ and $X_1(5l)$ and we explain how to compute the expansions
of the involved modular functions at every cusp.

\subsection{The modular curve $X_1(5)$}

In this section we recall
the definition of several 
classical level $5$ modular functions and we show how
to compute their expansions
at each of the four cusps of $X_1(5)$.
Let $b$ be an indeterminate and consider
 the  elliptic curve $E_b$ in Tate
normal form 
with equation 
\begin{equation}\label{eq:tatecurveb}
y^2+(1+b)xy+by=x^3+bx^2.
\end{equation} The point $P=(0,0)$ has order
$5$ and its multiples are  $2P=(-b,b^2)$, $3P=(-b,0)$, $4P=(0,-b)$.
Call $\PP^1_b$    the projective line
with parameter $b$. The modular invariant of $E_b$ is 
\begin{equation}\label{eq:j(b)}
j=j(b)=-\frac{(b^4+12b^3+14b^2-12b+1)^3}{b^5(b^2+11b-1)}
\end{equation}
Let   $$s=-\frac{11+5\sqrt{5}}{2}$$ and $\bar s$ 
be the two complex  roots of
$b^2+11b-1$.
We call $A_\infty$, $A_0$, $A_s$, $A_{\bar s}$ the points on $\PP^1_b$
corresponding to the values $\infty$, $0$, $s$ and $\bar s$ of the parameter
$b$.
The elliptic curve $E_b\rightarrow \PP^1_b-\{A_\infty, A_0, A_s, A_{\bar s}\}$
is the universal elliptic curve with one point of order $5$. So there exists a unique isomorphism 
between
the modular curve $X_1(5)=
\Gamma_1(5l)\backslash \gHb$ and $\PP^1_b$
 that is compatible with the 
moduli structure on either side. We want to compute this
isomorphism.  More precisely we
 compute the expansions of $b$ at every cusp of $X_1(5)$.
To this end we compare the  curve in Equation~(\ref{eq:tatecurveb})   
and the Tate curve 
\cite{tateanvers4} with equation
\begin{equation}\label{eq:tatecurveq}
y'^2+x'y'=x'^3 +a_4x'+a_6
\end{equation}
where 
\begin{eqnarray*}
a_4&=&-5\sum_{n\ge 1} \sigma_3(n)\,  q^n\\
a_6&=&-\sum_{n\ge 1} \frac{5\sigma_3(n)+7\sigma_5(n)}{12}\, q^n
\end{eqnarray*}
and $q$ is a formal parameter. We note that the coefficients
in the expansions above are integers and we have
\begin{eqnarray*}
a_4&=&\frac{1-E_4(q)}{48}\\
a_6&=&\frac{1-3E_4(q)+2E_6(q)}{1728}
\end{eqnarray*}
where
\begin{eqnarray*}
E_4(q)&=&1+240 \sum_{n\ge 1}\sigma_3(n)\, q^n\\
E_6(q)&=&1-504\sum_{n\ge 1}\sigma_5(n)\, q^n.
\end{eqnarray*}
The modular invariant 
of the Tate curve in Equation~(\ref{eq:tatecurveq}) is
\begin{equation}\label{eq:j(q)}
j(q)=\frac{1728E_4^3(q)}{E_4^3(q)-E_6^2(q)} = \frac{1}{q}+744+196884q+\cdots
\end{equation}
Any isomorphism between the two Weierstrass curves
in Equations~(\ref{eq:tatecurveb}) and (\ref{eq:tatecurveq})
must take the form
\begin{eqnarray}\label{eq:xx'}
 x&=&u^2x'+r\\
\nonumber  y&=&u^3y'+su^2x'+t.
\end{eqnarray}

Straightforward calculation gives the following
necessary and sufficient conditions
for the affine transform in (\ref{eq:xx'}) to induce an isomorphism
of Weierstrass curves:

\begin{eqnarray}\label{eq:isombq}
\frac{1728E_4^3(q)}{E_4^3(q)-E_6^2(q)}&=&-\frac{(b^4+12b^3+14b^2-12b+1)^3}{b^5(b^2+11b-1)}\\
\nonumber u^2&=&-\frac{E_4(q)}{E_6(q)}\times \frac{(b^2+1)(b^4+18b^3+74b^2-18b+1)}{(b^4+12b^3+14b^2-12b+1)}\\
\nonumber  r&=&\frac{u^2-b^2-6b-1}{12}  \\
\nonumber s&=& \frac{u-b-1}{2} \\
\nonumber t&=&  \frac{b^3+7b^2-(5+u^2)b+1-u^2}{24}
\end{eqnarray}

We can  simplify a bit these expressions. The
first one just means
$$j(q)=j(b).$$ 
From  
the classical  \cite[Proposition 7.1]{Schoof2} identities 
$$\left( \frac{qdj}{dq}  \right)^2=j(j-1728)E_4$$
\noindent 
and 
$$\left(
\frac{qdj}{dq}  \right)^3=-j^2(j-1728)E_6$$
\noindent 
 we deduce  
$$\left(
\frac{qdb}{dq}  \right)=
\frac{E_6(q)}{E_4(q)}\times \frac{b(b^2+11b-1)(b^4+12b^3+14b^2-12b+1)}{5(b^2+1)(b^4+18b^3+74b^2-18b+1)}.$$
So the expression for $u^2$ can be written
\begin{equation}\label{eq:u2b}
u^{2}= -\frac{b(b^2+11b-1)}{5\left(
\frac{qdb}{dq}  \right)}
\end{equation}
The expansion of $j$ as a series in $q$  has integer coefficients and can be 
computed using  Equation~(\ref{eq:j(q)}).
We deduce the expansion of $q$ as a series in $j^{-1}$
\begin{equation}\label{eq:q(j)}
q=j^{-1}+744j^{-2}+750420j^{-3}+\cdots
\end{equation}
It 
has integer coefficients
and one can compute it from the expansion of $j$ as a series in $q$
using any reasonable algorithm for the reversion of a power series: brute force linear algebra or the more efficient algorithms in 
\cite{BK} or the quasi-optimal algorithms in \cite{UK}.

{\it We first study the situation locally at $A_\infty$.}
A local parameter for $\PP^1_b$ at $A_\infty$
is $b^{-1}$.
 The expansion 
\begin{equation}\label{eq:jmu(b)inf}
j^{-1}=-b^{-5}+25b^{-6}+\cdots
\end{equation}
of $j^{-1}$  as a series 
in $b^{-1}$ has integer coefficients and can be computed 
using Equation~(\ref{eq:j(b)}) and standard
algorithms for polynomial arithmetic.
We substitute (\ref{eq:jmu(b)inf}) in (\ref{eq:q(j)}) and find
\begin{equation*}%\label{eq:q(b)inf}
q= -b^{-5}+25b^{-6}+\cdots
\end{equation*}
One more reversion gives the expansion of $b^{-1}$ as a series in $q^{\frac{1}{5}}$.
\begin{equation}\label{eq:bqinf}
b^{-1}=-q^\frac{1}{5} + 5q^\frac{2}{5}+\cdots.
\end{equation}
This expansion  defines  an embedding of the local field at 
$A_\infty$ inside the field of Puiseux series
$\CC\{\{q\}\}$ in the formal parameter $q$.
We deduce the expansion of $u$ from Equations~(\ref{eq:u2b}) 
and (\ref{eq:bqinf}) 
\begin{equation}\label{eq:uqinf}
u^2=q^{\frac{-2}{5}}+4q^{\frac{-1}{5}}+4-10q^{\frac{1}{5}}-
30q^{\frac{2}{5}}+\cdots 
\end{equation}
We also find
\begin{equation}\label{eq:rqinf}
r=-1-5q^{\frac{1}{5}}-10q^{\frac{2}{5}}+35q^{\frac{4}{5}}+
45q^{}+\cdots 
\end{equation}
The coordinates of the $5$-torsion $P$ on $E_b$
are $x_P=0$ and $y_P=0$.
The image of $P$ on the Tate curve has $x'$-coordinate
$$x'_P=-ru^{-2}=q^{\frac{2}{5}}+q^{\frac{3}{5}}+2q^{\frac{4}{5}}-
2q+\cdots $$
Since on  the Tate  curve we have
\begin{equation}\label{equation:expa}
x'(w, q)=\sum_{n\in \ZZ}\frac{wq^n}{(1-wq^n)^2}
-2\sum_{n\ge 1}\frac{nq^n}{1-q^n}
\end{equation}
\noindent we deduce that  the value of the
parameter $w$ at the $5$-torsion point $P$
is $$w(P)=q^{\pm \frac{2}{5}} \bmod <q>.$$ We may
take either  sign in the exponent above because we may choose  any of  the two
isomorphisms corresponding to either possible values for $u$. We
decide that 
\begin{equation}\label{eq:wPinfty}
w(P)=q^{\frac{2}{5}} \bmod <q>.
\end{equation}
So the limit curve when $q\rightarrow 0$ or equivalently when
$b\rightarrow \infty$ is a $5$-gon,  with a $5$-torsion
point lying  on the second component
after the one carrying the origin.
So let $$\gamma = \left(\begin{matrix}a & b\\ c & d  \end{matrix}\right)$$
be any matrix in $\SL_2(\ZZ)$ such that $c$ is $2$ modulo~$5$ and
$d$ is $0$ modulo~$5$. For example take
$$\gamma = \left(\begin{matrix}1 & 2\\ 2 & 5  \end{matrix}\right)$$

In particular, the cusp  $\gamma(\infty)=1/2$ has width $5$.
We set
$$q_\gamma(z)=\exp(2i\pi \gamma^{-1}(z)/5).$$
 If we replace $q^{1/5}$
by $q_\gamma$ in Equations~(\ref{eq:bqinf}),  (\ref{eq:uqinf})
and (\ref{eq:rqinf}) 
we  obtain the 
expansions of $b^{-1}$,  $u^2$ and $r$ 
at the cusp
$\gamma(\infty)=1/2$. 
Note in particular that $u^2$ is a modular function 
with weight $-2$ and level $5$. By construction, it has no zero and no
pole outside the cusps.
Similarly  $b$ is a modular function 
with weight $0$ and level $5$. By construction, it has no 
pole outside the cusps.

{\it We now study
the  situation locally at  $A_0$.} A local parameter at
$A_0$
is $b$. We find
 $j^{-1}=b^5+25b^6+\cdots$ and  $q=b^5+25b^6+\cdots$ and we
fix an embedding of the local field at $A_0$ inside $\CC\{\{q\}\}$ by setting
\begin{equation}\label{eq:bqzero}
b=q^{\frac{1}{5}}-5q^{\frac{2}{5}}+15q^{\frac{3}{5}}-30q^{\frac{4}{5}}+\cdots
\end{equation}
Using Equation~(\ref{eq:u2b})
we deduce 
\begin{equation}\label{eq:uqzero}
u^2=1-6q^{\frac{1}{5}}+19q^{\frac{2}{5}}-40q^{\frac{3}{5}}+55q^{\frac{4}{5}}+\cdots
\end{equation}
and
\begin{equation}\label{eq:rqzero}
r=-q^{\frac{1}{5}}+4q^{\frac{2}{5}}-10q^{\frac{3}{5}}+15q^{\frac{4}{5}}+\cdots
\end{equation}
So  the coordinate $x'_P$ of the $5$-torsion point $P$ is
$$x'_P=q^{\frac{1}{5}}+2q^{\frac{2}{5}}+3q^{\frac{3}{5}}+5q^{\frac{4}{5}}+
3q+O(q^{\frac{6}{5}})$$ and  the
parameter $w$ at $P$ can be taken to be  $w(P)=q^{\frac{1}{5}}\bmod <q>$ this time.
The limit curve when $q\rightarrow 0$ or equivalently when
$b\rightarrow 0$ is a $5$-gon,  with a $5$-torsion
point lying  on the first  component
after the one carrying the origin.
So let $$\gamma = \left(\begin{matrix}a & b\\ c & d  \end{matrix}\right)$$
be any matrix in $\SL_2(\ZZ)$ such that $c$ is $1$ modulo~$5$ and
$d$ is $0$ modulo~$5$. For example take
$$\gamma = \left(\begin{matrix}0 & -1\\ 1 & 0  \end{matrix}\right)$$
In particular, the cusp  $\gamma(\infty)=0$ has width $5$.
We set
$$q_\gamma(z)=\exp(2i\pi \gamma^{-1}(z)/5).$$
 If we replace $q^{1/5}$
by $q_\gamma$ in Equations~(\ref{eq:bqzero}),  (\ref{eq:uqzero})
and (\ref{eq:rqzero})
we  obtain the 
expansions of $b$, $u^2$ and $r$ 
at the cusp
$\gamma(\infty)=0$. 

{\it We now study the situation locally at $A_s$.}
A local parameter at
$A_s$
is $b-s$ and 
$$j^{-1}=(-\frac{1}{2}+\frac{11\sqrt{5}}{50})(b-s)+
(-\frac{45443}{125}+\frac{4064\sqrt{5}}{25}  )(b-s)^2+\cdots $$ 
and 
$$b-s=-\frac{125+55\sqrt{5}}{2}q
-(375+170\sqrt{5})q^2
-\frac{3375+1495\sqrt{5}}{2}q^3+\cdots $$
and
$$u^2=-\frac{25+11\sqrt 5}{2}-(200+90\sqrt 5)q-\frac{3575+1595\sqrt 5}{2}q^2+
\cdots$$
and
$$r=-\frac{7+3\sqrt 5}{2}-(100+45\sqrt 5)q-(1300+580\sqrt 5)q^2+
\cdots$$
and
$$x'_P=-\frac{1}{2}+\frac{\sqrt 5}{10}-\frac{5+\sqrt 5}{2}q+
\frac{-15+\sqrt 5}{2}q^2 +(-10+\sqrt 5)q^3+\cdots$$
We note  that the coordinate $x'(P)$ of the $5$-torsion point $P$ is
$$x'(P)=\frac{w}{(1-w)^2}+O(q)$$ where $w=\exp(\frac{4i\pi}{5})=\zeta_5^2$.
So the
parameter $w$ at $P$ can be taken to be  $w(P)= \zeta_5^{2}\, \bmod <q>$ this time.
The limit curve when $b\rightarrow s$ is thus
a $1$-gon equipped with
the $5$-torsion point $\zeta_5^2$ in its smooth locus $\Gm$.
Let $$\gamma = \left(\begin{matrix}a & b\\ c & d  \end{matrix}\right)$$
be any matrix in $\SL_2(\ZZ)$ such that $c$ is $0$ modulo~$5$ and
$d$ is $2$ modulo~$5$. For example take
$$\gamma = \left(\begin{matrix}3 & 1\\ 5 & 2  \end{matrix}\right)$$
In particular, the cusp  $\gamma(\infty)=3/5$ has width $1$.
We set
$$q_\gamma(z)=\exp(2i\pi \gamma^{-1}(z)).$$
 If we replace $q$
by $q_\gamma$ in the expansions above
we  obtain the 
expansions of $b$, $u^2$ and $r$ 
at the cusp
$\gamma(\infty)=0$.

{\it We finally  study the situation locally at $A_{\bar s}$.}
A local parameter at
$A_s$
is $b-\bar s$ and 
$$j^{-1}=(-\frac{1}{2}-\frac{11\sqrt{5}}{50})(b-\bar s)-
(\frac{45443}{125}+\frac{4064\sqrt{5}}{25}  )(b-\bar s)^2+\cdots$$ 
and 
$$b-\bar s=\frac{-125+55\sqrt{5}}{2}q
+(-375+170\sqrt{5})q^2
+\frac{-3375+1495\sqrt{5}}{2}q^3+\cdots$$
and
$$u^2=\frac{-25+11\sqrt 5}{2}+(-200+90\sqrt 5)q+\frac{-3575+1595\sqrt 5}{2}q^2+
\cdots$$
and
$$r=\frac{-7+3\sqrt 5}{2}+(-100+45\sqrt 5)q+(-1300+580\sqrt 5)q^2+
\cdots$$
and
$$x'_P=-\frac{1}{2}-\frac{\sqrt 5}{10}+\frac{-5+\sqrt 5}{2}q-
\frac{15+\sqrt 5}{2}q^2 -(10+\sqrt 5)q^3+\cdots$$
We note  that the coordinate $x'(P)$ of the $5$-torsion point $P$ is
$$x'(P)=\frac{w}{(1-w)^2}+O(q)$$ where $w=\exp(\frac{2i\pi}{5})=\zeta_5$.
So the
parameter $w$ at $P$ can be taken to be  $w(P)= \zeta_5 \, \bmod <q>$ this time.
The limit curve when $b\rightarrow \bar s$ is thus
a $1$-gon equipped with
the $5$-torsion point $\zeta_5$ in its smooth locus $\Gm$.
Let $$\gamma = \left(\begin{matrix}a & b\\ c & d  \end{matrix}\right)$$
be any matrix in $\SL_2(\ZZ)$ such that $c$ is $0$ modulo~$5$ and
$d$ is $1$ modulo~$5$. For example take
$$\gamma = \Id =  \left(\begin{matrix}1 & 0\\ 0 & 1  \end{matrix}\right)$$

In particular, the cusp  $\gamma(\infty)=\infty$ has width $1$.
We set
$$q_\Id(z)=\exp(2i\pi z).$$
 If we replace $q$
by $q_\Id$ in the expansions above
we  obtain the 
expansions of $b$, $u^2$ and $r$ 
at the cusp
$\infty$.

Altogether we have proved the following.

\begin{lem}[Computing expansions of $b$, $u^2$ and $r$]\label{lemma:bu2r}
There exists a deterministic algorithm that given an integer
$k\ge 1$ computes the $k$
first terms in the expansions of $b$, $u^2$ and $r$ at
each of the four cusps of $X_1(5)$,  at the expense of $k^\Theta$
elementary operations.
\end{lem}

We recall that the functions
$b$, $u^2$ and $r$ are defined in Equation~(\ref{eq:isombq}).
Further $b$ is a weight zero and level $5$
modular function having no pole
outside the cusps. 
The weight $-2$ and level $5$ modular function
$u^2$  has no zero and no pole outside the cusps. 
Finally $r$ is the sum of two level $5$ modular functions having
no pole outside the cusps. One of
weight $-2$ and one  of weight $0$. 
We also notice that  the coefficients in the expansions
of $b$, $u^2$ and $r$ lye in $\ZZ$ or $\ZZ[\frac{1+\sqrt 5}{2}]$.
Lemma~\ref{lem:PRM} implies that 
there exists a  positive constant $\Theta$   such that
for every integer $k\ge 1$,  the $k$-th coefficient in these
expansions has absolute value $\le \exp(\Theta \sqrt k)$.

\subsection{A plane model for $X_1(5l)$}\label{subsection:descripCl}

Let $l > 5$ be a prime.  We call   $X_l=X_1(5l)$ the
moduli of elliptic curves with 
one  point of order $5l$. The
genus of $X_l$ is $g_l=(l-2)^2$.
In this section we 
define and study
 a natural homogeneous singular plane  model $C_l$ for
this curve. 
In particular we    enumerate the geometric points on  $X_l$ above every 
singularity of $C_l$ and we  explain how to compute
 series expansions for affine
coordinates at every such branch.
Finally we 
recall how to compute the zeta function of the function
field $\Fp(X_l)$  for every prime integer
$p\not\in\{5,l\}$.

Let $b$ be an indeterminate and let
$E_b$ be the elliptic curve  in  Equation~(\ref{eq:tatecurveb}).
The field $\QQ(b)$ of rational fractions in $b$ is the function field
$\QQ(X_1(5))$ of the modular curve $X_1(5)$ over $\QQ$. The map
$$B_{5l,5,1} :X_l\rightarrow X_1(5)$$ introduced in Section~\ref{sec_modforms}  defines a degree $l^2-1$  extension $\QQ(X_l)/\QQ(b)$. We  construct 
an explicit
model for this extension.
The multiplication by $l$
isogeny $$[l] : E_b \rightarrow E_b$$
induces a degree $l^2$ rational function on $x$-coordinates:
$$x\mapsto \frac{N(x)}{D(x)}$$ where $N(x)$ is a monic degree $l^2$
polynomial in $\QQ(b)[x]$. Recursion formulae for division polynomials
(see \cite{enge} section 3.6) provide a quick algorithm for computing this
polynomial, and also show that the coefficients actually lie in
$\ZZ[b]$. The roots of $N(x)$ are the $x$-coordinates of the
points $Q$ on $E_b$
such that $[l]Q$ is $P=(0,0)$. If $l$ is congruent to $\pm 1$ modulo $5$ then $l P=\pm
P$ and $x$ divides $N(x)$. Otherwise $N(x)$ is divisible by $x-b$.
 Call $T_l(b,x)$ the
quotient of $N(x)$ by $x$ or $x-b$, accordingly.
 This  is 
 a monic  polynomial in $\ZZ[b][x]$ with 
degree
$l^2-1$ in $x$. As a polynomial in $x$ we have

$$T_l(b,x)=\sum_{0\le k \le
 l^2-1}a_{l^2-1-k}(b)x^k$$
\noindent where  $a_0(b)=1$.  We call
 $d$ be the total degree
of $T_l$.

Let $F$ be a field  extension of
$\QQ(b)$ where 
$T_l(b,x)\in \QQ(b)[x]$ has a root.
A suitable twist of the Tate  curve $E_b$
given by Equation~(\ref{eq:tatecurveb})
 has
a point of order $5l$ defined over $F$.
This proves that the function
field  extension  $\QQ(X_l)/\QQ(b)$
can be embedded in  $F/\QQ(b)$. 
Since the degree of
$T_l(b,x)$ in $x$ is equal to the degree of the
extension 
$\QQ(X_l)/\QQ(b)$ we deduce that
the polynomial
$T_l$ is irreducible in $\QQ (b)[x]$ and the quotient
field $\QQ(b)[x]/T_l$ is isomorphic to $\QQ(X_l)$. Since the latter
field is a regular extension of $\QQ(b)=\QQ(X_1(5))$ we deduce
that $T_l$ is absolutely  irreducible.

Let  $C_l\subset \PP^2$ be the  projective curve
 with homogeneous  equation
$T_l(\frac{\bbb}{\bz},\frac{\bx}{\bz})\bz^{d}$ in the variables
$\bbb$, $\bx$ and $\bz$.
The  map 
$B_{5l,5,1} : X_1(5l)\rightarrow X_1(5)$
is unramified
except at $b\in \{0,\infty,s,\bar s\}$.
So for every point $R$ on $X_l$ such that $b(R) \not \in
\{0,s,\bar s,\infty \}$,  the function $b -b (R)$ is a
uniformizing parameter at $R$. 
Let $\cU$ be the affine open set with equation $$\bz\bbb(\bbb^2-11\bbb\bz+\bz^2)\not =
0.$$ Every point  on  $C_l\cap \cU$  is smooth and all points on
$X_l$ above  points in
$C_l - \cU$ are cusps in the modular sense (i.e. the modular
invariant has a pole at these points). The smooth point $R=(b_R,x_R)$ 
on $C_l\cap \cU$ is the moduli of the curve $E_{b_R}$ equipped 
with the  unique $5l$-torsion point $Q$ having $x$-coordinate $x_R$
and such that $[l]Q=P\in E_{b_R}$.

\subsection{The singularities of $C_l$}\label{subsection:singCl}

We  study
the cusps of $X_1(5l)$
that are mapped onto  $A_\infty$ by $B_{5l,5,1}$.
Set  $\zeta_l=\exp(\frac{2i\pi}{l})$.
Let  $\alpha$ and $\beta$ be
integers such that $0\le \alpha,\beta\le l-1$. 
Let also  $\talpha$ and $\tbeta$ be
integers such that $0\le \talpha,\tbeta\le l-1$ and 
$$\tilde \alpha \equiv
\alpha /5 \bmod l$$ and $$\tilde \beta \equiv (\beta-2)/5 \bmod l.$$
We set
$w_Q=\zeta_l^\talpha q^{\frac{\tbeta}{l}}q^{\frac{2}{5l}}$ 
and observe that $$w_Q^5=\zeta_l^\alpha q^{\frac{\beta}{l}}\bmod <q>$$
and $$w_Q^l=q^{\frac{2}{5}}\bmod <q>=w_P$$
according to Equation~(\ref{eq:wPinfty}).
We denote by $Q$ the point on the Tate curve 
in   (\ref{eq:tatecurveq}), having 
 $w$-coordinate $w_Q$. The isomorphism given by
Equations~(\ref{eq:xx'}) maps $Q$ onto a $5l$-torsion point
on the curve $E_b$. This point is called $Q$ also.
The couple  $(E_b,Q)$ defines a point on $X_1(5l)$ that
is mapped onto $(E_b,P)$ by $B_{5l,5,1}$.
We substitute $w$ by $w_Q$ in 
expansion~(\ref{equation:expa})
and find
$$x'_{Q}=\zeta_l^\talpha q^{\frac{\tbeta}{l}}q^{\frac{2}{5l}}(1+O(q^{\frac{1}{5l}}))$$
\noindent if $0\le \tbeta \le \frac{l -1}{2}$ and
$$x'_{Q}=\zeta_l^{-\talpha}q^{\frac{l-\tbeta}{l}-\frac{2}{5l}}(1+O(q^{\frac{1}{5l}}))$$
\noindent if $\frac{l +1}{2}\le \tbeta \le l-1$.
Using Equation~(\ref{eq:xx'}) and the expansions in
Equations~(\ref{eq:uqinf}) and (\ref{eq:rqinf}) we 
find 
\begin{equation}\label{eq:polesx1}
x_{Q}+1=\zeta_l^\talpha q^{\frac{\tbeta}{l}+\frac{2}{5l}-\frac{2}{5}}(1+O(q^{\frac{1}{5l}}))
\end{equation}
\noindent if $0\le \tbeta \le \frac{l -1}{2}$  and

\begin{equation}\label{eq:polesx2}
x_{Q}+1=\zeta_l^{-\talpha} q^{\frac{l-\tbeta}{l}-\frac{2}{5l}-\frac{2}{5}}(1+O(q^{\frac{1}{5l}}))
\end{equation}
\noindent if $\frac{l +1}{2}\le \tbeta \le l-1$.

Let $$\gamma = \left(\begin{matrix}a & b\\ c & d  \end{matrix}\right)$$
be any matrix in $\SL_2(\ZZ)$ such that $c$ is $5\tbeta +2$ modulo~$5l$ and
$d$ is $5\talpha$ modulo~$5l$. 
In particular, the cusp  $\gamma(\infty)=a/c$ has width $w=\frac{5l}{\gcd(5l,c)}$. So $w$ is  $5$ if
$\beta=0$ and $5l$ otherwise.
We set
$$q_\gamma(z)=\exp(2i\pi \gamma^{-1}(z)/w).$$
 If we replace $q^{1/w}$
by $q_\gamma$ in Equation~(\ref{eq:polesx1})
or (\ref{eq:polesx2})  above
we  obtain the 
expansion of $x$ 
at the cusp
$\gamma(\infty)=a/c$.

The same method applies to cusps of $X_1(5l)$ that are mapped
onto $A_0$, $A_s$ or $A_{\bar s}$ by $B_{5l,5,1}$.

\begin{lem}[Computing expansions of $x$]\label{lemma:expx}
There exists a deterministic algorithm that given an integer
$k\ge 1$ computes the $k$
first terms in the expansions of $x$  at 
each of the  cusps of $X_1(5l)$, at the expense of $(kl)^\Theta$
elementary operations.
\end{lem}

We  notice that  the coefficients in these expansions
of $x$ lye in $\ZZ$ or $\ZZ[\frac{1+\sqrt 5}{2}]$
and 
there exists a  positive constant $\Theta$   such that
for every integer $k\ge 1$,  the $k$-th coefficient in these
expansions has absolute value $\le \exp((kl)^\Theta )$.

So we have a very accurate description of the singularities 
of $C_l$ since for every branch of $X_1(5l)$ above such 
a singularity we can compute expansions of both coordinates
$b$ and $x$ as series in the local parameter $q_\gamma$.

We shall
also need the following result due to Manin, Shokurov, Merel and
Cremona \cite{manin, merel, cremona, frey}.  

\begin{lem}[Manin, Shokurov, Merel, Cremona]\label{lemma:manin2}
For $l$ a prime and  $p\not\in \{5,l\}$ another prime, the
zeta function of $X_l \bmod p$ can be  computed in deterministic 
polynomial time
in $l$ and  $p$.
\end{lem}

We first  compute the action of the Hecke
operator $T_p$ on the space of Manin symbols for
the congruence group $\Gamma_1(5l)$ associated with  $X_l$.
Then, from the Eichler-Shimura identity $T_p=F_p+p<p>/F_p$ we deduce the characteristic polynomial
of the Frobenius $F_p$. \hfill $\Box$

\section{Power series}\label{sec_torsion_ps}

In this section we give  some notation and we
 state a few  useful  elementary properties 
of  power series in several variables. We are mainly interrested in relating
 the size
of coefficients in the expansions and the size of the values of the series where it converges. In the course of our calculations 
we shall encounter three kinds of power series.
Expansions of normalized  eigenforms have rather small coefficients, according to the Ramanujan conjecture.
To deal with such expansions we introduce in Definition~\ref{Type
of a power series} the {\it type} of a power series.
We shall also consider modular functions having no pole outside the cusps. The coefficients in the
expansions of such functions may be  larger, but they are controled by 
the Petersson and Rademacher's inequality as explained
in Lemma~\ref{lem:PRM} below.
Even more general modular functions may have quite big coefficients. To deal with this case, we introduce
in Definition~\ref{Exp-type of a power series} the exp-type of a power series.

\subsection{The type of a power series}\label{subsection:type}

Let $g\ge 1$ be an integer. The $L^\infty$  norm
of a vector  $\bbx=(x_1,\ldots ,x_g)\in \CC^g$ is $\max_k|x_k|$ 
and it is denoted
 $|\bbx |_\infty$. The $L^1$ norm
is  $|\bbx|_1=\sum_k|x_k|$ and the $L^2$ norm
is  $|\bbx|_2=\sqrt{\sum_k|x_k|^2}$.
We denote by  $\bB \bbx \bB$ the vector
 $(|x_1|, \ldots , |x_g|)$. If
$\bby=(y_1,\ldots ,y_g)$ is another vector
in $\CC^g$ we denote by   $\bbx \star \bby$\index{$\star$, the
componentwise product between two vectors} the
componentwise product 
$(x_1y_1,\ldots , x_gy_g)$.
We set  $\bzero=(0,\ldots ,0)\in \CC^g$ and  $\bun=(1,\ldots ,1)
\in \CC^g$ and $\bdeux=(2,\ldots,2)$.
If $\bn=(n_1,\ldots ,n_g)\in \NN^g$ we write $\bn!$ for the product
$n_1!n_2!\cdots n_g!$ and $\bbx^\bn$ for the product
$x_1^{n_1}\cdots x_g^{n_g}$.
We note
 $P(\bbx,\, \br)=\prod_{k=1}^g 
D(x_k, r_k) \subset \CC^g$\index{$P(\bbx, \, \br)$, the  polydisc
with center  $\bbx$ and polyradius  $\br$}, the  polydisc
with center  $\bbx$ and polyradius  $\br$. When $\br=(r,r,\ldots, r)$ we
just write $P(\bbx,r)$.
 If  $\bbx=(x_1,\ldots ,x_g)\in \RR^g$, we say that
 $\bbx \ge \bzero$ if and only if  $x_k \ge 0$ for every $k$. 
We say that  $\bbx > \bzero$ if and only if 
 $x_k > 0$ for every $k$. 
An entire  series  $f$ in the $g$ variables $x_1$, \ldots,
$x_g$ is  a formal sum
 $f=\sum_{\bk}f_\bk\bbx^\bk$ where the index
$\bk$ runs over $\NN^g$. 
\begin{defi}[Type of a power series in several variables]\index{Type
of a power series}\label{Type
of a power series}
Let  $A$ be a real number $\ge 1$ and
consider   $\bn =(n_1,\ldots ,n_g)\in \NN^g$ such that
$\bn \ge \bun$. We say that $f$ is of type  $(A,\, \bn)$ if
for every  $\bk \ge \bzero$
we have 
$$|f_\bk|\le A(\bk+\bun)^\bn=A\prod _{1\le m \le g}(k_m+1)^{n_m}.$$
\end{defi}
For every  $\bbz\in P(\bzero,\,1)$ we deduce an upper bound for the value of $f$ at $\bbz$.
\begin{eqnarray}\label{eq:majf}
\nonumber |f(\bbz)|\le \sum_{\bk \ge \bzero} A(\bk+\bun)^\bn|\bbz^\bk|&\le& A\prod_{1\le m\le g}\,\,\sum_{k_m\ge 0}(k_m+1)^{n_m}|z_m|^{k_m}\\\nonumber
&\le& \frac{\bn !A}{\prod_m(1-|z_m|)^{n_m+1}}\\
&=&\frac{\bn!A}{(\bun-\bB \bbz  \bB)^{\bn +\bun}}.
\end{eqnarray}
We check  that if $f$ is of type
$(A,\, \ba)$ and  $h$ is of type $(B,\, \bb)$, then
the product  $fh$ is of type  
\begin{equation}\label{eq:typprod}
(AB,\ba+\bb+\bun).
\end{equation}

\subsection{Refocusing a power series}\label{section:refoc}

Given an entire series $f$ of type $(A,\, \bn)$ 
and a vector  $$\bc = (c_1, \ldots , c_g)\in P(\bzero,1)$$ 
we set
\begin{eqnarray*}
F_\bc (\bby)&=&F_\bc (y_1, \ldots , y_g)\\
&=&f(\bc+\bby\star(\bun -\bB \bc \bB))=f((c_m +y_m (1-|c_m|))_m).
\end{eqnarray*}
We call $F_\bc$ the {\it refocused} series of $f$ at $\bc$.
According to \cite[Lemme 16]{Couveignes1} this is a series of type
 $(A_\bc, \, \bn+\bun)$
where 
\begin{equation}\label{eq:recen}
A_\bc = \bn !A\exp (g+|\bn|_1)2^{g+|\bn|_1}(\bun-\bB \bc \bB)^{-\bn-\bdeux}.
\end{equation}
In particular, it converges
for  $\bby \in P(\bzero, 1).$

\subsection{Bounding the remainder}\label{section:boundrem}

For any  integer  $u\ge 0$  
we denote by  $R_u(\bbx)$\index{$R_u(\bbx)$, the remainder
of order $u$} the remainder
of order $u$ of the series $f(\bbx)$. So
$$f(\bbx)=\sum_{|\bk|_1 \le u-1 }f_\bk\bbx^\bk + R_u(\bbx).$$
According to
\cite[Lemme 17]{Couveignes1}, if $f$ is of type
$(A, \bn)$ and if $r$ is a real in $]0,1[$ and if $\bbz \in P(\bzero, r)$  then
\begin{equation}\label{eq:majreste1}
|R_u(\bbz)|\le  B(u+1)^{(n+1)g}|\bbz|_\infty^u
\end{equation}
\noindent where  $n=|\bn|_\infty$ and  $$B=\frac{\bn! 2Ag}{(1-r)^{g+|\bn|_1}}.$$
Further, if $\kappa$ is a real in $]0,1[$ and if
$$u\ge \max (\frac{16(ng)^2}{(\log r)^2},\frac{2(\log \kappa -\log B)}{\log r})$$
\noindent then  
\begin{equation*}%\label{eq:majreste2}
|R_u(\bbz)|\le \kappa \text{ for } \bbz \in P(\bzero, r).
\end{equation*}

\subsection{The type of a quotient}\label{section:typequo}

Let  $f=\sum_{\bk}f_\bk\bbx^\bk$  
be  an entire  series in the $g\ge 2$ variables $x_1$, \ldots, $x_g$.
Assume that $f$ has type $(A,\bn)$ for some $A\ge 1$ and
$\bn\ge \bun$. 
Assume that $f$ is divisible by a polynomial
$P(z_1, \ldots, z_g)$. So there exists an entire  series
$h= \sum_\bk h_\bk \bbx^\bk$ such that $f=Ph$. We want to estimate
the size of coefficients in $h$. We shall only need the case
when $P=z_2-z_1$. So we restrict  to this special case.
Let $\bk=(k_1, \ldots, k_g)$. 
For every integer $m$ such that $2\le m\le g$ we set 
$$u_m=\frac{k_m+1}{k_m+n_m+1}.$$ We set
$$\hat u_1=\frac{k_1+1}{k_1+n_1+1}.$$
If $|\hat u_1-u_2|<|\frac{\hat u_1 +1}{2}-u_2|$
we set $u_1=\frac{\hat u_1+1}{2}$. Otherwise we set
$u_1=\hat u_1$. In any case
$$\frac{1}{|1-u_1|}\le \frac{2(k_1+n_1+1)}{n_1}$$
and
$$\frac{1}{u_1}\le \frac{k_1+n_1+1}{k_1+1}$$
and
$$\frac{1}{|u_2-u_1|}\le \frac{3(k_2+n_2+1)}{n_2}.$$
From Cauchy's  integral  $h_\bk$ is equal to 

$$\frac{1}{(2\pi i)^g}\int_{|\zeta_1|=u_1}
\dots
\int_{|\zeta_g|=u_g} \frac{f(\zeta_1, \zeta_2, \dots, \zeta_g)}{(\zeta_2-\zeta_1)\prod_{1\le m \le g}
\zeta_m^{k_m+1}}
d\zeta_1d\zeta_2\dots d\zeta_g.$$
Using Equation~(\ref{eq:majf})  we find
\begin{eqnarray*}
|f(\zeta_1, \zeta_2, \dots, \zeta_g)|&\le &
\frac{\bn !A}{\prod_{1\le m\le g} (1-u_m)^{n_m+1}}\\
&\le& \bn ! 2^{n_1+1}A\prod_{1\le m\le g}\left(\frac{k_m+n_m+1}{n_m}\right)^{n_m+1}.
\end{eqnarray*}
So
$$|h_\bk|\le \frac{3\bn ! 2^{n_1+1}A(k_2+n_2+1)}{n_2}   
\prod_{1\le m\le g}
\frac{\left(k_m+n_m+1\right)^{k_n+n_m+1}}{\left( k_m+1\right)^{k_m}
\left(n_m \right)^{n_m+1}}.$$

So $h$  has type $(12\bn ! 2^{n_1+|n|_1+g}A \exp(|n|_1), \bn+2\times \bun)$.

Now if we manage to divide $f$ by $K$ polynomials of the
form $z_{j_1}-z_{j_2}$ we obtain a series of type
\begin{equation}%\label{eq:typequo}
(A\exp(8gK(|\bn|_\infty+2K)^2),\bn+2K\times \bun).
\end{equation}
\subsection{The type of derivatives}\label{section:typederi}

If  $f$ is an entire  series in one variable
 of type  $(A,n)$, then the 
derivative $f'$
of $f$ is of type
\begin{equation*}
(A2^{n},n+1).
\end{equation*}
So the  the $d$-th derivative of $f$
is of type 
\begin{equation}\label{eq:majderiv}
(A2^{dn+\frac{d(d-1)}{2}},n+d).
\end{equation}

\subsection{The Petersson and Rademacher's inequality}\label{section:PRjmc}

The modular functions that appear in Lemma~\ref{lemma:bu2r} are not modular forms: they have poles at
the cusps. Since we plan to evaluate these functions 
at well chosen points $z$ in the Poincar{\'e} upper half 
plane, we must control the size of the
coefficients  in the expansions of these functions.

Let $f : \gH\rightarrow \CC$ an holomorphic periodic function
with integer period $e$. 
So 
%$e$ is the smallest positive integer such that 
$f(z+e)=f(z)$ for every $z\in \gH$. 
We assume that $f$ is meromorphic at $\infty$. So $f$
can be written as a series
$$f(z)=\bFF(q_e)=\sum_{\frac{k}{e}\ge v} a_{\frac{k}{e}}q_e^k$$
where $$q_e=\exp(\frac{2i\pi z}{e})$$ and $v\in \frac{1}{e}\ZZ$ is the
valuation  of $f$ at $\infty$ and the $a_{\frac{k}{e}}$ are
the coefficients in the Puiseux expansion of $f$ at $\infty$, and
the series $\bFF$ converges in the unit disk.

Conversely, to every Puiseux series
$$\sum_{\frac{k}{e}\ge v} a_{\frac{k}{e}}q^\frac{k}{e}$$
with radius of convergence $\ge 1$ we can associate an holomorphic
function $f : \gH \rightarrow \CC$ having  period 
%a divisor $e'$ of 
$e$, 
and meromorphic at $\infty$.

By abuse of notation we say that $f$ is a Puiseux series with radius
$\ge 1$ and period $e$.
The set of all such functions is a ring. An example of such 
a series is provided by Klein's modular invariant $\KJ(z)=\bJ(q)$
where $\bJ$ is the series in Equation~(\ref{eq:744}).
In view of Petersson and Rademacher's inequality  (\ref{eq:Petersson}) 
it is natural to state  the following lemma.

\begin{lem}[The  Petersson and Rademacher's property]\label{lem:PR}
Let 
$$f(z)=\sum_{\frac{k}{e}\ge v} a_{\frac{k}{e}}q^{\frac{k}{e}}$$
be a Puiseux series with radius of convergence $\ge 1$. Then the two
properties below are equivalent

\begin{enumerate}
\item There exist two  positive constants $K_1$ and $K_3$  such that
for every integer $k\ge K_3$  we have
\begin{equation}\label{eq:PR1}
\left|a_{\frac{k}{e}}\right|\le \exp \left( K_1\sqrt{\frac{k}{e}}\right).
\end{equation}
\item There exist two  positive constants $K_2$ and $K_4$  such that
for every $z = x+iy\in \gH $ such that $y^{-1}\ge K_4$
we have 
\begin{equation}\label{eq:PR2}
\left| f(z) \right| \le \exp \left(  \frac{K_2}{y} \right).
\end{equation}
\end{enumerate}

We say that such a Puiseux series is PR  (like Petersson and Rademacher).
The set of PR series is a ring which is integrally closed in the ring
of Puiseux series with radius of convergence $\ge 1$.
\end{lem}

This lemma is easily proven: one implication results from  Cauchy's
 formula and
the other implication is trivial.
 Using  Lemma~\ref{lem:PR} and the Petersson-Rademacher
inequality~(\ref{eq:Petersson}) we prove the 
following lemma.

\begin{lem}[Integral modular forms are PR]\label{lem:PRM}
Let $N\ge 1$ be an integer.
Let $f : \gH\rightarrow \CC$  be a modular function of weight
$0$ for the group $\Gamma_1(N)$. Assume that $f$ has no pole outside
the cusps. So $f$ is holomorphic on $\gH$, or equivalently
it belongs to the ring $\CC[Y_1(N)]$ of integral functions on $Y_1(N)$.
Then the expansion of $f$ at any cusp is PR.
\end{lem}

So the Petersson-Rademacher inequality for the Fourier
coefficients of $\KJ$ holds true for any integral function on $Y_1(N)$.
We also prove the following lemma concerning the discriminant form.

\begin{lem}[The discriminant and its inverse
 are PR]\label{lem:PRDelta}
Let $$\Delta(q)= q \prod_{n \ge 1} (1-q^n)^{24} = \sum_{n \ge 1} \tau(n) q^n$$
be the discriminant form. The inverse $\Delta^{-1}$ of $\Delta$
is a PR series in $q$.
\end{lem}

This results from the product formula for $\Delta$.
We deduce that Lemma~\ref{lem:PRM} extends to modular functions
of arbitrary weight.

\begin{lem}[Integral modular forms are PR]\label{lem:PRM2}
Let $N\ge 1$ be an integer.
Let $f : \gH\rightarrow \CC$  be a modular function of weight
$k\in \ZZ$ for the group $\Gamma_1(N)$. Assume that $f$ has no pole outside
the cusps. 
Then the expansion of $f$ at any cusp is PR.
\end{lem}

\subsection{The exp-type of a power series}\label{subsection:exptype}

We shall have to deal with series in one variable having bigger
coefficients than the ones introduced in Sections~\ref{subsection:type}
and \ref{section:PRjmc}.
%The series $\bJ$ introduced in equation (\ref{eq:744}) is a typical example.
The successive derivatives 
appearing in Lemma~\ref{lemma:algdep} are a good example. We no longer
care about convergence then. We just want to control the size
of the (logarithms) of the coefficients.

\begin{defi}[Exp-type of a power series in one
 variable]\index{Exp-type of a power series}\label{Exp-type of a power series}
Let  $n$ be
an integer $\ge 2$ and let  $A$ and  $B$ be   two
 real numbers $\ge 1$. We say that an entire
 series $f=\sum_{k\ge 0} f_kx^k$ is of exp-type  $(A, B, n)$ if
for every  $k \ge 0$
we have 

$$|f_k|\le \exp((Ak+B)^n).$$
\end{defi}

If $f_1$ is of exp-type $(A_1, B_1, n_1)$ and
$f_2$ of exp-type $(A_2, B_2, n_2)$ then the product
$f_1f_2$ is of exp-type 

\begin{equation}\label{eq:exptypeprod}
(A_1+A_2, B_1+B_2 +1, n)
\end{equation} where 
$n$ is the maximum of $n_1$ and $n_2$.

If $f$ is of exp-type $(A, B, n)$ and $k\ge 1$ is an integer, then
the $k$-th power $f^k$ is of exp-type $(kA, kB+k-1, n)$.
The derivative $f'=df/dx$ of $f$ is of exp-type 
\begin{equation}%\label{eq:exptypederiv}
(A,A+B+1, n).
\end{equation}

If $f$ is of exp-type $(A, B, n)$ and $g(x)=1/(1-xf(x))$, then
$g$ is of exp-type 
\begin{equation}\label{eq:exptypeinv}
(\sqrt{A+B+1}, 0,  2n).
\end{equation}

\section{Jacobians and Wronskians of power series}\label{section:JW}

In this section 
we state and prove an algebraic identity
relating 
Wronskian and Jacobian determinants. This identity
will be useful to control the local behaviour of the Jacobi
integration map.
We first state this identity  in its simplest 
and most natural form in Section~\ref{section:JWsimple}.
We then state  a more general identity in Section~\ref{section:JWgene}.
The proofs are given in Sections~\ref{section:JWP1} and
\ref{section:JWP2}.

\subsection{A special  case}\label{section:JWsimple}

We assume that  $g\ge 2$ is an integer and
we consider 
$g$ entire  series
 $f_1(x)$, $f_2(x)$, \ldots, $f_g(x)$ in one variable $x$, with coefficients
in $\CC$.
The {\it Wronskian} associated with
 $\bff = (f_1, \ldots, f_g)$  is the determinant
\begin{equation}\label{eq:Wfdef}
W_\bff(x)=\frac{1}{\prod_{1\le m\le g-1}m!}\times \left|  \begin{array}{ccc} f_1(x) & \dots & f_g(x)\\ f'_1(x) & \dots & f'_g(x)\\ \vdots & &\vdots \\
f^{(g-1)}_1(x) & \dots & f^{(g-1)}_g(x)\end{array} \right|.
\end{equation}
We may also introduce $g$
indeterminates $y_1$, $y_2$, \ldots, $y_g$ and define
the {\it Jacobian} associated with  $\bff$ to be the determinant
$$\cJ_\bff = \left|  \begin{array}{ccc} f_1(y_1) & \dots & f_g(y_1)\\ f_1(y_2) & \dots & f_g(y_2)\\ \vdots & &\vdots \\
f_1(y_{g}) & \dots & f_g(y_g)\end{array} \right|.$$\index{$\cJ_\bff $, the
Jacobian associated with the family $\bff$}

Let now  $D=\prod_{k<m}(y_m-y_k)$ be the reduced discriminant.
% and assume that  the $f_k(x)$ are entire series. 
Then
the Jacobian  $\cJ_\bff$ lies in the ring
$\CC[[y_1,\ldots,y_g]]$ and is divisible by
the reduced discriminant in this
ring. Further, the  quotient $\cJ_\bff/D$ is congruent
to $W_\bff(0)$ modulo the maximal 
ideal of  $\CC[[y_1,\ldots,y_g]]$:
\begin{equation}\label{eq:WJ}
\frac{\cJ_\bff}{D}\equiv  W_\bff(0)\bmod (y_1, y_2, \ldots, y_g)\CC[[y_1,\ldots,y_g]].
\end{equation}
A proof of this identity 
is given in the next  Sections~\ref{section:JWP1} and 
\ref{section:JWP2}.

Now assume that  all the series $f_k(x)$ have type $(A,n)$. Then the Jacobian
$\cJ_\bff$ has type 
\begin{equation*}%\label{eq:typejac}
(g!A^g,n\bun).
\end{equation*}
We deduce from Equation~(\ref{eq:majderiv}) and (\ref{eq:typprod})
that the Wronskian $W_\bff$ has type
\begin{equation}\label{eq:typewron}
(A^g\exp(\Theta ng^2+\Theta g^3),gn-1+\frac{g(g+1)}{2}).
\end{equation}

\subsection{A more general identity}\label{section:JWgene}

Let $g\ge 2$ be an integer. Let $\KK$ be any field with characteristic
zero. Let 
$f_1(x)$,
\ldots, $f_g(x)$ be $g$ entire  series in one variable $x$ having coefficients
in $\KK$.
Let $\bff=(f_1,\ldots,f_g)$ be the corresponding vector.
Let $1\le n\le g$  be an integer and consider the $n$ first derivatives
$$\bff^{(0)}, \, \bff^{(1)}, \ldots,
\bff^{(n-1)}$$
of $\bff$ with respect to the variable $x$. The exterior product
\begin{equation}\label{eq:partW}
W_{\bff,n}=\frac{\bff^{(0)}}{0!}\wedge\frac{\bff^{(1)}}{1!}\wedge  \dots\wedge \frac{\bff^{(n-1)}}{(n-1)!}\in 
\bigwedge^{n}\left( \KK[[x]] \right)^g
\end{equation}
is a sort of {\it partial Wronskian}  associated with the vector $\bff$.

Now let $S\ge 1$ be an integer and let 
$$n=m_1+m_2+\dots+m_S$$
be a partition of $n$ in $S$ parts. In particular $m_s$ is a positive
integer for every $1\le s\le S$.
Let $y_1$, $y_2$, \ldots, $y_S$ be $S$ distinct indeterminates and consider
the corresponding {\it partial Jacobian}
\begin{equation}\label{eq:partJ}
\cJ_{\bff,(m_s)_{1\le s\le S}}=W_{\bff,m_1}(y_1)\wedge W_{\bff,m_2}(y_2)\dots\wedge W_{\bff,m_S}(y_S) 
\end{equation}
in
$$  \bigwedge^{n} \left( \KK[[y_1,\ldots, y_S]] \right)^g.$$
Let 
$$D_{(m_s)_{1\le s\le S}}=\prod_{1\le s_1<s_2\le S}(y_{s_2}-y_{s_1})^{m_{s_1}m_{s_2}}$$ be the corresponding 
{\it partial weighted discriminant}. The following identity is a partial
generalization of Equation~(\ref{eq:WJ})

\begin{equation}\label{eq:WJpart}
\frac{\cJ_{\bff, (m_s)_{1\le s\le S}}}{D_{(m_s)_{1\le s\le S}}}\equiv  W_{\bff,n}(0)\bmod (y_1, y_2, \ldots, y_S).
\end{equation}

In particular $\cJ_{\bff, (m_s)_{1\le s\le S}}$
is divisible by the weighted discriminant 
  $D_{(m_s)_{1\le s\le S}}$ in the 
$\KK[[y_1,\ldots,y_S]]$-module  $\bigwedge^n \left( \KK[[y_1,\ldots,y_S]]\right)^g$. A proof of this identity 
 given in the next  Sections~\ref{section:JWP1} and 
\ref{section:JWP2}.

Now let $K\ge 1$ be an integer and let 
$$g=n_1+n_2+\dots+n_K$$
be a partition of the dimension $g$. So
  $n_k$ is a positive integer  for each $1\le k\le K$.
We consider $K$ vectors $\bff_1$, \ldots, $\bff_K$ in
$(\KK[[x]])^g$. 
We introduce $K$ indeterminates $x_1$, \ldots,
$x_K$ and following  Equation~(\ref{eq:partW}) 
we define the {\it total Wronskian} to
 be

\begin{equation}\label{eq:totalW}
W_{(\bff_k, n_k)_{1\le k\le K}} = \bigwedge_{1\le k\le K} 
W_{\bff_k,n_k}(x_k)\in \KK[[x_1, x_2, \ldots, x_K]].
\end{equation}

For every $1\le k\le K$ let $S_k$ be a positive integer and let
$$n_k=m_{k,1}+m_{k,2}+\dots+m_{k,S_k}$$
be a partition of $n_k$ into $S_k$ parts.
For every $1\le k\le K$ and $1\le s\le S_k$ we  introduce the new 
indeterminate $y_{k,s}$ and following Equation~(\ref{eq:partJ}) we define the {\it total Jacobian} to be

\begin{equation}\label{eq:totaljacob}
\cJ_{(\bff_k,(m_{k,s})_{1\le s\le S_k})_{1\le k\le K}}=
\bigwedge_{1\le k\le K}  \cJ_{\bff_k,  (m_{k,s})_{1\le s\le S_k}}
\end{equation}
in
$$\KK[[(y_{k,s})_{1\le k\le K;\, 1\le s\le S_k}   ]].$$

We set $$\bn=(n_k)_{1\le k\le K}$$ and 
$$\bm_k=(m_{k,s})_{1\le s\le S_k}$$ and 
$$\bm=(m_{k,s})_{1\le k\le K; \, 1\le s\le S_k}.$$ Both
$\bm$ and $\bn$ are partitions 
 of $g$. And $\bm$ is a refinement of $\bn$.
We define the discriminant relative to $\bm$ and $\bn$ to be
$$D_{\bm,\bn}=\prod_{1\le k\le K}D_{\bm_k}=
\prod_{1\le k\le K}\,\,\prod_{1\le s_1<s_2\le S_k}(y_{k,s_2}-y_{k,s_1})^{m_{k,s_1}m_{k,s_2}}.$$

Collecting  $K$ equations like
(\ref{eq:WJpart}) we obtain

\begin{equation}\label{eq:WJtot}
\frac{\cJ_{(\bff_k, \bm_k)_{1\le k\le K}}}{D_{\bm,\bn}}\equiv 
W_{(\bff_k,n_k)_{1\le k\le K}}(0)\bmod (y_{k,s})_{1\le k\le K; \, 1\le s\le S_k}.
\end{equation}

This generalization of Equation~(\ref{eq:WJ}) will be useful
in  Section~\ref{subsection:defons} when studying the 
Jacobi map. We now prove it using a formal analogue of the
Jacobi map.

\subsection{Proof of Equation (\ref{eq:WJpart}) in a special case}\label{section:JWP1}

In this section  assume that $S=n$ and $m_1=m_2=\dots=m_S=1$.
We write $\bff$ as a series in $x$ with coefficients in $\KK^g$
$$\bff=\Psi_{0}+\Psi_1 x+ \Psi_2 x^2+\dots \in \KK^g[[x]]$$
where $$\Psi_k=\frac{\bff^{(k)}}{k!}(0)\in \KK^g$$
for every $k\ge 0$.
We consider the formal integration 
\begin{equation*}%\label{eq:bFF}
\bFF = \int_0^x\bff(x)dx=\Psi_0 x+\frac{\Psi_1}{2}x^2+\frac{\Psi_2}{3}x^3+\dots \in \KK^g[[x]].
\end{equation*}
We introduce $S$ new indeterminates $y_1$, \ldots, $y_S$ and we  set
\begin{equation}\label{eq:defPhi}
\Phi(y_1, \ldots, y_S)=\bFF(y_1)+\dots+\bFF(y_S)\in \KK^g[[y_1,\ldots, y_S]].
\end{equation}
We denote by \begin{equation*}%\label{eq:mgot}
\mgot=(y_1, \ldots, y_S)\KK[[y_1, \ldots, y_S]]
\end{equation*}
the maximal ideal in $\KK[[y_1, \ldots, y_S]]$.
For every  $k\ge 1$ we call
$$\nu_k=y_1^k+\dots+y_S^k$$
the $k$-th Newton's power sum.
We check that
$$\Phi(y_1, \ldots, y_S)=\Psi_0 \nu_1+ \frac{\Psi_1}{2} \nu_2+
\dots+ \frac{\Psi_{S-1}}{S}\nu_S +R$$
where the remainder  $R$ is a vector in $\left(\KK[[y_1, \ldots, y_S]]\right)^g$ whose
coefficients are symmetric functions in the $(y_s)_{1\le s\le S}$ and belong to  $\mgot^{S+1}$.
So these coefficients
 belong to $\KK[[\nu_1, \ldots,\nu_S]]$ and for every 
$1\le k\le S$ the partial derivative 
$\frac{\partial R}{\partial \nu_k}$ is zero modulo $$\ngot = \mgot\cap \KK[[\nu_1, \ldots,\nu_S]]=   (\nu_1, \ldots, \nu_S) \KK[[\nu_1, \ldots,\nu_S]].$$
We deduce that for every $1\le k\le S$
\begin{equation}\label{eq:JW1}
\frac{\partial \Phi}{\partial \nu_k}\equiv \frac{\Psi_{k-1}}{k}\bmod \ngot.
\end{equation}
On the other hand, it is clear from the definition of $\Phi$ in 
Equation~(\ref{eq:defPhi})
that for every $1\le s\le S$
\begin{equation}\label{eq:JW2}
\frac{\partial \Phi}{\partial y_s}= \frac{d \bFF}{dx}(y_s)=\bff(y_s).
\end{equation}
Finally, for every $1\le k \le S$  and $1\le s\le S$ we have
$$\frac{\partial \nu_k}{\partial y_s} =k y_s^{k-1}$$
and the determinant
\begin{equation}\label{eq:JW3}
\left|\frac{\partial \nu_k}{\partial y_s}\right|_{k,s} =S!\prod_{1\le s_1<s_2\le S}(y_{s_2}-y_{s_1}).
\end{equation}

Equation~(\ref{eq:WJpart}) then follows from Equations~(\ref{eq:JW1}), (\ref{eq:JW2}) and
(\ref{eq:JW3}) applying the chain  rule for derivatives.

\subsection{Proof of Equation (\ref{eq:WJpart})  in general}\label{section:JWP2}

We introduce the $n$ indeterminates
$x_{s,j}$
for $1\le s\le S$ and $1\le j \le m_s$.
We  put the lexicographic order on these indeterminates and we apply Equation 
(\ref{eq:WJpart})
to the series $\bff$ and partition 
$n=1+1+\dots+1$. We obtain
\begin{eqnarray}\label{eq:genn1}
\scriptstyle \nonumber \bff (x_{1,1})&\scriptstyle  \wedge& 
\scriptstyle \dots \wedge \bff (x_{S,m_S})
=\left( W_{\bff , n}(0)+O(x_{1,1}, \ldots, x_{S,m_S})\right)\\
&\scriptstyle  \times & \scriptstyle \prod_{(1,1)\le (s_1,j_1)<(s_2,j_2)\le (S,m_S)}
(x_{s_2, j_2}-x_{s_1,j_1})
\end{eqnarray}
where $O(x_{1,1}, \ldots, x_{S,m_S})$ stands for any element in
the ideal 
generated by $x_{1,1}$,
\ldots, $x_{S,m_S}$ in  $\KK[[x_{1,1}, \ldots, x_{S,m_S}]]$ .

We introduce $S$ new indeterminates $y_1$,
\ldots, $y_S$. We also introduce $n$ indeterminates
$z_{s,j}$
for $1\le s\le S$ and $1\le j \le m_s$.
For every $1\le s\le S$ we consider $\bff(y_s+z)$ as
a series in $z$ with coefficients in $\KK[[y_s]]$.
We apply Equation~\ref{eq:WJpart}
to the series $\bff(y_s+z)$ and partition 
$m_s=1+1+\dots+1$. We obtain
\begin{eqnarray}\label{eq:genn2}
\scriptstyle  
\nonumber  \bff(y_s+z_{s,1})\wedge \dots \wedge \bff (y_s+z_{s,m_s})&
\scriptstyle =& \scriptstyle 
\left( W_{\bff , m_s}(y_s)+O(z_{s,1}, \ldots, z_{s,m_s})\right) \\
&\scriptstyle  \times&\scriptstyle   \prod_{1\le j_1 <j_2\le m_s}
(z_{s, j_2}-z_{s,j_1})
\end{eqnarray}
We now replace $x_{s,j}$ by $y_s+z_{s,j}$ in Equation~(\ref{eq:genn1}) and we obtain
\begin{eqnarray}\label{eq:genn3}
\scriptstyle \nonumber \bff (y_{1}+z_{1,1})&\scriptstyle \wedge &\scriptstyle \dots \wedge \bff (y_{S}+z_{S,m_S})
=\left( W_{\bff , n}(0)+O(y_1, \ldots, y_S, z_{1,1}, \ldots, z_{S,m_S})\right)\\
\scriptstyle & \scriptstyle \times &\scriptstyle \prod_{(1,1)\le (s_1,j_1)<(s_2,j_2)\le (S,m_S)}
(y_{s_2}-y_{s_1}+z_{s_2, j_2}-z_{s_1,j_1})
\end{eqnarray}
We now notice
that  the left hand side
of Equation~(\ref{eq:genn3}) is the  wedge 
product
of the left hand sides of the $S$ equations
like (\ref{eq:genn2}). Further, the  
discriminants on the right hand sides  of the $S$
equations like (\ref{eq:genn2})  divide the discriminant
in the right hand side of Equation~(\ref{eq:genn3}).
So we equate the right hand side
of Equation~(\ref{eq:genn3}) and the  wedge 
product
of the right hand sides of the $S$ equations
like (\ref{eq:genn2}). We then divide by the product of
the $S$ small discriminants. We then 
 reduce modulo the ideal generated by the variables
$z_{s,j}$ and we obtain
$$\scriptstyle {\bigwedge_{1\le s\le S} W_{\bff , m_s}(y_s) }=\left( {W_{\bff , n}(0)+O(y_1, \ldots, y_s)} \right)\times \prod_{1\le s_1<s_2\le S}(y_{s_2} - y_{s_1})^{m_{s_1}m_{s_2}}$$
as  was to be proven.

\section{A simple quantitative study of the Jacobi map}\label{section:quantjacob}

In this section we prove some upper and lower bounds for the
Jacobi map. Upper bounds are rather trivial but important to control
the complexity of the algorithms. Lower bounds are not very surprising
either. But they play an important role in the proof of Theorem~\ref{theorem:inverseJ}.

\subsection{Upper bounds for the Jacobi map} 

We first prove that the Jacobi map is Lipschitz
with constant $l^\Theta$. Indeed let $\gamma\in \Xi$ and call
$\phi_\gamma = \phi \circ \mu_\gamma$ 
the composition of the Jacobi map $\phi : X \rightarrow J(\CC)$ 
with the  modular parameterization 
$\mu_{\gamma} :  D(0,1)\rightarrow X$. So 
$$\phi_\gamma : q_\gamma \mapsto \left( \int^{q_\gamma } \omega 
\right)_{\omega\in \cB^1_{\rm DR}}$$
and  we don't need  to specify the origin of the integral here. Every 
$\omega = f(q)q^{-1}dq$ 
in $\cB^1_{\rm DR}$ can be written $h(q_\gamma)q_\gamma^{-1}dq_\gamma$
where $h(q_\gamma)$ is the expansion of the
 modular form $f$ at the cusp
$\gamma(\infty)$. It is a consequence of Ramanujan's conjectures
proven by Deligne and the explicit formulae by Asai
for pseudo-eigenvalues that the series $h(q_\gamma)q_{\gamma}^{-1}$ is
of type 

\begin{equation}\label{eq:typedeli}
(\Theta,\Theta).
\end{equation}

 Therefore if $q_1$ and $q_2$ belong to 
$D(0,\exp(-\pi/w_\gamma))\supset F_{w_\gamma}$ the integral
$\int_{q_1}^{q_2}h(q_\gamma)q_{\gamma}^{-1}dq_\gamma$ is bounded
in absolute value by $|q_2-q_1|$ times $l^{\Theta}$ according to 
Equation~(\ref{eq:majreste1}). So the Jacobi map is Lipschitz with
constant 

\begin{equation}\label{eq:lipjac}
\le l^\Theta.
\end{equation}

We now consider
some vector   $\hgamma = (\gamma_k)_{1\le k\le g} \in \Xi^g$ 
 and we call $$\phi_\hgamma : D(0,1)^g\rightarrow \CC^{\cB^1_{\rm DR}}/\cR$$ the composition
of the Jacobi map $\phi : X^g \rightarrow J(\CC)$ with the product of
the modular parameterizations $$\mu_{\hgamma} = \prod_{1\le k\le g} \mu_{\gamma_k} :D(0,1)^g\rightarrow X^g.$$

For every $1\le k\le g$ we set $w_k=w_{\gamma_k}$.
Let  $\bq=(q_1, \ldots, q_g)$ with
$q_k\in D(0,\exp(-\pi/w_k))$ for every $1\le k\le g$.
We study the map $\phi_\hgamma$ locally at $\bq$.
The tangent space to $D(0,1)^g \subset \CC^g$ at $\bq$ is identified with $\CC^g$ and we denote
by $(\delta_k)_{1\le k\le g}$ its canonical basis. 
The underlying $\RR$-vector space has basis
$(\delta_1, \ldots, \delta_g, i\delta_1, \ldots, i\delta_g)$. For $1\le k\le g$ we set $\delta_{k+g}=i\delta_k$.
Similarly we call  $(e_k)_{1\le k\le g}$  the canonical
basis of $\CC^g=\CC^{\cB^1_{\rm DR}}$ and set  $e_{k+g}=ie_k$ for $1\le k\le g$. So  $(e_k)_{1\le k\le 2g}$ is a  basis of the $\RR$-vector space underlying $\CC^g$.
Let $D_\bq\phi_\hgamma$ be the differential
of $\phi_\hgamma$ at $\bq$.
For $1\le k\le 2g$ let $\rho_k$ be the image of $\delta_k$ by this  
differential. 
The determinant  of $D_\bq\phi_\hgamma$ 
$$\frac{\rho_1\wedge \ldots \wedge \rho_{2g}}{e_1\wedge \ldots \wedge e_{2g}}$$
is the square of the absolute value of the Jacobian determinant
\begin{equation}\label{eq:jacdefi}
\cJ_\hgamma(\bq)=\left| \frac{\omega}{dq_{\gamma_k}}(q_k) \right|_{1\le k\le g, \, \omega\in \cB^1_{\rm DR}.}
\end{equation}
We shall need  an upper bound for the absolute value of 
$\cJ_\hgamma(\bq)$. 
As a series in the $g$ indeterminates $q_1$, $q_2$,
\ldots, $q_g$, the Jacobian $\cJ_\hgamma$ has type
\begin{equation}\label{eq:typejac2}
(\exp(\Theta g^2),\Theta\bun).
\end{equation}
This results from Equation~(\ref{eq:typedeli}) and  the 
definition
of $\cJ_\hgamma$  in Equation~(\ref{eq:jacdefi}).
If $\hgamma = (\gamma, \gamma, \ldots,
\gamma)$ is the repetition of $g$ times the same $\gamma$ in $\Xi$,
we  write $\cJ_\gamma$ for $\cJ_\hgamma$ and we 
denote by $W_\gamma (q)$ the  Wronskian associated with
$\cJ_\gamma$. This is a series in one variable $q$. We deduce
from Equations~(\ref{eq:typewron}) and  (\ref{eq:typedeli}) that
$W_\gamma$ has type
\begin{equation}\label{eq:typewron2}
(\exp(\Theta g^3),\Theta g^2).
\end{equation}

\subsection{Lower bounds for the Jacobi map}\label{sec:LBJM}

We now  bound from below  the
Wronskian $W_\gamma (q)$ and  the Jacobian 
$\cJ_\gamma(\bq)$ for some special values of  $q$ and $\bq$.
So we   assume that $\hgamma=(\gamma, \gamma, \ldots, \gamma)$
is the repetition of $g$ times the same $\gamma$ in $\Xi$ and 
we study 
 the Jacobi map $\phi : X^g \rightarrow J$
in the neighborhood of 
$(\gamma(\infty) , \ldots, \gamma(\infty))$ where  $\gamma(\infty)$ is the
cusp associated with $\gamma$. 
The parameter at the cusp $\gamma(\infty)$
is the  $q_\gamma$ from Equation~(\ref{eq:qgamma}).

We first treat the case when $\gamma=\Id$ and we write
$\cJ_\infty$ (resp. $W_\infty$) for $\cJ_\Id$ (resp. $W_\Id$).
We denote by  $$G_g=\prod_{1\le m\le g-1}m!$$ the denominator that
appears in the definition of the Wronskian in Equation~(\ref{eq:Wfdef}).
The expansions in $q=q_\Id$
of the  $\frac{\omega}{dq} $ for $\omega\in \cB^1_{\rm DR}$
are entire  series with algebraic integer coefficients; and they are permuted
by the absolute Galois group of $\QQ$. So the series 
$G_g^2\times W_\infty(q)^2$ 
has coefficients in $\ZZ$.
 
We set  $\eta=\frac{g(g+1)}{2}$.   
The product  $W_\infty (q)(dq)^\eta$ 
is a degree $\eta$  holomorphic
form on $X$.
Therefore it has  $2(g-1)\eta$ zeros counting multiplicities.
Since the $q$-valuation $v$  of $W_\infty (q)$
is  the multiplicity 
of the cusp
$\infty$ in the divisor of 
$W_\infty (q)(dq)^\eta$, we deduce that
$v \le 2(g-1)\eta$.
So the series $G_g^2\times W_\infty(q)^2$ has valuation $2v\le 2(g-1)g(g+1)$
and rational integer coefficients. We deduce from Equations~(\ref{eq:typewron2})
and (\ref{eq:typprod}) that the type of the series $W_\infty(q)^2$
is 
\begin{equation}\label{eq:typeW2}
(\exp(\Theta g^3),\Theta g^2).
\end{equation}
We write
$$W_\infty (q)^2=\frac{c}{G_g^2 }\times q^{2v}+R_{2v+1}(q)$$ where $c$ is a non-zero rational 
integer and $R_{2v+1}$ is the remainder of order $2v+1$. We can bound this
remainder using Equations~(\ref{eq:majreste1}) and (\ref{eq:typeW2}).
$$\left| R_{2v+1}(q)\right| \le \exp(\Theta l^6)|q|^{2v+1}.$$
So if $|q|\le \exp(-\Theta l^6)$ we have
  $|W_\infty (q)^2|\ge \frac{1}{2G_g^2}\times q^{2(g-1)g(g+1)}$  and

\begin{equation}\label{eq:infWinf}
|W_\infty (q)|\ge \frac{1}{2G_g}\times q^{(g-1)g(g+1)}\ge \exp(-\Theta l^5) \times q^{(g-1)g(g+1)}.
\end{equation} So we fix such 
a $q$. For example we take 
\begin{equation}\label{eq:leq}
q=10^{-\kappa_1 l^6}
\end{equation}
for some large enough positive constant $\kappa_1$.
We set  
\begin{equation}\label{eq:lebq}
\bq=q\bun =(q, \ldots, q)
\end{equation}
 and
 $\bbx = (x_1, \ldots, x_g)$ where $x_1$, \ldots, $x_g$ are new 
indeterminates. 
The  Jacobian $$\cJ_\infty(\bq+\bbx\star(\bun -\bB \bq \bB ))$$ 
is an entire series
in the  $g$ variables $x_1$, \ldots, $x_g$. 
This is indeed the Jacobian associated with the $g$
series  $$(f(q+x(1-|q|))/(q+x(1-|q|)))_{f(q)q^{-1}dq \in \cB^1_{\rm DR}}$$
in the variable $x$. 
Equation~(\ref{eq:WJ}) gives us the first non-zero term in the expansion
of this series at $\bbx=\bzero$:
\begin{eqnarray*}
\cJ_\infty(\bq+\bbx \star (\bun -\bB \bq \bB ) )=&&W_\infty(q)(1-|q|)^{\frac{g(g-1)}{2}}\prod_{k<m}(x_m-x_k) \\&+&R_{\frac{g(g-1)}{2}+1}(\bbx).
\end{eqnarray*}

The type of the 
Jacobian  $\cJ_\infty(\bq)$  as  a series
in  $\bq$  is given by Equation~(\ref{eq:typejac2}).
We deduce from Equation~(\ref{eq:recen}) that the refocused series
$\cJ_\infty(\bq+\bbx\star(\bun - \bB \bq \bB ))$ 
is a series in  $\bbx$ of type
 $(\exp(\Theta g^{2}), \Theta \bun )$.
Using Equation~(\ref{eq:majreste1}) we deduce that
for  $\bbx $ in  $P(\bzero ,\exp(-\pi))$
$$\left| R_{\frac{g(g-1)}{2}+1}(\bbx)\right| \le \exp({\Theta}g^2)|\bbx|_\infty^{\frac{g(g-1)}{2}+1}.$$

We set  $s=|\bbx|_\infty$  and we  assume
that  $\bbx$ takes the special form

\begin{equation}\label{eq:lebx}
\bbx= (\frac{s}{g},\frac{2s}{g}, \ldots, \frac{(g-1)s}{g}, s)
\end{equation}
and  $s\le \exp(-\pi)$.
Then

\begin{eqnarray*}
\left|W_\infty(q)(1-|q|)^{\frac{g(g-1)}{2}}\prod_{k<l}(x_l-x_k)\right| &\ge& \left| W_\infty(q)\right|
\left(\frac{s(1-|q|)}{g}\right)^{\frac{g(g-1)}{2}}\\
&\ge & \Theta^{-1}\left| W_\infty(q)\right| \left(\frac{s}{g}\right)^{\frac{g(g-1)}{2}}
\end{eqnarray*}
We take 

\begin{equation}\label{eq:les}
s=10^{-\kappa_2 l^{12}}
\end{equation}
for some large enough positive constant $\kappa_2$. Using 
Equations~(\ref{eq:leq}) and (\ref{eq:infWinf}) we obtain the following lower bound for
the Jacobian

\begin{equation}\label{eq:uninfJ}
|\cJ_\infty(\bq+\bbx  (1 -q  ) )|\ge \exp(-\Theta l^{16})
\end{equation}
when $q$, $\bq$, and $\bbx$ are given by 
Equations~(\ref{eq:leq}), (\ref{eq:lebq}), (\ref{eq:lebx}), and (\ref{eq:les}). In particular $\left|\bq+\bbx  (1 -q  )\right|_\infty$
 can be assumed to be $\le \exp(-2\pi)$. So
\begin{equation}\label{eq:infJinf}
\max_{\bq \in D(0,\exp(-2\pi))^g} \left|\cJ_\infty(\bq) \right|
\ge \exp(-\Theta l^{16}).
\end{equation}

In order to  bound from below  $\cJ_\gamma$ for any $\gamma \in \Xi$
we observe that $\cJ_\gamma$ and $\cJ$ are closely related:  If
$w$ is the width of the cusp $\gamma(\infty)$,  there
exists a $w$-th root of unity $\zeta_\gamma$ and an algebraic
number $\lambda_\gamma$ of absolute value $1$ (the product of all
pseudo-eigenvalues) such that  the following formal identity
in $\CC[[q_1, \ldots, q_g]]$ holds true
\begin{equation}\label{eq:transJJ}
\cJ_\gamma(q_1, \ldots, q_g)=\lambda_\gamma \cJ_\infty(\zeta_\gamma q_1,
\ldots, \zeta_\gamma q_g).
\end{equation}
So the lower bound in Equation~(\ref{eq:infJinf}) is also 
valid for every $\cJ_\gamma$.

\section{Equivalence of various norms}\label{sec_torsion_QJ}

The main algorithm 
  in this text (the one in Section~\ref{section:opjac})
  uses  a subroutine that computes the complex roots
of an analytic function  on a compact set.
 This problem is well conditioned according to
Lemma~\ref{lemma:anyzero} of Chapter~\ref{sec_couveignes_ZEROS},
 provided
we have a decent lower bound for the maximum of the function in question.
In our situation, the analytic functions  are derived from quadratic
differentials on $X$. We need
simple conditions for these functions
not to be uniformly small in absolute value in the neighborhood of any
cusp. The second inequality 
in Equation~(\ref{eq:equivnorm}) below  provides such  a condition.
In order to prove this inequality we 
study Jacobians associated
with weight $4$ cusp forms on $X$, locally at every cusp.

\subsection{Space of quadratic differentials}\label{section:SQD}

In this section we shall make use of parabolic
modular forms of weight $4$
on $X$. To every such form $f(q)$ one can associate a quadratic differential
$\omega= f(q)q^{-2}(dq)^2$. The divisor of $\omega$ is related to the
divisor of $f$ by the following relation
$$\Div (\omega)=\Div(f)-2\Cusps$$ where
$\Cusps$\index{$\Cusps$, the divisor of cusps on $X_1(5l)$} is the
sum of all cusps. Note that $X_1(5l)$ has no elliptic point.
The map $f(q)\mapsto f(q)q^{-2}(dq)^2$ defines a bijection
between the space of weight $4$ parabolic forms and the space
$\cH^2(\Cusps)$ of quadratic differentials with divisor
$\ge - \,  \Cusps$.  We denote by $g_2$ the dimension of the latter space.
This is $3g-3$ plus the degree of $\Cusps$ (the number of cusps).

We shall need a basis $\cB_{\rm quad}$\index{$\cB_{\rm quad}$, a basis for
$\cH^2(\Cusps)$}
for the space $\cH^2(\Cusps)$  or equivalently a basis for the space 
$S_4(\Gamma_1(5l))$ of weight
four  cusp forms.   We shall again use the standard basis made
of normalized newforms of level $5l$ together
with normalized newforms of level $l$ lifted
to level $5l$ by the two degeneracy maps.
We  also need the expansion of every form 
in $\cB_{\rm quad}$
at every cusp $\gamma(\infty)$ for $\gamma \in \Xi$.
More precisely, $f(q)q^{-2}(dq)^2$ should be rewritten as 
$h(q_\gamma)q_\gamma^{-2}(dq_\gamma)^2$ for 
every $\gamma$ in $\Xi$.
As for degree $1$ forms, and using the same methods,
we can
compute the expansion of all quadratic forms in $\cB_{\rm quad}$ at all
cusps in deterministic polynomial time $(klm)^\Theta$ where $k$
is the $q$-adic accuracy and $m$ the complex absolute
 accuracy of coefficients.

We now define several important norms on the space
 $\cH^2(\Cusps)$. 
If $\omega$ is a form in  $\cH^2(\Cusps)$,
we denote by $|\omega|_\infty$ the $L^\infty$ norm
in the basis $\cB_{\rm quad}$. To every cusp $\gamma(\infty)$
with $\gamma\in\Xi$ we associate  a norm  on 
$\cH^2(\Cusps)$. We  define  $|\omega|_\gamma$
to be the maximum of the modulus of the function
$\omega q_\gamma(dq_\gamma)^{-2}$ for $|q_\gamma|\le 1/2$. 
\begin{equation}\label{eq:defnormgamma}
|\omega|_\gamma=\max_{|q_\gamma|\le 1/2}
\left| \frac{q_\gamma\, \omega}{(dq_\gamma)^2}\right|.
\end{equation}
Any two such norms are of course equivalent: their ratios are bounded by a constant. More interestingly, the logarithm
of this constant factor
is polynomial in the level $5l$ of $X$:  for any $\gamma$ in $\Xi$ and any 
$\omega$
in $\cH^2(\Cusps)$ we have

\begin{equation}\label{eq:equivnorm}
l^{-\Theta}\times |\omega|_\gamma  \le |\omega|_\infty\le \exp(l^\Theta) \times |\omega|_\gamma.
\end{equation}

These inequalities will be proven in Section~\ref{subsection:equivnorm}.

\subsection{Jacobian of weight $4$ cusp forms}

Remind that we  have constructed in Section~\ref{section:SQD} a basis 
$\cB_{\rm quad}$
for the space $\cH^2(\Cusps)$ of quadratic differential  forms.
If $\hgamma \in \Xi^{g_2}$  and $\bq = (q_k)_{1\le k\le g_2}$
we define the quadratic Jacobian
\begin{equation}\label{eq:defJquad}
\cJ_\hgamma^{\rm quad}(\bq)=\left| \frac{q_{\gamma_k}\omega}{(dq_{\gamma_k})^2}(q_k) \right|_{ 1\le k\le g_2, \, \omega\in \cB_{\rm quad}.}
\end{equation}
It is a consequence of Ramanujan's conjectures
proven by Deligne and the explicit formulae by Asai
for pseudo-eigenvalues that the series $\frac{q_\gamma \omega}{(d q_\gamma)^2}$
is 
of type 
\begin{equation}\label{eq:typedeli2}
(\Theta,\Theta).
\end{equation}
for every $\omega$ in $\cB_{\rm quad}$ and every $\gamma$ in $\Xi$.
So the 
series in the  $g_2$  variables $q_1$, $q_2$, \ldots, $q_{g_2}$
defined by Equation~(\ref{eq:defJquad}) is of type
\begin{equation}\label{eq:typejac2quad}
(\exp(\Theta l^4),\Theta\bungd).
\end{equation}
Further, if $\hgamma = (\gamma,
\ldots, \gamma) \in \Xi^{g_2}$ is the repetition 
of $g_2$ times the same $\gamma$, we write $\cJ_\gamma^{\rm quad}$
for $\cJ_{\hgamma}^{\rm quad}$ and 
we denote by $W_\gamma^{\rm quad}$
the corresponding  Wronskian. This is a series in one variable $q$
and it is of type

\begin{equation}\label{eq:typewron2quad}
(\exp(\Theta l^6),\Theta l^4).
\end{equation}

We need a similar  estimate to Equation~(\ref{eq:infJinf}) for these  quadratic Jacobians
$\cJ_\gamma^{\rm quad}$.
We first treat the case when $\gamma=\Id$ and we write
$\cJ_\infty^{\rm quad}$ (resp. $W_\infty^{\rm quad}$) 
for $\cJ_\Id^{\rm quad}$ (resp. $W_\Id^{\rm quad}$).
 We denote by  $$G_{g_2}=\prod_{1\le m\le g_2-1}m!$$ the denominator 
in the definition of
$W_\infty^{\rm quad}$.

The expansions in $q=q_\Id$
of the  $\frac{q\omega}{(dq)^2} $ for $\omega\in \cH^2(\Cusps)$
are entire series with algebraic integer coefficients; and they are permuted
by the absolute Galois group of $\QQ$. So the series 
$G_{g_2}^2\times W_\infty^{\rm quad}(q)^2$ 
has coefficients in $\ZZ$.
 We set  $\eta_2=\frac{g_2(g_2+3)}{2}$.   
The product  $W_\infty^{\rm quad} (q)q^{-g_2}(dq)^{\eta_2}$ 
is a degree $\eta_2$  
form on $X$ and it is holomorphic outside  $\Cusps$. 
More precisely, it belongs to $\cH^{\eta_2}(g_2\Cusps)$.
Therefore it has  $$2(g-1)\eta_2+g_2\deg(\Cusps)\le \Theta l^6$$
 zeros counting multiplicities.
We deduce that the $q$-valuation $v_2$
of $W_\infty^{\rm quad} (q)$
is $\le \Theta l^6$.
So the series $G_{g_2}^2\times W_\infty^{\rm quad}(q)^2$ has valuation 
$\le \Theta l^6$
and rational integer coefficients. We deduce from 
Equations~(\ref{eq:typewron2quad})
and (\ref{eq:typprod}) that the type of the series $W_\infty^{\rm quad}(q)^2$
is 
\begin{equation}\label{eq:typeW2Q}
(\exp(\Theta l^6),\Theta l^4).
\end{equation}
We write
$$W_\infty ^{\rm quad}(q)^2=\frac{c}{G_{g_2}^2}q^{2v_2}+R_{2v_2+1}(q)$$ where 
$c$ is a non-zero rational
integer and $R_{2v_2+1}$ is the remainder of order $2v_2+1$. We can bound this
remainder using Equations~(\ref{eq:majreste1}) and (\ref{eq:typeW2Q}).
$$\left| R_{2v_2+1}(q)\right| \le \exp(\Theta l^6)|q|^{2v_2+1}.$$
So if $|q|\le \exp(-\Theta l^6)$ we have
  $|W_\infty^{\rm quad} (q)^2|\ge \frac{1}{2G_{g_2}^2}\left|q\right|^{2v_2}$  and
\begin{equation}\label{eq:infWinfquad}
|W_\infty^{\rm quad} (q)|\ge \exp(-\Theta l^5)\left|q\right|^{\Theta l ^6}.
\end{equation} So we fix such 
a $q$. For example we take 
\begin{equation}\label{eq:leqQ}
q=10^{-\kappa_3 l^6}
\end{equation}
where $\kappa_3$ is a large enough positive constant.
We set  
\begin{equation}\label{eq:lebqQ}
\bq=q\bungd =(q, \ldots, q)
\end{equation}
 and
 $\bbx = (x_1, \ldots, x_{g_2})$ where $x_1$, \ldots, $x_{g_2}$ are new 
indeterminates. 
The  Jacobian $$\cJ_\infty^{\rm quad}(\bq+\bbx\star(\bungd -\bB \bq \bB ))$$
is an entire series
in the  $g_2$ variables $x_1$, \ldots, $x_{g_2}$. 
This is indeed the Jacobian associated with the $g_2$
series  $$(f(q+x(1-|q|))/(q+x(1-|q|)))_{f(q)q^{-2}(dq)^2 \in \cB_{\rm quad}}$$ 
in the variable $x$. 
Equation~(\ref{eq:WJ}) gives us the first non-zero term in the expansion
of this series at $\bbx=\bzerogd$:

\begin{eqnarray*}
\cJ_\infty^{\rm quad}(\bq+\bbx \star (\bungd -\bB \bq \bB ) )=&&W_\infty^{\rm quad}(q)(1-|q|)^{\frac{g_2(g_2-1)}{2}}\prod_{k<l}(x_l-x_k)\\ &+&R_{\frac{g_2(g_2-1)}{2}+1}(\bbx).
\end{eqnarray*}

The type of the 
Jacobian  $\cJ_\infty^{\rm quad}(\bq)$  as  a series
in  $\bq$  is 
given by Equation~(\ref{eq:typejac2quad}).
We deduce from Equation~(\ref{eq:recen}) that the refocused series
$\cJ_\infty^{\rm quad}(\bq+\bbx\star(\bungd - \bB \bq \bB ))$ 
is a series in  $\bbx$ of type
 $(\exp(\Theta l^{4}), \Theta \bungd )$.
Using Equation~(\ref{eq:majreste1}) we deduce that
for  $\bbx $ in  $P(\bzerogd ,\exp(-\pi))$

$$\left| R_{\frac{g_2(g_2-1)}{2}+1}(\bbx)\right| \le \exp({\Theta}l^4)|\bbx|_\infty^{\frac{g_2(g_2-1)}{2}+1}.$$

We set  $s=|\bbx|_\infty$  and we  assume
that  $\bbx$ takes the special form

\begin{equation}\label{eq:lebxQ}
\bbx= (\frac{s}{g_2},\frac{2s}{g_2}, \ldots, \frac{(g_2-1)s}{g_2}, s)
\end{equation}
and  $s\le \exp(-\pi)$.
Then 

\begin{eqnarray*}
\scriptstyle \left|W_\infty^{\rm quad}(q)(1-|q|)^{\frac{g_2(g_2-1)}{2}}\prod_{k<l}(x_l-x_k)\right| & \scriptstyle \ge& \scriptstyle \left| W_\infty^{\rm quad}(q)\right|
\left(\frac{s(1-|q|)}{g_2}\right)^{\frac{g_2(g_2-1)}{2}}\\
&\scriptstyle \ge & \scriptstyle \Theta^{-1}\left| W_\infty^{\bf quad}(q)\right| \left(\frac{s}{g_2}\right)^{\frac{g_2(g_2-1)}{2}}
\end{eqnarray*}
We take 
\begin{equation}\label{eq:lesQ}
s=10^{-\kappa_4 l^{12}}
\end{equation}
\noindent where $\kappa_4$ is a large enough positive constant.
Using Equations~(\ref{eq:leqQ}) and (\ref{eq:infWinfquad}) we obtain the following lower bound for
the quadratic Jacobian
\begin{equation}\label{eq:uninfJQ}
|\cJ_\infty^{\rm quad}(\bq+\bbx  (1 -q  ) )|\ge \exp(-\Theta l^{16})
\end{equation}
when $q$, $\bq$, and $\bbx$ are given by Equations~(\ref{eq:leqQ}), (\ref{eq:lebqQ}), (\ref{eq:lebxQ}), and (\ref{eq:lesQ}). 
In particular $\left|\bq+\bbx  (1 -q  )\right|_\infty$
 can be assumed to be $\le 1/2$. So

\begin{equation}\label{eq:infJinfQ}
\max_{\bq \in D(0,1/2)^{g_2}} \left|\cJ_\infty^{\rm quad}(\bq) \right|
\ge \exp(-\Theta l^{16}).
\end{equation}

In order to  bound  from  below $\cJ_\gamma^{\rm quad}$ for any $\gamma \in \Xi$
we observe that if
$w$ is the width of the cusp $\gamma(\infty)$,  there
exists a $w$-th root of unity $\zeta_\gamma$ and an algebraic
number $\lambda_\gamma^{\bf quad}$ of absolute value $1$ (the product of all
pseudo-eigenvalues) such that  the following formal identity
in $\CC[[q_1, \ldots, q_g]]$ holds true
\begin{equation}%\label{eq:transJJQ}
\cJ_\gamma^{\rm quad}(q_1, \ldots, q_{g_2})=\lambda_\gamma^{\rm quad}
 \cJ_\infty^{\rm quad}(\zeta_\gamma q_1,
\ldots, \zeta_\gamma q_{g_2}).
\end{equation}

So the lower bound in Equation~(\ref{eq:infJinfQ}) is also 
valid for every $\cJ_\gamma^{\rm quad}$.

\subsection{Equivalence of norms on $\cH^2(\Cusps)$}\label{subsection:equivnorm}

In Section~\ref{section:SQD} we have defined various norms
on the space of quadratic differential forms $\cH^2(\Cusps)$.
If 
$$\omega = f(q)q^{-2}(dq)^2$$
 is a quadratic differential form 
in $\cH^2(\Cusps)$, the norm $|\omega|_\infty$ is the 
$L^\infty$ norm associated with  the basis $\cB_{\rm quad}$.
For every $\gamma \in \Xi$,  the norm $|\omega|_\gamma$ is defined
by Equation~(\ref{eq:defnormgamma}). We write
$$\omega = f(q)q^{-2}(dq)^2=h(q_\gamma)q_\gamma^{-2}(dq_\gamma)^2.$$
We must prove both inequalities in Equation~(\ref{eq:equivnorm}).
The first inequality is a trivial consequence of
Equations~(\ref{eq:typedeli2}) and (\ref{eq:majreste1}).

We denote by $\cM_\gamma^{\rm quad}$ the matrix occurring in the definition
of the jacobian $\cJ_\gamma^{\rm quad}$.
So 

\begin{eqnarray*}
\cM_\gamma^{\rm quad}(q_1, q_2, \ldots, q_{g_2})&=&\left( \frac{q_{\gamma}\omega}{(dq_{\gamma})^2}(q_k) \right)_{ 1\le k\le g_2, \, \omega\in \cB_{\rm quad}}
\\&=& \left( \frac{h(q_k)}{q_k} \right)_{ 1\le k\le g_2, \, h(q_\gamma)q_\gamma^{-2}(dq_\gamma)^2 \in \cB_{\rm quad}.}
\end{eqnarray*}

In particular $\cJ_\gamma^{\rm quad}$ is the determinant of 
$\cM_\gamma^{\rm quad}$. 

Now let $q$,  $\bq$ and $\bbx$  be given by Equations~(\ref{eq:leqQ}), (\ref{eq:lebqQ}), (\ref{eq:lebxQ}), and (\ref{eq:lesQ}).
We set
$$\br = \bq+\bbx (1-q) =(r_1, r_2, \ldots, r_{g_2})\in D(0,1/2)^{g_2}$$
and denote by $\cM_\gamma^{\rm quad}(\br)$
the evaluation of $\cM_\gamma^{\rm quad}$ at $\br$. 
The entries in  $\cM_\gamma^{\rm quad}(\br)$ are bounded
above by $l^\Theta$ in absolute value.
Using Equation~(\ref{eq:uninfJQ}) we deduce
that the entries in the inverse matrix of $\cM_\gamma^{\rm quad}(\br)$ 
are bounded in absolute value by $\exp(\Theta l^{16})$.

Let $\omega$ be a form in $\cH^2(\Cusps)$ and let
$\bc \in \CC^{g_2}$  be the coordinate vectors
of $\omega = h(q_\gamma)q_\gamma^{-2}(dq_\gamma)^2$ in the basis $\cB_{\rm quad}$.
For every $1\le k\le g_2$ set $v_k=h(r_k)/r_k$ and let
$\bv = (v_1, \ldots, v_{g_2})$ be the corresponding vector. We
have 
$$\bv^t = \cM_\gamma^{\rm quad}(\br) \times \bc^t$$
where $\bv^t$ is the transposed  vector of $\bv$ and
 $\bc^t$ is the transposed  vector of $\bc$.
So
$$|\omega|_\infty = |\bc|_\infty \le 
 g_2\exp(\Theta l^{16}) \times |\bv|_\infty$$
and this is $\le \exp(\Theta l^{16}) \times |\omega|_\gamma $
from  the definition of $|\omega|_\gamma $ given in 
Equation~(\ref{eq:defnormgamma}).

\section{An elementary operation in the jacobian}\label{section:opjac}

An important prerequisite for the  explicit computation
in the Jacobian $J$ is to  be able to compute the
linear space associated with some divisor on $X$.
In this section, we describe
an algorithm to solve the  following elementary problem:
given $3g-4$ points $P_1$, $P_2$, \ldots, $P_{3g-4}$ in
$X(\CC)$, find $g$ points $Q_1$, \ldots, $Q_g$ in $X(\CC)$ such
that 
$$Q_1+\cdots +Q_{g}\sim 2\cK - (P_1+\cdots +P_{3g-4})$$
where $\sim$ stands for linear equivalence of divisors and  
$\cK$\index{$\cK$, the canonical divisor}
is the canonical class.
This elementary problem will be used as a building block for 
explicit arithmetic operations in the jacobian $J$ of $X$.
We observe that the solution is not always unique.  However,
the image of $Q_1+\cdots +Q_{g}$ by the Jacobi integration map
$\phi' : \Sym^g X\rightarrow J$ is well defined. When doing numerical
approximations, it will be convenient to measure the error in 
$J(\CC)=\CC/\Lambda$ in terms of  the distance $d_J$ 
defined in Equation~(\ref{eq:distJ}).

We shall solve the above problem in two steps. 
We set $$P=P_1+\cdots +P_{3g-4}$$ and we
 first look for a differential quadratic form $\omega$
in the linear space $\cH^2(-P)\subset \cH^2(\Cusps)$ using our explicit knowledge
of the latter space and linear algebra algorithms.
We then compute the divisor $(\omega)$ of $\omega$ and output the (effective)  difference
$(\omega)-P$.
We now provide
 details for  these  two steps.

We denote by $T=P+\Cusps$ the divisor obtained by adding the cusps to
$P$. The degree of $T$ is $g_2-1$.
We write $T=T_1+T_2+\dots+T_{g_2-1}$ where
 $T_k=(\gamma_k,q_k)$ for every $1\le k\le g_2-1$. Let $\epsilon_1 =
\exp(-m_1)$  be a positive real number. We assume  that $m_1\ge \Theta l$.
We find an  $\epsilon_1$-simple
divisor 
 $T'=T'_1+T'_2+\dots+T'_{g_2-1}$  
such that for every $1\le k\le g_2-1$
we have $T'_k=(\gamma_k,q'_k)$
and $$|q'_k-q_k|\le \Theta (g_2-1)\epsilon_1.$$
We look for a quadratic form  in $\cH^2(\Cusps)$ having
divisor $\ge -\, \Cusps+T'$. The space of such forms
can be described as  the kernel
of a matrix $\cM$. Each of the  $g_2-1$ lines of $\cM$ 
corresponds to a point $T'_k$ for some $1\le k\le g_2-1$. The
$g_2$ columns of $\cM$ correspond to the $g_2$ forms
in the basis $\cB_{\rm quad}$. The  entry
of $\cM$ at the line corresponding
to the point $T'_k = (\gamma_k,q'_k)$ and column corresponding
to the form $\omega$ in $\cB_{\rm quad}$ is obtained in the following
way: we consider the expansion of $\omega$ in the variable $q_{\gamma_k}$
$$\omega=h(q_{\gamma_k})q_{\gamma_k}^{-2}(dq_{\gamma_k})^2$$
and we evaluate the function $h(q_{\gamma_k})q_{\gamma_k}^{-1}$
at the value $q'_k$ of $q_{\gamma_k}$ corresponding to $T'_k$.
The entries  in $\cM$ are bounded by $l^\Theta$ in absolute value according
to Equations~(\ref{eq:typedeli}) and  (\ref{eq:majreste1}).

We can't compute $\cM$ exactly. Instead of that, we fix a positive
real $\epsilon_2=\exp(-m_2)$ and we  compute a matrix $\cM'$ with 
decimal entries  in $\ZZ[i,1/10]$ such that 
the difference $\cM'-\cM$ has $L^\infty$ norm $\le \epsilon_2$.
The entries  in this matrix $\cM'$ can be chosen to have numerators
and denominators bounded in absolute value by $\exp(\Theta(l+m_2))$.
We find a non-zero vector in the kernel of $\cM'$, 
 having coefficients  in $\ZZ[i]$ and 
 bounded in absolute value by $\exp(\Theta l^2(l+m_2))$.
We divide this vector by its largest coefficient an obtain
a vector $v=(v_k)_{1\le k\le g_2}$ with
$L^\infty$ norm equal to $1$. This vector may not lie in the 
kernel of $\cM$ but $\cM v$  has coefficients 
 bounded by $g_2\exp(-m_2)$.
We call $\nu$ the quadratic differential form in $\cH^2(\Cusps)$
having  coordinate vector $v$ in the basis
$\cB_{\rm quad}$. By definition we have
\begin{equation}\label{eq:norm1}
\left| \nu \right|_\infty=1
\end{equation}
 Using 
Lemma~\ref{lemma:anyzero} of Chapter~\ref{sec_couveignes_ZEROS} together
with Equations~(\ref{eq:equivnorm}) and (\ref{eq:norm1}) we show that
if  $m_2\ge l^\Theta$, then 
for every 
 $1\le  k \le g_2-1$ the form $\nu$ has a
zero $T''_k=(\gamma_k,q''_k)$ such that 
$$|q'_k-q''_k|\le \exp(-\Theta^{-1} \sqrt{m_2}).$$
If $m_2\ge \Theta m_1^2$ then  these $g_2-1$ zeros must be pairwise distinct,  because each of them is close to some  $T'_k$ and
the latter points form an $\epsilon_1$-simple divisor.  The divisor of $\nu$ can be written 
$$(\nu)=T''-\Cusps+Q'$$
where $T''=\sum_{1\le k\le g_2-1}T''_k$ and $Q'$ is a degree $g$
effective divisor. We rewrite $(\nu)$ as
$$(\nu)=P+Q'+\sum_{1\le k\le g_2-1}\left( T''_k-T_k\right).$$

The image of the error term $\sum_{1\le k\le g_2-1}\left( T''_k-T_k\right)$
by $\phi'$ is small in the torus $J(\CC)=\CC/\cR$. More precisely

$$d_J(0,\phi'(\sum_{1\le k\le g_2-1}\left( T''_k-T_k\right)))
\le \exp(-m_1/\Theta)$$
provided $m_1\ge l^\Theta$ and $m_2\ge m_1^\Theta$.

So $Q'$ is a good approximation for the solution $Q$ to the original
problem. Using the algorithm in Theorem~\ref{theorem:findingzeros}
 of Chapter~\ref{sec_couveignes_ZEROS}, we compute an approximation
of the divisor of $\nu$ and output the corresponding approximation
$Q''$ of $Q'$.

\begin{lem}[An elementary operation]\label{lemma:flip}
There is  a deterministic algorithm that on input a degree $3g-4$ effective divisor 
$P=P_1+\cdots +P_{3g-4}$ on  $X_1(5l)$,  returns a degree $g$ effective divisor 
$Q=Q_1+\cdots+Q_{g}$ such that
$$Q_1+\cdots +Q_{g}\sim 2\cK - (P_1+\cdots +P_{3g-4})$$
where $\cK$ is the canonical class on $X_1(5l)$. The running time is 
$(lm)^\Theta$ where
$5l$ is the level and $m$ the required absolute   accuracy of the result.
\end{lem}

Remind that the accuracy of the result in the above
statement is measured in the torus $\CC/\Lambda$ using the distance
$d_J$ introduced in Equation~(\ref{eq:distJ}).

\section{Arithmetic operations in the Jacobian}\label{sec_torsion_arith}

We fix a degree $g$ effective 
divisor $\Omega$ on $X$.
 We also need an effective
 degree $g-4$ auxiliary divisor $\Pi$. For example, we may
choose a point  $O$ as origin for the Jacobi integration map (e.g. $O$ could be the
cusp at infinity), and set $\Omega=g O$ and $\Pi=(g-4)O$. An element in
 $\Pic^0(X)$ is given as the class of a divisor $Q-\Omega$ where $Q$
is a degree $g$  effective divisor.  Let   $R$ be another degree $g$ effective
divisor. 
In order to add the class of  $Q-\Omega$ and the class of
$R-\Omega$ we apply Lemma~\ref{lemma:flip}  twice.
We first apply it to the divisor  $Q+R+\Pi$. This is indeed a degree  $3g-4$ effective divisor.
We obtain a degree $g$ effective divisor  $T$ such that $T\sim 2\cK-Q-R-\Pi$.
We again apply Lemma~\ref{lemma:flip} to the divisor $T+\Omega+\Pi$ this time. And we obtain 
a degree $g$ effective divisor $U$ such that  $U+\Omega\sim Q+R$. So the class of  $U-\Omega$ is the sum
of the classes of $Q-\Omega$ and  $R-\Omega$.

In order to compute the opposite of the class $Q-\Omega$, we apply Lemma~\ref{lemma:flip} to the divisor
$2\Omega+\Pi$ and obtain a degree $g$ effective divisor  $R$, linearly  equivalent to  $2\cK-\Pi-2\Omega$.
We apply Lemma~\ref{lemma:flip} to the divisor  $R+Q+\Pi$ and obtain a degree $g$ effective divisor 
$T$ such that  $T-\Omega\sim -(Q-\Omega)$.

\begin{thm}[Arithmetic operations in $J_1(5l)$]\label{theorem:arithJ}
Addition and subtraction in the jacobian of  $X_1(5l)$ can be computed
 in deterministic time
$(lm)^\Theta$ where $5l$ is the level and $m$ the required
absolute  accuracy of the result.
\end{thm}

Again, the accuracy of the result is measured in the torus $\CC/\Lambda$ using the distance
$d_J$ introduced in Equation~(\ref{eq:distJ}). In particular the error belongs to a \underline{group},
and when  chaining operations in the jacobian, the successive errors  add to each other: the error on the result
is the sum of the errors 
on either input
 plus the error introduced in the current calculation. 
This observation
is particularly useful in conjunction with the fast exponentiation 
algorithm of Section~\ref{sec:comp}: if
we multiply a divisor $Q- \Omega$ by a positive integer $N$, assuming
that every elementary 
operation introduces and error $\le \epsilon$, then
the error on the final result is $\le \Theta \times \epsilon \times N\log N$ so
the loss of accuracy is $\le \Theta \log N$.

\begin{thm}[Fast exponentiation in $J_1(5l)$]\label{theorem:exp}
There is a deterministic  algorithm that on input two degree $g$ effective divisors
$\Omega$ and $Q$ on $X_1(5l)$ and a positive integer $N$ outputs
a degree $g$ effective divisor $R$ such that 
$$R-\Omega\sim N(Q-\Omega).$$

The algorithm runs in time $(lm\log N)^\Theta$ where $5l$ is the level and $m$ the required absolute 
accuracy of the result.
\end{thm}

\section{The inverse Jacobi problem}\label{sec_torsion_invjac}

In this section we are 
given a degree $g$ effective origin divisor $\Omega$ on $X$ and
an element $x$ in  $\CC^{\cB^1_{\rm DR}}$, and 
we want to solve the inverse Jacobi problem
for $$x+\Lambda \in \CC^{\cB^1_{\rm DR}}/\Lambda=J(\CC).$$
So we look for a degree $g$ effective divisor $P=P_1+\cdots+P_g$
on $X$ such that $\phi' (P - \Omega)=x+\Lambda$.
We note that the solution might not be unique. 

The main idea is the following: we start from a family of $2g$
classes $b_1+\Lambda$, $b_2+\Lambda$, \ldots, $b_{2g}+\Lambda$
in $J$ for which the inverse Jacobi problem is already solved: 
for every $1\le k\le 2g$, we know
a divisor $B_k-\Omega$ such that $\phi'(B_k-\Omega)=b_k+\Lambda$. We
try to approximate $x$ by an integer combination 
$\sum_{1\le k\le 2g} N_kb_k$. This should not be too difficult if
the $b_k$ are very small and $\RR$-linearly independent: 
we compute the coordinates
of $x\in \CC^g$ in the $\RR$-basis made of the $b_k$ and we round each of
these coordinates
to the closest integer.
Once we have found the $N_k$ we note that the divisor $\sum_kN_kB_k$
would be  a nice solution to the problem if it were a difference
between two effective degree $g$ divisors. This is not the case of course,
but using the algorithms in Theorems~\ref{theorem:exp} and 
\ref{theorem:arithJ} we find a degree $g$ effective  divisor $P$ such that
$P-\Omega$ is linearly  equivalent to $\sum_kN_kB_k$. We output $P$ and we 
are done.

There remains to explain how to find the $b_k$ and the corresponding $B_k$.
For every
$1\le k\le 2g$, we set $B_k=R'_k-R_k$ where $R_k$ and $R'_k$ are two points on $X$
that are very close. 
More precisely, we choose the $g$ first points $R_k=(\gamma_k,q_k)$  for $1\le k\le g$, and we 
set $R_{k+g}=R_k$. We also choose a  positive integer $\Upsilon=\exp(-\chi)$.
We assume  that   $\chi \ge l^\Theta$ so $\Upsilon$ is small. For $1\le k\le g$,  we set 
$R'_k=(\gamma_k,q_k+\Upsilon)$ and $R'_{k+g}=(\gamma_k, q_k+i\Upsilon)$. We set $$b_k=\left(\int_{R_k}^{R'_k}\omega\right)_{\omega \in {\cB^1_{\rm DR}}}\in \CC^{\cB^1_{\rm DR}}.$$
These integrals can be  computed efficiently using the same method as for period integrals. 
We now want to quantify the condition that these $b_k$ should be $\RR$-linearly independent.  So let $(e_k)_{1\le k\le g}$ be the canonical
basis of $\CC^g$ and set  $e_{k+g}=ie_k$ for $1\le k\le g$. So  $(e_k)_{1\le k\le 2g}$ is a  basis of the $\RR$-vector space underlying $\CC^g$.
We shall need a lower bound for the determinant 
\begin{equation}\label{eq:detbe}
\frac{b_1\wedge \ldots \wedge b_{2g}}{e_1\wedge \ldots \wedge e_{2g}}.
\end{equation}
 Since $\Upsilon$
is going to be small, we derive such a lower bound from the local study of the Jacobi integration map.
We call $\hgamma \in \Xi^g$ the vector $(\gamma_k)_{1\le k\le g}$ and we call $$\phi_\hgamma : D(0,1)^g\rightarrow \CC/\cR$$ the composition
of the Jacobi map $\phi : X^g \rightarrow J(\CC)$ with the product of
the modular parameterizations $$\mu_{\hgamma} = \prod_{1\le k\le g} \mu_{\gamma_k} :D(0,1)^g\rightarrow X^g.$$

For every $1\le k\le g$ we set $w_k=w_{\gamma_k}$.
Let  $\bq=(q_1, \ldots, q_g)$ with
$q_k\in D(0,\exp(-\pi/w_k))$ for every $1\le k\le g$.
We study the map $\phi_\hgamma$ locally at $\bq$.
The tangent space to $D(0,1)^g \subset \CC^g$ at $\bq$ is identified with $\CC^g$ and we denote
by $(\delta_k)_{1\le k\le g}$ its canonical basis. 
The underlying $\RR$-vector space has basis
$(\delta_1, \ldots, \delta_g, i\delta_1, \ldots, i\delta_g)$. For $1\le k\le g$ we set $\delta_{k+g}=i\delta_k$.
Let $D_\bq\phi_\hgamma$ be the differential
of $\phi_\hgamma$ at $\bq$.
For $1\le k\le 2g$ let $\rho_k$ be the image of $\delta_k$ by this  differential. Assuming $\chi\ge l^\Theta$ and using 
Equations~(\ref{eq:typedeli})  and   (\ref{eq:majreste1}) we 
prove that the coefficients of $\rho_k$
in the basis $(e_m)_{1\le m\le 2g}$ have absolute value
\begin{equation}\label{eq:boundgamma}
\le l^\Theta.
\end{equation}

The determinant  of $D_\bq\phi_\hgamma$ 
$$\frac{\rho_1\wedge \ldots \wedge \rho_{2g}}{e_1\wedge \ldots \wedge e_{2g}}$$
is the square of the absolute value of the Jacobian determinant

\begin{equation}\label{eq:jacdefi2}
\cJ_\hgamma(\bq)=\left| \frac{\omega}{dq_{\gamma_k}}(q_k) \right|_{1\le k\le g, \, \omega\in \cB^1_{\rm DR}.}
\end{equation}

We  denote $\lambda$ the opposite of the logarithm of the absolute value of the above determinant, and we call it
the {\it illconditioning} of $\bq$. We shall see that the inverse Jacobi problem is well conditioned  unless the illconditioning
is large.

We first observe that we can  bound from below the norm of every $\rho_k$
in terms of $\lambda$. Indeed, the determinant of 
$D_\bq\phi_\hgamma$  is  bounded from above  by the product of
the $L^2$ norms $\prod_{1\le k\le 2g} |\rho_k|_2$ so

$$\exp(-2\lambda)\le |\rho_k|_2 \times  \prod_{j\not =k} |\rho_j|_2
\le |\rho_k|_2 \times \exp(\Theta l^3)$$
using Equation~(\ref{eq:boundgamma}). 
So 
\begin{equation}\label{eq:mingamma}
|\rho_k|_2\ge \exp(-2\lambda-\Theta l^3).
\end{equation}
Our next 
 concern is to  bound from below  the determinant of Equation~(\ref{eq:detbe}) in terms of the illconditioning
$\lambda$. For every $1\le k\le 2g$ we notice that $\Upsilon \rho_k$ is the first order approximation
of $b_k$. We deduce from Equations~(\ref{eq:typedeli})  and   (\ref{eq:majreste1}) that 
$$\left|  b_k-\Upsilon \rho_k \right|  \le l^\Theta\Upsilon^2.$$ Using the lower bound (\ref{eq:mingamma}) we deduce that
$$\frac{\Upsilon}{2} |\rho_k|_2\le  |b_k|_2\le \frac{3\Upsilon}{2}|\rho_k|_2\le l^{\Theta}\Upsilon $$
provided $\chi \ge l^\Theta+2\lambda$.

Using multilinearity of the  determinant  we can bound the difference 
$$ \frac{b_1\wedge \ldots \wedge b_{2g}}{e_1\wedge \ldots \wedge e_{2g}} - \Upsilon^{2g}\frac{\rho_1\wedge \ldots \wedge \rho_{2g}}{e_1\wedge \ldots \wedge e_{2g}}$$
by 
$$2^g\left(\max_{1\le k\le 2g}|b_k|_2\right)^{2g-1} \max_{1\le k\le 2g}|b_k-\Upsilon\rho_k|_2\le 2^gl^{\Theta(2g-1)}l^{\Theta}\Upsilon^{2g+1}$$
and this is less than half of $\Upsilon^{2g}|\frac{\rho_1 \wedge \cdots \wedge \rho_{2g}}{e_1\wedge \cdots \wedge e_{2g}}|=
\exp(-2\lambda)\Upsilon^{2g}$ as soon as  $$\chi \ge  \Theta l^3 +2\lambda.$$

We deduce that
\begin{equation*}
\left| \frac{b_1\wedge \cdots \wedge b_{2g}}{e_1\wedge \cdots \wedge e_{2g}} \right| \ge \frac{1}{2}\exp(-2\lambda-2g\chi).
\end{equation*}

So we have a lower bound for the determinant
of the transition matrix between the basis
$(e_k)_{1\le k \le 2g}$  and the basis  $(b_k)_{1\le k \le 2g}$. 
Further the entries in this matrix are
bounded by $l^{\Theta} \Upsilon$ in absolute value.
Therefore the entries in the inverse  matrix are bounded 
by 

$$2\exp(2\lambda+2g\chi)l^{\Theta g}\Upsilon^{2g-1}\le \exp(2\lambda+\chi+ \Theta l^3)$$ in absolute value. We thus can compute this inverse matrix in time polynomial
in $l$, $\lambda$, $\chi$ and the required absolute  accuracy.

In Section~\ref{section:latticeperiods} we have constructed  a basis $\cB_{\rm per}$ for the lattice of periods,
consisting of vectors with coordinates bounded by $\exp(l^{\Theta})$ in absolute value  in the basis $(e_k)_{1\le k \le 2g}$.
The coordinates of these periods in the basis
$(b_k)_{1\le k\le 2g}$ are bounded by  $$\exp(l^{\Theta}+2\lambda+\chi)$$
in absolute value.

Every point in the fundamental parallelogram associated with the basis $\cB_{\rm per}$ (i.e. having coordinates in $[0,1]$ in this basis)
has coordinates  $\le \exp(2\lambda+\chi+l^{\Theta} )$ in absolute value in the basis $(b_k)_{1\le k \le 2g}$. When we replace the latter coordinates
by the closest integer, the induced error is bounded by  $l^{\Theta}\Upsilon$ 
according to the $L^2$ norm for the canonical basis $(e_k)_{1\le k \le 2g}$.

According to
Equation~(\ref{eq:uninfJ})
 there exists
 a vector $\bq$ with illconditioning $$\lambda \le l^\Theta.$$ This  
 finishes the proof of the following theorem.

\begin{thm}[Inverse Jacobi problem]\label{theorem:inverseJ}
The exists a  deterministic algorithm that takes as  input 

\begin{itemize}
\item a prime integer $l$,
\item an element $x$ in  the tangent space $\CC^{\cB^1_{\rm DR}}$ to $J_1(5l)$ at the origin (where $\cB^1_{\rm DR}$ is the basis
of $H^1_{\rm DR}$ made of normalized newforms of level $5l$ together
with normalized newforms of level $l$ lifted
to level $5l$ by the two degeneracy maps),
\item a degree $g$ effective (origin)  divisor $\Omega$ on $X_1(5l)$,
\end{itemize}
and returns an approximation of the
 degree $g$ effective divisor $$P=P_1+\cdots+P_g$$
on $X_1(5l)$ such that $\phi' (P - \Omega)=x+\Lambda$.

The running time is $\left( l\times 
\log \left( 2+\left| x \right|_\infty\right) \times m \right)^\Theta$ where $5l$ is the level, $m$ is the required absolute  accuracy of the result
and $\log \left(2+\left| x \right|_\infty\right)$ 
is the size of $x$ i.e.
the logarithm of its  $L^\infty$ norm 
in the canonical  basis of $\CC^{\cB^1_{\rm DR}}$.

\end{thm}

We insist  that 
the absolute accuracy in the above theorem is measured in the space $\CC^{\cB^1_{\rm DR}}$ using the $L^\infty$
norm; or equivalently in the jacobian  $J=J_1(5l)$  using the distance $d_J$.
However, when $x$ belongs to the Ramanujan subspace $W_f$
and assuming the degree $g$  origin divisor $\Omega$ is the 
 cuspidal
divisor $D_0$ manufactured  in Section~\ref{sec_constr_D}, then there
is a unique effective degree $g$ divisor $Q$ such that
$\phi'(Q-\Omega)=x+\Lambda$.  We call it
the {\it Ramanujan divisor} associated with $x$.
In that case, we can and must control
the error in $X^g$.
This error can be expressed as the distance between
$Q=\sum_{1\le k\le g}Q_k$ and
the output divisor
$Q'=\sum_{1\le k\le g}Q'_k$. Assume that $Q_k=(\gamma_k,q_k)$.
The   distance between
$Q$ and $Q'$ is  defined to be the minimum
over all permutations
$\sigma$ of $\{1,2, \ldots, g\}$ of the quantity
$$\max_{1\le k\le g} \left|q_k-q_{\gamma_k}(Q'_{\sigma(k)})\right|.$$

Sections~\ref{section:stabalg} to \ref{section:BXg}
will be mainly devoted to the proof  of the theorem below.

\begin{thm}[Approximating $V_f$ over the complex numbers]\label{theorem:torsion_main} 
There exists a  deterministic algorithm that takes as  input 
an even integer $k > 2$,  a prime integer $l>6(k-1)$, a finite field $\FF$
with characteristic $l$,  a ring epimorphism 
$f : \TT(1,k)\rightarrow  \FF$,
and a cuspidal divisor $\Omega$ on $X_1(5l)$
 like the divisor $D_0$  constructed in Section \ref{sec_constr_D},
and computes complex approximations for every element in
$W_f\subset J_1(5l)$,  the image of $V_f\subset J_1(l)$ by
$B_{5l,l,1}^*$. Here $V_f\subset J_1(l)$ is defined by 
Equation~(\ref{eq:defVltor}) and we  assume that 
the image of the Galois representation $\rho_f$ associated with $f$
contains $\SL  (V_f)$. The algorithm returns 
 for every element  $x$ in $W_f$ a complex approximation
of the unique   degree $g$   divisor  $Q_x=\sum_{1\le n\le g}Q_{x,n}$ such that 
$Q_x-\Omega$ lies in the class represented by $x$. Every point $Q_{x,n}$ is given as 
a couple $(\gamma,q)$ where $\gamma\in \Xi$ and $q$ is an 
approximation
of the value at $Q_{x,n}$ 
of the local analytic parameter $q_\gamma$ defined by 
Equation~(\ref{eq:qgamma}).
The running time of the algorithm is $\le (m\times \# V_f)^\Theta$ for some absolute constant $\Theta$.
Here $\# V_f$ is the cardinality of the Galois representation $V_f$ and $m$ is the required absolute accuracy 
of the returned approximations  for the $q$ associated with every $Q_{x,n}$.
\end{thm}

There are 
two important differences between Theorem~\ref{theorem:inverseJ}
and Theorem~\ref{theorem:torsion_main}. While Theorem~\ref{theorem:inverseJ} controls the error in $J_1(5l)$, Theorem~\ref{theorem:torsion_main} controls the error in $X_1(5l)^g$. Unfortunately, Theorem~\ref{theorem:torsion_main}
only applies to special divisors like the $Q_x-\Omega$. For these divisors,  one  can prove
that  the inverse
Jacobi problem is reasonably well conditioned. Results in  
Section~\ref{sec_constr_D} prove that
$Q_x$ is well defined;  and using 
Proposition~\ref{prop_arakelov_contrib_1} in Section~\ref{sec_final_estimates} one can show that computing $Q_x$ from $x$ is a well conditioned problem.
This will be the purpose of the next three sections.

Another remark concerning notation. We denote by $\Omega$ the divisor $D_0$ introduced in 
section \ref{sec_constr_D}. And we write $Q_x$ rather than
 $D_x$. The only reason for this slight change in notation
is that many things are already called $D$ in this chapter, and we want to avoid any possible confusion.

\section{The algebraic conditioning}\label{section:stabalg}

An important feature of Theorems~(\ref{theorem:arithJ}),
 (\ref{theorem:exp}), and  (\ref{theorem:inverseJ}) is that the  error
in all these statements is measured in the torus $J(\CC)$. We have seen
that this helps controlling the accumulation of errors when we chain
computations. However, when solving the inverse Jacobi
problem, we want to control the error
in $\Sym^g X$,  at least for the final result of the  computation. 
So we need a theoretical estimate for  the
 error in 
$\Sym^g X$ in terms of the error in  $J(\CC)$. This will be the main concern
of this and the following two   sections.

So assume that we are given a degree $g$ effective origin divisor 
$\Omega$ on $X$ and
a vector  $x$ in  $\CC^{\cB^1_{\rm DR}}$, and 
we look for an effective degree 
$g$ divisor $P$ such that 
$\phi' (P - \Omega)=x+\Lambda$, where $\Lambda$ is the lattice
of periods. We write 
\begin{equation}\label{eq:defP}
P=n_1P_1+n_2P_2+\cdots+n_KP_K
\end{equation}
where the $(P_k)_{1\le k\le K}$ are pairwise distinct 
points on $X$ and the 
$(n_k)_{1\le k\le K}$ are positive integers such that 
$$n_1+n_2+\cdots+n_K=g.$$

We shall make two assumptions.

The {\it first assumption} is rather essential: we  assume that the divisor
$P$ is {\it non-special} or equivalently that :
\begin{equation}\label{eq:nonspecial}
\dim(\Lambda(P))=1.
\end{equation}
 A first interresting consequence of this
first assumption is that the answer to the inverse Jacobi
problem is unique; and the error can be defined as the distance
to the unique solution.

The {\it second assumption} we make is more technical. We assume that
Klein's
 modular fonction $\KJ$ does not take the values $0$ or $1728$ at any
of the points $(P_k)_{1\le k\le K}$. So
$$\KJ (P_k)\not \in \{0,1728\}.$$

Removing
the second assumption would only result in an heavier presentation.
By contrast, the first assumption plays
 a crucial role in the forthcoming calculations.
Its meaning  is that the Jacobi map
$\phi' : \Sym^g X \rightarrow J$ is a local diffeomorphism
at $P-\Omega$. 
Our first task is to reformulate this first assumption in a more algebraic
setting. We look for an algebraic variant of the
Jacobi determinant in Equation~(\ref{eq:jacdefi2}).
Knowing that such an algebraic quantity is non-zero we will  bound 
it from below in Section~\ref{sec_torsion_heights}.

We assume  that we are given a basis $\cB^1_{\ZZ}$\index{$\cB^1_{\ZZ}$, a basis for $\cH^1=H^1_{\rm DR}(X_1(5l))$ consisting of forms with rational integer
coefficients in their $q$-expansion} of
$H^1_{\rm DR}(X_1(5l))$ such that for every $\omega = f(q) q^{-1}dq$
in  $\cB^1_{\ZZ}$, the associated modular form $f(q)$ has rational
integer coefficients and type $(\exp(\Theta l^4),\Theta)$. The
existence of such a basis is granted  by Lemma~\ref{lemma:baserat} below. 
We stress that we don't try   to compute such a basis. We are just
happy to know that it exists.

We also need an algebraic  uniformizing parameter 
$t_Q$ at every
point $Q$ on $X$ such that $\KJ (Q)\not \in \{0,1728\}$. 
When $Q$ is not a cusp either,  the differential 
$d\KJ$ of Klein's function $\KJ$ has no pole
nor  zero  at $Q$.
So $\KJ-\KJ(Q)$
is a uniformizing parameter
at $Q$. So we set $t_Q=\KJ-\KJ(Q)$ in that case.
When $Q=\gamma(\infty)$ is a cusp and
 $\gamma\in \Xi$,  we set $$\KJ_\gamma=\KJ\circ W^{-1}_\gamma$$
where $W_\gamma$ is given in Equation~(\ref{eq:Wgamma}). We deduce
from Equation~(\ref{eq:Wgamma2}) that the automorphism
$W_\gamma^{-1}$ maps
the cusp $\gamma(\infty)$ to the cusp $\infty$. So $\KJ_\gamma$
has a simple pole at $\gamma(\infty)$ and $\KJ^{-1}_\gamma$ is a 
uniformizing parameter at $\gamma(\infty)$. So we set
$t_Q=\KJ^{-1}_\gamma$ in that case.
We notice that $\KJ_\gamma$
only depends on the width of $\gamma(\infty)$. We call
 $\KJ_1=\KJ$, $\KJ_5$, $\KJ_l$
and $\KJ_{5l}$ the four corresponding functions.
Alltogether,  $t_Q$
is one of the  following functions: $\KJ-\KJ(Q)$,
$1/\KJ$, $1/\KJ_5$, $1/\KJ_l$ or  $1/\KJ_{5l}$.

Let  $\omega$ be a form in $\cB^1_{\ZZ}$ 
and let $k$ be an integer such that 
$1\le k\le K$ where $K$ is the number of distinct
points in the
divisor $P$ of Equation
(\ref{eq:defP}). Let $t_k=t_{P_k}$ be the algebraic uniformizing
parameter at $P_k$ and consider the Taylor expansion of
$\omega/dt_k$ at $P_k$
\begin{equation}\label{eq:Taylor}
\scriptstyle \omega/dt_k=\mW_{k,0}^{\, \omega}+\mW_{k,1}^{\, \omega}\times \frac{t_k}{1!}+
\mW_{k,2}^{\, \omega} \times \frac{t_k^2}{2!}+\cdots+
\mW_{k,n_k-1}^{\, \omega} \times \frac{t_k^{n_k-1}}{(n_k-1)!}+
O(t_k^{n_k}).
\end{equation}

We only need the first $n_k$
terms in this expansion,  where
$n_k$ is the multiplicity of $P_k$ in the divisor $P$.
We form the matrix
\begin{equation}\label{eq:wronalg}
\cM_P^{\rm alg}=\left(\mW_{k,m}^{\, \omega}  \right)_{1\le k\le K,\, 
0\le m\le n_k-1 ;  \, \omega\in \cB^1_{\ZZ}.}
\end{equation}

The lines in $\cM_P^{\rm alg}$ are indexed by  pairs $(k,m)$
where $1\le k\le K$ and $0\le m\le n_k-1$.  The columns in
$\cM_P^{\rm alg}$ are indexed by  forms $\omega$  in
$\cB^1_{\ZZ}$.
We call $\cM_P^{\rm alg}$ the {\it algebraic
Wronskian matrix }\index{$\cM_P^{\rm alg}$, the algebraic Wronskian matrix}
at $P$. 
%This matrix is something in between the usual
%Wronskian and the Jacobian matrix. For example, if 
%the divisor
%$P$ is simple, then $K=g$ and $n_1=n_2=\cdots=n_g=1$ and
%$\cM_P^{\rm alg}$ is a Jacobian matrix. On the other hand, if 
%$K=1$ and $n_1=g$, then $P=gP_1$ and $\cM_P^{\rm alg}$ is a
%usual Wronskian matrix. 
The determinant of $\cM_P^{\rm alg}$
will play an important role in the sequel. We call it the
{\it algebraic conditioning}\index{Algebraic conditioning}.

An important feature of the matrix $\cM_P^{\rm alg}$
is that its entries
are algebraic functions evaluated at the points $P_k$ for
$1\le k\le K$. Indeed let $t_k=t_{P_k}\in \{\KJ-\KJ(P_k), 1/\KJ, 1/\KJ_5, 1/\KJ_l, 
 1/\KJ_{5l}\}$ be the chosen algebraic uniformizing  
parameter  at the point $P_k$ and set
$\omega^{(0)}=\omega/d t_k$. For every integer $m\ge 0$
set  $\omega^{(m+1)}=d\omega^{(m)}/dt_k$. Then
$$\mW_{k,m}^{\, \omega} =\omega^{(m)}(P_k).$$

Lemma~\ref{lemma:algdep} below provides more detailed information
about the algebraic dependency between the derivatives
$\omega^{(m)}$ and Klein's function $\KJ$.

It is an important 
 consequence of our first assumption in Equation~(\ref{eq:nonspecial})
that the determinant of the algebraic
Wronskian matrix is non-zero
$$\det \cM_{P}^{\rm alg}\not = 0.$$

In the next Section~\ref{sec_torsion_heights} we shall
derive a lower bound for this determinant, using the theory
of heights.
To finish this section, there remains to state and prove the two
Lemmas~\ref{lemma:baserat} and \ref{lemma:algdep}. We first construct
the basis $\cB^1_{\ZZ}$.

\begin{lem}[A rational basis]\label{lemma:baserat}
There is  a basis $\cB^1_{\ZZ}$ of
$H^1_{\rm DR}(X_1(5l))$ such that for every $\omega = f(q) q^{-1}dq$
in  $\cB^1_{\ZZ}$, the associated modular form $f(q)$ has rational
integer coefficients and type $$(\exp(\Theta l^4),\Theta).$$
The transition matrix from the basis $\cB^1_{\rm DR}$ to
the basis  $\cB^1_{\ZZ}$  has algebraic integer entries bounded
in absolute value by $\exp(\Theta l^4)$ and its determinant is the square root
of a non-zero rational integer.
\end{lem}

We construct $\cB^1_\ZZ$ from 
$\cB^1_{\rm DR}$ using a descent process.

Let $\omega $ be a differential form in the basis
$\cB^1_{\rm DR}$ and let $$f(q) = \omega q(dq)^{-1}=\sum_{j\ge 1}f_jq^j$$
 be the corresponding  modular form. Let $\ZZ_f$  be the ring generated
by the coefficients of $f(q)$. As a $\ZZ$-module, $\ZZ_f$ is generated by 
the $f_j$ for $j\le  4l^2$. Let $\KK_f$ be the fraction
field of $\ZZ_f$. Let $\LL$ be a strict subfield of $\KK_f$.
Let $\ba=(a_j)_{1\le j \le 4l^2}$ be a vector
with  rational
integer coefficients.  The associated linear combination
$\sum_{1\le j\le 4l^2}a_jf_j$ belongs to $\LL$ if and only if 
$\ba$ belongs to a  submodule  %$\SSS_\LL $ 
of $\ZZ^{4l^2}$
with rank $<4l^2$. The degree $d_f$ of $\KK_f$ over $\QQ$ is bounded above
by $2g$. So the number of strict subfields of $\KK_f$
is $ <  2^{2g}$. So there exists rational integers
$(a_j)_{1\le j\le 4l^2}$ such that $0\le a_j<2^{2g}$ and 
$\theta = \sum_{1\le j\le 4l^2}a_jf_j$ does not belong to any strict 
subfield of
$\KK_f$. This $\theta$ is an algebraic integer that generates
$\KK_f$ over $\QQ$. For $0\le k\le d_f-1$ we set 
$$\Tr(\theta^kf)=\sum_{j\ge 1}\Tr(\theta^kf_j)q^j$$
where $\Tr : \ZZ_f \rightarrow \ZZ$ is the trace map.
Since $|f_j|\le \Theta (j+1)^\Theta$ and $|\theta|\le \exp(\Theta l^2)$ we deduce
that $|\Tr(\theta^kf_j)|\le \exp(\Theta l^4)(j+1)^\Theta$.
So the series $\Tr(\theta^kf)$ for
$0\le k\le d_f-1$ have type $(\exp(\Theta l^4),\Theta)$.

We do the same construction for every Galois orbit in $\cB^1_{\rm DR}$.
We collect all the forms  thus obtained. This makes
a basis $\cB^1_{\ZZ}$ of $H^1_{\rm DR}(X_1(5l))$ consisting of
forms $f(q)q^{-1}dq$ where $f(q)$ is a series with integer coefficients
and of type $(\exp(\Theta l^4),\Theta)$. \hfill $\Box$

\bigskip 

Now let's prove some quantitative statement about the algebraic
dependency
between  the successive 
 derivatives $\omega^{(m)}$ and $\KJ$.

\begin{lem}[An algebraic relation]\label{lemma:algdep}
Let $t$ be one of the  functions $\KJ, 1/\KJ, 1/\KJ_5, 1/\KJ_l, 
 1/\KJ_{5l}$. Let $\omega$ be a form in $\cB^1_{\ZZ}$ and let
$m\ge 0$ be an integer. Set $\omega^{(0)} =\omega/dt$
and $\omega^{(m)}=d^m\omega^{(0)}/(dt)^m$.
There exists a non-zero irreducible
polynomial $E(X,Y)\in\ZZ[X,Y]$ such that
$$E(\omega^{(m)}, t )=0.$$
Its degree
in either variable is 
$$\le (lm)^\Theta,$$ and 
its coefficients are bounded in absolute value by $$\exp((lm)^\Theta).$$
\end{lem}

\bigskip 

All the functions involved belong to the field $\QQ(X_1(5l))$ of
 modular functions having Puiseux expansion in $\QQ\{\{q\}\}$
at the cusp $\infty$.
This field is a regular extension of $\QQ(\KJ)$ corresponding to
the standard model of $X_1(5l)$ over $\QQ$.

We assume  that $t=\KJ$ since the other cases are quite similar.
The differential
$d \KJ$ has $4(l^2-1)$ zeros of multiplicity $2$ (the points
in the fiber of $\KJ$ above $0$) and
$6(l^2-1)$ zeros of multiplicity $1$ 
(the points
in the fiber of $\KJ$ above $1728$).
 So $\omega^{(0)}=\omega/d\KJ$
has degree $\le 14l^2$ and less than $10l^2$ poles.
When we differentiate $\omega^{(0)}$ we increase by one the multiplicity
of each pole. So $d\omega^{(0)}$ has the same poles as
$\omega^{(0)}$ and the total multiplicity of these poles
is less
than $24l^2$. When we divide by $d\KJ$ we don't add poles but we
increase the multiplicities by $1$ (in the fiber of $\KJ$
above $1728$) or $2$ (in the fiber of $\KJ$
above $0$). The degree of the polar divisor of $\omega^{(1)}$
is thus
 $\le 38 l^2$. We go on like that and we prove that 
the  degree of $\omega^{(m)}$ is $\le (14+24m)l^2\le 24(m+1)l^2$.
The degree of $\KJ$ is $12(l^2-1)$. So there is an irreducible
polynomial $E(x,y)$ in $\ZZ[x,y]$ such that
$E(\omega^{(m)},\KJ^{-1})=0$ and   $\deg_x E\le  12(l^2-1)$ and
$\deg_y E\le 24(m+1)l^2$.

In order to bound the coefficients in $E(x,y)$ we consider the
expansions of $\KJ^{-1}$ and $\omega^{(m)}$ at the cusp $\infty$.
Remember that $$\KJ  (q)=\frac{1}{q}+744+\sum_{k\ge 1}c(k)q^k$$
so 
$$-q^2\frac{d\KJ}{d q }=1-q\sum_{k\ge 1}kc(k)q^k.$$
From Equation~(\ref{eq:Petersson}) we deduce that
$\sum_{k\ge 1}kc(k)q^k$ has exp-type $(\Theta,0,2)$ in the sense
of Section~\ref{subsection:exptype}.
Using Equation~(\ref{eq:exptypeinv}) we deduce
that $\left(-q^2\frac{d\KJ}{d q }\right)^{-1}$ has exp-type
$(\kappa_1,0,4)$ where $\kappa_1\ge 1$ is an absolute  constant.
We write $$\omega^{(0)}=-q^2\frac{\omega}{dq}
\left( -q^2\frac{d\KJ}{dq}\right)^{-1}.$$
The series $\omega/dq$ has type $(\exp(\Theta l^4), \Theta)$
and exp-type $(\Theta,\Theta l^2,2)$. Using 
Equation~(\ref{eq:exptypeprod}) we deduce that
 $\omega^{(0)}$
has exp-type $(\kappa_2, l^2\kappa_3 ,4)$ for
some absolute constants $\kappa_2\ge 1$
and $\kappa_3\ge 1$. A simple iteration shows that $\omega^{(m)}$
has exp-type $$(\kappa_2+m\kappa_1, l^2\kappa_3+m\kappa_2+\frac{m(m-1)}{2}\kappa_1+2m,4).$$

So if $a$ is an integer such that $0\le a\le 12(l^2-1)$ then
$\left(\omega^{(m)}\right)^a$ has exp-type
$$(\Theta l^2(m+1), \Theta l^2(l+m+1)^2,4).$$

On the other hand, $\KJ^{-1}$ has exp-type $(\Theta,0,4)$ and
if $b$ is an integer such that $0\le b\le 24(m+1)l^2$ then
$\KJ^{-b}$ has exp-type
$$(\Theta(m+1)l^2,\Theta(m+1)l^2,4).$$

So all the monomials $\left(\omega^{(m)}\right)^a\KJ^{-b}$
arising in equation $E(\omega^{(m)}, \KJ^{-1})=0$ have
exp-type
$$(\Theta l^2(m+1), \Theta l^2(l+m+1)^2, 4)$$ and the coefficients
in their $q$-expansions up to order 
$$\deg(\omega^{(m)})\times \deg (\KJ^{-1})\le \Theta(m+1)l^4$$
are rational integers  bounded in absolute value by
$$\exp(\Theta l^{24}(m+1)^8).$$

Since the coefficients in $E(x,y)$ are solutions of the
 homogeneous
system given by these truncated $q$-expansions,
 they are bounded in absolute value by
$$\exp(\Theta l^{28}(m+1)^9).$$\hfill $\Box$

\section{Heights}\label{sec_torsion_heights}

In this section we recall basis facts about heights of algebraic
numbers and we deduce upper and lower bounds for the determinant
of the algebraic Wronskian matrix in Equation~(\ref{eq:wronalg})
when the divisor $P$ is a Ramanujan divisor.

Let $\Qb\subset \CC$ be the algebraic closure of $\QQ$
in $\CC$.
Let $\alpha \in \Qb$ be an algebraic number.
The {\it degree} $d_\alpha$
of $\alpha$
is the degree of the field extension $\QQ(\alpha)/\QQ$. 
Let $$f(x)=a_{d_\alpha} x^{d_\alpha}
+a_{d_\alpha-1}x^{d_\alpha -1}+\dots+a_0$$ be the unique
irreducible  polynomial in $\ZZ[x]$ such that $f(\alpha)=0$ and $a_{d_\alpha}
>0$.
We say that $a_{d_\alpha}$ is the {\it denominator} of $\alpha$
and we denote it $\dgot_\alpha$\index{$\dgot_\alpha$, the denominator of $\alpha$}.

Let $\KK$ be a number  field containing $\alpha$.
The multiplicative height of $\alpha$  with respect to $\KK$ is 
$$H_\KK(\alpha)=\prod_{\sigma}
\max(1,|\sigma(\alpha)|)\prod_{v}\max(1,|\alpha|_v)$$\index{$H_\KK(\alpha)$, the multiplicative height of $\alpha$}
where the $\sigma$ in the first product runs over the set of embeddings of $\KK$ into $\CC$ and 
the $v$ in the second product runs over the non-archimedean places of $\KK$.
In the special case $\KK=\QQ(\alpha)$ we have
$$H_{\QQ(\alpha)}(\alpha)=\dgot_\alpha \prod_{1\le k\le d_\alpha}
\max(1,|\alpha_k|)$$
where the $\alpha_k$ are the $d_\alpha$ roots of $f(x)$.

The {\it logarithmic height} of $\alpha$ with respect to $\KK$ is 
$$h_\KK(\alpha)=\log H_\KK(\alpha)$$\index{$h_\KK(\alpha)$, the logarithmic height of $\alpha$}
and the absolute (logarithmic) height of $\alpha$ is
$$h(\alpha)=\frac{h_\KK(\alpha)}{\deg(\KK/\QQ)}=\frac{h_{\QQ(\alpha)}(\alpha)}{d_\alpha}.$$\index{$h(\alpha)$, the absolute
logarithmic height of $\alpha$}
Knowing the degree $d_\alpha$ and absolute
height $h(\alpha)$ of a non-zero algebraic
number $\alpha$, we deduce the following upper and lower bounds
 for the absolute value $|\alpha|$
$$\exp\left(-d_\alpha \times h(\alpha)\right)\le |\alpha |
\le \exp \left(d_\alpha \times h(\alpha)\right).$$

Let $F(x)$ be a degree $d_F$ polynomial in $\ZZ[x]$
and assume that all coefficients in $F(x)$ are bounded by
$H_F$ in absolute value. Let $\KK$ be a number field and let $\alpha\in \KK$ be an algebraic number.
We set $\beta=F(\alpha)$. 

If $\sigma$ is any embedding of $\KK$ into  $\CC$ we 
have  $\sigma(\beta)=F(\sigma(\alpha))$ so
\begin{eqnarray*}
\max(1,|\sigma(\beta)|)&\le &{(d_F+1)} H_F \max(1,|\sigma(\alpha)|)^{d_F}.
\end{eqnarray*}

\medskip 

Now, let $v$ be a non-archimedean valuation of $\KK$. We have
\begin{eqnarray*}
\max(1,|\beta|_v)&\le &\max(1,|\alpha|_v)^{d_F}.
\end{eqnarray*}

Forming the product over all $\sigma$'s and all $v$'s
we find that the  absolute 
 logarithmic height of $\beta$ is 
\begin{equation}\label{eq:heightg(a)}
\le d_F\times h(\alpha) +\log(d_F+1)+\log H_F.
\end{equation}

Now assume that $\alpha$ and $\beta$ belong to a degree $d$ extension $\KK$
of
$\QQ$ and let $E(x,y)\in \ZZ[x,y]$ be a polynomial such that
$E(\alpha,y)\not =0$ and $E(\alpha,\beta)=0$. Assume that all coefficients in
$E(x,y)$ are bounded by $H_E$ in absolute value. Call $d_x$ (resp. $d_y$)
the degree of $E(x,y)$ with respect to the variable $x$ (resp. $y$).

We write $$E(x,y)=\sum_{0\le k\le d_y}E_k(x)y^k.$$
We deduce from
inequality (\ref{eq:heightg(a)}) 
that every  $E_k(\alpha)$ has absolute height

\begin{equation}\label{eq:hEk}
\le d_x\times h(\alpha)+\log (d_x+1)+\log H_E.
\end{equation}

Let $K$ be the largest $k$ such that 
$E_k(\alpha)\not =0$. 
Set $$F(y)=E(\alpha,y)=\sum_{0\le k\le K}E_k(\alpha)y^k.$$
If $\sigma$ is any embedding of $\KK$ into  $\CC$ we call 
$${}^\sigma\! F(y)=E(\sigma(\alpha),y)=
\sum_{0\le k\le K}E_k(\sigma(\alpha))y^k$$ the polynomial 
obtained by applying $\sigma$ to all the coefficients in $F(y)$.
Applying  Landau's inequality to ${}^\sigma\! F(y)$
we find that 
$$\max(1,|\sigma(\beta)|) \le  \frac{\sqrt{d_y+1}\times \max_{0\le k\le K}(|E_k(\sigma(\alpha))|)}{\left| E_K(\sigma(\alpha))\right|}$$
and this is 
$$\le 
{\sqrt{d_y+1}\times (d_x+1)H_E\times \max(1,|\sigma(\alpha)|)^{d_x}}\times {\max (1,\left| E_K(\sigma(\alpha))\right|^{-1})}.$$

\medskip 

Now, let $v$ be a non-archimedean valuation of $\KK$. Applying Gauss'
lemma to $F(y)$ we find that
\begin{eqnarray*}
\max(1,|\beta|_v)&\le &|E_K(\alpha)|_v^{-1}\times \max(1,|\alpha|_v)^{d_x}\\
&\le& 
 \max(1,|E_K(\alpha)|_v^{-1})\times \max(1,|\alpha|_v)^{d_x}.
\end{eqnarray*}

Forming the product over all $\sigma$'s and all $v$'s
we find that the  absolute 
 logarithmic height of $\beta$ is bounded from above by

$$\frac{\log(d_y+1)}{2}+\log(d_x+1)+\log H_E+d_x\times h(\alpha)
+h(E_K(\alpha)).$$

Using Equation~(\ref{eq:hEk}) we deduce that 

\begin{equation*}
h(\beta)\le \frac{\log(d_y+1)}{2}+2\log(d_x+1)+2\log H_E+2d_x\times h(\alpha).
\end{equation*}

\begin{lem}[Relating heights of algebraic numbers]\label{lemma:heightsrel}
Let $\alpha$ and $\beta$ belong to a degree $d$ extension of
$\QQ$ and let $E(x,y)\in \ZZ[x,y]$ be a polynomial such that
$E(\alpha,y)\not =0$ and $E(\alpha,\beta)=0$. Assume that all coefficients in
$E(x,y)$ are bounded by $H_E$ in absolute value. Call $d_x$ (resp. $d_y$)
the degree of $E(x,y)$ with respect to the variable $x$ (resp. $y$).
Then the absolute heights of $\alpha$ and $\beta$ are related by
the following inequality

$$h(\beta)\le 
\frac{\log(d_y+1)}{2}+2\log(d_x+1)+2\log H_E+2d_x\times h(\alpha).$$
\end{lem}

\bigskip 

We now can bound the heights of the entries in the algebraic Wronskian of
Section~\ref{section:stabalg}, at least in the cases we are interrested
in. Let 
$$x+\Lambda\in W_f\subset J(\CC)=\CC^{\cB^1_{\rm DR}}/\Lambda$$
be a non-zero  vector in the Ramanujan subspace $W_f$ 
of Theorem~\ref{theorem:torsion_main}.
Assume that 
the degree $g$  origin  divisor $\Omega$ on $X_1(5l)$ is the cuspidal
divisor $D_0$ manufactured  in Section~\ref{sec_constr_D}. 
Let 
$P$ be the unique    degree 
$g$ divisor  such that 
$$\phi' (P - \Omega)=x+\Lambda.$$

As proven
in Section~\ref{sec_constr_D}, the
  two assumptions of Section~\ref{section:stabalg} are satisfied in that
case:
the divisor
$P=n_1P_1+n_2P_2+\cdots+n_KP_K$ is {\it non-special} and 
Klein's
 modular fonction $\KJ$ does not take the values $0$ or $1728$ at any
of the points $(P_k)_{1\le k\le K}$. 
Proposition~\ref{prop_arakelov_contrib_1} in Section~\ref{sec_final_estimates}
implies that the absolute height of every $\KJ(P_k)$ is $\le l^\Theta$.
Using Lemma~\ref{lemma:algdep}
 and Lemma~\ref{lemma:heightsrel}  we deduce
that every entry in the algebraic  Wronskian matrix of Equation~(\ref{eq:wronalg}) has absolute  height $\le l^\Theta$.
Further any such entry generates an extension  of $\QQ$
of degree $\le (\# V_f)^\Theta$.
 So   every entry
in the algebraic Wronskian has 
denominator and absolute value bounded
by $\exp((\# V_f)^\Theta)$. So both the denominator and absolute value of the determinant
are $\le \exp((\# V_f)^\Theta)$. 
So 
\begin{equation}\label{eq:minalgstab}
h(\det (\cM_P^{\rm alg}))\le (\# V_f)^\Theta
\end{equation}
in that case.

Note also that,  as an algebraic
integer,  the algebraic conditioning as degree 
at most twice  $\# V_f$ because the square of it lies in the
definition field of $x+\Lambda$ and the latter field is a degree 
$\# V_f-1$ extension of $\QQ$. 

We deduce that for all primes $p$ but a finite number bounded by
$(\# V_f)^\Theta$, the divisor $P$ remains non-special when we reduce modulo
any place $\pgot$ above $p$.
Indeed, we can assume that $p\not \in \{5, l\}$
so $X_1(5l)$ has good reduction modulo
$p$. We also can assume that $\KJ(P_k)\not \in \{ 0, 1728\}
\bmod \pgot$ because both the degree and 
the absolute height of $\KJ(P_k)$ are  $\le (\# V_f)^\Theta$. We also
can assume that the $P_{k}$ remain pairwise distinct
modulo $\pgot$ for the same reason: we just need
to  exclude less than $(\# V_f)^\Theta$ primes  $p$.
We  can also  assume that every $P_k$ which is not a cusp does not
reduce modulo $\pgot$ onto a cusp.
We also assume that $p$ is larger  than  the genus
$g$ of $X_1(5l)$ so that the Taylor expansion in
Equation~(\ref{eq:Taylor}) remains valid modulo $\pgot$.  So the algebraic Wronskian matrix reduces
modulo $\pgot$ to the algebraic Wronskian matrix. Excluding 
a few more primes, but no more than $(\# V_f)^\Theta$, we can assume
that the algebraic  conditioning does not vanish modulo $\pgot$.
So the divisor $P$ remains non-special.

\begin{lem}[Reduction modulo $p$ of a Ramanujan divisor]\label{lem:jmcuni}
Call $\Omega$  the cuspidal divisor  $D_0$
on $X_1(5l)$ introduced
in Section~\ref{sec_constr_D} and let $W_f$ be the Ramanujan subspace $W_f$ 
of Theorem~\ref{theorem:torsion_main}. Let 
$$x\in W_f\subset J_1(5l)$$
and let $P$ be the unique divisor
on $X_1(5l)$ such that $P-\Omega$ lies in the class defined
by $x$. Then $P$ is non-special and for all primes $p$
but a finite number bounded by $(\# V_f)^\Theta$, the divisor
$P$  remains non-special modulo any place $\pgot$ above $p$.
\end{lem}

\medskip

We stress that the above lemma is very similar
to Theorem~\ref{thm_exist_good_places}.

\section{Bounding the error in $X^{\protect\lowercase{g}}$}\label{section:BXg}

In this section we fill the gap between Theorems~\ref{theorem:inverseJ} and~\ref{theorem:torsion_main}.
For $x+\Lambda \in W_f\subset J(\CC)$ and $P$ 
the corresponding Ramanujan divisor, we relate the error on $P$ and 
the error on $x$. We need some control on  the Jacobi integration map locally at $P$. 
A first step  in this direction is the upper bound for the height of the algebraic stability, given
in Equation~(\ref{eq:minalgstab}).
In Section~\ref{section:anacon} we deduce a lower bound for the determinant of the differential of the Jacobi
map at $P$. This bound implies that the Jacobi map is non-singular on a reasonably large 
neighborhood of $P$, as we show in Section~\ref{subsection:defons}.
In  Section~\ref{sec_torsion_relate}  we deduce that the inverse of the Jacobi map
is well defined and Lipschitz (with reasonably small constant) on a (reasonably large) neighborhood of $x+\Lambda$.
We shall need the
identities  between Wronskians and Jacobians of power series
proven in Section~\ref{section:JW}.
The reason for these algebraic complications is the following: we need an equation for the singular locus in $X^g$
of the Jacobi integration map. The space $X^g$ is stratified by the diagonals. The strata correspond to partitions of $\{1, 2, 3,  \ldots, g\}$.
We obtain a different equation for the singular locus on every stratum. These various equations
are related by algebraic identities between Jacobians and Wronskians.

\subsection{The analytic conditioning of  a divisor}\label{section:anacon}

Let 
\begin{equation}\label{eq:iciP}
P=n_1P_1+\dots+n_KP_K
\end{equation}
 be a degree $g$ effective divisor
on $X$ where the $n_k$ are positive integers for $1\le k\le K$.
For every $k$ 
we write $$P_k=(\gamma_k, q_k)$$ where $\gamma_k\in \Xi$ and
$q_k\in F_{w_k}\subset D(0,\exp(-\pi/w_k))$,
where  $w_k$ is the width of the cusp $\gamma_k(\infty)$.
Let 
$$\phi_{(n_k)_{1\le k\le K}} : \prod_{1\le k\le K}\Sym^{n_k}X
\rightarrow J(\CC)=\CC^{\cB^1_{\rm DR}}/\Lambda$$
be the relevant Jacobi integration map
in this context. We stress that this map is different from
the maps introduced before. Its initial set is a sort of
{\it semi-symmetric product} i.e. something between
the Cartesian product $X^g$ and the full symmetric product
$\Sym^gX$.
In order to write down the differential of this map at the divisor
$P$,  we need a local system of coordinates on 
$\prod_{1\le k\le K}\Sym^{n_k}X$. We shall use partial
Newton power sums. 
%\medskip 
Let $$R=(\gamma , q_\gamma(R))$$ be a point on $X$.
For every $n\ge 1$ and $m\ge 1$ we define the
$m$-th power sum $\nu_{R,n,m}$
on $\Sym^n X$ to be the function that takes the value
$$\left(q_\gamma(Q_1)-q_\gamma(R)\right)^m+
\left(q_\gamma(Q_2)-q_\gamma(R)\right)^m+
\dots+\left(q_\gamma(Q_n)-q_\gamma(R)\right)^m$$
at $\{Q_1, \ldots, Q_n\}$.

\medskip

Now let us come back to the  divisor $P$ in Equation~(\ref{eq:iciP}).
The functions $$(\nu_{P_k,n_k,m})_{1\le k\le K;\, 1\le m\le n_k}$$
form a local system of coordinates on 
$\prod_{1\le k\le K}\Sym^{n_k}X$
at the divisor $P$.
The matrix of the differential at $P$
of the integration map
$\phi_{(n_k)_{1\le k\le K}}$ in the bases
$(d\nu_{P_k,n_k,m})_{1\le k\le K;\, 1\le m\le n_k}$
and ${\cB^1_{\rm DR}}$ is

\begin{equation}\label{eq:mpana}
\cM_P^{ana}=\frac{1}{\prod_{1\le k\le K}\prod_{1\le m\le n_k}m!}\left( \hmW_{k,m}^{\, \omega} \right)_{1\le k\le K; \, 0\le m
\le n_k-1}
\end{equation}
where $\hmW_{k,m}^{\, \omega}$
is the $m$-th derivative of $\omega/d q_{\gamma_k}$ with respect
to $q_{\gamma_k}$, evaluated at $P_k=(\gamma_k ,q_k)$. 
The determinant of this matrix is called the {\it analytic conditioning}\index{Analytic conditioning}
of the divisor $P$.
This analytic conditioning and the algebraic conditioning defined in
Section~\ref{section:stabalg} after Equation~(\ref{eq:wronalg})  differ
by simple factors we should not be afraid of.

\medskip

Now, we assume that we are in the context of
Theorem~\ref{theorem:torsion_main}. We call
$\Omega$  the  cuspidal divisor $D_0$ constructed in Section~\ref{sec_constr_D} and we assume that the class of $P-\Omega$
corresponds to a point $x$ in the Ramanujan subspace $W_f\subset J$.
Equation~(\ref{eq:minalgstab}) implies that the algebraic conditioning of $P$
is $\ge \exp(-(\# V_f)^\Theta)$. In order  to obtain a similar 
lower bound for the analytic conditioning, we must  bound from below  the
complementary factors. 

The factor $1/\prod_{1\le k\le K}\prod_{1\le m\le n_k}m!$
 is $\ge \exp(-\Theta l^7)$. 

Lemma~\ref{lemma:baserat} implies that the factor coming from the change 
of bases from $\cB^1_{\rm DR}$
to $\cB^1_{\ZZ}$ is $\ge \exp(-l^\Theta)$.

There are also  factors  due to the 
change of coordinates.
Assume for example that $P_k$ is not a cusp.
So  the algebraic parameter  at $P_k$ is $\KJ-\KJ(P_k)$. The analytic 
parameter  at $P_k$ is $q_{\gamma_k} - q_k$. The extra factor in the
analytic conditioning is thus 
$$\left(\frac{d\KJ}{dq_{\gamma_k}}(P_k)\right)^{\frac{n_k(n_k+1)}{2}}$$
Lemma~\ref{lemma:djmin} below ensures that there exists
a constant $K_1$ such that
this factor is $\ge \exp(-(\# V_f)^{K_1})$ provided
there exists a constant $K_2$
such that  $\left|\KJ(P_k)\right|$ and 
  $\left|\KJ(P_k)-1728\right|$ are both $\ge \exp(-(\# V_f)^{K_2})$.
But this latter condition is met  because
$\KJ(P_k)$ is not $\{0,1728\}$, and both its degree and  logarithmic
height are $\le (\# V_f)^\Theta$.

Assume now that $P_k$ is a cusp. Then the algebraic parameter
at $P$ is the function $1/\KJ_{w_k}$ introduced in Section~\ref{section:stabalg} where $w_k\in \{1,5,l,5l\}$ is
the width of  the cusp $P_k=\gamma_k(\infty)$. And the
derivative $d(1/\KJ_{w_k})/dq_{\gamma_k}$ is just $1$.
So the analytic conditioning at a Ramanujan $l$-torsion
divisor  is
\begin{equation}\label{eq:anacon}
\ge \exp(-(\# V_f)^\Theta).
\end{equation}

\medskip

To finish this section, there remains to state and prove Lemma~\ref{lemma:djmin}.

\smallskip

\begin{lem}[Lower bounds for $d\bJ/dx$]\label{lemma:djmin}
For every  positive real number  $K_1$ there exists
a positive real number $K_2$ such that
the following statement is true.

Let $L\ge 2$ be an integer and let $\bJ(x)$ be Klein's
series given in Equation~(\ref{eq:744}). 
Let $q\in D(0,1)$
be a complex number in the unit disk such that 
$1-|q|\ge L^{-K_1}$ and  $|\bJ(q)|\ge \exp(-L^{K_1})$ and 
$|\bJ(q)-1728|\ge \exp(-L^{K_1})$.
Then $|\frac{d\bJ}{dx}(q)|\ge \exp(-L^{K_2})$.
\end{lem}

We first note that $\frac{d\bJ}{dx}(q)$
only vanishes if  $\bJ(q)\in \{0,1728\}$.
So we just want to prove that if $\frac{d\bJ}{dx}$ is small
at $q$ then $q$ is close to a zero of it.
We would like to apply Lemma~\ref{lemma:anyzero}
 to the series $-x^2\bJ'(x)=-x^2
\frac{d\bJ}{dx}$. But this series
is not of type $(A,n)$ because its coefficients are a bit
too large.
So we set $R=(1+|q|)/2$ and we set $\KJJ(x)=-x^2R^2\times \bJ'(Rx)$.
From the hypothesis in the lemma  there exists a  positive
real $K_3$ such that $\log R^n\le -n\times L^{-K_3}$. From 
Equation~(\ref{eq:Petersson}) there exists an absolute
 constant $\kappa$ such that
the coefficient of $x^n$ in $-x^2\times \bJ'(x)$  is
$\le \exp(\kappa \sqrt n)$. This implies that the series
$\KJJ(x)$ has type $(\exp(L^{K_4}), K_4)$ for some $K_4\ge1$.
We apply Lemma~\ref{lemma:anyzero} of Chapter~\ref{sec_couveignes_ZEROS} to the series $\KJJ(x)$ at
$q/R$ and we are done.
\hfill $\Box$

\subsection{The neighborhood of a non-special divisor}\label{subsection:defons}

The singular locus of the Jacobi integration map
$\phi : X^g \rightarrow J$ is a strict closed subset.
So every non-special effective degree $g$
divisor has a neighborhood in $\Sym^g X$
consisting of  non-special divisors. In this section, we 
provide a quantified version of this statement.
We first need to define simple neighborhoods of
the divisor $P$ in Equation~(\ref{eq:iciP}).
Let $\epsilon$ be a positive real number.  We call $\cP_\epsilon$
the set of degree $g$ effective divisors $P'$ that can be written
$P'=\sum_{1\le k\le K}\sum_{1\le m\le n_k}P'_{k,m}$
where $P'_{k,m}=(\gamma_k,q'_{k,m})$ and $|q'_{k,m}-q_k|\le \epsilon$.

We assume that $P$
is non-special. 
We expect that if $\epsilon$ is small enough, then $P'$
is  non-special as well. In order to write down the 
analytic conditioning of $P'$ we must take
 multiplicities into account. So we rewrite 
$P'$ as 
$$P'=\sum_{1\le k\le K}\sum_{1\le s\le S_k}m_{k,s}P'_{k,s},$$
where $\bm_k=(m_{k,s})_{1\le s\le S_k}$ is a partition
of $n_k$ into $S_k$ non-empty parts. In particular
$$n_k=m_{k,1}+m_{k,2}+\dots+m_{k,S_k}.$$

For every $1\le k\le K$ we call $$x_k=q_{\gamma_k}-q_k$$
the local analytic parameter at $P_k$. We call
$\bff_k$ the vector in $(\CC[[x_k]])^{\cB^1_{\rm DR}}$ defined by

\begin{equation}\label{eq:taylorfk}
\bff_k=
\left(\Taylor\left(\frac{\omega}{dq_{\gamma_k}},q_k\right)  \right)_{\omega\in \cB^1_{\rm DR}}
\in (\CC[[x_k]])^{\cB^1_{\rm DR}}
\end{equation}
where $\Taylor\left(\frac{\omega}{dq_{\gamma_k}},q_k\right)$
is the Taylor expansion of $\frac{\omega}{dq_{\gamma_k}}$
at $q_k$ in the parameter $x_k$.

For every $1\le k\le K$ and $1\le s\le S_k$ we write
$$P'_{k,s}=(\gamma_k,q'_{k,s})$$
and we set  $y_{k,s}=q'_{k,s}-q_k$.

The analytic conditioning of $P$ is the Wronskian
in Equation~(\ref{eq:totalW}) evaluated at $(0,0,\ldots,0)\in \CC^K$ and divided by 
$\bn != \prod_{1\le k\le K} n_k!$

The analytic conditioning of $P'$ is the Jacobian
in Equation~(\ref{eq:totaljacob}) divided by
$\bm !=\prod_{1\le k\le K}\prod_{1\le s\le S_k}m_{k,s}!$ More precisely, the
Jacobian in Equation~(\ref{eq:totaljacob}) is $\bm !$ times the Taylor expansion of
the analytic conditioning in the parameters $(y_{k,s})_{1\le k\le K; 1\le s
\le S_k}$.

These two conditionings are related by Equation~(\ref{eq:WJtot}). The conditioning at $P'$ is divisible
by the determinant $D_{\bm, \bn}$ relating the two partitions.
And the quotient specializes to the conditioning at $P$ times $\bn !/\bm !$

\smallskip

Let $S=\sum_{1\le k\le K} S_k$. Then 
$\cJ_{(\bff_k, \bm_k)_{1\le k\le K}}$ is a series
in the $S$ variables $(y_{k,s})_{1\le k\le K; \, 1\le s\le S_k}$.

\smallskip 

We assume that $q_k$ and $q'_{k,m}$ belong to
$D(0,\exp(-\pi/w_k))$.

\smallskip 

Using Equation~(\ref{eq:typedeli}) about
the size of coefficients in modular forms,
the type of refocused series as described 
in Section~\ref{section:refoc},
the type of derivatives given in Section~\ref{section:typederi}, 
the type of quotient series given in Section~\ref{section:typequo}, 
we can prove that the series
$$\frac{\cJ_{(\bff_k, \bm_k)_{1\le k\le K}}}{D_{\bm,\bn}}$$
has type $(\exp(l^\Theta),l^\Theta\bun )$.

We  denote by $\lambda$ the opposite of the logarithm of the absolute value of 
the analytic conditioning of $P$. 
We call it the  {\it analytic illconditioning} of $P$.
It generalizes the illconditioning introduced
in Section~\ref{sec_torsion_invjac}. The only difference is that
here we take multiplicities into account.
Using Equation~(\ref{eq:WJtot}) and inequality
(\ref{eq:majreste1}) on  bounding the remainder, we now
prove that if 
\begin{equation}\label{eq:nonspe}
-\log \epsilon \ge l^\Theta+\lambda
\end{equation}
 then
$P'$ is a non-special divisor also.  According
to Equation~(\ref{eq:anacon}) the analytic illconditioning
$\lambda$ is $\le (\# V_f)^\Theta$. So 
any divisor  in the neighborhood
$\cP_\epsilon$ is  non-special, provided 
$\epsilon \le \exp(-(\# V_f)^\Theta)$.

\subsection{Relating the direct and inverse error}\label{sec_torsion_relate}

We go on with the notation in the previous Section~\ref{subsection:defons}. We have a divisor $P$ as
in Equation~(\ref{eq:iciP}). 
We denote by $\cS_{(n_k)_{1\le k\le K}}$ or just $\cS$
the map $$\scriptstyle (R_1, \ldots,
R_g)\mapsto (\{R_{1}, \ldots, R_{n_1}\},
\{R_{{n_1+1}}, \ldots, R_{n_1+n_2}\}, \ldots, \{R_{
n_1+n_2+\dots+n_{K-1}+1}, \ldots, R_{g}\})$$
and we check that the following diagram commutes:

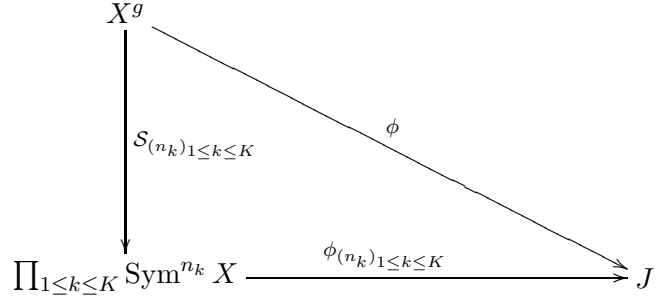
\begin{equation*}
\xymatrix{X^g \ar@{->}[ddd]^{\cS_{(n_k)_{1\le k\le K}}} \ar@{->}[dddrrrrr]^{\phi} &&&&&\\
&&&&&\\
&&&&&\\
\prod_{1\le k\le K}\Sym^{n_k}X \ar@{->}[rrrrr]^{\phi_{(n_k)_{1\le k\le K}}}&&&&& 
J}
\end{equation*}

\medskip 

We can see the divisor  $P$ 
as a point on the semi-symmetric product $\prod_{1\le k\le K}\Sym^{n_k}X.$
We call $U\in X^g$ the unique $g$-uple such that $\cS(U)=P$.

\smallskip 

We assume that  $P$ is a non-special divisor.
The map $\phi_{(n_k)_{1\le k\le K}}$ is thus a local diffeomorphism
at $P$. However, the 
map $\phi$ is not a local diffeomorphism at $U$ because
$\cS$ is not, unless the partition $\bn$ is
$(1,1,  \ldots, 1)$. We shall allow ourselves to
write $\phi'$ instead of 
 $\phi_{(n_k)_{1\le k\le K}}$ in this section.

The previous Section~\ref{subsection:defons}
provides an explicit analytic description of the
maps $\phi$ (resp. $\phi'$) at $U$ (resp. $P$).
We have local coordinates 
$(z_{k,m})_{1\le k\le K;\, 1\le m\le n_k}$ at
$U\in X^g$.
If $R$ is a $g$-uple  in $X^g$
we write $$R=(R_{k,m})_{1\le k\le K;\, 1\le m\le n_k}$$
and we set 
$$z_{k,m}(R)=q_{\gamma_k}(R_{k,m})-q_k.$$

For every $1\le k\le K$ we define the series
$\bFF_k\in \left( \CC[[x_k]]\right)^{\cB^1_{\rm DR}}$ to be the formal
integral
$$\bFF_k(x_k)=\int_0^{x_k} \bff_k(x_k) dx_k=\bFF_{k,1}\times x_k
+\bFF_{k,2}\times x_k^2+\dots$$
where $\bff_k$ is the vector defined in Equation~(\ref{eq:taylorfk}).
This is the Taylor expansion of the Jacobi
map at the point $P_k$.
The Taylor expansion of $\phi$
at $U$ is

$$\Taylor(\phi, U)=\sum_{1\le k\le K}\sum_{1\le m\le n_k }
\bFF_k(z_{k,m})$$
and it lies in $$\left(\CC[[(z_{k,m})_{1\le k\le K;\, 1\le m\le n_k} ]]\right)^{\cB^1_{\rm DR}}.$$

We stress that for every $\omega$ in $\cB^1_{\rm DR}$, the corresponding
coordinate in the above Taylor expansion is a series of type

\begin{equation}\label{eq:typeT}
(l^\Theta,\Theta\times \bun)
\end{equation}
 because it is mostly the expansion of a modular
form at a point which is not too close to the boundary of 
the unit disk.

We now split this Taylor expansion  in two pieces. 
We write 
$$\Taylor(\phi, U)=T_1+T_2$$
where $T_1$ is a sort of  principal part
\begin{eqnarray}\label{eq:decoup}
\nonumber  T_1&=&\sum_{1\le k\le K}\, \, \sum_{1\le m\le n_k }\sum_{1\le j\le n_k} 
\bFF_{k,j}\times z_{k,m}^j\\
\nonumber &=&\sum_{1\le k\le K } \sum_{1\le j\le n_k } \bFF_{k,j}
\sum_{1\le m\le n_k } z_{k,m}^j\\
\nonumber&=&\sum_{1\le k\le K } \sum_{1\le j\le n_k } \bFF_{k,j}\times 
\nu_{P_k,n_k,j}(\{R_{k,1}, R_{k,2}, \ldots, R_{k,n_k}\})\\
&=&\sum_{1\le k\le K } \sum_{1\le j\le n_k } \bFF_{k,j}\times 
\nu_{P_k,n_k,j}(\cS(R)),
\end{eqnarray}
and $T_2=\Taylor(\phi, U)-\phi_1$ is the corresponding remainder.
It is clear from Equation~(\ref{eq:decoup}) that 
$T_1$ can be written
$$T_1=T'_1\circ \cS$$
where $T'_1$ is the first order term (the differential) in the Taylor
expansion of $\phi'$ at $P=\cS(U)$. We write 
$T_2=T'_2\circ \cS$ where $T'_2$ is the corresponding remainder.
Our goal now is to prove that in many circumstances the 
principal term $T'_1$ dominates the remainder $T'_2$.
In order to control how 
 close is $R$ to $U$ (or equivalently how close is $\cS(R)$ to $P$)
we set

\begin{eqnarray*}
\epsilon (R) &=& \max_{1\le k\le K;\, 1 \le m \le n_k}\left|
z_{k,m}(R_{k,m})\right|\\
&=&\max_{1\le k\le K;\, 1 \le m \le n_k}\left|
q_{\gamma_k}(R_{k,m}) -q_k\right|
\end{eqnarray*}

and
 $$\eta (\cS(R)) =\max_{1\le k\le K;\, 1 \le j \le n_k}\left|\nu_{P_k, n_k, j}(\{R_{k,1},
R_{k,2}, \ldots, R_{k,n_k}\})\right|$$
and we assume that  $\eta(\cS(R))$  is
$\le \exp(-l^\Theta)$.

\smallskip

We first  bound the remainder  $T_2(R)=T'_2(\cS(R))$ from above
in terms of $\eta (\cS(R))$.
Using Buckholtz inequality~(\ref{eq:buckholtz})
we can bound the coordinates $z_{k,m}(R)=q_{\gamma_k}(R_{k,m})-q_k$
for every $1\le k\le K$ and $1\le m\le n_k$ in terms
of the bound $\eta (\cS(R))$ on Newton power  functions:
$$|z_{k,m}(R)|\le 5 \eta(\cS(R))^{1/n_k}$$
so 
$$|z_{k,m}(R)|^{n_k+1}\le 5^{g+1}\times \eta (\cS(R)) ^{(n_k+1)/n_k}\le 
5^{g+1}\times \eta (\cS(R))^{1+1/g}.$$

Using the type estimate in Equation~(\ref{eq:typeT}) and
Equation~(\ref{eq:majreste1})
we deduce that
\begin{equation}\label{eq:T2}
\left| T'_2(\cS(R)) \right|_\infty \le \exp(l^\Theta)\times \eta (\cS(R))^{1+1/g}.
\end{equation}

\medskip 

We now  bound from below  $T_1(R)=T'_1(\cS(R))$ in terms of $\eta (\cS(R))$.
This time we need some hypothesis on the divisor $P$. We assume
that  we are in the context of
Theorem~\ref{theorem:torsion_main}. We call
$\Omega$  the  cuspidal divisor $D_0$
constructed in Section~\ref{sec_constr_D} and we assume that the class of $P-\Omega$
corresponds to a point $x$ in the Ramanujan subspace $W_f\subset J$.
We note that $T'_1$ is a linear map and its matrix in the bases
$(\nu_{P_k,n_k,j})_{1\le k\le K;\, 1\le j\le n_k}$ and $\cB^1_{\rm DR}$ is the matrix
$\cM_P^{ana}$ of Equation~(\ref{eq:mpana}). 
The coefficients in this matrix are $\le l^\Theta$ in absolute value.
The  determinant of this matrix 
is the analytic  conditioning. 
According to Equation~(\ref{eq:anacon})
the analytic conditioning at such 
a  Ramanujan $l$-torsion
divisor  is
$\ge \exp(-(\# V_f)^\Theta)$.
We deduce
that
\begin{equation}\label{eq:T1}
\left| T'_1(\cS(R)) \right|_\infty \ge \exp(-(\# V_f)^\Theta)
\times \eta (\cS(R)).
\end{equation}

\smallskip 

We deduce from Equations~(\ref{eq:T1}) and  (\ref{eq:T2})
that $$\left| T'_1(\cS(R))\right|_\infty \ge 2\left| T'_2(\cS(R))\right|_\infty$$ as soon as $\epsilon (R)\le \exp(-(\# V_f)^\Theta)$.

\smallskip 

Using Lemma~\ref{lemma:defsur} below we deduce that there exist
two 
absolute  positive constants $\kappa_1$  and $\kappa_2$
such that
if $$\left|\phi (R)-\phi'(P)\right|_\infty\le \exp(-(\# V_f)^{\kappa_1})$$ then
there exists a divisor  $S$ in   $\prod_{1\le k\le K}\Sym^{n_k}X$
such that $$\eta(S)\le \exp((\# V_f)^{\kappa_2})\times \left|\phi(R)-\phi'(P)\right|_\infty$$
and $\phi'(S)=\phi(R)$. 

Using again that the conditioning
of $P$ is $\ge \exp(-(\# V_f)^\Theta)$ together
with the estimates  in Section~\ref{subsection:defons}, e.g. 
 Equation~(\ref{eq:nonspe}), 
we deduce that if the constant $\kappa_1$  has been chosen big enough, then 
 this divisor $S$ is non-special. Therefore
$S=\cS(R)$ and we have proven that
$$\eta (\cS(R))\le \exp((\# V_f)^{\kappa_2})\times \left|\phi (R)-\phi' (P)\right|_\infty$$
therefore
$$\epsilon (R)\le \exp((\# V_f)^{\kappa_3})\times \left|\phi (R)-\phi' (P)\right|_\infty$$ for some absolute positive constant $\kappa_3$.

 This relation between the error
in $J$ and the error in $X^g$ 
finishes the proof of Theorem~\ref{theorem:torsion_main}.

\medskip 

To finish this section there remains to state and prove  
Lemma~\ref{lemma:defsur}.

\begin{lem}[Perturbation of a non-singular linear map]\label{lemma:defsur}
Let $g\ge 1$ be a positive integer and let
$\bar P(\bzero,\bun)\in \CC^g$ be the closed polydisk with polyradius
$\bun$ and centered at the origin. Let
$T : \bar P(\bzero,\bun) \rightarrow \CC^g$
be a continuous function and let
$T_1 : \CC^g  \rightarrow \CC^g$
be a linear function.
Set $T_2=T-T_1$ and assume that for every $R=(R_1, \ldots, R_g)$
in $\bar P(\bzero,\bun)$ we have
$|T_1(R)|_\infty> 2|T_2(R)|_\infty$. Let $K>0$ be a positive real
number such that  for every $R$ in $\CC^g$
$$\left|T_1(R)\right|_\infty \ge K\times \left| R\right|_\infty.$$
Then  the image by $T$ of the
polydisk $\bar P(\bzero,\bun)$ contains
the polydisk $\bar P(\bzero,\frac{K}{2}\times \bun)$ 
with polyradius $K/2\times \bun$.
\end{lem}

For every real number $t\in [0,1]$ we denote
by $\cS_t$ the image by $T_1+tT_2$
of the 
$L^\infty$-sphere 
with radius $1$. 
 Then $\cS_t$
is contained in $$\CC^g-\bar P(\bzero,\frac{K}{2}\times \bun )$$
for every $t\in [0,1]$. So $\cS_1$
is homologous
 to $\cS_0$ in
$\CC^g-  \bar P(\bzero,\frac{K}{2}\times \bun ) $
and its class in
$$H_{2g-1}(\CC^g-  \bar P(\bzero,\frac{K}{2}\times \bun ),\ZZ )$$ is non-zero.
So  for 
every $x$ in $\bar P(\bzero,\frac{K}{2}\times \bun )$, 
the class
of $\cS_1$ in $H_{2g-1}(\CC^g-\{x\},\ZZ)$ is non-zero either. Assume
now that $x$ does not belong to the image by $T$
of the polydisk $\bar P(\bzero,\bun)$. Then
$\cS_1$ is the boundary of $T(\bar P(\bzero,\bun))\subset \CC^g-\{x\}$ 
and its class in $H_{2g-1}(\CC^g-\{x\},\ZZ)$ is trivial.
A contradiction. \hfill $\Box$

\section{Final result of this chapter}\label{section:finalC}

We now state the final result of this chapter.

\begin{thm}[Approximating $V_f$ over the complex numbers]\label{theorem:torsion_main2} 
There exists a  deterministic algorithm that takes as  input 
an even integer $k > 2$,  a prime integer $l>6(k-1)$, a finite field $\FF$
with characteristic $l$,  a ring epimorphism 
$f : \TT(1,k)\rightarrow  \FF$,
and a cuspidal divisor $\Omega$ on $X_1(5l)$
like the divisor $D_0$  constructed in Section~\ref{sec_constr_D},
and computes complex approximations for every element in
$W_f\subset J_1(5l)$,  the image of $V_f\subset J_1(l)$ by
$B_{5l,l,1}^*$. Here $V_f\subset J_1(l)$ is defined by 
Equation~(\ref{eq:defVltor}) and we  assume that 
the image of the Galois representation $\rho_f$ associated with $f$
contains $\SL  (V_f)$. The algorithm returns 
 for every element  $x$ in $W_f$ 
the unique degree $g$ effective
divisor  $Q_x$ such that 
$Q_x-\Omega$ lies in the class represented by $x$. 
More precisely, the algorithm returns the cuspidal part 
$Q_x^{\rm cusp}$ of $Q_x$ and a complex approximation of its
finite part $Q_x^{\rm fin} = \sum_{1\le n\le d_x}Q_{x,n}$.
   Every point $Q_{x,n}$ is given by complex approximations
of its affine coordinates
$(b(Q_{x,n}), x(Q_{x,n}))$ in the  plane model $C_l$ of $X_1(5l)$ given
in Section~\ref{subsection:descripCl}. 
The running time of the algorithm is $\le (m\times \# V_f)^\Theta$ for some absolute constant $\Theta$.
Here $\# V_f$ is the cardinality of the Galois representation $V_f$
and $m$ is the required absolute accuracy.
\end{thm}

There are  only two differences between this Theorem~\ref{theorem:torsion_main2}  and the previous Theorem~\ref{theorem:torsion_main}. Firstly we claim that we can separate
the cuspidal and the   finite part of $Q_x$. Secondly we return algebraic
coordinates $b$ and $x$ for the points $Q_{x,i}$ rather than analytic
ones.

Indeed,  Theorem~\ref{theorem:torsion_main} gives
us for every point $P$ 
in the support of $Q_x$
an 
analytic coordinate $(\gamma, q)$. We want to decide if the point $P$
in question is equal to $\gamma(\infty)$. Let  $\KJ$ be Klein's modular
function.
We know that  both the degree and  the logarithmic
height of $\KJ (P)$ are $\le (\# V_f)^\Theta$. So there exists an absolute
constant $\kappa_1$ such that if $1/\KJ (P)$ is 
$\le \exp(-(\# V_f)^{\kappa_1})$ then it is zero. But 
the  Petersson-Rademacher inequality tells us that 
there exists an $\epsilon >0$ such that if $x$ is a complex
number bounded by $\epsilon$ in absolute value then $$\left|\bJ(x)\right|
\ge \frac{1}{2|x|}$$
where $\bJ(x)$ is  Klein's
series given in Equation~(\ref{eq:744}). 
Using Equation~(\ref{eq:KJJ})
we deduce that there exists an absolute
constant $\kappa_2$ such that if $q$ is 
$\le \exp(-(\# V_f)^{\kappa_2})$ then $q=0$ and 
the point $P=(\gamma,q)$
is the  cusp $\gamma(\infty)$. 
This explains why we can recognize cusps.

Now assume that $P=(\gamma,q)$ is a point in the support of some
$Q_x$ and assume that $P$ is not a cusp. Then
$|q| \ge \exp(-(\# V_f)^{\kappa_2})$ and
 in order to compute $b(P)$ and $x(P)$ we just substitute
 $q_\gamma$ by $q$ in the expansions for $x$
and $b$ at $\gamma(\infty)$ computed in Section~\ref{sec:modfunjmc}.

%JMCfin
\chapter{Computing $V_{\protect\lowercase{f}}$ 
modulo $\protect\lowercase{p}$}\label{sec_couveignes_modp}

\author{J.-M. Couveignes}

\bigskip

\bigskip

% author: Jean-Marc
%JMCdebut

In this chapter we  address the problem of computing 
in  the group of $l^k$-torsion 
rational points in  the jacobian variety 
of algebraic curves  over  finite fields, with an application to
 computing modular representations.

Let $p$ be a prime and let $\Fp=\ZZ/\! p\ZZ$ be the  field 
with $p$ elements. Let 
$\Fbar_p$ be an  algebraic closure
of $\Fp$. For any power $q$ of  $p$ we call $\Fq\subset \Fbar_p$ the field
with $q$ elements. 
Let $\AA^2\subset \PP^2$ be the affine and projective planes over $\Fq$.
Let $C\subset \PP^2$ be  a plane projective geometrically integral 
curve  over $\Fq$. Let 
 $X$ be  its smooth projective model and let $J$ be  the jacobian variety of
$X$. 
We note  $g$  the genus of $X$ and $d$ the degree of $C$.
We assume that we are given the numerator of the zeta function of the
function
field $\Fq(X)$. So we know the characteristic polynomial of the
Frobenius endomorphism $F_q$ of  $J$. This is a monic degree $2g$ polynomial
$\chi(x)$ with integer coefficients.

Let $l\not = p $ be a prime integer and let $k\ge 1$ be an integer.
 We
look for a {\it nice generating set} for the group $J[l^k](\Fq)$ 
of  $l^k$-torsion points in $J(\Fq)$. By {\it nice}  we mean that the
generating set $(g_i)_{1\le i\le I}$ should induce  a decomposition of
$J[l^k](\Fq)$ as  a direct product $\prod_{1\le i\le I} <g_i>$
of cyclic subgroups with non-decreasing orders. 
Given such a generating set and an $\Fq$-endomorphism of $J$,  we
also want to describe the action
of this endomorphism on  $J[l^k](\Fq)$ by an $I\times I$ integer matrix.

 By an algorithm in this paper we usually mean a
probabilistic Las Vegas algorithm. 
 In some places we shall give
deterministic algorithms or probabilistic Monte-Carlo algorithms, but
this will be stated explicitly. See Section~\ref{sec:comp} for a reminder
of computational complexity theory.
The main  reason for using probabilistic
Turing machines is that we shall need to construct 
generating 
sets for  the Picard group of curves over finite fields. Solving such a problem in
the deterministic  world is out of reach at this time.
See Section~\ref{sec:comp}.

In Section~\ref{section:classics} we recall how to compute in the 
Picard group
$J(\Fq)$. Section~\ref{section:picking}  gives a naive algorithm for picking 
random elements in this
group. Pairings are useful when looking for relations between
divisor classes. So we recall how to compute pairings in 
Section~\ref{section:pairings}. Section~\ref{section:divisiblegroups} is
concerned with characteristic subspaces for the action of Frobenius
inside the $l^\infty$-torsion of $J(\Fbar_p)$. In 
Section~\ref{section:kummer} we look for  a convenient surjection from $J(\Fq)$ onto
its
$l^k$-torsion subgroup. We use the Kummer exact sequence and the
structure
of the ring generated by the Frobenius endomorphism. 
In Section~\ref{section:relations} we give
an  algorithm that, on input a degree $d$ plane
projective curve  over $\Fq$, plus some information on its
singularities, and the zeta function of its function field, 
   returns a nice generating set
for the group of $l^k$-torsion points inside $J(\Fq)$ in
probabilistic  polynomial time in $\log q$, $d$ and $l^k$. 
In Section~\ref{section:ramanujan} we apply the general algorithms 
in
Section~\ref{section:relations}  to the
modular curve $X_1(5l)$  in
 order to compute explicitly a  modular
representations  $V_f$ modulo $l$.
Such a representation
 modulo $l$ can be realized   as a
subgroup $W_f$ 
 inside the $l$-torsion of $J_1(5l)/\QQ$. The idea is to
compute the reduction modulo $p$ of the
group scheme  $W_f$ as a subgroup of   $J_1(5l)/\Fp$,    for many small primes
$p$. One can then lift using the Chinese Remainder Theorem, as will be
explained
in Section~\ref{sec_prob_algorithm}. 

%newmodif

\begin{remark} The symbol $\Theta$ in this chapter  stands for
a positive effective absolute constant. So any statement containing
this symbol becomes true if the symbol is replaced in every occurrence by some 
large enough real number.
\end{remark}

\section{Basic algorithms for plane curves}\label{section:classics}

In this section, we  recall elementary results about computing 
in the Picard group of an algebraic curve over a finite field. See
\cite{hache, volcheck, Diem} for a more detailed treatment.

\subsection{Finite fields}
We should first explain how finite fields are represented.
The prime field $\Fp$ is just $\ZZ/p\ZZ$ so elements in it
are represented as integers in $[0,p-1[$.
The base field $\Fq$ is given 
as $\Fp[x]/f(x)$ where $f(x)$ is an 
 irreducible unitary polynomial   with degree
$a$ in  $\Fp[x]$ where $p$ is the characteristic and
$q=p^a$.  A finite  extension of $\Fq$ is given
as $\Fq[y]/h(y)$ where $h(y)$
is a unitary irreducible polynomial in $\Fq[y]$. We shall never 
use two extensions of $\Fq$ simultaneously.
Remind polynomial
factoring in $\Fq[x]$
is Las Vegas probabilistic polynomial time in $\log q$ and the degree of the
 polynomial to be
factored \cite[Chapter 14]{VZG}.  

\subsection{Plane projective curves and their smooth model}\label{sec_couveignes_plane_curves}

We now explain how curves are represented in
this
paper. 
To start with, a  projective  plane curve $C$ over $\Fq$ is  given by  
a degree $d$ homogeneous
polynomial $E(\bx,\by ,\bz)$ in the three
variables $\bx$, $\by$ and $\bz$, with coefficients in $\Fq$.
The curve $C$ is assumed to be absolutely integral.
By a {\it point} on $C$ we mean a geometric point: an element
of $C({{\Fbar_q }})$ where ${{\Fbar_q }}=\Fbar_p$ is the algebraic closure
of $\Fp$ fixed in the introduction.
Any  ${{\Fbar_q }}$-point on $C$ can be represented by its affine or projective
coordinates. 

Let $X$ be a smooth model of $C$ and let $X\rightarrow C$ be
the  desingularization map.
If 
$P\in X({{\Fbar_q }})$ is  a geometric point on $X$
 above a singular point $S$
on $C$, we say that $P$ is a {\it singular branch}. 
The {\it conductor} $\gotC$ is an effective  divisor
on $X$. It is the closed subscheme of $X$ defined by the sheaf 
$\Ann_{\cO_C}(\cO_X/\cO_C)$. Every multiplicity in $\gotC$ is even.
Some authors call $\gotC$  the {\it adjunction divisor}.
Its support consists 
of all singular branches. The conductor  expresses
the local behaviour of the map
$X\rightarrow C$.
See \cite[IV.1]{Serre}, \cite{gorenstein}.
We have  $\deg(\gotC)=2\delta$ where
$$\delta = \frac{(d-1)(d-2)}{2} -g$$
is the difference between the arithmetic genus 
of $C$ and
the  geometric genus of $X$. Since $\delta\le {(d-1)(d-2)}/{2}$, 
the support of $\gotC$
contains  at most  ${(d-1)(d-2)}/{2}$ geometric points 
in $X({{\Fbar_q }})$. So the field of definition of any
singular branch on $X$  is an extension of $\Fq$ with 
degree $\le  {(d-1)(d-2)}/{2}$.
A modern reference for singularities of plane curves
is \cite[Section 5.8]{casas}.

The smooth model $X$ of $C$ is not given as a projective variety.
Indeed, we shall only need  a nice local
description of $X$ above
 every singularity of $C$. This means that we need a 
list of all singular points on $C$, and a 
list (a labelling) of all  points
in $X({{\Fbar_q }})$ lying  above every singularity of $C$ (the 
singular branches), and  a uniformizing
parameter at every such branch. We also need the  Laurent series expansions of
affine  plane coordinates in terms of all these uniformizing
parameters.

%Singularities can be detected by computing and factoring discriminants. This is probabilistic polynomial in $d$ and $\log q$.

More precisely,  let
$P\in X({{\Fbar_q }})$ be a geometric point above a singular point $S$, and 
let $v$ be the corresponding
valuation.  The field of definition of $P$
 is an extension
field $\FF_P$ of $\Fq$ with  degree  $\le {(d-1)(d-2)}/{2}$. 
Let $x_S$ and $y_S$
be affine coordinates that vanish at the  singular point $S$ on $C$. We need
a local parameter $t_P$ at $P$ and  expansions $x_S=\sum_{k\ge
  v(x_S)}a_kt_P^k$ and $y_S=\sum_{k\ge
  v(y_S)}b_kt_P^k$ with coefficients in $\FF_P$. 

Because these expansions are not finite, we just assume that we are given
a black box  that on input a positive integer $n$
returns  the first $n$ terms in all these expansions. In all the cases
we shall be interested in, this black box will contain a Turing
machine that answers in time polynomial in $n$ and $\log q$ and the genus $g$.
This is the case for
curves with ordinary multiple points for example. We have
shown in Section~\ref{sec:modfunjmc} 
that this is also the case for the standard plane model of 
   modular curves $X_1(5l)$. Using general normalization and factorization algorithms \cite[Section 2]{Diem}
one may show that this is indeed the case for all plane curves without any restriction, but this is beyond the scope of this
text.

We may also assume that  we are given the conductor $\gotC$ of $C$ as a
combination of singular branches with even  coefficients.  
The following algorithms still work if  the conductor
is  replaced by any  divisor
$\gotD$  that 
is greater than the conductor and has polynomial degree  in $d$. Such
a divisor can be found easily: the singular
branches on $X$ are supposed to be known already, and the multiplicities
are bounded above by ${(d-1)(d-2)}/{2}$.

\subsection{Divisors, forms,  and functions}

Smooth ${{\Fbar_q }}$-points on $C$ are  represented by their affine or projective
coordinates. Labelling for the  branches above singular points is
given in the description of $X$. 
So we know how to represent divisors on $X$.
For any   integer $h\ge 0$   the $\Fq$-linear space
$H^0(\PP^2,\cO_{\PP^2}(h))$
of 
 degree $h$ homogeneous polynomials in $\bx$,
$\by$, and $\bz$ has  dimension  ${(h+1)(h+2)}/{2}$. A
basis for it is made of all monomials of the form $\bx^a\by^b\bz^c$
with $a,b,c \in \NN$ and  $a+b+c=h$. 
We denote by
$\cO_{X}(h)$  the pullback of $\cO_{\PP^2}(h)$ to $X$. 
Let  $F$ be a degree $h$   form on $\PP^2$ having non-zero pullback  $F_X$ on  $X$. Let $\Delta= (F_X)$ be the divisor of this restriction. 
The map $f\mapsto
\frac{f}{F_X}$ is a bijection from $H^0(X,\cO_{X}(h))$ to the linear space $H^0(X,\cO_X(\Delta))$ associated with $\Delta$.

We assume that we are given
 a divisor  $\gotD$  bigger than the conductor
$\gotC$. We assume that the degree of $\gotD$ is $\le d^\Theta$. 
We have explained in the previous Section~\ref{sec_couveignes_plane_curves}
 how to find
such a divisor.
The dimension of $H^0(X, \cO_X(h)(-\gotD))$ is at least
$dh+1-g-\deg({\gotD})$ and is equal to this number when it
exceeds $g-1$. This is the case if  $$h\ge d+\frac{\deg\gotD}{d}.$$ The
dimension of 
$H^0(X, \cO_X(h)(-\gotD))$ is     greater than $2g$ if 
$$h\ge \frac{3d}{2}+\frac{\deg\gotD}{d}.$$ 

We take $h$ to be the smallest integer fulfilling this condition.
The composite map $X{\onto} C {\into} \PP^2$ induces a map $$\rho_h : H^0(\PP^2,\cO_{\PP^2}(h))\rightarrow H^0(X,\cO_X(h)).$$ The image
of $\rho_h$ contains   $H^0(X, \cO_X(h)(-\gotD))$. This is known as Noether's residue theorem \cite[Theorem 7]{gorenstein}.
It will be convenient to 
describe $H^0(X,\cO_X(h)(-\gotD))$ as a quotient
$$0\rightarrow \Ker \rho_h
\rightarrow \rho_h^{-1} (H^0(X,\cO_X(h)(-\gotD) ))  \stackrel{\rho_h}{\rightarrow} H^0(X,\cO_X(h)(-\gotD)) \rightarrow 0.$$

We need 
linear equations for  $$\rho_h^{-1} (H^0(X,\cO_X(h)(-\gotD) ))  \subset H^0(\PP^2,\cO_{\PP^2}(h)).$$ We consider a
generic homogeneous form
$$F(\bx,\by,\bz)=\sum_{a+b+c=h}\epsilon_{a,b,c}\bx^a\by^b\bz^c$$ of 
degree $h$
in  $\bx$, $\by$  and $\bz$.
For every branch $P$  above a singular point $S\in C$ 
%with homogeneous coordinates $$[X,Y,Z]$$ \noindent 
(assuming for example  that  $S$ has non-zero
$\bz$-coordinate)
we replace
in $F(\frac{\bx}{\bz},\frac{\by}{\bz},1)$  the affine coordinates 
$x=\frac{\bx}{\bz}$ and $y=\frac{\by}{\bz}$
by their   expansions as  series in  the local parameter $t_P$ at this
branch. We ask  the resulting series in $t_P$ to have
valuation at least the multiplicity of $P$ 
in the divisor  $\gotD$.
Every singular branch thus  produces linear equations
in the $\epsilon_{a,b,c}$.  The
collection of all such equations defines the subspace $\rho_h^{-1} (H^0(X,\cO_X(h)(-\gotD) ))$.

A basis for the  subspace $\Ker \rho_h$ of $\rho_h^{-1} (H^0(X,\cO_X(h)(-\gotD) ))$
consists of all  $\bx^a\by^b\bz^c E(\bx,\by,\bz)$ with 
$a+b+c=h-d$. We fix a supplementary space $M_C$ to $\Ker \rho_h$
in $\rho_h^{-1} (H^0(X,\cO_X(h)(-\gotD) ))$ and we assimilate $H^0(X,\cO_X(h)(-\gotD))$  to it.

Given
a homogeneous form in three variables  one can  compute its  divisor on $X$ 
using resultants and the given expansions of affine coordinates in terms of
the
local parameters at every singular branch.
A function is given as a quotient of two forms.

\subsection{The Brill-Noether algorithm}\label{subsubsection:brill}

The linear  space
$M_C$ computed in the previous paragraph is isomorphic to
$H^0(X,\cO_X(h)(-\gotD))$ via the map $\rho_h$. This space
allows us to  compute in the group $J (\Fq)$ of $\Fq$-points
in the jacobian of $X$. We fix an effective $\Fq$-divisor
$\Omega$ with
degree
$g$ on $X$. This $\Omega$ will play the role of origin:
a point  $\alpha \in J (\Fq)$ is
represented by a divisor $A-\Omega$ in the corresponding linear
equivalence
class, where $A$ is an effective $\Fq$-divisor on $X$ with degree $g$.
Given another point  $\beta \in J (\Fq)$ by a similar divisor
$B-\Omega$, we can compute the space $H^0(X, \cO_X(h)(-\gotD-A-B))$ which is
non-trivial  and pick a non-zero form $F_1$ in it.  The divisor of $F_1$
is
$(F_1)=A+B+\gotD+R$ where $R$ is an effective divisor with  degree
$hd-2g-\deg(\gotD)$. The linear space $H^0(X,\cO_X(h)(-\gotD-R-\Omega))$
has dimension at least $1$. We pick a non-zero form $F_2$ in it. It
has divisor $(F_2)=\gotD+R+\Omega+D$ where $D$ is effective with
degree $g$. And $D-\Omega$ is linearly equivalent to
$A-\Omega+B-\Omega$.

In order to invert the class $\alpha$  of $A-\Omega$ we pick a
non-zero form  $F_1$ in $H^0(X,\cO_X(h)(-\gotD-2\Omega))$.  The divisor of $F_1$
is
$(F_1)=2\Omega+\gotD+R$ where $R$ is an effective divisor with  degree
$hd-2g-\deg(\gotD)$. The linear space $H^0(X,\cO_X(h)(-\gotD-R-A))$
has dimension at least $1$. We pick a non-zero form $F_2$ in it. It
has divisor $(F_2)=\gotD+R+A+B$ where $B$ is effective with
degree $g$. And $B-\Omega$ is linearly equivalent to
$-(A-\Omega)$.

\begin{lem}[Arithmetic operations in the jacobian]\label{lemma:arithmeticoperations}
Let $C/\Fq$ be a degree $d$ plane projective absolutely integral
curve. Let $g$ be the geometric genus of $C$.
Assume that we are given the smooth model $X$ of $C$ and a 
   $\Fq$-divisor with degree
$g$  on $X$, denoted $\Omega$. We assume that  $\Omega$
is given as a difference between two effective divisors with degrees bounded
by  $d^\Theta$. This $\Omega$ serves as an origin.
Arithmetic operations in the Picard group  $\Pic^0(X /\Fq)$ can be performed
in   time  polynomial  in $\log q$ and $d$. This includes
addition, subtraction and comparison of divisor classes.
\end{lem}
  
If $\Omega$ is not effective, we use Lemma~\ref{lemma:explicitRR} below to
compute a non-zero function $f$  in $H^0(X,\cO_X(\Omega))$ and we
write $\Omega'=(f)+\Omega$. This is an effective divisor with degree
$g$. We replace $\Omega$ by $\Omega'$ and finish as in the paragraph
before Lemma~\ref{lemma:arithmeticoperations}
\hfill $\Box$

\medskip 

We now recall
the principle  of the Brill-Noether algorithm for computing
complete linear series. Remind functions in $\Fq(X)$ are represented as
quotients of forms.

\begin{lem}[Brill-Noether]\label{lemma:explicitRR}
There exists an algorithm that on input 
 a degree $d$ plane projective absolutely integral 
curve $C/\Fq$ and the  smooth model $X$ of $C$
 and  two effective $\Fq$-divisors $A$ and $B$ on
 $X$,  computes a basis for $H^0(X,\cO_X(A-B))$  in time  polynomial 
in $d$ and $\log q$ and the degrees of $A$ and $B$.
\end{lem}

We assume that $\deg(A)\ge \deg(B)$, otherwise $H^0(X,\cO_X(A-B))=0$. Let
$a$ be the degree of $A$. We let $h$ be the smallest integer
such that $$hd-g+1 > a+\deg\gotD.$$
The space $H^0(X,\cO_X(h)(-\gotD-A))$ is non-zero. It
 is contained in the 
image of the  map 
$$\rho_h : H^0(\PP^2,\cO_{\PP^2}(h))\rightarrow H^0(X,\cO_X(h))$$
 so that we can represent it
  as a subspace of $H^0(\PP^2,\cO_{\PP^2}(h))$. We pick a non-zero  form
$f$  in   $H^0(X, \cO_X(h)(-\gotD-A))$ and compute its divisor
$(f)=\gotD+A+D$.  
The space $H^0(X,\cO_X(h)(-\gotD-B-D))$  is contained in the 
image of the  map $\rho_h$
 so that we can represent it
  as a subspace of $H^0(\PP^2,\cO_{\PP^2}(h))$. 
We compute forms $\gamma_1$, $\gamma_2$,
  \ldots , $\gamma_k$ in $H^0(\PP^2,\cO_{\PP^2}(h))$ such that their 
images by $\rho_h$ provide
  a basis for  $H^0(X,\cO_X(h)(-\gotD-B-D))$. A basis for 
$H^0(X, \cO_X(A-B))$ is made
  of the functions $\frac{\gamma_1}{f}$, $\frac{\gamma_2}{f}$, \ldots,
  $\frac{\gamma_k}{f}$. 
\hfill $\Box$

We deduce an explicit moving lemma for divisors.

\begin{lem}[Moving divisor lemma I]\label{lemma:moving}
There exists an algorithm that on input 
 a degree $d$ plane projective absolutely integral
curve $C/\Fq$ and the  smooth model $X$ of $C$
 and  a degree zero  divisor $D=D^+-D^-$ and an
effective divisor $A$ with degree $<q$ on 
 $X$ computes a divisor $E=E^+-E^-$ linearly equivalent to $D$ and
disjoint to $A$   in time  polynomial 
in $d$ and $\log q$ and the degrees of $D^+$,  and
$A$. Further
the degree of $E^+$ and $E^-$ can be taken to be
$\le 2gd$.
\end{lem}

Let $O$ be an $\Fq$-rational divisor on $X$ such that
$1\le \deg (O)\le d$
 and disjoint to $A$.  We may take $O$ to be
a well chosen fiber of some plane coordinate function on $X$.
We  compute the linear space
$H^0(X,\cO_X(D^+-D^-+2g O))$. The subspace
of $H^0(X,\cO_X(D^+-D^-+2g O))$
consisting  of functions $f$  
 such that
$(f)+D^+-D^-+2g O$ is not disjoint to $A$ is contained  in a union of
at most $\deg(A) < q$ hyperplanes. We conclude invoking 
Lemma~\ref{lemma:inequ} below.\hfill $\Box$

There remains to state and prove the

\begin{lem}[Solving inequalities]\label{lemma:inequ}
Let $q$ be  a prime power, $d\ge 2$  and $n\ge  1$ two  integers and
let $H_1$, \ldots , $H_n$ be hyperplanes inside  $V=\FF_q^d$,
each  given by a linear equation. Assume that 
$n<q$. There exists a deterministic
algorithm
that finds a vector in $U=V-\bigcup_{1\le k\le n}H_k$ in time polynomial in $\log q$, $d$ and $n$.
\end{lem}

This is proved by lowering  the dimension $d$. For $d=2$ we pick
any affine line $L$ in $V$ not containing the origin. We  observe that
there are at least $q-n$ points in 
$U\cap L=L-\bigcup_{1\le k\le n}L\cap H_k$. We enumerate points in $L$ until we
find one which is not in any $H_k$. This requires at most $n+1$
trials.

Assume now that  $d$ is bigger than $2$. Hyperplanes in $V$ are parametrized
by the projective space $\PP(\hat V)$ where $\hat V$ is the dual of
$V$. We enumerate points in $\PP(\hat V)$ until we find a hyperplane
$K$ distinct from every  $H_k$. We compute a basis for $K$ and an
equation for every $H_k\cap K$ in this basis. This way,  we have lowered
the dimension by $1$.\hfill $\Box$

We can strengthen a bit the moving divisor algorithm by removing the
condition that $A$ has degree $<q$.  Indeed, in case this condition is
not met, we call
$\alpha$ the smallest integer such that
$q^\alpha > \deg (A)$ and we set $\beta =\alpha+1$. We apply 
Lemma~\ref{lemma:moving} after base change to the field with $q^\alpha$
elements and find a divisor $E_\alpha$. We call $e_\alpha$ the norm of
$E_\alpha$
from $\FF_{q^\alpha}$ to $\FF_q$. It is equivalent to $\alpha D$. We
similarly construct a divisor $e_\beta$ that is equivalent to $(\alpha+1)D$.
We
return the divisor $E=e_\beta-e_\alpha$. We observe that we can 
take $\alpha \le
1+\log_q \deg (A)$ so the degree of 
the positive part $E^+$ of $E$ is $\le 6g d(\log_q(\deg (A))+1)$.

\begin{lem}[Moving divisor lemma II]\label{lemma:moving2}
There exists an algorithm that on input 
 a degree $d$ plane projective absolutely integral
curve $C/\Fq$ and the  smooth model $X$ of $C$
 and  a degree zero  $\Fq$-divisor $D=D^+-D^-$ and an
effective divisor $A$ on 
 $X$ computes a divisor $E=E^+-E^-$ linearly equivalent to $D$ and
disjoint to $A$   in time  polynomial 
in $d$ and $\log q$ and the degrees of $D^+$,  and
$A$. Further
the degree of $E^+$ and $E^-$ can be taken to be
$\le 6gd (\log_q(\deg (A))+1)$.
\end{lem}

\section{A first approach to picking random divisors}\label{section:picking}

Given a finite field $\Fq$ and a plane projective absolutely
integral  curve $C$ over $\Fq$  with projective smooth  model
 $X$, we call $J$ the jacobian of $X$
and we consider  
two related problems:  picking a random element 
in  $J (\Fq)$  with (close to) uniform distribution and finding a
generating set for (a large subgroup of) $J (\Fq)$. 
Let $g$ be the genus of $X$.
We  assume that we are given  a degree $1$ divisor $O=O^+-O^-$ where $O^+$
and $O^-$ are effective, $\Fq$-rational and have degree bounded by 
$\Theta g^\Theta$ for some positive constant $\Theta$.

We know from \cite[Theorem 2]{MST} that   
 the   group $\Pic^0(X/\Fq)$ is generated by the classes 
$[\pgot -\deg(\pgot)O]$ where $\pgot$ runs over the set
of prime divisors of
degree $\le 1+2\log_q(4g-2)$.  For the convenience of the reader 
we quote this result as a lemma.

\begin{lem}[M\"uller, Stein, Thiel]\label{lemma:Stein}
Let $K$ be an algebraic  function field of one variable
over $\Fq$. Let  $N\ge 0$ be an integer.
Let $g$ be the genus of $K$.
Let $\chi : \Div (K)\rightarrow \CC^*$ be a  character
of finite order which is non-trivial when restricted to
$\Div^0$. Assume  
that $\chi(\gotB)=1$ for every prime divisor $\gotB$ of degree
$\le N$. Then 
$$N<{2\log_q(4g-2)}.$$
\end{lem}

If  $q<4g^2$,  the number of prime divisors of degree $\le 1+2\log_q(4g-2)$
is bounded by $\Theta g^\Theta$. So we  can compute easily
 a small generating set
for $J(\Fq)$.
In the rest of this section, we will assume that the size 
 $q$ of the
field is greater than or equal to $4g^2$. This condition
ensures the
existence of a $\Fq$-rational point.

Picking efficiently and provably random elements in $J (\Fq)$ with
uniform distribution  seems
difficult to us. We first give here an algorithm for efficiently constructing
random
divisors with a distribution that  is far from uniform but still
sufficient to construct a generating set for a large subgroup of
$J (\Fq)$. Once given generators, picking random elements becomes
much easier.

Let $r$ be the smallest prime integer bigger than $30$, $2g -2$ and $d$.
We observe $r$ is less than $\max (4g-4,2d,60)$.
The set  $\cP(r,q)$ of $\Fq$-places with  degree $r$ 
on $X$ has cardinality
$$\#\cP(r,q)= \frac{\# X (\FF_{q^{r}}) -
  \# X (\FF_{q})}{r}.$$

So 
$$(1-10^{-2})\frac{q^{r}}{r}\le \# \cP(r,q)  \le (1+10^{-2})\frac{q^{r}}{r}.$$

Indeed, $\left|\# X (\FF_{q^{r}})-q^{r}-1 \right|\le
2g q^{\frac{r}{2}}$ and $\left|\# X (\FF_{q})-q-1 \right|\le
2g q^{\frac{1}{2}}$.
 So   $$\left | \# \cP(r,q) -
\frac{q^{r}}{r}\right| \le
\frac{4g+3}{r}q^{\frac{r}{2}}\le 8q^{\frac{r}{2}}$$
and $8r q^{\frac{-r}{2}}\le r2^{3-\frac{r}{2}}\le
10^{-2}$ since  $r \ge 31$.

Since we are given a degree $d$  plane model $C$ for the curve
$X$,
we have a degree $d$ map $x : X \rightarrow \PP^1$. This
is the composition of the desingularization map
$X\rightarrow C$ with the  restriction to $C$ of the
rational map $[\bx, \by, \bz]\rightarrow [\bx, \bz]$.
Since $d<r$, the  function $x$ maps $\cP(r,q)$
to the set  $\cU(r,q)$ 
of monic prime polynomials of degree $r$ over
$\Fq$. The cardinality of   $\cU(r,q)$ is
$\frac{q^{r}-q}{r}$ so 
$$(1-10^{-9})\frac{q^{r}}{r}\le \# \cU(r ,q)  \le
\frac{q^{r}}{r}.$$

The fibers of the map $x : \cP(r,q)\rightarrow \cU(r,q)$ have
cardinality between $0$ and $d$.
We can pick a random element in $\cU(r,q)$ with uniform
distribution in the following way: we pick a random monic polynomial
 of
degree $r$ with coefficients in $\Fq$, with uniform
distribution. We  check whether  it is
irreducible. If it is, we output it. Otherwise we start again. This is
polynomial time in $r$ and $\log q$.

Given a random  element in $\cU(r,q)$ with uniform distribution,
we can compute the fiber of $x : \cP(r,q)\rightarrow \cU(r,q)$
above it and, provided  this fiber
 is non-empty,  pick a random element in it
with uniform distribution. If the
fiber is empty, we pick another element in $\cU(r,q)$ until we
find a non-empty fiber. At least one in every  $d\times(0.99)^{-1}$ fibers
is non-empty.
 We
thus define  a distribution $\mu$ on $\cP(r,q)$ and prove the following.

\begin{lem}[A very rough measure]\label{lemma:random}
There is a unique measure $\mu$ on $\cP(r,q)$ such that
all non-empty fibers of the map $x : \cP(r,q)\rightarrow \cU(r,q)$
have the same measure, and all points in a given fiber
have the same measure.
There exists a probabilistic algorithm that picks a random element
in $\cP(r,q)$ with distribution $\mu$ in time polynomial in $d$
and $\log q$. For every subset $Z$ of $\cP(r,q)$ the measure
$\mu (Z)$ is related to  the uniform measure
$\frac{\#Z}{\#\cP(r,q)}$ by

$$\frac{\#Z}{d\#\cP(r,q)}\le   \mu (Z) \le \frac{d\#Z}{\#\cP(r,q)}.$$
\end{lem}

Now let $\cD(r,q)$ be the set  of effective $\Fq$-divisors with  degree $r$ 
on $X$.
Since we have assumed that  $q\ge 4g^2$ we know that   $X$ has at least one
$\Fq$-rational point. 
Let  $\Omega$ be  a degree $r$ effective divisor on
$X/\Fq$. We  associate to every $\alpha$ in $\cD(r,q)$ the class of
$\alpha-\Omega$ in $J (\Fq)$. This defines a surjection $$\phi_r :
\cD(r,q)\rightarrow J (\Fq)$$ with all
its fibers having  cardinality $\#\PP^{r-g}(\Fq)$.
So the set  $\cD( r,q)$  has cardinality
$\frac{q^{r-g+1}-1}{q-1}\# J(\Fq)$. 
 So
$$\#\cP(r,q)  \le   \#\cD(r, q)  \le 
q^{r-g}\frac{1-\frac{1}{q^{r-g+1}}}{1-\frac{1}{q}}q^{g}(1+\frac{1}{\sqrt
q})^{2g}.$$

Since   $q\ge 4g^2$ we have   $$\# \cD(r, q)\le 2eq^{r}.$$
Assume that $G$ is a finite group
and $\psi$  an epimorphism of groups $$\psi : J (\Fq)\rightarrow
G.$$  We look for some divisor  $\Delta\in \cD(r,q)$ such that
$$\psi(\phi_r(\Delta))\not= 0\in G.$$ Since all the fibers of $\psi\circ \phi_r$
have the same cardinality, the fiber above $0$ has at most
$$\frac{2eq^{r}}{\#G}$$ elements.
So the number of prime divisors $\Delta\in \cP(r,q)$ such that
$\psi(\phi_r(\Delta))$ is not $0$ is at least
$$q^{r}(\frac{0.99}{r}-\frac{2e}{\# G}).$$ We assume that $\#G$ is
at least $12r$. Then at least half of the divisors in
$\cP(r,q)$ are not mapped onto $0$ by $\psi\circ \phi_r$.
The $\mu$-measure of the subset consisting of these elements
is at least $\frac{1}{2d}$.
So if we pick a random $\Delta$ in $\cP(r,q)$ with $\mu$-measure as in
Lemma~\ref{lemma:random}, the probability of 
success is at least $\frac{1}{2d}$. If we make $2d$ trials, the probability
of success is $\ge 1-\exp(-1)\ge \frac{1}{2}$.

\begin{lem}[Finding non-zero classes]\label{lemma:montecarlo}
There exists
  a probabilistic (Monte-Carlo) algorithm that takes as
input

\begin{enumerate}
\item   a degree $d$ and geometric genus $g$  plane projective 
absolutely integral  curve $C$ over
$\Fq$, such that $q\ge 4g^2$,

\item the smooth model $X$ of $C$,

\item  a degree $1$ divisor $O=O^+-O^-$ where $O^+$
and $O^-$ are effective, $\Fq$-rational and have degree bounded by 
$\Theta g^\Theta$ for some positive constant $\Theta$,

\item an epimorphism   $\psi : \Pic^0(X/\Fq)\rightarrow G$  (that need not be computable) such that the cardinality of
$G$ is at least $\max(48g,24d,720)$,
\end{enumerate}
\noindent 
and outputs  a  sequence of $2d$
elements in $\Pic^0(X/\Fq)$  such that at least one of them is not in
the kernel of $\psi$ with probability $\ge \frac{1}{2}$. The algorithm
is polynomial time in $d$ and $\log q$.
\end{lem}

As a special case we take $G=G_0=J (\Fq)$ and $\psi=\psi_0$ the identity.
Applying Lemma~\ref{lemma:montecarlo} we find  a sequence of elements
in $J (\Fq)$ out of which one at least is non-zero (with high
probability). We take  $G_1$
to be the quotient of $G$ by the subgroup generated by these elements and
$\psi_1$ the quotient map. Applying the lemma again we construct
another sequence of elements in $J (\Fq)$ out of which one at
least is not in $G_0$ (with high probability). We go on like that and produce a sequence of
subgroups in $J (\Fq)$
that increase with constant probability until the index in
$J (\Fq)$ becomes smaller than $\max(48g,24d,720)$. Note that
every step in this method is probabilistic: it succeeds with
some probability, that can be made very high (exponentially close
to $1$) while keeping a polynomial overall complexity.

\begin{lem}[Finding an almost  generating set]\label{lemma:generators}
There exists a probabilistic (Monte-Carlo) algorithm that takes as
input

\begin{enumerate}
\item   a degree $d$ and geometric genus $g$  plane projective 
absolutely integral  curve $C$ over
$\Fq$, such that $q\ge 4g^2$,

\item the smooth model $X$ of $C$,

\item  a degree $1$ divisor $O=O^+-O^-$ where $O^+$
and $O^-$ are effective, $\Fq$-rational and have degree bounded by 
$\Theta g^\Theta$ for some positive constant $\Theta$,
\end{enumerate}
\noindent
and outputs  a  sequence of 
elements in $\Pic^0(X/\Fq)$   that generate a subgroup of index at most 
$$\max(48g,24d,720)$$
\noindent 
with probability $\ge \frac{1}{2}$. The algorithm
is polynomial time in $d$ and $\log q$.
\end{lem}

Note that we do not catch the whole group $J(\Fq)$ of rational points but
a subgroup $\cG$ with  index at most
$\iota = \max ( 48g, 24d, 720)$. This is a small
but annoying gap. In the sequel we shall try
to compute the $l$-torsion of the group $J(\Fq)$
of rational points.
Because of the small gap in the above lemma, we may miss
some $l$-torsion points if $l$ is smaller than 
$\iota$. However, let $k$ be an integer
such that $l^k > \iota$. And let $x$ be a point
of order $l$
in $J(\Fq)$. Assume that there exists 
a point  $y$ in $J(\Fq)$ such that
$x=l^{k-1}y$. The group $<y>$ generated by $y$ and the group
$\cG$ have non-trivial intersection because the 
product of their orders is bigger than the order of
$J(\Fq)$. Therefore  $x$ belongs to $\cG$.

Our strategy for computing $J(\Fq)[l]$ will be to find
a minimal field extension $\FQ$ of $\Fq$ such that all points
in $J(\Fq)[l]$ are divisible by $l^{k-1}$ in $J(\FQ)$. We then
shall apply the above lemma to $J(\FQ)$. To finish
with,  we shall have to compute  $J(\Fq)$ as a subgroup
of $J(\FQ)$. To this end, we shall use the Weil pairing.

\section{Pairings}\label{section:pairings}

Let $n$ be a prime to $p$ integer and $J$ a jacobian variety over $\Fq$. 
 The Weil
pairing relates the full $n$-torsion subgroup $J({{\Fbar_q }})[n]$ with
itself. It can be defined using Kummer theory and is geometric in
nature.
The Tate-Lichtenbaum-Frey-R{\"u}ck pairing  is more arithmetic
and relates
the $n$-torsion $J(\Fq)[n]$ in the group of $\Fq$-rational points
and the quotient $J(\Fq)/nJ(\Fq)$. In this section, we quickly
review the definitions and algorithmic properties of these pairings,
following
work by Weil, Lang, Menezes,
Okamoto, Vanstone, Frey and R{\"u}ck.

We first recall the definition of
Weil pairing following \cite{langab}.
 Let
$k$ be an algebraically  closed field with characteristic $p$.
For every abelian variety $A$ over $k$, we denote by $Z_0(A)_0$
the group of $0$-cycles with degree $0$ and by $S : Z_0(A)_0\rightarrow A$ 
the summation map, that associates to every $0$-cycle of degree $0$
the corresponding sum in $A$.
Let $V$ and $W$ be two projective non-singular  integral
varieties over $k$, and let 
$\alpha : V \rightarrow A$ and $\beta : W \rightarrow B$ be the
 canonical
maps into their Albanese varieties. Let $D$ be a correspondence on
$V\times W$.
Let $n\ge 2$ be a prime to $p$ integer. Let
$\agot$ (resp. $\bgot$) be a $0$-cycle of degree $0$ on $V$
(resp. $W$) and let $a=S(\alpha(\agot))$ (resp. $b=S(\beta(\bgot))$) 
be the associated
point in $A$ (resp. $B$). Assume that $na=nb=0$. 
The  Weil
pairing $e_{n,D}(a,b)$ is defined 
in  \cite[VI, \S 4, Theorem 10]{langab}.
It is an $n$-th root of unity in $k$. It
is   linear in
$a$, $b$ and $D$.

Assume now that $V=W=X$ is a smooth projective  integral 
  curve over $k$ and assume that $A=B=J$ is its
jacobian and $$\alpha=\beta=\phi  :X \rightarrow J$$
is the Jacobi map (once  an origin on $X$ has been chosen). If we take $D$ to
be the diagonal on $X\times X$ we define a pairing $e_{n,D}(a,b)$
that
will be denoted $e_n(a,b)$ or $e_{n,X}(a,b)$. It does not depend on
the origin for the Jacobi map. It is non-degenerate.

The  jacobian $J$ comes with a  principal polarization i.e. an isomorphism
 $\lambda : J \rightarrow \hat J$
between $J$ and its dual $\hat J$. If
$\alpha$ is an  endomorphism $\alpha : J\rightarrow J$,
we denote by ${}^t\alpha$
 its transpose  ${}^t\alpha : \hat J
\rightarrow \hat J$. If $D$ is a divisor on $J$ that
is algebraically equivalent to zero, the image by ${}^t \alpha$
of the linear equivalence class of $D$ is the linear equivalence
class of the inverse image $\alpha^{-1}(D)$. See
\cite[V, \S 1]{langab}. The Rosati dual of $\alpha$ is defined to
be $\alpha^*=\lambda ^{-1}\circ {}^t\alpha \circ \lambda$. The 
map $\alpha \rightarrow \alpha^*$ is an involution, and 
$\alpha^*$ is the adjoint of $\alpha$ for the Weil pairing
\begin{equation}\label{equation:adjoint}
e_{n,X}(a,\alpha(b))=e_{n,X}(\alpha^*(a),b)
\end{equation}
\noindent according to \cite[VII, \S 2, Proposition 6]{langab}. 

If $Y$ is another smooth projective integral 
 curve over $k$ and $K$ its jacobian and
$ f  : X\rightarrow Y$ a non-constant map with degree $d$, and
$f^* : K\rightarrow J$ the associated map between jacobians,  then
for $a$ and $b$ of order dividing $n$ in $K$ one has
$$e_{n,X}(f^*(a),f^*(b))=e_{n,Y}(a,b)^d.$$

%\medskip 

The Frey-R{\"u}ck pairing can be constructed from the Lichtenbaum version
of Tate's
pairing \cite{lichtenbaum} as was shown in \cite{freyruck}. Let $q$ be a
power of $p$. 
Let  again $n\ge 2$ be a  prime to $p$    integer 
and let $X$  be  a smooth projective absolutely
integral  curve over  $\Fq$. Let $g$ be the genus of $X$.
We assume that $n$ divides $q-1$. 
%We assume $\cC$ has a rational point over $\Fq$. 
Let $J$ be the jacobian of $X$. The Frey-R{\"u}ck
pairing $$\{,\}_n: J (\Fq)[n]\times J (\Fq)/nJ (\Fq)\rightarrow
\Fqs/(\Fqs)^n$$
is defined as follows.  We take a class of order dividing $n$ in
$J (\Fq)$. Such a class can be represented by an
$\Fq$-divisor $D$ with degree $0$. We take  a class in $J (\Fq)$ and pick a degree zero
$\Fq$-divisor $E$  in this class, that we assume to be disjoint to $D$. The
pairing evaluated at the classes $[D]$ and $[E] \bmod n$ is $$\{[D],[E]\bmod n \}_n=f(E)\bmod {(\Fqs)^n}$$ where $f$
is any function with divisor $nD$.
This is a non-degenerate pairing. 

\smallskip 

We now explain how one can compute the Weil pairing, following work by
Menezes, Okamo\-to, Vanstone, Frey and R{\"u}ck. The
Tate-Lichtenbaum-Frey-R{\"u}ck
pairing can be computed similarly.
%The Weil pairing is computed  as follows. 
As usual, we assume that we are given a
degree $d$
plane model $C$ for $X$. 
Assume that $\agot$ and $\bgot$ have
disjoint
support (otherwise we may replace $\agot$ by some linearly
equivalent
divisor using the explicit moving Lemma~\ref{lemma:moving}.)
We compute  a function $f_\agot$  with divisor $n\agot$. 
We similarly compute  a function $f_\bgot$  with divisor $n\bgot$.
Then $$e_n(a,b)=\frac{f_\bgot(\agot)}{f_\agot(\bgot)}.$$ This algorithm is
polynomial in the degree $d$  of $C$ and the order $n$ of the divisors,
provided
the initial divisors $\agot$ and $\bgot$ are given as differences between effective
divisors with  polynomial degree  in $d$.

Using an idea that appears in a paper by Menezes, Okamoto and Vanstone
\cite{MOV} in the context of elliptic curves, and in \cite{freyruck}
for general curves, one can make this algorithm polynomial in $\log n$ in the
following way.
We write $\agot = \agot_0=\agot_0^+-\agot_0^-$ where $\agot_0^+$ and $\agot_0^-$ are effective
divisors. Let $f_\agot$  be the function computed in the above simple minded
algorithm. One has $(f)= n\agot_0^+-n\agot_0^-$. We want to express $f_\agot$ as a
product of small degree functions. We use a variant of fast
exponentiation. Using Lemma~\ref{lemma:moving} we compute a divisor 
$\agot_1= \agot_1^+-\agot_1^-$ and a function
$f_1$ such that $\agot_1$ is disjoint to $\bgot$ and
$(f_1)=\agot_1-2\agot_0$ and such that the degrees of $\agot_1^+$
and 
$\agot_1^-$ are
$\le 6gd (\log_q(\deg (\bgot))+1)$.
We go on  and compute, for $k\ge 1$ an integer, a divisor 
$\agot_k= \agot_k^+-\agot_k^-$ and a function
$f_k$ such that $\agot_k$ is disjoint to $\bgot$ and
$(f_k)=\agot_k-2\agot_{k-1}$ and such that the degrees of $\agot_k^+$
and 
$\agot_k^-$ are
$\le 6gd ( \log_q(\deg (\bgot)) +1)$.
We write the  base $2$ expansion of $n=\sum_k\epsilon_k2^k$  with
$\epsilon_k\in \{0,1\}$.
We compute the function $\Phi$ with divisor $\sum_k \epsilon_k\agot_k$.
We claim that the function $f_\agot$ can be written as a product of the
$f_k$, for $k\le \log_2n$,  and $\Phi$  with suitable integer exponents
bounded by $ n$ in absolute value. Indeed we write $F_1=f_1$, $F_2=f_2f_1^2$,
$F_3=f_3f_2^2f_1^4$ and so on. We have $(F_k)=\agot_k-2^k\agot$ and
$${\Phi^{}}{\prod_kF_k^{-\epsilon_k}}$$ has divisor $n\agot$ so is
the $f_\agot$ we were looking for.

\begin{lem}[Computing the Weil pairing]\label{lemma:weilp}
There exists an algorithm that on input a 
 prime to $q$     integer $n\ge
2$  and 
a degree $d$  absolutely integral 
 plane  projective curve $C$
over $\Fq$  and its  smooth
model $X$ 
and  two  $\Fq$-divisors on
 $X$,   denoted $\agot=\agot^+-\agot^-$ and $\bgot=\bgot^+-\bgot^-$, with
degree $0$,  and  order dividing $n$ in the jacobian,  computes the
Weil
pairing $e_n(\agot,\bgot)$ in time  polynomial 
in  $d$,  $\log q$,  $\log n$ and the degrees of $\agot^+$, $\agot^-$,
 $\bgot^+$, $\bgot^-$, the positive and negative parts of $\agot$ and
$\bgot$.
\end{lem}

\begin{lem}[Computing   Tate-Lichtenbaum-Frey-R{\"u}ck
    pairings]\label{lemma:tatep}
There exists an algorithm that on input an integer $n\ge
2$  dividing $q-1$ and 
a degree $d$  absolutely integral   plane  projective curve $C$
over $\Fq$  and its  smooth
model $X$ 
and  two $\Fq$-divisors on
 $X$,   denoted $\agot=\agot^+-\agot^-$ and $\bgot=\bgot^+-\bgot^-$, with
degree $0$,  and   such that the class of $\agot$ has order dividing $n\ge 2$ in the jacobian,  computes the
Tate-Lichtenbaum-Frey-R{\"u}ck 
pairing $\{\agot,\bgot\}_n$ in time polynomial
in  $d$,  $\log q$,  $\log n$ and the degrees of $\agot^+$, $\agot^-$,
 $\bgot^+$, $\bgot^-$, the positive and negative parts of $\agot$ and
$\bgot$.
\end{lem}

\section{Divisible groups}\label{section:divisiblegroups}

Let $\Fq$ be a finite field with characteristic $p$ and
let  $X$ be a projective smooth absolutely integral  algebraic
curve over $\Fq$. Let $g$ be the genus of $X$ and let $l\not = p$ be a
prime integer.  We assume that $g\ge 1$.
Let $J$ be the jacobian
of $X$ and let  $\End(J/\Fq)$ be the ring of
endomorphisms of $J$ over $\Fq$. Let $F_q$
be the Frobenius endomorphism.
In this section we study the action of $F_q$ on $l^k$-torsion
points of $J$. We first consider the whole $l^k$-torsion
group. We then restrict to some well chosen subgroups
where this action is more amenable.

Let  $\chi(x)$  be the characteristic polynomial  of 
$F_q \in  \End(J/\Fq)$.
The Rosati dual  to $F_q$  is
$q/F_q$. 
Let $$\cO = \ZZ[x]/\chi(x)$$ and
$\cO_l=\ZZ_l[x]/\chi(x)$.  We set $$\varphi_q=x\bmod
\chi(x)\in \cO.$$
Mapping  $\varphi_q$ onto $F_q$ defines an epimorphism from the ring
$\cO$  onto $\ZZ[F_q]$. 
In order to control the degree of the field of 
definition of $l^k$-torsion points we shall bound the order
of $\varphi_q$ in $(\cO/l^k\cO)^*$.
We set  $$\cU_1=(\cO/l\cO)^*=(\Fl[x]/\chi(x))^*.$$
Let the prime
factorization of $\chi(x)\bmod l$ be  $\prod_i\chi_i(x)^{e_i}$ with
$\deg(\chi_i)=f_i$. The order of $\cU_1$ is
$\prod_il^{(e_i-1)f_i}(l^{f_i}-1)$. Let $\gamma$ be the
smallest integer such that  $l^{\gamma}$ is bigger than or equal
to  $2g$. Then the exponent of the group $\cU_1$ divides 
$$A_1=l^{\gamma}\prod_i(l^{f_i}-1).$$
We set
$$B_1=\prod_i(l^{f_i}-1)$$
and $$C_1=l^{\gamma}.$$ There is a unique polynomial $M_1(x)\in \ZZ[x] $
with degree $< 2g$  such that $$\frac{\vp_q^{A_1}-1}{l}=M_1(\vp_q)\in
\cO.$$

Now  for every positive 
integer $k$,  the element  $\vp_q$ belongs to the unit
group $$\cU_k=(\cO/l^k\cO)^*$$   of the quotient algebra
$\cO/l^k\cO=\ZZ[x]/(l^k,\chi(x))$.
The prime
factorization of $\chi(x)\bmod  l$ is lifted modulo $l^k$ as
$\prod_i\Xi_i(x)$ with $\Xi_i$ monic and 
$\deg(\Xi_i)=e_if_i$,  and  the order of $\cU_k$ is
$\prod_il^{f_i(ke_i-1)}(l^{f_i}-1)$. The exponent of the latter
group divides
$$A_k=A_1l^{k-1}.$$ So we 
set $$B_k=B_1=\prod_i(l^{f_i}-1)$$
and $$C_k=C_1l^{k-1}=l^{k-1+\gamma}.$$ There is a unique polynomial
$M_k(x)
\in \ZZ[x] $
with degree $<  \deg(\chi)$  such that 
$$\frac{\vp_q^{A_k}-1}{l^k}=M_k(\vp_q)\in
\cO.$$

For every integer $N\ge 2$ we can compute $M_k(x)\bmod N$ from
$\chi(x)$ in probabilistic polynomial time
in $\log q$, $\log l$, $\log N$, $k$, $g$. Indeed   we first factor
$\chi(x)\bmod l$ then compute the $\chi_i$ and the $e_i$ and $f_i$. We compute 
$$x^{A_k} \bmod (\chi(x), l^kN)$$
using the fast exponentiation algorithm presented
in Section~\ref{sec:comp}. We remove $1$ and divide by $l^k$.

\begin{lem}[Frobenius and  $l$-torsion]\label{lemma:frobtorsion}
 Let $k$ be a
positive integer and $l\not = p$ a prime. Let $\chi(x)$ be the
characteristic polynomial of the Frobenius $F_q$ of $J/\Fq$.
Let $e_i$ and $f_i$ be the multiplicities and inertiae in the prime
decomposition
of $\chi(x)\bmod l$. Let $\gamma$ be the smallest integer such that
$l^\gamma$ is bigger than or equal to  $2g$. Let  
$B=\prod_i(l^{f_i}-1)$. Let 
 $C_k=l^{k-1+\gamma}$ and let $A_k=BC_k$. 
The $l^k$-torsion in $J$ splits completely over the degree $A_k$
extension of $\Fq$. There is a degree
$<  2g$ polynomial $M_k(x)\in \ZZ[x]$ such that
$$F_q^{A_k}=1+l^kM_k(F_q).$$
For every integer $N$ one  can compute such a $M_k(x)\bmod N$ from
$\chi(x)$ in probabilistic 
polynomial time 
in $\log q$, $\log l$, $\log N$, $k$, $g$.
\end{lem}

In order to state sharper results  it is convenient to introduce
$l$-divisible subgroups inside  the $l^\infty$-torsion of a
jacobian $J$,
that may
or may not correspond to subvarieties. We now see how to define such
subgroups and control their rationality properties.

\begin{lem}[Divisible group]\label{definition:divisible}
Let  $\Pi :
J[l^\infty]\rightarrow J[l^\infty]$ be  a group
homomorphism whose restriction to its image $\GG$  is a bijection. 
Multiplication by $l$ is then  a surjection from $\GG$ onto itself.
We denote by $\GG[l^k]$ the $l^k$-torsion in $\GG$.  There is an integer $w$ such that $\GG[l^k]$ is a
free $\ZZ/l^k\ZZ$ module of rank $w$ for every $k$.
We assume that  $\Pi$ commutes with the Frobenius endomorphism  $F_q$.
We then say $\GG$ is {\rm the
divisible group associated with $\Pi$}. From Tate's theorem \cite{tate}  $\Pi$ is  induced by some 
endomorphism in $\End (J/\Fq)\otimes_\ZZ\ZZ_l$ and  we can define
$\Pi^*$  the Rosati dual of 
$\Pi$  and denote by  $\GG^*=\IM(\Pi^*)$ the associated
divisible group,  that we  call  the adjoint  of $\GG$.
\end{lem}

\begin{remark}
The dual $\GG^*$ does not only depend on $\GG$. It may depend
on $\Pi$ also.  This will not be a problem for us.
\end{remark}

\begin{remark}
We may equivalently 
 define $\Pi^*$ as the dual of $\Pi$ for the Weil pairing.
See Equation~(\ref{equation:adjoint}).
\end{remark}

We now give an example of divisible group.
Let $F(x)=F_1(x)$ and $G(x)=G_1(x)$ be two
monic coprime polynomials in $\Fl[x]$ such that 
$$\chi (x)=F_1(x)G_1(x)\bmod l.$$
According to 
Bezout's theorem we have two polynomials $H_1(x)$ and $K_1(x)$ in
$\Fl[x]$  such that $$F_1H_1+G_1K_1=1$$ and $\deg(H_1)<\deg(G_1)$ and $\deg(K_1)<\deg(F_1)$. 

From Hensel's lemma, for every
positive integer $k$ there exist four polynomials $F_k(x)$, $G_k(x)$,
$H_k(x)$ and $K_k(x)$ in $(\ZZ/l^k\ZZ)[x]$ such that $F_k$ and
$G_k$ are monic and
$$\chi(x)=F_k(x)G_k(x)\bmod{l^k}$$ and
$$F_kH_k+G_kK_k=1\bmod{l^k}$$  and $\deg(H_k)<\deg(G_1)$ and
$\deg(K_k)<\deg(F_1)$ and  $F_1=F_k\bmod l$, $G_1=G_k\bmod l$,
$H_1=H_k\bmod l$, $K_1=K_k\bmod l$.

The sequences $(F_k)_k$, $(G_k)_k$, $(H_k)_k$, $(K_k)_k$  converge in
$\ZZ_l[x]$ to $F_0$, $G_0$, $H_0$, $K_0$. 
If we substitute $F_q$ for $x$ 
 in $F_0H_0$ we obtain  a map $$\Pi_G : J[l^\infty]\rightarrow
J[l^\infty]$$ and similarly, if we substitute  $F_q$  for
$x$ in
$G_0K_0$ we obtain  a map $\Pi_F$.  It is clear that $\Pi_F^2=\Pi_F$
and $\Pi_G^2=\Pi_G$ and $\Pi_F+\Pi_G=1$ and $\Pi_F\Pi_G=0$. We call
$\GG_F=\IM(\Pi_F)$ and $\GG_G=\IM(\Pi_G)$
the associated supplementary $l$-divisible groups.

\begin{defi}[Characteristic subspaces]\label{definition:characteristic}
For every non-trivial monic factor $F(x)$ of $\chi(x)\bmod l$ such that the
cofactor 
$G=\chi/F\bmod l$ is prime to $F$, we write  $\chi=F_0G_0$ the
corresponding factorization in  $\ZZ_l[x]$. The $l$-divisible group
$$\GG_F=\IM(\Pi_F)$$ is
called the $F_0$-torsion in $J[l^\infty]$ and is denoted 
$J[l^\infty,F_0]$. It is the characteristic subspace of $F_q$
associated with the factor $F$. If $F=(x-1)^e$ is the largest  power
of $x-1$ dividing $\chi(x)\bmod l$ we abbreviate
$\GG_{(x-1)^e}=\GG_1$. If  $F=(x-q)^e$ then 
we write  similarly $\GG_{(x-q)^e}=\GG_{q}$.
\end{defi}

We notice that there exists a unit $u$ in $\End (J/\Fq)\otimes_\ZZ\ZZ_l$ such that
the Rosati dual $\Pi_1^*$ of $\Pi_1$  is 

$$\Pi_1^*=\Pi_q \circ u.$$

Therefore  $$\GG_q=\GG_1^*$$ and
the restriction of the Weil pairing to $\GG_1[l^k]\times 
\GG_{q}[l^k]$  is non-degenerate for any integer  $k\ge 1$.

We now compute fields of definitions for torsion points inside
such divisible groups.
The action of $F_q$ on the $l^k$-torsion $\GG_F[l^k]=J[l^k,
F_0]$ inside $\GG_F$
factors through the  ring $\cO_l/(l^k,F_0(\varphi_q))=\ZZ_l[x]/(l^k,F_0)$. We
deduce the following.

\begin{lem}[Frobenius and  $F_0$-torsion]\label{lemma:frobtorsionGF}
 Let $k$ be a
positive integer and $l\not = p $ a prime. Let $\chi(x)$ be the
characteristic polynomial of the Frobenius $F_q$ of $J$. Let $\chi =FG\bmod l$ with $F$ and $G$
monic coprime.
Let $e_i$ and $f_i$ be the multiplicities and inertiae in the prime
decomposition of $F(x)\bmod l$. Let $\gamma$ be the smallest integer such that
$l^\gamma$ is bigger than or equal to  $2g$. Let  
$B(F)=\prod_i(l^{f_i}-1)$. Let $C_k(F)=l^{k-1+\gamma}$ and $A_k(F)=B(F)C_k(F)$. 
The $l^k$-torsion in $\GG_F$ splits completely over the degree $A_k(F)$
extension of $\Fq$. There is a degree
$<  \deg(F)$ polynomial $M_k(x)\in \ZZ_l [x]$ such that
$$\Pi_F \circ F_q^{A_k(F)}=\Pi_F+l^k\Pi_F\circ M_k(F_q).$$
For every power  $N$ of $l$,  one  can compute such an $M_k(x)$
modulo $N$ from
$\chi(x)$ and $F(x)$ in probabilistic 
polynomial time 
in $\log q$, $\log l$, $\log N$, $k$, $g$.
\end{lem}

If we take for $F$ the largest power of $x-1$ dividing $\chi(x)\bmod l$  in
the above lemma, we can take    $B(F)=1$ so
$A_k(F)$ is an $l$ power $\le 2gl ^{k}$.

If we take for  $F$ the largest power of $x-q$ dividing $\chi(x)\bmod l$ in
the above lemma, we have
 $B(F)= l -1$ so
$A_k(F)$ is  $\le 2g(l-1)l ^{k}$.

So the characteristic  spaces associated with the eigenvalues $1$ and
$q$ split completely over small degree extensions of $\Fq$.

\section{The Kummer map}\label{section:kummer}

Let $X$ be a smooth projective absolutely  integral   curve over
$\Fq$ of  genus $g$ and
$J$ the jacobian of $X$.
Let $n\ge 2$ be an integer dividing $q-1$. We assume that $g\ge 1$.
In this section, we construct a convenient surjection from
$J(\Fq)$ to $J(\Fq)[n]$.
If $P$ is in $J(\Fq)$ we take some 
$R\in J({{\Fbar_q }})$ such that $nR=P$ and form the $1$-cocycle
$({}^\sigma R -R)_\sigma$ in $H^1(\Fq, J [n])$. Using the Weil pairing we
deduce an element 
$$S \mapsto (e_n({}^\sigma R -R,S))_\sigma$$
 in
$$\Hom
(J [n](\Fq),H^1(\Fq, \mmu_n))=\Hom(J [n](\Fq),\Hom(\Gal(\Fq),\mmu_n)).$$

The map that sends $P\bmod nJ (\Fq)$ to $S\mapsto  (e_n({}^\sigma R -R,S))_\sigma$ is injective because the
Frey-R{\"u}ck pairing is non-degenerate. We observe that
$\Hom(\Gal(\Fq),\mmu_n)$
is isomorphic to $\mmu_n$ because  giving
an
 homomorphism from $\Gal(\Fq)$ to
$\mmu_n$ is equivalent to giving
the image of the Frobenius generator
$F_q$. We obtain a bijection $$T_{n,q} :
J (\Fq)/nJ (\Fq) \rightarrow 
\Hom(J [n](\Fq),\mmu_n)$$
 that we call the {\it Tate
map}. It maps  $P$ onto $S \mapsto e_n({}^{F_q} R -R,S)$. If
$J[n]$ splits completely over $\Fq$ we set $K_{n,q}(P)={}^{F_q} R -R$
and define a bijection $$K_{n,q} : J(\Fq)/nJ(\Fq) \rightarrow
J[n](\Fq)=J[n]$$ that we call the {\it Kummer map}.

\begin{defi}[The Kummer map] Let $J/\Fq$ be a jacobian and $n\ge 2$ a
prime to $p$
integer. Assume that $J[n]$ splits completely over $\Fq$. For $P$
in $J(\Fq)$ we choose any $R$ in $J({{\Fbar_q }})$ such that
$nR=P$ and we set $K_{n,q}(P)= {}^{F_q} R -R$.
This defines a bijection
$$K_{n,q} : J(\Fq)/nJ(\Fq) \rightarrow
J[n](\Fq)=J[n].$$
\end{defi}

We now assume that $$n=l^k$$ is a  power of some prime
integer $l\not = p $. We still   make the (strong~!) assumption that $J [n]$
splits completely over $\Fq$. 
We want to compute the Kummer map $K_{n,q}$
explicitly. Let $P$ be an
$\FF_{q}$-rational
point in $J$. Let $R$ be such that $nR=P$. %We observe that
				%$F_q^{n}$ fixes $R$. 
Since $F_q-1$ annihilates  $J [n]$, 
there is an $\Fq$-endomorphism $\kappa$
of $J$ such that $F_q-1=n\kappa$.  We note that $\kappa$ belongs to 
$\ZZ[F_q]\otimes_\ZZ\QQ=\QQ[F_q]$  and therefore commutes with $F_q$. We have 
$$\kappa
(P) = (F_q-1)(R)=K_{n,q}(P)$$ and $\kappa (P)$ is $\FF_{q}$-rational.
So we can compute $K_{n,q}(P)$ without computing $R$. We don't need
to  divide $P$ by $n$.

The Kummer map will show  very useful but its definition  requires that
$J [n]$ splits completely
over $\Fq$. If this is not the case,
we must base change to some extension
of $\Fq$.

Let $\chi(x)$ be the characteristic polynomial of $F_q$ and
let $B=\prod_i(l^{f_i}-1)$ where the $f_i$ are the degrees of the irreducible 
factors of $\chi(x)\bmod l$. Let $l^\gamma$ be the smallest power of $l$ that is
bigger than or equal to $2g$.
Let $C_k=l^{\gamma+k-1}$ and $A_k=BC_k$.   Set $Q=q^{A_k}$.
  From Lemma~\ref{lemma:frobtorsion} 
there is a polynomial $M_k(x)$ such that $$F_Q=1+l^kM_k(F_q).$$
So for $P$ an
$\FF_{Q}$-rational
point in $J$ and  $R$ such that $nR=P$,  the Kummer map  $K_{n,Q}$
applied to $P$ is  $$M_k(F_q)(P) = (F_Q-1)(R)=K_{n,Q}(P)$$ and this is
an  $\FF_{Q}$-rational point.

\begin{lem}[Computing the Kummer map]\label{lemma:kummer1}
Let $J/\Fq$ be a jacobian. Let $g\ge 1$ be its dimension. 
Let  $l\not= p$ be a
prime integer and $n=l^k$ a power of $l$. 
Let $\chi(x)$ be the characteristic polynomial of $F_q$ and
let $B=\prod_i(l^{f_i}-1)$ where the $f_i$ are the degrees of the irreducible
factors  of $\chi(x)\bmod l$. 
Let $l^\gamma$ be the smallest power of $l$ that is
bigger than or equal to $2g$.
Let $C_k=l^{\gamma+k-1}$ and $A_k=BC_k$.   Set $Q=q^{A_k}$ and  observe
that  $n$ divides $Q-1$ because $J[n]$ splits completely over $\FF_Q$.
There exists an endomorphism $\kappa \in \ZZ[F_q]$ of $J$ such that $n\kappa
=F_Q-1$
 and  for every $\FF_{Q}$-rational point $P$ and
any
$R$ with $nR=P$ one has  $\kappa
(P) = (F_{Q}-1)(R)=K_{n,Q}(P)$. This endomorphism $\kappa$ induces a bijection
between $J (\FF_{Q})/nJ (\FF_{Q})$ and $J [n](\FF_Q)=J [n]$. Given $\chi(x)$
and a positive 
integer $N$ one can compute $\kappa \bmod N$  as a polynomial in $F_q$
with coefficients in $\ZZ/N\ZZ$ in probabilistic polynomial time in
$g$, $\log l$,  $\log q$, $k$, $\log N$.
\end{lem}

This lemma is not of much use in practice because the 
 field  $\FQ$  is too big.
On the other hand, we may not be interested in the whole $n$-torsion in
$J$ but just a small piece in it, namely the $n$-torsion
of a given divisible group.

So let $l\not =
p$ be a prime integer and $\GG$ an $l$-divisible group in
$J[l^\infty]$. Let $$\Pi \in \End (J/\Fq)\otimes_\ZZ\ZZ_l$$
be a projection onto $\GG$. 
So $$\Pi : J[l^\infty]\rightarrow \GG.$$
We assume that $\Pi^2=\Pi$. 
Let   $n=l^k$ and let $Q$ be a power of $q$ such that $\GG[n]$
splits completely over $\FF_Q$. 
%We want to compute the Kummer map $K_{n,Q}$ restricted to
%$\GG(\FQ)$. 
Let $P$ be an
$\FQ$-rational
point in $\GG$. Let $R\in \GG({{\Fbar_q }})$ be such that $nR=P$.
We set $$K_{\GG,n,Q} (P)= {}^{F_Q}R-R$$ and define an isomorphism
$$K_{\GG,n,Q} : \GG(\FF_Q)/n\GG(\FF_Q) \rightarrow
\GG(\FF_Q)[n]=\GG[n].$$

In order to make this construction explicit, we now assume that there
exists some 
$\kappa \in \ZZ_l[F_q]$ such that 
$$(F_Q-1-n\kappa )\Pi=0.$$
Lemma~\ref{lemma:frobtorsionGF} provides us with such a $Q$ and such a
$\kappa$ when 
$\GG=J[l^\infty,F_0]$ is some
characteristic subspace.

We now can  compute this new  Kummer map $K_{\GG,n,Q}$. Let $P$ be an
$\FQ$-rational
point in $\GG$. Let $R\in \GG$ be such that $nR=P$. 
From $$(F_Q-1-n\kappa )\Pi(R)=0=(F_Q-1-n\kappa )(R)$$ we deduce that
$K_{\GG, n,Q}(P)=\kappa(P)$. Hence the  following lemma.

\begin{lem}[The Kummer map for a divisible group]\label{lemma:kummer2}
Let $J/\Fq$ be a jacobian. Let $g$ be its dimension.
Let $l\not = p$ be a
prime integer and $n=l^k$ a power of $l$. We assume that $g\ge 1$.
Let $\chi(x)$ be the characteristic polynomial of $F_q$. Assume
that $\chi(x)=F(x)G(x) \bmod l$ with $F$ and $G$ monic coprime polynomials in
$\Fl[x]$ and let $\GG_F$ be the associated divisible group.
Let $B=\prod_i(l^{f_i}-1)$ where the $f_i$ are the degrees of the 
irreducible factors  of $F(x)\bmod l$. 
 Let $l^\gamma$ be the smallest power of $l$ that is
bigger than or equal to $2g$.
Let $C_k=l^{k-1+\gamma}$ and $A_k=BC_k$. Set $Q=q^{A_k}$. From
Lemma~\ref{lemma:frobtorsionGF} there exists an endomorphism  $\kappa \in \ZZ_l[F_q]$ such
that   $$\Pi_F(n\kappa- F_Q+1)=0$$
 and  for every $\FF_{Q}$-rational point $P\in \GG_F$ and
any
$R\in \GG_F$ with $nR=P$ one has  $$\kappa
(P) = (F_{Q}-1)(R)=K_{\GG, n,Q}(P).$$ This endomorphism $\kappa$ induces a bijection
between $\GG_F (\FF_{Q})/n\GG_F (\FF_{Q})$ and $\GG_F [n](\FQ)=\GG_F [n]$. Given $\chi(x)$
and $F(x)$ and  a power 
$N$ of $l$,  one can compute $\kappa \bmod N$  as a polynomial in $F_q$
with coefficients in $\ZZ/N\ZZ$ in probabilistic polynomial time in
$g$, $\log l$,  $\log q$, $k$, $\log N$. Any
 $N\ge n(4Q)^g$ suffices
 for the purpose of computing $\kappa(P)$.
\end{lem}

\section{Linearization of torsion classes}\label{section:relations}

Let $C$ be a  degree $d$   plane projective absolutely integral  
curve $C$ over
$\Fq$ with geometric genus $g\ge 1$,
and assume that  we are given  the smooth model $X$ of $C$. 
We also assume that we are given  a degree $1$ divisor $O=O^+-O^-$ where $O^+$
and $O^-$ are effective, $\Fq$-rational and have degree bounded by 
 $\Theta g^\Theta$ for some 
 constant  $\Theta$.

Let 
$J$ be the  jacobian of $X$.
We  assume that $l\not = p$ is a prime integer that divides
$\#J(\Fq)$. Let $n=l^k$ be a power of $l$. We want
to describe $J(\Fq)[l^k]$ by generators and relations.

If  $x_1$, $x_2$, \dots, $x_I$ are  elements in  a finite commutative group $G$ we let $\cR$
be the kernel of the map $\xi : \ZZ^I\rightarrow G$ defined by 
$\xi(a_1,\cdots,a_I)=\sum_ia_ix_i$. We call $\cR$ the {\it lattice of
relations} between the $x_i$. 
We first give a very general and rough algorithm for computing
relations in any finite commutative group.

\begin{lem}[Finding relations in  blackbox groups]\label{lemma:idiot}
Let $G$ be a finite and commutative group and let $x_1$, $x_2$, \ldots,
$x_I$ be elements in $G$. A basis for the lattice of relations between
the $x_i$ can be
computed at the expense of $3I\#G$ operations (or comparisons) in $G$.
\end{lem}

We first compute and store all the multiples of $x_1$. So  we list $0$, $x_1$,
$2x_1$, \dots until we find the first multiple $e_1x_1$ that is equal to zero.
This gives us the relation $r_1=(e_1,0,\ldots,0)\in \cR$. This first
step requires at most $o=\#G$ operations in $G$  and $o$ comparisons.

We then compute successive multiples of $x_2$ until we find the first
one
$e_2x_2$ that is in $L_1=\{0,x_1,\ldots,(e_1-1)x_1\}$. This gives us a
second
relation $r_2$. The couple $(r_1,r_2)$ is a basis for the
 lattice of relations between $x_1$ and $x_2$. Using this lattice, we
 compute the list $L_2$ of elements in the group generated by $x_1$
 and $x_2$. This second
step requires at most $2o$ operations and $e_1e_2\le o$ comparisons.

We then compute successive multiples of $x_3$ until we find the first
one
$e_3x_3$ that is in $L_2$. This gives us a
third
relation $r_3$. The triple $(r_1,r_2,r_3)$ is a basis for the
 lattice of relations between $x_1$, $x_2$ and $x_3$. Using this lattice, we
 compute the list $L_3$ of elements in the group generated by $x_1$,
 $x_2$ and $x_3$. This third
step requires at most $2o$ operations and $o$ comparisons. And we go
on like this. \hfill
$\Box$

We stress that the algorithm above
 is far from efficient and will  not be very useful
 unless the group $G$ is very small.

\medskip 

We  now come back to the computation of generators and relations
for $J(\Fq)[l^k]$. 
Let $B=l-1$. Let $l^\gamma$ be
the smallest power of $l$  that  is bigger than or equal to $2g$
and let $A_k= Bl^{\gamma + k-1}$. We set $Q_k=q^{A_k}$. 
Let $F$ be  the largest
power of $x-1$ dividing the characteristic polynomial
$\chi(x)$ of $F_q$.  Definition~\ref{definition:characteristic}
and Lemma~\ref{lemma:kummer2}  provide us with
 two surjective maps $$\Pi_1 :
J(\FF_{Q_k})[l^\infty]\rightarrow \GG_1(\FF_{Q_k})$$ and
$$K_{\GG_1,l^k,Q_k} : \GG_1(\FF_{Q_k}) \rightarrow \GG_1[l^k].$$

If we now take for  $F$ the largest
 power of $x-q$ dividing $\chi(x)$, Definition~\ref{definition:characteristic}
and Lemma~\ref{lemma:kummer2}  give  two surjective maps $$\Pi_{q} :
J(\FF_{Q_k})[l^\infty]\rightarrow \GG_{q}(\FF_{Q_k})$$ and
$$K_{\GG_{q},l^k,Q_k} : \GG_{q}(\FF_{Q_k}) \rightarrow \GG_{q}[l^k].$$

Remember that $\GG_q=\GG_1^*$ and
the restriction of the Weil pairing to $\GG_1[l^k]\times 
\GG_{q}[l^k]$  is non-degenerate. We use this pairing to build
a presentation for $\GG_1[l^k]$ and $\GG_q[l^k]$ simultaneously.
The motivation for this approach  is that  generators for $\GG_1[l^k]$ provide
relations for $\GG_q[l^k]$,  and conversely.

If $Q_k\ge 4g^2$, we use  Lemma~\ref{lemma:generators} to
  produce a sequence
$\gamma_1$, \dots, $\gamma_I$ of
elements in  $J(\FF_{Q_k})$ that generate (with high probability)
a subgroup of index at
most
$$\iota = \max ( 48g, 24d, 720).$$

If $Q_k\le 4g^2$  we use Lemma~\ref{lemma:Stein} to produce
a sequence 
$\gamma_1$, \dots, $\gamma_I$ of
elements in  $J(\FF_{Q_k})$ that generate it.

 Let $N$ be the largest divisor of $\#J(\FF_{Q_k})$ which is prime to
$l$. 

We set $$\alpha_i=K_{\GG_1,l^k,Q_k}(\Pi_1(N\gamma_i))$$ and 
$$\beta_i=K_{\GG_{q},l^k,Q_k}(\Pi_{q}(N\gamma_i)).$$

The group $\cA_k$ generated by the $\alpha_i$ has index at most $\iota$
in
$\GG_1[l^k]$.
 The group $\cB_k$ generated by the $\beta_i$ has index at most $\iota$
in
$\GG_{q}[l^k]$.
Let $l^\delta$ be smallest power of $l$ that is bigger than
$\iota$ and assume that $k>\delta$.  Then $$\GG_1[l^{k-\delta}]
\subset \cA_k.$$

We now explain how to compute the lattice of relations  between
given elements $\rho_1$, \dots , $\rho_J$ in $\GG_1[l^k]$. We denote by
$\cR$ this lattice. We recall that 
the restriction of the Weil pairing  to $\GG_1[l^k]\times 
\GG_{q}[l^k]$  is a non-degenerate pairing
$$e_{l^k} : \GG_1[l^k]\times 
\GG_{q}[l^k]\rightarrow \mmu_{l^k}.$$

We fix an isomorphism 
between the group $$\mmu_{l^k}({{\Fbar_q }})=
\mmu_{l^k}(\FF_{Q_k})$$ of $l^k$-th roots of unity and
$\ZZ/l^k\ZZ$. 
Having  chosen the preimage of $1\bmod l^k$, computing
this isomorphism is a problem called {\it discrete logarithm}.
We can  compute this discrete logarithm
 by exhaustive search at the expense of $O(l^k)$ operations in $\FF_{Q_k}$.
There exist more efficient algorithms, but we don't need them for our
complexity
estimates.

We regard  the matrix $(e_{l^k}(\rho_j, \beta_i))$ as a matrix with
$I$ rows, $J$ columns and coefficients in $\ZZ/l^k\ZZ$. This
matrix defines a morphism from $\ZZ^J$ to $(\ZZ/l^k\ZZ)^I$ whose
kernel
is a lattice $\cR'$ that contains $\cR$. The index of $\cR$ in $\cR'$
is at most $\iota$. Indeed $\cR'/\cR$ is isomorphic to the orthogonal
subspace to  $\cB_k$ inside $<\rho_1,\ldots,\rho_J>\subset \GG_1[l^k]$. So it  has order $\le \iota$.
We then compute a basis of $\cR'$. This boils down to computing
the kernel of an $I\times (J+I)$ integer matrix with
entries bounded by $l^k$. This can be done
by putting this 
matrix in Hermite normal form (see \cite[2.4.3]{Cohen}). The complexity is polynomial in 
$I$, $J$ and $k\log l$. See \cite{Havas},  \cite[2.4.3]{Cohen} and \cite{vdk}.

Once given a basis of $\cR'$, the sublattice $\cR$ can be  computed
using Lemma~\ref{lemma:idiot} at the expense of $\le 3J\iota$ operations.

We apply this method to the generators $(\alpha_i)_i$ of $\cA_k$. Once
given the lattice $\cR$ of relations between the $\alpha_i$ it is a
matter of linear algebra to find a basis $(b_1,\dots,b_w)$  for
$\cA_k[l^{k-\delta}]=\GG_1[l^{k-\delta}]$.
The latter group is a rank $w$ free module over $\ZZ/l^{k-\delta}\ZZ$ and is acted
on by the $q$-Frobenius $F_q$. For every $b_j$ we can compute the
lattice
of relations between $F_q(b_j)$, $b_1$, $b_2$, \dots, $b_w$ and deduce
the matrix of $F_q$ with respect to
 the basis $(b_1,\dots,b_w)$. From this matrix 
we deduce
a nice generating set  for the kernel  of $F_q-1$ in $\GG_1[l^{k-\delta}]$. This
kernel is $J[l^{k-\delta}](\Fq)$. 
We deduce the following.

\begin{thm}[Computing the $l^k$-torsion  in the Picard group]\label{theorem:mainjacobi}
There is a probabilistic Monte-Carlo  algorithm that on  input 

\begin{enumerate}
\item   a degree $d$ and geometric genus $g$  plane projective 
absolutely integral  curve $C$ over
$\Fq$, 

\item the smooth model $X$ of $C$,

\item  a degree $1$ divisor $O=O^+-O^-$ where $O^+$
and $O^-$ are effective, $\Fq$-rational and have degree bounded by 
$\Theta g^\Theta$ for some positive constant $\Theta$,

\item a prime $l$ different from the characteristic $p$ of $\Fq$
  and a power $n=l^k$ of $l$,
\item the zeta function of $X$,
\end{enumerate}
outputs 
a set $g_1$, \dots, $g_I$ of divisor classes in   the
Picard group of $X/\Fq$,   such that
the $l^k$-torsion $\Pic(X/\Fq)[l^k]$ is the direct product
of the $<g_i>$,  and the orders of the $g_i$ form a non-decreasing
sequence.
Every class $g_i$ is given by a divisor $G_i-gO$ in the class, where $G_i$ is
a degree $g$  effective $\Fq$-divisor on $X$.

The algorithm runs in probabilistic polynomial time in $d$, 
$\log q$ and $l^k$. It outputs the correct answer with probability
$\ge \frac{1}{2}$. Otherwise, it may return either nothing or a strict
subgroup of $\Pic(X/\Fq)[l^k]$.

If one is given a degree zero $\Fq$-divisor $D=D^+-D^-$ of order
dividing $l^k$,
one can compute the coordinates of the class of $D$ in the basis
$(g_i)_{1\le i \le I}$ in polynomial time in $d$, $\log q$, $l^k$
and
the degree of $D^+$. These coordinates are integers $x_i$ such that $\sum_{1\le i\le I}x_ig_i=[D]$.
\end{thm}

\section{Computing $V_{\protect\lowercase{f}}$ 
modulo $\protect\lowercase{p}$}\label{section:ramanujan}

In this section, we apply the general algorithm
given in Section~\ref{section:relations}
to the plane curve $C_l$ constructed 
in Section~\ref {subsection:descripCl} and 
we compute Ramanujan divisors modulo $p$.
So we assume that we  are given an even integer $k > 2$,  
a prime integer $l>6(k-1)$, a finite field $\FF$
with characteristic $l$, and a ring epimorphism 
$f : \TT(1,k)\rightarrow  \FF$. More precisely, we 
are given the images $f(T_i)$ for $i\le k/12$.
We want to compute
the associated Galois representation
$V_f\subset J_1(l)$ or rather
its image $W_f\subset J_1(5l)$
by $B_{5l,l,1}^* : J_1(l) \rightarrow J_1(5l)$.
We will assume that the image of $\rho_f$ contains $\SL  (V_f)$. 
We set   $X_l=X_1(5l)$ and $J_l=J_1(l)$ and we denote by $g_l$ the genus
of $X_l$. 

Let $p\not\in \{5,l\}$ a prime integer.
We explain how to compute divisors
on $X_l/\Fp$ associated to every element $x$ in $W_f/\Fp\subset J/\Fp$.
The definition field $\Fq$ for such divisors  can be predicted from the
characteristic polynomial of the Frobenius endomorphism
$F_p$ acting on $V_f$. So the strategy is to
pick random $\Fq$-points in the
$l$-torsion of the  jacobian
$J_l$ and to project them onto $W_f$ using Hecke operators.

The covering map $B_{5l,l,1} : X_l\rightarrow X_1(l)$ has degree
$24$. We call it $\pi$.
It induces
two morphisms $\pi ^* : J_1(l)\rightarrow J_l$ and $\pi_* :
J_l\rightarrow J_1(l)$ such that the composite map 
$\pi_*\circ \pi^*$ is   multiplication by $24$ in  $J_1(l)$.
We denote by $\cA_l \subset J_l$ the image of 
$\pi^*$. This is a subvariety of $J_l$ isogenous to
$J_1(l)$. 
The restriction of $\pi^*\circ\pi_*$ to $\cA_l$ is multiplication by $24$.
The maps $\pi^*$  and $\pi_*$ induce Galois equivariant
bijections between
the
$N$-torsion subgroups $J_1(l)[N]$ and $\cA_l[N]$ for every integer $N$
which is prime to
$6$.

Using Theorem~\ref{thm_red_to_wt_2_Delta}
we derive from $f$ a finite  set $(t_1, \ldots, t_r)$ of elements
in  $\TT(l,2)$
with $r=(l^2-1)/6$ and 
\begin{equation}\label{eq:defVltor2}
V_f = \bigcap_{1\le i \le r} \ker\left(t_i, J_1(l)(\Qbar)[l]\right).
\end{equation}
and  $W_f\subset \cA_l\subset J_l$ is the image of $V_f$ by $\pi^*$. 

We  choose an integer    $s$
such that $24s$ is congruent
to $1$ modulo $l$. For every integer $n\ge 2$ we note $T_n \in \TT(l,2)$
the $n$-th Hecke operator with weight $2$ and level $l$. We can see
$T_n$ as endomorphism of $J_1(l)$. We
 set  $\hT_n=[s]\circ \pi^*\circ T_n\circ
\pi_*$. We  notice that $$\hT_n\circ \pi^*=\pi^*\circ T_n$$
on $J_1(l)[l]$. This
way, the map 
$\pi^* : J_1(l)\rightarrow
J_l$ induces   a Galois equivariant bijection of $\TT(l,2)$-modules
between   $J_1(l)[l]$ and  $\cA_l[l]$. And $W_f=\pi^*(V_f)$
 is the subspace in $\cA_l[l]$  cut out by all
$\hT_n-\tau(n)$. 
We notice that $\pi^*$, $\pi_*$, $T_n$, and $\hT_n$ can be seen as
correspondences
as well as morphisms between  jacobians. The following lemma states that
the Hecke action on divisors can be efficiently computed.

\begin{lem}[Computing the Hecke action]\label{lemma:hecke}
Let  $l$ and  $p$ be primes such that $p\not \in \{5,l\}$.
 Let $n\ge 2$
be an integer.  Let $q$ be a
power of $p$ and let $D$ be an effective $\Fq$-divisor of degree
$\deg (D)$ on $X_l \bmod p$. The divisor $\pi^*\circ \pi_*(D)$ 
  can be computed in  polynomial time in
  $l$, $\deg(D)$ and $\log q$. The divisor 
$\pi^*\circ T_n\circ \pi_*(D)$
  can be computed in  polynomial time in
  $l$, $\deg(D)$, $n$ and $\log q$.
\end{lem}

If $n$ is prime to $l$, we define
the Hecke operator $T(n,n)$ as an element in the ring of correspondences 
on $X_1(l)$ tensored by $\QQ$. See \cite[VII, \S 2 ]{langmf}.
From \cite[VII, \S 2, Theorem 2.1]{langmf} we have 
$T_{l^i}=(T_l)^i$ and $T_{n^{i}}=T_{n^{i-1}}T_n-nT_{n^{i-2}}T(n,n)$
if $n$ is prime and $n\not =l$.
And of course $T_{n_1}T_{n_2}=T_{n_1n_2}$ if $n_1$ and $n_2$ are
coprime.
So it suffices to explain how to compute $T_l$ and also 
$T_n$ and $T(n,n)$ for $n$ prime and
$n \not =l$.

Let  $x=(E,u)$ be  a point on $Y_1(l)\subset X_1(l)$
 representing an elliptic
curve $E$ with one $l$-torsion point $u$. 
Let $n$ be an integer. The Hecke operator $T_n$  maps
$x$ onto the sum of all $(E_I,I(u))$, where $I : E\rightarrow
E_I$ runs over the set of all  isogenies of degree $n$ from
$E$ such that $I(u)$ still has order $l$.  
If $n$ is prime to
$l$, the Hecke operator $T(n,n)$  maps $x$ onto
$\frac{1}{n^2}$ times $(E,nu)$. 
So we can compute the action of these  Hecke correspondences on points
$x=(E,u)$ using V{\'e}lu's formulae \cite{velu}.

There remains to treat the case of cusps.

We call $\sigma_{ \beta}$ for $1\le  \beta\le
\frac{l-1}{2}$ the cusp on $X_1(l)$ corresponding to
the $l$-gon equipped with an $l$-torsion point on the $\beta$-th component.
The corresponding Tate curves  $\CC^*/q$ have an  $l$-torsion
point $w=\zeta_l^\star q^{\frac{ \beta}{l}}$ where the star runs
over the set of all
 residues modulo $l$. There are $l$ Tate curves at every
 such cusp.

We call  $\mu_{ \alpha}$ for
$1\le {\alpha}\le \frac{l -1}{2}$ the cusps on $X_1(l)$
corresponding to a $1$-gon equipped with
the $l$-torsion point $\zeta_l^\alpha$ in its smooth locus $\Gm$.
The Tate curve at $\mu_\alpha$
is the Tate curves $\CC^*/q$ with $l$-torsion
point $w=\zeta_l^{ \alpha}$. One single Tate curve here:
no ramification.

For $n$ prime and  $n\not = l$ we have

$$T_n(\sigma_{\beta})=\sigma_{\beta}+n\sigma_{n\beta}$$
\noindent  and
$$T_n(\mu_{\alpha})=n\mu_{\alpha}+\mu_{n\alpha},$$
\noindent 
where $n \alpha$ in $\mu_{n \alpha}$  (resp. $n\beta$ in $\sigma_{n\beta}$)
 should be understood as a class in
$(\ZZ/l\ZZ)^*/\{ 1,-1\}$.

Similarly
$$T_l(\sigma_{\beta})=\sigma_{\beta}+2l\sum_{1\le
\alpha\le \frac{l -1}{2}}\mu_{\alpha}$$
\noindent  and
$$T_l(\mu_{\alpha})=l\mu_{\alpha}.$$

And of course, if $n$ is prime to $l$, then   $$T(n,n)(\sigma_\beta)=
\frac{1}{n^2}\sigma_{n\beta}$$ and $$T(n,n)(\mu_\alpha)=
\frac{1}{n^2}\mu_{n\alpha}.$$

All together, one can compute the effect of $T_n$ on
cusps for all $n$. 
For the sake of completeness, we also give the action of the diamond
operator $\langle  n\rangle $ on cusps. If $n$ is prime to $l$ then
$\langle n\rangle (\sigma_\beta)=\sigma_{n\beta}$ and $\langle n \rangle (\mu_\alpha)=\mu_{n\alpha}$.

\hfill $\Box$

We can now state the  following theorem.

\begin{thm}[Computing $V_f$ modulo $p$]\label{theorem:computingltorsionmodp}
There is a probabilistic (Las Vegas)  algorithm that 
takes as input 
an even integer $k > 2$,  a prime integer $l>6(k-1)$, a finite field $\FF$
with characteristic $l$,  a ring epimorphism 
$f : \TT(1,k)\rightarrow  \FF$,
 a cuspidal divisor $\Omega$ on $X_1(5l)$
 as constructed in Section~\ref{sec_constr_D},
and a prime $p\not\in \{5,l\}$, 
and computes the reduction modulo $p$
of  every element in
$W_f\subset J_1(5l)$. Here $V_f\subset J_1(l)$ is defined by 
Equation~(\ref{eq:defVltor2}) and $W_f$
  the image of $V_f\subset J_1(l)$ by
$B_{5l,l,1}^*$, and 
we  assume that 
the image of the Galois representation $\rho_f$ associated with $f$
contains $\SL  (V_f)$. The algorithm returns 
 for every element  $x$ in $W_f$ 
a  degree $g$ effective
divisor  $Q_{x,p}$ on $X_1(5l)/\Fp$ such that 
$Q_{x,p}-\Omega$ lies in the class represented by $x$ modulo $p$. 
The running time of the algorithm is $\le ( p\times \# V_f)^\Theta$ 
for some absolute constant $\Theta$.
\end{thm}

\begin{remark}
It has been proven
in Section~\ref{sec_constr_D}
that 
the  divisor $Q_x$ on $X_1(5l)/\QQ$ 
associated to $x\in W_f\subset J_1(5l)(\Qb)[l]$ is non-special.
According to 
Lemma~\ref{lem:jmcuni}
and~Theorem~\ref{thm_exist_good_places}, for every prime
$p$ but a finite number bounded by $l^\Theta$, 
the divisor $Q_x$ remains  non-special modulo
$p$ for every $x$ in $W_f$, and is equal to the divisor 
$Q_{x,p}$ returned by the algorithm above.
\end{remark}

To prove Theorem~\ref{theorem:computingltorsionmodp}
 we notice
that Section~\ref{subsection:descripCl} 
 gives us a  plane model for $X_l\bmod p$ and a resolution of its
singularities. From Lemma~\ref{lemma:manin2} we obtain the zeta function of
$X_l \bmod p$.  
The characteristic polynomial of $F_p$ acting on the $2$-dimensional
 $\FF$-vector
space  $V_f$
is $X^2-f(T_p)X+p^{k-1}\bmod l$.
Since we know $f(T_i)$  for $2\le i \le k/12$ we
deduce $f(T_p)$ using Manin-Drinfeld-Shokurov  theory.
 Knowing  the characteristic polynomial of $F_p$, we deduce
the order of $F_p$ acting on $V_f\bmod p$. We deduce some small enough
 splitting field  $\Fq$ for $V_f \bmod p$.
We then apply Theorem~\ref{theorem:mainjacobi} and obtain
a basis for the $l$-torsion in the Picard group of $X_l
/\Fq$. The same theorem allows us to compute the matrix of the
endomorphism $\pi^*\circ\pi_*$
in this basis. We deduce  a $\Fl$-basis 
for the image $\cA[l](\Fq)$ of $\pi^*\circ\pi_*$.
Using Theorem~\ref{theorem:mainjacobi} again, we now write down the matrices of the Hecke operators $\hT_n$ in this
basis
for all $n \le (l^2-1)/6$. It is then a matter of linear algebra to compute a
basis
for the intersection of the kernels of all $(t_i)_{1\le i\le (l^2-1)/6}$ in
$\cA[l](\Fq)$. The algorithm is Las Vegas rather than Monte-Carlo
because we 
can check the result,  the group  $W_f$ having known cardinality $(\#\FF)^2$.
\hfill $\Box$

%JMCfin
\chapter{Computing the residual Galois representations}
\label{sec_comp_mod_l_rep} 

\author{B. Edixhoven}

\bigskip

\bigskip

% authors Bas and Jean-Marc

In this chapter we first combine the results of
Chapters~\ref{chap_bnd_height} and~\ref{sec_couveignes_TORSION} in
order to work out the strategy of Chapter~\ref{chap_first_descr} in the
setup of Section~\ref{sec_setup_tau}. This gives the main result,
Theorem~\ref{thm_comp_rep_mod_l}: a deterministic polynomial time
algorithm, based on computations with complex numbers. The crucial
transition from approximations to exact values is done in
Section~\ref{sec_compute_Qzl-alg}, and the proof
of Theorem~\ref{thm_comp_rep_mod_l} is finished in
Section~\ref{sec_extract_rho}. In Section~\ref{sec_prob_algorithm} we
replace the complex computations with the computations over finite
fields from Chapter~\ref{sec_couveignes_modp}, and give a
probabilistic (Las Vegas type) polynomial time variant of the
algorithm in Theorem~\ref{thm_comp_rep_mod_l}.

\section{Main result}\label{sec_main_result}

% And explain that if $l>6(k-1)$ then \rho is  reducible or has image
% containing SL_2(\FF), Theorem~\ref{thm_large_image}.

For positive integers $k$ and $n$ we have defined, in
Section~\ref{sec_galreps}, $\TT(n,k)$ as the $\ZZ$-algebra in
$\End_\CC(S_k(\Gamma_1(n)))$ generated by the Hecke operators $T_m$
($m\geq 1$) and the~$\ld a\rd$ ($a$ in~$(\ZZ/n\ZZ)^\times$).
Theorem~\ref{thm_gen_Hecke} says that $\TT(n,k)$ is generated as
$\ZZ$-module by the Hecke operators $T_i$ with $1\leq i\leq
k{\cdot}[\SL_2(\ZZ):\Gamma_1(n)]/12$. In particular, $\TT(1,k)$ is
generated as $\ZZ$-module by the $T_i$ with $i\leq k/12$. For each $k$
and~$n$, each surjective ring morphism $\TT(n,k)\to\FF$ gives rise to
a Galois representation $\rho_m\colon\Gal(\Qbar/\QQ)\to\GL_2(\FF)$.
Theorem~\ref{thm_large_image} says that if $n=1$ and the
characteristic $l$ of $\FF$ satisfies $l>6(k-1)$, then $\rho_m$ is
reducible, or has image containing $\SL_2(\FF)$.

\begin{thm}\label{thm_comp_rep_mod_l}
There is a deterministic algorithm that on input a positive
integer~$k$, a finite field~$\FF$, and a surjective ring morphism
$f\colon\TT(1,k)\to\FF$ such that the associated Galois
representation $\rho\colon \Gal(\Qbar/\QQ)\to\GL_2(\FF)$ is reducible
or has image containing $\SL_2(\FF)$, computes $\rho$ in time
polynomial in $k$ and~$\#\FF$. The morphism $f$ is given by the images of\/
$T_1,\ldots,T_{\lfloor k/12\rfloor}$. More explicitly, the algorithm gives:
\begin{enumerate}
\item a Galois extension $K$ of\/ $\QQ$, given as a $\QQ$-basis $e$
and the products~$e_ie_j$ (i.e., the $a_{i,j,k}$ in $\QQ$ such that
$e_ie_j=\sum_k a_{i,j,k}e_k$ are given);
\item a list of the elements~$\sigma$ of~$\Gal(K/\QQ)$, where each
$\sigma$ is given as its matrix with respect to~$e$;
\item an injective morphism~$\rho$ from $\Gal(K/\QQ)$
  into~$\GL_2(\FF)$, making $\FF^2$ into a semi-simple
  representation of $\Gal(\Qbar/\QQ)$,
\end{enumerate}
such that $K$ is unramified outside~$l$, with $l$ the
characteristic of\/ $\FF$, and such that for all prime numbers
$p$ different from~$l$ we have:
\[
\trace(\rho(\Frob_p))=f(T_p)\quad\text{and}\quad
\det(\rho(\Frob_p))=p^{k-1}\quad \text{in\/ $\FF$.}
\]
\end{thm}
\begin{rem}
Of course, we do not only prove \emph{existence} of such an algorithm,
but we actually \emph{describe} one in the proof. We cannot claim that
we really \emph{give} such an algorithm because we did not make all
constants in our estimates explicit, e.g., those in
Theorem~\ref{thm_arakelov_contrib_2} and
Theorem~\ref{thm_exist_good_places}.
\end{rem}
\begin{proof}
As this proof is rather long, we divide it into sections.

\section{Reduction to irreducible representations}\label{sec_red_to_irr}
Let $k$, $\FF$, $l$ and~$f$ be as in
Theorem~\ref{thm_comp_rep_mod_l}. By definition, the associated
representation $\rho\colon\Gal(\Qbar/\QQ)\to\GL_2(\FF)$ is
semi-simple, unramified outside~$\{l\}$ and has
$\det\circ\rho=\chi_l^{k-1}$, with
$\chi_l\colon\Gal(\Qbar/\QQ)\to\FF_l^\times$ the mod~$l$ cyclotomic
character. Hence $\rho$ is reducible if and only if it is of the form
$\chi_l^i\oplus\chi_l^j$, for some $i$ and $j$ in $\ZZ/(l{-}1)\ZZ$
with $i{+}j=k{-}1$ in $\ZZ/(l{-}1)\ZZ$. The following result gives us
an effective way to decide if $\rho$ is reducible, and to determine
$\rho$ in that case. More precisely, using the standard algorithms
based on modular symbols
%%% reference????, to the beginning of the introduction??
this proposition reduces the proof of~Theorem~\ref{thm_comp_rep_mod_l}
to the case where $\rho$ is irreducible.
%% explain this in more detail??

\begin{prop}
In this situation, if\/ $l=2$, then $\rho\cong 1\oplus 1$, and if\/
$l=3$, then $\rho\cong 1\oplus\chi_3$. Assume now that $l\geq5$. Let
$i$ and $j$ be in $\ZZ/(l{-}1)\ZZ$ such that $i+j=k{-}1$. Then
$\rho$ is isomorphic to $\chi_l^i\oplus\chi_l^j$ if and only if for
all prime numbers $p\neq l$ with $p\leq (l^2{-}1)/12$ we have $f(T_p) =
p^i{+}p^j$ in\/~$\FF$.
\end{prop}
\begin{proof}
The statements about $l{=}2$ and $l{=}3$ are proved in Th\'eor\`eme~3
of~\cite{Serre10} (see also Theorem~3.4 of~\cite{Edixhoven3}).  As
$\TT(1,k)$ is not zero (it has $\FF$ as quotient), $k$ is even.

The idea for the rest of the proof is to use suitable cuspidal
eigenforms over $\FF_l$ whose associated Galois representations give
all the $\chi_l^i{\oplus}\chi_l^j$ with $i{+}j$ odd, and then to apply
Proposition~\ref{prop_test_equality_res_repr}. 

For $a$ in $\ZZ$ even such that $4\leq a\leq l{-}3$ or $a=l{+}1$, let
$\ol{E}_a$ be the element of $M_a(1,\FF_l)$ with $a_1(\ol{E}_a){=}1$
and $T_p(\ol{E}_a)=(1{+}p^{a-1})\ol{E}_a$ for all primes~$p$. These
$\ol{E}_a$ can be obtained by reducing the Eisenstein series $E_a$
modulo~$l$, after multiplication by $B_k/2k$ (see
Example~\ref{exam_eisenstein}). We cannot use $E_{l-1}$ because its
reduction has constant $q$-expansion~$1$, and so it cannot be
normalised as we need. The $\ol{E}_a$ are eigenforms, and we have
$\rho_{\ol{E}_a}=1{\oplus}\chi_l^{a-1}$. We note that the
$\chi_l^{a-1}$ give all powers of $\chi_l$ except~$\chi_l^{-1}$. But
we do have $1\oplus\chi_l^{-1} =
\chi_l^{-1}\otimes\rho_{\ol{E}_{l+1}}$. We also note that, for
each~$a$, $\ol{E}_a$ and $\theta^{l-1}\ol{E}_a$ give the same Galois
representation. Therefore, all $\chi_l^i{\oplus}\chi_l^j$ with $i{+}j$
odd are associated with suitable cupsidal eigenforms. The proof is then
finished by invoking Proposition~\ref{prop_test_equality_res_repr}.
\end{proof}

\section{Reduction to torsion in Jacobians}\label{sec_red_to_jac_2}
Let $k$, $\FF$, $l$ and~$f$ be as in
Theorem~\ref{thm_comp_rep_mod_l}, such that the representation $\rho$
attached to $f$ is irreducible. The following proposition is a
special case of Theorem~3.4 of~\cite{Edixhoven3}.
\begin{prop}
In this situation, there is a $k'$ in $\ZZ_{>0}$ and a surjective ring
morphism $f'\colon\TT(1,k')\to\FF$ and an $i$ in $\ZZ/(l{-}1)\ZZ$ such
that $2\leq k'\leq l{+}1$ and $\rho\cong\rho'\otimes\chi_l^i$, where
$\rho'\colon\Gal(\Qbar/\QQ)\to\GL_2(\FF)$ is the Galois representation
attached to~$f'$.
\end{prop}
Such a $k'$, $f'$ and~$i$ can be computed in time polynomial in $k$
and~$l$ as follows. First, one computes the $\FF_l$-algebras
$\FF_l\otimes\TT(1,k')$ for $2\leq k'\leq l{+}1$. Then, for each $i$
in $\ZZ/(l{-}1)\ZZ$ and for each~$k'$, one checks if there exists an
$\FF_l$-linear map from $\FF_l\otimes\TT(1,k')$ to $\FF$ that sends,
for all $m\leq (l^2{-}1)/12$ with $l$ not dividing~$m$, $T_m$ to
$m^{-i}f(T_m)$, and if so, if it is an algebra morphism. The previous
proposition guarantees that such $k'$, $f'$ and~$i$ do
exist. Proposition~\ref{prop_test_equality_res_repr} guarantees that
$\rho\cong\rho'\otimes\chi_l^i$.

Using standard algorithms for linear algebra over~$\QQ$, the computation
of $\rho$ is reduced to that of~$\rho'$. %%%%%%%%%%%%%%%%%%%%%%%%%%%%%
%% Say something about complexity here??? %%%%%%%%%%%%%%%%%%%%%%%%%%%%
By Theorem~\ref{thm_red_to_wt_2_Delta}, $\rho'$ is realised in
$J_1(l)(\Qbar)[l]$ as the intersection of the kernels of a set of
elements of $\FF_l\otimes\TT(l,2)$ that can be computed in time
polynomial in $k$ and~$l$. Hence, the proof of
Theorem~\ref{thm_comp_rep_mod_l} is reduced to the case where $\rho$
is irreducible, and is realised in $J_1(l)(\Qbar)[l]$ as described in
Theorem~\ref{thm_red_to_wt_2_Delta}.

\section{Computing the $\QQ(\zeta_{\protect\lowercase{l}})$-algebra
  corresponding to~$V$}\label{sec_compute_Qzl-alg}

We recall the situation.  The representation $\rho$ of
$\Gal(\Qbar/\QQ)$ attached to the surjective ring morphism
$f_k\colon\TT(1,k)\to\FF$ has image containing~$\SL_2(\FF)$, and is
also attached to a surjective ring morphism
$f_2\colon\TT(2,l)\to\FF$. In particular, $\rho$ is realised on the
two-dimensional $\FF$-vector space $V$ in $J_1(l)(\Qbar)[l]$
consisting of all elements annihilated by~$\ker(f_2)$. We note that
$l>5$. As in Section~\ref{sec_setup_tau} we let $X_l$ be the modular
curve $X_1(5l)$, over~$\QQ$, and we embed $V$ in the Jacobian $J_l$ of
$X_l$ via pullback by the standard map from $X_l$ to~$X_1(l)$. We take
a cuspidal divisor $D_0$ on $X_{l,\QQ(\zeta_l)}$ as in
Theorem~\ref{thm_exist_D}. We have, for each $x\in V$, a unique
effective divisor $D_x$ of degree $g_l$ (the genus of~$X_l$)
on~$X_{l,\Qbar}$, such that in $J_l(\Qbar)$ we have $x=[D_x-D_0]$. We
write each $D_x$ as $D_x^\fin+D_x^\cusp$, where $D_x^\cusp$ is
supported on the cusps and $D_x^\fin$ is disjoint from the cusps. We
write $D_x=\sum_{i=1}^{g_l}Q_{x,i}$, with $Q_{x,i}$ in $X_l(\Qbar)$,
such that $D_x^\fin=\sum_{i=1}^{d_x}Q_{x,i}$.

Theorem~\ref{theorem:torsion_main} says that we have an analytic
description of $V$ inside $J_l(\CC)$, and, that for every embedding
$\sigma\colon\QQ(\zeta_l)\to\CC$, complex approximations
$D_{\sigma,x}$ of the $D_x$ can be computed in polynomial time in
$\#\FF$ and the required accuracy (the number of digits on the right
of the decimal point). Each such approximation is given as a sum of
$g_l$ points $Q_{\sigma,x,i}$ in~$X_l(\CC)$, the numbering of which by
$i\in\{1,\ldots,g_l\}$ is completely arbitrary, i.e., unrelated to
each other when $\sigma$ varies. Similarly, complex approximations of
all $b_l(Q_{x,i})$ or $(1/b_l)(Q_{x,i})$ (one of which two has
absolute value~$<2$) and of all $x'_l(Q_{x,i})$ or $(1/x'_l)(Q_{x,i})$
can be computed in polynomial time in $\#\FF$ and the required
accuracy. We denote such approximations by $b_l(Q_{\sigma,x,i})$, etc.

We compute such approximations, and also of $j(Q_{x,i})$ or
$(1/j)(Q_{x,i})$, for all $\sigma$, $x$ and~$i$, with accuracy a
sufficiently large absolute constant
times~$l^{15}{\cdot}(\#\FF)^6$. Here, and the rest of this section, we
will use the $O$-notation without making the implied ``absolute''
constants explicit.

Using these approximations we will first decide for which
$(\sigma,x,i)$ the point $Q_{\sigma,x,i}$ approximates a cusp or
not. For $x=0$ we have $D_x=D_0$ and so all $Q_{0,i}$ are
cusps. Recall that the cusps are precisely the poles of the rational
function~$j$. Hence a necessary for $Q_{\sigma,x,i}$ to approximate a
cusp is that $(1/j)(Q_{\sigma,x,i})$ is small. Let $x$ in $V$ be
non-zero, and $i$ in $\{1,\ldots,g_l\}$ such that $j(Q_{x,i})\neq
0$. By Proposition~\ref{prop_arakelov_contrib_1}, the expression
for~$j$ in~$b$ in Proposition~\ref{prop_Y_1_5}, and
Lemma~\ref{lem_sum_and_product}, we have
$h((1/j)(Q_{x,i}))=O(l^{12})$. The degree of $(1/j)(Q_{x,i})$
over~$\QQ$ is at most $l^2{\cdot}(\#\FF)^2$. By
Lemma~\ref{lem_lower_bd_nonzero_elmt}, we have, for all
$\sigma\colon\Qbar\to\CC$, that if $(1/j)(Q_{x,i})\neq0$, then
$|\sigma((1/j)(Q_{x,i}))|\geq \exp(-O(l^{14}(\#\FF)^2))$. We conclude
that the $Q_{\sigma,x,i}$ for which
$|\sigma((1/j)(Q_{\sigma,x,i}))|<\exp(-O(l^{14}(\#\FF)^2))$ are the
ones that approximate cusps. This gives us the correct value of the
integers~$d_x$, and, after renumbering the~$Q_{\sigma,x,i}$,
approximations $D_{\sigma,x}^\fin=\sum_{i=1}^{d_x}Q_{\sigma,x,i}$ of
the $D_x^\fin$.

The next step is to get an integer $n$ with $0\leq n\leq
l^4{\cdot}(\#\FF)^4$ such that the function $f_l:=b_l+nx_l'$ separates
the various $Q_{x,i}$ with $x\in V$ and $i\leq d_x$ that are
distinct. We do this using just one embedding $\sigma$ of
$\QQ(\zeta_l)$ in~$\CC$. For $x$ in $V$ and $i\leq d_x$ we have, by
Proposition~\ref{prop_arakelov_contrib_1},
$h((b_l(Q_{x,i}),x_l'(Q_{x,i})))=O(l^{12})$. For $x$ and $y$ in $V$,
the field over which they are both defined has degree at most
$l{\cdot}(\#\FF)^3$ over~$\QQ$. By
Lemma~\ref{lem_lower_bd_nonzero_elmt}, we conclude that
$Q_{\sigma,x,i}$ and $Q_{\sigma,y,j}$ approximate the same point if
and only if
$|b_l(Q_{\sigma,x,i})-b_l(Q_{\sigma,y,j})|<\exp(-O(l^{13}{\cdot}(\#\FF)^3))$
and
$|x_l'(Q_{\sigma,x,i})-x_l'(Q_{\sigma,y,j})|<\exp(-O(l^{13}{\cdot}(\#\FF)^3))$.
We observe that the required approximations can indeed be computed
within the required time because the height bounds from
Proposition~\ref{prop_arakelov_contrib_1} imply that
$|b_l(Q_{\sigma,x,i})|<O(l^{14}{\cdot}(\#\FF)^2))$, i.e., also on the
left of the decimal point there are not too many digits. An integer
$n$ as above does not give a suitable $f_l$ if and only if there are
$Q_{\sigma,x,i}$ and $Q_{\sigma,y,j}$ that approximate different
points and still
$|f_l(Q_{\sigma,x,i})-f_l(Q_{\sigma,y,j})|<\exp(-O(l^{13}{\cdot}(\#\FF)^3)))$.
Trying the possible $n$ one by one until we have a suitable one gives
us a function $f_l$ as desired.

Now that we have our function~$f_l$, we continue, as explained at the
end of Section~\ref{sec_setup_tau}, by computing an integer $m$ with
$0\leq m \leq l^2{\cdot}(\#\FF)^4$ such that the function:
\[
a_m\colon V\to \Qbar, \quad 
x\mapsto \prod_{i=1}^{d_x}(m-f_l(Q_{x,i}))
\]
is injective, and hence a generator of the $\QQ(\zeta_l)$-algebra
$A_{\QQ(\zeta_l)}$ corresponding to~$V$. As in the previous step, we
do this using just one embedding $\sigma$ of $\QQ(\zeta_l)$
in~$\CC$. We estimate the loss in accuracy in computing a product
$\prod_{i\leq d_x}(m-f_l(Q_{\sigma,x,i}))$. For $x$ in $V$ and $i\leq
d_x$ we have $h(f_l(Q_{x,i}))=O(l^{12})$. Hence, in the product, we
have $|f_l(Q_{\sigma,x,i})|=O(l^{14}{\cdot}(\#\FF)^2)$. As there are
at most $l^2$ factors, the loss of accuracy is at most
$O(l^{16}{\cdot}(\#\FF)^2)$ digits. Hence, from our approximations
$f_l(Q_{\sigma,x,i})$ we get approximations $a_{\sigma,m}(x)$ of the
$a_m(x)$ with accuracy $O(l^{15}{\cdot}(\#\FF)^3)$. For all
candidates~$m$ and all $x$ in~$V$, that $h(a_m(x))=O(l^{14})$, and the
degree of $a_m(x)$ over~$\QQ$ is at most $(\#\FF)^2$. We conclude that
a candidate $m$ is not suitable if and only if there are distinct $x$
and $y$ in~$V$ with
$|a_{\sigma,m}(x)-a_{\sigma,m}(y)|<\exp(O(-l^{14}{\cdot}(\#\FF)^2))$. Trying
one by one gives us a suitable~$m$.

We denote by $a=a_{D_0,f_l,m}$ the generator of $A_{\QQ(\zeta_l)}$ that
we obtained by finding suitable $n$ and~$m$. We will now compute the
minimal polynomial of~$a$ over~$\QQ(\zeta_l)$:
\begin{eqn}\label{eqn_min_pol_a_4}
P=\prod_{x\in V}(T-a(x)), \quad P=\sum P_jT^j, \quad P_j\in\QQ(\zeta_l).
\end{eqn}
As $h(a(x))=O(l^{14})$, and, for $j$ in $\{0,\ldots,\# V\}$, $P_j$ is,
up to a sign, an elementary symmetric polynomial in the~$a(x)$,
Lemma~\ref{lem_bound_sym} gives that
$h(P_j)=O(l^{14}{\cdot}(\#\FF)^2)$. We write the $P_j$ in the
$\QQ$-basis $(1,\zeta_l,\ldots,\zeta_l^{l-2})$ of~$\QQ(\zeta_l)$:
\[
P_j=\sum_{i<l-1}P_{j,i}\zeta_l^i, \quad P_{j,i}\in\QQ.
\]
Then, for each~$j$, $(P_{j,0},\ldots,P_{j,l-2})$ is the unique
solution in $\Qbar$ of the system of linear equations, indexed by the
$\sigma$ in $\Gal(\QQ(\zeta_l)/\QQ)$:
\[
\sum_{i<l-1}\sigma(\zeta_l^i)P_{j,i} = \sigma(P_j).
\]
Applying Lemma~\ref{lem_bound_h_det} and Cramer's rule give:
\begin{eqn}\label{eqn_h_Pji}
h(P_{j,i}) \leq l\log l + l{\cdot}h(P_j) = O(l^{15}{\cdot}(\#\FF)^2).
\end{eqn}
So, in order to deduce the $P_{j,i}$ from approximations as in
Proposition~\ref{prop_cont_frac} the accuracy we need is
$O(l^{15}{\cdot}(\#\FF)^2)$. We estimate the loss of accuracy in the
evaluation of the product in~(\ref{eqn_min_pol_a_4}). We already know
that $h(a(x))=O(l^{14})$, and that the degree of $a(x)$ over~$\QQ$
is~$(\#\FF)^2$. Therefore, for all embeddings $\sigma$ of
$\QQ(\zeta_l)$ into~$\CC$, we have
$|a_{\sigma}(x)|=\exp(O(l^{14}{\cdot}(\#\FF)^2))$. As there are
$(\#\FF)^2$ factors in~(\ref{eqn_min_pol_a_4}), the loss of accuracy
is at most $O(l^{14}{\cdot}(\#\FF)^4))$. We conclude that our
approximations $P_\sigma=\prod_{x\in V}(T-a_\sigma(x))$ at all
$\sigma$ are accurate enough to get good enough approximations of the
$P_{j,i}$ such that Proposition~\ref{prop_cont_frac} gives us the
exact values of the~$P_{j,i}$. So, finally, we know $A_{\QQ(\zeta_l)}$
explicitly as:
\[
A_{\QQ(\zeta_l)} = \QQ(\zeta_l)[T]/(P) 
= \bigoplus_{0\leq i<(\#\FF)^2}\QQ(\zeta_l){\cdot}T^i.
\]
We remark that the definition of~$a$ directly implies that
$a(0)=1$. Hence $P$ has a factor $T-1$. Under our assumption that the
image of $\rho$ contains $\SL_2(\FF)$ the polynomial $P/(T-1)$ is
irreducible over~$\QQ(\zeta_l)$.

\section{Computing the vector space
  structure}\label{sec_comp_v-sp-str} The addition map $+\colon
V\times V\to V$ corresponds to a morphism of $\QQ$-algebras $+^*\colon
A\to A\otimes_\QQ A$, called co-adition. We will now explain how to
compute the co-addition over~$\QQ(\zeta_l)$, i.e., the morphism of
$\QQ(\zeta_l)$-algebras:
\[
+^*\colon A_{\QQ(\zeta_l)} \lto
A_{\QQ(\zeta_l)}\otimes_{\QQ(\zeta_l)} A_{\QQ(\zeta_l)} = 
\QQ(\zeta_l)[U,V]/(P(U),P(V)).
\]
To give this morphism is equivalent to give the image of our generator
$a$ of the previous section. This image can be written uniquely as a
polynomial in $U$ and $V$ of degree less than $(\#\FF)^2$ in each
variable. Hence, there are unique $\mu_{i,j}$ in $\QQ(\zeta_l)$, for
$i$ and $j$ in $\{0,\ldots,(\#\FF)^2{-}1\}$, such that for all $x$ and
$y$ in~$V$ we have, in~$\Qbar$:
\begin{eqn}\label{eqn_add_V}
a(x+y) = \sum_{i,j}\mu_{i,j}{\cdot}a(x)^ia(y)^j.
\end{eqn}
We view (\ref{eqn_add_V}) as an inhomogeneous system of $(\#\FF)^4$
linear equations in the~$\mu_{i,j}$. Then our bound
$h(a(x))=O(l^{14})$, together with Cramer's rule and
Lemma~\ref{lem_bound_h_det} give
$h(\mu_{i,j})=O(l^{14}{\cdot}(\#\FF)^6)$. Writing
$\mu_{i,j}=\sum_{k<l-1}\mu_{i,j,k}\zeta_l^k$, we have
$h(\mu_{i,j,k})=O(l^{15}{\cdot}(\#\FF)^6)$. Hence our approximations
$a_\sigma(x)$ are sufficiently precise to deduce the exact values of
the~$\mu_{i,j,k}$.

We also want to compute the multiplication map $\FF\times V\to
V$. That is, for each $\lambda$ in $\FF$ we want to know the map
$(\lambda{\cdot})^*$ from $A_{\QQ(\zeta_l)}$ to itself that it
induces. For $\lambda=0$ this is the map that sends $a$ to~$1$. Let
now $\lambda$ be in~$\FF^\times$. Then there are unique $\alpha_i$ in
$\QQ(\zeta_l)$, for $i$ in $\{0,\ldots,(\#\FF)^2{-}1\}$ such that
$(\lambda{\cdot})^*(a)=\sum_i\alpha{\cdot}a^i$. These $\alpha_i$
uniquely determined by the following equalities in~$\Qbar$, for all
$x$ in~$V$:
\begin{eqn}\label{eqn_mult_V}
a(\lambda{\cdot}x) = \sum_i \alpha_i{\cdot}a(x)^i.
\end{eqn}
Arguments as above for the addition show that our approximations
$a_\sigma(x)$ allow us to get the exact values of
the~$\alpha_i$.

\section{Descent to $\QQ$}\label{sec_descent_to_Q}
At this moment, we finally have to pay the price for working with a
divisor~$D_0$ on~$X_{l,\QQ(\zeta_l)}$ and not on~$X_l$ itself. We have
$A_{\QQ(\zeta_l)}=\QQ(\zeta_l)\otimes A$, hence we have a semi-linear
action of $\Gal(\QQ(\zeta_l)/\QQ)$ on the
$\QQ(\zeta_l)$-algebra~$A_{\QQ(\zeta_l)}$: for $\tau$ in
$\Gal(\QQ(\zeta_l)/\QQ)$, $\lambda$ in $\QQ(\zeta_l)$, and $x$ in
$A_{\QQ(\zeta_l)}$ we have $\tau(\lambda
x)=\tau(\lambda){\cdot}\tau(x)$. The $\QQ$-algebra $A$ is precisely
the subset of $A_{\QQ(\zeta_l)}$ of elements of that are fixed by this
action.

In order to understand what the action of $\Gal(\QQ(\zeta_l)/\QQ)$
does with our generator $a$ of $A_{\QQ(\zeta_l)}$, we must include the
divisor $D_0$ into its notation: we will write $a_{D_0}$ for it. Then,
for each $\tau$ in $\Gal(\QQ(\zeta_l)/\QQ)$ we have
$\tau(a_{D_0})=a_{\tau D_0}$, where $a_{\tau D_0}$ is defined as $a$,
but with the divisor $D_0$ replaced by~$\tau D_0$. For each
$\sigma\colon\QQ(\zeta_l)\to\CC$, and for each $x$ in $V\subset
J_l(\CC)$, we have the approximations $a_{\sigma, D_0}(x)$ of
$a_{D_0}(x)$, and $a_{\sigma\tau, D_0}(x)$ of~$a_{\tau D_0}(x)$.

We take a generator $\tau$ of
$\Gal(\QQ(\zeta_l)/\QQ)=\FF_l^\times$. There are unique $c_i$ in
$\QQ(\zeta_l)$, for $i$ in $\{0,\ldots,(\#\FF)^2{-}1\}$, such that
$\tau(a) = \sum_i c_i{\cdot}a^i$. These $c_i$ are uniquely determined
by the following system of equalities in~$\Qbar$, indexed by the $x$
in~$V$:
\begin{eqn}\label{eqn_def_ci_1}
a_{\tau D_0}(x) = \sum_{0\leq i<(\#\FF)^2} c_i{\cdot}a_{D_0}(x)^i.
\end{eqn}
Lemma~\ref{lem_bound_h_det} implies that for all $i$ we have
$h(c_i)=O(l^{14}{\cdot}(\#\FF)^6)$, and hence, writing
$c_i=\sum_{j<l-1}c_{i,j}\zeta_l^j$,
$h(c_{i,j})=O(l^{15}{\cdot}(\#\FF)^6)$. We conclude that our
approximations $a_{\sigma, D_0}(x)$ of
$a_{D_0}(x)$, and $a_{\sigma\tau, D_0}(x)$ of~$a_{\tau D_0}(x)$, for
all $x$ and $\sigma$, are sufficiently accurate to get the exact
values of the~$c_i$.

Linear algebra over~$\QQ$ gives us then $A$, in terms of a $\QQ$-basis
with multiplication table, and with the maps $+^*\colon A\to A\otimes
A$ and $(\lambda{\cdot})^*\colon A\to A$ that correspond to the
$\FF$-vector space structure on~$V$.

\section{Extracting the Galois representation}\label{sec_extract_rho}
We finish our computation of the Galois representation $\rho$ as
indicated in Chapter~\ref{chap_first_descr}. We view $V\times V$ as
$\Hom_\FF(\FF^2,V)$. This gives a right-action by $\GL_2(\FF)$ on
$V\times V$, hence a left-action on $A\otimes A$. This action can be
expressed in the co-addition and the $\FF^\times$-action. We let $B$
be the $\QQ$-algebra corresponding to the subset $\Isom_\FF(\FF^2,V)$
of $\Hom_\FF(\FF^2,V)$. To compute $B$, as a factor of~$A\otimes A$,
we compute its idempotent, i.e., the element of $A\otimes A$ that is
$1$ on $\Isom_\FF(\FF^2,V)$ and $0$ on its complement, as follows. In
$A_{\QQ(\zeta_l)}=\QQ(\zeta_l)[T]/((T-1)P_l)$ we have the idempotent
$a_1=P_1/P_1(1)$ which, as function on~$V$, is the characteristic
function of~$\{0\}$. Then $a_1$ is an element of~$A$. Let $a_2=1-a_1$
in~$A$. Then $a_2$ is the characteristic function of $V-\{0\}$. We let
$a_3$ be the element of $A\otimes A$ obtained by taking the product of
the $g{\cdot}(a_2\otimes 1)$, where $g$ ranges
through~$\GL_2(\FF)$. Then $a_3$ is the characteristic function of
$\Isom_\FF(\FF^2,V)$. We compute $B=(A\otimes A)/(1-a_3)$ by linear
algebra over~$\QQ$, in terms of a basis with a multiplication table,
and with the $\GL_2(\FF)$-action.

We factor the algebra $B$ as a product of fields, using a polynomial
time factoring algorithm over~$\QQ$ (see~\cite{LLL}
and~\cite{Lenstra_Hendrik_1}). Each factor $K$ of $B$ then gives us an
explicit realisation of~$\rho$, as explained in
Chapter~\ref{chap_first_descr}: let $G\subset \GL_2(\FF)$ be the
stabiliser of a chosen factor~$K$; then $G=\Gal(K/\QQ)$ and the
inclusion $\emph{is}$ a representation from $G$ to $\GL_2(\FF)$. This
finishes the proof of Theorem~\ref{thm_comp_rep_mod_l}.
\end{proof}
\begin{rem}
In the factorisation of the algebra $B$ above, our assumption that
$\im\rho$ contains $\SL_2(\FF)$ implies that the idempotents lie in
the sub-algebra $B^G$ of invariants by the subgroup $G$ of
$\GL_2(\FF)$ consisting of the $g$ with $\det(g)$ a $k{-}1$th power
in~$\FF_l^\times$. This subalgebra is a product of copies of~$\QQ$.

If $\FF=\FF_l$, factoring $B$ can be avoided by twisting $\rho$ by a
suitable power of~$\chi_l$. Indeed, $\rho':=\rho\otimes\chi_l^{1-k/2}$
has image $\GL_2(\FF_l)$, and $\rho$ can then be obtained as
$\rho'\otimes\chi_l^{k/2-1}$.
\end{rem}

\section{A probabilistic variant}\label{sec_prob_algorithm}
In this section we give a probabilistic Las Vegas type algorithm,
based on the results of Chapter~\ref{sec_couveignes_modp}, that
computes the representation $\rho$ as in
Theorem~\ref{thm_comp_rep_mod_l}, in probabilistic running time
polynomial in $k$ and~$\#\FF$. A nice feature of computations over
finite fields is that there is no loss of accuracy, as in the previous
sections where computations with complex numbers were used. On the
other hand, information obtained modulo varying primes is a bit harder
to combine, and here the point of view of Galois theory that we have
taken, relating sets with Galois action to algebras, is very
convenient.

Let $f\colon \TT(1,k)\to\FF$ be as in
Theorem~\ref{thm_comp_rep_mod_l}, as well as $l$ and the
representation
$\rho\colon\Gal(\Qbar/\QQ)\to\GL_2(\FF)$. Sections~\ref{sec_red_to_irr}
and~\ref{sec_red_to_jac_2} apply without any change. So we can now put
ourselves in the situation as in the beginning of
Section~\ref{sec_compute_Qzl-alg}. Then $\im\rho$ contains
$\SL_2(\FF)$, $\rho$ is realised on a two-dimensional $\FF$-vector
space in $J_l(\Qbar)[l]$, and for each $x$ in~$V$ there is a unique
effective divisor $D_x$ of degree $g_l$ on $X_{l,\Qbar}$ such that
$x=[D_x-D_0]$ in~$J_l(\Qbar)$. For each $x$ in~$V$, we write
$D_x=D_x^\fin+D_x^\cusp$ as before, with
$D_x^\fin=\sum_{i=1}^{d_x}Q_{x,i}$ and $D_x=\sum_{i=1}^{g_l}Q_{x,i}$.

Using Theorem~\ref{theorem:computingltorsionmodp}, we try to compute
the reductions $D_{x,\FF_q}$ over suitable extensions of the residue
fields of $\ZZ[\zeta_l]$ at successive prime numbers~$p$ up to
$c{\cdot}l^{15}{\cdot}(\#\FF)^6$, with $c$ a suitable absolute
constant, skipping $5$ and~$l$.  These fields $\FF_q$ have degree at
most $l{\cdot}(\#\FF)^3$ over their prime field~$\FF_p$.  If some
$D_{x,\FF_q}$ is not unique, this will be detected by our
computations, and we throw the corresponding prime~$p$ away. By
Theorem~\ref{thm_exist_good_places}, at most
$O(l^{12}{\cdot}(\#\FF)^2)$ primes~$p$ are thrown away. For the
$V$-good primes $p\leq B$, with $V$-good defined as in
Theorem~\ref{thm_exist_good_places}, we then have computed
all~$D_{x,\FF_q}$.

We split each such $D_{x,\FF_q}$ in a cuspidal part
$D_{x,\FF_q}^\cusp$ and a non-cuspidal part~$D_{x,\FF_q}^\fin$. By
Theorem~\ref{thm_exist_good_places}, there are at most
$O(l^{12}{\cdot}(\#\FF)^2)$ primes~$p$ where at some $\FF_q$, the sum
$\sum_{x\in V}\deg D_{x,\FF_q}^\fin$ is less than $\sum_{x\in
  V}d_x$. (In fact, as $\im\rho$ contains $\SL_2(\FF)$, all $d_x$ for
$x\neq0$ are equal, but our argument does not need this.) This means
that we have computed the unordered list of~$d_x$. We discard the
primes~$p$ where for some $\FF_q$ the sum $\sum_{x\in V}\deg
D_{x,\FF_q}^\fin$ is less than $\sum_{x\in V}d_x$. For the remaining
primes, we have computed the $D_{x,\FF_q}^\fin$ for all~$x$ in~$V$.

We want an integer $n$ such that $f_l:=b_l+nx_l'$ separates
the~$D_x^\fin$. Let $p$ be the smallest $V$-good prime with
$p>l^4{\cdot}(\#\FF)^4$. We compute the $D_{x,\FF_q}^\fin$ over some
$\FF_q$ at~$p$. We view the effective divisor $D_{\FF_q}:=\sum_{x\in
  V}D_{x,\FF_q}^\fin$ as a closed subscheme of~$X_{l,\FF_q}$. Then
$f_l$ as above embeds $D_{\FF_q}$ into $\AA^1_{\FF_q}$ if and only if
it is injective on the geometric points of $D_{\FF_q}$ and has
non-zero derivative at the multiple points of~$D_{\FF_q}$. As the
degree of $D_{\FF_q}$ is at most $l^2{\cdot}(\#\FF)^2$, this excludes
at most $l^4{\cdot}(\#\FF)^4$ elements of $\FF_p$, and the algorithms
of Chapter~\ref{sec_couveignes_modp} let us compute these in
polynomial time. We choose $n$ in $\{0,\ldots,l^4{\cdot}(\#\FF)^4\}$ such
that $f_l$ embeds $D_{\FF_q}$ into~$\AA^1_{\FF_q}$. Then, by
Nakayama's lemma, $f_l$ embeds $\sum_{x\in V}D_x^\fin$  into
$\AA^1_\Qbar$, and hence separates the~$D_x^\fin$.

The next step  is to compute an integer $m$ in
$\{0,l^2{\cdot}(\#\FF)^4\}$ such that the function:
\[
a_m\colon V\to \Qbar, \quad 
x\mapsto \prod_{i=1}^{d_x}(m-f_l(Q_{x,i}))
\]
is injective, and hence a generator of the $\QQ(\zeta_l)$-algebra
$A_{\QQ(\zeta_l)}$ corresponding to~$V$. Let $p$ and $\FF_q$ be as in
the preceding paragraph. For $x\in V$, write
$D_{x,\FF_q}^\fin=\sum_{i=1}^{d_x}Q_{x,i,\Fbar_q}$. We take $m$ in
$\{0,\ldots,l^2{\cdot}(\#\FF)^4\}$ such that the elements
$\prod_{i=1}^{d_x}(m-f_l(Q_{x,i,\Fbar_q}))$ of $\FF_q$, for $x\in V$, are all
distinct. Then $a:=a_m$ has the desired property.

The minimal polynomial $P$ in $\QQ(\zeta_l)[T]$ of $a$
over~$\QQ(\zeta_l)$ is given as:
\[
P = \prod_{x\neq 0}(T-a(x)), \quad P=\sum_j P_jT^j, \quad P_j\in\QQ(\zeta_l).
\]
We have seen in (\ref{eqn_h_Pji}) that, when writing $P_j =
\sum_{i<l-1}P_{j,i}\zeta_l^i$, we have
$h(P_{j,i})=O(l^{15}{\cdot}(\#\FF)^2)$. Our construction shows that
all $P_{j,i}$ are integral at all primes~$p$ that are $V$-good. From
the $D_{x,\FF_q}^\fin$ that we have computed, we get the images of the
$P_{j,i}$ in~$\FF_p$, for all $p$ that are $V$-good, up to
$c{\cdot}l^{15}{\cdot}(\#\FF)^2$. Proposition~\ref{prop_x_from_reductions}
then gives us the~$P_{j,i}$. So, at this point, we have computed the
$\QQ(\zeta_l)$-algebra $A_{\QQ(\zeta_l)}$.

To compute the $\FF$-vector space structure, we proceed as in
Section~\ref{sec_comp_v-sp-str}. For each $V$-good prime $p$ up to our
bound $c{\cdot}l^{15}{\cdot}(\#\FF)^6$, we compute at each $\FF_q$
over $p$ the images in $\FF_q$ of the $\mu_{i,j}$ from the linear
system~(\ref{eqn_add_V}), over~$\FF_q$, and then the $\mu_{i,j,k}$
in~$\FF_p$. Proposition~\ref{prop_x_from_reductions} and the height
bound $h(\mu_i,j,k)=O(l^{15}{\cdot}(\#\FF)^6)$ give us
the~$\mu_{i,j,k}$. The $\alpha_i$ as in~(\ref{eqn_mult_V}) can be
computed in similarly, as well as the $c_i$ as
in~(\ref{eqn_def_ci_1}). From here on, the computation then proceeds
as in Section~\ref{sec_extract_rho}.

\chapter{Computing coefficients of modular forms}\label{chap_comp_coefs}

\author{B. Edixhoven}

\bigskip

\bigskip

% authors Bas and Jean-Marc, if he wants.

In this chapter we apply the our main result on the computation of
Galois representations attached to modular forms of level one to the
computation of coefficients of modular forms. In
Section~\ref{sec_comp_tau} we treat the case of the discriminant
modular form, i.e., the computation of Ramanujan's $\tau$-function at
primes. In Section~\ref{sec_comp_tn} we deal with the more general
case of forms of level one and arbitrary weight~$k$, reformulated as
the computation of Hecke operators $T_n$ as $\ZZ$-linear combinations
of the $T_i$ with $i<k/12$. In Section~\ref{sec_quad_forms} we give an
application to theta functions of even, unimodular positive definite
quadratic forms over~$\ZZ$.

\section{Computing $\tau(\protect\lowercase{p})$ in time polynomial 
in $\log\protect\lowercase{p}$}
\label{sec_comp_tau}
We recall that Ramanujan's $\tau$-function is defined by the following
identity of formal power series with integer coefficients:
\[
x\prod_{n\geq1}(1-x^n)^{24} = \sum_{n\geq 1}\tau(n)x^n.
\]

\begin{thm}\label{thm_computation_tau}
There exists a deterministic algorithm that on input a prime
number~$p$ gives~$\tau(p)$, in running time polynomial in~$\log p$.
\end{thm}
\begin{proof}
Deligne has proved in \cite{Deligne1} and~\cite{Deligne2} that for all
prime numbers $p$ we have $|\tau(p)| < 2p^{11/2}$. Therefore, it
suffices to compute $\tau(p)$ mod~$l$ for all primes $l<x$, if the
product of these~$l$ is at least~$4p^{11/2}$. Analytic number theory
(see for example~\cite{Tenenbaum}, I.2.6, Corollary~10.1) tells us
that we can take $x=O(\log p)$, hence the proof is reduced to showing
that there is a deterministic algorithm that computes $\tau(p)$
mod~$l$ for prime numbers $p$ and~$l$ in time polynomial in $\log p$
and~$l$. Of course, the slightly weaker but much more elementary bound
$|\tau(p)|=O(n^6)$ in~\cite[Cor.~2.1.6]{Miyake1} also suffices for our
purposes.

We take an algorithm as in Theorem~\ref{thm_comp_rep_mod_l}, and we
apply it with $k=12$. We have $\TT(1,12)=\ZZ$, and for each $n$ in
$\ZZ_{\geq 1}$, the element $T_n$ of $\TT(1,k)$ is the
integer~$\tau(n)$. We must now show that from the output of the
algorithm we can compute $\tau(p)$ mod~$l$ deterministically, in time
polynomial in $\log p$ and~$l$. For $l$ prime, we let
$\rho_l\colon\Gal(K_l/\QQ)\into\GL_2(\FF_l)$ denote the Galois
representation attached to~$\Delta$. As the discriminant of the ring
of integers $A$ of~$K_l$ is a power of~$l$, Theorem~1.4
of~\cite{Buchmann-Lenstra} gives the existence of a deterministic
polynomial time algorithm that, given~$K_l$, produces~$A$, given by a
$\ZZ$-basis. The maximal order $A$ is preserved by the action
of~$\Gal(K_l/\QQ)$. Then $\Gal(K_l/\QQ)$ acts on the étale
$\FF_p$-algebra $\ol{A}:=A/pA$, and $\Hom(\ol{A},\Fbar_p)$ is a
$\Gal(K_l/\QQ)$-torsor. Moreover, $\ol{A}$ is the product of its
finitely many residue fields~$\ol{A}/m$, where $m$ ranges through the
maximal ideals of~$\ol{A}$.

We let $\Frob$ denote the absolute Frobenius endomorphism of~$\ol{A}$;
it sends $a$ to~$a^p$, it is an automorphism and it induces the
absolute Frobenius automorphism on each of the residue fields. The
matrix of $\Frob$ can be computed in time polynomial in~$l$ and~$\log
p$.

The Frobenius element $\sigma_m$ attached to a maximal ideal~$m$
of~$\ol{A}$ is the unique element $\sigma$ of $\Gal(K_l/\QQ)$ that
fixes $m$ and induces the absolute Frobenius on~$\ol{A}/m$. For
varying~$m$, the $\sigma_m$ form the Frobenius conjugacy class
(at~$p$) in~$\Gal(K_l/\QQ)$.

For each $\sigma$ in $\Gal(K_l/\QQ)$ we let $\ol{A}_\sigma$ be the
quotient of~$\ol{A}$ by the ideal generated by the image
of~$\Frob-\sigma$. Such an $\ol{A}_\sigma$ can be computed in
polynomial time. The $\sigma$ in the Frobenius conjugacy class are
precisely those $\sigma$ for which $\ol{A}_\sigma$ is non-zero. We can
try the $\sigma$ one by one until we have found a $\sigma$ in the
Frobenius conjugacy class. Then we apply the map $\rho_l$ given by
Theorem~\ref{thm_comp_rep_mod_l} to it, and get an element
$\rho_l(\sigma)$ of~$\GL_2(\FF_l)$. The trace of $\rho_l(\sigma)$ is
then $\tau(p)$ mod~$l$.
\end{proof}

%% In the proof above, should we move the explanation of how to
%% compute Frobenius elements without factoring to an earlier
%% section??? 

\section{Computing $T_{\protect\lowercase{n}}$ for 
large $\protect\lowercase{n}$ and large weight}\label{sec_comp_tn}

In this section we prove the following two results, for which we first
recall some notation.  For positive integers $k$ and $N$ we have
defined, in Section~\ref{sec_galreps}, $\TT(N,k)$ as the $\ZZ$-algebra
in $\End_\CC(S_k(\Gamma_1(N)))$ generated by the Hecke operators $T_n$
($n\geq 1$) and the~$\ld a\rd$ ($a$ in~$(\ZZ/N\ZZ)^\times$). By
Theorem~\ref{thm_gen_Hecke}, $\TT(N,k)$ is generated as $\ZZ$-module
by the Hecke operators $T_i$ with $1\leq i\leq
k{\cdot}[\SL_2(\ZZ):\Gamma_1(N)]/12$. Just before
(\ref{eqn_TS_pairing}) we have defined $M_k(\Gamma_1(N),\ZZ)$ as the
sub $\ZZ$-module of $M_k(\Gamma_1(N))$ consisting of the $f$ with all
$a_i(f)$ in~$\ZZ$.

\begin{thm}\label{thm_comp_tn}
One can compute, on input two positive integers $k$ and~$n$, and the
factorisation of $n$ into prime factors, the element $T_n$ of the
Hecke algebra $\TT(1,k)$, by computing the Galois representations
attached to sufficiently many maximal ideals of~$\TT(1,k)$ as in
Theorem~\ref{thm_comp_rep_mod_l}. The computation gives $T_n$ as
$\ZZ$-linear combination of the $T_i$ with $i\leq k/12$. This
algorithm is deterministic. For fixed $k$, it has running time
polynomial in $\log n$. If the generalised Riemann hypothesis (GRH)
holds for all number fields that are quotients of
$\QQ\otimes\TT(1,k)$'s, then the algorithm has running time polynomial
in $k$ and~$\log n$.
\end{thm}

\begin{cor}\label{cor_comp_an}
Assume GRH. There exists a deterministic algorithm that on input the
weight $k\geq0$ and the coefficients $a_i(f)$ with $0\leq i\leq k/12$
of a modular form $f$ in $M_k(\SL_2(\ZZ),\ZZ)$, and a positive
integer~$n$ together with its factorisation into primes, computes
$a_n(f)$ in running time polynomial in $k$, $\log n$, and the maximum
of the $\log(1+|a_i(f)|)$ with $i\leq k/12$.
\end{cor}
The principle of the proof of Theorem~\ref{thm_comp_tn} is first to
reduce the computation of $T_n$ to that of the $T_p$ for the primes
$p$ that divide~$n$, using the identities implicit
in~(\ref{eqn_dirichlet_hecke}). For $p$ prime, $T_p$ is computed from
its images in sufficiently many residue fields~$\TT(1,k)/m$, using the
LLL-algorithm for lattice reduction. Such an image is computed as the
trace of a Frobenius element at~$p$ of the Galois representation
$\rho_m$ attached to~$m$. The Galois representation is computed as in
Theorem~\ref{thm_comp_rep_mod_l}. The problem in doing all this is to
keep the residue fields small, because the computation of $\rho_m$
takes time polynomial in $k$ and~$\#(\TT(1,k)/m)$. Here it makes a big
difference if one assumes GRH or not. Assuming GRH, there are
suffiently many $m$ with $\TT(1,k)/m$ of size polynomial in $k$ and
$\log p$ so that $T_p$ can be reconstructed from its images in
these~$\TT(1,k)/m$. If we do not assume GRH, then we cannot rule out
the possibility that all ``small'' primes (small in terms of the
discriminant of~$\TT(1,k)$) are completely inert in~$\TT(1,k)$. Before
we give the proof of Theorem~\ref{thm_comp_tn}, where the details of
the algorithm are given, and the analysis of the running time, we
state and prove some preliminary results. Corollary~\ref{cor_comp_an}
will be deduced from Theorem~\ref{thm_comp_tn}, using some elementary
properties of the Eisenstein series~$E_k$, and the interpretation of
the $\ZZ$-module $S_k(\SL_2(\ZZ),\ZZ)$ as the $\ZZ$-dual
of~$\TT(1,k)$.

We start with a simple result that is well known.
\begin{prop}\label{prop_basis_T}
Let $k\geq 4$ be an even integer. Let $n=(k-14)/12$ if $k\equiv 2$
mod~$12$, and $n=\lfloor k/12\rfloor$ otherwise. Then $\TT(1,k)$ is
free of rank~$n$ as $\ZZ$-module, and the $T_i$ with $1\leq i\leq n$
form a $\ZZ$-basis for~$\TT$.
\end{prop}
\begin{proof}
Theorem~\ref{thm_gen_Hecke} tells us that $\TT:=\TT(1,k)$ is generated
as $\ZZ$-module by the $T_i$ with $1\leq i\leq k/12$, but we will not
use this. What we do use is that the pairing:
\[
S_k(\SL_2(\ZZ),\ZZ) \times \TT \lto \ZZ, \quad (f,t)\mapsto a_1(tf)
\]
is perfect (see~(\ref{eqn_TS_pairing})). What me must show is then
that $S_k(\SL_2(\ZZ),\ZZ)$ is free of rank~$n$ as $\ZZ$-module, and
that the maps $a_i\colon f\mapsto a_i(f)$, with $1\leq i\leq n$, form
a $\ZZ$-basis of the dual of~$S_k(\SL_2(\ZZ),\ZZ)$. As $\Delta$ is $q$
times a unit in~$\ZZ[[q]]$, we have:
\[
S_k(\SL_2(\ZZ),\ZZ) = \Delta{\cdot}M_{k-12}(\SL_2(\ZZ),\ZZ).
\]
According to Swinnerton-Dyer (\cite[\S3]{Swinnerton-Dyer1}):
\begin{eqn}
\label{eqn_mk_e4e6delta}
M(\SL_2(\ZZ),\ZZ) = \ZZ[E_4,\Delta] \oplus E_6{\cdot}\ZZ[E_4,\Delta],
\end{eqn}
with $E_4$ and $\Delta$ algebraically independent, and
$E_6^2=E_4^3-1728\Delta$. It follows that $S_k(\SL_2(\ZZ),\ZZ)$ is
free of rank~$n$. Suitable monomials in~$E_4$, $\Delta$ and~$E_6$ show
that $(a_1,\ldots,a_n)$ is a $\ZZ$-basis
of~$S_k(\SL_2(\ZZ),\ZZ)^\vee$.
\end{proof}
We note that we do not absolutely need the previous proposition,
because from the set of generators $T_i$ with $1\leq i\leq k/12$ as
given by Theorem~\ref{thm_gen_Hecke} one can also compute a
$\ZZ$-basis in time polynomial in~$k$, and this is what one will
probably do in the case of arbitrary level and weight. But in this case
of level one we have chosen to be more explicit.

We consider the Hecke algebras $\TT(1,k)$ as lattices in the
$\RR$-algebras $\TT(1,k)_\RR:=\RR\otimes\TT(1,k)$. As all $T_m$ are
selfadjoint as operators on $S_k(\SL_2(\ZZ))$ with respect to the
Petersson inner product (see Section~\ref{sec_modforms}), all their
eigenvalues are real. As the level is~$1$, all normalised eigenforms
are newforms; let us write them as $f_1,\ldots,f_n$, and view them as
ring morphisms $f_i\colon\TT(1,k)\to\RR$. Then the map:
\[
f\colon \TT(1,k)_\RR \lto \RR^n, \quad t\mapsto(f_1(t),\ldots,f_n(t)),
\]
is an isomorphism of $\RR$-algebras. The trace form on the
$\ZZ$-algebra $\TT(1,k)$ induces the trace form of the
$\RR$-algebra~$\RR^n$, i.e., the standard inner product on~$\RR^n$. We
equip $\TT(1,k)_\RR$ with the standard volume form, i.e., the one for
which a unit cube has volume~$1$. Our first goal is now to get a bound
for the absolute discriminant of~$\TT(1,k)$, or, equivalently, over
the volume of the quotient~$\TT(1,k)_\RR/\TT(1,k)$.

\begin{prop}\label{prop_bound_disc_T}
Let $k\geq 12$ be an even integer. Then we have:
\[
\log\Vol\left(\TT(1,k)_\RR/\TT(1,k)\right) 
= \frac{1}{2}\log|\discr\TT(1,k)| \leq \frac{k^2}{24}\log k.
\]
\end{prop}
\begin{proof}
As $T_1,\ldots,T_n$ is a $\ZZ$-basis of~$\TT(1,k)$, we have:
\[
\Vol(\TT(1,k)_\RR/\TT(1,k)) = |\det(f(T_1),\ldots,f(T_n))|.
\]
Using Deligne's bound $|a_p(f_j)| \leq 2p^{(k-1)/2}$ of~\cite{Deligne2}
as we used the Weil bounds in the proof of
Lemma~\ref{lem_bnd_coefficients} we get:
\begin{eqn}\label{eqn_bnd_ai_of_fj}
\begin{aligned}
|(f(T_i))_j| = |a_i(f_j)| & \leq \sigma_0(i)i^{(k-1)/2} \\
& \leq 2{\cdot}i^{1/2}i^{(k-1)/2}= 2{\cdot}i^{k/2}. 
\end{aligned}
\end{eqn}
Hence the square of the length of $f(T_i)$ is at most~$4ni^k$. This
implies:
\[
|\det(f(T_1),\ldots,f(T_n))| \leq \prod_{i=1}^n
\left(2{\cdot}i^{k/2}{\cdot}\sqrt{n}\right) 
= 2^n{\cdot}(n!)^{k/2}{\cdot}n^{n/2}. 
\]
Hence:
\begin{align*}
\log\Vol(\TT(1,k)_\RR/\TT(1,k)) & \leq 
n\log 2 + \frac{k}{2}\log(n!) + \frac{n}{2}\log n \\
& \leq n\log 2 + \frac{kn}{2}\log n -\frac{kn}{2} + \frac{k}{2} + 
\frac{n}{2}\log n,
\end{align*}
where we have used that $\log n! \leq n\log n -n+1$. Simple estimates
for $k>24$ and direct checks in the remaining cases give the result
(the only non-trivial case being $k=24$, in which case the
discriminant equals $2^6 3^2 144169$.
\end{proof}
The next ingredient to be used in our proof of
Theorem~\ref{thm_comp_tn} comes from analytic number theory: an
\emph{effective} prime number theorem for number fields, under the
assumption of GRH. The effectivity just alluded to means that for all
real numbers~$x$ in a specified interval such as $(2,\infty)$ an
estimate for the number of prime ideals of norm at most~$x$ of an
arbitrary number field~$K$ is given, whereas the usual prime number
theorem is, for each~$K$, a different asymptotic statement. 

\begin{thm}[Weinberger]\label{thm_weinberger}
For $K$ a number field and $x$ a real number let $\pi(x,K)$ denote the
number of maximal ideals $m$ of the ring of integers $O_K$ of $K$
with~$\#(O_K/m)\leq x$. For $x>2$ in~$\RR$ let $\li x = \int_2^x (1/\log
y)dy$. Then there exists $c_1$ in~$\RR$ such that for every number
field~$K$ for which GRH holds, and for every $x>2$ one has:
\[
|\pi(x,K)-\li x| \leq 
c_1\sqrt{x}\log\left(\left|\discr(O_K)x^{\dim_\QQ K}\right|\right).
\]
\end{thm}
Weinberger states this result in \cite{Weinberger2} (it is the lemma
on page~181, and says that it is proved in more generality
in~\cite{Weinberger1} (probably on page~328, just after~(4.4)). 
\begin{cor}\label{cor_pi1xKlarge}
There exist $c_2$ and $c_3$ in~$\RR$ such that for every number field
$K$ for which GRH holds and for every~$x$ in~$\RR$ such that:
\begin{align*}
x & > c_2{\cdot}(\log|\discr(O_K)|)^2{\cdot}(\log(1+\log|\discr(O_K)|))^2, 
\quad\text{and}\\
x & > c_3{\cdot}(\dim_\QQ K)^2{\cdot}(1+\log\dim_\QQ K)^4
\end{align*}
we have:
\[
\pi(x,K) \geq \frac{1}{2}\frac{x}{\log x}.
\]
\end{cor}
\begin{proof}
Let $K$ be a number field for which GRH holds, and let
$x\in\RR_{>2}$. We write $d_K$ for $|\discr(O_K)|$ and $n_K$ for
$\dim_\QQ K$. Then Theorem~\ref{thm_weinberger} says that:
\begin{align*}
\pi(x,K) & \geq \li(x) - c_1x^{1/2}\log(d_Kx^{n_K}) \\
& \geq \frac{x-2}{\log x} - c_1x^{1/2}\log(d_Kx^{n_K}) \\
& = \frac{x}{\log x}\left(
1-\frac{2}{x} - \frac{c_1(\log x)(\log d_K)}{x^{1/2}}
-\frac{c_1n_K(\log x)^2}{x^{1/2}}\right).
\end{align*}
In order to estimate the last term of the previous line, one uses the
substitutions $x=y^4$ and $y=n_K^{1/2}z$. For the one but last term,
one uses $x=y^2$ and $y=2c_1(\log d_K)z$ if $d_K\neq 1$.
\end{proof}
In order to find an element of $\TT(1,k)$ from sufficiently many of
its images modulo maximal ideals, we will need a lower bound on the
length of a shortest non-zero vector in the intersection of the
maximal ideals. The following lemma gives a general, well known lower
bound in the context of orders in number fields.
\begin{lem}\label{lem_lower_bnd_mu1}
Let $K$ be a number field, $n$ its dimension as $\QQ$-vector space,
$A\subset K$ an order, i.e., a subring of finite index in the ring of
integers of~$K$, and $I\subset A$ a non-zero ideal. We equip
$K_\RR:=\RR\otimes K$ with the inner product induced from the standard
inner product on $\CC^n$, where $K_\RR$ is embedded in $\CC^n$ via all
distinct $\sigma_i\colon K\to\CC$. We consider $A$ and $I$ as a
lattices in~$K_\RR$. Then we have:
\[
\mu_1(I) \geq \sqrt{n}{\cdot}\#(A/I)^{1/n},
\]
where $\mu_1(I)$ is the length of a shortest element of $I-\{0\}$.
\end{lem}
\begin{proof}
Let $x\in I$ be non-zero. Then $A{\cdot}x\subset I$, and we have:
\[
|x_1|\cdots|x_n| = |\rN_{K/\QQ}(x)| 
= \#(A/A{\cdot}x) \geq \#(A/I).
\]
The inequality between geometric and arithmetic mean gives:
\[
\left(|x_1|^2\cdots|x_n|^2\right)^{1/n}
\leq \frac{|x_1|^2+\cdots+|x_n|^2}{n} = \frac{\|x\|^2}{n}.
\]
Combining the last two inequalities gives:
\[
\|x\| \geq \sqrt{n}{\cdot}|\rN_{K/\QQ}(x)|^{1/n} \geq \#(A/I)^{1/n},
\]
which finishes the proof.
\end{proof}
The next result summarises the standard approach for using the
LLL-algorithm for the ``closest or nearest vector problem''. For more
details on lattice reduction we refer to~\cite{LLL},
\cite{Lenstra_Hendrik_2} and~\cite[\S2.6]{Cohen}.
\begin{prop}\label{prop_cvp_lll}
Let $n\geq 0$ and let $L$ be a free $\ZZ$-module of finite rank~$n$,
equipped with a positive symmetric bilinear form $b\colon L\times
L\to\ZZ$. We view $L$ as a lattice in the $\RR$-vector space $L_\RR$
equipped with the inner product given by~$b$, and for $x$ in $L_\RR$
we put $\|x\|:=(b(x,x))^{1/2}$. Let $L'$ be a submodule of finite
index of~$L$, and let $\mu_1(L')$ be the length of a shortest non-zero
element of~$L'$. Let $t$ be an element of~$L$ such that:
\[
\|t\| < 2^{-(n+1)/2} \mu_1(L').
\]
Let $e=(e_1,\ldots,e_n)$
be an ``LLL-reduced basis'' of~$L'$: if $e^*=(e_1^*,\ldots,e_n^*)$
denotes the orthogonal $\RR$-basis of $L_\RR$ obtained from $e$ by
letting $e^*_i$ be the orthogonal projection of $e_i$ to the
orthogonal complement of the subspace of $L_\RR$ generated by
$\{e_j\,|\, j<i\}$ (i.e., by the Gram-Schmidt orthogonalisation
process), and $\mu_{i,j}:=b(e_i,e_j^*)/b(e_j^*,e_j^*)$, then we have:
\begin{align*}
|\mu_{i,j}| & \leq \frac{1}{2}\quad \text{for $1\leq j<i\leq n$, and}\\
\|e_i^*\|^2 & \geq
\left(\frac{3}{4}-\mu_{i,i-1}^2\right)\|e_{i-1}^*\|^2
\quad \text{for $1<i\leq n$}.
\end{align*}
Then $t$ is the shortest element of $t+L'$, and for any $x$ in $t+L'$
we recover $t$ as follows:
%%%%%% Make a nice picture!!!!!!!
\begin{itemize}
\item put $x_n:=x$;
\item for $i$ going down from $n$ to $1$ let 
$x_{i-1}:=x_i-[b(x_i,e_i^*)/b(e_i^*,e_i^*)]e_i$, where, for $y$ in~$\QQ$,
$[y]$ denotes the largest of the (one or two) integers nearest to~$y$; 
\item then $t=x_0$.
\end{itemize}
\end{prop}
\begin{proof}
We claim that the orthogonal block:
\[
B:=\left\{\sum_i \lambda_i e_i^*\,|\,-1/2\leq \lambda_i<1/2\right\}
\]
is a fundamental domain for $L'$ acting on $L_\RR$ by
translations. Indeed, for $x$ in~$L_\RR$ and the $x_i$ as above, $x_0$
is in~$B$ and $x-x_0$ is in~$L'$, and moreover, $B$ and $L_\RR/L'$
have the same volume, namely~$\prod_i\|e_i^*\|$. By the defining
properties of an LLL-reduced basis we have, for $i$ in
$\{2,\ldots,n\}$, that $\|e_i^*\|^2\geq (1/2)\|e_{i-1}^*\|^2$. Also,
as $e_1^*=e_1$, we have $\|e_1^*\|\geq\mu_1(L')$. It follows that for
all~$i$:
\[
\|e_i^*\| \geq 2^{-(i-1)/2}\|e_1\| \geq 2^{-(i-1)/2}\mu_1(L'),
\]
and, in particular:
\[
\|e_i^*\| \geq 2^{-(n-1)/2}\mu_1(L').
\]
Hence, for any $x$ in~$L_\RR$ not in~$B$ we have:
\[
\|x\| \geq \min_i\|e_i^*\|/2 \geq 2^{-(n+1)/2}\mu_1(L').
\]
\end{proof}
We can now finally prove Theorem~\ref{thm_comp_tn}. We split it in
three parts: description of the algorithm, proof of its correctness,
and running time analysis.
\begin{proof}[of Theorem~\ref{thm_comp_tn}]
\subsection*{Description of the algorithm}
Let $k$ and~$n$ be given, with the factorisation of $n=\prod_p
p^{v_p(n)}$ in prime factors. Let $\TT:=\TT(1,k)$, and let $r$ be its
rank.  Proposition~\ref{prop_basis_T} tells us that $r\leq k/12$ and
that $(T_1,\ldots,T_r)$ is a $\ZZ$-basis of~$\TT$. Using modular
symbols algorithms (see Chapter~8 of~\cite{Stein}), one computes the
$\ZZ$-algebra structure of~$\TT$, i.e., one computes $T_iT_j$ as
linear combinations of the~$T_l$, with $1\leq i,j,l\leq r$. The
identity~(\ref{eqn_dirichlet_hecke}) gives us:
\[
T_n = \prod_{p|n}T_{p^{v_p(n)}}, \quad 
T_{p^i}= T_pT_{p^{i-1}}-p^{k-1}T_{p^{i-2}},
\]
where $p$ is prime and  $i\geq 2$. Hence the computation of $T_n$ is
reduced to that of the $T_p$ for the primes $p$ dividing~$n$. 

Let now $p$ be a prime dividing~$n$. Using a factorisation algorithm
as in~\cite{LLL}, we factor the $\QQ$-algebra $\TT_\QQ=\prod_i K_i$ as
a product of fields. In fact, in all cases that we know of, the number
of factors $K_i$ is at most one (see~\cite{Farmer-James}). For
each~$i$, let $A_i$ be the image of $\TT$ in~$K_i$. The $A_i$ are
computed as $\ZZ$-algebra, and the surjections $\TT\to A_i$ are
described by the images of the $T_j$ for $j\leq r$. These morphisms
embed $\TT$ into~$\prod_i A_i$. The computation of $T_p$ is reduced to
that of its images in the~$A_i$.

Let $K$ be one of the factors~$K_i$, let $A$ be~$A_i$, let
$n_K:=\dim_\QQ K$ and let $d_A:=|\discr(A)|$. We define:
\begin{eqn}\label{eqn_bound_for_ls}
B_A := c{\cdot}((\log d_A)^2+n_K^2 + n_K{\cdot}k{\cdot}\log p),
\end{eqn}
for $c$ a suitable absolute constant. For all primes $l$ with:
\[
6(k-1)<l<B_A{\cdot}(\log B_A)^4
\]
we compute the maximal ideals of~$A/lA$, and order them by their
norms: $\#(A/m_1A)\leq \#(A/m_2 A)$, etc. Then we take $j$ minimal
such that:
\begin{eqn}\label{eqn_l_bnd_ind_ideal}
\#(A/m_1\cdots m_j) > 
\left(2^{(n_K+1)/2}{\cdot}2{\cdot}p^{(k-1)/2}\right)^{n_K}.
\end{eqn}
For each of the $m_i$ with $i\leq j$, we compute the Galois representation
$\rho_{m_i}$ as in Theorem~\ref{thm_comp_rep_mod_l}, and a Frobenius
element $\rho_{m_i}(\Frob_p)$ at~$p$ as described in
Section~\ref{sec_comp_tau}. Then we have the images of $T_p$ in all
$\TT/m_i$ for $1\leq i\leq j$. 

We compute a $\ZZ$-basis $e'=(e'_1,\ldots,e'_{n_K})$ of~$m_1\cdots
m_j$, starting from the $\ZZ$-basis of~$A$ that is part of the
description of~$A$, adding the congruence condition modulo the~$m_i$'s
one by one. The LLL-algorithm gives us an LLL-reduced basis $e$
of~$m_1\cdots m_j$.

Let $\ol{T_p}$ denote the image of $T_p$ in $\prod_{i\leq j}A/m_i$. We
compute a preimage $T_p'$ in~$A$ of~$\ol{T_p}$, adding the congruence
conditions one by one. We compute $T_p$ itself with the algorithm of
Proposition~\ref{prop_cvp_lll}.

\subsection*{Correctness of the algorithm}
We will now show that the computation works, and that it gives the
correct result, if the constant $c$ in~(\ref{eqn_bound_for_ls}) is
large enough.

The first thing that we have to show is that there are distinct
maximal ideals $m_1,\ldots,m_j$ of $A$, whose residue characteristics
are between $6(k-1)$ and $B_A{\cdot}(\log B_A)^4$, not equal to~$p$,
and such that (\ref{eqn_l_bnd_ind_ideal}) holds. We note that for each
$l$ in that range there is at least one maximal ideal $m$ in $A/lA$,
and that $\#(A/m)$ is at least~$l$. According to Corollary~10.1 in
Section~I.2.6 of~\cite{Tenenbaum}, for sufficiently large real
numbers~$x$, the sum $\sum_{l\leq x}\log(l)$ is at least~$x/2$. So,
indeed, taking all $l$ between $6(k-1)$ and a constant times
$n_K{\cdot}k{\cdot}\log p$ implies that the required $m_i$ exist. We
observe that the choice of $B_A$ is sufficiently large for this.

The second point where an argument is needed is that the $\rho_{m_i}$
are either reducible, or have image containing $\SL_2(A/m_i)$, so that
they can be computed by the algorithm of
Theorem~\ref{thm_comp_rep_mod_l}. But this is guaranteed by
Theorem~\ref{thm_large_image}. 

The third point is that in the situation where we invoke the algorithm
of Proposition~\ref{prop_cvp_lll}, the assumptions of that proposition
hold. That means that we must check that $\|T_p\|$ is strictly less
than $2^{-(n_K+1)/2}{\cdot}\mu_1(m_1\cdots m_j)$. By Deligne's bound
of~\cite{Deligne2}, we know that $\|T_p\|\leq
\sqrt{n_K}{\cdot}2{\cdot}p^{(k-1)/2}$. Lemma~\ref{lem_lower_bnd_mu1}
says that $\mu_1(m_1\cdots m_j)\geq\sqrt{n_K}{\cdot}(\#(A/m_1\cdots
m_j))^{1/n_K}$. Hence the required inequality follows
from~(\ref{eqn_l_bnd_ind_ideal}). This finishes the proof of the
correctness of the algorithm.

\subsection*{Running time analysis}
We show the two claims on the running time. We also indicate at the
appopriate places that the algorithm is deterministic.

Let $k$, and $n=\prod_p p^{v_p(n)}$ be given. The computation of the
Hecke algebra $\TT:=\TT(1,k)$, i.e., of the products $T_iT_j$ as
$\ZZ$-linear combinations of the $T_l$, with $1\leq i,j,l\leq r$,
using modular symbols as in Chapter~8 of~\cite{Stein}, is
deterministic and is done in time polynomial in~$k$. Multiplication of
two elements of $\TT$ can be done in time polynomial in $k$ and the
maximum of the heights of the coordinates of the elements with respect
to the $\ZZ$-basis
$T=(T_1,\ldots,T_r)$. Lemma~\ref{lem_bound_h_sol_lin_eq} and
(\ref{eqn_bnd_ai_of_fj}) give that the heights of the coordinates of
$T_m$ (with $m\in\ZZ_{>0}$) with respect to $T$ are
$O(k^3(\log(km)))$. Hence the computation of $T_n$ from the $T_p$ for
$p$ dividing~$n$ is done in time polynomial in $k$ and~$\log n$. The
number of primes dividing~$n$ is at most $\log_2 n$, and each of them
is at most~$n$. Let $p$ be a prime dividing~$n$.

Factorisation as in~\cite{LLL} is a deterministic polynomial time
algorithm, hence we do get the factorisation $\TT_\QQ=\prod_i K_i$ and
the surjections $\TT\to A_i\subset K_i$ in time polynomial in~$k$. 

The number of factors $K_i$ is at most~$r\leq k/12$.  Let $A$ be one
of the~$A_i$. As the morphism $\TT\to\prod_i A_i$ is injective, we
have:
\[
\prod_i \discr(A_i) \leq \discr(\TT) \leq (k^2\log k)/12.
\]
Therefore, the number $B_A$ as in~(\ref{eqn_bound_for_ls}) is
$O(k^2{\cdot}((k\log k)^2+\log n))$. 

For $l$ a prime number, the computation of the maximal ideals of
$A/lA$ can be done deterministically in time polynomial in $l$
and~$n_K$ as follows (compare with Algorithm~3.4.10 of~\cite{Cohen},
``Berlekamp for small primes''). Take $i\in\ZZ_{\geq 1}$ such that
$l^i\geq n_K$. The product $I$ of the maximal ideals of $A/lA$ is the
subspace of $x$ with $x^{l^i}=0$. This can be computed with linear
algebra over~$\FF_l$. One computes, again with linear algebra, the
sub-algebra $(A/lA)'$ of elements $x$ with $x^l=x$. This algebra is
isomorphic to $\FF_l^d$ for some~$d$. If $d>1$, then for any $x$ in
$(A/lA)'$ that is not in~$\FF_l$ there is an $a$ in $\FF_l$ such that
$(x-a)^{l-1}$ is a non-trivial idempotent in~$(A/lA)'$. Repeating this
procedure splits $(A/lA)'$ completely. The maximal ideal corresponding
to an elementary idempotent $\iota$ of $(A/lA)'$ is then $(1-\iota)I$.

We conclude that the computation of all maximal ideals of the $A/lA$,
for all primes $l$ with $6(k-1)<l<B_A{\cdot}(\log B_A)^4$ can be done
deterministically in time polynomial in $k$ and $\log n$. We can order
them by their norms during their computation. Then we have
$m_1,\ldots,m_j$ such that (\ref{eqn_l_bnd_ind_ideal}) holds.  As
$\#(A/m_i)\geq 2$ for all~$i$, we have:
\[
j\leq n_K{\cdot}((n_K+3)/2 + (k-1)(\log_2 p)/2) = O(k^2\log n).
\]
Let $m$ be one of the $m_i$, where $i\leq j$. Then $\#(A/m)\leq
l^{n_K}$, where $l$ is the largest prime with $l<B_A{\cdot}(\log
B_A)^4$. A rough estimate gives that there is an absolute constant
such that:
\[
\#(A/mA)\leq (c{\cdot}k^5\log n)^{k/6}.
\]
The running time for computing the Galois representation $\rho_m$, is
therefore bounded by a fixed power of $(c{\cdot}k^5\log n)^{k/6}$. For
fixed $k$, this is of polynomial size in~$\log n$. Assuming GRH for
$K$, we get a much smaller estimate for $\#(A/m)$ as follows. We use
the fact that if $a$ is in $\RR_{>e}$, and $x>2a\log a$, then $x/\log
x > a$ (to prove this, use that $x\mapsto x/\log x$ is increasing on
$\RR_{>e}$, and that $e/\log e>2$).  We let $x=B_A{\cdot}(\log
B_A)^4$. Then $x$ satisfies the two hypotheses of
Corollary~\ref{cor_pi1xKlarge} (assuming the constant $c$
of~(\ref{eqn_bound_for_ls}) large enough with respect to $c_2$
and~$c_3$), and therefore (here we use GRH!)  $\pi(x,K)\geq x/2\log
x$. But then $\pi(x,K)\geq B_A/2$, i.e., there are at least $B_A/2$
maximal ideals of $O_K$ with $\#(O_K/m)\leq x$. Now $O_K/A$ is finite,
of order at most $d_A^{1/2}$, hence there are at most
$(1/2){\cdot}\log_2 d_A$ distinct primes dividing $\#(O_K/A)$, and
hence at most $(n_K/2){\cdot}\log_2 d_A$ maximal ideals of $A$ where
$A$ and $O_K$ differ. Similarly, there are at most $6(k-1)n_K$ maximal
ideals of~$A$ with residue characteristic at most~$6(k-1)$. The
definition of $B_A$ in equation~(\ref{eqn_bound_for_ls}) implies that:
\[
\frac{x}{2\log x} -\frac{n_K}{2\log_2 d_A} - 6kn_K >
n_K\left(\frac{n_K+1}{2} + 1 + \frac{k-1}{2}\log_2 p\right),
\]
assuming the absolute constant $c$ large enough. Let $I$ be the
product of the maximal ideals $m$ in $A$ with $\#(A/m)\leq x$ and with
residue characteristic at least~$6(k-1)$. Using that for all $m$
containing $I$ we have $\#(A/m)\geq 2$ we get:
\[
\#(A/I) > \left(2^{(n_K+1)/2}{\cdot}2{\cdot}p^{(k-1)/2}\right)^{n_K}.
\]
The definition of $j$, see~(\ref{eqn_l_bnd_ind_ideal}), implies that
for every $i\leq j$, we have $\#(A/m_i)\leq x$. As we have already
noticed above, $x$ is of size polynomial in $k$ and~$\log
n$. Summarising: if we assume GRH for~$K$, then running time for the
computation of $\rho_m$ is polynomial in $k$ and~$\log n$.

The computation of a Frobenius element $\rho_m(\Frob_p)$,
given~$\rho_m$, is deterministic and has running time polynomial in
$k$ and $\log n$; see the proof of Theorem~\ref{thm_computation_tau}. 

At this point we have computed the images of $T_p$ in the $A/m_i$, for
all $i\leq j$. The computation of a $\ZZ$-basis $e'$ of $m_1\cdots
m_j$ can be done in time polynomial in $k$ and $\log n$. One starts
with the basis of~$A$ that is part of its description, and computes
successively $\ZZ$-bases of $m_1$, $m_1m_2$, etc. At the $i$th step,
the maximum of the absolute values of the coordinates of the elements
of our temporary basis with respect to the basis of~$A$ gets at most
$n_Kl_i$ times larger, where $l_i$ is the characteristic
of~$A/m_i$. It follows that the absolute values of the coordinates of
the $e'_i$ with respect to the basis of~$A$ are at
most~$n_K^j{\cdot}\,\#(A/m_1\cdots m_j)$, which is of size polynomial in
$k$ and~$\log n$.

The LLL -algorithm gives us an LLL-reduced basis $e$ of $m_1\cdots
m_j$ in time polynomial in $k$ and $\log n$, see~\cite{LLL},
\cite{Lenstra_Hendrik_2} or~\cite[\S2.6]{Cohen}. 

A preimage $T_p'$ in $A$ of the image $\ol{T_p}$ of $T_p$ in
$\prod_{i\leq j}A/m_i$ can be computed in time polynomial in $k$ and
$\log n$ as follows. One lifts the image of $T_p$ in $A/m_1$ to an
element $T_{p,1}$ of~$A$, with small coordinates with respect to the
$\ZZ$-basis of~$A$, then one adjusts $T_{p,1}$ with an element
of~$m_1$ with small coordinates with respect to the $\ZZ$-basis that
was already computed, to get the correct image in $A/m_2$, etc. The
size of the coordinates of $T_p'$ with respect to the $\ZZ$-basis
of~$A$ is polynomial in~$k$ and~$\log n$.

The algorithm of Proposition~\ref{prop_cvp_lll} then computes $T_p$
for us, in time polynomial in $k$ and~$\log n$.
\end{proof}

\begin{proof}[of Corollary~\ref{cor_comp_an}]
Let $k\geq0$, and let $f$ be in $M_k(\SL_2(\ZZ),\ZZ)$. We are given
the integers $k$ and the $a_i(f)$ for $0\leq i\leq k/12$, and a
positive integer $n$ together with its factorisation into primes. We
are to compute~$a_n(f)$. 

If $M_k(\SL_2(\ZZ),\ZZ)$ is zero, then so is~$f$, and we have
$a_n(f)=0$. If $k=0$, then $f$ is a constant power series, hence
$a_n(f)=0$. So, we may and do assume that $k\geq4$, and that $k$ is
even.

Let $E_k$ be the Eisenstein form of weight~$k$ and level one
normalised as in Example~\ref{exam_eisenstein}. Note that the
coefficients $a_i(E_k)$ are rational numbers, not necessarily
integers. Let $g:=f-a_0(f){\cdot}E_k$ in $S_k(\SL_2(\ZZ),\QQ)$. Then
we have:
\[
f = a_0(f){\cdot}E_k + g, \quad\text{hence}\quad 
a_n(f) = a_0(f)a_n(E_k) + a_n(g).
\]
As the factorisation of $n$ into primes is given, we can compute the
coefficient:
\[
a_n(E_k)=-\frac{2k}{B_k}\sigma_{k-1}(n) = 
-\frac{2k}{B_k}\prod_{p|n}\left(1+p^{k-1}+\cdots+p^{(k-1)v_p(n)}\right)
\]
in time polynomial in~$k$ and~$\log n$. Hence it remains to
compute~$a_n(g)$.

Let $\TT$ denote the Hecke algebra $\TT(1,k)$ acting
on~$S_k(\SL_2(\ZZ),\ZZ)$. Recall that its rank is at most~$k/12$. Via
the perfect pairing from~(\ref{eqn_TS_pairing}):
\[
S_k(\SL_2(\ZZ),\ZZ) \times \TT \lto \ZZ, \quad (h,t)\mapsto a_1(th)
\]
we view $g$ as the element in $\Hom_{\text{$\ZZ$-Mod}}(\TT,\QQ)$
that sends $t$ to~$a_1(tg)$. We then have, by~(\ref{eqn_a1Tnf}):
\[
a_r(g) = g(T_r) \quad \text{for all $r\geq 1$}.
\]
The $T_i$ with $1\leq i\leq \rank(\TT)$ form a $\ZZ$-basis for~$\TT$, by
Proposition~\ref{prop_basis_T}. The products $T_iT_j$ of these basis
elements can be computed in time polynomial in~$k$.

The coefficients $a_i(E_k)$ for $0\leq i\leq k/12$ can be computed in
time polynomial in~$k$ from the standard formulas as given in
Example~\ref{exam_eisenstein}. The coefficients $a_i(f)$ for $0\leq
i\leq k/12$ are given. Hence the $g(T_i)=a_i(g)$ with $1\leq i\leq
\rank(\TT)$ are computed in time polynomial in~$k$ and the maximum of
the $\log(1+|a_i(f)|)$ with $0\leq i\leq
k/12$. Theorem~\ref{thm_comp_tn} says that $T_n$ can be computed as
$\ZZ$-linear combination of the $T_i$ with $1\leq i\leq\rank(\TT)$, in
time polynomial in~$k$ and~$\log n$. Applying $g$, viewed as element
of $\Hom_{\text{$\ZZ$-Mod}}(\TT,\QQ)$, to this linear combination
gives $a_n(g)=g(T_n)$, in the required time.
\end{proof}

\section{An application to quadratic forms}
\label{sec_quad_forms}
In this section we apply our results on the computation of
coefficients of modular forms to theta functions of even, unimodular
positive definite quadratic forms over~$\ZZ$. According to
Corollary~\ref{cor_theta_even_unimod}, such a theta function is a
modular form on $\SL_2(\ZZ)$, of weight equal to half the dimension of
the space of the quadratic form.

In view of Example~\ref{exam_Leech}, Theorem~\ref{thm_computation_tau}
has the following consequence.
\begin{thm}\label{thm_leech}
There is a deterministic algorithm that on input an integer~$m>0$
together with its factorisation into primes gives the representation
number $r_L(m)$ where $L$ is the Leech lattice, in running time
polynomial in~$\log m$.
\end{thm}
Let us now turn to arbitrary even unimodular lattices. We can then apply
Corollary~\ref{cor_comp_an} to prove the following result.
\begin{thm}\label{thm_even_unimod_latt}
Assume GRH. There is a deterministic algorithm that, on input the
rank~$n_L$ and the integers $r_L(i)$ for $1\leq i\leq n_L/24$ of an
even unimodular lattice $(L,b)$, and an integer $m>0$ together with
its factorisation into primes, computes $r_L(m)$ in running time
polynomial in~$n_L$ and~$\log(m)$.
\end{thm}
\begin{proof}
Let $n_L$ and the $r_L(i)$ for $1\leq i\leq n_L/24$ be given. As $L$
is even and unimodular, $n_L$ is even (in fact, it is a multiple
of~$8$, see~\cite[VII,~\S6]{Serre5}).  Let $k:=n_L/2$. Then $\theta_L$
is in $M_k(\SL_2(\ZZ),\ZZ)$. The coefficients $a_i(\theta_L)=r_L(2i)$
with $0\leq i\leq k/12$ are given to us. The fact that the open balls
in $L_\RR$ centered at the $x\in L$ and with radius~$2^{-1/2}$ do not
overlap gives that for $i$ with $0\leq i\leq k/12$ we have
$\log(1+r_L(i))\leq k\log(1+k/6)$. Corollary~\ref{cor_comp_an} then
means that $r_L(m)=a_{m/2}(\theta_L)$ can be computed in the required
time.
\end{proof}
Let us point out that Theorem~\ref{thm_even_unimod_latt} can be
applied to the orthogonal direct sum~$(L,b)$ of even unimodular
lattices $(L_i,b_i)$ for which the $r_{L_i}(j)$ can be computed in
time polynomial in $n_{L_i}$ and~$j$. Indeed:
\[
\theta_{L_1}\cdots\theta_{L_r} = 
\sum_{x\in L}q^{(b_1(x_1,x_1)+\cdots+b_r(x_r,x_r))/2} = 
\sum_{x\in L}q^{b(x,x)/2} = \theta_L.
\]
The coefficients $a_i(\theta_L)$ with $0\leq i\leq n_L/24$ can be
computed in time polynomial in~$n_L$ by computing the product of the
images of the~$\theta_{L_i}$ in~$\ZZ[[q]]/(q^{n_L/24+1})$.

%% -- Conclusion, prospects
\chapter*{Epilogue}\label{subsec_future}
\addcontentsline{toc}{chapter}{Epilogue}

% authors Bas and Jean-Marc...  

Theorems~\ref{thm_comp_rep_mod_l} and~\ref{thm_comp_tn} will certainly
be generalised to spaces of cuspforms of arbitrarily varying level and
weight.  This is currently being worked out for the probabilistic
variants, and at least for square free levels (and for level two times
a square free number for reasons that will become clear below). We
describe some consequences of this current work. The reader will
understand what the consequences are of deterministic generalisations
to arbitrary levels and weights.

\subsection*{Computation of Galois representations}
Peter Bruin will give, in his forthcoming PhD thesis (Summer 2010), a
probabilistic algorithm that on input positive integers $k$ and $N$,
with $N$ square free, and a surjective morphism
$f\colon\TT(2N,k)\to\FF$ to a finite field, computes $\rho_f$ with
probabilistic running time polynomial in $k$, $N$ and~$\#\FF$. An
important new ingredient in his method, suggested by Couveignes, is to
be more flexible concerning the choice of the divisor $D_0$ of
Section~\ref{sec_setup_tau}. His algorithm computes, for each $x$
in~$V$, the smallest integer~$d_x$ such that
$h^0(X_1(2N)_\Qbar,\calL_x(d_x{\cdot}\infty))=1$, and then represents
$x$ by the unique effective divisor $D_x$ of degree $d_x$ on
$X_1(2N)_\Qbar$ such that $x=[D_x-d_x{\cdot}\infty]$ in
$J_1(2N)(\Qbar)$. This generalises the probabilistic variant of
Theorem~\ref{thm_comp_rep_mod_l} that is described in
Section~\ref{sec_prob_algorithm}. 

\subsection*{Computing coefficients of modular forms}
Concerning Theorem~\ref{thm_comp_tn}, Peter Bruin will obtain a
probabilistic algorithm that, assuming GRH, on input positive integers
$k$, $N$ and~$n$, with $N$ squarefree, together with the factorisation
of $n$ into prime factors, the element $T_n$ of the Hecke algebra
$\TT(2N,k)$ in probabilistic running time polynomial in $k$, $N$
and~$\log n$.

\subsection*{Lattices, sums of squares}
This last result, applied to theta functions of lattices as in
Section~\ref{sec_quad_forms}, will have the following remarkable
consequence that, assuming GRH, there is a probabilistic algorithm
that on input positive integers $k$ and $n$, together with the
factorisation of $n$ into prime factors, computes the number:
\[
r_{\ZZ^{2k}}(n) = \#\{x\in\ZZ^{2k}\;|\;x_1^2+\cdots+x_{2k}^2=n\}
\]
in time polynomial in $k$ and~$\log n$.

For some small values of $k$ there are well known explicit formulas
for the representation numbers $r_{\ZZ^{2k}}(n)$, owing their
existence to the fact that there are no non-CM cuspidal eigenforms on
$\Gamma_1(4)$ of weight~$k$. In order to give these formulas, let
$\chi\colon\ZZ\to\CC$ be the map obtained from the character
$(\ZZ/4\ZZ)^\times=\{1,-1\}\subset\CC^\times$, by extending it by zero
to $\ZZ/4\ZZ$ and composing the result with the reduction map
$\ZZ\to\ZZ/4\ZZ$. Then we have:
\[
\begin{aligned}
r_{\ZZ^2}(n) & = 4\sum_{d|n}\chi(d), \\
r_{\ZZ^4}(n) & =8\sum_{2\nmid d|n}d + 16\sum_{2\nmid d|(n/2)}d, \\
r_{\ZZ^6}(n) & = 16\sum_{d|n}\chi(n/d)d^2 -4\sum_{d|n}\chi(d)d^2, \\
r_{\ZZ^8}(n) & = 16\sum_{d|n}d^3 - 32\sum_{d|(n/2)}d^3 +
256\sum_{d|(n/4)}d^3, \\
r_{\ZZ^{10}}(n) & = (4/5)\sum_{d\mid n}\chi(d)d^4 + 
(64/5)\sum_{d\mid n}\chi(n/d)d^4  \\
& \qquad + (8/5)\sum_{d\in\ZZ[i],\, |d|^2=n}d^4.
\end{aligned}
\]
For the history of these formulas, featuring, among others, Fermat,
Euler, Lagrange, Legendre, Gauss, Jacobi and Liouville, we refer the
reader to~\cite{Milne2} and to Chapter~20 of~\cite{Hardy-Wright}.

In her forthcoming masters thesis (Summer 2010), Ila Varma will show
that there are no other values of~$k$ for which the theta function
$\theta_{\ZZ^{2k}}=\theta_{\ZZ}^{2k}$, where $\ZZ^{2k}$ is equipped
with the standard inner product, is a linear combination of Eisenstein
series and cuspforms coming from Hecke characters.

We conclude that, even in the absence of formulas as above, one will
be able to compute the numbers $r_{\ZZ^{2k}}(n)$ as fast as if such
formulas existed. Often, theta functions are considered to be modular
forms of which the coefficients are easy to compute, and hence useful
for computing Hecke operators and coefficients of eigenforms. However,
it seems that for coefficients $a_n(f)$ with $n$ large, this will be
the other way around, from now on.

\subsection*{Point counting on modular curves}
Another consequence of Peter Bruin's results mentioned above is that,
again assuming GRH, there will be a probabilistic algorithm that on
input a positive square free number~$n$ and a finite field~$\FF_q$
computes the number $\# X_1(n)(\FF_q)$ in time polynomial in~$n$
and~$\log q$. Indeed, this is a matter of computing the element $T_p$
(where $p$ is the prime dividing~$q$) in the Hecke algebra $\TT(n,2)$
acting on the space $S_2(\Gamma_1(n),\ZZ)$.

\subsection*{Point counting in a more general context}
The methods that we have used in this book can also be tried outside
the context of modular forms. Let us consider, for example, a smooth
surface $S$ in~$\PP^3_\QQ$ of degree at least~$5$, say. Letting $l$ be
a prime number, one has the cohomology groups $\rH^i(S_{\Qbar,\et},
\FF_l) $ for $0 \leq i \leq 4$, being finite dimensional
$\FF_l$-vector spaces with $\Gal(\Qbar/\QQ)$-action. It seems
reasonable to suspect that, again, there is an algorithm that on input
a prime $l$ computes these cohomology groups, with their
$\Gal(\Qbar/\QQ)$-action, in time polynomial in~$l$.  Once such an
algorithm is known, one also has an algorithm that, on input a prime
$p$ of good reduction of~$S$, gives the number $\# S(\FF_p)$ of
$\FF_p$-valued points of $S$ in time polynomial in~$\log p$. This
result would be of interest because the known $p$-adic algorithms for
finding such numbers have running time exponential in~$\log p$.

In this case, we choose a Lefschetz fibration from a blow-up of $S$ to
$\PP^1_\QQ$, and use the derived direct images of the constant sheaf
$\FF_{l,S}$ under this fibration to compute the \'etale cohomology
of~$S_\Qbar$ with $\FF_l$-coefficients. The most complicated
contribution then comes from the first derived image $\calF_l$, which
is a locally constant sheaf of $\FF_l$-vector spaces of dimension~$r$,
say, with $r$ independent of~$l$, on the open part $U$ of $\PP^1_\QQ$
over which the fibration is smooth. This open part $U$ is independent
of~$l$ as well, and it is the analog of the open part of the $j$-line
over which all modular curves are unramified.

For each $l$ let $V_l := \ul{\Isom}_U(\FF_l^r, \calF_l)$. These $V_l$
play the role of the covers $X_1(l)$ of the $j$-line, as, by
definition, they trivialise the sheaves~$\calF_l$. Each cover $V_l \to
U$ is finite Galois with group $G=\GL_r(\FF_l)$, and $\rH^1(
U_{\Qbar,\et}, \calF_l) $ is closely related to
$\rH^1(V_{l,\Qbar,\et},\FF_\ell^r)$ which sits in the $l$-torsion of
the Jacobian of the smooth projective model $\ol{V_l}$ of $V_l$. It is
our hope that methods as in this book (height bounds, approximations)
can show that there is a polynomial algorithm for computing these
cohomology groups.

In \cite{Edix-de-Jong-Schepers} a first step in this program is taken,
by proving, in the function field case, that, for varying~$l$, the
height of $\ol{V_l}$ is bounded by a fixed power of~$l$.

\subsection*{Modular forms of half integral weight}
Much to our regret, we have nothing to say about modular forms of half
integer weight. Nevertheless, it would be very interesting to be able
to compute coefficients of such modular forms, as they encode
interesting arithmetic information, such as class numbers.

%\printindex

\end{document}